\documentclass[12pt, a4paper]{report}

\usepackage[utf8]{inputenc}
\usepackage[T1]{fontenc}
\usepackage{amsmath, amsthm, amssymb} % For advanced mathematics
\usepackage{mathtools}
\usepackage{geometry} % For page layout
\usepackage{graphicx} % To include graphics
\usepackage{xcolor}
\usepackage{tikz}
\usepackage{tikz-cd}
\usepackage{setspace}
\usepackage{enumitem}
\usepackage{array}
\usepackage{booktabs}
\usepackage{longtable}
\usepackage{colortbl} % light rules between the rows of tables
\usepackage{caption}
\usepackage{listings} % Code listings (Chapter~\ref{chap:implementation}, Appendix~\ref{chap:code})
\usepackage{needspace} % keeps a file heading in Appendix~\ref{chap:code} with its first lines
\usepackage{sectsty}
\allsectionsfont{\raggedright}
\usepackage[numbers,sort&compress]{natbib}
\usepackage{hyperref} % For hyperlinks and metadata (loaded late)

\usetikzlibrary{positioning, arrows.meta, calc, decorations.markings}

\newenvironment{itemizedot}{\begin{itemize} }{\end{itemize}}

\hypersetup{
    colorlinks=true,
    linkcolor=blue,
    citecolor=blue!60!black,
    filecolor=magenta,
    urlcolor=cyan!50!black,
    pdftitle={Automating the Derivation of Advanced Algorithms in Numerical Linear Algebra using Associative Algebras with Involution},
    pdfauthor={Ran Gutin},
    pdfsubject={PhD Thesis},
    pdfkeywords={Numerical Linear Algebra, Star-Algebras, Quivers, Path Algebras, Automatic Differentiation},
    bookmarks=true,
}

\theoremstyle{plain}
\newtheorem{theorem}{Theorem}[chapter]
\newtheorem{lemma}[theorem]{Lemma}
\newtheorem{proposition}[theorem]{Proposition}
\newtheorem{corollary}[theorem]{Corollary}
\theoremstyle{definition}
\newtheorem{definition}[theorem]{Definition}
\newtheorem{example}[theorem]{Example}
\newtheorem{convention}[theorem]{Convention}
\theoremstyle{remark}
\newtheorem{remark}[theorem]{Remark}
\theoremstyle{plain}
\newtheorem{mainthm}{Theorem}

\newcommand{\staralg}[0]{$*$-algebra}
\newcommand{\staralgs}[0]{$*$-algebras}
\newcommand{\starring}[0]{$*$-ring}
\newcommand{\rad}[0]{\operatorname{rad}} % Jacobson radical

\newcommand{\epsilondiv}{\epsilon_{\mathrm{div}}}
\newcommand{\epsilonform}{\epsilon_{\mathrm{form}}}
\newcommand{\RR}{\mathbb{R}}
\newcommand{\CC}{\mathbb{C}}
\newcommand{\HH}{\mathbb{H}}
\newcommand{\QQ}{\mathbb{Q}}
\newcommand{\ZZ}{\mathbb{Z}}

\newcommand{\KK}{\mathbb{K}}

\newcommand{\RP}{\mathbb{RP}}
\newcommand{\PSL}{\mathrm{PSL}}
\newcommand{\SL}{\mathrm{SL}}
\newcommand{\GL}{\mathrm{GL}}
\newcommand{\PU}{\mathrm{PU}}
\newcommand{\SO}{\mathrm{SO}}
\newcommand{\Or}{\mathrm{O}}
\newcommand{\Un}{\mathrm{U}}
\newcommand{\Sp}{\mathrm{Sp}}
\newcommand{\tr}{\operatorname{tr}}
\newcommand{\Arg}{\operatorname{Arg}}
\newcommand{\sign}{\operatorname{sgn}}
\newcommand{\st}{\operatorname{st}}

\newcommand{\LieQuad}{\mathrm{L}}
\newcommand{\bmat}[1]{\begin{pmatrix}#1\end{pmatrix}}
\newcommand{\Herm}{\operatorname{Herm}}
\newcommand{\twoR}{{}^{2}\mathbb{R}}

\DeclareMathOperator{\Hom}{Hom}
\DeclareMathOperator{\End}{End}

\DeclareMathOperator{\Rep}{Rep}
\DeclareMathOperator{\diag}{diag}



\colorlet{codegray}{black!50} % a named colour: `black!50' is misread inside `postbreak'
\lstdefinestyle{dsqcode}{
  basicstyle=\ttfamily\scriptsize\setstretch{1},
  frame=none,
  xleftmargin=0pt,
  xrightmargin=0pt,
  aboveskip=\smallskipamount,
  belowskip=\smallskipamount,
  numbers=left,
  numberstyle=\tiny\color{codegray},
  numbersep=8pt,
  stepnumber=5,
  firstnumber=1,
  breaklines=true,
  postbreak=\mbox{\textcolor{codegray}{$\hookrightarrow$}\space},
  extendedchars=true,
  literate={·}{{$\cdot$}}1 {⁻}{{$^{-}$}}1 {¹}{{$^{1}$}}1 {…}{{\ldots}}1,
}
\lstdefinestyle{dsqtext}{style=dsqcode, language={}, numbers=none, xleftmargin=\parindent, aboveskip=\medskipamount, belowskip=\medskipamount}

\tikzset{
  qvertex/.style={circle, draw, thick, inner sep=1pt, minimum size=26pt, font=\small},
  qarr/.style={-{Latex[length=2.5mm]}, thick},
  qinv/.style={<->, densely dashed, gray!70!black},
  qbar/.style={<->, densely dotted, thick, blue!60!black},
}

\newcommand{\rowrule}{\arrayrulecolor{black!30}\midrule[0.3pt]\arrayrulecolor{black}}

\begin{document}
%======================================================================
%   TITLE PAGE
%======================================================================
\pagenumbering{arabic}
\begin{titlepage}
\centering
\vspace*{1.5cm}
{\LARGE\bfseries Automating the Derivation of Advanced Algorithms in Numerical Linear Algebra using Associative Algebras with Involution\par}
\vspace{2.5cm}
{\Large Ran Gutin\par}
\vspace{2.5cm}
{\large Department of Computing\\[1ex] Imperial College London\par}
\vspace{2.5cm}
{\large A thesis submitted for the degree of\\[1ex] Doctor of Philosophy (PhD)\par}
\vfill
{Submitted September 2025\\[1ex] Revised following examination, October 2026\par}
\end{titlepage}
\setcounter{page}{2} % continuous Arabic numbering: the title page is page 1

%======================================================================
%   ABSTRACT
%======================================================================
\chapter*{Abstract}
\addcontentsline{toc}{chapter}{Abstract}
Standard algorithms of numerical linear algebra are written against seven operations: addition, subtraction, multiplication, division, the adjoint and the $*$-square root of scalars, and $(i,j)$ indexing of matrices. This thesis overloads them, as forward-mode automatic differentiation overloads arithmetic, so that the scalars may lie in a \emph{$*$-algebra}, an associative algebra with an involution. One program then becomes different classical algorithms over different $*$-algebras, and several classical decompositions become one: the Hermitian eigendecomposition, a tridiagonal, block-diagonal canonical form under unitary similarity.

Over a $*$-algebra, indexing must say what the entries of a matrix are. We prove an analogue of Gabriel's theorem: every finite-dimensional real $*$-algebra consists of matrices over a quotient of the path $*$-algebra of its \emph{doubly-signed quiver}. This is a directed graph whose vertices record whether their scalars are real, complex or quaternionic, and how the involution acts. Each doubly-signed quiver encodes a classical problem as unitary similarity of Hermitian matrices. The spectral theorems, the singular value decomposition, the Takagi factorisation, Sylvester's law of inertia and the real Jordan form arise in this way.

Division and $*$-square roots can fail in the $*$-algebra itself, so they take place in a completed localisation, which inverts the arrows. For the real singular value decomposition, Householder tridiagonalisation becomes Golub--Kahan bidiagonalisation. Over the split-complex numbers (over which the Hermitian eigendecomposition is the real Jordan form), the QR algorithm becomes Rutishauser's LR algorithm, two steps at a time. There the Cholesky factorisation becomes Doolittle's $LU$ factorisation. These hold for suitable $*$-square roots, whenever the steps are defined.

Infinitesimals make limits of algorithms exact computations: the Gram--Schmidt process yields the Bruhat decomposition, and Gaussian elimination, with pivots chosen by orders of magnitude, yields Dijkstra's shortest-path algorithm to leading order. A Python library implements the seven operations and these algorithms.

\cleardoublepage

%======================================================================
%   PROVENANCE AND ACKNOWLEDGEMENTS
%======================================================================
\chapter*{Declarations}
\addcontentsline{toc}{chapter}{Declarations}
\section*{Originality}
I declare that the work in this thesis is my own, and that the work of others is acknowledged and referenced where it is used. I had substantial help from AI assistants. Gemini, made by Google, helped me to understand large amounts of background material. Claude, made by Anthropic, and ChatGPT, made by OpenAI, helped me to carry out experiments, to write and test software, and to edit and revise the text much faster than I could have done alone. I take full responsibility for the content of this thesis.

\section*{Copyright}
The copyright of this thesis rests with the author. Unless otherwise indicated, its contents are licensed under a Creative Commons Attribution 4.0 International Licence (CC~BY 4.0). The licence is at \url{https://creativecommons.org/licenses/by/4.0/}.

\chapter*{Statement of Provenance}
\addcontentsline{toc}{chapter}{Statement of Provenance}
Parts of this thesis build on papers and preprints written by the author during his doctoral studies. Their relationship to the thesis is as follows.
\begin{itemize}
  \item Theorem~\ref{thm:spectral-power-series} and Corollaries~\ref{cor:svd-power-series} and~\ref{cor:twisted-hermitian} in Chapter~\ref{chap:infinitesimals} generalise, from the dual numbers to formal power series, results of~\cite{Gutin2022DualNumbers}: its spectral theorems for symmetric and for Hermitian matrices, and its T-SVD $M=U\Sigma V^T$ with $U^TU=V^TV=I$. Its $*$-SVD is not treated in this thesis.
  \item Section~\ref{subsec:2R} and Proposition~\ref{prop:2R-algorithms}, on the split-complex numbers $\twoR$, develop the preprint~\cite{Gutin2021DoubleNumbers}, which observed, under the name of double numbers, that over this algebra the $LDL^*$ decomposition yields the $LU$ decomposition of a real matrix and an algorithm for the singular value decomposition becomes the LR algorithm. Proposition~\ref{prop:2R-algorithms} makes these statements precise, with the choice of $*$-square roots and the conditions under which the steps are defined.
  \item The parity-equivariant spectral theorem of Section~\ref{sec:parity} makes precise a unifying observation made in~\cite{Gutin2024Clifford} about canonical forms over Clifford algebras.
  \item Appendix~\ref{chap:gutin-2022} is based on the preprint~\cite{Gutin2022Visualisation}. Its pictures are those of the preprint, redrawn; their derivation from the space of inner products on $\RR^2$, and Section~\ref{sec:qr-lr} on the LR algorithm, are new to this thesis. Its warm-up visualisation for positive semidefinite matrices first appeared in~\cite{Gutin2022PSD}. The appendix is not required for the remainder of the thesis.
\end{itemize}
Further related work by the author is~\cite{Gutin2021AllRings,Gutin2023DoubleComplex,Gutin2023Conjecture}. Results that are not new are attributed to their sources at the point where they are stated.

\chapter*{Acknowledgements}
\addcontentsline{toc}{chapter}{Acknowledgements}
I am grateful to my supervisor, Professor Abbas Edalat, for his guidance and support throughout this work. I thank my examiners, L. Feng and H. Wiklicky, for their careful reading of the thesis and for detailed comments that improved it. I also thank Gregory Gutin for helpful discussions.

\cleardoublepage

%======================================================================
%   TABLE OF CONTENTS
%======================================================================
\tableofcontents
\cleardoublepage

\chapter*{List of Notation}
\addcontentsline{toc}{chapter}{List of Notation}
The main symbols, with the place where they are defined. Where a letter has a second meaning that occurs in several places, both are listed. Letters are also used locally, within a proof or an example, with meanings that are introduced there.

\begin{singlespace}
{\small
\renewcommand{\arraystretch}{1.3}
\begin{longtable}{@{}>{\raggedright\arraybackslash}p{0.2\linewidth}>{\raggedright\arraybackslash}p{0.56\linewidth}>{\raggedright\arraybackslash}p{0.18\linewidth}@{}}
\toprule
\textbf{Symbol} & \textbf{Meaning} & \textbf{Defined in} \\
\midrule
\endhead
\multicolumn{3}{@{}l}{\emph{Algebras and involutions}} \\
$x^*$, $X^*$ & the adjoint of $x$, given by the involution of a $*$-algebra; for matrices over it, the conjugate transpose & Def.~\ref{def:starring}, Conv.~\ref{conv:notation} \\
$X^T$, $X^H$ & transpose; conjugate transpose of a real or complex matrix & Conv.~\ref{conv:notation} \\
$\Herm_n(A)$, $\Un_n(A)$ & Hermitian and unitary $n\times n$ matrices over $A$ & Conv.~\ref{conv:notation} \\
$\twoR$ & the split-complex numbers $\RR\oplus\RR$, with $(x,y)^*=(y,x)$ & Conv.~\ref{conv:notation}, Sec.~\ref{subsec:2R} \\
$X^\omega$, $J$ & the symplectic involution $J^{-1}X^TJ$, with $J=\bigl(\begin{smallmatrix}0&I_n\\-I_n&0\end{smallmatrix}\bigr)$; in Section~\ref{sec:parity}, $J$ is $\diag(I_p,-I_q)$ & Ex.~\ref{ex:staralgs}, Sec.~\ref{sec:hermitian-sqrt} \\
$\varepsilon$ & a positive infinitesimal; the generator of the dual numbers $\RR[\varepsilon]/(\varepsilon^2)$ & Ch.~\ref{chap:infinitesimals}, Ex.~\ref{ex:staralgs} \\
$\epsilon$ & the sign $\pm1$ of an $\epsilon$-Hermitian form & Def.~\ref{def:hermitian} \\
$\rad(A)$ & the Jacobson radical & Sec.~\ref{subsec:radicals} \\
$S$, $R$ & a $*$-stable semisimple subalgebra and the radical, in $A=S\oplus R$ with $R=\rad(A)$; in a factorisation $H=QR$, $R$ is the triangular factor & Thm.~\ref{thm:wedderburn-malcev}, Lemma~\ref{lem:star-wedderburn} \\
$e_v$, $1_c$ & the idempotent of a vertex $v$; of an orbit $c$ of vertices & Def.~\ref{def:path-algebra}, Def.~\ref{def:matrix-like} \\
$c_k$ & a position: the $k$-th index carrying the idempotent $1_c$ & Sec.~\ref{sec:indexing}, Cor.~\ref{cor:positions} \\
\addlinespace
\multicolumn{3}{@{}l}{\emph{Quivers}} \\
$Q=(Q_0,Q_1,s,t)$ & a quiver; in a factorisation $H=QR$, $Q$ is the unitary factor & Def.~\ref{def:quiver} \\
$kQ$, $R_Q$ & the path algebra and its arrow ideal & Def.~\ref{def:path-algebra} \\
$r$ & the cycle rank of a quiver & Sec.~\ref{sec:tits_form} \\
$\dagger$; $\epsilon(v)$, $\epsilon(\alpha)$ & the involution and the signs of a signed quiver & Def.~\ref{def:signed-quiver} \\
$v^{\top}$, $\alpha^{\top}$; $\bar v$, $\bar\alpha$ & the transpose and the conjugation of a doubly-signed quiver & Def.~\ref{def:doubly-signed} \\
$\epsilondiv$, $\epsilonform$ & the two signs of a vertex of a doubly-signed quiver & Def.~\ref{def:doubly-signed}, Table~\ref{tab:vertex-types} \\
$(D,\epsilonform)$, $(\CC,\mathrm u)$ & labels of the ten vertex types & Sec.~\ref{sec:dsq} \\
$G_v$, $C_v$; $\kappa_\alpha$, $\mu_\alpha$ & the matrices and the signs of the block model & Table~\ref{tab:vertex-blocks}, \eqref{eq:sign-conditions} \\
$\iota$ & the conjugation of a complexified $*$-algebra, whose fixed points form the real algebra & Thm.~\ref{thm:classification}, Def.~\ref{def:path-star-algebra} \\
$\RR Q$, $S_Q$, $M_Q$, $B_Q$ & the path $*$-algebra, its parts of degree $0$ and $1$, and $\RR Q/\RR Q_{\geq2}$ & Def.~\ref{def:path-star-algebra} \\
$\mathcal O(Q)$ & the orbits of the vertices under the transpose and the conjugation & Sec.~\ref{sec:matrices-with-signs} \\
$M_{n,\eta}(B)$; $\eta_c$, $n_c$ & matrices with signs; the form and the number of positions at an orbit $c$ & Sec.~\ref{sec:matrices-with-signs}, Def.~\ref{def:matrix-like} \\
$E$, $E_0$, $E_1$ & the matrix algebra of the encoding theorem and its parts of degree $0$ and $1$ & Thm.~\ref{thm:encoding} \\
\addlinespace
\multicolumn{3}{@{}l}{\emph{Localisation and completion}} \\
$kQ_{\mathrm{loc}}$ & the localisation of $kQ$ at its arrows & Def.~\ref{def:localisation} \\
$\Lambda$ & the completed localisation; in Section~\ref{subsec:cat-symmetric}, a diagonal matrix of eigenvalues & Def.~\ref{def:completion} \\
$v(x)$, $x^\circ$ & the valuation and the lowest-degree part of an element of $\Lambda$ & Def.~\ref{def:completion} \\
$k\{X\}$, $D_r$ & the free field; the coefficient skew field of the embedding theorem & Sec.~\ref{sec:embedding} \\
\addlinespace
\multicolumn{3}{@{}l}{\emph{Infinitesimals}} \\
$\RR\langle\varepsilon\rangle$ & the real closed field of algebraic Puiseux series & Sec.~\ref{sec:puiseux} \\
$k((t))$, $k[[t]]$ & formal Laurent series and formal power series & Sec.~\ref{sec:puiseux} \\
$\st(x)$ & the standard part of a finite element & Sec.~\ref{sec:standard-part} \\
${}^*\RR$ & the hyperreal numbers of nonstandard analysis & Sec.~\ref{subsec:robinson} \\
\addlinespace
\multicolumn{3}{@{}l}{\emph{Groups}} \\
$\Or(n)$, $\Un(n)$, $\Sp(n)$, $\Sp(2n,\RR)$, $\Or^*(2n)$, $\Or(n,\CC)$, $\Sp(2n,\CC)$, $\GL_n(D)$ & the unitary groups of the ten vertex types & Table~\ref{tab:vertex-types} \\
\bottomrule
\end{longtable}
}
\end{singlespace}
\cleardoublepage
% (continuous Arabic numbering from the title page)
%======================================================================
%   CHAPTER 1: INTRODUCTION
%======================================================================
\chapter{Introduction}
\label{chap:intro}

The standard algorithms of numerical linear algebra, such as Householder reduction, the QR algorithm and the Cholesky and $LU$ factorisations, are written against seven operations: addition, subtraction, multiplication, division, the adjoint $x\mapsto x^*$ and the $*$-square root of scalars, and $(i,j)$ indexing of matrices. This thesis overloads these seven operations, as forward-mode automatic differentiation overloads arithmetic~\cite{GriewankWalther2008}, so that the scalars may come from an associative algebra with an involution, called a \emph{\staralg} (Definition~\ref{def:staralg}), instead of from $\RR$ or $\CC$; the involution is the adjoint.

The outcome is that matrix decompositions are unified, and so are several of the algorithms that compute them. The spectral theorems, the singular value decompositions (SVDs) over $\RR$, $\CC$ and $\HH$, the Takagi factorisation, Sylvester's law of inertia and the real Jordan form are instances of one problem, unitary similarity of Hermitian matrices, over different \staralgs{} (Theorem~\ref{thm:A}). And an algorithm written once, which we call a \emph{generic algorithm}, becomes different classical algorithms over different \staralgs{}: Householder tridiagonalisation becomes the bidiagonalisation of Golub and Kahan, the QR algorithm becomes forms of Rutishauser's LR algorithm, and the Cholesky factorisation becomes Doolittle's $LU$ factorisation (Theorem~\ref{thm:B}). This is the sense of the title: each of these algorithms is derived from another by running one program over different scalars, with no step done by hand.

Most of the thesis is therefore about one problem and the algorithms that solve it: the eigendecomposition of a Hermitian matrix, computed by Householder tridiagonalisation and the QR algorithm. Here a Hermitian eigendecomposition is a block-diagonal form under unitary similarity whose blocks cannot be split further (compare Theorem~\ref{thm:hermitian-ks}). Over $\RR$, $\CC$ and $\HH$ the blocks are $1\times1$ and the form is diagonal; over other $*$-algebras they can be larger, as for the real Jordan form over the split-complex numbers $\twoR=\RR\oplus\RR$ (Section~\ref{subsec:2R}). For the decompositions of Theorem~\ref{thm:A}(3) the blocks have normal forms and can be taken tridiagonal (Section~\ref{sec:catalogue}), so that the Hermitian eigendecomposition is a tridiagonal, block-diagonal canonical form.

The other factorisations follow from this one by limits (Section~\ref{sec:chain}). The QR factorisation, as distinct from the QR algorithm, is a limit of the singular value decomposition (Proposition~\ref{prop:qr-from-svd}). The Bruhat decomposition of Theorem~\ref{thm:B}(4) and the $LU$ factorisation are limits of the QR factorisation, computed over a field with an infinitesimal (\emph{auto-limits}; Theorem~\ref{thm:bruhat} and Proposition~\ref{prop:gs-limit}). Dijkstra's shortest-path algorithm is the leading-order part of Gaussian elimination, with the pivots chosen by orders of magnitude (Proposition~\ref{prop:dijkstra}). What does not reduce to the Hermitian eigendecomposition is the machinery that sets up the reduction: the structure theorem for real $*$-algebras (Chapter~\ref{chap:quivers}), division and $*$-square roots in the completed localisation (Chapter~\ref{division-problem}) and the transfer principle (Chapter~\ref{chap:infinitesimals}).

Addition, subtraction, multiplication and the adjoint carry over to every \staralg{}; the algorithms use the adjoint to form conjugate transposes. The other three operations do not carry over directly, and giving them a meaning is the technical content of the thesis. Indexing must say what the entries of a matrix over a \staralg{} are. By a theorem of Gabriel, every basic finite-dimensional algebra over an algebraically closed field is a quotient of the path algebra of its quiver~\cite[Theorem~II.3.7]{AssemSimsonSkowronski2006}; for finite-dimensional real \staralgs{}, an analogue of this theorem for \emph{doubly-signed quivers} answers the question (Chapter~\ref{chap:quivers}). Division and the $*$-square root can fail in the \staralg{} itself, because generic algorithms divide by zero divisors and normalise vectors whose norm vanishes; we carry them out in a \emph{completed localisation} of the path $*$-algebra, in which the arrows of the quiver are inverted (Chapter~\ref{division-problem}); for short, the completed localisation of the quiver. The simplest instance is the field of Laurent series, which is the completed localisation of the quiver with one loop. Its variable $\varepsilon$ is an infinitesimal. Running an algorithm over the field of algebraic Puiseux series in $\varepsilon$ and keeping the constant terms of its finite outputs gives the limits of those outputs as a real parameter $\varepsilon$ tends to zero (Proposition~\ref{prop:autolimit}).

Section~\ref{sec:proposal} lists the seven operations and what each needs over a \staralg{}, and Section~\ref{sec:outcomes} states the two outcomes as Theorems~\ref{thm:A} and~\ref{thm:B}. Section~\ref{sec:motivation} describes the two examples from which the approach grew. Section~\ref{sec:limitations} explains why existing unifying frameworks do not provide transformations of programs, and Section~\ref{sec:related} surveys related work. Section~\ref{sec:contributions} describes the contributions chapter by chapter.

\section{Seven Overloaded Operations}
\label{sec:proposal}
An algorithm of numerical linear algebra, written for real or complex matrices, uses the following operations, six on its scalars and one on its matrices:
\begin{enumerate}[label=(O\arabic*),nosep]
  \item addition, $x+y$;
  \item subtraction, $x-y$;
  \item multiplication, $xy$;
  \item division, $x/y=xy^{-1}$;
  \item the \emph{adjoint}, $x\mapsto x^*$;
  \item the \emph{$*$-square root}: for a self-adjoint $s$, some $x$ with $x^*x=s$;
  \item $(i,j)$ indexing: access to the entry $X_{ij}$ of a matrix $X$.
\end{enumerate}
The adjoint is complex conjugation over $\CC$ and the identity over $\RR$, which is why it is rarely counted as an operation; the algorithms use it to form conjugate transposes $X^*$ and inner products $x^*y=\sum_ix_i^*y_i$. Over $\RR$ and $\CC$, the $*$-square roots of $s\geq0$ are the numbers $\sqrt s$ multiplied by a sign or a phase. An algorithm may also compare the sizes of scalars, where it chooses a sign or a pivot. Comparison is not counted among the seven operations: it never produces a scalar, and only selects which admissible $*$-square root or pivot the other operations then use. Listing~\ref{lst:householder} shows Householder tridiagonalisation written against these operations.

\begin{lstlisting}[float=tbp,label={lst:householder},caption={Householder tridiagonalisation written against the seven operations, with one comparison to choose the sign of the root, in Python with NumPy-style slicing; \texttt{X.H} is the conjugate transpose, formed with the adjoint. It simplifies the function \texttt{householder} of Appendix~\ref{chap:code}, which also reports breakdowns and skips columns that are already reduced.}]
def householder(H, sqrt):
    """Tridiagonalise a Hermitian H: return T, U with H = U T U*."""
    n = H.shape[0]
    T, U = H.copy(), identity(n)
    for k in range(n - 2):
        x = T[k+1:, k]                   # (i, j) indexing
        gamma = sqrt(x.H @ x)            # *-square root: gamma* gamma = x* x
        if norm(x[0] + gamma) > norm(x[0] - gamma):
            gamma = -gamma               # comparison: avoid cancellation
        w = x.copy()
        w[0] = x[0] - gamma              # subtraction
        s = w.H @ x                      # adjoint, multiplication and addition
        P = identity(n)
        P[k+1:, k+1:] -= (w / s) @ w.H   # division, in the corner of column k
        T = P @ T @ P.H
        U = U @ P.H
    return T, U
\end{lstlisting}

To \emph{overload} the operations is to give each of them a meaning for scalars from a \staralg{} $A$, so that the same program runs on matrices with entries in $A$. In forward-mode automatic differentiation, $A$ is the algebra of dual numbers $\RR[\varepsilon]/(\varepsilon^2)$ with the trivial involution, and the program computes a derivative along with each value. A square root or an inverse that the listing needs may not exist over $A$; Proposition~\ref{prop:householder-any-root} says when each step can be carried out. When every step can, the listing computes a tridiagonal $T$ and a unitary $U$ over $A$ with $H=UTU^*$, because the proof that the reflectors $P=I-ws^{-1}w^*$ are unitary uses only identities that hold in every \staralg{}. Table~\ref{tab:operations} summarises what each operation means over $A$, what makes it difficult, and where the difficulty is addressed.

\begin{table}[htbp]
\centering
\caption{The seven operations over a $*$-algebra $A$: what they mean, what makes them difficult, and where the difficulty is addressed.}
\label{tab:operations}
\small
\renewcommand{\arraystretch}{1.25}
\begin{tabular}{@{}>{\raggedright\arraybackslash}p{0.13\linewidth}>{\raggedright\arraybackslash}p{0.25\linewidth}>{\raggedright\arraybackslash}p{0.25\linewidth}>{\raggedright\arraybackslash}p{0.29\linewidth}@{}}
\toprule
\textbf{Operation} & \textbf{Meaning over $A$} & \textbf{Difficulty} & \textbf{Addressed by} \\
\midrule
$+$, $-$, $\times$ & the operations of $A$ & representing the elements of $A$ & paths with coefficient blocks (Theorem~\ref{thm:structure}, Proposition~\ref{prop:block-model}) \\
\rowrule
$/$ & $x/y=xy^{-1}$, with $y$ invertible between idempotents (Definition~\ref{def:invertible-between}) & zero divisors: no non-zero element of the radical is invertible & the completed localisation $\Lambda$ (Chapter~\ref{division-problem}) \\
\rowrule
adjoint & $x\mapsto x^*$, the involution of $A$ & it is not the identity: over $\twoR$ it exchanges the components, and in a path $*$-algebra it reverses paths & the block model \eqref{eq:block-star}; on matrices, the conjugate transpose \eqref{eq:signed-adjoint} \\
\rowrule
$*$-square root & for self-adjoint $s$, some $x$ with $x^*x=s$ & not unique, and need not exist & any root in the corner of the next position with $w^*x$ invertible (Proposition~\ref{prop:householder-any-root}); roots in $\Lambda$ (Section~\ref{sec:lambda-sqrt}) \\
\rowrule
$(i,j)$ indexing & entries $X_{c_k,d_l}\in1_cA1_d$ at positions $c_k$, $d_l$ & which idempotents, and in which order & the orbits of the vertices of a doubly-signed quiver (Corollary~\ref{cor:positions}); a cyclic order (Proposition~\ref{prop:golub-kahan}) \\
\bottomrule
\end{tabular}
\end{table}

\paragraph{Addition, subtraction and multiplication.} These are the operations of $A$. Its elements are represented as matrices over a quotient of the path $*$-algebra of a doubly-signed quiver, that is, over linear combinations of paths (Theorem~\ref{thm:structure} and Proposition~\ref{prop:block-model}); this representation also serves indexing and division.

\paragraph{The adjoint.} The adjoint $x\mapsto x^*$ is the involution of $A$: it reverses products, $(xy)^*=y^*x^*$, and its square is the identity. On matrices, the adjoints of the entries, with the indices exchanged, give the conjugate transpose, which algorithms for Hermitian matrices use at every step.

\paragraph{Indexing.} Over $A$, the entries of a matrix are refined into elements of corners $1_cA1_d$ cut out by self-adjoint idempotents, and a matrix is indexed by \emph{positions}, pairs $c_k$ of an idempotent $1_c$ and a counter $k$ (Section~\ref{sec:indexing}). Because the idempotents are self-adjoint, the adjoint acts on matrices as the conjugate transpose, $(X^*)_{ij}=(X_{ji})^*$, up to signs for indefinite forms (Remark~\ref{rem:indefinite}); this is what an algorithm written for Hermitian matrices needs. Which idempotents to use is answered by the analogue of Gabriel's theorem: every finite-dimensional real \staralg{} is a matrix algebra over a quotient of the path $*$-algebra of a doubly-signed quiver, and the positions are grouped by the orbits of its vertices under the transpose and the conjugation (Corollary~\ref{cor:positions}). The positions are listed in a cyclic order through the orbits. The order matters: for the quiver of the singular value decomposition, it is what turns tridiagonalisation into bidiagonalisation (Proposition~\ref{prop:golub-kahan}).

\paragraph{Division.} Every finite-dimensional real \staralg{} other than $\RR$, $\CC$ and $\HH$ contains zero divisors, and the generic algorithms divide by them, even when division is taken between idempotents (Definition~\ref{def:invertible-between}). For example, over the four-dimensional $*$-algebra of the singular value decomposition (Section~\ref{sec:unification-example}), whose radical squares to zero, the norms $x^*x$ of the columns that Householder's algorithm reduces vanish. In the path $*$-algebra of its quiver they are multiples $c\,ba$ or $c\,ab$ of paths of length two, which have the $*$-square roots $\pm\sqrt c\,a$ or $\pm\sqrt c\,b$ but no inverses. We therefore run the algorithms in the \emph{completed localisation} $\Lambda$ of this path $*$-algebra, obtained by inverting every arrow between its end points and completing in the path degree (Chapter~\ref{division-problem}). In $\Lambda$ every path is invertible between its end points. For all the classical examples, each corner of $\Lambda$ is, after complexification, a space of matrices over a field, and elimination inverts every invertible element (Proposition~\ref{prop:lambda-complete}). A result holds in the original algebra whenever the computed factors, such as $U$ and $T$ in $H=UTU^*$, contain no inverted arrows (Theorem~\ref{thm:reduction}).

\paragraph{The $*$-square root and comparisons.} A $*$-square root is not unique (Lemma~\ref{lem:sqrt-unique}). Householder reflectors work with any $*$-square root that lies in the corner of the next position and for which the denominator of the reflector is invertible (Proposition~\ref{prop:householder-any-root}), so the phase condition of the classical algorithm is not needed, and in $\Lambda$ the roots are computed degree by degree (Section~\ref{sec:lambda-sqrt}). Comparisons, such as the choice of sign in Listing~\ref{lst:householder}, need a measure of the size of a scalar. Over a real closed field with an infinitesimal, they are comparisons in an ordered field, and such fields also supply square roots (Chapter~\ref{chap:infinitesimals}); in $\Lambda$, the implementation compares norms of coefficients and, to choose pivots, degrees and leading coefficients (Section~\ref{sec:lambda-inverses}).

\section{Two Outcomes: Unified Decompositions and Unified Algorithms}
\label{sec:outcomes}
We now state the two outcomes, with references to the precise results. A quiver is a directed graph. A \emph{doubly-signed quiver} $Q$ (Definition~\ref{def:doubly-signed}) is a quiver with two involutions, a \emph{transpose} that reverses arrows and a \emph{conjugation} that preserves them, together with signs; at each vertex they record whether the scalars are real, complex or quaternionic, and how the involution acts there. Its path $*$-algebra $\RR Q$ is graded by the length of paths: the vertices have degree $0$ and the arrows degree $1$.

\begin{mainthm}[Unified decompositions]\label{thm:A}
\leavevmode
\begin{enumerate}
  \item Every finite-dimensional real $*$-algebra is $*$-isomorphic to an algebra of block matrices over a quotient $\RR Q/I$ of the path $*$-algebra of a doubly-signed quiver $Q$, with one block of rows and columns for each orbit of vertices and a symmetric form on each block that enters the adjoint. The quiver $Q$, the sizes of the blocks and the signatures of the forms are invariants of the $*$-algebra (Theorem~\ref{thm:structure}).
  \item Let $Q$ be a doubly-signed quiver, let $I\subseteq\RR Q$ be an ideal spanned by homogeneous elements of degree at least $2$ and stable under the involution, and consider matrices over $\RR Q/I$ indexed by positions. Two Hermitian matrices of degree $1$, that is, with entries that are linear combinations of arrows, are unitarily similar if and only if they are related by a unitary matrix of degree $0$. The unitary matrices of degree $0$ form a product of classical groups, one for each orbit of vertices, and the Hermitian matrices of degree $1$ form a direct sum of classical matrix spaces, one for each orbit of arrows (Theorem~\ref{thm:encoding}).
  \item For the quivers of Table~\ref{tab:catalogue}, the problem of (2) is the one solved by the spectral theorem for real symmetric, real skew-symmetric, complex Hermitian or quaternionic Hermitian matrices; by the singular value decomposition over $\RR$, $\CC$ or $\HH$; by the Takagi factorisation or the Youla normal form; by Sylvester's law of inertia; or by the real Jordan form (Section~\ref{sec:catalogue}).
\end{enumerate}
\end{mainthm}

Part~(1) says how matrices over a \staralg{} are indexed: by positions, grouped by the orbits of the vertices of $Q$. Part~(2) says which problem $Q$ encodes: a product of classical groups acting on a direct sum of classical matrix spaces. Part~(3) says that the familiar decompositions are normal forms for some of these problems. So a doubly-signed quiver encodes a problem, and a decomposition is a normal form for it. Section~\ref{subsec:recipe} reads part~(2) backwards, as a recipe for encoding a new problem. Theorem~\ref{thm:A} does not say that every quiver has a normal form; that depends on its representation type (Section~\ref{sec:wild-behaviour}). Part~(1) rests on Theorem~\ref{thm:classification}; parts~(2) and~(3) use only the explicit construction of the path $*$-algebra in Section~\ref{sec:path-star-algebras}.

The second outcome concerns the following programs, each written once against the seven operations: Householder tridiagonalisation (Listing~\ref{lst:householder}); the Gram--Schmidt process, and the QR algorithm, whose step $H=QR\mapsto RQ$ computes $Q$ and $R$ by it; and the Cholesky factorisation $H=R^*R$. They are run over the completed localisation $\Lambda$, with the graded $*$-square roots of Section~\ref{sec:lambda-sqrt}. Part~(3) is over the split-complex numbers $\twoR$. They are the path $*$-algebra of a doubly-signed quiver without arrows, whose two vertices are real and exchanged by the transpose (a single vertex orbit of type $(\RR,0)$), and their Hermitian matrices $(X,X^T)$, of degree~$0$, encode similarity of real matrices $X$ (Section~\ref{subsec:2R}).

\begin{mainthm}[Unified algorithms]\label{thm:B}
\leavevmode
\begin{enumerate}
  \item Householder tridiagonalisation behaves as follows.
  \begin{enumerate}
    \item Over the doubly-signed quiver of the real singular value decomposition, with the positions alternating between its two vertices, $u_1,v_1,u_2,v_2,\dots$, and for the matrix that encodes a real $m\times n$ matrix $W$, $m\geq n$, it is defined at every step, and it is the bidiagonalisation of $W$ by Golub and Kahan, up to signs (Proposition~\ref{prop:golub-kahan}).
    \item If every component of the complexified quiver is an oriented cycle or a linearly oriented path, then, whenever its steps are defined, it transforms a Hermitian matrix of degree~$1$, with all signs at the positions positive, only by the classical groups of Theorem~\ref{thm:A}(2), and so reduces it to tridiagonal form within the encoded problem (Proposition~\ref{prop:householder-graded}).
  \end{enumerate}
  \item Let $H$ encode an invertible square matrix $W$ over the quiver of the singular value decomposition over $\RR$, $\CC$ or $\HH$. List the positions alternately from its two orbits of vertices, and use principal graded $*$-square roots. Then two steps of the QR algorithm on $H$ are one step of the unshifted QR algorithm on each of the Gram matrices $W^*W$ and $WW^*$, and a single step is one step of the Cholesky LR algorithm on each of them (Proposition~\ref{prop:qr-steps}).
  \item Over $\twoR$, let $X$ be an invertible real matrix, and use the $*$-square roots $(c,1)$ of the self-adjoint elements $(c,c)$. Then the Cholesky factorisation of $(X,X^T)$ exists exactly when the leading principal minors of $X$ are non-zero, and is Doolittle's $LU$ factorisation of $X$. One step of the QR algorithm exists exactly when those of $X^2$ are non-zero, and is two steps of Rutishauser's LR algorithm on $X$ when these exist (Proposition~\ref{prop:2R-algorithms}). The QR step is the same over the similarity quiver, which encodes this problem in degree~$1$ (Section~\ref{subsec:cat-similarity}).
  \item Over the field $\RR((\varepsilon))$ of Laurent series, which is the completed localisation of the quiver with one loop, the Gram--Schmidt process with the inner product $\sum_i\varepsilon^{i-1}x_iy_i$, followed by taking standard parts, computes a Bruhat decomposition of every invertible real matrix, and the orders of the norms determine its permutation (Theorem~\ref{thm:bruhat}).
\end{enumerate}
\end{mainthm}

\begin{table}[htbp]
\centering
\caption{One program, many algorithms: what programs written once against the seven operations compute over different scalars. For a quiver, the scalars are the completed localisation of its path $*$-algebra, and the matrices are the Hermitian matrices of degree~$1$ that it encodes. The statements, with their hypotheses, are those of Theorem~\ref{thm:B} and the results that it cites; the row for Householder tridiagonalisation over $\twoR$ holds whenever every step is defined (Section~\ref{subsec:2R}). Table~\ref{tab:algorithms} gives further cases.}
\label{tab:one-program}
\small
\renewcommand{\arraystretch}{1.25}
\begin{tabular}{@{}>{\raggedright\arraybackslash}p{0.27\linewidth}>{\raggedright\arraybackslash}p{0.22\linewidth}>{\raggedright\arraybackslash}p{0.43\linewidth}@{}}
\toprule
\textbf{Program, written once} & \textbf{Scalars} & \textbf{Classical algorithm obtained} \\
\midrule
Householder tridiagonalisation & $\RR$, $\CC$, $\HH$ & Householder tridiagonalisation \\ \addlinespace[3pt]
 & the real SVD quiver & Golub--Kahan bidiagonalisation \\ \addlinespace[3pt]
 & $\twoR$ & tridiagonalisation by similarity \\
\rowrule
QR step, by the Gram--Schmidt process & $\RR$, $\CC$, $\HH$ & QR step \\ \addlinespace[3pt]
 & the SVD quivers & two steps: a QR step on $W^*W$ and on $WW^*$; one step: a Cholesky LR step on each \\ \addlinespace[3pt]
 & the Takagi quiver & two steps: a QR step on $WW^*$; one step: a Cholesky LR step on $WW^*$, by a unitary congruence of $W$ \\ \addlinespace[3pt]
 & $\twoR$ & two steps of Rutishauser's LR algorithm \\
\rowrule
Cholesky factorisation & $\twoR$ & Doolittle $LU$ factorisation \\
\rowrule
Gram--Schmidt process, weighted by $\varepsilon^{i-1}$ & $\RR((\varepsilon))$ & Bruhat decomposition \\
\bottomrule
\end{tabular}
\end{table}

Parts~(2) and~(3) concern different $*$-algebras, and different relations between the QR and LR algorithms. As algebras, both are $\RR\oplus\RR$ in degree~$0$, but for the real singular value decomposition the involution fixes the two vertex idempotents, so that the unitary group is $\Or(n)\times\Or(m)$ and the matrix $W$ sits on the arrows, whereas in $\twoR$ it exchanges them, so that the unitary group is $\GL_n(\RR)$ and $X$ sits in degree~$0$. In~(2) the LR steps are Cholesky steps on positive definite Gram matrices, two of which make one QR step, as for every positive definite matrix~\cite{Xu1998}; in~(3) they are steps of the LR algorithm on an invertible real matrix, which need not be symmetric. It is~(3) that unifies the QR algorithm with Rutishauser's LR algorithm.

In each case the classical algorithm arises from the overloaded operations alone: the programs contain no code for bidiagonalisation, for the LR algorithm or for elimination. This is a statement about how the algorithms are related, not about their cost: the specialised implementations, such as those in LAPACK, are faster and better studied, and the generic programs are not proposed as replacements for them. Nor is it a statement about rounding: the identifications are exact, and the numerical stability of the generic programs is not studied (Section~\ref{sec:open-problems}). Example~\ref{ex:golub-kahan-numbers} follows Householder's algorithm through a $3\times2$ example by hand, and Table~\ref{tab:one-program} collects these results by program.

The generic algorithms can also break down where the encoded problem has a solution. Over the congruence quiver, which encodes Sylvester's law of inertia, every column whose entries are linear combinations of arrows is isotropic, $x^*x=0$, so Householder reduction stops at the first non-zero column; the generic algorithm there is the $LDL^*$ factorisation, which needs no norms (Section~\ref{subsec:cat-ldl}). Over $\twoR$, the QR step on an invertible $X$ breaks down exactly when a leading principal minor of $X^2$ vanishes (Proposition~\ref{prop:2R-algorithms}). And when the quiver is wild, as for a vertex with two loops, there is no normal form of the kind found in the catalogue (the number of parameters grows quadratically, not linearly), and a Jacobi-type method needs quotients that do not exist in $\Lambda$ but do exist in a non-commutative skew field, such as Cohn's free field (Section~\ref{sec:wild-behaviour}).

\section{Two Examples}
\label{sec:motivation}
The approach grew out of two examples. The first shows that the limit of an algorithm can be computed by overloading its operations; the second, that overloading can turn one classical algorithm into another, and where it can fail.

\subsection{Example: From QR to LU by an Algebraic Limit}
\label{subsec:qr-lu}
The QR and LU decompositions are among the most important factorisations in numerical analysis~\cite{GolubVanLoan2013}. Compare the Gram--Schmidt process, which computes a QR factorisation, with Gaussian elimination, which computes an LU factorisation~\cite{GolubVanLoan2013}:
\begin{align}
    u_k &= a_k - \sum_{j=1}^{k-1} \frac{\langle u_j, a_k\rangle}{\langle u_j, u_j\rangle}\, u_j \quad &(\text{Gram--Schmidt for QR}), \label{eq:gs} \\
    u_k^T &= a_k^T - \sum_{j=1}^{k-1} l_{kj}\, u_j^T \quad &(\text{Gaussian elimination for LU}). \label{eq:ge}
\end{align}
Here $a_k$ denotes the $k$-th column of a square matrix $A$ in \eqref{eq:gs} and the $k$-th row of $A$ in \eqref{eq:ge}. The two recurrences have the same shape: each step subtracts from the current vector suitable multiples of the vectors already produced. They differ only in how the multipliers are chosen.

This resemblance can be made precise. Fix $\varepsilon>0$ and replace the Euclidean inner product by the weighted inner product
\begin{equation}\label{eq:weighted-ip}
  \langle x, y\rangle_\varepsilon \;=\; \sum_{i=1}^{n} \varepsilon^{\,i-1} x_i y_i ,
  \qquad\text{with quadratic form}\qquad q_\varepsilon(x)=\sum_{i=1}^n \varepsilon^{\,i-1}x_i^2 .
\end{equation}
For small $\varepsilon$, the first coordinate dominates every inner product, then the second, and so on. Running the (unnormalised) Gram--Schmidt process \eqref{eq:gs} with respect to $\langle\cdot,\cdot\rangle_\varepsilon$ and letting $\varepsilon\to 0^+$ produces the factorisation that Gaussian elimination computes.

\begin{proposition}\label{prop:gs-limit}
Let $A\in M_n(\RR)$ have non-zero leading principal minors. Apply the recurrence \eqref{eq:gs}, with the inner product \eqref{eq:weighted-ip}, to the columns of $A$, and denote the resulting vectors by $u_1(\varepsilon),\dots,u_n(\varepsilon)$. As $\varepsilon\to 0^+$, the matrix $[u_1(\varepsilon)\,\cdots\,u_n(\varepsilon)]$ converges to the lower triangular factor $\hat L$ of the Crout factorisation $A=\hat L R$, in which $R$ is unit upper triangular, and the Gram--Schmidt coefficients converge to the entries of $R$.
\end{proposition}
\begin{proof}[Sketch of proof]
Induction on $k$ shows that the $i$-th entry of $u_k(\varepsilon)$ is $O(\varepsilon^{k-i})$ for $i<k$. Because $u_j$ is $\langle\cdot,\cdot\rangle_\varepsilon$-orthogonal to $u_1,\dots,u_{j-1}$, the coefficient $\langle u_j,a_k\rangle_\varepsilon/\langle u_j,u_j\rangle_\varepsilon$ equals $\langle u_j,a_k^{(j)}\rangle_\varepsilon/\langle u_j,u_j\rangle_\varepsilon$, where $a_k^{(j)}$ is $a_k$ with its components along $u_1,\dots,u_{j-1}$ removed. Both inner products are dominated by their $j$-th terms, which are of order $\varepsilon^{j-1}$, so the coefficient tends to the ratio of the $j$-th entries of $a_k^{(j)}$ and $u_j$. Subtracting this multiple of $u_j$ annihilates the $j$-th entry, which is exactly the elimination step of Crout's method. (Crout's method eliminates by columns; applied to $A^T$ it is the row elimination \eqref{eq:ge}.)
\end{proof}

Theorem~\ref{thm:bruhat} in Chapter~\ref{chap:infinitesimals} removes the assumption on the minors: for every invertible $A$, the limit is a Bruhat decomposition, whose permutation is read off from the orders of magnitude of the norms. Two caveats apply. The process must be left unnormalised: if the vectors are normalised, as in the usual QR factorisation, the normalising factors $q_\varepsilon(u_k)^{-1/2}$ grow like $\varepsilon^{-(k-1)/2}$ and there is no limit. And classical Gram--Schmidt is rarely used to compute QR factorisations in practice, because it can lose orthogonality for ill-conditioned matrices~\cite{Higham2002,GolubVanLoan2013}. Householder reflections~\cite{GolubVanLoan2013} avoid both, and their limit is a genuine elimination step. Let $D_\varepsilon=\diag(1,\varepsilon,\dots,\varepsilon^{n-1})$ and let $a_1$ be the first column of $A$, with $a_{11}\neq 0$. The $D_\varepsilon$-orthogonal reflection that maps $a_1$ to a multiple of $e_1$ is
\[
P_\varepsilon \;=\; I-\frac{2\,v v^{T}D_\varepsilon}{v^{T}D_\varepsilon v},\qquad v=a_1+\sign(a_{11})\,q_\varepsilon(a_1)^{1/2}\,e_1 .
\]
As $\varepsilon\to 0^+$ we have $q_\varepsilon(a_1)^{1/2}\to|a_{11}|$, $D_\varepsilon v\to 2a_{11}e_1$ and $v^TD_\varepsilon v\to 4a_{11}^2$, hence
\begin{equation}\label{eq:householder-limit}
P_\varepsilon \;\longrightarrow\; G \;=\; \diag(-1,1,\dots,1)\,\bigl(I-\ell\, e_1^{T}\bigr),
\qquad \ell=\Bigl(0,\tfrac{a_{21}}{a_{11}},\dots,\tfrac{a_{n1}}{a_{11}}\Bigr)^{T}.
\end{equation}
The matrix $I-\ell e_1^T$ is the Gauss transformation that eliminates the first column of $A$ below the diagonal, and $G$ differs from it only by the sign of the first row; like each $P_\varepsilon$, it is an involution. Applying the same reasoning to the trailing submatrix, on which the weights are again graded, shows that Householder QR with respect to $\langle\cdot,\cdot\rangle_\varepsilon$ tends to Gaussian elimination without pivoting, up to signs. Limits of this kind, of Householder QR applied to a problem with heavily weighted rows, also appear in the analysis of the weighting method for least squares problems with equality constraints, where they yield elimination methods~\cite{CoxHigham1999}.

\paragraph{From limits to infinitesimals.}
Computing such limits by hand is laborious, and it has to be repeated for every algorithm. The key observation is that the limit can be computed \emph{algebraically}, by running the algorithm once over a larger field, that is, by overloading its operations. Every step of Householder QR is a sequence of field operations and square roots applied to quantities depending on $\varepsilon$. Suppose we execute the algorithm with $\varepsilon$ a \emph{formal positive infinitesimal}: an element of an ordered field $\KK\supset\RR$ satisfying $0<\varepsilon<1/m$ for every positive integer $m$. A suitable field is the field of real algebraic Puiseux series (Section~\ref{sec:puiseux}), whose elements are convergent series $\sum_{q} c_q\varepsilon^{q}$ with rational exponents. All intermediate quantities are then exact elements of $\KK$, and at the end we keep only the constant terms (the \emph{standard parts}) of the outputs. Proposition~\ref{prop:autolimit} shows that, whenever the outputs are finite, these standard parts are exactly the limits as $\varepsilon\to 0^+$ of the outputs of the same algorithm executed in real arithmetic. We call a limit computed in this way an \emph{auto-limit}, by analogy with automatic differentiation. Thus the analytic operation ``take the limit of the algorithm'' is replaced by the algebraic operation ``run the algorithm over a larger field and truncate'', just as automatic differentiation replaces symbolic differentiation by arithmetic in the dual numbers.

Two properties of $\KK$ are needed for this to work. First, square roots of positive elements must exist in $\KK$, since the algorithm computes norms; this holds because $\KK$ is a \emph{real closed field} (Definition~\ref{def:rcf}). Second, the facts used to prove the algorithm correct over $\RR$ (for example, that $P_\varepsilon$ is orthogonal for $\langle\cdot,\cdot\rangle_\varepsilon$, or that norms of non-zero vectors are positive) must remain true over $\KK$. This is guaranteed by the \emph{transfer principle} of Tarski (Theorem~\ref{thm:tarski-transfer}). Chapter~\ref{chap:infinitesimals} develops both points. In the present example the infinitesimal $\varepsilon$ lives in a commutative field, so the limit can also be computed directly, as in the proof of Proposition~\ref{prop:gs-limit}. The algebraic viewpoint is more convenient once $\varepsilon$ is combined with the more general $*$-algebras of the next example and of later chapters, but there too the limit can be found by asymptotic analysis, by computer algebra with $\varepsilon$ as a symbol or, less reliably, by numerical evaluation at small $\varepsilon$.

\subsection{Example: Similarity Transformations via the \texorpdfstring{$*$}{*}-Algebra \texorpdfstring{$\twoR$}{2R}}
\label{subsec:2R}
Not every relationship between algorithms is a limit. Consider the problem of reducing a real square matrix to tridiagonal form by a similarity transformation. Every real square matrix is similar to a tridiagonal matrix, as its real Jordan form shows (Section~\ref{subsec:cat-similarity}). For symmetric matrices, orthogonal reduction to tridiagonal form precedes the computation of eigenvalues; general matrices are usually reduced to Hessenberg form instead, because reduction to tridiagonal form by similarity can break down~\cite{GolubVanLoan2013}. Over $\RR$, the natural canonical form under similarity is the real Jordan form, obtained from the primary decomposition of the matrix.

Let $\twoR:=\RR\oplus\RR$ with componentwise operations and with the involution that exchanges the two components, $(x,y)^*=(y,x)$. This is the algebra of split-complex numbers $p+qj$, with $j^2=1$ and $j^*=-j$, written in the basis of its idempotents $\tfrac12(1\pm j)$. A matrix over $\twoR$ is a pair $(X,Y)$ of real matrices, and the involution of $\twoR$ induces the conjugate transpose $(X,Y)^*=(Y^T,X^T)$. It follows that
\begin{itemize}
  \item the Hermitian matrices, $(X,Y)^*=(X,Y)$, are exactly the pairs $(X,X^T)$ with $X\in M_n(\RR)$ arbitrary;
  \item the unitary matrices, $U^*U=I$, are exactly the pairs $U=(P,P^{-T})$ with $P\in\GL_n(\RR)$;
  \item for such $U$ and $H=(X,X^T)$, one has $UHU^*=(PXP^{-1},\,(PXP^{-1})^T)$.
\end{itemize}
Thus \emph{unitary similarity of Hermitian matrices over $\twoR$ is exactly similarity of arbitrary real matrices}. If we take an algorithm for the unitary tridiagonalisation of Hermitian matrices written against the seven operations, such as Listing~\ref{lst:householder}, and execute it with these operations overloaded for $\twoR$, the result is an algorithm that reduces arbitrary real matrices to tridiagonal form by similarity. Here $*$-square roots always exist: a self-adjoint element $(c,c)$ has the $*$-square root $z=(c,1)$, since $z^*z=(1,c)(c,1)=(c,c)$, whatever the sign of $c$.

In exact arithmetic the generic algorithm can, however, break down: it may encounter a non-zero vector $v=(x,y)$ with $v^*v=(y^Tx,y^Tx)=0$, by which it would have to divide. This mirrors the ``serious breakdown'' of the unsymmetric Lanczos process~\cite{GolubVanLoan2013}. Over $\twoR$ the QR algorithm in turn becomes Rutishauser's LR algorithm, two steps at a time, and the Cholesky factorisation becomes Doolittle's $LU$ factorisation (Proposition~\ref{prop:2R-algorithms}).

The analogue of the spectral theorem in its diagonal form fails here: a matrix with non-real eigenvalues, or one that is not diagonalisable, is not similar to a real diagonal matrix. What remains is the Hermitian eigendecomposition in the sense used throughout this thesis: a unitary similarity to a direct sum of blocks, none of which is unitarily similar to a direct sum of smaller blocks (compare Theorem~\ref{thm:hermitian-ks}). Over $\twoR$ these blocks are the Hermitian matrices $(J,J^T)$ for the real Jordan blocks $J$, each of which can be taken tridiagonal (Section~\ref{subsec:cat-similarity}), so the real Jordan form is the Hermitian eigendecomposition over $\twoR$.

This example shows the two themes of the thesis in miniature: a classical algorithm transported to a new $*$-algebra solves a different classical problem, and the points where it fails (division by elements that are not invertible) are precisely where further theory is required.

\section{Limitations of Existing Unifying Frameworks}
\label{sec:limitations}
Several theories relate matrix decompositions to each other, but they do not, as they stand, provide transformations of programs.

\emph{Lie groups and symmetric spaces.} A Lie group is a group that is also a smooth manifold, such as $\GL_n(\RR)$, the orthogonal group $\Or(n)$ or the symplectic group $\Sp(2n,\RR)$. Many matrix decompositions are instances of structural decompositions of Lie groups: the polar decomposition and the SVD are instances of the Cartan and $KAK$ decompositions, and the QR decomposition is an instance of the Iwasawa decomposition~\cite{Helgason1978,EdelmanJeong2023}. A symmetric space is a quotient $G/K$ of a Lie group by the fixed points of an involutive automorphism; algebraic treatments exist~\cite{Loos1969}. Edelman and Jeong~\cite{EdelmanJeong2023} derive fifty-three matrix factorisations systematically from the Cartan ($KAK$) decompositions of the classical symmetric spaces. These theories describe the \emph{sets} and \emph{groups} that appear in a decomposition, and they are expressed in terms of Lie algebras and root systems. Numerical code, however, is expressed in terms of arithmetic on scalars. The theories therefore explain why decompositions exist, but in general they do not tell us how to turn one program into another. An exception is the family of Jacobi-type methods, which have been formulated on compact Lie algebras so that one algorithm covers several eigenvalue problems~\cite{KleinsteuberHelmkeHuper2004}; see Section~\ref{sec:related}.

\emph{Scalar product spaces.} Mackey, Mackey and Tisseur~\cite{MackeyMackeyTisseur2006} unify many structured factorisations by working in a \emph{scalar product space}, that is, a real or complex vector space with a non-degenerate bilinear or sesquilinear form; the matrices that are self-adjoint, skew-adjoint or unitary with respect to the form give rise to Jordan algebras, Lie algebras and automorphism groups. This is the closest existing framework to ours. Its scalars remain $\RR$ or $\CC$, however, and the scalar product is an extra datum. In our setting the form is absorbed into the involution of a matrix $*$-algebra, and the scalars may themselves be quaternions, dual numbers, or elements of any finite-dimensional real $*$-algebra. Conversely, matrices over many of these $*$-algebras can be represented by larger real or complex matrices on which the involution acts as the adjoint for a form: a quaternionic matrix by a complex matrix of twice the size, with the conjugate transpose~\cite{Zhang1997}, and a matrix $(X,Y)$ over $\twoR$ by $\diag(X,Y)$, with the adjoint for the indefinite form $\bigl(\begin{smallmatrix}0&I\\I&0\end{smallmatrix}\bigr)$.

\emph{Mismatch of objectives.} Decompositions are used differently in representation theory and in numerical computation. The \emph{Bruhat decomposition} of $\GL_n$ writes an invertible matrix as $M=LPU$, with $L$ lower triangular, $U$ upper triangular and $P$ a permutation matrix. This resembles LU factorisation with pivoting, but it is qualitatively different: the permutation $P$ is uniquely determined by $M$ (it records the ranks of the leading, upper-left, submatrices of $M$) and it depends discontinuously on $M$. In computable analysis every computable function on real inputs is continuous~\cite{Weihrauch2000}, so $P$ cannot be computed from arbitrarily accurate approximations of $M$; it is therefore unclear what it would mean to compute the Bruhat decomposition numerically. From exact input it can be computed directly, by Gaussian elimination with a fixed pivot rule (Section~\ref{sec:bruhat}), and also by running the Gram--Schmidt process over a field of Laurent series (Theorem~\ref{thm:bruhat}).\footnote{Care is needed for decompositions whose factors are not unique, such as the SVD, where the singular vectors are determined only up to symmetries. There the appropriate notion is that of a computable multi-valued function~\cite{Weihrauch2000}.} By contrast, LU factorisation with partial pivoting chooses $P$ by comparisons of approximate magnitudes, and this choice is not required to be unique.\footnote{The Jordan normal form suffers from the same difficulty as the Bruhat decomposition: its block structure is uniquely determined by the matrix and varies discontinuously with it. In numerical practice one replaces it by decompositions with less rigidly determined factors, such as the Schur decomposition. Our framework covers both kinds: the real Jordan form is among the decompositions of Theorem~\ref{thm:A}(3), and Theorem~\ref{thm:B}(4) computes the Bruhat decomposition from exact input, while the algorithms of Theorem~\ref{thm:B}(1)--(3) compute decompositions of the latter kind, such as tridiagonal and bidiagonal forms.}

More generally, the purpose of a decomposition in representation theory is usually to \emph{classify} objects up to isomorphism, typically via the Krull--Schmidt theorem (Theorem~\ref{thm:krull-schmidt}) and its analogues for forms (Theorem~\ref{thm:hermitian-ks}). The purpose of a decomposition in numerical analysis is weaker: to find a basis in which certain computations become easy. Despite this difference, the computationally useful decompositions are related in the way the classifying ones are: each is a normal form under unitary similarity of Hermitian matrices over some $*$-algebra (Theorem~\ref{thm:A}). It is this relation that the thesis exploits.

\section{Related Work and the State of the Art}
\label{sec:related}
We now situate the thesis within existing work. No previous framework, to our knowledge, overloads the operations of numerical algorithms for finite-dimensional $*$-algebras, including non-semisimple ones, and uses the representation theory of quivers to identify the problem that the transported algorithm solves.

\paragraph{Structured numerical linear algebra.}
Structure-preserving algorithms for matrices that are symmetric, skew-symmetric, Hamiltonian, symplectic, or self-adjoint with respect to an indefinite form are surveyed in~\cite{BunseGerstnerByersMehrmann1992}; the scalar product spaces of Mackey, Mackey and Tisseur~\cite{MackeyMackeyTisseur2006} were discussed in Section~\ref{sec:limitations}. Quaternionic matrices and their SVD are treated in~\cite{Zhang1997}; Hamiltonian square roots of skew-Hamiltonian matrices are studied in~\cite{FassbenderMackeyMackeyXu1999}, and we use them in Section~\ref{sec:hermitian-sqrt}.

\paragraph{Classical relations between algorithms.}
Some of the relations that Theorem~\ref{thm:B} establishes are classical in special cases. Many of the decompositions are also, classically, the Hermitian eigendecomposition of a larger block matrix, such as $\bigl(\begin{smallmatrix}0&W^*\\W&0\end{smallmatrix}\bigr)$ for the singular value decomposition of $W$, and the thesis points these out where they arise. For a square upper bidiagonal matrix $W$, a perfect shuffle of the rows and columns of $\bigl(\begin{smallmatrix}0&W^T\\W&0\end{smallmatrix}\bigr)$ gives a symmetric tridiagonal matrix with zero diagonal, whose eigenvalues are plus and minus the singular values of $W$~\cite{GolubKahan1965,DemmelKahan1990}; this is the encoding of Section~\ref{sec:unification-example} with the labels of the arrows removed. Two steps of the Cholesky LR algorithm on a positive definite matrix are one step of the QR algorithm~\cite{Xu1998}. The framework explains these relations by the algebra of the scalars, and extends them, under the hypotheses of Theorem~\ref{thm:B}, to several of the other $*$-algebras of Chapter~\ref{chap:indexing}.

\paragraph{Lie-theoretic unification.}
The Lie-theoretic frameworks of Section~\ref{sec:limitations} unify decompositions through the groups involved~\cite{EdelmanJeong2023,KleinsteuberHelmkeHuper2004}. Our approach is complementary: we unify through the \emph{scalars}, so that the same program, with overloaded arithmetic, implements the different cases. The two viewpoints meet in Chapter~\ref{chap:indexing}, where the unitary matrices of degree~$0$ turn out to form products of the classical groups (Theorem~\ref{thm:encoding}).

\paragraph{Operator algebras in numerical analysis.}
Hagen, Roch and Silbermann~\cite{HagenRochSilbermann2001} use $C^*$-algebras to study approximation methods for operators on infinite-dimensional Hilbert spaces. A typical result (see also~\cite{BoettcherSilbermann1999}) is that the finite section method for a Toeplitz operator is stable if and only if a certain element of a quotient $C^*$-algebra of matrix sequences is invertible. There, $C^*$-algebras serve as a tool for proving stability and convergence of discretisations. In this thesis, $*$-algebras play a different role: they are the scalars over which a finite-dimensional algorithm is run. Moreover, our algebras are generally not $C^*$-algebras, since a finite-dimensional $C^*$-algebra is semisimple (Remark~\ref{rem:cstar}), whereas the non-semisimple algebras, such as the dual numbers, are essential to our applications. Section~\ref{sec:operator} discusses the possibility of extending our framework to operator algebras.

\paragraph{Automatic differentiation and perturbation theory.}
Forward-mode automatic differentiation evaluates a program over the dual numbers~\cite{GriewankWalther2008}. For matrix functions, the Fr\'echet derivative of $f$ at $A$ in direction $E$ can be read off from $f$ applied to the block matrix $\bigl(\begin{smallmatrix}A&E\\0&A\end{smallmatrix}\bigr)$~\cite{Higham2008}, which is the same construction viewed through a matrix representation of the dual numbers. Analytic perturbation theory describes how eigenvalues and eigenvectors of a symmetric matrix depend on a parameter; by Rellich's theorem they can be chosen analytically~\cite{Kato1995}, and analytic singular value decompositions can be computed numerically~\cite{BunseGerstnerByersMehrmannNichols1991}. Decompositions of matrices over the dual numbers were studied by the author in~\cite{Gutin2022DualNumbers}. Chapter~\ref{chap:infinitesimals} recovers and extends its spectral theorems and its singular value decomposition with respect to the transpose (the T-SVD) by a single algebraic argument.

\paragraph{Representation theory of quivers and of forms.}
A quiver is a directed graph, and the representations of quivers describe modules over finite-dimensional algebras~\cite{AssemSimsonSkowronski2006,Schiffler2014,DerksenWeyman2017}. Gabriel's classification of the quivers of finite representation type~\cite{Gabriel1972} and Kac's theorem~\cite{Kac1980} relate representation type to Dynkin diagrams and root systems (this is not the theorem of Gabriel whose analogue is Theorem~\ref{thm:A}(1)), and Drozd's theorem~\cite{Drozd1980} shows that every finite-dimensional algebra over an algebraically closed field is tame or wild. Species~\cite{DlabRingel1976} extend quivers to non-algebraically closed fields. Representations with additional orthogonal or symplectic structure are described by symmetric quivers~\cite{DerksenWeyman2002} and by Shmelkin's signed quivers~\cite{Shmelkin2006}. Sergeichuk reduced the classification of systems of forms and linear maps to the classification of representations of quivers together with Hermitian forms over division algebras~\cite{Sergeichuk1988}; see also~\cite{HornSergeichuk2006}. The decomposition of Hermitian forms into orthogonally indecomposable summands was studied in additive categories by Quebbemann, Scharlau and Schulte~\cite{QuebbemannScharlauSchulte1979}, and in the context of linear algebra with duality by Lorand~\cite{Lorand2020}. Our doubly-signed quivers combine the signed quivers with a second signing that records real division algebras, and the novelty lies in using this combinatorial data to transform algorithms.

\paragraph{Real algebraic geometry and infinitesimals.}
Tarski's theorem~\cite{Tarski1951} shows that the theory of real closed fields admits quantifier elimination, so that first-order statements transfer between real closed fields~\cite{BochnakCosteRoy1998,Marker2002}. Computational real algebraic geometry routinely uses formal infinitesimals to perturb inputs into general position~\cite{BasuPollackRoy2006}. Robinson's nonstandard analysis~\cite{Robinson1966} provides infinitesimals with a transfer principle for a much richer language, but for the single extension $\RR\subset{}^*\RR$ and at the cost of constructivity. We use the former and compare it with the latter in Chapter~\ref{chap:infinitesimals}.

\paragraph{Localisation of path algebras.}
The localisation of a path algebra at its arrows is a universal localisation in the sense of Cohn~\cite{Cohn2006}, and it is the algebra of the free groupoid of the quiver~\cite{GabrielZisman1967}. The free field of Cohn~\cite{Cohn2006,CohnReutenauer1999} is the universal skew field of fractions of a free algebra, and Schofield~\cite{Schofield1985} studied homomorphisms from rings to skew fields and to simple Artinian rings. Leavitt path algebras~\cite{AbramsAraSilesMolina2017} and graph $C^*$-algebras~\cite{Raeburn2005} adjoin partial inverses instead. Each edge $e$ acquires a ``ghost'' edge $e^*$, with $e^*e$ the idempotent of its range, but $ee^*$ is the idempotent of its source only when no other edge leaves that source; in general an edge is a partial isometry, not an invertible element. For a cycle the two constructions agree, $L_k(C_n)\cong M_n(k[x,x^{-1}])$, but for a vertex with two loops the Leavitt path algebra is the Leavitt algebra $L(1,2)$, in which the loops have left inverses and are not invertible, whereas the localisation is the group algebra of the free group. The ghost edges are also new generators, whereas the involution of a doubly-signed quiver maps arrows to arrows of the same quiver.

\paragraph{Position of this thesis.}
The thesis combines these strands. From structured linear algebra and Lie theory it takes the catalogue of decompositions to be unified; from representation theory it takes symmetric and signed quivers, which it extends to real $*$-algebras by a second signing; from model theory it takes the transfer principle that justifies running algorithms over fields with infinitesimals; and from the localisation of rings it takes the means of dividing by elements that are zero divisors in the original algebra, while the theory of skew fields shows how far such means can go.

\section{Contributions and Outline}
\label{sec:contributions}
The main contributions of the thesis are the seven overloaded operations and the two outcomes of Section~\ref{sec:outcomes}. The chapters are as follows.

\emph{Chapter~\ref{chap:staralg}} collects the algebraic background: $*$-algebras, radicals and semisimple algebras, Hermitian modules, and quivers and path algebras. It states the Krull--Schmidt theorem and its analogue for Hermitian modules, and it determines exactly when the uniqueness part of the latter holds for quadratic forms.

\emph{Chapter~\ref{chap:matrixops}} gives the seven operations their meaning over a $*$-algebra, and shows that Householder reflectors work with any $*$-square root for which the denominator of the reflector is invertible (Proposition~\ref{prop:householder-any-root}). Its example is the first part of Theorem~\ref{thm:B}: over the quiver of the singular value decomposition, Householder tridiagonalisation becomes the bidiagonalisation of Golub and Kahan (Proposition~\ref{prop:golub-kahan}).

\emph{Chapter~\ref{chap:infinitesimals}} turns limits of algorithms into exact computations over a field with infinitesimals, using Tarski's transfer principle and Puiseux series. Its main results are the Bruhat decomposition by the Gram--Schmidt process over a field of Laurent series, which is part~(4) of Theorem~\ref{thm:B} (Theorem~\ref{thm:bruhat}), and Dijkstra's shortest-path algorithm as the leading-order part of Gaussian elimination (Proposition~\ref{prop:dijkstra}).

\emph{Chapter~\ref{chap:quivers}} introduces doubly-signed quivers, which adapt Shmelkin's signed quivers~\cite{Shmelkin2006} to the real numbers, and constructs their path $*$-algebras (Proposition~\ref{prop:block-model}). Its main result is Theorem~\ref{thm:A}(1), an analogue of Gabriel's theorem: every finite-dimensional real $*$-algebra is an algebra of block matrices over a quotient of the path $*$-algebra of its doubly-signed quiver (Theorem~\ref{thm:structure}).

\emph{Chapter~\ref{chap:indexing}} proves the encoding theorem, which is Theorem~\ref{thm:A}(2) (Theorem~\ref{thm:encoding}), and applies it to a catalogue of quivers, which gives Theorem~\ref{thm:A}(3) (Table~\ref{tab:catalogue}). It then proves the remaining parts of Theorem~\ref{thm:B} (Section~\ref{sec:encoded-algorithms}), and relates the behaviour of generic algorithms to the representation type of the quiver (Section~\ref{sec:wild-behaviour}).

\emph{Chapter~\ref{division-problem}} constructs the algebra in which the algorithms run: the \emph{completed localisation} of a path $*$-algebra, in which every arrow is inverted between its end points. It shows how inverses and $*$-square roots are computed there, and when a result holds in the original algebra (Theorem~\ref{thm:reduction}).

\emph{Chapter~\ref{chap:implementation}} describes a Python library that implements the seven operations and the algorithms written against them, and reports experiments in which they reproduce the classical algorithms to within rounding errors.

\emph{Chapter~\ref{chap:conclusion}} summarises the results, with the status of each, and lists limitations and open problems.

Appendix~\ref{chap:gutin-2022}, which is based on the preprint~\cite{Gutin2022Visualisation}, visualises the dynamics of the QR algorithm on $2\times 2$ real matrices, in a way inspired by Yaglom's accounts of the geometry of low-dimensional algebras~\cite{Yaglom1968-ComplexNumbers,Yaglom1979-SimpleNonEuclid}. It draws a matrix in the hyperbolic plane of inner products on $\RR^2$, and shows the QR algorithm becoming the LR algorithm as the inner product degenerates, as in Proposition~\ref{prop:gs-limit}. Appendix~\ref{chap:code} contains the source code of the library of Chapter~\ref{chap:implementation}; its test suite is part of the source distribution submitted with the thesis.

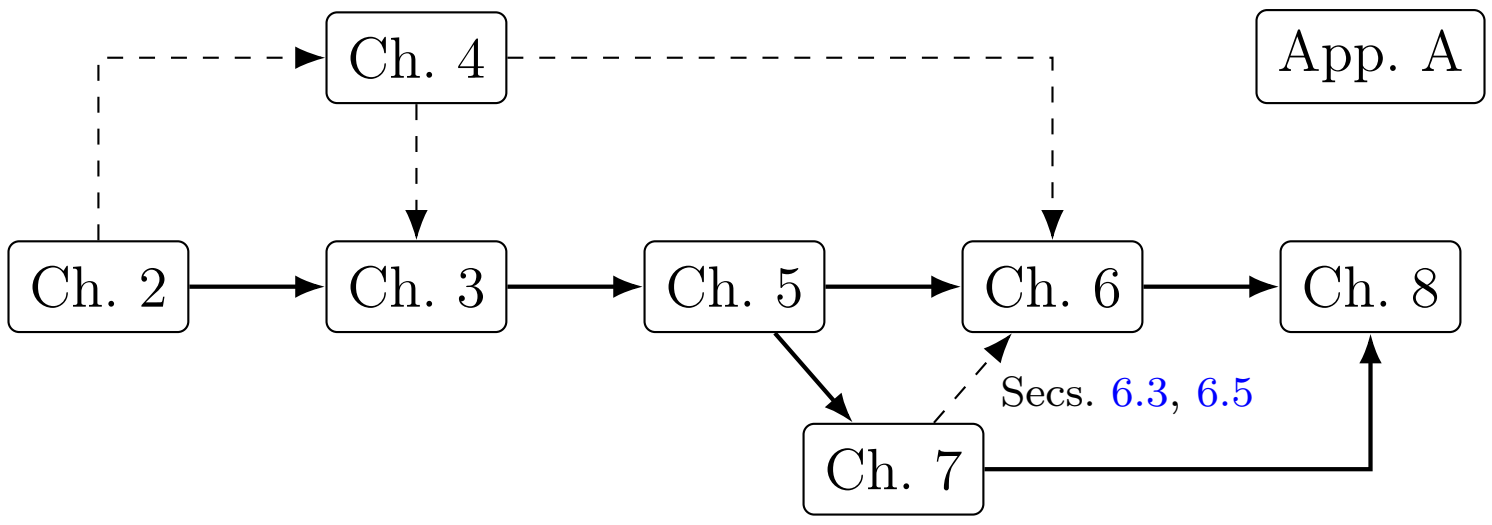
\begin{figure}[ht]
\centering
\begin{tikzpicture}[
  ch/.style={draw, rounded corners=2pt, minimum width=10mm, minimum height=6mm, font=\small},
  dep/.style={-{Latex[length=2mm]}, thick},
  side/.style={-{Latex[length=2mm]}, dashed}]
\node[ch] (c2) {Ch.~2};
\node[ch, right=9mm of c2] (c3) {Ch.~3};
\node[ch, right=9mm of c3] (c5) {Ch.~5};
\node[ch, right=9mm of c5] (c6) {Ch.~6};
\node[ch, right=9mm of c6] (c8) {Ch.~8};
\node[ch, below=9mm of $(c5)!0.5!(c6)$] (c7) {Ch.~7};
\node[ch, above=9mm of c3] (c4) {Ch.~4};
\node[ch, above=9mm of c8] (ca) {App.~A};
\draw[dep] (c2) -- (c3);
\draw[dep] (c3) -- (c5);
\draw[dep] (c5) -- (c6);
\draw[dep] (c6) -- (c8);
\draw[dep] (c5) -- (c7);
\draw[dep] (c7) -| (c8);
\draw[side] (c7) -- node[right=4pt, font=\scriptsize, pos=0.3] {Secs.~\ref{sec:encoded-algorithms}, \ref{sec:wild-behaviour}} (c6);
\draw[side] (c4) -- (c3);
\draw[side] (c2) |- (c4);
\draw[side] (c4) -| (c6);
\end{tikzpicture}
\caption{Dependencies between the chapters. Solid arrows show the main line. Dashed arrows show the side branch of Chapter~\ref{chap:infinitesimals}, and the use of Chapter~\ref{division-problem} in Sections~\ref{sec:encoded-algorithms} and~\ref{sec:wild-behaviour}. No chapter depends on Appendix~\ref{chap:gutin-2022}.}
\label{fig:reading-guide}
\end{figure}

\paragraph{Reading guide.} Figure~\ref{fig:reading-guide} shows how the chapters depend on each other. Chapter~\ref{chap:staralg} collects background and can be consulted as needed. The main line runs through Chapters~\ref{chap:matrixops}, \ref{chap:quivers}, \ref{chap:indexing} and~\ref{division-problem}. Sections~\ref{sec:encoded-algorithms} and~\ref{sec:wild-behaviour} use results of Chapter~\ref{division-problem}; Section~\ref{sec:encoded-algorithms} begins with a summary of what it uses, so Chapter~\ref{division-problem} can also be read after Chapter~\ref{chap:indexing}. Chapter~\ref{chap:infinitesimals} is a side branch: it is needed for the limits of algorithms in Section~\ref{subsec:qr-lu}, for the transfer principle with which Proposition~\ref{prop:symplectic-sqrt} passes from $\RR$ to real closed fields, for Theorem~\ref{thm:B}(4), and for the power-series results to which the catalogue of Chapter~\ref{chap:indexing} refers. Chapter~\ref{chap:implementation} describes the implementation of Chapters~\ref{chap:matrixops} to~\ref{division-problem}, and no chapter depends on Appendix~\ref{chap:gutin-2022}.

%======================================================================
%   CHAPTER 2: STAR-ALGEBRAS AND HERMITIAN MODULES
%======================================================================
\chapter{\texorpdfstring{$*$}{*}-Algebras, Modules and Quivers}
\label{chap:staralg}

This chapter fixes notation and recalls the algebraic objects on which the thesis is built: $*$-algebras, which supply the scalars for the seven operations of Section~\ref{sec:proposal}; the structure of finite-dimensional algebras (radicals, semisimple and Artinian rings, the theorems of Wedderburn, basic algebras and idempotents, tensor products); modules over them, and modules equipped with Hermitian forms; and quivers and path algebras, which describe finite-dimensional algebras combinatorially. Most of the material is standard; the references are~\cite{Lam2001,Knus1991,Scharlau1985,KMRT1998,AssemSimsonSkowronski2006}. The main point of the chapter is Section~\ref{sec:krull-schmidt}, which states the Krull--Schmidt theorem and its analogue for Hermitian modules, and determines precisely when the uniqueness part of the latter can hold for quadratic forms.

\section{Notation and Conventions}
\label{sec:notation}
\begin{convention}\label{conv:notation}
Throughout the thesis the following conventions are in force.
\begin{enumerate}
  \item Rings are associative and have an identity element; ring homomorphisms preserve the identity. Unless stated otherwise, $k$ denotes a field of characteristic different from $2$, and algebras are finite-dimensional over $k$.
  \item The involution of a $*$-ring or $*$-algebra is always written $x\mapsto x^*$, using the (centred) mathematical asterisk, and $x^*$ is called the \emph{adjoint} of $x$. Complex conjugation of a scalar is written $z\mapsto\bar z$, and the transpose of a matrix $X$ is written $X^T$. For matrices of real or complex numbers we write $X^H=\bar X^T$.
  \item If $A$ is a $*$-algebra, the \emph{conjugate transpose} of $X=(x_{ij})\in M_n(A)$ is the matrix $X^*$ with entries $(X^*)_{ij}=(x_{ji})^*$. With this involution $M_n(A)$ is again a $*$-algebra. We write $\Herm_n(A)=\{X\in M_n(A): X^*=X\}$ and $\Un_n(A)=\{U\in M_n(A): U^*U=I\}$.
  \item The symbols $\dagger$, $(\cdot)^{\top}$ and $\overline{(\cdot)}$ are reserved for involutions of quivers (Chapter~\ref{chap:quivers}); they are never used for the involution of an algebra. When two different involutions of the same algebra have to be compared, the second one receives its own symbol, such as $X^\omega$ for the symplectic involution (Section~\ref{sec:hermitian-sqrt}) or $x^\natural$ for the involution induced on a quotient (Lemma~\ref{lem:star-wedderburn}).
  \item $\twoR$ denotes the $*$-algebra $\RR\oplus\RR$ with involution $(x,y)^*=(y,x)$.
\end{enumerate}
\end{convention}

\section{\texorpdfstring{$*$}{*}-Algebras}
\label{sec:staralg-def}
Many algebras that occur in analysis, geometry and physics carry a natural ``adjoint'' or ``conjugation'' operation. Examples are complex conjugation, the transpose and conjugate transpose of matrices, and the adjoint of an operator on a Hilbert space. All of these reverse the order of products and square to the identity, and abstracting these two properties leads to the notion of an involution.

\begin{definition}\label{def:starring}
A \emph{\starring} (or \emph{involutive ring}) is a ring $A$ together with a map $A\to A$, $x\mapsto x^*$, called an \emph{involution}, such that for all $x,y\in A$:
\begin{enumerate}
    \item $(x^*)^* = x$;
    \item $(x + y)^* = x^* + y^*$ and $(xy)^* = y^*x^*$.
\end{enumerate}
\end{definition}

It follows from the axioms that $1^*=1$, since $1^*=1^*(1^*)^*=(1^*\,1)^*=(1^*)^*=1$. Thus the involution is an anti-automorphism of $A$ that is its own inverse.

\begin{definition}[$*$-algebra]\label{def:staralg}
Let $R$ be a commutative \starring{} with involution $r \mapsto r'$. A \emph{\staralg{} over $R$} is an associative $R$-algebra $A$ which is also a \starring{}, such that
\[(rx)^* = r'x^* \qquad\text{for all } r\in R,\ x\in A.\]
A \emph{$*$-homomorphism} $f: A \to B$ between $*$-algebras is an algebra homomorphism such that $f(a^*) = f(a)^*$ for all $a \in A$.
\end{definition}

If $R=\CC$ with complex conjugation, the involution of $A$ is \emph{conjugate-linear}: $(\lambda x)^*=\bar\lambda x^*$. If the involution of $R$ is trivial, as for $R=\RR$, the involution of $A$ is $R$-linear. Every \starring{} is a \staralg{} over $\ZZ$ with the trivial involution. Unless said otherwise, the $*$-algebras in this thesis are finite-dimensional over $k=\RR$ (or over a real closed field), with $k$-linear involution.

\begin{example}\label{ex:staralgs}
The following $*$-algebras will be used repeatedly.
\begin{enumerate}
  \item $\CC$ with $z^*=\bar z$, viewed as a two-dimensional $\RR$-algebra.
  \item $M_n(\RR)$ with $X^*=X^T$; $M_n(\CC)$ with $X^*=X^H$.
  \item The quaternions $\HH$ with quaternion conjugation $(a+bi+cj+dk)^*=a-bi-cj-dk$.
  \item The \emph{dual numbers} $\RR[\varepsilon]/(\varepsilon^2)$, with either the trivial involution $\varepsilon^*=\varepsilon$ or the involution $\varepsilon^*=-\varepsilon$.
  \item The algebra $\twoR$ of Convention~\ref{conv:notation}, whose involution exchanges the two summands.
  \item $M_{2n}(k)$ with the \emph{symplectic involution} $X^*=J^{-1}X^TJ$, where $J=\bigl(\begin{smallmatrix}0&I_n\\-I_n&0\end{smallmatrix}\bigr)$.
  \item For any $*$-algebra $A$, the matrix algebra $M_n(A)$ with the conjugate transpose of Convention~\ref{conv:notation}.
\end{enumerate}
The algebras in (2), (3) and (6) are \emph{semisimple}, that is, they have no non-zero nilpotent ideals (Definition~\ref{def:semisimple}); the dual numbers are not, since $\RR\varepsilon$ is a nilpotent ideal, and $\twoR$ is semisimple but not simple.
\end{example}

\begin{remark}[$C^*$-algebras]\label{rem:cstar}
A \emph{$C^*$-algebra} is a complex Banach algebra with a conjugate-linear involution whose norm satisfies $\|x^*x\|=\|x\|^2$ for all $x$. Every finite-dimensional $C^*$-algebra is $*$-isomorphic to a finite direct sum $\bigoplus_i M_{n_i}(\CC)$ with the conjugate-transpose involution~\cite[Theorem III.1.1]{Davidson1996}; in particular it is semisimple. The $C^*$-identity implies that $x^*x=0$ forces $x=0$, and that a self-adjoint nilpotent element $h$ is zero, since $\|h\|^{2^m}=\|h^{2^m}\|$ for all $m$. The dual numbers violate this in either of their involutions: $\varepsilon^*\varepsilon=\pm\varepsilon^2=0$ although $\varepsilon\neq 0$. More generally, every finite-dimensional $*$-algebra $A$ with a non-zero radical (the largest nilpotent ideal of $A$; Section~\ref{subsec:radicals}) admits no $C^*$-norm: the radical is stable under the involution (Lemma~\ref{lem:radical-stable}), so for $0\neq x\in\rad(A)$ the element $h=x^*x$ is self-adjoint and nilpotent, hence $h=0$ and then $x=0$. Since non-semisimple algebras are central to this thesis, the theory of $C^*$-algebras does not apply directly to our setting; we return to this point in Section~\ref{sec:operator}.
\end{remark}

\section{The Structure of Finite-Dimensional Algebras}
\label{sec:alg-background}

\subsection{Radicals, Semisimple Algebras and Artinian Rings}
\label{subsec:radicals}
Let $A$ be a ring. The \emph{(Jacobson) radical} $\rad(A)$ is the intersection of all maximal left ideals of $A$; it is a two-sided ideal, and it also equals the intersection of all maximal right ideals~\cite[\S4]{Lam2001}. If $A$ is a finite-dimensional algebra over a field, then $\rad(A)$ is the largest nilpotent two-sided ideal of $A$: there is an $N$ with $\rad(A)^N=0$, and every nilpotent two-sided ideal is contained in $\rad(A)$. Intuitively, the radical consists of the ``infinitesimal'' elements of $A$, which are invisible to the simple modules.

\begin{definition}\label{def:semisimple}
A finite-dimensional algebra $A$ is \emph{semisimple} if $\rad(A)=0$. It is \emph{simple} if its only two-sided ideals are $0$ and $A$.
\end{definition}

\begin{definition}\label{def:artinian}
A ring $R$ is \emph{left Artinian} if every descending chain $I_1\supseteq I_2\supseteq\cdots$ of left ideals of $R$ is eventually constant. A \emph{simple Artinian ring} is a non-zero left Artinian ring whose only two-sided ideals are $0$ and $R$.
\end{definition}

Every finite-dimensional algebra is Artinian, since a strictly descending chain of subspaces has length at most the dimension. The ring $\ZZ$ is not Artinian, because $2\ZZ\supsetneq4\ZZ\supsetneq8\ZZ\supsetneq\cdots$. The following theorem, which is the theorem of Wedderburn used in Chapter~\ref{division-problem}, describes semisimple and simple Artinian rings completely.

\begin{theorem}[Wedderburn--Artin]\label{thm:wedderburn-artin}
A ring $R$ is left Artinian with $\rad(R)=0$ if and only if
\[ R\cong M_{n_1}(D_1)\times\cdots\times M_{n_r}(D_r) \]
for some integers $n_i\geq 1$ and division rings (skew fields) $D_i$. The pairs $(n_i,D_i)$ are uniquely determined up to isomorphism of the $D_i$ and permutation. In particular, $R$ is simple Artinian if and only if $R\cong M_n(D)$ for a single $n$ and a single division ring $D$.
\end{theorem}
See~\cite[Theorems~3.5 and~3.10]{Lam2001}. Over $k=\RR$ the finite-dimensional division algebras are $\RR$, $\CC$ and $\HH$ (Frobenius' theorem), so every finite-dimensional semisimple real algebra is a product of matrix algebras over $\RR$, $\CC$ and $\HH$. We shall also need the following classical result about automorphisms of simple algebras. An algebra over $k$ is \emph{central simple} if it is simple and its centre is $k$; examples are $M_n(k)$ and, for $k=\RR$, $M_n(\HH)$.

\begin{theorem}[Skolem--Noether]\label{thm:skolem-noether}
Let $A$ be a finite-dimensional central simple algebra over a field $k$. Every $k$-algebra automorphism $\varphi$ of $A$ is inner: there is an invertible $g\in A$ with $\varphi(x)=gxg^{-1}$ for all $x\in A$.
\end{theorem}
See~\cite{Jacobson1989,KMRT1998}. Applied to the composite of an anti-automorphism with the transpose, it shows that every anti-automorphism of $M_n(k)$ has the form $X\mapsto GX^TG^{-1}$; this is how the types of involutions are determined in Chapter~\ref{chap:quivers}.

\begin{example}\label{ex:radicals}
\begin{enumerate}
  \item The dual numbers $\RR[\varepsilon]/(\varepsilon^2)$ have radical $\RR\varepsilon$ and $A/\rad(A)\cong\RR$.
  \item In the algebra $T_n(k)$ of lower triangular $n\times n$ matrices, the radical consists of the strictly lower triangular matrices, and $T_n(k)/\rad(T_n(k))\cong k^n$, the diagonal.
  \item $M_n(\RR)$, $\HH$, $\CC$ and $\twoR$ are semisimple; $M_n(\RR)$ and $\HH$ are simple, whereas $\twoR$ is not.
\end{enumerate}
\end{example}

\begin{lemma}\label{lem:radical-stable}
If $A$ is a finite-dimensional $*$-algebra, then $\rad(A)^*=\rad(A)$.
\end{lemma}
\begin{proof}
The involution is an anti-automorphism, so $\rad(A)^*$ is a two-sided ideal, and it is nilpotent because $(\rad(A)^*)^N=(\rad(A)^N)^*=0$. Hence $\rad(A)^*\subseteq\rad(A)$. Applying the involution once more gives the reverse inclusion.
\end{proof}

The semisimple quotient $A/\rad(A)$ can be realised inside $A$. The following result is the associative analogue of the Levi decomposition of a Lie algebra.

\begin{theorem}[Wedderburn--Malcev]\label{thm:wedderburn-malcev}
Let $A$ be a finite-dimensional algebra over a field $k$ of characteristic $0$ (more generally, assume that $A/\rad(A)$ is separable over $k$). Then there is a subalgebra $S\subseteq A$ such that $A=S\oplus\rad(A)$ as vector spaces; in particular $S\cong A/\rad(A)$. Any two such subalgebras are conjugate by an element of the form $1+r$ with $r\in\rad(A)$. If $A$ is a $*$-algebra, $S$ can be chosen with $S^*=S$.
\end{theorem}
See~\cite[\S11.6]{Pierce1982} for the first two statements and~\cite{Taft1957} for the last. We call $A=S\oplus\rad(A)$ a \emph{Wedderburn--Malcev decomposition} of $A$. When $\rad(A)^2=0$, the existence of a $*$-stable $S$ has an elementary proof, which we give in Lemma~\ref{lem:star-wedderburn}.

\subsection{Basic Algebras and Idempotents}
\label{subsec:basic}
In ordinary linear algebra the entries of a matrix are scalars, and a matrix algebra $M_n(k)$ is built from the one-dimensional algebra $k$. For a general finite-dimensional algebra, the role of the scalars is played by its \emph{basic algebra}: the smallest algebra from which it can be rebuilt by forming matrices. We first give the definition, which uses only the semisimple quotient, and then explain how idempotents make the idea precise.

\begin{definition}\label{def:basic}
A finite-dimensional algebra $A$ is \emph{basic} if $A/\rad(A)$ is a product of division algebras, that is, if all the integers $n_i$ in Theorem~\ref{thm:wedderburn-artin} for $A/\rad(A)$ equal $1$.
\end{definition}

For example, $\RR$, $\CC$, $\HH$, the dual numbers, $\twoR$ and $T_n(k)$ are basic, whereas $M_n(k)$ for $n\geq2$ is not. Its basic algebra is $k$, from which it is rebuilt as $M_n(k)$.

An element $e\in A$ is an \emph{idempotent} if $e^2=e$. Idempotents $e,f$ are \emph{orthogonal} if $ef=fe=0$, and a non-zero idempotent is \emph{primitive} if it is not the sum of two non-zero orthogonal idempotents. A \emph{complete set of primitive orthogonal idempotents} is a set $\{e_1,\dots,e_n\}$ of pairwise orthogonal primitive idempotents with $e_1+\cdots+e_n=1$. In $M_n(k)$ the diagonal matrix units $E_{11},\dots,E_{nn}$ form such a set. Given such a set, every $a\in A$ decomposes as
\begin{equation}\label{eq:peirce}
 a=\sum_{i,j} e_iae_j,\qquad e_iae_j\in e_iAe_j,
\end{equation}
the \emph{Peirce decomposition}, which generalises the decomposition of a matrix into its entries. Two idempotents $e,f$ are \emph{equivalent} if the left modules $Ae$ and $Af$ are isomorphic; equivalently, there are $u\in fAe$ and $v\in eAf$ with $vu=e$ and $uv=f$. A finite-dimensional algebra is basic if and only if, in a complete set of primitive orthogonal idempotents, no two distinct idempotents are equivalent~\cite[\S I.6]{AssemSimsonSkowronski2006}. If $e$ is the sum of one primitive idempotent from each equivalence class, then $AeA=A$ and the algebra $eAe$ is basic. The algebras $A$ and $eAe$ are then \emph{Morita equivalent}: their categories of modules are equivalent~\cite{AssemSimsonSkowronski2006,AuslanderReitenSmalo1995}. In this sense every finite-dimensional algebra behaves like a ``matrix algebra over its basic algebra''.

For $*$-algebras we need a variant that takes the involution into account, since a corner $fAf$ is a $*$-algebra only when $f$ is self-adjoint. The involution maps idempotents to idempotents, but it need not fix any primitive idempotent. For instance, in $M_2(\RR)$ with the symplectic involution $X^*=J^{-1}X^TJ$, the only self-adjoint elements are the scalar matrices, and the involution exchanges $E_{11}$ and $E_{22}$; this $*$-algebra cannot be cut down to a smaller one. We therefore count only self-adjoint idempotents.

\begin{definition}\label{def:basic-star}
A finite-dimensional $*$-algebra $A$ is \emph{basic as a $*$-algebra} if every self-adjoint idempotent of $A/\rad(A)$ is central. If $f\in A$ is a self-adjoint idempotent with $AfA=A$ and $fAf$ is basic as a $*$-algebra, then $fAf$ is called a \emph{basic $*$-algebra of $A$}.
\end{definition}

Without the word ``self-adjoint'', the condition would say that $A/\rad(A)$ is a product of division algebras, which is Definition~\ref{def:basic}, because $M_n(D)$ has a non-central idempotent exactly when $n\geq2$. Thus $M_n(\RR)$ with the transpose has the basic $*$-algebra $E_{11}M_n(\RR)E_{11}\cong\RR$, whereas $M_2(\RR)$ with the symplectic involution is itself basic as a $*$-algebra, although it is not basic as an algebra. The difference between the two is exactly the difference between the vertex signs $+1$ and $-1$ of a signed quiver (Chapter~\ref{chap:quivers}). It would not suffice to ask that the involution permute a complete set of primitive idempotents: in $M_2(\RR)$ with the involution $X\mapsto\eta X^T\eta$, $\eta=\diag(1,-1)$, the involution exchanges the primitive idempotents $\frac12\bigl(\begin{smallmatrix}1&\pm1\\\pm1&1\end{smallmatrix}\bigr)$, but $E_{11}$ is a non-central self-adjoint idempotent, and the basic $*$-algebra is $E_{11}M_2(\RR)E_{11}\cong\RR$. Over $\RR$, a $*$-algebra $A$ is basic as a $*$-algebra exactly when $A/\rad(A)$ is a product of copies of the ten $*$-algebras of Table~\ref{tab:vertex-types}, and every finite-dimensional real $*$-algebra has a basic $*$-algebra (Lemma~\ref{lem:star-simple} and Theorem~\ref{thm:structure}).

\subsection{Tensor Products}
\label{subsec:tensor}
\begin{definition}\label{def:tensor}
Let $A$ and $B$ be algebras over $k$. Their \emph{tensor product} $A\otimes_kB$ is the vector space tensor product with multiplication determined by $(a\otimes b)(a'\otimes b')=aa'\otimes bb'$. If $A$ and $B$ are $*$-algebras, $A\otimes_kB$ is a $*$-algebra with $(a\otimes b)^*=a^*\otimes b^*$.
\end{definition}

The tensor product of matrix algebras is the \emph{Kronecker product}: the map $X\otimes Y\mapsto X\otimes_{\mathrm{Kr}}Y$, where $X\otimes_{\mathrm{Kr}} Y$ is the block matrix $(x_{ij}Y)_{ij}$, is an isomorphism $M_m(k)\otimes_kM_n(k)\cong M_{mn}(k)$. More generally $M_n(A)\cong M_n(k)\otimes_kA$, and if $A$ is a $*$-algebra this is a $*$-isomorphism when $M_n(k)$ carries the transpose and $M_n(A)$ the conjugate transpose. Two further examples are used in Chapter~\ref{chap:quivers}: the complexification $B\otimes_\RR\CC$ of a real algebra $B$, for instance $\CC\otimes_\RR\CC\cong\CC\times\CC$ and $\HH\otimes_\RR\CC\cong M_2(\CC)$.

\section{Modules and Hermitian Modules}
\label{sec:modules}
Linear algebra is built on vector spaces over fields. Modules are the analogue of vector spaces when the scalars form a ring~\cite{DummitFoote2003,Lam2001}.

\begin{definition}[Module]\label{def:module}
Let $R$ be a ring. A \emph{left $R$-module} is an abelian group $(M, +)$ together with a map $R \times M \to M$, $(r, m) \mapsto rm$, such that for all $r, s \in R$ and $m, n \in M$:
\[
 r(m+n) = rm + rn,\qquad (r+s)m = rm + sm,\qquad (rs)m = r(sm),\qquad 1m = m.
\]
A \emph{right $R$-module} is defined analogously with a map $M \times R \to M$, $(m,r)\mapsto mr$, satisfying in particular $m(rs) = (mr)s$.
\end{definition}

A module is \emph{indecomposable} if it is non-zero and is not the direct sum of two non-zero submodules. It is \emph{Noetherian} (respectively \emph{Artinian}) if every ascending (respectively descending) chain of submodules is eventually constant, and it has \emph{finite length} if it is both. Every finitely generated module over a finite-dimensional algebra has finite length.

Let $A$ be a $*$-algebra and $M$ a left $A$-module. The \emph{dual module} $M^\vee=\Hom_A(M,A)$ consists of the left $A$-linear maps $f:M\to A$. It becomes a left $A$-module through the involution, by setting $(a\cdot f)(m)=f(m)\,a^*$.

We now generalise the Hermitian inner product on $\CC^n$ to modules over $*$-algebras. The key notion is that of a sesquilinear form, ``sesqui'' (Latin for ``one and a half'') referring to linearity in one argument and conjugate-linearity in the other.

\begin{definition}[Sesquilinear and $\epsilon$-Hermitian forms]\label{def:hermitian}
Let $A$ be a $*$-algebra and $M$ a left $A$-module. A \emph{sesquilinear form} on $M$ is a biadditive map $h: M \times M \to A$ such that for all $m, n\in M$ and $a\in A$:
\[
 h(am, n) = a\, h(m, n),\qquad h(m, an) = h(m, n)\, a^*.
\]
Let $\epsilon\in\{1,-1\}$. The form is \emph{$\epsilon$-Hermitian} if $h(n,m)=\epsilon\, h(m,n)^*$ for all $m,n\in M$; it is \emph{Hermitian} if $\epsilon=1$ and \emph{skew-Hermitian} if $\epsilon=-1$. It is \emph{non-degenerate} if $M$ is finitely generated and projective and the map $M\to M^\vee$, $n\mapsto h(\cdot,n)$, is an isomorphism. A pair $(M,h)$ consisting of a module and a non-degenerate $\epsilon$-Hermitian form is called an \emph{$\epsilon$-Hermitian module}, or simply a \emph{Hermitian module} when $\epsilon=1$.
\end{definition}

For $A=\CC$ with complex conjugation this is the usual notion of a Hermitian form on $\CC^n$, except that positivity is not required. If $(M,h)$ is Hermitian then $h(m,m)=h(m,m)^*$ is self-adjoint for every $m$, generalising the fact that $\langle v,v\rangle$ is real for a complex inner product.

\begin{example}\label{ex:hermitian-modules}
\begin{enumerate}
  \item If $A=\RR$ with the trivial involution, Hermitian modules are real vector spaces with a non-degenerate symmetric bilinear form, and skew-Hermitian modules are symplectic vector spaces.
  \item If $A=\CC$ with complex conjugation, a Hermitian form on $\CC^n$ is $h(x,y)=y^HSx$ for a Hermitian matrix $S$; skew-Hermitian forms are $i$ times Hermitian forms.
  \item Let $A=\twoR$ and let $e_1=(1,0)$ and $e_2=(0,1)$, so that $e_1^*=e_2$. An $A$-module is a pair $M=M_1\oplus M_2$ of real vector spaces, with $M_i=e_iM$. Since $h(e_1m,e_1n)=e_1h(m,n)e_1^*=e_1e_2\,h(m,n)=0$, a sesquilinear form pairs $M_1$ only with $M_2$, and an $\epsilon$-Hermitian form is determined by a bilinear pairing $M_1\times M_2\to\RR$. Non-degeneracy says that this pairing identifies $M_2$ with the dual space $M_1^{\vee}$. The isometries are the maps $(P,P^{-T})$ with $P\in\GL(M_1)$. This is the structure that turned unitary similarity into ordinary similarity in Section~\ref{subsec:2R}.
\end{enumerate}
\end{example}

A Hermitian form brings geometric notions into the purely algebraic setting of modules. Let $(M,h)$ be an $\epsilon$-Hermitian module.
\begin{itemize}
    \item Two elements $m, n \in M$ are \emph{orthogonal}, written $m \perp n$, if $h(m, n) = 0$. Because $h$ is $\epsilon$-Hermitian, $m \perp n$ if and only if $n \perp m$.
    \item The \emph{orthogonal complement} of a submodule $W \subseteq M$ is $W^\perp = \{ m \in M : h(m, w) = 0 \text{ for all } w \in W \}$.
    \item If $M = M_1 \oplus M_2$ with $M_1 \subseteq M_2^\perp$, the decomposition is an \emph{orthogonal sum}, written $(M,h) = (M_1,h_1) \perp (M_2,h_2)$, where $h_i$ is the restriction of $h$.
    \item An \emph{isometry} $(M,h)\to(M',h')$ is an isomorphism $g:M\to M'$ with $h'(gm,gn)=h(m,n)$. The isometries of $(M,h)$ to itself form the \emph{unitary group} $\Un(M,h)$.
    \item The \emph{adjoint} of $T\in\End_A(M)$ is the unique $T^\sharp\in\End_A(M)$ with $h(Tm,n)=h(m,T^\sharp n)$; it exists by non-degeneracy. $T$ is \emph{self-adjoint} if $T^\sharp=T$. For $M=A^n$ and the standard form $h(x,y)=\sum_i x_iy_i^*$, the adjoint is the conjugate transpose and $\Un(M,h)=\Un_n(A)$.
\end{itemize}

The decompositions unified in this thesis are, in this language, orthogonal decompositions of a Hermitian module into submodules invariant under a self-adjoint endomorphism. For example, the spectral theorem for a Hermitian matrix $S\in M_n(\CC)$ asserts that $(\CC^n,h)$, with $h$ the standard form, is the orthogonal sum of the eigenspaces of $S$.

A particularly simple class of $\epsilon$-Hermitian modules is formed by the hyperbolic modules. If $P$ is a finitely generated projective $A$-module, the \emph{hyperbolic module} $\mathbb{H}_\epsilon(P)$ is $P\oplus P^\vee$ with the form
\[
  h\bigl((p,f),(q,g)\bigr) = g(p) + \epsilon\, f(q)^* .
\]
This form is sesquilinear, $\epsilon$-Hermitian and non-degenerate. The hyperbolic module depends only on $P$ and not on any form on $P$. For example, for $A=P=\RR$ and $\epsilon=1$, the hyperbolic plane $\mathbb{H}_1(\RR)$ has the Gram matrix $\bigl(\begin{smallmatrix}0&1\\1&0\end{smallmatrix}\bigr)$, which is congruent to $\bigl(\begin{smallmatrix}1&0\\0&-1\end{smallmatrix}\bigr)$ via the change of basis $\bigl(\begin{smallmatrix}1&1\\1&-1\end{smallmatrix}\bigr)$, up to scaling.

\section{The Krull--Schmidt Theorem and Its Hermitian Analogue}
\label{sec:krull-schmidt}
None of the results in this section is new; they are recalled because the classification problems addressed later are of this kind, and because the precise statements clarify what can and cannot be expected of numerical decompositions.

\begin{theorem}[Krull--Schmidt]\label{thm:krull-schmidt}
Let $R$ be a ring and $M$ a non-zero $R$-module of finite length. Then
\begin{enumerate}
  \item $M\cong M_1\oplus\cdots\oplus M_r$ for some indecomposable modules $M_1,\dots,M_r$;
  \item if also $M\cong N_1\oplus\cdots\oplus N_s$ with every $N_j$ indecomposable, then $r=s$ and there is a permutation $\pi$ of $\{1,\dots,r\}$ such that $M_i\cong N_{\pi(i)}$ for all $i$.
\end{enumerate}
\end{theorem}
For a proof see~\cite[\S19]{Lam2001} or~\cite{AssemSimsonSkowronski2006}. Existence follows by induction on the length. Uniqueness rests on \emph{Fitting's lemma}: every endomorphism of an indecomposable module of finite length is either invertible or nilpotent. It follows that the endomorphism ring of such a module is a local ring.

\begin{remark}\label{rem:ks-uniqueness}
The uniqueness in Theorem~\ref{thm:krull-schmidt} is uniqueness \emph{up to isomorphism}: the summands themselves, as submodules, are generally not unique. For example, $k^2=k e_1\oplus k e_2=k(e_1+e_2)\oplus ke_2$. This is the essential difference between classification in representation theory and decomposition in numerical analysis, where one wants an actual basis or an actual collection of submodules.
\end{remark}

We now ask for an analogue of the Krull--Schmidt theorem for Hermitian modules. A non-degenerate $\epsilon$-Hermitian module $(M,h)$ is \emph{orthogonally indecomposable} if $M\neq 0$ and $(M,h)$ is not an orthogonal sum of two non-zero $\epsilon$-Hermitian modules. Two questions arise: does every $(M,h)$ decompose as an orthogonal sum of orthogonally indecomposable modules, and is such a decomposition unique up to isometry of the summands? The answers are, respectively, ``yes'' and ``it depends on the field''~\cite{QuebbemannScharlauSchulte1979,Sergeichuk1988,Lorand2020}.

\begin{proposition}\label{prop:orth-existence}
Let $A$ be a finite-dimensional $*$-algebra. Every non-degenerate $\epsilon$-Hermitian $A$-module is an orthogonal sum of finitely many orthogonally indecomposable $\epsilon$-Hermitian modules.
\end{proposition}
\begin{proof}
We argue by induction on the length of $M$. If $(M,h)$ is orthogonally indecomposable there is nothing to prove. Otherwise $(M,h)=(M_1,h_1)\perp(M_2,h_2)$ with $M_1,M_2\neq 0$. The summands are again non-degenerate, because $M^\vee=M_1^\vee\oplus M_2^\vee$ and the map $n\mapsto h(\cdot,n)$ respects the decomposition; and both have smaller length than $M$. The induction hypothesis applies to each summand.
\end{proof}

Uniqueness is more delicate. It fails already in the simplest case $A=k$ with the trivial involution, where Hermitian modules are symmetric bilinear spaces. In characteristic different from $2$ every such space has an orthogonal basis, so the orthogonally indecomposable modules are the one-dimensional forms $\langle a\rangle$, $a\in k^\times$, and $\langle a\rangle\cong\langle b\rangle$ if and only if $b\in a\,k^{\times 2}$. Uniqueness of orthogonal decompositions therefore means uniqueness of diagonalisations up to squares. Recall that a field is \emph{quadratically closed} if every element is a square, and \emph{Euclidean} if it is ordered and every positive element is a square (for instance $\RR$, or any real closed field).

\begin{proposition}\label{prop:quadratic-uniqueness}
Let $k$ be a field of characteristic different from $2$. The following are equivalent.
\begin{enumerate}
  \item Every non-degenerate symmetric bilinear space over $k$ has a unique orthogonal decomposition into one-dimensional summands, up to isometry and permutation of the summands.
  \item $k$ is quadratically closed or Euclidean.
\end{enumerate}
\end{proposition}
\begin{proof}
(2)$\Rightarrow$(1). If $k$ is quadratically closed, every one-dimensional form is isometric to $\langle 1\rangle$, and the number of summands is the dimension. If $k$ is Euclidean, every one-dimensional form is isometric to $\langle 1\rangle$ or $\langle -1\rangle$, and the numbers of each are determined by Sylvester's law of inertia, whose usual proof works over any ordered field in which positive elements are squares~\cite{Lam2005}.

(1)$\Rightarrow$(2). We use the elementary fact that if $c\in k^\times$ is represented by the form $\langle a,b\rangle$, then $\langle a,b\rangle\cong\langle c,abc\rangle$~\cite{Lam2005}. Suppose first that $-1$ is a square. Then $\langle 1,1\rangle\cong\langle1,-1\rangle$ is a hyperbolic plane, which represents every $c\in k^\times$, so $\langle1,1\rangle\cong\langle c,c\rangle$; by (1), $c$ is a square. Hence $k$ is quadratically closed. Suppose next that $-1$ is not a square. The hyperbolic plane $\langle1,-1\rangle$ again represents every $c$, so $\langle 1,-1\rangle\cong\langle c,-c\rangle$, and (1) forces $c\in k^{\times2}$ or $c\in -k^{\times 2}$. Thus $k^\times=k^{\times2}\cup(-k^{\times2})$. Moreover, if $c=x^2+y^2\neq 0$ then $\langle1,1\rangle\cong\langle c,c\rangle$, and (1) forces $c\in k^{\times 2}$, so sums of squares are squares; in particular $-1$ is not a sum of squares, so $k$ is formally real. Declaring the non-zero squares to be the positive elements then defines an ordering in which every positive element is a square, so $k$ is Euclidean.
\end{proof}

Proposition~\ref{prop:quadratic-uniqueness} shows that uniqueness fails over most fields of characteristic different from $2$, for example over $\QQ$ (where $\langle1,1\rangle\cong\langle2,2\rangle$ but $\langle1\rangle\not\cong\langle2\rangle$), over finite fields of odd characteristic (where $\langle1,1\rangle\cong\langle u,u\rangle$ for a non-square $u$), and over $\RR(t)$. In characteristic $2$ the theory is different again: symmetric bilinear forms need not be diagonalisable (alternating forms are not), and quadratic forms are no longer determined by their bilinear forms; we exclude this case throughout. For the fields used in this thesis, uniqueness holds.

\begin{theorem}[Hermitian Krull--Schmidt theorem]\label{thm:hermitian-ks}
Let $k$ be an algebraically closed field of characteristic different from $2$, or a real closed field (for example $\RR$). Let $A$ be a finite-dimensional $*$-algebra over $k$ and $\epsilon\in\{1,-1\}$. Then every non-degenerate $\epsilon$-Hermitian $A$-module is an orthogonal sum of orthogonally indecomposable $\epsilon$-Hermitian modules, and the summands are unique up to isometry and permutation.
\end{theorem}

Existence is Proposition~\ref{prop:orth-existence}. For uniqueness we refer to~\cite{QuebbemannScharlauSchulte1979,Sergeichuk1988} (see also~\cite{Wall1976,Knus1991,Lorand2020}) and only outline the argument, which is a chain of reductions.
\begin{enumerate}
  \item \emph{Reduction to the semisimple case.} Let $J=\rad(A)$ be the radical of $A$, its largest nilpotent ideal (Section~\ref{subsec:radicals}), which is stable under the involution (Lemma~\ref{lem:radical-stable}). The form $h$ induces a non-degenerate form on $M/JM$ over $A/J$. Orthogonal decompositions of $M/JM$ lift to orthogonal decompositions of $M$, uniquely up to isometry, by lifting idempotents (elements $e$ with $e^2=e$; Section~\ref{subsec:basic}) of $\End_A(M)$ modulo its radical compatibly with the adjoint involution~\cite{QuebbemannScharlauSchulte1979}.
  \item \emph{Reduction to simple $*$-stable factors.} A semisimple $*$-algebra is a product of simple algebras that are either stable under the involution or exchanged in pairs $B\times B^{\mathrm{op}}$. Over an exchanged pair, every $\epsilon$-Hermitian module is hyperbolic, $\mathbb{H}_\epsilon(P)$ for a $B$-module $P$, and its decompositions are governed by the ordinary Krull--Schmidt theorem for $P$.
  \item \emph{Reduction to division algebras.} A simple $*$-stable algebra is $M_n(D)$ for a division algebra $D$ (Theorem~\ref{thm:wedderburn-artin}). Hermitian Morita theory~\cite[Ch.~I, \S9]{Knus1991} identifies $\epsilon$-Hermitian modules over $M_n(D)$ with $\epsilon'$-Hermitian forms over $D$, with respect to a suitable involution of $D$ and some $\epsilon'=\pm1$. Concretely, for a rank-one idempotent $e\in M_n(D)$, a module $M$ over $M_n(D)$ corresponds to the $D$-vector space $eM$, and the form on $M$ corresponds to an $\epsilon'$-Hermitian form on $eM$. For example, for $M_{2m}(\RR)$ with the symplectic involution, Hermitian modules correspond to skew-symmetric bilinear forms over $\RR$.
  \item \emph{Division algebras.} Over a division algebra with involution, non-alternating forms are diagonalisable and alternating ones are sums of hyperbolic planes. Uniqueness thus reduces to the uniqueness of diagonalisations. When $k$ is real closed, $D$ is $k$, $k(\sqrt{-1})$ or the quaternions over $k$, and uniqueness follows from Sylvester's law of inertia and its Hermitian and quaternionic analogues; when $k$ is algebraically closed, $D=k$ and the forms are classified by their dimension.
\end{enumerate}
By Proposition~\ref{prop:quadratic-uniqueness}, the last step, and with it the uniqueness statement of Theorem~\ref{thm:hermitian-ks}, fails for $A=k$ unless $k$ is quadratically closed or Euclidean. This is why the thesis works over $\RR$, over real closed fields, and over algebraically closed fields.

\section{Quivers and Path Algebras}
\label{sec:quivers-basic}
Quivers provide the combinatorial language in which we describe finite-dimensional algebras. We introduce them here and develop the variants with involution in Chapter~\ref{chap:quivers}.

\begin{definition}[Quiver]\label{def:quiver}
A \emph{quiver} $Q=(Q_0,Q_1,s,t)$ consists of a finite set $Q_0$ of \emph{vertices}, a finite set $Q_1$ of \emph{arrows}, and two maps $s,t:Q_1\to Q_0$ assigning to each arrow $\alpha$ its \emph{source} $s(\alpha)$ and \emph{target} $t(\alpha)$. We write $\alpha: i\to j$ if $s(\alpha)=i$ and $t(\alpha)=j$. Loops ($i=j$) and multiple arrows are allowed.
\end{definition}

In other words, a quiver is a finite directed multigraph. A \emph{path} of length $\ell\geq1$ is a sequence of arrows $p=\alpha_\ell\cdots\alpha_2\alpha_1$ with $t(\alpha_i)=s(\alpha_{i+1})$; it starts at $s(p)=s(\alpha_1)$ and ends at $t(p)=t(\alpha_\ell)$. We compose from right to left, as for functions. For every vertex $v$ there is a path $e_v$ of length $0$ with $s(e_v)=t(e_v)=v$.

\begin{definition}[Path algebra]\label{def:path-algebra}
The \emph{path algebra} $kQ$ of a quiver $Q$ over a field $k$ is the vector space with basis the set of all paths in $Q$, including the paths $e_v$ of length zero. The product of two paths $p$ and $q$ is their concatenation $pq$ if $t(q)=s(p)$, and $0$ otherwise; this is extended bilinearly to all of $kQ$.
\end{definition}

The elements $e_v$ form a complete set of primitive orthogonal idempotents, and every path satisfies $p=e_{t(p)}\,p\,e_{s(p)}$. Thus $e_jkQe_i$ is spanned by the paths from $i$ to $j$, and the Peirce decomposition \eqref{eq:peirce} of an element of $kQ$ is its decomposition according to start and end points. The path algebra is finite-dimensional if and only if $Q$ has no oriented cycles. The \emph{arrow ideal} $R_Q$ is the two-sided ideal spanned by the paths of length at least one. If $Q$ has no oriented cycles then $R_Q=\rad(kQ)$. For every quiver, the quotient
\[ kQ/R_Q^2 \]
is finite-dimensional, with basis the vertices and the arrows; in it, the product of any two arrows is zero. It is a basic algebra whose radical $R_Q/R_Q^2$ squares to zero. Conversely, by a theorem of Gabriel, every basic finite-dimensional algebra over an algebraically closed field is a quotient $kQ/I$ of the path algebra of a uniquely determined quiver $Q$ by an ideal with $R_Q^m\subseteq I\subseteq R_Q^2$ for some $m$~\cite[Theorem~II.3.7]{AssemSimsonSkowronski2006}. In particular, the basic algebras with radical square zero are exactly the algebras $kQ/R_Q^2$.

\begin{example}\label{ex:path-algebras}
\begin{enumerate}
  \item For the linearly oriented quiver $1\xrightarrow{\alpha_1}2\xrightarrow{\alpha_2}\cdots\xrightarrow{\alpha_{n-1}}n$, the map $e_i\mapsto E_{ii}$, $\alpha_i\mapsto E_{i+1,i}$ extends to an isomorphism of $kQ$ with the algebra $T_n(k)$ of lower triangular matrices; the path from $i$ to $j>i$ maps to $E_{ji}$. (With the opposite convention for composing paths one obtains upper triangular matrices.)
  \item For the quiver with one vertex and one loop $x$, the path algebra is the polynomial ring $k[x]$, and $kQ/R_Q^2=k[x]/(x^2)$ is the algebra of dual numbers.
  \item For the quiver with two vertices $u,v$ and two arrows $a:u\to v$ and $b:v\to u$, the algebra $kQ/R_Q^2$ is four-dimensional, with basis $e_u,e_v,a,b$ and $ab=ba=0$. It reappears in Section~\ref{sec:unification-example}.
\end{enumerate}
\end{example}

\begin{definition}[Representation]\label{def:representation}
A \emph{representation} $V$ of $Q$ over $k$ assigns a finite-dimensional vector space $V_i$ to each vertex $i$ and a linear map $V_\alpha:V_{s(\alpha)}\to V_{t(\alpha)}$ to each arrow $\alpha$. Its \emph{dimension vector} is $(\dim V_i)_{i\in Q_0}$. A morphism $V\to W$ is a family of linear maps $f_i:V_i\to W_i$ with $f_{t(\alpha)}V_\alpha=W_\alpha f_{s(\alpha)}$ for all $\alpha$.
\end{definition}

Representations of $Q$ are the same as $kQ$-modules, the vertex spaces being $V_i=e_iV$~\cite{AssemSimsonSkowronski2006,Schiffler2014}. Representations of dimension vector $\mathbf d$, up to isomorphism, are the orbits of the group $\prod_i\GL(d_i)$ acting on the space $\bigoplus_{\alpha}\Hom(k^{d_{s(\alpha)}},k^{d_{t(\alpha)}})$ by change of basis at every vertex. For the loop quiver of Example~\ref{ex:path-algebras}(2), a representation is a square matrix up to similarity, whose classification is the Jordan normal form. For the quiver with two vertices and a single arrow $u\to v$, a representation is a rectangular matrix up to the action $X\mapsto PXQ^{-1}$, whose classification is given by the rank. Requiring the base changes to be orthogonal instead of arbitrary turns the latter problem into the singular value decomposition; this is the role of the involutions introduced in Chapter~\ref{chap:quivers}.

\paragraph{A preview of signed and doubly-signed quivers.}
A \emph{signed quiver} (Definition~\ref{def:signed-quiver}) is a quiver equipped with an involution that reverses arrows, together with signs $\pm1$ on the vertices and arrows that the involution fixes; it describes a basic $*$-algebra with radical square zero over an algebraically closed field. A \emph{doubly-signed quiver} (Definition~\ref{def:doubly-signed}) carries in addition a second involution, which preserves arrows, together with a second set of signs on its fixed vertices; it describes $*$-algebras over the real numbers, and the second signs record whether the scalars at a vertex are real, complex or quaternionic. Figures~\ref{fig:signed-quiver-example} and~\ref{fig:dsq-examples} show examples.
%======================================================================
%   CHAPTER 3: THE SEVEN OPERATIONS OVER *-ALGEBRAS
%======================================================================
\chapter{The Seven Operations over \texorpdfstring{$*$}{*}-Algebras}
\label{chap:matrixops}

This chapter gives the seven operations of Section~\ref{sec:proposal} their meaning over a $*$-algebra. Addition, subtraction, multiplication and the adjoint need little discussion: they are the operations of the algebra, and on matrices the adjoint gives the conjugate transpose (Convention~\ref{conv:notation}). The other three do: $(i,j)$ indexing, which we base on self-adjoint idempotents and \emph{positions} (Section~\ref{sec:indexing}); division, which we define between idempotents (Section~\ref{sec:division}); and the $*$-square root, which is not unique, although Householder reflectors work with any $*$-square root in the corner of the next position for which the denominator of the reflector is invertible (Section~\ref{sec:hermitian-sqrt}). Section~\ref{sec:unification-example} then shows, in the simplest non-trivial example, how the singular value decomposition arises as a spectral theorem over a non-semisimple $*$-algebra, why the generic algorithms must leave that algebra, and how, once they do, Householder tridiagonalisation becomes the Golub--Kahan bidiagonalisation; this is the first case of Theorem~\ref{thm:B}. The algebraic background that it uses, from radicals and idempotents to quivers and path algebras, is collected in Chapter~\ref{chap:staralg}.
\section{Matrix Indexing in Finite-Dimensional \texorpdfstring{$*$}{*}-Algebras}
\label{sec:indexing}
One of the most basic operations on a matrix is to access its entry $M_{ij}$ using two indices. This operation generalises to every finite-dimensional $*$-algebra, and it still needs only two indices. The entries, however, are no longer scalars: they are elements of the corners $eAf$ cut out by idempotents $e,f$. This section explains which idempotents to use. The choice matters for algorithms: Section~\ref{sec:unification-example} shows that the order of the indices decides which classical algorithm a generic algorithm becomes. We first treat the case in which $A$ really is a matrix algebra over a smaller algebra.

\begin{proposition}\label{prop:matrix-structure}
Let $A$ be a finite-dimensional algebra with a complete set of orthogonal idempotents $e_1,\dots,e_n$ that are pairwise equivalent. Choose $u_i\in e_iAe_1$ and $v_i\in e_1Ae_i$ with $v_iu_i=e_1$ and $u_iv_i=e_i$ (with $u_1=v_1=e_1$). Then
\[
  \varphi: A\to M_n(e_1Ae_1),\qquad \varphi(a)_{ij}=v_i\,a\,u_j,
\]
is an algebra isomorphism. If $A$ is a $*$-algebra, the $e_i$ are self-adjoint and $v_i=u_i^*$, then $\varphi$ is a $*$-isomorphism, where $M_n(e_1Ae_1)$ carries the conjugate transpose.
\end{proposition}
\begin{proof}
The elements $E_{ij}:=u_iv_j$ satisfy $E_{ij}E_{kl}=u_i(v_ju_k)v_l=\delta_{jk}u_ie_1v_l=\delta_{jk}E_{il}$, because $v_ju_k\in e_1Ae_je_kAe_1$ vanishes for $j\neq k$ and equals $e_1$ for $j=k$; moreover $\sum_iE_{ii}=\sum_ie_i=1$. Using $1=\sum_k u_kv_k$ we get
\[
 \varphi(ab)_{ij}=v_iab\,u_j=\sum_k (v_iau_k)(v_kbu_j)=\sum_k\varphi(a)_{ik}\varphi(b)_{kj},
\]
so $\varphi$ is multiplicative, and $(b_{ij})\mapsto\sum_{i,j}u_ib_{ij}v_j$ is its inverse. In the $*$-case, $\varphi(a^*)_{ij}=u_i^*a^*u_j=(u_j^*au_i)^*=(\varphi(a)_{ji})^*$.
\end{proof}

Proposition~\ref{prop:matrix-structure} applies, for example, to $M_n(\RR[\varepsilon]/(\varepsilon^2))$, the matrix algebras over the dual numbers, and it identifies an algebra with a single equivalence class of primitive idempotents as a matrix algebra over its basic algebra $e_1Ae_1$.

\paragraph{Peirce indexing.} For a general finite-dimensional $*$-algebra $A$ we can use the Peirce decomposition \eqref{eq:peirce} directly. Let $E=\{e_1,\dots,e_n\}$ be a complete set of primitive orthogonal idempotents that is closed under the involution, and let $\sigma$ be the permutation with $e_i^*=e_{\sigma(i)}$. For $a\in A$ put $a_{ij}:=e_iae_j\in e_iAe_j$. Then
\begin{equation}\label{eq:generalised-indexing}
  (ab)_{ij}=\sum_{k}a_{ik}b_{kj},\qquad (a^*)_{ij}=e_ia^*e_j=(e_{\sigma(j)}\,a\,e_{\sigma(i)})^*=\bigl(a_{\sigma(j)\sigma(i)}\bigr)^*.
\end{equation}
Thus every element of $A$ is a generalised matrix $(a_{ij})$ whose $(i,j)$ entry lies in $e_iAe_j$; multiplication is matrix multiplication, and the involution is ``transpose, then re-index by $\sigma$, then apply the involution entrywise''. For $M_2(\RR)$ with the symplectic involution, $E=\{E_{11},E_{22}\}$ and $\sigma$ exchanges the two indices; formula \eqref{eq:generalised-indexing} reproduces $X^*=J^{-1}X^TJ=\bigl(\begin{smallmatrix}x_{22}&-x_{12}\\-x_{21}&x_{11}\end{smallmatrix}\bigr)$, the signs coming from the involution on the corners. For $\twoR$, $E=\{(1,0),(0,1)\}$, the off-diagonal corners vanish, and \eqref{eq:generalised-indexing} reduces to the exchange of the two components.

\paragraph{Indexing by self-adjoint idempotents.} If the indexing idempotents are self-adjoint, then $\sigma$ is the identity, and \eqref{eq:generalised-indexing} becomes
\[ (a^*)_{ij}=(a_{ji})^* : \]
the involution acts on generalised matrices exactly as the conjugate transpose acts on ordinary ones. This is the property that an algorithm written for matrices with the conjugate transpose needs in order to run without change. A set $E$ as above can be coarsened to self-adjoint idempotents by replacing $e_i$ and $e_{\sigma(i)}$, when they differ, by their sum, which is an idempotent because $e_i$ and $e_{\sigma(i)}$ are orthogonal, and which is self-adjoint. We call the resulting idempotents $1_c$, one for each orbit $c$ of $\sigma$, the \emph{orbit idempotents}. The price of the coarsening is that the corners $1_cA1_d$ are larger than the corners $e_iAe_j$. For a basic real $*$-algebra, the orbits correspond to the orbits of the vertices of its doubly-signed quiver under the transpose and the conjugation (Chapter~\ref{chap:quivers}); for a matrix algebra over it, the positions introduced next account for the matrix part.

\paragraph{Positions.} Matrices over $A$ are indexed in the same way. Let $1_c$ ($c\in\mathcal O$) be orthogonal self-adjoint idempotents of $A$ with $\sum_c1_c=1$. In $M_n(A)$ the elements $E_{kk}1_c$, for $1\leq k\leq n$ and $c\in\mathcal O$, are again orthogonal self-adjoint idempotents with sum $1$. We call the pairs $(c,k)$ \emph{positions} and write $c_k$ for $(c,k)$. A matrix $X=(x_{kl})\in M_n(A)$ becomes the generalised matrix with entries
\[ X_{c_k,d_l}=1_c\,x_{kl}\,1_d\in1_cA1_d , \]
it is multiplied as a matrix, and its adjoint is its conjugate transpose. Different orbits may carry different numbers $n_c$ of positions. The matrices are then the elements of a corner $eM_N(A)e$, as in the rectangular problems of Chapter~\ref{chap:indexing}, and they form the $*$-algebra $\End(V,\eta)\otimes_{\RR^{\mathcal O}}A$ of Definition~\ref{def:matrix-like} with $\dim V_c=n_c$ and all signs $\eta_i=+1$. By Theorem~\ref{thm:structure}, every finite-dimensional real $*$-algebra has the form $\End(V,\eta)\otimes_{\RR^{\mathcal O}}\RR Q/I$, with signs $\eta_i=\pm1$ that can be taken to be $+1$ except at orbits of signature type, and a matrix algebra over it has the same form (Corollary~\ref{cor:positions}). Position indexing is therefore the data layout of every algorithm in this thesis; signs $\eta_i=-1$ change the adjoint to \eqref{eq:signed-adjoint} (Remark~\ref{rem:indefinite}).

The positions are listed in a \emph{cyclic order}: we fix an order $c^1,\dots,c^m$ of the orbits and take $c^1_1,c^2_1,\dots,c^m_1,c^1_2,c^2_2,\dots$, skipping an orbit once its positions are used up. For example, with $3$ positions for an orbit $u$ and $5$ for an orbit $v$, the order is $u_1,v_1,u_2,v_2,u_3,v_3,v_4,v_5$. For the singular value decomposition, Householder's algorithm requires the positions to alternate between the two orbits, and this turns tridiagonalisation into bidiagonalisation (Proposition~\ref{prop:golub-kahan} and the remarks after it).

An element of $A$, or a matrix over $A$, is therefore stored as an array indexed by pairs of positions, with entries in the corner spaces. How the entries themselves are represented, as series over the paths of a quiver, is the subject of Chapters~\ref{division-problem} and~\ref{chap:implementation}.

\section{Division and Zero Divisors}
\label{sec:division}
Division is the operation most likely to fail when an algorithm is transported to a new algebra. Let $A$ be a finite-dimensional algebra, and let $A=S\oplus\rad(A)$ be a Wedderburn--Malcev decomposition (Theorem~\ref{thm:wedderburn-malcev}). Since $\rad(A)$ is nilpotent, its elements are ``infinitesimal'', and invertibility is decided by the semisimple part.

\begin{proposition}\label{prop:invertibility}
For $x\in A$ the following are equivalent:
\begin{enumerate}
  \item $x$ is invertible in $A$;
  \item the image of $x$ in $A/\rad(A)$ is invertible;
  \item the left multiplication map $L_x:A\to A$, $y\mapsto xy$, has non-zero determinant.
\end{enumerate}
\end{proposition}
\begin{proof}
(1)$\Rightarrow$(2) is clear. For (2)$\Rightarrow$(1), let $y\in A$ be such that the images of $xy$ and $yx$ in $A/\rad(A)$ equal $1$. Then $xy=1-r$ with $r\in\rad(A)$. Since $r^N=0$ for some $N$, $1-r$ is invertible with inverse $1+r+\cdots+r^{N-1}$, so $x$ has the right inverse $y(1-r)^{-1}$. Similarly $x$ has a left inverse, and therefore $x$ is invertible. For (1)$\Leftrightarrow$(3): if $x$ is invertible then $L_x$ is invertible with inverse $L_{x^{-1}}$. Conversely, if $L_x$ is bijective there is $y$ with $xy=1$; then $L_xL_y=\mathrm{id}$, hence $L_yL_x=\mathrm{id}$ by finite-dimensionality, and evaluating at $1$ gives $yx=1$.
\end{proof}

Proposition~\ref{prop:invertibility} gives an explicit formula for inverses. Let $n=\dim_kA$ and let $\chi_x(t)=\det(tI-L_x)=t^n+c_{n-1}t^{n-1}+\cdots+c_1t+c_0$ be the characteristic polynomial of $L_x$, so that $c_0=(-1)^n\det(L_x)$. Write $\chi_x(t)=t\,p_x(t)+c_0$.

\begin{proposition}\label{prop:inverse-formula}
If $\det(L_x)\neq0$, then
\[
  x^{-1}=\frac{(-1)^{n+1}}{\det(L_x)}\;p_x(x).
\]
\end{proposition}
\begin{proof}
The map $x\mapsto L_x$ is an injective algebra homomorphism $A\to\End_k(A)$, since $L_{xy}=L_xL_y$ and $L_x(1)=x$. By the Cayley--Hamilton theorem $\chi_x(L_x)=0$, hence $\chi_x(x)=0$ in $A$, that is, $x\,p_x(x)=-c_0=(-1)^{n+1}\det(L_x)$. Since $p_x(x)$ commutes with $x$, this gives a two-sided inverse.
\end{proof}

The formula requires only ring operations and one division by a scalar, so it can be implemented for every finite-dimensional algebra once a basis and its structure constants are known. (The sign $(-1)^{n+1}$ equals $-1$ when $\dim_kA$ is even and $+1$ when it is odd.)

\begin{example}[Division of dual numbers]\label{ex:dual-division}
Let $A=\RR[\varepsilon]/(\varepsilon^2)$ and $x=a+b\varepsilon$. In the basis $(1,\varepsilon)$ we have $x\cdot1=a+b\varepsilon$ and $x\cdot\varepsilon=a\varepsilon$, so
\[
  L_x=\begin{pmatrix}a&0\\ b&a\end{pmatrix},\qquad \chi_x(t)=(t-a)^2=t^2-2at+a^2,\qquad \det(L_x)=a^2,\qquad p_x(t)=t-2a .
\]
Here $n=2$, and Proposition~\ref{prop:inverse-formula} gives, for $a\neq0$,
\[
  (a+b\varepsilon)^{-1}=-\frac{x-2a}{a^2}=\frac{a-b\varepsilon}{a^2}=\frac1a-\frac{b}{a^2}\,\varepsilon .
\]
Consequently, for $c+d\varepsilon\in A$,
\[
  \frac{c+d\varepsilon}{a+b\varepsilon}=(c+d\varepsilon)\,\frac{a-b\varepsilon}{a^2}=\frac{c}{a}+\frac{ad-bc}{a^2}\,\varepsilon .
\]
Indeed, $(a+b\varepsilon)\bigl(\tfrac ca+\tfrac{ad-bc}{a^2}\varepsilon\bigr)=c+\tfrac{ad-bc}{a}\varepsilon+\tfrac{bc}{a}\varepsilon=c+d\varepsilon$. If $c+d\varepsilon=f(t)+f'(t)\varepsilon$ and $a+b\varepsilon=g(t)+g'(t)\varepsilon$, the $\varepsilon$-coefficient of the quotient is $(f'g-fg')/g^2$, which is the quotient rule; this is how automatic differentiation differentiates quotients.

If $a=0$, then $x=b\varepsilon$ is a zero divisor and has no inverse. Division of one zero divisor by another is not well defined: the equation $\varepsilon y=\varepsilon$ is solved by $y=1+t\varepsilon$ for every $t\in\RR$. Nevertheless, a computation often points to one value of such a quotient, its \emph{intended value}. For instance, $(\varepsilon b)/(\varepsilon c)$ is meant to be $b/c$, the value obtained by cancelling $\varepsilon$ as if it were invertible; Chapter~\ref{division-problem} makes this precise.
\end{example}

The last observation is not a curiosity. Consider the Jacobi eigenvalue algorithm~\cite{GolubVanLoan2013} applied to the symmetric matrix $H=\bigl(\begin{smallmatrix}1+\varepsilon a&\varepsilon b\\ \varepsilon b&1+\varepsilon c\end{smallmatrix}\bigr)$ over the dual numbers. The rotation angle $\theta$ that diagonalises $H$ is determined by $\tan 2\theta=2h_{12}/(h_{11}-h_{22})=2\varepsilon b/(\varepsilon(a-c))$, a quotient of two zero divisors. The intended value is $\tan2\theta=2b/(a-c)$: with this choice the rotation diagonalises the first-order part $\bigl(\begin{smallmatrix}a&b\\ b&c\end{smallmatrix}\bigr)$ and hence $H$ itself. First-order perturbation theory finds this value directly: it diagonalises the perturbation within each eigenspace of the unperturbed matrix~\cite{Kato1995}, as Proposition~\ref{prop:dual-spectral} does over the dual numbers. It needs to know which eigenvalues coincide, and an iterative algorithm does not know this in advance. Inside such an algorithm, the quotient is the simplest instance of the \emph{division problem}, which is solved in Chapter~\ref{division-problem}. There, $\varepsilon$ is lifted from the dual numbers $k[x]/(x^2)$ to the path algebra $k[x]$ of the loop quiver, where $x$ is no longer nilpotent; $x$ is then made invertible, and the computation is carried out in the field $k((x))$ of Laurent series, where $(xb)/(xc)=b/c$; finally, the result is reduced modulo $x^2$.

\paragraph{Division between idempotents.} In an algebra that is not a division algebra, the useful notion of invertibility is invertibility \emph{between} two idempotents.

\begin{definition}\label{def:invertible-between}
Let $A$ be a ring and $f,g\in A$ idempotents. An element $p\in gAf$ is \emph{invertible from $f$ to $g$} if there is $q\in fAg$ with $qp=f$ and $pq=g$. The element $q$ is then unique, and we write $q=p^{-1}$.
\end{definition}

Uniqueness holds because $q=qg=q(pq')=(qp)q'=fq'=q'$ for any two such elements $q,q'$. For example, in $M_2(k)$ the matrix unit $E_{21}$ is invertible from $E_{11}$ to $E_{22}$, with inverse $E_{12}$, although it is nilpotent as an element of $M_2(k)$. On the other hand, a non-zero element $p$ of $\rad(A)$ is never invertible between idempotents, because $f=qp$ would be an idempotent in the ideal $\rad(A)$, hence nilpotent and equal to $0$, and then $p=pf=0$. Radical elements can therefore only become invertible in a larger ring; Chapter~\ref{division-problem} constructs one.

The division used by the algorithms of this thesis is defined relative to a fixed complete set of orthogonal idempotents $\{e_i\}$ that is permuted by the involution, for example the vertices of a quiver. For $y\in A$ let $f_y$ be the sum of the $e_j$ with $ye_j\neq0$ and $g_y$ the sum of the $e_i$ with $e_iy\neq0$, so that $y\in g_yAf_y$. We write $y^{-1}$ for the inverse of $y$ from $f_y$ to $g_y$ if it exists, and $x/y:=xy^{-1}$. For $y=0$ we have $f_y=g_y=0$, and the definition gives $0^{-1}=0$. For an invertible element $y$ of $A$, $y^{-1}$ is the ordinary inverse, since then $ye_j\neq0$ for all $j$. Since the involution permutes the $e_i$, we have $f_{y^*}=g_y^*$ and $g_{y^*}=f_y^*$, and applying the involution to $y^{-1}y=f_y$ and $yy^{-1}=g_y$ shows that $(y^*)^{-1}=(y^{-1})^*$.

The quotient satisfies $(x/y)\,y=xf_y$, so it solves $zy=x$ exactly when $xf_y=x$, that is, when $x$ vanishes outside the corner on which $y$ is invertible. For zero divisors this can fail. In $\twoR$ with $e_1=(1,0)$ and $e_2=(0,1)$, the element $(c,0)$ with $c\neq0$ has $f=g=e_1$ and the inverse $(c^{-1},0)$, but $(p,q)/(c,0)=(p/c,0)$ does not solve $z\,(c,0)=(p,q)$ when $q\neq0$; relative to the single self-adjoint idempotent $1=e_1+e_2$, the element $(c,0)$ has no inverse at all. The implementation of Chapter~\ref{chap:implementation} therefore checks $(x/y)\,y=x$ for every quotient and reports a breakdown when it fails. Householder reflectors need more: by Proposition~\ref{prop:householder-any-root}, their denominator must be invertible in the corner $fAf$ of its column, where $f$ is the idempotent of the column's position. In the Gram--Schmidt process, on the other hand, $0^{-1}=0$ lets a column that depends on the previous ones, and whose residual is therefore zero, produce a zero column of $Q$; Chapter~\ref{chap:implementation} shows that $Q$ is then a partial isometry.

\section{\texorpdfstring{$*$}{*}-Square Roots}
\label{sec:hermitian-sqrt}
In standard linear algebra the length of a vector is computed as
\[ \|v\| = \Bigl(\sum_{i=1}^n v_i^* v_i\Bigr)^{1/2} = (v^* v)^{1/2}. \]
The map $x \mapsto x^* x$ is the \emph{Hermitian square}, and computing $\|v\|$ requires the inverse operation, a \emph{$*$-square root}: given a self-adjoint $s$, find $x$ with $x^*x=s$. The root $x$ need not be self-adjoint. Norms are needed by Householder reflections, Givens rotations and Gram--Schmidt orthonormalisation, so the existence of $*$-square roots is an important question for our programme. The answer depends strongly on the involution.
\begin{itemize}
  \item For $\CC$, or $M_n(\CC)$ with the conjugate transpose, Hermitian squares are exactly the non-negative real numbers, respectively the positive semidefinite matrices; $*$-square roots of other self-adjoint elements do not exist.
  \item For $\twoR$, the self-adjoint elements are $(c,c)$ with $c\in\RR$, and every such element is a Hermitian square: $(c,c)=z^*z$ for $z=(c,1)$, whatever the sign of $c$.
  \item For the dual numbers with $\varepsilon^*=\varepsilon$, every element is self-adjoint and $(a+b\varepsilon)^*(a+b\varepsilon)=a^2+2ab\varepsilon$; the Hermitian squares are $0$ and the elements $c+d\varepsilon$ with $c>0$. With $\varepsilon^*=-\varepsilon$, the self-adjoint elements are the real numbers, and the Hermitian squares are the non-negative reals.
\end{itemize}
The symplectic involution behaves like $\twoR$: every self-adjoint element is a Hermitian square.

\begin{proposition}\label{prop:symplectic-sqrt}
Let $R$ be a real closed field (for example $\RR$), let $J=\bigl(\begin{smallmatrix}0&I_m\\-I_m&0\end{smallmatrix}\bigr)$, and equip $M_{2m}(R)$ with the symplectic involution $X^\omega=J^{-1}X^TJ$.
\begin{enumerate}
  \item For every $M\in M_m(R)$ the matrix $N=\bigl(\begin{smallmatrix}0&M\\-I&0\end{smallmatrix}\bigr)$, which in Kronecker product notation is $N=\bigl(\begin{smallmatrix}0&1\\0&0\end{smallmatrix}\bigr)\otimes M-\bigl(\begin{smallmatrix}0&0\\1&0\end{smallmatrix}\bigr)\otimes I$, satisfies
  \[ N^\omega N=\begin{pmatrix}M^T&0\\0&M\end{pmatrix}. \]
  \item Every $H\in M_{2m}(R)$ with $H^\omega=H$ is a Hermitian square: $H=Y^\omega Y$ for some $Y\in M_{2m}(R)$.
\end{enumerate}
\end{proposition}
\begin{proof}
(1) We have $N^T=\bigl(\begin{smallmatrix}0&-I\\M^T&0\end{smallmatrix}\bigr)$, hence $N^TJ=\bigl(\begin{smallmatrix}I&0\\0&M^T\end{smallmatrix}\bigr)$ and $N^\omega=J^{-1}N^TJ=\bigl(\begin{smallmatrix}0&-M^T\\I&0\end{smallmatrix}\bigr)$. Multiplying, $N^\omega N=\bigl(\begin{smallmatrix}M^T&0\\0&M\end{smallmatrix}\bigr)$.

(2) The condition $H^\omega=H$ says that $JH$ is skew-symmetric, that is, $H$ is \emph{skew-Hamiltonian}. Every real skew-Hamiltonian matrix is similar, via a symplectic matrix $S$ (one with $S^\omega S=I$), to a block-diagonal matrix $D=\diag(M,M^T)$~\cite{FassbenderMackeyMackeyXu1999,VanLoan1984}. By (1), applied to $M^T$ in place of $M$, there is $N'$ with $N'^\omega N'=D$. Then $H=SDS^{-1}=SN'^\omega N'S^\omega=(N'S^\omega)^\omega(N'S^\omega)$, since $(S^\omega)^\omega=S$. This proves (2) for $R=\RR$. For fixed $m$, statement (2) is a first-order sentence in the language of ordered fields, so it holds over every real closed field by Tarski's transfer principle (Theorem~\ref{thm:tarski-transfer}).
\end{proof}

\paragraph{Roots are not unique.} If $x$ is a $*$-square root of $s$ and $u$ satisfies $u^*u=1$, then $ux$ is another. Conversely, if $x$ is invertible from $f$ to $g$, every other root $x'\in Af$ has the form $ux$ with $u^*u=g$, and $u$ is unitary from $g$ to $g'$ when $x'$ is invertible from $f$ to $g'$.

\begin{lemma}\label{lem:sqrt-unique}
Let $f,g$ be self-adjoint idempotents of a $*$-algebra $A$, let $x\in gAf$ be invertible from $f$ to $g$, and let $x'\in Af$ satisfy $x'^*x'=x^*x$. Then $u=x'x^{-1}$ satisfies $u^*u=g$ and $x'=ux$. If, moreover, $x'$ is invertible from $f$ to a self-adjoint idempotent $g'$, then $u$ is unitary from $g$ to $g'$: $u^*u=g$ and $uu^*=g'$.
\end{lemma}
\begin{proof}
We have $ux=x'x^{-1}x=x'f=x'$ and
\[ u^*u=(x^{-1})^*\,x'^*x'\,x^{-1}=(x^{-1})^*x^*xx^{-1}=(xx^{-1})^*(xx^{-1})=g . \]
If $x'$ is invertible from $f$ to $g'$, then $u\in g'Ag$ is invertible from $g$ to $g'$, with inverse $xx'^{-1}$. Multiplying $u^*u=g$ on the right by $u^{-1}$ gives $u^*=u^*g'=u^{-1}$, hence $uu^*=g'$.
\end{proof}

A $*$-square root is therefore a choice, and an algorithm should depend on it as little as possible. Classically, a Householder reflector maps a column $x$ to a multiple of the first unit column $\hat e_1$ with a prescribed phase, $\gamma=-e^{i\theta}\|x\|$ where $e^{i\theta}$ is the phase of $x_1$. Over a $*$-algebra no phase condition is needed: any $*$-square root of $x^*x$ will do, provided that the denominator of the reflector is invertible. We state this in the generality of position indexing, with self-adjoint idempotents $E_1,\dots,E_m$ and the identity $I=\diag(E_1,\dots,E_m)$.

\begin{proposition}\label{prop:householder-any-root}
Let $A$ be a $*$-algebra, let $E_1,\dots,E_m$ and $f$ be self-adjoint idempotents of $A$, and let $x=(x_1,\dots,x_m)^T$ be a column with $x_i\in E_iAf$. Let $\gamma\in E_1Af$ satisfy $\gamma^*\gamma=x^*x$, put $w=x-\hat e_1\gamma$ and $s=w^*x\in fAf$, and suppose that $s$ is invertible from $f$ to $f$. Then $P=I-ws^{-1}w^*$ satisfies
\[ P^*P=PP^*=I\qquad\text{and}\qquad Px=\hat e_1\gamma . \]
\end{proposition}
\begin{proof}
Since $\gamma^*\gamma=x^*x$,
\[ w^*w=x^*x-x_1^*\gamma-\gamma^*x_1+\gamma^*\gamma=(x^*x-\gamma^*x_1)+(x^*x-x_1^*\gamma)=s+s^*. \]
The inverse $s^{-1}$ lies in $fAf$, and $(s^*)^{-1}=(s^{-1})^*$. Hence $(s^*)^{-1}(w^*w)s^{-1}=(s^*)^{-1}+s^{-1}$, and, since $Iw=w$,
\[ P^*P=I-ws^{-1}w^*-w(s^*)^{-1}w^*+w\bigl((s^*)^{-1}+s^{-1}\bigr)w^*=I . \]
The same computation with the factors in the other order, using $s^{-1}(w^*w)(s^*)^{-1}=(s^*)^{-1}+s^{-1}$, gives $PP^*=I$. Finally, $Px=x-ws^{-1}(w^*x)=x-wf=x-w=\hat e_1\gamma$.
\end{proof}

Over $\CC$, the classical choice of $\gamma$ makes $s$ real and $P$ Hermitian; by the proposition it is needed only for $P$ to be Hermitian, not for $P$ to be unitary. What remains of the classical choice is the condition that $s$ be invertible from $f$ to $f$. Over $\RR$, $\CC$ and $\HH$ with their standard involutions, it holds for every $x\neq0$ if the sign of $\gamma$ is chosen to avoid cancellation in $x_1-\gamma$, as in the classical algorithm; this choice also keeps a nearly reduced column accurate in floating-point arithmetic. In corners with zero divisors, the choice of the root itself can matter. In $\twoR$, with $(p,q)^*=(q,p)$, the column $x=\bigl((1,-1),(\sqrt2,\sqrt2)\bigr)^T$ has $x^*x=1$; the roots $\gamma=\pm(1,1)$ give the zero divisors $s=(0,2)$ and $s=(2,0)$, so they yield no reflector, whereas the root $\gamma=(2,\tfrac12)$ gives the invertible $s=(\tfrac12,3)$.

\paragraph{Norms in the radical.} If $B$ is a $*$-algebra whose radical $R$ satisfies $R^2=0$, then every column $x$ with entries in $R$ has
\[ x^*x=\sum_ix_i^*x_i\in R^2=0, \]
because $R^*=R$ (Lemma~\ref{lem:radical-stable}). Such columns have zero norm. A Householder reflector or a Gram--Schmidt step for them would need an element of the radical, such as $s=w^*x$ when $\gamma$ lies in the radical, to be invertible from $f$ to $f$, which is impossible by Section~\ref{sec:division}. Section~\ref{sec:unification-example} shows that this happens for every column of the matrices that encode the singular value decomposition, and how lifting to a path algebra removes the obstruction.

Proposition~\ref{prop:symplectic-sqrt} gives $*$-square roots in one class of semisimple algebras. For the algebras of this thesis, square roots are computed directly in a path $*$-algebra and its completed localisation, degree by degree (Section~\ref{sec:lambda-sqrt}). In the algebra $K$ of a single vertex (Table~\ref{tab:vertex-types}), $-1$ is a Hermitian square exactly when $K$ is not of signature type: if it is not, $-1=uu^*$ for some $u$ by Lemma~\ref{lem:star-simple}(1), and if $K$ is $\RR$, $\CC$ or $\HH$ with its standard involution, then $y^*y=|y|^2\geq0$ for every $y$.

\section{A First Unification: The SVD as a Spectral Theorem}
\label{sec:unification-example}
We now show, in the simplest non-trivial case, how a classical decomposition arises as the analogue of the Hermitian eigendecomposition over a non-semisimple $*$-algebra, and how the generic algorithm for the latter becomes the classical algorithm for the former. Let $B=\RR Q/R_Q^2$ be the four-dimensional algebra of Example~\ref{ex:path-algebras}(3), with basis $e_u,e_v,a,b$, where $a=e_vae_u$ is an arrow from $u$ to $v$, $b=e_ube_v$ is an arrow from $v$ to $u$, and $ab=ba=0$. Define an involution on $B$ by $e_u^*=e_u$, $e_v^*=e_v$ and $a^*=b$. In the language of Chapter~\ref{chap:quivers}, $B$ is the $*$-algebra of the doubly-signed quiver of Figure~\ref{fig:dsq-examples}(a); here we work with it directly.

A matrix $H\in M_n(B)$ can be written uniquely as $H=e_uX+e_vY+aW+bZ$ with $X,Y,W,Z\in M_n(\RR)$, where $e_uX$ denotes the matrix with entries $x_{ij}e_u$. Since $(aW)^*$ has entries $(w_{ji}a)^*=w_{ji}b$, we have $H^*=e_uX^T+e_vY^T+bW^T+aZ^T$. Hence
\[
  H\in\Herm_n(B)\iff X=X^T,\quad Y=Y^T,\quad Z=W^T .
\]
Similarly, a direct computation using $e_ua=0$, $ae_u=a$, $e_va=a$, $ae_v=0$ (and the analogous relations for $b$) shows that $U=e_uP+e_vR+aS+bT$ is unitary if and only if
\[
  P^TP=I,\qquad R^TR=I,\qquad T=-PS^TR .
\]
Let $U_0=e_uP+e_vR$ denote the semisimple part of $U$. If all entries of $H$ lie in the radical, that is, $H=aW+bW^T$, then every term of $UHU^*$ that involves $S$ or $T$ is a product of two radical elements and vanishes, and
\begin{equation}\label{eq:svd-action}
  UHU^*=U_0HU_0^*=a\,(RWP^T)+b\,(RWP^T)^T .
\end{equation}
Thus \emph{unitary similarity of Hermitian matrices over $B$ with entries in the radical is the same as the action $W\mapsto RWP^T$ of $\Or(n)\times\Or(n)$ on real matrices}. A unitary diagonalisation of $H$, with diagonal entries $a\sigma_i+b\sigma_i$, is precisely a singular value decomposition $W=R^T\diag(\sigma_i)P$. Theorem~\ref{thm:encoding} generalises this computation to every doubly-signed quiver.

When $X$ and $Y$ are non-zero the situation changes. Taking $P=R=I$ and $T=-S^T$, one finds that the $a$-component of $UHU^*$ is $W+SX-YS$. If $X=\diag(x_j)$ and $Y=\diag(y_i)$, then the $(i,j)$ entry of $SX-YS$ is $s_{ij}(x_j-y_i)$, so every entry $w_{ij}$ with $x_j\neq y_i$ can be removed. In general, therefore, the canonical form consists of eigendecompositions of $X$ and $Y$ together with singular value decompositions of the blocks of $W$ that connect equal eigenvalues of $X$ and $Y$. The SVD appears in pure form exactly when $X=Y=cI$ for some $c\in\RR$, that is, when all entries of $H-c\cdot1$ lie in the radical; the case $c=0$ is the case $H=aW+bW^T$ above.

This has an algorithmic consequence. An algorithm such as the Jacobi method, transported to $M_n(B)$, computes rotation angles from quotients of entries of $H$. When the semisimple parts of the relevant entries coincide, as they do in the pure SVD case, these quotients are quotients of zero divisors, such as $(aw_{ij})/(aw_{jj})$, whose intended value is $w_{ij}/w_{jj}$. As in Example~\ref{ex:dual-division}, such quotients are meaningless in $B$.

\paragraph{Householder's algorithm.} The obstruction is even more basic for Householder's algorithm, because the norms vanish. We index $H=aW+bW^T$ by positions (Section~\ref{sec:indexing}), using the self-adjoint idempotents $e_u$ and $e_v$, and we allow $W\in M_{m\times n}(\RR)$ to be rectangular: there are $n$ positions $u_1,\dots,u_n$ and $m$ positions $v_1,\dots,v_m$, and
\[ H_{v_i,u_j}=w_{ij}\,a,\qquad H_{u_j,v_i}=w_{ij}\,b , \]
all other entries being zero. The part of the column of $u_1$ below its diagonal is a column $x$ with entries $w_{i1}a$ at the positions $v_i$, and $x^*x=\sum_iw_{i1}^2\,ba$, which is $0$ in $B$. In the path algebra $\RR Q$, however, $ba\neq0$, and $x^*x=c\,ba$ with $c=\sum_iw_{i1}^2$ has the $*$-square root $\gamma=\sqrt c\,a$, since $a^*a=ba$. Householder's algorithm therefore runs over $\RR Q$ once the elements $s=w^*x$, which are multiples of $ba$ or $ab$, can be inverted. To this end we adjoin inverses of the arrows: elements $a^{-1}\in e_u\Lambda_0e_v$ and $b^{-1}\in e_v\Lambda_0e_u$ with $a^{-1}a=e_u$, $aa^{-1}=e_v$, $b^{-1}b=e_v$ and $bb^{-1}=e_u$, and we extend the involution by $(a^{-1})^*=b^{-1}$. Chapter~\ref{division-problem} constructs the resulting algebra $\Lambda_0$ for every quiver (Section~\ref{sec:localisation}).

\begin{proposition}\label{prop:golub-kahan}
Let $W\in M_{m\times n}(\RR)$ with $m\geq n$, and let $H=aW+bW^T$ be indexed by the positions $u_1,v_1,u_2,v_2,\dots$ in cyclic order. Run Householder tridiagonalisation over $\Lambda_0$, taking for the $*$-square root of $x^*x=c\,ba$ (respectively $c\,ab$) the element $\gamma=\pm\sqrt c\,a$ (respectively $\pm\sqrt c\,b$), with the sign that avoids cancellation, and inverting each denominator $s$ in the corner of its column. Then:
\begin{enumerate}
  \item every step is defined, and every reflector has the form $P=e_uP_u+e_vP_v$, that is, $P$ acts on the positions of $u$ by a real orthogonal matrix $P_u$ and on those of $v$ by $P_v$, and has no other entries;
  \item the reflectors are, in order, the left and right Householder reflectors of the Golub--Kahan bidiagonalisation of $W$~\cite{GolubKahan1965}, with the same choice of signs, except where the part of a column or row beyond its first entry vanishes but the first entry does not: there the generic step multiplies one position by $-1$, whereas the classical algorithm, as implemented in LAPACK, applies no reflector;
  \item the result is $T=U^*HU=aW'+bW'^T$ with $U=e_uU_u+e_vU_v$ and $W'=U_v^TWU_u$ upper bidiagonal, and $W'$ is the Golub--Kahan bidiagonal matrix of $W$ up to the signs of its entries.
\end{enumerate}
\end{proposition}
\begin{proof}
Suppose that, before some step, the current matrix is $a\tilde W+b\tilde W^T$ for a real matrix $\tilde W=(\tilde w_{ij})$ whose first $k-1$ columns and first $k-1$ rows are already in bidiagonal form, and that the step concerns the column of $u_k$. Since $Q$ has no loops, the entries of this column at positions of $u$ vanish, so the part $x$ below the diagonal has the entries $\tilde w_{ik}a$ at the positions $v_i$ with $i\geq k$ and zero entries elsewhere. Hence $x^*x=c\,ba$ with $c=\sum_{i\geq k}\tilde w_{ik}^2$, and $\gamma=\pm\sqrt c\,a$ lies in $e_v\Lambda_0e_u$, which is the corner of the next position, $v_k$, as Proposition~\ref{prop:householder-any-root} requires. If $c=0$, then $w=0$ and $P=I$. Otherwise $w=x-\hat e_{v_k}\gamma$, where $\hat e_{v_k}$ is the unit column of the position $v_k$, has the entries $\hat w_ia$ at the positions $v_i$, where $\hat w=(\tilde w_{kk}\mp\sqrt c,\,\tilde w_{k+1,k},\dots,\tilde w_{mk})^T$ is the classical Householder vector of the column $(\tilde w_{kk},\dots,\tilde w_{mk})^T$, and $s=w^*x=c'\,ba$ with $c'=\hat w^T(\tilde w_{kk},\dots,\tilde w_{mk})^T=c+\sqrt c\,|\tilde w_{kk}|>0$ for the sign that avoids cancellation. In $\Lambda_0$ we have $s^{-1}=c'^{-1}a^{-1}b^{-1}$, and the entry of $ws^{-1}w^*$ at the positions $(v_i,v_j)$ is
\[ \hat w_ia\;c'^{-1}a^{-1}b^{-1}\;\hat w_jb=c'^{-1}\hat w_i\hat w_j\,e_v , \]
while all its other entries vanish. Thus $P$ is the identity at the positions of $u$ and the classical reflector $I-\hat w\hat w^T/c'$ at the positions $v_k,\dots,v_m$. Conjugation by $P$ replaces $\tilde W$ by $P_v\tilde W$, which is the left Householder step of the Golub--Kahan bidiagonalisation on column $k$. A step at the column of $v_k$ is the same computation with $a$ and $b$, $u$ and $v$, and rows and columns exchanged: its reflector acts on the positions $u_{k+1},\dots,u_n$ and replaces $\tilde W$ by $\tilde WP_u$, which is the right Householder step on row $k$. The positions alternate between $u$ and $v$, so the steps come in the order of the Golub--Kahan algorithm, and tridiagonalisation, which treats the columns of the first $m+n-2$ positions, includes all its steps because $m\geq n$. When the entries $\tilde w_{k+1,k},\dots,\tilde w_{mk}$ vanish but $\tilde w_{kk}$ does not, the sign that avoids cancellation gives $\gamma=-\tilde w_{kk}a$, and $P$ multiplies the position $v_k$ by $-1$, whereas the classical algorithm applies no reflector; this proves (1) and (2). All reflectors have the form $e_uP_u+e_vP_v$, hence so does $U$, and $U^*(aW)U=a\,U_v^TWU_u$, which gives (3).
\end{proof}

The proposition shows several things at once. The generic tridiagonalisation, written once for Hermitian matrices, becomes the Golub--Kahan bidiagonalisation without any code for it. The unitary $U$ has the form $e_uU_u+e_vU_v$ with $U_u\in\Or(n)$ and $U_v\in\Or(m)$, so it lies in the group $\Or(n)\times\Or(m)$ that acts on such matrices by Theorem~\ref{thm:encoding}, and the result $T$ has entries in the radical, so $U$ and $T$ lie in $B$ even though the intermediate quantities $s^{-1}$ do not; Theorem~\ref{thm:reduction} makes this precise.

And the order of the positions matters. If the part $x$ of the column of a position of $u$ below the diagonal is non-zero, then $x^*x=c\,ba$ with $c>0$, and no $*$-square root of it lies in the corner $e_u\Lambda_0e_u$: that corner is the ring $\RR[g,g^{-1}]$ with $g=ba=g^*$, in which $\gamma^*\gamma=\gamma^2$, and $c\,g$ is not a square. By Proposition~\ref{prop:householder-any-root}, the next position must therefore belong to $v$, and symmetrically for columns of $v$. For almost every $W$ these parts are non-zero at every step, and the orders that Householder's algorithm accepts are therefore exactly the two cyclic orders, starting with $u$ or with $v$.

With the plain indexing of $M_n(B)$, in which every position carries the idempotent $1$, the corner condition is empty, but the result is then only a tridiagonal matrix. In the experiments of Chapter~\ref{chap:implementation}, the case $m<n$ and the complex and quaternionic singular value decompositions behave in the same way.

With the labels $a$ and $b$ removed, $H$ is the matrix $\bigl(\begin{smallmatrix}0&W^T\\W&0\end{smallmatrix}\bigr)$, listed in the cyclic order $u_1,v_1,u_2,v_2,\dots$ instead of $u_1,\dots,u_n,v_1,\dots,v_m$, that is, after a perfect shuffle of its rows and columns. That this permuted matrix is tridiagonal, with zero diagonal, when $W$ is upper bidiagonal, and that its non-zero eigenvalues are $\pm$ the non-zero singular values of $W$, is classical; it relates the singular values of bidiagonal matrices to symmetric tridiagonal eigenproblems~\cite{GolubKahan1965,DemmelKahan1990}. What the algebra adds is that the labels force this order, through the corner condition of Proposition~\ref{prop:householder-any-root}, and that the same program runs unchanged over the other $*$-algebras of Chapter~\ref{chap:indexing}.

\begin{example}[A $3\times2$ matrix]\label{ex:golub-kahan-numbers}
Let
\[ W=\begin{pmatrix}3&1\\4&2\\0&0.3\end{pmatrix},\qquad H=\begin{pmatrix}0&3b&0&4b&0\\3a&0&a&0&0\\0&b&0&2b&0.3b\\4a&0&2a&0&0\\0&0&0.3a&0&0\end{pmatrix}, \]
with the rows and columns of $H$ indexed by the positions $u_1,v_1,u_2,v_2,v_3$ in this cyclic order. Householder tridiagonalisation treats the columns of the first three positions.

\emph{Column of $u_1$.} Below the diagonal it is $x=(3a,0,4a,0)^T$, at the positions $v_1,u_2,v_2,v_3$, and $x^*x=25\,ba$. The graded root is $5a$, and the sign that avoids cancellation gives $\gamma=-5a$, which lies in the corner of the next position, $v_1$. Then $w=(8a,0,4a,0)^T$, $s=w^*x=40\,ba$ and $s^{-1}=\tfrac1{40}\,a^{-1}b^{-1}$. In $ws^{-1}w^*$ the inverse letters cancel: its entry at the positions $(v_i,v_j)$ is $\tfrac1{40}\hat w_i\hat w_j\,e_v$ with $\hat w=(8,4,0)^T$, and its other entries vanish. So $P$ is the identity at the positions of $u$ and the classical reflector
\[ P_v=I-\tfrac1{40}\hat w\hat w^T=\begin{pmatrix}-0.6&-0.8&0\\-0.8&0.6&0\\0&0&1\end{pmatrix} \]
at those of $v$, and the new matrix $PHP^*$ encodes $P_vW=\bigl(\begin{smallmatrix}-5&-2.2\\0&0.4\\0&0.3\end{smallmatrix}\bigr)$.

\emph{Column of $v_1$.} Below the diagonal it is $(-2.2b,0,0)^T$, at the positions $u_2,v_2,v_3$. Only its first entry is non-zero, which is the exceptional case of Proposition~\ref{prop:golub-kahan}(2). Here $\gamma=2.2\,b$, $s=9.68\,ab$, and $P$ multiplies the position $u_2$ by $-1$, so that the encoded matrix becomes $\bigl(\begin{smallmatrix}-5&2.2\\0&-0.4\\0&-0.3\end{smallmatrix}\bigr)$.

\emph{Column of $u_2$.} Below the diagonal it is $(-0.4a,-0.3a)^T$, at the positions $v_2,v_3$, with $x^*x=0.25\,ba$ and $\gamma=0.5\,a$. Then $s=0.45\,ba$, and $P$ acts on $v_2,v_3$ by the classical reflector $\bigl(\begin{smallmatrix}-0.8&-0.6\\-0.6&0.8\end{smallmatrix}\bigr)$.

The result is $T=U^*HU$ with $U=e_uU_u+e_vU_v$, and $T$ encodes the upper bidiagonal matrix
\[ W'=U_v^TWU_u=\begin{pmatrix}-5&2.2\\0&0.5\\0&0\end{pmatrix}, \]
whose singular values $5.4663\ldots$ and $0.4573\ldots$ are those of $W$. The bidiagonalisation of Golub and Kahan, with the sign conventions of LAPACK, gives the entries $-5$, $-2.2$ and $-0.5$; the two results differ by the sign of the second column, as Proposition~\ref{prop:golub-kahan} predicts. Neither $U$ nor $T$ contains an inverse letter, so $H=UTU^*$ holds over the four-dimensional algebra $B$ itself (Theorem~\ref{thm:reduction}). Listing~\ref{lst:session} runs this example with the library of Chapter~\ref{chap:implementation}.
\end{example}
%======================================================================
%   CHAPTER 4: INFINITESIMALS AND TRANSFER
%======================================================================
\chapter{Infinitesimals and Transfer Principles}
\label{chap:infinitesimals}

This chapter treats ordered fields with infinitesimals as scalars. Over a real closed field, the seven operations of Section~\ref{sec:proposal} and the comparisons satisfy the same first-order statements as over $\RR$ (Theorem~\ref{thm:tarski-transfer}), and running an algorithm once over the field of algebraic Puiseux series in an infinitesimal $\varepsilon$, then taking standard parts of its finite outputs, gives their limits as a real parameter $\varepsilon$ tends to zero (Proposition~\ref{prop:autolimit}). The Bruhat decomposition of Theorem~\ref{thm:B}(4) is computed in this way (Theorem~\ref{thm:bruhat}).

\section{Motivation and Overview}
\label{sec:inf-motivation}
Infinitesimals enter this thesis in three ways.
\begin{enumerate}
  \item \emph{Automatic limits.} Section~\ref{subsec:qr-lu} showed that Gaussian elimination is the limit, as a weight $\varepsilon$ tends to $0$, of Householder QR with respect to a weighted inner product. Running an algorithm once over a field containing a formal infinitesimal $\varepsilon$, and then discarding infinitesimal parts, computes such limits mechanically (Proposition~\ref{prop:autolimit}). This is the auto-limit of Section~\ref{subsec:qr-lu}, named by analogy with automatic differentiation. As with derivatives, there are other ways to find these limits: analysis by hand, as in the proof of Proposition~\ref{prop:gs-limit}; computer algebra, with $\varepsilon$ kept as a symbol; and, less reliably, running the algorithm numerically at a small value of $\varepsilon$.
  \item \emph{Perturbation theory.} A matrix whose entries are power series in $\varepsilon$ describes a family of matrices depending on a parameter. Decompositions over power series rings describe how decompositions vary with the parameter, and truncating modulo $\varepsilon^2$ recovers decompositions over the dual numbers~\cite{Gutin2022DualNumbers}.
  \item \emph{Division.} The radical of a finite-dimensional algebra consists of nilpotent elements, which behave like infinitesimals but cannot be inverted. Chapter~\ref{division-problem} replaces them by \emph{invertible} ones, by inverting the arrows of a path algebra and completing. For the quiver with one loop, this gives the field of Laurent series (Section~\ref{sec:pythagorean}), in which the auto-limits of this chapter can be computed.
\end{enumerate}

To use infinitesimals rigorously we need two things. First, an ordered field $\KK\supsetneq\RR$ containing infinitesimals in which the operations required by numerical algorithms, in particular square roots of positive elements, are available. Second, a \emph{transfer principle}, guaranteeing that the theorems which justify an algorithm over $\RR$ (for instance, the existence of an orthonormal eigenbasis for a symmetric matrix) remain true over $\KK$. Without a transfer principle every such theorem would have to be re-proved for each new field.

There are two well-known ways of obtaining these. \emph{Robinson's nonstandard analysis}~\cite{Robinson1966} produces a field ${}^*\RR$ with infinitesimals to which \emph{every} first-order statement about $\RR$ transfers, but ${}^*\RR$ is obtained non-constructively and cannot be computed with. \emph{Tarski's theorem} on real closed fields~\cite{Tarski1951} gives a transfer principle only for statements about polynomial equations and inequalities, but it applies to explicit, computable fields such as the field of real algebraic Puiseux series. All the statements we need are of this restricted kind, so the results of this thesis rely only on Tarski's theorem. We nevertheless describe Robinson's approach (Section~\ref{subsec:robinson}), both because it is the best-known theory of infinitesimals and because it clarifies the relationship between nilpotent infinitesimals, such as those of the dual numbers, and invertible ones (Section~\ref{sec:standard-part}). In short, \emph{nonstandard analysis in Robinson's sense is not needed for the results of this thesis}; what we use instead is the constructive transfer principle of Theorem~\ref{thm:tarski-transfer}. The statements that we transfer, such as the spectral theorem, could also be proved directly over an arbitrary real closed field.

The chapter is organised as follows. Section~\ref{sec:rcf} recalls real closed fields. Section~\ref{sec:transfer} states the two transfer principles precisely. Section~\ref{sec:puiseux} introduces Puiseux series, the fields with infinitesimals that we use. Section~\ref{sec:standard-part} defines the standard part map, proves that it computes limits of algorithms, and explains how nilpotent infinitesimals arise as quotients. Section~\ref{sec:pythagorean} shows that fields of Laurent series, which are the completed localisations of the quiver with one loop, are Pythagorean, and Section~\ref{sec:bruhat} uses them to compute the Bruhat decomposition as an auto-limit of the Gram--Schmidt process. Sections~\ref{sec:power-series} and~\ref{sec:parity} apply these tools to prove spectral theorems and singular value decompositions over formal power series rings. Section~\ref{sec:shortest-paths} shows that Dijkstra's shortest-path algorithm is the leading-order part of Gaussian elimination over Puiseux series, and Section~\ref{sec:chain} assembles the chain of limits that leads from the Hermitian eigendecomposition, through the singular value decomposition and the QR factorisation, to the Bruhat and $LU$ factorisations and to Dijkstra's algorithm.

\section{Real Closed Fields}
\label{sec:rcf}
\begin{definition}\label{def:rcf}
A field $R$ is \emph{real closed} if it is an ordered field in which every positive element has a square root and every polynomial of odd degree has a root.
\end{definition}

By the Artin--Schreier theory~\cite[\S1.2]{BochnakCosteRoy1998}, a field $R$ is real closed if and only if any of the following equivalent conditions holds:
\begin{itemize}[leftmargin=1.25em]
  \item $R$ is not algebraically closed, but $R(\sqrt{-1})$ is;
  \item $R$ admits an ordering, and no proper algebraic extension of $R$ admits an ordering extending it;
  \item $R$ admits an ordering, and every polynomial over $R$ factors into linear and quadratic factors.
\end{itemize}
A real closed field has a unique ordering, in which the positive elements are exactly the non-zero squares. Examples are $\RR$, the field of real algebraic numbers, and the field of real Puiseux series (Section~\ref{sec:puiseux}). The theorem of Frobenius extends to real closed fields: the finite-dimensional division algebras over a real closed field $R$ are $R$, $R(\sqrt{-1})$ and the quaternions over $R$.

Two weaker notions will also be useful. An ordered field is \emph{Euclidean} if every positive element is a square, and a field is \emph{Pythagorean} if every sum of squares is a square. Real closed fields are Euclidean, and Euclidean fields are Pythagorean.

The properties of real closed fields used in practice are these: square roots of positive elements exist, which is needed for orthogonalisation, polar decompositions and SVDs; every polynomial splits into linear and quadratic factors, which is needed for eigenvalue problems; and the transfer principle of the next section holds.

\section{First-Order Logic and Transfer Principles}
\label{sec:transfer}

\subsection{First-Order Languages}
\label{subsec:fol}
We recall only what is needed; see~\cite{Marker2002} for details. A \emph{signature} is a set of constant, function and relation symbols. The signature of \emph{ordered rings} is $\mathcal L_{\mathrm{or}}=\{0,1,+,-,\cdot,<\}$. \emph{Terms} are built from variables and constants using the function symbols; in $\mathcal L_{\mathrm{or}}$ they are polynomials with integer coefficients. \emph{Formulas} are built from equations and relations between terms using $\neg,\wedge,\vee,\Rightarrow$ and the quantifiers $\forall x$ and $\exists x$, where $x$ ranges over \emph{elements} of the structure. A \emph{sentence} is a formula without free variables. A structure in which a sentence is true is called a \emph{model} of it.

The restriction that quantifiers range over elements, and not over subsets, functions or natural numbers, is what ``first-order'' means. Statements about matrices of a \emph{fixed} size are first-order, because a matrix is a tuple of elements. For example, the spectral theorem for symmetric $2\times 2$ matrices over an ordered field is the sentence
\begin{multline*}
\forall a\,\forall b\,\forall c\;\exists p\,\exists q\,\exists\lambda\,\exists\mu\;\bigl(p^2+q^2=1\;\wedge\; a=p^2\lambda+q^2\mu\\ \wedge\; b=pq(\lambda-\mu)\;\wedge\; c=q^2\lambda+p^2\mu\bigr),
\end{multline*}
which says that $\bigl(\begin{smallmatrix}a&b\\b&c\end{smallmatrix}\bigr)=R\,\diag(\lambda,\mu)\,R^T$ for the rotation $R=\bigl(\begin{smallmatrix}p&-q\\q&p\end{smallmatrix}\bigr)$. For each $n$ there is a similar sentence $\varphi_n$ expressing the spectral theorem for $n\times n$ symmetric matrices. The statement ``$\varphi_n$ holds for all $n$'' is not first-order, so transfer principles are applied to each $\varphi_n$ separately.

\subsection{The Tarski--Seidenberg Transfer Principle}
\label{subsec:tarski}
\begin{theorem}[Tarski]\label{thm:tarski-transfer}
The theory of real closed fields admits quantifier elimination in the signature $\mathcal L_{\mathrm{or}}$: every formula is equivalent, in all real closed fields, to a quantifier-free formula with the same free variables. Consequently:
\begin{enumerate}
  \item \emph{(Transfer)} A sentence of $\mathcal L_{\mathrm{or}}$ holds in one real closed field if and only if it holds in every real closed field.
  \item \emph{(Transfer with parameters)} If $R\subseteq R'$ are real closed fields and $\varphi(x_1,\dots,x_m)$ is a formula of $\mathcal L_{\mathrm{or}}$, then for all $a_1,\dots,a_m\in R$, $\varphi(a_1,\dots,a_m)$ holds in $R$ if and only if it holds in $R'$.
\end{enumerate}
\end{theorem}
See~\cite{Tarski1951},~\cite[Proposition~5.2.3]{BochnakCosteRoy1998} and~\cite[Corollary~3.3.16]{Marker2002}. Quantifier elimination is effective: there are algorithms that compute the quantifier-free equivalent of a given formula~\cite{BasuPollackRoy2006}. In this sense the transfer principle is constructive. We use it as follows: to show that a theorem about $n\times n$ matrices, true over $\RR$, also holds over a real closed field $R$ containing infinitesimals, it suffices to write the theorem as a first-order sentence in $\mathcal L_{\mathrm{or}}$.

\subsection{Robinson's Nonstandard Analysis}
\label{subsec:robinson}
Robinson's nonstandard analysis rests on the compactness theorem of first-order logic, which together with G\"odel's completeness theorem is one of the two fundamental theorems of the subject.

\begin{theorem}[Compactness]\label{thm:compactness}
A set of first-order sentences has a model if and only if each of its finite subsets has a model.
\end{theorem}

Consider the signature $\Sigma$ that contains a function symbol for every function $\RR^m\to\RR^n$ (including constants, as functions $\RR^0\to\RR$) and a relation symbol for every subset of $\RR^n$, for all $m,n$. This signature is uncountable, which is allowed. Let $T$ be the set of all $\Sigma$-sentences that are true in $\RR$. Add a new constant symbol $\delta$ and the sentences $0<\delta<1/m$ for every positive integer $m$. Every finite subset of the enlarged set of sentences has a model, namely $\RR$ itself with $\delta$ interpreted as a sufficiently small positive real number. By compactness, the whole set has a model ${}^*\RR$. In ${}^*\RR$ the element $\delta$ is positive but smaller than every $1/m$, so it is a non-zero infinitesimal, and ${}^*\RR\neq\RR$.

\begin{theorem}[Transfer principle of nonstandard analysis]\label{thm:robinson-transfer}
There is an ordered field ${}^*\RR\supsetneq\RR$ in which every function and relation on $\RR$ has a natural extension, such that a $\Sigma$-sentence with parameters from $\RR$ holds in $\RR$ if and only if it holds in ${}^*\RR$.
\end{theorem}
Its language is much richer than that of Theorem~\ref{thm:tarski-transfer}: it applies to all functions (exponentials, trigonometric functions, arbitrary maps), not only to polynomials. On the other hand, it concerns the single extension $\RR\subset{}^*\RR$, whereas Theorem~\ref{thm:tarski-transfer} applies to every pair of real closed fields, including the fields of Puiseux series used below. Its main use in mathematics is to replace arguments involving many nested quantifiers, the ``epsilon management'' of analysis, by arguments with infinitesimals~\cite{Tao2007}. The construction is, however, non-constructive: it requires the compactness theorem (or an ultrafilter), and there is no explicit, computable description of the elements of ${}^*\RR$.

\subsection{Comparison of the Two Approaches}
\label{subsec:comparison}
Each approach is strong in one respect and weak in another.
\begin{itemize}
    \item Robinson's approach preserves \emph{every} first-order property of the real numbers, but it does not provide a constructive or computable field of infinitesimals.
    \item The approach based on Tarski's theorem transfers only properties expressible with polynomial equations and inequalities. In exchange it applies to concrete, computable fields with infinitesimals, such as the real algebraic Puiseux series. Moreover, real algebraic Puiseux series converge, so they can be interpreted as functions of a real parameter; this suggests borrowing ideas from complex analysis, such as the substitution of a root of unity in the proof of Theorem~\ref{thm:spectral-power-series}, although that proof is purely algebraic (Section~\ref{sec:power-series}).
\end{itemize}
For extending linear algebra to the perturbative setting, the second approach provides exactly the right balance, and we use it throughout. All arithmetic operations of numerical linear algebra, as well as square roots, are semialgebraic, so the restriction to polynomial statements is harmless for our purposes.

\section{Puiseux Series}
\label{sec:puiseux}
Let $k$ be a field of characteristic $0$ and $t$ an indeterminate.
\begin{itemize}[leftmargin=1.25em]
  \item The field of \emph{formal Laurent series} is $k((t))$, whose elements are series $\sum_{j\geq j_0}c_jt^j$ with $j_0\in\ZZ$; the ring of \emph{formal power series} $k[[t]]$ consists of the series with $j_0\geq0$.
  \item The field of \emph{formal Puiseux series} is $k\{\{t\}\}=\bigcup_{N\geq1}k((t^{1/N}))$, whose elements are Laurent series in a root $t^{1/N}$ of $t$.
  \item If $k=\RR$ or $k=\CC$, a Puiseux series is \emph{convergent} if it converges on some punctured disc $0<|t|<\rho$; the convergent Puiseux series form a subfield.
\end{itemize}
If $k$ is real closed, $k\{\{t\}\}$ is ordered by the sign of the leading coefficient; in this ordering $t$ is a positive infinitesimal.

\begin{theorem}[Newton--Puiseux]\label{thm:newton-puiseux}
If $k$ is algebraically closed of characteristic $0$, then $k\{\{t\}\}$ is algebraically closed. If $k$ is real closed, then $k\{\{t\}\}$ is real closed. The same holds for convergent Puiseux series over $\CC$ and $\RR$ respectively.
\end{theorem}
See~\cite[\S2.6]{BasuPollackRoy2006} and~\cite{BochnakCosteRoy1998}. The real closure of the rational function field $\RR(t)$, ordered so that $t$ is a positive infinitesimal, is the field $\RR\langle t\rangle$ of \emph{algebraic} Puiseux series, consisting of the Puiseux series that are algebraic over $\RR(t)$. These series converge, and $\RR\langle t\rangle$ is isomorphic to the field of germs at $0^+$ of continuous semialgebraic functions on intervals $(0,\rho)$~\cite[\S2.6 and Chapter~3]{BasuPollackRoy2006}; its elements can be represented exactly and manipulated algorithmically. More general real closed fields with infinitesimals are the Hahn fields $k((t^{\QQ}))$ of series with well-ordered support in $\QQ$~\cite{Ribenboim1999}, which contain the Puiseux series.

Formal Puiseux series supply a clean algebraic, ordered and real closed base field; convergent Puiseux series provide the analytic counterpart, in which each element is a function of the real parameter $t$.

\section{Infinitesimals and the Standard Part}
\label{sec:standard-part}

\subsection{Finite Elements, Infinitesimals and Asymptotic Notation}
Let $F$ be an ordered field containing $\RR$. An element $x\in F$ is \emph{finite} if $|x|<m$ for some integer $m$, and \emph{infinitesimal} if $|x|<1/m$ for every positive integer $m$. We write $F_{\mathrm{fin}}$ for the set of finite elements and $\mathfrak m$ for the set of infinitesimals.

\begin{proposition}\label{prop:standard-part}
$F_{\mathrm{fin}}$ is a subring of $F$, $\mathfrak m$ is its unique maximal ideal, and every $x\in F_{\mathrm{fin}}$ is infinitely close to a unique real number $\st(x)$, the \emph{standard part} of $x$. The map $\st:F_{\mathrm{fin}}\to\RR$ is a surjective ring homomorphism with kernel $\mathfrak m$.
\end{proposition}
\begin{proof}
Closure of $F_{\mathrm{fin}}$ and $\mathfrak m$ under the ring operations follows from the corresponding facts for bounds. If $x\in F_{\mathrm{fin}}\setminus\mathfrak m$ then $1/m<|x|<m'$ for some integers, so $1/x$ is finite; hence every element outside $\mathfrak m$ is a unit and $\mathfrak m$ is the unique maximal ideal. For $x\in F_{\mathrm{fin}}$, the set $\{r\in\RR:r<x\}$ is non-empty and bounded above, so it has a supremum $s\in\RR$ by completeness of $\RR$; for every $m$, $s-1/m<x\leq s+1/m$, so $x-s\in\mathfrak m$. Uniqueness holds because two distinct reals are not infinitely close. That $\st$ is a homomorphism follows since $\mathfrak m$ is an ideal.
\end{proof}

Asymptotic notation has a precise meaning in this setting: for $x\in F$ we set
\[ O(x)=x\,F_{\mathrm{fin}},\qquad o(x)=x\,\mathfrak m, \]
so that $O(1)=F_{\mathrm{fin}}$ and $o(1)=\mathfrak m$. These are $F_{\mathrm{fin}}$-submodules of $F$, and the usual rules such as $O(x)O(y)=O(xy)$ and $O(x)+o(x)=O(x)$ hold.

\subsection{Limits of Algorithms as Standard Parts}
We now justify the ``auto-limit'' of Section~\ref{subsec:qr-lu}. By an \emph{algorithm over ordered fields} we mean a finite sequence of instructions, each of which either applies one of the operations $+,-,\times$, division by a non-zero element, or the square root of a non-negative element to previously computed quantities, or branches according to the sign of a previously computed quantity. The algorithm may loop, provided that each run on the inputs under consideration executes finitely many instructions.

\begin{proposition}[Auto-limit]\label{prop:autolimit}
Let $\mathcal A$ be an algorithm over ordered fields and let its inputs $x_1(\varepsilon),\dots,x_m(\varepsilon)$ be elements of the field $\RR\langle\varepsilon\rangle$ of algebraic Puiseux series, for example polynomials in $\varepsilon$. Suppose that $\mathcal A$, executed over $\RR\langle\varepsilon\rangle$, terminates without error and produces outputs $y_1,\dots,y_r\in\RR\langle\varepsilon\rangle$. Then there is $\varepsilon_0>0$ such that for every real $\varepsilon'\in(0,\varepsilon_0)$, executing $\mathcal A$ in real arithmetic on the inputs $x_1(\varepsilon'),\dots,x_m(\varepsilon')$ follows the same branches and produces the outputs $y_1(\varepsilon'),\dots,y_r(\varepsilon')$. In particular, if $y_j$ is finite, then
\[ \lim_{\varepsilon'\to0^+}y_j(\varepsilon')=\st(y_j). \]
\end{proposition}
\begin{proof}
Identify $\RR\langle\varepsilon\rangle$ with the field of germs at $0^+$ of continuous semialgebraic functions (Section~\ref{sec:puiseux}). The field operations and square roots of germs are computed pointwise, and a non-zero germ has constant sign on some interval $(0,\rho)$. Each of the finitely many executed instructions therefore agrees with its pointwise counterpart on some interval $(0,\rho_i)$, and each branch is decided in the same way there; take $\varepsilon_0=\min_i\rho_i$. Finally, a finite element $y=\sum_{q\geq0}c_q\varepsilon^q$ of $\RR\langle\varepsilon\rangle$ is a convergent series without negative exponents, so $y(\varepsilon')\to c_0=\st(y)$.
\end{proof}

The proposition explains the transition from limits to infinitesimals in Section~\ref{subsec:qr-lu}: instead of analysing the limit of an algorithm as a parameter tends to zero, we execute the algorithm once, exactly, over $\RR\langle\varepsilon\rangle$ and take standard parts. Two remarks are in order. First, the limit of the outputs of an algorithm need not equal the output of the algorithm applied to the limiting input; the latter may not even be defined, as when Householder QR with respect to $\langle\cdot,\cdot\rangle_\varepsilon$ is applied at $\varepsilon=0$, where the ``inner product'' is degenerate. Second, Theorem~\ref{thm:tarski-transfer} guarantees that the algorithm is correct over $\RR\langle\varepsilon\rangle$ whenever its correctness over $\RR$ is a first-order statement, for each fixed matrix size.

\subsection{Nilpotent Infinitesimals as Quotients}
Theories of \emph{nilpotent} infinitesimals, such as the dual numbers, make it easy to automate the computation of derivatives by operator overloading~\cite{GriewankWalther2008}, and they support intuitive geometric derivations of formulas for areas, volumes and lengths. They are also the basis of Synthetic Differential Geometry~\cite{Kock2006}, which develops calculus with nilpotent infinitesimals using intuitionistic logic. The following proposition shows that, for the purposes of this thesis, nilpotent infinitesimals can be obtained from invertible ones by passing to a quotient, so that no separate logical framework is required.

\begin{proposition}\label{prop:dual-quotient}
Let $F\supseteq\RR$ be an ordered field and $\delta\in F$ a positive infinitesimal. Then $A_\delta:=\RR+O(\delta)$ is a subring of $F_{\mathrm{fin}}$, $o(\delta)$ is an ideal of $A_\delta$, and the map $x+a\delta\mapsto x+\st(a)\,\varepsilon$ ($x\in\RR$, $a\in F_{\mathrm{fin}}$) induces an isomorphism
\[ A_\delta/o(\delta)\;\cong\;\RR[\varepsilon]/(\varepsilon^2). \]
\end{proposition}
\begin{proof}
$A_\delta$ is closed under multiplication because $(x+a\delta)(x'+a'\delta)=xx'+(xa'+x'a+aa'\delta)\delta$. The map is well defined: if $x+a\delta=x'+a'\delta$ then $x-x'\in O(\delta)\cap\RR=\{0\}$ and $a=a'$. It is a ring homomorphism because $\delta^2=\delta\cdot\delta\in o(\delta)$ maps to $\varepsilon^2=0$; it is surjective; and its kernel consists of the elements $a\delta$ with $\st(a)=0$, that is, of $o(\delta)$.
\end{proof}

If $f:\RR\to\RR$ is differentiable at $x\in\RR$ and $f$ extends to $F$ (for instance, $f$ semialgebraic and $F$ real closed, or $F={}^*\RR$), then $f(x+\delta)=f(x)+f'(x)\delta+o(\delta)$ for every infinitesimal $\delta\neq0$, so $f(x+\delta)\in A_\delta$ and its image in the dual numbers is $f(x)+f'(x)\varepsilon$. This is precisely the value computed by forward-mode automatic differentiation. In this sense the dual numbers, and similarly the algebras $\RR[\varepsilon]/(\varepsilon^{m})$ with $A_\delta=\RR+\RR\delta+\cdots+O(\delta^{m-1})$, are absorbed into the theory of ordered fields with infinitesimals. The same idea, reversed, underlies the completed localisation of Chapter~\ref{division-problem}: a nilpotent element is lifted to an invertible infinitesimal, the computation is carried out in a ring in which it is invertible (for the quiver with one loop, a field of Laurent series), and the result is reduced again.

\begin{remark}[Historical note]
Infinitesimals were used informally by Leibniz, Euler and others. In the nineteenth century Cauchy, Bolzano and above all Weierstrass based analysis on the notion of limit, formulated with inequalities and nested quantifiers. Robinson showed in the 1960s that a consistent theory of infinitesimals can be obtained from model theory~\cite{Robinson1966}. Fields such as $\RR(\varepsilon)$, which contain infinitesimals but support little analysis, had long been known; the point of the transfer principles above is that they allow statements to pass between $\RR$ and its extension in both directions.
\end{remark}
\section{Pythagorean Fields and Laurent Series}
\label{sec:pythagorean}
The field of Laurent series is not real closed ($t$ has no square root in $k((t))$), but it has the property needed for Gram--Schmidt orthogonalisation.

\begin{proposition}\label{prop:laurent-pythagorean}
Let $k$ be a real closed field. Then the field $k((t))$ of formal Laurent series, ordered by the sign of the leading coefficient, is Pythagorean. A sum of squares $s=\sum_i x_i^2$ of elements of $k((t))$ has even valuation, and $s$ is zero only if all $x_i$ are. If $x\in k((t))^n$ satisfies $\sum_ix_i^2=1$, then $x\in k[[t]]^n$. The same statements hold for convergent Laurent series over $\RR$.
\end{proposition}
\begin{proof}
Let $v$ be the smallest valuation among the non-zero $x_i$, and let $c_i$ be the coefficient of $t^v$ in $x_i$. The coefficient of $t^{2v}$ in $s$ is $\sum_ic_i^2>0$, and all lower coefficients vanish. Hence $s=c\,t^{2v}(1+h)$ with $c>0$ and $h\in t\,k[[t]]$. Then $s=(\sqrt c\,t^v(1+h)^{1/2})^2$, where $(1+h)^{1/2}$ is given by the binomial series, which converges $t$-adically, and converges analytically if $h$ does. If $\sum_ix_i^2=1$ then $2v=0$ for the smallest valuation $v$, so no $x_i$ has a negative power of $t$.
\end{proof}

Thus over $k((t))$ the norm $\|x\|=(\sum_ix_i^2)^{1/2}$ of a vector is defined, it vanishes only for $x=0$, and normalised vectors have entries in $k[[t]]$. Consequently the Gram--Schmidt process can be applied without restriction over $k((t))$.

The field $\RR((t))$ is the completed localisation (Definition~\ref{def:completion}, with $k=\RR$) of the path algebra $\RR[t]$ of the quiver with one loop $t$: the localisation adjoins $t^{-1}$, and the completion admits series in $t$. With $t^*=t$, it is the real part of the completed localisation of the doubly-signed quiver with one vertex of type $(\RR,+1)$ and one loop fixed by the transpose with sign $+1$ (Section~\ref{sec:completion}). So the computations of this chapter in $\RR((t))$ are computations in $\Lambda$, and the implementation of Chapter~\ref{chap:implementation} performs them there. Since $\RR((t))$ is Pythagorean but not real closed, an algorithm that needs square roots of elements of odd order, such as $t$ itself, must be run over a larger field, for example the Puiseux series; in the algorithms below this is avoided by writing the parameter as $\varepsilon=t^2$.

\section{The Bruhat Decomposition as an Auto-Limit}
\label{sec:bruhat}
Proposition~\ref{prop:gs-limit} showed that the Gram--Schmidt process, with respect to the weighted inner product \eqref{eq:weighted-ip}, tends to Gaussian elimination when the leading principal minors of the matrix are non-zero. Without this assumption the limit is a Bruhat decomposition. Every $M\in\GL_n(\RR)$ can be written as $M=LwU$, with $L$ lower triangular, $U$ upper triangular and $w$ a permutation matrix, and $w$ is determined by $M$; we write $we_j=e_{w(j)}$. This is the Bruhat decomposition $\GL_n=\bigsqcup_wB_-wB$, where $B$ and $B_-$ are the upper and lower triangular groups, which follows from the decomposition $\GL_n=\bigsqcup_wBwB$~\cite{Humphreys1975} applied to $w_0M$, where $w_0$ reverses the order of the coordinates and $w_0B_-=Bw_0$. It can be computed directly by Gaussian elimination: working through the columns from left to right, take as pivot in each column its first non-zero entry in a row not yet used as a pivot row, and use it to eliminate the entries below it; the pivot of column~$j$ then lies in row~$w(j)$. We now obtain the decomposition in another way, as an auto-limit (Proposition~\ref{prop:autolimit}), in fields of Laurent series.

\begin{theorem}\label{thm:bruhat}
Let $M\in\GL_n(\RR)$, and let $M=LwU$ be a Bruhat decomposition. Apply the Gram--Schmidt recurrence \eqref{eq:gs}, with the inner product $\langle x,y\rangle_\varepsilon=\sum_i\varepsilon^{i-1}x_iy_i$, to the columns of $M$ over the field $\RR((\varepsilon))$, in which $\varepsilon$ is a positive infinitesimal. This gives vectors $u_1,\dots,u_n$ and a unit upper triangular matrix $R=(r_{kj})$ of coefficients with $M=[u_1\cdots u_n]\,R$.
\begin{enumerate}
  \item $\langle u_j,u_j\rangle_\varepsilon=\varepsilon^{w(j)-1}\bigl(c_j+O(\varepsilon)\bigr)$ with $c_j=(L_{w(j)w(j)}U_{jj})^2>0$. In particular, the orders of the norms determine $w$, and $c_j$ depends only on $M$.
  \item The entries of the $u_j$ and of $R$ are finite. Put $P'=\st[u_1\cdots u_n]$ and $U'=\st(R)$. Then $M=P'U'$, $U'$ is unit upper triangular, and $L'=P'w^{-1}$ is lower triangular and invertible. Thus $M=L'wU'$ is a Bruhat decomposition.
  \item Let $t=\sqrt\varepsilon$, and let $\diag(1,t,\dots,t^{n-1})\,M=\hat Q\hat R$ be the QR factorisation over $\RR((t))$, with positive diagonal in $\hat R$. Then $\hat Q$ has finite entries, and $\st(\hat Q)$ is the signed permutation matrix whose column $j$ is $\sign(L'_{w(j)w(j)})\,e_{w(j)}$.
\end{enumerate}
By Proposition~\ref{prop:autolimit}, applied in the field of algebraic Puiseux series, which contains all the quantities computed here, these standard parts are the limits as $\varepsilon\to0^+$ of the same computations in real arithmetic.
\end{theorem}

\begin{proof}
For sets $S,T$ of indices, write $M_{S,T}$ for the corresponding submatrix of $M$, and put $[j]=\{1,\dots,j\}$, $W_j=\{w(1),\dots,w(j)\}$ and $\sigma_j=\sum_{k\leq j}(w(k)-1)$. Since $U$ is upper triangular, the first $j$ columns of $M$ are the columns $Le_{w(1)},\dots,Le_{w(j)}$ of $L$, multiplied on the right by the leading $j\times j$ block $U_{[j]}$ of $U$. So for $|S|=j$,
\[ \det M_{S,[j]}=\pm\det L_{S,W_j}\,\det U_{[j]} . \]
If $s_1<\dots<s_j$ and $q_1<\dots<q_j$ are the elements of $S$ and $W_j$, then $\det L_{S,W_j}=0$ unless $s_k\geq q_k$ for all $k$: otherwise, if $s_m<q_m$, the rows $s_1,\dots,s_m$ of $L_{S,W_j}$ vanish in the columns $q_m,\dots,q_j$, because $L$ is lower triangular, and a $j\times j$ matrix with a zero block of size $m\times(j-m+1)$ is singular. Moreover $\det L_{W_j,W_j}=\prod_{k\leq j}L_{w(k)w(k)}\neq0$. Let $G=M^TD_\varepsilon M$ be the Gram matrix of the columns of $M$, where $D_\varepsilon=\diag(1,\varepsilon,\dots,\varepsilon^{n-1})$, and let $G_j$ be its leading $j\times j$ submatrix. By the Cauchy--Binet formula,
\[ \det G_j=\sum_{|S|=j}\varepsilon^{\sum_{s\in S}(s-1)}\bigl(\det M_{S,[j]}\bigr)^2 . \]
Every set $S$ that contributes satisfies $s_k\geq q_k$ for all $k$, so its exponent is at least $\sigma_j$, with equality only for $S=W_j$. Hence $\det G_j=\varepsilon^{\sigma_j}(\gamma_j+O(\varepsilon))$ with $\gamma_j=\prod_{k\leq j}(L_{w(k)w(k)}U_{kk})^2>0$.

(1) The vectors $u_j$ are orthogonal for $\langle\cdot,\cdot\rangle_\varepsilon$, so $G=R^T\Delta^2R$ with $\Delta^2=\diag(\langle u_j,u_j\rangle_\varepsilon)$. Hence $\langle u_j,u_j\rangle_\varepsilon=\det G_j/\det G_{j-1}=\varepsilon^{w(j)-1}(c_j+O(\varepsilon))$, with $c_j=\gamma_j/\gamma_{j-1}$.

(2) Since $R^T\Delta^2$ is lower triangular, the first $k$ rows of $G$ are $(R^T\Delta^2)_{[k],[k]}R_{[k],[n]}$, so for $k<j$ we get $\det G_{[k],[k-1]\cup\{j\}}=\det G_k\cdot\det R_{[k],[k-1]\cup\{j\}}=\det G_k\cdot r_{kj}$. By the Cauchy--Binet formula, $\det G_{[k],[k-1]\cup\{j\}}=\sum_{|S|=k}\varepsilon^{\sum_{s\in S}(s-1)}\det M_{S,[k]}\det M_{S,[k-1]\cup\{j\}}$, which is $O(\varepsilon^{\sigma_k})$ by the argument above, whereas $\det G_k=\varepsilon^{\sigma_k}(\gamma_k+O(\varepsilon))$. So $r_{kj}$ is finite for all $k<j$, $R$ is finite, and so is $[u_1\cdots u_n]=MR^{-1}$. Taking standard parts in $M=[u_1\cdots u_n]R$ gives $M=P'U'$, where $U'$ is unit upper triangular and $P'=MU'^{-1}$ is invertible. By (1), $\varepsilon^{i-1}u_{j,i}^2\leq\langle u_j,u_j\rangle_\varepsilon=O(\varepsilon^{w(j)-1})$, so $\st(u_{j,i})=0$ for $i<w(j)$: column $j$ of $P'$ vanishes above row $w(j)$. This column is column $w(j)$ of $L'=P'w^{-1}$, so $L'$ is lower triangular, and it is invertible.

(3) The factorisation exists because $\RR((t))$ is Pythagorean (Proposition~\ref{prop:laurent-pythagorean}). Put $T=\diag(1,t,\dots,t^{n-1})$ and $\Delta=\diag(\langle u_j,u_j\rangle_\varepsilon^{1/2})$. Then
\[ TM=\bigl(T[u_1\cdots u_n]\Delta^{-1}\bigr)(\Delta R), \]
and the first factor has orthonormal columns, so it is $\hat Q$. By Proposition~\ref{prop:laurent-pythagorean}, the entries of $\hat Q$ are finite, and $\st(\hat Q)$ is orthogonal because $\st$ is a ring homomorphism. By (1), $\hat Q_{ij}=t^{i-1}u_{j,i}/\Delta_{jj}$ with $\Delta_{jj}=t^{w(j)-1}(\sqrt{c_j}+O(t^2))$. For $i>w(j)$ this is infinitesimal, since $u_{j,i}$ is finite, and for $i=w(j)$ its standard part is $\st(u_{j,w(j)})/\sqrt{c_j}$, which has the sign of $L'_{w(j)w(j)}$. So column $j$ of $\st(\hat Q)$ vanishes below row $w(j)$. By induction on $w(j)$, the column with $w(j)=m$ is orthogonal to the columns $\pm e_1,\dots,\pm e_{m-1}$ found before, so it is $\pm e_m$.
\end{proof}

For $w=1$ the theorem is Proposition~\ref{prop:gs-limit}. Two further remarks. First, parts (1) and (2) use only field operations, in the field $\RR((\varepsilon))$, which is the completed localisation of the quiver with one loop $\varepsilon$; part (3) also needs the square roots $t^{i-1}$ of the weights, and $\RR((t))$ is the completed localisation of the quiver with one loop $t$. So the whole computation takes place in completed localisations, and the implementation of Chapter~\ref{chap:implementation} carries it out there; for $n=4$ it recovers $w$, $L'$ and $U'$ in all $24$ Bruhat cells (Table~\ref{tab:experiments}). Second, the computation needs $M$ exactly: the orders of the norms, and hence $w$, change under arbitrarily small perturbations of $M$, in agreement with the discontinuity of the Bruhat decomposition discussed in Section~\ref{sec:limitations}.

\section{Matrix Decompositions over Formal Power Series}
\label{sec:power-series}

\subsection{The Spectral Theorem}
The following theorem generalises the result of~\cite{Gutin2022DualNumbers}, which concerns matrices over the dual numbers $k[\varepsilon]/(\varepsilon^2)$; that result is recovered by reducing modulo $\varepsilon^2$.

\begin{theorem}[Spectral theorem over formal power series]\label{thm:spectral-power-series}
Let $k$ be a real closed field and let $S(\varepsilon)$ be a symmetric matrix with entries in $k[[\varepsilon]]$. Then there are an orthogonal matrix $U(\varepsilon)$ and a diagonal matrix $\Lambda(\varepsilon)$, both with entries in $k[[\varepsilon]]$, such that
\[ S(\varepsilon) = U(\varepsilon)\,\Lambda(\varepsilon)\,U(\varepsilon)^T. \]
If $k=\RR$ and the entries of $S(\varepsilon)$ are convergent power series, then $U(\varepsilon)$ and $\Lambda(\varepsilon)$ can be chosen with convergent entries.
\end{theorem}

For $k=\RR$ and convergent $S$, Theorem~\ref{thm:spectral-power-series} is a form of Rellich's theorem from analytic perturbation theory~\cite[Chapter~II]{Kato1995}. The proof below is purely algebraic; it applies to formal power series over any real closed field, and it illustrates the use of the transfer principle.

\begin{proof}
Let $K=k\{\{\varepsilon\}\}$ be the field of formal Puiseux series, which is real closed (Theorem~\ref{thm:newton-puiseux}) and contains $k[[\varepsilon]]$.

\emph{Step 1: the eigenvalues lie in $K$.} For each $n$, the spectral theorem for $n\times n$ symmetric matrices is a first-order sentence (Section~\ref{subsec:fol}), true over $\RR$. By Theorem~\ref{thm:tarski-transfer} it holds over $K$. Hence $S(\varepsilon)$ is orthogonally diagonalisable over $K$, and all roots of its characteristic polynomial lie in $K$.

\emph{Step 2: the eigenvalues are finite.} Over $\RR$, every eigenvalue $\lambda$ of a symmetric matrix satisfies $\lambda^2\leq\sum_{i,j}s_{ij}^2$; this inequality is first-order, so it transfers to $K$. The right-hand side lies in $k[[\varepsilon]]$, so it is smaller than $c=1+\sum_{i,j}s_{ij}(0)^2\in k$, and since $c\geq1$ we get $|\lambda|<c$. An element of $K$ that is bounded in absolute value by an element of $k$ has no term with a negative exponent, since such a term would dominate. Hence no eigenvalue has a term with a negative exponent.

\emph{Step 3: the eigenvalues have no fractional exponents.} Let $\lambda$ be an eigenvalue, let $d\geq1$ be minimal with $\lambda\in k((\varepsilon^{1/d}))$, put $s=\varepsilon^{1/d}$ and write $\lambda=f(s)=\sum_jf_js^j$. Suppose, for a contradiction, that $d\geq2$. Let $C=k(\sqrt{-1})$, which is algebraically closed, and let $\omega\in C$ be a primitive $2d$-th root of unity, so that $\omega^d=-1$. The substitution $s\mapsto\omega s$ defines a ring automorphism of $C((s))$ that maps $\varepsilon=s^d$ to $-\varepsilon$. Applying it to the identity $\det(S(\varepsilon)-f(s)I)=0$ yields
\[ \det\bigl(S(-\varepsilon)-f(\omega s)\,I\bigr)=0 . \]
So $f(\omega s)=\sum_jf_j\omega^j\varepsilon^{j/d}$ is an eigenvalue of the symmetric matrix $S(-\varepsilon)$, whose entries again lie in $k[[\varepsilon]]$. By Step 1 applied to $S(-\varepsilon)$, this eigenvalue lies in $K$, that is, all its coefficients $f_j\omega^j$ lie in $k$. Since $\omega^j\in k$ only when $\omega^j=\pm1$, that is, when $d$ divides $j$, we conclude that $f_j=0$ unless $d\mid j$. Then $\lambda\in k((\varepsilon))$, contradicting the minimality of $d$. Together with Step 2, every eigenvalue lies in $k[[\varepsilon]]$.

\emph{Step 4: eigenvectors.} For each eigenvalue $\lambda\in k[[\varepsilon]]$, the matrix $S(\varepsilon)-\lambda I$ has entries in the field $k((\varepsilon))$, so Gaussian elimination over $k((\varepsilon))$ yields a basis of its kernel with entries in $k((\varepsilon))$.

\emph{Step 5: orthonormalisation.} By Proposition~\ref{prop:laurent-pythagorean}, $k((\varepsilon))$ is an ordered Pythagorean field, so the Gram--Schmidt process applies within each eigenspace and produces an orthonormal basis of it with entries in $k((\varepsilon))$.

\emph{Step 6: orthogonality between eigenspaces.} Eigenvectors of a symmetric matrix for distinct eigenvalues are orthogonal over $\RR$; this is first-order and transfers to $K$.

\emph{Step 7: integrality.} The resulting eigenvectors are unit vectors over $k((\varepsilon))$, so by Proposition~\ref{prop:laurent-pythagorean} their entries lie in $k[[\varepsilon]]$. Assembling them into the columns of $U(\varepsilon)$ gives the required decomposition.

For the convergent case, replace formal Puiseux and Laurent series by convergent ones throughout. Convergent Puiseux series form a real closed field (Theorem~\ref{thm:newton-puiseux}), convergent Laurent series form a Pythagorean field (Proposition~\ref{prop:laurent-pythagorean}), and the substitution in Step 3 preserves convergence.
\end{proof}

\subsection{The Singular Value Decomposition}
\begin{corollary}[SVD over formal power series]\label{cor:svd-power-series}
Let $k$ be a real closed field and $B\in M_{p\times q}(k[[\varepsilon]])$. Then there are orthogonal matrices $Q\in\Or_p(k[[\varepsilon]])$ and $P\in\Or_q(k[[\varepsilon]])$ and a diagonal matrix $\Sigma\in M_{p\times q}(k[[\varepsilon]])$, whose non-zero diagonal entries have positive leading coefficients, such that $B=Q\Sigma P^T$. If $k=\RR$ and $B$ is convergent, the factors can be chosen convergent.
\end{corollary}
\begin{proof}
Apply Theorem~\ref{thm:spectral-power-series} to $B^TB=P\Lambda P^T$, with $\Lambda=\diag(\mu_1,\dots,\mu_q)$ and columns $p_i$ of $P$. Then $\mu_i=\|Bp_i\|^2$ is a sum of squares in $k[[\varepsilon]]$, so by Proposition~\ref{prop:laurent-pythagorean} it has a square root $\sigma_i\in k[[\varepsilon]]$ with non-negative leading coefficient. Order the columns of $P$ so that $\sigma_1,\dots,\sigma_r$ are the non-zero ones; then $Bp_i=0$ for $i>r$, since $\|Bp_i\|^2=\mu_i=0$. For $i\leq r$ put $q_i=Bp_i/\sigma_i$. These are orthonormal vectors over $k((\varepsilon))$, since $q_i^Tq_j=p_i^TB^TBp_j/(\sigma_i\sigma_j)=\delta_{ij}$, and by Proposition~\ref{prop:laurent-pythagorean} their entries lie in $k[[\varepsilon]]$. Extend $q_1,\dots,q_r$ to an orthonormal basis $q_1,\dots,q_p$ of $k((\varepsilon))^p$ by the Gram--Schmidt process; again all entries lie in $k[[\varepsilon]]$. With $Q=[q_1\cdots q_p]$ and $\Sigma$ the $p\times q$ diagonal matrix with entries $\sigma_1,\dots,\sigma_r$ followed by zeros, we have $BP=Q\Sigma$.
\end{proof}

For $k=\RR$ and convergent $B$, this is the existence of an analytic singular value decomposition, which can also be computed numerically~\cite{BunseGerstnerByersMehrmannNichols1991}; reducing modulo $\varepsilon^2$ recovers the SVD of dual-number matrices of~\cite{Gutin2022DualNumbers}.

\section{A Parity-Equivariant Spectral Theorem}
\label{sec:parity}
In~\cite{Gutin2024Clifford} we observed that many important decompositions of numerical linear algebra can be viewed as variants of the spectral theorem over various $*$-algebras, equivariant with respect to a group of symmetries. In this section we make this observation precise for one case, in which the symmetry is the parity automorphism $\varepsilon\mapsto-\varepsilon$ of the power series ring.

Throughout this section, $k$ is a real closed field, $\mathcal R=k[[\varepsilon]]$, and $\gamma:\mathcal R\to\mathcal R$ is the automorphism $\gamma(f)(\varepsilon)=f(-\varepsilon)$, extended entrywise to matrices. The fixed ring is $\mathcal F=\mathcal R^\gamma=k[[\varepsilon^2]]$. An element or matrix $X$ is \emph{even} if $\gamma(X)=X$ and \emph{odd} if $\gamma(X)=-X$; every $X$ is uniquely the sum of an even and an odd part. Fix $p,q\geq0$, let $J=\diag(I_p,-I_q)$, and let $\tau=J\gamma$ be the semilinear involution $v\mapsto J\gamma(v)$ of $\mathcal R^{p+q}$. A vector $v$ is \emph{$\tau$-homogeneous} of parity $\pm1$ if $\tau v=\pm v$.

Step 3 of the proof below uses Galois descent~\cite[\S18]{KMRT1998} for a quadratic extension, in the following elementary form.

\begin{lemma}\label{lem:descent}
Let $L$ be a field of characteristic different from $2$, $\theta$ an automorphism of $L$ of order $2$, and $K=L^\theta$ its fixed field. Let $\theta$ act on $L^n$ coordinatewise. Then every $\theta$-stable $L$-subspace $W\subseteq L^n$ has an $L$-basis consisting of $\theta$-fixed vectors.
\end{lemma}
\begin{proof}
Choose $a\in L\setminus K$ and put $c=a-\theta(a)\neq0$, so that $\theta(c)=-c$. For $w\in W$ we have
\[ w=\tfrac12\,(w+\theta w)+\tfrac12\,c^{-1}\bigl(c\,(w-\theta w)\bigr), \]
and both $w+\theta w$ and $c(w-\theta w)$ are $\theta$-fixed elements of $W$. Hence the $\theta$-fixed vectors span $W$ over $L$, and a basis can be extracted from them.
\end{proof}

\begin{theorem}[Parity-equivariant spectral theorem]\label{thm:parity}
Let $S$ be a symmetric $(p+q)\times(p+q)$ matrix over $\mathcal R$ satisfying the equivariance condition $J\gamma(S)J=S$; equivalently, the diagonal blocks of $S$ (of sizes $p$ and $q$) are even and the off-diagonal blocks are odd. Then there exist an orthogonal matrix $U$ over $\mathcal R$ and a diagonal sign matrix $J'$ with the same number of entries $+1$ as $J$, such that $J\gamma(U)=UJ'$ and $U^TSU$ is block diagonal, each block being either
\begin{itemize}
  \item a $1\times1$ block $(\alpha)$ with $\alpha\in\mathcal F$, or
  \item a $2\times2$ block $\begin{pmatrix}\alpha&\beta\\ \beta&\alpha\end{pmatrix}$ with $\alpha\in\mathcal F$ and $\beta\in\varepsilon\mathcal F$ odd, on which $J'$ restricts to $\diag(1,-1)$.
\end{itemize}
\end{theorem}
\begin{proof}
By Theorem~\ref{thm:spectral-power-series}, $S=\sum_\lambda\lambda P_\lambda$, where $\lambda$ runs over the distinct eigenvalues in $\mathcal R$ and $P_\lambda$ is the orthogonal projector onto the $\lambda$-eigenspace; the $P_\lambda$ have entries in $\mathcal R$, since $P_\lambda=\sum_iu_iu_i^T$ for orthonormal eigenvectors $u_i$ over $\mathcal R$. Note that $\langle\tau u,\tau v\rangle=\gamma\langle u,v\rangle$ for the standard bilinear form, because $J^TJ=I$.

\emph{Step 1: $\tau$ permutes the eigenspaces.} The equivariance condition says that $\tau S=S\tau$. If $Sv=\lambda v$ then $S\tau v=\tau(Sv)=\tau(\lambda v)=\gamma(\lambda)\tau v$. Hence $\tau$ maps the $\lambda$-eigenspace onto the $\gamma(\lambda)$-eigenspace, and the eigenvalues fall into orbits $\{\lambda,\gamma(\lambda)\}$ of size one (when $\lambda$ is even) or two.

\emph{Step 2: orbits of size two.} Let $\lambda\neq\gamma(\lambda)$ and write $\lambda=\alpha+\beta$ with $\alpha$ even and $\beta\neq0$ odd. Let $v_1,\dots,v_m$ be an orthonormal basis over $\mathcal R$ of the $\lambda$-eigenspace. Then $\tau v_1,\dots,\tau v_m$ is an orthonormal basis of the $\gamma(\lambda)$-eigenspace, and $\langle v_i,\tau v_j\rangle=0$ because eigenvectors for distinct eigenvalues are orthogonal. Put
\[ u_i^{\pm}=\tfrac{1}{\sqrt2}\,(v_i\pm\tau v_i). \]
Then $\tau u_i^\pm=\pm u_i^\pm$, the vectors $u_i^\pm$ are orthonormal (using $\langle\tau v_i,\tau v_j\rangle=\gamma(\delta_{ij})=\delta_{ij}$), and
\[ Su_i^+=\tfrac{1}{\sqrt2}\bigl(\lambda v_i+\gamma(\lambda)\tau v_i\bigr)=\alpha u_i^++\beta u_i^-,\qquad Su_i^-=\alpha u_i^-+\beta u_i^+ . \]
So on $\mathrm{span}(u_i^+,u_i^-)$ the matrix of $S$ is $\bigl(\begin{smallmatrix}\alpha&\beta\\\beta&\alpha\end{smallmatrix}\bigr)$, and $\tau$ acts as $\diag(1,-1)$.

\emph{Step 3: orbits of size one.} Let $\lambda=\gamma(\lambda)$, so $\lambda\in\mathcal F$, and let $W$ be the $\lambda$-eigenspace over the fraction field $k((\varepsilon))$. By Step 1, $W$ is $\tau$-stable, and $\tau$ acts on $k((\varepsilon))^{p+q}$ as $J$ composed with the coordinatewise automorphism $\gamma$. Applying Lemma~\ref{lem:descent} after the change of coordinates $v\mapsto(v_1,\dots,v_p,\varepsilon v_{p+1},\dots,\varepsilon v_{p+q})$, which turns $\tau$ into the coordinatewise action of $\gamma$, shows that $W$ has a basis of $\tau$-fixed vectors. We orthonormalise it by the Gram--Schmidt process and claim that homogeneity is preserved. If $u,w$ are homogeneous of parities $\pi_u,\pi_w$, then $\gamma\langle u,w\rangle=\langle\tau u,\tau w\rangle=\pi_u\pi_w\langle u,w\rangle$, so $\langle u,w\rangle$ is homogeneous of parity $\pi_u\pi_w$, and $\langle u,u\rangle$ is even. Hence $w-\frac{\langle u,w\rangle}{\langle u,u\rangle}u$ is homogeneous of parity $\pi_w$. Finally, $\langle w,w\rangle\in k((\varepsilon^2))$ is positive, so by the argument of Proposition~\ref{prop:laurent-pythagorean} its square root is $\sqrt c\,\varepsilon^{m}(1+h)^{1/2}$ with $h$ even, which is homogeneous of parity $(-1)^m$; dividing by it preserves homogeneity. We obtain an orthonormal basis of $W$ consisting of homogeneous unit vectors, whose entries lie in $\mathcal R$ by Proposition~\ref{prop:laurent-pythagorean}. On these vectors $S$ acts as the scalar $\lambda\in\mathcal F$.

\emph{Step 4: assembly.} Let $U$ be the orthogonal matrix whose columns are the vectors constructed in Steps 2 and 3, grouped by blocks, and let $J'$ be the diagonal matrix of their parities. Then $\tau U=UJ'$, that is, $J\gamma(U)=UJ'$, and $U^TSU$ has the stated block form. Evaluating $J\gamma(U)=UJ'$ at $\varepsilon=0$ gives $J'=U(0)^TJU(0)$, so $\tr J'=\tr J$, and $J'$ has the same number of entries $+1$ as $J$.
\end{proof}

The $2\times2$ parity blocks play the role of the $2\times2$ blocks in the real normal forms of normal and skew-symmetric matrices: they occur exactly when the eigenvalues are exchanged by the symmetry. Two consequences are worth recording.

\begin{corollary}[SVD of odd matrices]\label{cor:odd-svd}
Let $B\in M_{p\times q}(\mathcal R)$ be odd. Then $B=Q\Sigma P^T$ with $P\in\Or_q(\mathcal F)$ and $Q\in\Or_p(\mathcal F)$ even and orthogonal, and $\Sigma$ diagonal with odd entries in $\varepsilon\mathcal F$.
\end{corollary}
\begin{proof}
Since $B$ is odd, $B=\varepsilon C$ with $C$ having entries in $\mathcal F=k[[\varepsilon^2]]$, which is again a formal power series ring over $k$, in the variable $\varepsilon^2$. Corollary~\ref{cor:svd-power-series} applied to $C$ over $\mathcal F$ gives $C=Q\Sigma_0P^T$ with $P,Q,\Sigma_0$ over $\mathcal F$, and then $B=Q(\varepsilon\Sigma_0)P^T$.
\end{proof}

Corollary~\ref{cor:odd-svd} is the decomposition behind the self-adjoint matrix $\bigl(\begin{smallmatrix}0&B^T\\B&0\end{smallmatrix}\bigr)$, which satisfies the hypothesis of Theorem~\ref{thm:parity} with $J=\diag(I_q,-I_p)$ precisely when $B$ is odd; its eigenvalues $\pm\sigma_i$ are exchanged by $\gamma$, and Theorem~\ref{thm:parity} produces the parity blocks $\bigl(\begin{smallmatrix}0&\sigma_i\\\sigma_i&0\end{smallmatrix}\bigr)$.

\begin{corollary}[Hermitian matrices over $(\mathcal R,\varepsilon^*=-\varepsilon)$]\label{cor:twisted-hermitian}
Equip $\mathcal R$ with the involution $f^*=\gamma(f)$, so that $\varepsilon^*=-\varepsilon$, and let $H\in M_n(\mathcal R)$ satisfy $H^*=H$, that is, $\gamma(H)^T=H$. Then there is $U\in M_n(\mathcal R)$ with $U^*U=I$ such that $U^*HU$ is block diagonal with $1\times1$ blocks $(\alpha)$, $\alpha\in\mathcal F$, and $2\times2$ blocks $\bigl(\begin{smallmatrix}\alpha&\beta\\-\beta&\alpha\end{smallmatrix}\bigr)$ with $\alpha\in\mathcal F$ and $\beta\in\varepsilon\mathcal F$.
\end{corollary}
\begin{proof}[Sketch of proof]
Let $C=k(\sqrt{-1})$ and $i=\sqrt{-1}$. The substitution $\varepsilon=i\delta$ embeds $\mathcal R$ into $C[[\delta]]$, and under this embedding the involution $\gamma$ becomes complex conjugation of the coefficients, $\sigma(\sum c_j\delta^j)=\sum\bar c_j\delta^j$. Hence $H$ becomes a complex Hermitian matrix over $C[[\delta]]$. Applying Theorem~\ref{thm:spectral-power-series} to its realification $\bigl(\begin{smallmatrix}\mathrm{Re}\,H&-\mathrm{Im}\,H\\ \mathrm{Im}\,H&\mathrm{Re}\,H\end{smallmatrix}\bigr)$ shows that it is unitarily diagonalisable over $C[[\delta]]$ with eigenvalues in $k[[\delta]]$. The image of $\mathcal R$ is the fixed ring of the automorphism $\theta=\sigma\circ(\delta\mapsto-\delta)$ of $C[[\delta]]$, and $\theta$ permutes the eigenvalues as $\lambda(\delta)\mapsto\lambda(-\delta)$. The argument of Theorem~\ref{thm:parity}, with $\tau$ replaced by $\theta$ and with the anti-fixed unit $i$ used to convert anti-fixed vectors into fixed ones, produces a $\theta$-fixed unitary matrix, that is, a unitary matrix over $\mathcal R$. An orbit $\{\lambda(\delta),\lambda(-\delta)\}$ of size two with $\lambda=\alpha+\beta'$ ($\alpha$ even, $\beta'$ odd in $\delta$) yields the block $\bigl(\begin{smallmatrix}\alpha&i\beta'\\-i\beta'&\alpha\end{smallmatrix}\bigr)$, and rewriting $\delta=-i\varepsilon$ shows that $\alpha\in\mathcal F$ and $\beta:=i\beta'\in\varepsilon\mathcal F$.
\end{proof}

For example, $H=\bigl(\begin{smallmatrix}0&\varepsilon\\-\varepsilon&0\end{smallmatrix}\bigr)$ is Hermitian for $\varepsilon^*=-\varepsilon$ but cannot be diagonalised: its eigenvalues, $\pm\sqrt{-1}\,\varepsilon$, do not lie in $\mathcal R$. It is already a parity block. Modulo $\varepsilon^2$, Corollary~\ref{cor:twisted-hermitian} describes Hermitian matrices over the dual numbers with the involution $\varepsilon^*=-\varepsilon$; Chapter~\ref{chap:indexing} shows that the corresponding classical problem is the normal form of real skew-symmetric matrices under orthogonal similarity.

\section{Shortest Paths as an Auto-Limit}
\label{sec:shortest-paths}
Auto-limits also apply to algorithms whose inputs carry orders of magnitude. Let $G$ be a directed graph with vertex set $V=\{1,\dots,n\}$, without loops and with at most one arrow from $i$ to $j$, carrying a positive rational weight $d_{ij}$ on each arrow $i\to j$; put $d_{ij}=\infty$ if there is no such arrow. (Dijkstra's algorithm also allows the weight $0$, which we exclude.) The \emph{length} $\ell(\gamma)$ of a walk $\gamma$ is the sum of the weights of its arrows. Write $\mathrm{dist}(i,j)$ for the least length of a path from $i$ to $j$, with $\mathrm{dist}(i,i)=0$ and $\mathrm{dist}(i,j)=\infty$ if there is no path, and $N(i,j)$ for the number of paths of that length. Encoding each weight $d$ by the quantity $\varepsilon^d$, of order of magnitude $d$, put
\[ M=I-A,\qquad A=\bigl(\varepsilon^{d_{ij}}\bigr)_{i,j\in V},\qquad\text{that is,}\quad m_{ij}=\delta_{ij}-\varepsilon^{d_{ij}}, \]
over the field $\RR\{\{\varepsilon\}\}$ of real Puiseux series (Section~\ref{sec:puiseux}), with $\varepsilon^\infty=0$.

A non-zero $x\in\RR\{\{\varepsilon\}\}$ has a \emph{leading term} $c\,\varepsilon^q$, where $q$ is the least exponent of $x$, its \emph{order}, and $c\neq0$. For series with non-negative coefficients nothing cancels: the order of a sum is the least order of the summands, whose leading coefficients are added if several have this order, and orders add under multiplication. Taking leading terms therefore turns sums into minima of exponents and products into sums, and maps such series to the min-plus, or tropical, semiring~\cite{MaclaganSturmfels2015}. That shortest paths can be found by Gaussian elimination over this semiring is classical~\cite{Carre1971,GondranMinoux2008}. The semiring has no subtraction, and elimination there replaces the inverse of $1-a$ by the closure $a^*=1\oplus a\oplus a^2\oplus\cdots$. Over $\RR\{\{\varepsilon\}\}$, subtraction and division are available, like the other operations of Section~\ref{sec:proposal}, and the elimination becomes ordinary Gaussian elimination, followed by taking leading terms.

For $S\subseteq V$ and $i,j\in V$, let $w_S(i,j)$ be the sum of $\varepsilon^{\ell(\gamma)}$ over the walks $\gamma$ from $i$ to $j$ that have at least one arrow and whose interior vertices, all but the first and the last, lie in $S$. If the weights have common denominator $r$, a walk with $k$ arrows has length at least $k\min d_{ij}$, so only finitely many walks have length below any bound, and $w_S(i,j)\in\RR((\varepsilon^{1/r}))$. Its coefficients are non-negative integers. If $w_S(i,j)\neq0$, its order is the least length $\mathrm{dist}_S(i,j)$ of such a walk, and its leading coefficient is the number $N_S(i,j)$ of these walks of that length; for $i\neq j$ they are paths, because removing a closed sub-walk shortens a walk.

\begin{proposition}\label{prop:dijkstra}
In this notation:
\begin{enumerate}
  \item $M$ is invertible, and $M^{-1}=I+\bigl(w_V(i,j)\bigr)_{i,j\in V}$. If $\mathrm{dist}(i,j)<\infty$, the leading term of $(M^{-1})_{ij}$ is $N(i,j)\,\varepsilon^{\mathrm{dist}(i,j)}$; otherwise $(M^{-1})_{ij}=0$.
  \item For every ordering of the vertices, Gaussian elimination without pivoting applies to $M$. Once the vertices of a set $S$ have been eliminated, the Schur complement is $\bigl(\delta_{ij}-w_S(i,j)\bigr)_{i,j\notin S}$, and every pivot has leading term $1$.
  \item Let $s\in V$, and solve $M^Tz=e_s$ by Gaussian elimination, taking as the next pivot a remaining vertex $u$ whose current right-hand side entry $b_u$ has the least order. Then $b_u$ and $z_u$ both have the leading term $N(s,u)\,\varepsilon^{\mathrm{dist}(s,u)}$ if $\mathrm{dist}(s,u)<\infty$, and both vanish otherwise. Eliminating $u$ gives each remaining $b_j$ the leading term of $b_j+\varepsilon^{d_{uj}}b_u$.
\end{enumerate}
\end{proposition}
\begin{proof}
(1) The $(i,j)$ entry of $A^k$ is the sum of $\varepsilon^{\ell(\gamma)}$ over the walks from $i$ to $j$ with $k$ arrows, and its order is at least $k\min d_{ij}$. So $\sum_{k\geq0}A^k$ converges, it is the inverse of $I-A$, and its $(i,j)$ entry is $\delta_{ij}+w_V(i,j)$. The rest follows from the description of $w_V(i,j)$ above.

(2) Let $T=V\setminus S$. By (1) applied to the graph on $S$, $M_{S,S}$ is invertible, and $(M_{S,S}^{-1})_{pq}$ is $\delta_{pq}$ plus the sum over the walks from $p$ to $q$ inside $S$. A walk from $i\in T$ to $j\in T$ with interior in $S$ is either an arrow $i\to j$, or an arrow $i\to p$ into $S$, followed by a walk from $p$ to $q$ inside $S$, possibly of length zero, and an arrow $q\to j$. So the $(i,j)$ entry of the Schur complement $M_{T,T}-M_{T,S}M_{S,S}^{-1}M_{S,T}$, which is $\delta_{ij}-\varepsilon^{d_{ij}}-\sum_{p,q\in S}\varepsilon^{d_{ip}}(M_{S,S}^{-1})_{pq}\,\varepsilon^{d_{qj}}$, is $\delta_{ij}-w_S(i,j)$. Its diagonal entries $1-w_S(i,i)$ are the pivots, and $w_S(i,i)$ has positive order.

(3) Initially $b=e_s$, so $s$ is eliminated first, and $b_s$ and $z_s=1+w_V(s,s)$ have leading term $1$. By the same computation as in (2), once a set $S\ni s$ has been eliminated, $b_j=w_S(s,j)$ for $j\notin S$. Let $u$ be the next pivot and $S'=S\cup\{u\}$. We show, by induction on $|S|$, that every shortest path from $s$ to $u$ has its interior in $S$ and that $\mathrm{dist}(s,p)\leq\mathrm{dist}(s,u)$ for $p\in S$. If $\mathrm{dist}(s,u)=\infty$, then $b_u=0$ and $z_u=0$. Otherwise, let $\gamma$ be a shortest path from $s$ to $u$ and $v$ its first vertex outside $S$. The part of $\gamma$ up to $v$ has its interior in $S$, so $\mathrm{dist}_S(s,v)\leq\mathrm{dist}(s,u)$, with strict inequality if $v\neq u$. Since $b_u$ has the least order, $\mathrm{dist}_S(s,u)\leq\mathrm{dist}_S(s,v)$. Hence $v=u$, $\mathrm{dist}_S(s,u)=\mathrm{dist}(s,u)$ and $N_S(s,u)=N(s,u)$, so $b_u=w_S(s,u)$ has the leading term $N(s,u)\,\varepsilon^{\mathrm{dist}(s,u)}$, and so has $z_u=(M^{-1})_{su}$, by (1).

After the elimination of $u$, $b_j=w_{S'}(s,j)$. The walks counted in $w_S(s,j)+\varepsilon^{d_{uj}}w_S(s,u)$ are those that avoid $u$, and those that visit $u$ once and then take the arrow $u\to j$. Every other walk counted in $w_{S'}(s,j)$ is longer than one of these. If it visits $u$ more than once, removing the closed walk between its first and last visits to $u$ gives a shorter walk that visits $u$ once. If it visits $u$ once and then takes an arrow $u\to p$ with $p\in S$, replace its part up to $p$, of length at least $\mathrm{dist}(s,u)+d_{up}$, by a shortest path from $s$ to $p$; by induction, this path has its interior in $S$ and length at most $\mathrm{dist}(s,u)$, so the result is a shorter walk that avoids $u$. Since nothing cancels, $b_j$ has the leading term of $w_S(s,j)+\varepsilon^{d_{uj}}w_S(s,u)$. Its order is at least that of $b_u$, so the next pivot $u'$ satisfies $\mathrm{dist}(s,u')\geq\mathrm{dist}(s,u)$, which completes the induction.
\end{proof}

In leading terms, part~(3) is Dijkstra's algorithm~\cite{Dijkstra1959}. The orders of the entries $b_j$ are its tentative distances; the pivot rule chooses a remaining vertex of least tentative distance, which is then final; and $b_j\mapsto b_j+\varepsilon^{d_{uj}}b_u$ is its update $\mathrm{dist}_j\leftarrow\min(\mathrm{dist}_j,\mathrm{dist}_u+d_{uj})$, in which the leading coefficients also count shortest paths. Part~(3) also accounts for what Dijkstra's algorithm leaves out: in this order of elimination, the fill-in of the Schur complements, the pivots and the back substitution change no leading term. In a fixed order, Gaussian elimination with back substitution still gives the distances, by part~(1), but the fill-in is then needed; in leading terms, this is the elimination method of Carr\'e~\cite{Carre1971}. All the quantities computed lie in the field $\RR\langle\varepsilon\rangle$ of algebraic Puiseux series, so Proposition~\ref{prop:autolimit} applies: for every sufficiently small real $\varepsilon>0$, Gaussian elimination in real arithmetic, with the remaining vertex of largest $b_u$ as the next pivot, makes choices allowed in part~(3), and, if $\mathrm{dist}(s,u)<\infty$,
\[ \mathrm{dist}(s,u)=\lim_{\varepsilon\to0^+}\frac{\log z_u(\varepsilon)}{\log\varepsilon}. \]

By part~(2), $PMP^T$ lies in the big cell $B_-B$ of the Bruhat decomposition (Section~\ref{sec:bruhat}) for every permutation matrix $P$, so no pivoting is needed for the factorisation to exist. The pivots of part~(3) are chosen for another reason. Like the permutation in Theorem~\ref{thm:bruhat}, they are read off from orders of magnitude, here those of the right-hand side, and they make the leading terms available without fill-in or back substitution.

For an undirected graph, that is, when $d_{ij}=d_{ji}$, the matrix $M$ is symmetric. It is positive definite, because it differs from $I$ by infinitesimals, so the elimination computes $M=LDL^T$, and $M=(LD^{1/2})(LD^{1/2})^T$ is its Cholesky factorisation; the pivots in $D$ and their square roots have leading term~$1$. For a directed graph, $M$ is not symmetric, but the elimination is still the Cholesky factorisation of a Hermitian matrix. By Proposition~\ref{prop:2R-algorithms}(1), whose proof uses only field operations and so holds over $\RR\{\{\varepsilon\}\}$, the Doolittle factorisation of $M$ is the Cholesky factorisation of the Hermitian matrix $(M,M^T)$ over $\twoR\otimes\RR\{\{\varepsilon\}\}$. The involution of $\twoR$ exchanges $M$ and $M^T$, that is, it reverses every arrow. So Dijkstra's algorithm, like Rutishauser's LR algorithm, can be obtained from an algorithm for Hermitian matrices over a $*$-algebra, here followed by an auto-limit. This reinterprets a classical algorithm and gives no new or faster way to find shortest paths: Dijkstra's algorithm, which works with the orders alone, is much cheaper than Gaussian elimination over $\RR\{\{\varepsilon\}\}$ or in real arithmetic.

\section{A Chain of Limits}
\label{sec:chain}
The factorisations of this chapter reduce to the Hermitian eigendecomposition by a chain of limits. The singular value decomposition of a matrix $W$ is a Hermitian eigendecomposition: over the $*$-algebra of the quiver of the singular value decomposition (Theorem~\ref{thm:A}), or classically, of the matrix $\bigl(\begin{smallmatrix}0&W^*\\W&0\end{smallmatrix}\bigr)$. The next proposition obtains the QR factorisation $W=QR$ (the factorisation, not the QR algorithm) as a limit of the singular value decomposition, with the columns of $W$ weighted by decreasing orders of magnitude. No finite algorithm computes a singular value decomposition, so this limit is taken through the identification of $\RR\langle\varepsilon\rangle$ with germs of semialgebraic functions, as in the proof of Proposition~\ref{prop:autolimit}, rather than as the auto-limit of an algorithm. With the rows weighted instead, the QR factorisation in turn has the Bruhat decomposition as an auto-limit (Theorem~\ref{thm:bruhat}), and the $LU$ factorisation when the leading principal minors are non-zero (Proposition~\ref{prop:gs-limit}). Finally, Dijkstra's shortest-path algorithm is the leading-order part of Gaussian elimination (Proposition~\ref{prop:dijkstra}). So the chain runs from the Hermitian eigendecomposition, through the singular value decomposition and the QR factorisation, to the Bruhat and $LU$ factorisations and to Dijkstra's algorithm.

\begin{proposition}\label{prop:qr-from-svd}
Let $W\in\RR^{m\times n}$ have linearly independent columns $w_1,\dots,w_n$, and let $W=QR$ be its QR factorisation, where $Q$ has orthonormal columns $q_1,\dots,q_n$ and $R$ is upper triangular with positive diagonal. Put $D_\varepsilon=\diag(1,\varepsilon,\dots,\varepsilon^{n-1})$, and let $WD_\varepsilon=U\Sigma V^T$ be a singular value decomposition over $\RR\langle\varepsilon\rangle$, which exists by Theorem~\ref{thm:tarski-transfer}, in which $U$ has orthonormal columns $u_1,\dots,u_n$, $V$ is orthogonal, $\Sigma=\diag(\sigma_1,\dots,\sigma_n)$ with $\sigma_1\geq\dots\geq\sigma_n>0$, and the signs are chosen so that $u_k^Tw_k>0$. Then $\sigma_k=\varepsilon^{k-1}r_{kk}\bigl(1+O(\varepsilon^2)\bigr)$, the matrices $U$ and $U^TW=\Sigma V^TD_\varepsilon^{-1}$ are finite, and
\[ \st(U)=Q,\qquad \st\bigl(\Sigma V^TD_\varepsilon^{-1}\bigr)=R. \]
In particular, as a real parameter $\varepsilon$ tends to $0^+$, the left singular vectors of $WD_\varepsilon$, with these signs, tend to the columns of $Q$.
\end{proposition}
\begin{proof}
For $1\leq k\leq n$, let $\Lambda^k$ denote the $k$-th compound matrix, whose columns are indexed by the subsets $J\subseteq\{1,\dots,n\}$ with $k$ elements. Since $\Lambda^k$ is multiplicative, $\Lambda^k(WD_\varepsilon)=\Lambda^k(U)\,\Lambda^k(\Sigma)\,\Lambda^k(V)^T$ is a singular value decomposition, so the largest singular value of $\Lambda^k(WD_\varepsilon)$ is $s_k=\sigma_1\cdots\sigma_k$, with left singular vector $u_1\wedge\dots\wedge u_k$. The column indexed by $J$ is $\varepsilon^{e_J}w_J$, where $w_J=\bigwedge_{j\in J}w_j$ and $e_J=\sum_{j\in J}(j-1)$. The least exponent, $e_k=0+1+\dots+(k-1)$, occurs only for $J=[k]=\{1,\dots,k\}$, every other $J$ has $e_J\geq e_k+1$, and $\|w_{[k]}\|=r_{11}\cdots r_{kk}$, the volume spanned by $w_1,\dots,w_k$. Bounding $s_k$ below by the norm of the column $[k]$ and above by the Frobenius norm gives $s_k^2=\varepsilon^{2e_k}(r_{11}\cdots r_{kk})^2\bigl(1+O(\varepsilon^2)\bigr)$, and $\sigma_k=s_k/s_{k-1}$ gives the formula for $\sigma_k$. If $x$ is a unit left singular vector for $s_k$, then
\[ s_k^2=\bigl\|\Lambda^k(WD_\varepsilon)^Tx\bigr\|^2=\sum_J\varepsilon^{2e_J}(x\cdot w_J)^2=\varepsilon^{2e_k}(x\cdot w_{[k]})^2+O(\varepsilon^{2e_k+2}), \]
so $(x\cdot\hat w)^2=1+O(\varepsilon^2)$ for the unit vector $\hat w=w_{[k]}/\|w_{[k]}\|=q_1\wedge\dots\wedge q_k$, and $\st(x)=\pm\hat w$. Taking $x=u_1\wedge\dots\wedge u_k$ and using that $\st$ is a ring homomorphism, $\st(u_1)\wedge\dots\wedge\st(u_k)=\pm q_1\wedge\dots\wedge q_k$, so the orthonormal vectors $\st(u_1),\dots,\st(u_k)$ span the same subspace as $q_1,\dots,q_k$. By induction on $k$, $\st(u_k)$ is a unit vector in this subspace orthogonal to $q_1,\dots,q_{k-1}$, so $\st(u_k)=\pm q_k$; the sign is $+$, because $u_k^Tw_k>0$ and $q_k^Tw_k=r_{kk}>0$. Hence $\st(U)=Q$ and $\st(U^TW)=Q^TW=R$, and $U^TW=\Sigma V^TD_\varepsilon^{-1}$ because $U^TWD_\varepsilon=\Sigma V^T$. Finally, the singular values have distinct orders, so for small real $\varepsilon>0$ they are distinct and the decomposition with these signs is unique; the identification of $\RR\langle\varepsilon\rangle$ with germs of semialgebraic functions (Section~\ref{sec:puiseux}) gives the last statement.
\end{proof}

The proposition is elementary; its proof uses minors, as does that of Theorem~\ref{thm:bruhat}. Together with that theorem and Propositions~\ref{prop:gs-limit} and~\ref{prop:dijkstra}, it shows that the factorisations of this chapter, and Dijkstra's algorithm, are reached from the Hermitian eigendecomposition by exact encodings and limits, one step at a time; the steps weight different things, rows or columns, and do not combine into a single limit. The relation is structural, not computational: in practice the singular value decomposition is computed by algorithms built from QR factorisations, not the other way round.
%======================================================================
%   CHAPTER 5: DOUBLY-SIGNED QUIVERS AND REAL *-ALGEBRAS
%======================================================================
\chapter{Doubly-Signed Quivers and Real \texorpdfstring{$*$}{*}-Algebras}
\label{chap:quivers}

Chapter~\ref{chap:matrixops} stored the elements of a $*$-algebra as matrices indexed by positions $c_k$, with entries in the corners between the orbit idempotents $1_c$. This chapter explains where the orbits and the corners come from. Its main result, Theorem~\ref{thm:structure}, which is part~(1) of Theorem~\ref{thm:A}, is an analogue of Gabriel's theorem for real $*$-algebras: every finite-dimensional real $*$-algebra $A$ is $*$-isomorphic to
\[ \End(V,\eta)\otimes_{\RR^{\mathcal O}}\RR Q/I , \]
a matrix algebra, with prescribed sizes and signs, over the quotient of the path $*$-algebra $\RR Q$ of a doubly-signed quiver $Q$ by an ideal of relations $I$. The first factor supplies the positions and the second the entries. The quiver, the sizes and the signatures of the forms are invariants of $A$, and all continuous parameters lie in $I$.

Two algebras with the quiver $Q$ have distinct roles in this thesis. The quotient $\RR Q/\RR Q_{\geq2}$, in which all products of arrows vanish, is the algebra whose Hermitian matrices encode the classical decompositions (Chapter~\ref{chap:indexing}). The path $*$-algebra $\RR Q$ itself is where computations take place, once its arrows have been inverted (Chapter~\ref{division-problem}).

The chapter has two parts. Sections~\ref{sec:quivers} to~\ref{sec:dsq} introduce four kinds of quivers, each with more structure than the one before: quivers; \emph{symmetric quivers}, which carry an involution that reverses arrows; \emph{signed quivers}, which carry in addition a sign on each vertex and arrow that the involution fixes; and \emph{doubly-signed quivers}, which carry a second involution and a second family of signs, which record whether the scalars at a vertex are real, complex or quaternionic. Sections~\ref{sec:classification} to~\ref{sec:structure} are about real $*$-algebras: those whose radical squares to zero (Theorem~\ref{thm:classification}), the $*$-simple ones (Lemma~\ref{lem:star-simple}), the path $*$-algebra of a doubly-signed quiver (Definition~\ref{def:path-star-algebra}) and the structure theorem.

\section{Quivers}
\label{sec:quivers}
Quivers, path algebras and representations were introduced in Section~\ref{sec:quivers-basic}. This section recalls what a quiver records about an algebra, and how its representation type, which indicates whether canonical forms can exist, is read off from it.

\paragraph{What a quiver records.} A finite-dimensional algebra $A$ has two layers of discrete data. The first is its semisimple quotient $A/\rad(A)$, a product of matrix algebras over division algebras (Theorem~\ref{thm:wedderburn-artin}). The second is the bimodule $\rad(A)/\rad(A)^2$ over $A/\rad(A)$, which is a direct sum of simple bimodules with multiplicities. Together they determine $A/\rad(A)^2$, and over an algebraically closed field they are recorded, up to Morita equivalence, by a quiver: a vertex for each simple factor and, from $i$ to $j$, as many arrows as the multiplicity of the simple bimodule between the factors $i$ and $j$. Over such a field, the basic algebras whose radical squares to zero are exactly the algebras $kQ/R_Q^2$ (Section~\ref{sec:quivers-basic}), so they are classified by discrete data. The next layer is not discrete: an algebra $kQ/I$ with $R_Q^3\subseteq I\subseteq R_Q^2$ is determined by the subspace $I/R_Q^3$ of the span of the paths of length two, and such subspaces vary in continuous families. By Gabriel's theorem (Section~\ref{sec:quivers-basic}), relations are all that is needed beyond the quiver.

\subsection{Representation Type and the Tits Form}
\label{sec:tits_form}
Whether a decomposition can have a canonical form depends on the \emph{representation type} of the quiver involved. The Tits form gives a simple criterion, which can be understood by counting degrees of freedom. In this subsection $k$ is algebraically closed.

Fix a dimension vector $\mathbf d=(d_i)_{i\in Q_0}$. The representations of $Q$ with dimension vector $\mathbf d$ form the affine space
\[ \Rep(Q,\mathbf d)=\bigoplus_{\alpha\in Q_1}\Hom\bigl(k^{d_{s(\alpha)}},k^{d_{t(\alpha)}}\bigr),\qquad \dim\Rep(Q,\mathbf d)=\sum_{\alpha:i\to j}d_id_j, \]
and the isomorphism classes are the orbits of the group $\GL(\mathbf d)=\prod_i\GL(d_i)$, of dimension $\sum_id_i^2$, acting by change of basis. The group dimension is the number of degrees of freedom available for choosing bases, and the dimension of $\Rep(Q,\mathbf d)$ is the number of degrees of freedom that a representation has.

\begin{definition}[Tits form]\label{def:tits}
The \emph{Tits form} of $Q$ is the quadratic form
\[ q_Q(\mathbf d)=\sum_{i\in Q_0}d_i^2-\sum_{\alpha:i\to j}d_id_j=\dim\GL(\mathbf d)-\dim\Rep(Q,\mathbf d). \]
\end{definition}

\begin{proposition}\label{prop:tits-infinite}
If $q_Q(\mathbf d)\leq0$ for some non-zero dimension vector $\mathbf d$, then $Q$ has infinitely many isomorphism classes of representations of dimension vector $\mathbf d$.
\end{proposition}
\begin{proof}
The scalar matrices $\lambda I\in\GL(\mathbf d)$ act trivially, so every orbit has dimension at most $\dim\GL(\mathbf d)-1<\dim\Rep(Q,\mathbf d)$. A finite union of such orbits cannot fill $\Rep(Q,\mathbf d)$.
\end{proof}

In numerical terms: if $q_Q(\mathbf d)<0$, the group has too few degrees of freedom to bring an arbitrary representation into a simple normal form such as a diagonal one; the isomorphism classes then depend on at least $1-q_Q(\mathbf d)$ continuous parameters. What distinguishes the representation types is how this number grows with the dimension vector. The precise classification is as follows. A connected quiver is called
\begin{itemize}
    \item \emph{of finite type} if it has only finitely many indecomposable representations up to isomorphism. By Gabriel's classification of quivers of finite type~\cite{Gabriel1972}, this holds if and only if its underlying graph is a Dynkin diagram of type $A_n$, $D_n$, $E_6$, $E_7$ or $E_8$, if and only if $q_Q$ is positive definite; the indecomposable representations then correspond to the positive roots via their dimension vectors.
    \item \emph{tame} if it is not of finite type, but for each dimension vector the indecomposable representations form finitely many one-parameter families. This holds if and only if the underlying graph is an extended Dynkin diagram $\tilde A_n$, $\tilde D_n$, $\tilde E_6$, $\tilde E_7$ or $\tilde E_8$, if and only if $q_Q$ is positive semidefinite but not definite~\cite{DlabRingel1976}. The loop quiver, whose representations are classified by the Jordan normal form, is an example.
    \item \emph{wild} otherwise, in which case $q_Q$ is not positive semidefinite. Classifying the representations of a wild quiver contains the problem of classifying pairs of matrices up to simultaneous similarity, and is considered intractable~\cite{Kac1980,Drozd1980}.
\end{itemize}
The \emph{cycle rank} of a quiver is $|Q_1|-|Q_0|+c$, where $c$ is the number of connected components of its underlying graph; it is the number of independent cycles of that graph, counting loops and multiple arrows as cycles. The Dynkin diagrams are trees, and among the extended Dynkin diagrams only $\tilde A_n$, a single cycle, is not a tree, so a connected quiver of cycle rank at least $2$ is wild. Nearly all quivers are wild. For example, the Kronecker quiver with two arrows $u\rightrightarrows v$ has $q_Q(\mathbf d)=(d_u-d_v)^2\geq0$ and is tame (its representations are matrix pencils, classified by the Kronecker canonical form), whereas with three arrows $q_Q(1,1)=-1$ and the quiver is wild. Likewise, the quiver with one vertex and two loops has $q_Q(d)=-d^2$ and is wild. The Tits form also has a homological interpretation, which we do not use: $q_Q(\mathbf{dim}\,V)=\dim_k\End(V)-\dim_k\mathrm{Ext}^1(V,V)$ for every representation $V$~\cite{AssemSimsonSkowronski2006}.

\section{Symmetric Quivers}
\label{sec:symmetric-quivers}
An involution of an algebra $A$ permutes the simple components of $A/\rad(A)$, which are the vertices of its quiver, and, because it reverses products, it maps an arrow from $i$ to $j$ to an arrow from the image of $j$ to the image of $i$. A quiver that describes a $*$-algebra must record this.

\begin{definition}[Symmetric quiver]\label{def:symmetric-quiver}
A \emph{symmetric quiver} is a quiver $Q$ together with a map $\dagger:Q_0\sqcup Q_1\to Q_0\sqcup Q_1$ that maps vertices to vertices and arrows to arrows, satisfies $(x^\dagger)^\dagger=x$, and reverses arrows: if $\alpha:i\to j$ then $\alpha^\dagger:j^\dagger\to i^\dagger$.
\end{definition}

Symmetric quivers were introduced by Derksen and Weyman~\cite{DerksenWeyman2002}. Their definition also requires every arrow from a vertex $i$ to $i^\dagger$ to be fixed; we impose this only for signed quivers (Definition~\ref{def:signed-quiver}), whose signs allow a pair of exchanged arrows from $i$ to $i^\dagger$ to be replaced by two fixed ones. Each vertex is either fixed by the involution or exchanged with another vertex, and the same holds for arrows. An arrow $\alpha:i\to j$ with $\alpha^\dagger=\alpha$ satisfies $j=i^\dagger$; it is a loop if $i$ is fixed.

The path algebra of a symmetric quiver is a $*$-algebra, with the $k$-linear involution
\[ e_v^*=e_{v^\dagger},\qquad (\alpha_\ell\cdots\alpha_1)^*=\alpha_1^\dagger\cdots\alpha_\ell^\dagger . \]
This involution preserves $R_Q^2$, so $kQ/R_Q^2$ is a $*$-algebra as well. These algebras are easy to compute with, since the paths form a basis and the involution is path reversal composed with $\dagger$. After complexification, the path $*$-algebra of every doubly-signed quiver has the same form, except that each path carries a coefficient block of size at most two and the involution carries signs (Proposition~\ref{prop:block-model}).

For example, a loop $x$ at a fixed vertex, with $x^\dagger=x$, gives the polynomial ring $k[x]$ with $x^*=x$. Two fixed vertices $u,v$ with an arrow $a:u\to v$ and its partner $a^\dagger:v\to u$ give the four-dimensional algebra $kQ/R_Q^2$ of Example~\ref{ex:path-algebras}(3) with $a^*=a^\dagger$, the algebra of the singular value decomposition (Section~\ref{sec:unification-example}).

\section{Signed Quivers}
\label{sec:signed_quivers}
A symmetric quiver does not determine the involution of the algebras it describes. By the Skolem--Noether theorem (Theorem~\ref{thm:skolem-noether}), a $k$-linear involution of $M_n(k)$, for $k$ algebraically closed of characteristic different from $2$, has the form $X\mapsto GX^TG^{-1}$ with $G^T=\pm G$: it is the adjoint involution of a symmetric or of a skew-symmetric form, and it is said to be of \emph{orthogonal} or of \emph{symplectic type}~\cite{KMRT1998}. Moreover, for a suitable choice of the arrows, the involution of the algebra maps each arrow fixed by $\dagger$ to plus or minus itself. Signs on the fixed vertices and arrows record these two kinds of information.

\begin{definition}[Signed quiver]\label{def:signed-quiver}
A \emph{signed quiver} is a symmetric quiver $(Q,\dagger)$ together with a sign $\epsilon(x)\in\{+1,-1\}$ for every vertex and every arrow $x$ that is fixed by the involution, $x^\dagger=x$. It is convenient to extend $\epsilon$ to all of $Q_0\sqcup Q_1$ by setting $\epsilon(x)=0$ when $x^\dagger\neq x$; then $\epsilon:Q_0\sqcup Q_1\to\{-1,0,+1\}$ is characterised by
\[ \epsilon(x)^2=\begin{cases}1&\text{if }x^\dagger=x,\\0&\text{if }x^\dagger\neq x.\end{cases} \]
We also impose the following \emph{normalisation}: every arrow from a vertex $v$ to $v^\dagger$ is fixed by the involution.
\end{definition}

The normalisation removes a redundancy. If $\alpha:v\to v^\dagger$, then $\alpha^\dagger$ also goes from $v$ to $v^\dagger$, and if $\alpha^\dagger\neq\alpha$ the two arrows span the same space as $\alpha+\alpha^\dagger$ and $\alpha-\alpha^\dagger$, which are fixed by the involution with signs $+1$ and $-1$. So a pair of exchanged arrows from $v$ to $v^\dagger$ can always be replaced by two fixed arrows with opposite signs.

Signed quivers describe the $*$-algebras $A$ with $\rad(A)^2=0$ and a $k$-linear involution over an algebraically closed field $k$ of characteristic different from $2$; over $\RR$ they do not suffice (Section~\ref{sec:dsq}). The signs have the following meaning. A vertex $v$ with $v^\dagger=v$ corresponds to a simple component of $A/\rad(A)$ that is stable under the involution, and the vertex sign $\epsilon(v)=+1$ or $-1$ records whether the involution is of orthogonal or of symplectic type there. A pair of vertices $v\neq v^\dagger$ corresponds to two simple components exchanged by the involution, as in $\twoR$. A fixed arrow $\alpha$ goes from a vertex $v$ to $v^\dagger$, and the involution maps it to $\epsilon(\alpha)\alpha$, so that $\alpha$ behaves like a symmetric ($+1$) or skew-symmetric ($-1$) matrix. Arrows with $\alpha^\dagger\neq\alpha$ come in pairs.

\begin{remark}\label{rem:shmelkin}
The name is taken from Shmelkin~\cite{Shmelkin2006}, who records the same data in a folded form: a quiver with signs in $\{-1,0,+1\}$ on its vertices and arrows, from which the symmetric quiver with its signs is obtained by adding a partner for each unsigned vertex, other than his \emph{twins}, and for each unsigned arrow~\cite[Section~3]{Shmelkin2006}. Our signed quivers are the unfolded form, in which the involution is part of the data. His sign on a loop $\alpha$ at a fixed vertex $v$ records whether the bilinear form $(x,y)\mapsto\langle V_\alpha x,y\rangle$ of a $*$-representation (Definition~\ref{def:star-rep}) is symmetric or skew-symmetric; this is $\epsilon(v)\epsilon(\alpha)$ in our notation, so the two conventions differ at loops on vertices of sign $-1$.
\end{remark}

Figure~\ref{fig:signed-quiver-example} shows four signed quivers. In the diagrams, vertices are drawn as circles labelled with their name and, if fixed by the involution, their sign; dashed grey lines connect vertices or arrows exchanged by $\dagger$; arrows fixed by $\dagger$ are labelled with their sign.

\begin{figure}[ht]
\centering
\begin{minipage}{0.45\linewidth}\centering
\begin{tikzpicture}
  \node[qvertex] (v) {$v$\,\scriptsize$(+)$};
  \path[qarr] (v) edge [loop above, min distance=12mm] node[above] {$x\ (+)$} (v);
\end{tikzpicture}\\[1mm]
{\small (a) Symmetric loop: symmetric matrices under orthogonal similarity.}
\end{minipage}\hfill
\begin{minipage}{0.45\linewidth}\centering
\begin{tikzpicture}
  \node[qvertex] (v) {$v$\,\scriptsize$(+)$};
  \path[qarr] (v) edge [loop above, min distance=12mm] node[above] {$h\ (-)$} (v);
\end{tikzpicture}\\[1mm]
{\small (b) Skew loop: skew-symmetric matrices under orthogonal similarity.}
\end{minipage}

\vspace{5mm}
\begin{minipage}{0.45\linewidth}\centering
\begin{tikzpicture}
  \node[qvertex] (u) {$u$};
  \node[qvertex, right=22mm of u] (w) {$u^\dagger$};
  \path[qarr] (u) edge node[above] {$a\ (+)$} (w);
  \draw[qinv] (u.south) to[bend right=35] node[below]{\scriptsize$\dagger$} (w.south);
\end{tikzpicture}\\[1mm]
{\small (c) Exchanged vertices with a symmetric arrow: symmetric matrices under congruence.}
\end{minipage}\hfill
\begin{minipage}{0.45\linewidth}\centering
\begin{tikzpicture}
  \node[qvertex] (u) {$u$\,\scriptsize$(+)$};
  \node[qvertex, right=22mm of u] (v) {$v$\,\scriptsize$(+)$};
  \path[qarr] (u) edge[bend left=25] node[above] {$a$} (v);
  \path[qarr] (v) edge[bend left=25] node[below] {$a^\dagger$} (u);
\end{tikzpicture}\\[1mm]
{\small (d) A pair of exchanged arrows: rectangular matrices under $W\mapsto RWP^T$ (the SVD).}
\end{minipage}
\caption{Four signed quivers and the matrix problems they encode (Chapter~\ref{chap:indexing}). Signs of fixed vertices and arrows are shown in parentheses; the dashed line in (c) indicates that $\dagger$ exchanges $u$ and $u^\dagger$. The problems named are the ones obtained over $\RR$, when each fixed vertex of sign $+1$ carries a real inner product, as in the doubly-signed quivers of Figure~\ref{fig:dsq-examples}.}
\label{fig:signed-quiver-example}
\end{figure}

\subsection{\texorpdfstring{$*$}{*}-Representations}
\label{sec:star-reps}
For a signed quiver, the analogue of a representation is a representation together with compatible forms. Write $\tilde\epsilon(v)=\epsilon(v)$ if $v^\dagger=v$ and $\tilde\epsilon(v)=1$ otherwise, and choose signs $\kappa_\alpha$ for the arrows with $\kappa_\alpha=\epsilon(\alpha)$ if $\alpha^\dagger=\alpha$ and
\[ \kappa_\alpha\kappa_{\alpha^\dagger}=\tilde\epsilon(s(\alpha))\,\tilde\epsilon(t(\alpha))\qquad\text{if }\alpha^\dagger\neq\alpha , \]
for instance $\kappa_\alpha=1$ for one arrow of each exchanged pair. This is the first condition of \eqref{eq:sign-conditions} on the signs of the block model of Section~\ref{sec:path-star-algebras}, in which $g_v=\tilde\epsilon(v)$.

\begin{definition}\label{def:star-rep}
A \emph{$*$-representation} of a signed quiver $(Q,\dagger,\epsilon)$ over $k$ is a representation $V$ of $Q$ (Definition~\ref{def:representation}) together with a non-degenerate bilinear form $\langle\cdot,\cdot\rangle$ on $V=\bigoplus_vV_v$ such that
\begin{enumerate}
  \item $\langle V_u,V_v\rangle=0$ unless $v=u^\dagger$, and $\langle y,x\rangle=\langle x,y\rangle$ for $x\in V_u$, $y\in V_{u^\dagger}$ when $u\neq u^\dagger$;
  \item for each vertex $v$ with $v^\dagger=v$, the restriction of the form to $V_v$ is symmetric if $\epsilon(v)=1$ and skew-symmetric if $\epsilon(v)=-1$;
  \item for each arrow $\alpha$ and all $x\in V_{s(\alpha)}$, $y\in V_{t(\alpha)^\dagger}$: $\langle V_\alpha x,y\rangle=\kappa_\alpha\,\langle x,V_{\alpha^\dagger}y\rangle$.
\end{enumerate}
An isomorphism of $*$-representations is an isomorphism of representations that is an isometry.
\end{definition}

Condition (3) says that $V_{\alpha^\dagger}$ is the adjoint of $V_\alpha$, up to the sign $\kappa_\alpha$. The condition on $\kappa_\alpha\kappa_{\alpha^\dagger}$ is forced: writing (3) for $\alpha$ and for $\alpha^\dagger$ and using the symmetry of the forms gives $\tilde\epsilon(s(\alpha))\tilde\epsilon(t(\alpha))\langle y,V_\alpha x\rangle=\kappa_\alpha\kappa_{\alpha^\dagger}\langle y,V_\alpha x\rangle$, so any other choice would force $V_\alpha=0$, for example for an arrow from a vertex of sign $+1$ to one of sign $-1$. The choice of $\kappa_\alpha$ within an exchanged pair does not matter, since replacing $V_{\alpha^\dagger}$ by $-V_{\alpha^\dagger}$ changes it. For example, for the signed quiver of Figure~\ref{fig:signed-quiver-example}(a), a $*$-representation is a vector space with a non-degenerate symmetric form and a self-adjoint endomorphism; over $\RR$ and with a positive definite form, its classification is the spectral theorem. For Figure~\ref{fig:signed-quiver-example}(d), it consists of two inner product spaces and a linear map between them together with its adjoint; its classification is the singular value decomposition. Chapter~\ref{chap:indexing} proves a general version of these statements in the language of matrices over $*$-algebras. For the doubly-signed quivers of Section~\ref{sec:dsq} one additionally requires a real structure, compatible with the conjugation, on the complexified representation.

The count of degrees of freedom of Section~\ref{sec:tits_form} applies to $*$-representations as well, with $\GL(\mathbf d)$ replaced by the product $G$ of the groups at the vertices (an orthogonal or a symplectic group at each fixed vertex, according to its sign, and a general linear group for each pair of exchanged vertices) and $\Rep(Q,\mathbf d)$ replaced by the space of $*$-representations with given forms. Every orbit has dimension at most $\dim G$, so there are infinitely many isomorphism classes as soon as $\dim G$ is smaller than the dimension of that space. An analogous trichotomy, governed by root systems of signed quivers, is developed in~\cite{Shmelkin2006,DerksenWeyman2002}.

\section{Doubly-Signed Quivers}
\label{sec:dsq}
Over $\RR$, a simple algebra is $M_n(D)$ with $D\in\{\RR,\CC,\HH\}$. To describe a real algebra $B$ by a quiver, we pass to its complexification $B_\CC=B\otimes_\RR\CC$, on which complex conjugation acts by $\iota(b\otimes z)=b\otimes\bar z$. This conjugation is a conjugate-linear \emph{automorphism}, not an anti-automorphism, so it maps an arrow from $i$ to $j$ to an arrow from the conjugate of $i$ to the conjugate of $j$. Each simple component of $B$ gives rise to simple components of $B_\CC$ as follows: $M_n(\RR)$ gives $M_n(\CC)$, fixed by $\iota$; $M_n(\HH)$ gives $M_{2n}(\CC)$, again fixed by $\iota$; and $M_n(\CC)$ gives $M_n(\CC)\times M_n(\CC)$, two components exchanged by $\iota$. On a component fixed by $\iota$, the Skolem--Noether theorem shows that $\iota(X)=G\bar XG^{-1}$ with $G\bar G=\pm I$, the sign being $+1$ for real and $-1$ for quaternionic type. A doubly-signed quiver records this conjugation and its sign, alongside the involution and the signs of a signed quiver.

\begin{definition}[Doubly-signed quiver]\label{def:doubly-signed}
A \emph{doubly-signed quiver} is a signed quiver $Q$, whose involution $x\mapsto x^{\top}$ we call the \emph{transpose} and whose sign function we denote by $\epsilonform$, together with
\begin{itemize}
  \item a map $\overline{(\cdot)}:Q_0\sqcup Q_1\to Q_0\sqcup Q_1$, the \emph{conjugation}, that maps vertices to vertices and arrows to arrows, satisfies $\bar{\bar x}=x$, and preserves arrows: if $\alpha:i\to j$ then $\bar\alpha:\bar i\to\bar j$;
  \item a sign $\epsilondiv(v)\in\{+1,-1\}$ for every vertex $v$ with $\bar v=v$, extended by $\epsilondiv(v)=0$ when $\bar v\neq v$;
\end{itemize}
subject to the following rules:
\begin{enumerate}
    \item the transpose and the conjugation commute; their composite is written $(\cdot)^H=\overline{(\cdot)}^{\,\top}$ and called the \emph{conjugate transpose};
    \item $\epsilondiv(v^{\top})=\epsilondiv(v)$ for every vertex $v$;
    \item $\epsilonform(\bar x)=\epsilonform(x)$ for every vertex or arrow $x$;
    \item if both end points of an arrow $\alpha:i\to j$ are fixed by the conjugation, then $\bar\alpha=\alpha$ when $\epsilondiv(i)=\epsilondiv(j)$, and $\bar\alpha\neq\alpha$ otherwise;
    \item every arrow $\alpha$ from a vertex $v$ to $v^H$ is fixed by the conjugate transpose, that is, $\bar\alpha=\alpha^{\top}$.
\end{enumerate}
\end{definition}

\begin{figure}[ht]
\centering
\begin{minipage}{0.47\linewidth}\centering
\begin{tikzpicture}
  \node[qvertex] (u) {\shortstack{$u$\\[-1pt]\scriptsize$\RR,+$}};
  \node[qvertex, right=24mm of u] (v) {\shortstack{$v$\\[-1pt]\scriptsize$\RR,+$}};
  \path[qarr] (u) edge[bend left=25] node[above] {$a$} (v);
  \path[qarr] (v) edge[bend left=25] node[below] {$a^{\top}$} (u);
\end{tikzpicture}\\[1mm]
{\small (a) Real SVD. Both vertices are fixed by the transpose and by conjugation; $a$ and $a^{\top}$ are exchanged by the transpose and fixed by conjugation.}
\end{minipage}\hfill
\begin{minipage}{0.47\linewidth}\centering
\begin{tikzpicture}
  \node[qvertex] (u) {\shortstack{$u$\\[-1pt]\scriptsize$\HH,-$}};
  \node[qvertex, right=24mm of u] (v) {\shortstack{$v$\\[-1pt]\scriptsize$\HH,-$}};
  \path[qarr] (u) edge[bend left=25] node[above] {$a$} (v);
  \path[qarr] (v) edge[bend left=25] node[below] {$a^{\top}$} (u);
\end{tikzpicture}\\[1mm]
{\small (b) Quaternionic SVD. As in (a), but both vertices are quaternionic with $\epsilondiv=\epsilonform=-1$.}
\end{minipage}

\vspace{5mm}
\begin{minipage}{0.47\linewidth}\centering
\begin{tikzpicture}
  \node[qvertex] (u) {$u$};
  \node[qvertex, right=30mm of u] (ub) {$\bar u$};
  \node[qvertex, below=16mm of u] (v) {$v$};
  \node[qvertex, below=16mm of ub] (vb) {$\bar v$};
  \path[qarr] (u) edge[bend right=25] node[left] {$a$} (v);
  \path[qarr] (v) edge[bend right=25] node[right] {$\bar a^{\top}$} (u);
  \path[qarr] (ub) edge[bend left=25] node[right] {$\bar a$} (vb);
  \path[qarr] (vb) edge[bend left=25] node[left] {$a^{\top}$} (ub);
  \draw[qbar] (u) -- node[above]{\scriptsize $\bar{\ }$, $\top$} (ub);
  \draw[qbar] (v) -- node[below]{\scriptsize $\bar{\ }$, $\top$} (vb);
\end{tikzpicture}\\[1mm]
{\small (c) Complex SVD. Each complex vertex appears as a pair exchanged both by conjugation and by the transpose (unitary type); the arrows form two $2$-cycles exchanged by conjugation.}
\end{minipage}\hfill
\begin{minipage}{0.47\linewidth}\centering
\begin{tikzpicture}
  \node[qvertex] (v) {$v$};
  \node[qvertex, right=24mm of v] (vb) {$\bar v$};
  \path[qarr] (v) edge[bend left=25] node[above] {$\alpha\ (+)$} (vb);
  \path[qarr] (vb) edge[bend left=25] node[below] {$\bar\alpha\ (+)$} (v);
  \draw[qbar] (v.west) to[out=-90,in=-90,looseness=1.3] node[below]{\scriptsize $\bar{\ }$, $\top$} (vb.east);
\end{tikzpicture}\\[1mm]
{\small (d) Takagi factorisation. A unitary-type complex vertex with an arrow from $v$ to $\bar v=v^{\top}$ that is fixed by the transpose with sign $+1$; conjugation exchanges $\alpha$ and $\bar\alpha$.}
\end{minipage}
\caption{Doubly-signed quivers for four classical decompositions, drawn after complexification. Vertex labels give the division algebra ($\RR$, $\CC$ or $\HH$) and $\epsilonform$ for vertices fixed by the respective involutions; dotted blue lines join vertices exchanged by conjugation $\bar{\ }$ and by the transpose $\top$. Chapter~\ref{chap:indexing} derives the corresponding matrix problems.}
\label{fig:dsq-examples}
\end{figure}

Thus a doubly-signed quiver carries two involutions: the transpose, which reverses arrows, and the conjugation, which preserves them and plays the role of complex conjugation. It carries two sign functions: $\epsilondiv$ on the vertices fixed by the conjugation, which distinguishes $\RR$ from $\HH$, and $\epsilonform$ on the vertices and arrows fixed by the transpose, which distinguishes symmetric from skew-symmetric behaviour. The scalars at a vertex $v$ are real if $\epsilondiv(v)=+1$, complex if $\epsilondiv(v)=0$ (in which case $v$ and $\bar v$ together form one vertex, with complex scalars, of the real algebra), and quaternionic if $\epsilondiv(v)=-1$.

Rules~(1) to~(3) express the facts that conjugation and transposition commute, that transposition does not change the type of division algebra, and that conjugation does not change the symmetry type. Rule~(4) reflects how conjugation acts on the arrows between two vertices fixed by it. Between two vertices with the same $\epsilondiv$ it acts as a real structure on the space of arrows, so a basis of fixed arrows can be chosen, and we always choose one; between a real and a quaternionic vertex its square is $-1$, so it fixes no arrow, and the arrows come in exchanged pairs. Rule~(5) is a normalisation of the same kind: the conjugate transpose acts on the arrows from $v$ to $v^H$ as a real structure, so these arrows can be chosen fixed. Both facts are part of Theorem~\ref{thm:classification}(3). When $v$ is fixed by the conjugation, $v^H=v^{\top}$, and rule~(5) follows from the normalisation of Definition~\ref{def:signed-quiver} and rule~(4).

\begin{table}[htbp]
\centering
\caption{The vertex types of a doubly-signed quiver, the corresponding $*$-simple real $*$-algebras of minimal size, and the unitary groups $\Un_n$ of the matrix $*$-algebras over them. Indefinite unitary groups such as $\Or(p,q)$, $\Un(p,q)$ and $\Sp(p,q)$ arise from the same vertex types when the matrix algebra carries an indefinite form (Remark~\ref{rem:indefinite}).}
\label{tab:vertex-types}
\small
\renewcommand{\arraystretch}{1.2}
\begin{tabular}{@{}lll@{}}
\toprule
\textbf{Vertex data} & \textbf{$*$-simple $*$-algebra (minimal size)} & \textbf{Unitary group of $M_n$} \\
\midrule
$\epsilondiv=+1$, $\epsilonform=+1$ & $\RR$, trivial involution & $\Or(n)$ \\
$\epsilondiv=+1$, $\epsilonform=-1$ & $M_2(\RR)$, symplectic involution & $\Sp(2n,\RR)$ \\
$\epsilondiv=-1$, $\epsilonform=-1$ & $\HH$, quaternion conjugation & $\Sp(n)$ \\
$\epsilondiv=-1$, $\epsilonform=+1$ & $\HH$, involution of orthogonal type & $\Or^*(2n)$ \\
$\epsilondiv=0$, $\epsilonform=+1$ & $\CC$, trivial involution & $\Or(n,\CC)$ \\
$\epsilondiv=0$, $\epsilonform=-1$ & $M_2(\CC)$, symplectic involution & $\Sp(2n,\CC)$ \\
$\epsilondiv=0$, $v^{\top}=\bar v$ & $\CC$, complex conjugation & $\Un(n)$ \\
$\epsilondiv=+1$, $v^{\top}\neq v$ & $\RR\times\RR^{\mathrm{op}}$, exchange involution & $\GL_n(\RR)$ \\
$\epsilondiv=-1$, $v^{\top}\neq v$ & $\HH\times\HH^{\mathrm{op}}$, exchange involution & $\GL_n(\HH)$ \\
$\epsilondiv=0$, $v^{\top}\notin\{v,\bar v\}$ & $\CC\times\CC^{\mathrm{op}}$, exchange involution & $\GL_n(\CC)$ \\
\bottomrule
\end{tabular}
\end{table}

The vertex types of a doubly-signed quiver correspond exactly to the ten families of classical groups over $\RR$ (Table~\ref{tab:vertex-types}); this is the classification of involutions on simple real algebras~\cite{KMRT1998,Weyl1939,Helgason1978}. For the doubly-signed quiver of a real $*$-algebra (Theorem~\ref{thm:classification}), the unitary group at each vertex is the group in the table. Note the pattern: when $\epsilondiv=\epsilonform$ the group is compact ($\Or(n)$, $\Sp(n)$), whereas when $\epsilondiv=-\epsilonform$ it is not, except for $\Or^*(2)$, which is a circle group.

In labels we write the type of a vertex as a pair $(D,\epsilonform)$, with the division algebra $D=\RR,\CC,\HH$ in place of $\epsilondiv=+1,0,-1$, and with $\epsilonform=0$ for the exchange types, in which $v^{\top}\notin\{v,\bar v\}$. The type with $v^{\top}=\bar v$ is written $(\CC,\mathrm{u})$, for unitary, and in diagrams the signs $\pm1$ are abbreviated to $\pm$. For example, $(\RR,+1)$ is the type of $\RR$ with the trivial involution, and $(\RR,0)$ that of $\RR\times\RR^{\mathrm{op}}$.
\section{Real \texorpdfstring{$*$}{*}-Algebras with Radical Square Zero}
\label{sec:classification}
We now turn to real $*$-algebras, beginning with those whose radical squares to zero. We show that each of them has a doubly-signed quiver, that the signs carry the information described in Section~\ref{sec:dsq}, that every doubly-signed quiver arises in this way, and that a basic such algebra is determined by its doubly-signed quiver. Two lemmas are needed first: a $*$-stable version of the Wedderburn--Malcev theorem, which has an elementary proof when the radical squares to zero, and a lemma on lifting unitaries.

\begin{lemma}\label{lem:star-wedderburn}
Let $A$ be a finite-dimensional $*$-algebra over a field of characteristic $0$ with $R=\rad(A)$ and $R^2=0$. Then $A$ has a Wedderburn--Malcev decomposition $A=S\oplus R$ with $S^*=S$.
\end{lemma}
\begin{proof}
By Lemma~\ref{lem:radical-stable}, $R^*=R$, so the involution descends to an involution $x\mapsto x^\natural$ of $\bar A=A/R$. By Theorem~\ref{thm:wedderburn-malcev} there is an algebra homomorphism $\varphi:\bar A\to A$ that is a section of the quotient map. Define $\psi(x)=\varphi(x^\natural)^*$. Then $\psi$ is again a homomorphic section: $\psi(xy)=\varphi(y^\natural x^\natural)^*=(\varphi(y^\natural)\varphi(x^\natural))^*=\psi(x)\psi(y)$, and $\psi(x)\equiv (x^\natural)^\natural=x$ modulo $R$. Let $D=\psi-\varphi$, which takes values in $R$. Since $R^2=0$,
\[ D(xy)=\psi(x)\psi(y)-\varphi(x)\varphi(y)=D(x)\varphi(y)+\varphi(x)D(y). \]
Put $j=\tfrac12(\varphi+\psi)=\varphi+\tfrac12D$. Then
\[ j(x)j(y)=\varphi(xy)+\tfrac12\bigl(D(x)\varphi(y)+\varphi(x)D(y)\bigr)+\tfrac14D(x)D(y)=\varphi(xy)+\tfrac12D(xy)=j(xy), \]
so $j$ is a homomorphic section, and $S=j(\bar A)$ is a subalgebra complementary to $R$. Finally $j(x)^*=\tfrac12(\varphi(x)^*+\varphi(x^\natural))=j(x^\natural)$, because $\psi(x^\natural)=\varphi(x)^*$; hence $S^*=S$.
\end{proof}

If $N$ is the radical of a finite-dimensional real $*$-algebra and $h\in1+N$ is self-adjoint, we put
\[ h^{-1/2}=\sum_{j\geq0}\binom{-1/2}{j}(h-1)^j . \]
The sum is finite because $N$ is nilpotent, $h^{-1/2}$ is a self-adjoint polynomial in $h$, and $h^{-1/2}hh^{-1/2}=1$, since this is an identity of formal power series in $h-1$.

\begin{lemma}\label{lem:unitary-lifting}
Let $A$ be a finite-dimensional real $*$-algebra with radical $N$.
\begin{enumerate}
  \item Every unitary of $A/N$ is the image of a unitary of $A$.
  \item If $f,f'\in A$ are self-adjoint idempotents whose images in $A/N$ are conjugate by a unitary of $A/N$, then $f'=ufu^*$ for a unitary $u\in A$. In particular $fAf\cong f'Af'$ as $*$-algebras.
\end{enumerate}
\end{lemma}
\begin{proof}
(1) If $v\in A$ lifts the unitary $\bar v$, then $h=v^*v\in1+N$ is self-adjoint, and $u=vh^{-1/2}$ lifts $\bar v$ and satisfies $u^*u=h^{-1/2}hh^{-1/2}=1$.

(2) By (1), the unitary of $A/N$ that conjugates the image of $f$ to that of $f'$ lifts to a unitary $u\in A$, and $g=ufu^*$ is a self-adjoint idempotent with $g\equiv f'$ modulo $N$. Then $z=f'g+(1-f')(1-g)\equiv1$ is invertible (Proposition~\ref{prop:invertibility}), $f'z=f'g=zg$, and $z^*z=gf'g+(1-g)(1-f')(1-g)$ commutes with $g$. Hence $w=z(z^*z)^{-1/2}$ is unitary, and since $(z^*z)^{-1/2}$ commutes with $g$, we get $wg=zg(z^*z)^{-1/2}=f'z(z^*z)^{-1/2}=f'w$, that is, $wgw^*=f'$. So $f'=(wu)f(wu)^*$, and $x\mapsto(wu)x(wu)^*$ restricts to a $*$-isomorphism $fAf\to f'Af'$.
\end{proof}

Let $B$ be a finite-dimensional real $*$-algebra with $\rad(B)^2=0$, and let $B_\CC=B\otimes_\RR\CC$ be its complexification. On $B_\CC$ there are two commuting structures: the $\CC$-linear extension of the involution of $B$, which we continue to write as $x\mapsto x^*$ and which is an anti-automorphism, and the conjugation $\iota=\mathrm{id}\otimes\overline{(\cdot)}$, which is a conjugate-linear automorphism. By Lemma~\ref{lem:star-wedderburn} applied to $B$, there is a $*$-stable Wedderburn--Malcev decomposition $B=S_B\oplus\rad(B)$, and then $B_\CC=S\oplus R$ with $S=S_B\otimes\CC$ and $R=\rad(B)\otimes\CC=\rad(B_\CC)$, where $S$ is stable under both $*$ and $\iota$. Write $S=\prod_iS_i$ with $S_i\cong M_{n_i}(\CC)$ simple, and let $1_i$ be the unit of $S_i$. Then $R=\bigoplus_{i,j}1_jR1_i$, and $1_jR1_i$ is an $S_j$--$S_i$-bimodule, hence isomorphic to a direct sum of $m_{ji}$ copies of the simple bimodule $M_{n_j\times n_i}(\CC)$.

\begin{definition}\label{def:quiver-of-B}
The \emph{quiver of $B_\CC$} has a vertex for each simple component $S_i$ and $m_{ji}$ arrows from $i$ to $j$.
\end{definition}

\begin{theorem}[Doubly-signed quivers of real $*$-algebras]\label{thm:classification}
Let $B$ be a finite-dimensional real $*$-algebra with $\rad(B)^2=0$, and let $Q$ be the quiver of $B_\CC$.
\begin{enumerate}
  \item The involution $*$ of $B_\CC$ induces an involution $(\cdot)^{\top}$ of $Q$ that reverses arrows, and the conjugation $\iota$ induces an involution $\overline{(\cdot)}$ of $Q$ that preserves arrows; these commute.
  \item For a vertex $i$ with $i^{\top}=i$, the restriction of $*$ to $S_i\cong M_{n_i}(\CC)$ is $X\mapsto GX^TG^{-1}$ with $G^T=\epsilonform(i)G$ and $\epsilonform(i)\in\{\pm1\}$. For a vertex with $\bar i=i$, the restriction of $\iota$ is $X\mapsto C\bar XC^{-1}$ for some $C$ with $C\bar C=\epsilondiv(i)I$ and $\epsilondiv(i)\in\{\pm1\}$; the simple component of $B/\rad(B)$ underlying $S_i$ is a matrix algebra over $\RR$ if $\epsilondiv(i)=1$ and over $\HH$ if $\epsilondiv(i)=-1$. Vertices with $\bar i\neq i$ come from simple components of $B/\rad(B)$ that are matrix algebras over $\CC$.
  \item On the multiplicity space of the arrows from $i$ to $i^{\top}$, the involution $*$ induces a linear involution. The dimensions of its $\pm1$-eigenspaces do not depend on the coordinates used to define it, and these eigenspaces define arrows fixed by the transpose with signs $\epsilonform(\alpha)=\pm1$. On the multiplicity space of the arrows between two vertices $i,j$ fixed by conjugation, $\iota$ induces a conjugate-linear map whose square is $\epsilondiv(i)\epsilondiv(j)$. On the multiplicity space of the arrows from $i$ to $i^H$, the composite $\iota\circ*$ induces a conjugate-linear map whose square is $1$, so these arrows can be chosen fixed by the conjugate transpose.
  \item With these data, $Q$ is a doubly-signed quiver $Q(B)$, which does not depend, up to isomorphism, on the choices made in its construction. It is unchanged when $B$ is replaced by a matrix algebra $M_n(B)$ with the conjugate-transpose involution, or by a corner $fBf$, where $f$ is a self-adjoint idempotent with $BfB=B$; for example, by a basic $*$-algebra of $B$ (Definition~\ref{def:basic-star}).
\end{enumerate}
Conversely, every doubly-signed quiver $Q$ arises in this way, from the basic $*$-algebra $B_Q$ of Definition~\ref{def:path-star-algebra}, whose semisimple quotient is a product of algebras from Table~\ref{tab:vertex-types}, and $B$ is determined by $Q(B)$ up to the replacements in (4), in the following sense: two real $*$-algebras with radical square zero and isomorphic doubly-signed quivers have $*$-isomorphic basic $*$-algebras. In particular, two basic real $*$-algebras with radical square zero and isomorphic doubly-signed quivers are $*$-isomorphic. If $A$ is a complex algebra with a $\CC$-linear involution, the same construction without the conjugation $\iota$ yields a signed quiver.
\end{theorem}

\begin{proof}
Throughout, $N=\rad(B)$ and $R=\rad(B_\CC)$. Since $R^2=0$, the product of an element of $B_\CC$ with an element of $R$ depends only on the image of the former in $B_\CC/R$. So $R$ is a bimodule over $B_\CC/R$, and the quotient map identifies $S$ with $B_\CC/R$, compatibly with $*$ and $\iota$. The constructions below use $S$ only through this identification, so they do not depend on the choice of $S_B$.

\emph{(1)} The involution $*$ is an anti-automorphism of $B_\CC$ that preserves $S$ and $R$. So it permutes the simple components $S_i$, defining $i\mapsto i^{\top}$, and maps $R_{ji}=1_jR1_i$ onto $R_{i^{\top}j^{\top}}$, reversing arrows. Likewise $\iota$ permutes the $S_i$, defining $i\mapsto\bar i$, and maps $R_{ji}$ onto $R_{\bar j\bar i}$, preserving arrows. These maps of the vertices are involutions, and they commute because $*$ and $\iota$ commute. The corresponding involutions of the arrows are constructed in Step~4 below.

\emph{(2)} If $i^{\top}=i$, then $*$ restricts to a $\CC$-linear anti-automorphism of $S_i\cong M_n(\CC)$ of order two. Composing with the matrix transpose gives an automorphism, which is inner by the Skolem--Noether theorem (Theorem~\ref{thm:skolem-noether}); so $X^*=GX^TG^{-1}$ for some invertible $G$. Since $*$ has order two, $G(G^T)^{-1}$ is central, $G^T=\lambda G$, and $G=(G^T)^T=\lambda^2G$ forces $\lambda=\pm1$. The case of $\iota$ is the same with complex conjugation in place of the transpose: $\iota(X)=C\bar XC^{-1}$ with $C\bar C=\mu I$. Conjugating this identity gives $\bar CC=\bar\mu I$, while $\bar C=\mu C^{-1}$ gives $\bar CC=\mu I$; hence $\mu$ is real, and rescaling $C$ by a positive real number gives $\mu=\pm1$. The fixed algebra of $\iota$ on $S_i$ is $M_n(\RR)$ if $\mu=1$ and $M_{n/2}(\HH)$ if $\mu=-1$~\cite{KMRT1998}. If $\bar i\neq i$, the fixed points of $\iota$ on $S_i\times S_{\bar i}$ form a copy of $M_n(\CC)$.

\emph{Step 1: coordinates.} Let $F$ be a $*$-simple factor of $S_B$. By the proof of Lemma~\ref{lem:star-simple}, whose existence part does not use the present theorem, $F\cong M_{r,\eta}(K)$ for some $r\geq1$, some signs $\eta$ and one of the ten algebras $K$ of Table~\ref{tab:vertex-types}. By the inspection of Table~\ref{tab:vertex-blocks} in the proof of Proposition~\ref{prop:block-model}(3), which does not use the present theorem either, $K$ is $*$-isomorphic to the algebra of fixed points of the map $\iota$ of \eqref{eq:block-star} in $\prod_vM_{d_v}(\CC)$, with the involution of \eqref{eq:block-star}, where $v$ runs through the complex vertices of one row of that table. Hence $F\otimes_\RR\CC\cong M_r(K\otimes_\RR\CC)\cong\prod_vM_r\bigl(M_{d_v}(\CC)\bigr)$. This matches the simple components of $F\otimes_\RR\CC$ with the vertices $v$ of the row, compatibly with the transpose and the conjugation. It also gives coordinates $S_i\cong M_{n_i}(\CC)$, with $n_i=rd_v$ for the vertex $v$ matched with $i$, in which, for $X\in S_i$,
\[ X^*=G_iX^TG_i^{-1}\in S_{i^{\top}},\qquad \iota(X)=C_i\overline XC_i^{-1}\in S_{\bar i}. \]
Here $G_i=\eta\otimes G_v$ and $C_i=I_r\otimes C_v$ are Kronecker products, and $\eta$ also denotes the diagonal matrix of the signs. Since $G_{v^{\top}}=G_v$, $C_{\bar v}=C_v$ and $\eta$ is real, \eqref{eq:vertex-data} gives
\[ G_{i^{\top}}=G_i,\qquad C_{\bar i}=C_i,\qquad G_i^T=g_iG_i,\qquad C_i\overline{C_i}=h_iI,\qquad C_{i^{\top}}\overline{G_i}\,C_i^T=G_{\bar i}, \]
where $g_i=g_v$ and $h_i=h_v$. In particular, the vertex data of (2) are those of the row: $\epsilonform(i)=g_i$ if $i^{\top}=i$, and $\epsilondiv(i)=h_i$ if $\bar i=i$. The numbers $g_i$ and $h_i$ are constant on orbits, and $g_i=1$, so that $G_i$ is symmetric, if $i^{\top}\neq i$.

\emph{Step 2: multiplicity spaces.} For vertices $i$ and $j$, let $M_{ji}=M_{n_j\times n_i}(\CC)$, an $S_j$--$S_i$-bimodule through these coordinates. The algebra $S_j\otimes S_i^{\mathrm{op}}\cong M_{n_jn_i}(\CC)$ is simple, $M_{ji}$ is its only simple module, and the endomorphisms of $M_{ji}$ are the scalars. Hence the evaluation map $M_{ji}\otimes U_{ji}\to R_{ji}$, $Y\otimes u\mapsto u(Y)$, is an isomorphism, where $U_{ji}=\Hom_{S_j\text{-}S_i}(M_{ji},R_{ji})$ is the multiplicity space; so $m_{ji}=\dim U_{ji}$. On $M_{ji}$ we define the reference maps
\[ \phi(Y)=G_iY^TG_j^{-1}\in M_{i^{\top}j^{\top}},\qquad \psi(Y)=C_j\overline YC_i^{-1}\in M_{\bar j\bar i}. \]
Writing elements of $S$ by their coordinates, we have $\phi(aYb)=b^*\phi(Y)a^*$ and $\psi(aYb)=\iota(a)\psi(Y)\iota(b)$ for $a\in S_j$ and $b\in S_i$. Hence, for $u\in U_{ji}$, the maps $\tau(u)=*\circ u\circ\phi^{-1}$ and $\gamma(u)=\iota\circ u\circ\psi^{-1}$ are again bimodule maps. They define a linear map $\tau:U_{ji}\to U_{i^{\top}j^{\top}}$ and a conjugate-linear map $\gamma:U_{ji}\to U_{\bar j\bar i}$ with
\[ u(Y)^*=\tau(u)\bigl(\phi(Y)\bigr),\qquad \iota\bigl(u(Y)\bigr)=\gamma(u)\bigl(\psi(Y)\bigr); \]
that is, under the evaluation isomorphisms, $*$ is $\phi\otimes\tau$ and $\iota$ is $\psi\otimes\gamma$. The identities of Step~1 give
\[ \phi(\phi(Y))=g_ig_j\,Y,\qquad \psi(\psi(Y))=h_ih_j\,Y,\qquad \psi(\phi(Y))=\phi(\psi(Y)), \]
as in the proof of Proposition~\ref{prop:block-model}(1) for a path of length one with signs $1$. Since $*$ and $\iota$ are commuting maps of order two, these identities give, for all $u\in U_{ji}$,
\begin{equation*}
\tau(\tau(u))=g_ig_j\,u,\qquad \gamma(\gamma(u))=h_ih_j\,u,\qquad \gamma(\tau(u))=\tau(\gamma(u)).\tag{R}
\end{equation*}
In particular, $\tau$ is a linear involution of $U_{i^{\top}i}$. If $\bar i=i$ and $\bar j=j$, then $\gamma$ is a conjugate-linear map of $U_{ji}$ with square $\epsilondiv(i)\epsilondiv(j)$. The map $\gamma\circ\tau$, which is induced by $\iota\circ*$, is a conjugate-linear map of $U_{i^Hi}$ with square $1$. This proves (3), apart from the choice of bases, which we make in Step~4.

\emph{Step 3: the signs are invariants.} The involution $\tau$ of $U_{i^{\top}i}$ does not depend on the coordinates, because the reference map $\phi$ on $M_{i^{\top}i}$ does not. If $i^{\top}=i$, then $M_{ii}=S_i$ as a bimodule, and $\phi$ is the involution of $S_i$. If $i^{\top}\neq i$, then $G_i$ is symmetric. After the coordinates of $S_{i^{\top}}$ are changed by $X\mapsto G_i^{-1}XG_i$, the involution maps $X\in S_i$ to $X^T\in S_{i^{\top}}$, and $\phi$ becomes the transpose of $M_{i^{\top}i}=M_{n_i}(\CC)$. Any two systems of coordinates in which $*$ maps $X\in S_i$ to $X^T\in S_{i^{\top}}$ differ by $X\mapsto PXP^{-1}$ on $S_i$ and $X\mapsto P^{-T}XP^T$ on $S_{i^{\top}}$. This changes $Y\in M_{i^{\top}i}$ into $P^{-T}YP^{-1}$, which commutes with the transpose. Hence the eigenspaces of $\tau$ on $U_{i^{\top}i}$, and their dimensions, depend only on $B$. (If $n_i$ is even, as at the vertices of type $(\HH,0)$, one could also use coordinates in which $*$ maps $X\in S_i$ to $\hat JX^T\hat J^{-1}$ for a skew-symmetric $\hat J$, with the reference map $Y\mapsto\hat JY^T\hat J^{-1}$. This would change the signs of all arrows from $i$ to $i^{\top}$. Step~1 excludes such coordinates, because there $G_i^T=g_iG_i$ with $g_i=1$; this is the convention of Table~\ref{tab:vertex-blocks}.)

\emph{Step 4: arrows.} We call bases $\mathcal B_{ji}$ of the spaces $U_{ji}$ \emph{adapted} if they come with involutions $\alpha\mapsto\alpha^{\top}$ and $\alpha\mapsto\bar\alpha$ of $\mathcal B=\bigsqcup_{i,j}\mathcal B_{ji}$, mapping $\mathcal B_{ji}$ onto $\mathcal B_{i^{\top}j^{\top}}$ and onto $\mathcal B_{\bar j\bar i}$ respectively, with the following properties. The two involutions commute and satisfy the normalisation of Definition~\ref{def:signed-quiver} and rules (4) and (5) of Definition~\ref{def:doubly-signed}, and there are signs $\kappa_\alpha,\mu_\alpha\in\{\pm1\}$ with
\[ \tau(\alpha)=\kappa_\alpha\,\alpha^{\top},\qquad \gamma(\alpha)=\mu_\alpha\,\bar\alpha\qquad(\alpha\in\mathcal B), \]
and $\mu_\alpha=1$ whenever $\bar\alpha=\alpha$. The arrows of $Q$ from $i$ to $j$ are then the elements of $\mathcal B_{ji}$, with these involutions. If $\alpha^{\top}=\alpha$, its sign $\epsilonform(\alpha)=\kappa_\alpha$ is the eigenvalue of $\tau$ at $\alpha$. This is the construction of statements (1) and (3).

Adapted bases exist. The group $\Gamma=\{1,\top,{}^-,H\}\cong(\ZZ/2)^2$ acts on pairs of vertices by $(i,j)^{\top}=(j^{\top},i^{\top})$, $\overline{(i,j)}=(\bar i,\bar j)$ and $(i,j)^H=(j^H,i^H)$, where $(i,j)$ stands for the arrows from $i$ to $j$. In each orbit of $\Gamma$ we choose a pair $(i,j)$ and a basis of $U_{ji}$ according to the stabiliser $\Gamma_{(i,j)}$ of the pair, and we transport the basis to the other pairs of the orbit by $\tau$ and $\gamma$:
\begin{enumerate}
  \item[(a)] $\Gamma_{(i,j)}=\{1\}$. Take any basis of $U_{ji}$. For each of its elements $\alpha$, put $\alpha^{\top}=\tau(\alpha)$, $\bar\alpha=\gamma(\alpha)$ and $\bar\alpha^{\top}=\overline{\alpha^{\top}}=\tau(\gamma(\alpha))=\gamma(\tau(\alpha))$. The two involutions act on these four elements as $\Gamma$ acts on the four pairs of the orbit.
  \item[(b)] $\Gamma_{(i,j)}=\{1,\top\}$, that is, $j=i^{\top}$ and $\bar i\neq i$. Take a basis of eigenvectors $\alpha=\alpha^{\top}$ of the involution $\tau$ of $U_{i^{\top}i}$, and put $\bar\alpha=\gamma(\alpha)$. By (R), $\tau(\bar\alpha)=\kappa_\alpha\bar\alpha$, so $\bar\alpha^{\top}=\bar\alpha$ and $\epsilonform(\bar\alpha)=\epsilonform(\alpha)$.
  \item[(c)] $\Gamma_{(i,j)}=\{1,{}^-\}$, that is, $\bar i=i$, $\bar j=j$ and $j\neq i^{\top}$. If $\epsilondiv(i)=\epsilondiv(j)$, then $\gamma$ is a real structure on $U_{ji}$; take a basis of fixed vectors $\alpha=\bar\alpha$. Otherwise $\gamma^2=-1$, and $\gamma$, together with the complex scalars, makes $U_{ji}$ a right vector space over $\HH$. A basis $\beta_1,\dots,\beta_k$ over $\HH$ gives the basis $\beta_1,\gamma(\beta_1),\dots,\beta_k,\gamma(\beta_k)$ over $\CC$, and we put $\bar\beta_l=\gamma(\beta_l)$. In both cases, put $\alpha^{\top}=\tau(\alpha)$ and $\overline{\alpha^{\top}}=\bar\alpha^{\top}$ for the elements $\alpha$ of the basis.
  \item[(d)] $\Gamma_{(i,j)}=\{1,H\}$, that is, $j=i^H$ and $\bar i\neq i$. Then $\gamma\circ\tau$ is a real structure on $U_{i^Hi}$; take a basis of fixed vectors $\alpha$. By (R), $\tau(\alpha)=\gamma(\gamma(\tau(\alpha)))=\gamma(\alpha)$, and we put $\alpha^{\top}=\bar\alpha=\tau(\alpha)$, whose transpose and conjugate are both $\alpha$.
  \item[(e)] $\Gamma_{(i,j)}=\Gamma$, that is, $j=i^{\top}$ and $\bar i=i$. Then $\tau$ is an involution and $\gamma$ is a real structure of $U_{i^{\top}i}$, because $h_ih_{i^{\top}}=1$. They commute, so $\gamma$ preserves the eigenspaces of $\tau$, and each eigenspace has a basis of fixed vectors of $\gamma$. We take these vectors, with $\alpha^{\top}=\bar\alpha=\alpha$.
\end{enumerate}
In each case, (R) shows that $\tau$ and $\gamma$ map each element of the resulting bases to plus or minus its image under the corresponding involution, with $\mu_\alpha=1$ when $\bar\alpha=\alpha$. The normalisation and rules (4) and (5) hold by construction. Indeed, the arrows from $i$ to $i^{\top}$ lie over pairs of type (b) or (e); the arrows between two vertices fixed by the conjugation, over pairs of type (c) or (e); and the arrows from $i$ to $i^H$, over pairs of type (d) or (e). In type (e), $\epsilondiv(i^{\top})=\epsilondiv(i)$ by Step~1, and $i^H=i^{\top}$.

\emph{(4)} The remaining rules of Definition~\ref{def:doubly-signed} follow from the fact that $*$ commutes with $\iota$. Since $*$ maps $S_i$ anti-isomorphically onto $S_{i^{\top}}$ and commutes with $\iota$, the component $S_i$ is $\iota$-stable if and only if $S_{i^{\top}}$ is, with the same type, so $\epsilondiv(i^{\top})=\epsilondiv(i)$. Similarly $\epsilonform(\bar i)=\epsilonform(i)$ for the vertices. For an arrow with $\alpha^{\top}=\alpha$, also $\bar\alpha^{\top}=\bar\alpha$, and (R) gives $\kappa_\alpha\mu_\alpha\bar\alpha=\gamma(\tau(\alpha))=\tau(\gamma(\alpha))=\mu_\alpha\kappa_{\bar\alpha}\bar\alpha$, so $\epsilonform(\bar\alpha)=\epsilonform(\alpha)$, for any adapted bases. So $Q$ is a doubly-signed quiver.

Up to isomorphism, a doubly-signed quiver is determined by its vertices, with their involutions and signs, by the numbers of arrows from $i$ to $j$, and by the numbers of arrows from $i$ to $i^{\top}$ with sign $+1$ and with sign $-1$. Indeed, by the normalisation and rules (3) to (5), the arrows over an orbit of $\Gamma$ are determined up to isomorphism by these numbers for one pair $(i,j)$ of the orbit. In the notation of the five cases above: in case (a), each arrow from $i$ to $j$ lies in a free orbit of the two involutions; in case (b), the arrows from $i$ to $i^{\top}$ are fixed by the transpose, and the conjugation maps them to arrows with the same signs; in case (c), they are all fixed by the conjugation or come in exchanged pairs, according to whether $\epsilondiv(i)=\epsilondiv(j)$; in case (d), they are fixed by the conjugate transpose; and in case (e), they are fixed by both involutions. For $Q(B)$, the vertex data depend only on $B$ by (2), the numbers of arrows are the integers $m_{ji}$, and the numbers of arrows from $i$ to $i^{\top}$ of each sign are the dimensions of the eigenspaces of $\tau$, which depend only on $B$ by Step~3. So $Q(B)$ does not depend, up to isomorphism, on the choices of $S_B$, of the coordinates and of the adapted bases.

Replacing $B$ by $M_n(B)$ replaces $S_i$ by $M_n(S_i)$ and $R_{ji}$ by $M_n(R_{ji})$. With the matrices $I_n\otimes G_i$ and $I_n\otimes C_i$, which are those of Step~1 for $M_n(B)$, the reference maps are the transpose and the entrywise conjugation of $M_n(\CC)$ tensored with $\phi$ and $\psi$. The involution and the conjugation of $M_n(R_{ji})$ are the same maps tensored with $*$ and $\iota$. So the spaces $U_{ji}$ and the maps $\tau$ and $\gamma$ do not change, and neither do the vertex data; hence $Q(M_n(B))\cong Q(B)$.

Now let $f$ be a self-adjoint idempotent with $BfB=B$. Its image in $B/N\cong S_B$ is a self-adjoint idempotent $f_0\in S_B$, and $f=uf_0u^*$ for a unitary $u\in B$ by Lemma~\ref{lem:unitary-lifting}(2); so we may assume that $f\in S_B$. Write $f=\sum_if_i$ with $f_i\in S_i$. Then $f_i^*=f_{i^{\top}}$ and $\iota(f_i)=f_{\bar i}$, and $f_i\neq0$ for every $i$ because $BfB=B$. The algebra $fBf$ has the radical $fNf$, whose square is zero, and the simple components $f_iS_if_i$, with the same transpose and conjugation. Moreover $f_jR_{ji}f_i\cong f_jM_{ji}f_i\otimes U_{ji}$, where $f_jM_{ji}f_i$ is the simple $f_jS_jf_j$--$f_iS_if_i$-bimodule. Since $\phi(f_jYf_i)=f_{i^{\top}}\phi(Y)f_{j^{\top}}$ and $\psi(f_jYf_i)=f_{\bar j}\psi(Y)f_{\bar i}$, the restrictions of $\phi$ and $\psi$ are reference maps for $fBf$ with the same properties, and they induce the same maps $\tau$ and $\gamma$ on the same spaces $U_{ji}$. At the arrows from $i$ to $i^{\top}$, they are the reference maps of Step~3 for $fBf$. If $i^{\top}=i$, the restriction of $\phi$ is the involution of $f_iS_if_i$. If $i^{\top}\neq i$, the coordinates of Step~3 can be chosen so that $f_i$ and $f_{i^{\top}}=f_i^T$ are both the diagonal matrix whose first $\operatorname{rank}f_i$ entries are $1$ and whose other entries are $0$; then the restriction of the transpose is the transpose of the upper left blocks. Finally, the vertex data do not change. If $i^{\top}=i$, the symmetric or skew-symmetric form defined by $G_i$ restricts to a form of the same type on the image of $f_i$, which is non-degenerate because $f_i^*=f_i$. If $\bar i=i$, the real or quaternionic structure $x\mapsto C_i\bar x$ preserves the image of $f_i$, because $\iota(f_i)=f_i$. Since $Q(B)$ is determined up to isomorphism by these data, $Q(fBf)\cong Q(B)$.

\emph{The converse.} Let $Q$ be a doubly-signed quiver, with admissible signs. By Proposition~\ref{prop:block-model}, $B_Q$ is a real $*$-algebra with radical $M_Q$, $M_Q^2=0$ and $S_Q\cong\prod_cK_c$, so that it is basic by Lemma~\ref{lem:star-simple}(3). Its complexification is $S_\CC\oplus M_\CC$ with the maps \eqref{eq:block-star}. These parts of the proof of Proposition~\ref{prop:block-model} do not use the present theorem. The block coordinates are coordinates as in Step~1, with $r=1$ and $\eta=1$. The copies of $M_{d_j\times d_i}(\CC)$ in $M_\CC$ that belong to the arrows $\alpha:i\to j$ show that these arrows form a basis of $U_{ji}$, and \eqref{eq:block-star} gives $\tau(\alpha)=\kappa_\alpha\alpha^{\top}$ and $\gamma(\alpha)=\mu_\alpha\bar\alpha$. Since $Q$ satisfies the rules of Definition~\ref{def:doubly-signed} and the signs are admissible, this basis is adapted, and the doubly-signed quiver that it defines is $Q$. Hence $Q(B_Q)=Q$.

\emph{Uniqueness.} Let $B$ be basic. Every self-adjoint idempotent of $S_B\cong B/N$ is central. So, by Lemma~\ref{lem:star-simple}(3), $r=1$ in Step~1 for every $*$-simple factor of $S_B$, and we may take $\eta=1$, since $M_{1,(-1)}(K)=M_{1,(1)}(K)$. Choose adapted bases, let $Q=Q(B)$ be the resulting doubly-signed quiver, and let $\kappa_\alpha,\mu_\alpha$ be the signs. By (R),
\[ \tau(\tau(\alpha))=\kappa_\alpha\kappa_{\alpha^{\top}}\alpha=g_ig_j\alpha,\qquad \gamma(\gamma(\alpha))=\mu_\alpha\mu_{\bar\alpha}\alpha=h_ih_j\alpha, \]
and $\gamma(\tau(\alpha))=\kappa_\alpha\mu_{\alpha^{\top}}\overline{\alpha^{\top}}$ equals $\tau(\gamma(\alpha))=\mu_\alpha\kappa_{\bar\alpha}\bar\alpha^{\top}$, for $\alpha\in\mathcal B_{ji}$. This gives \eqref{eq:sign-conditions}. Since moreover $\kappa_\alpha=\epsilonform(\alpha)$ when $\alpha^{\top}=\alpha$ and $\mu_\alpha=1$ when $\bar\alpha=\alpha$, the signs are admissible. Let $B_Q$ be built with these signs. Define $\Psi:B_\CC\to S_\CC\oplus M_\CC$ on $S$ as the coordinate isomorphism $S\to S_\CC=\prod_iM_{d_i}(\CC)$ of Step~1 (all vertices of a row of Table~\ref{tab:vertex-blocks} have the same $d_v$, $G_v$ and $C_v$), and on $R$ as the linear map with $\alpha(Y)\mapsto Y\alpha$ for $\alpha\in\mathcal B_{ji}$ and $Y\in M_{ji}$. Then $\Psi$ is an isomorphism of algebras, because it is an isomorphism of bimodules on $R$ and $R^2=M_\CC^2=0$. It commutes with $*$ and $\iota$: on $S$ by Step~1, and on $R$ because, for instance, $\alpha(Y)^*=\kappa_\alpha\,\alpha^{\top}\bigl(\phi(Y)\bigr)$ is mapped to $\kappa_\alpha G_iY^TG_j^{-1}\alpha^{\top}=(Y\alpha)^*$. So $\Psi$ restricts to a $*$-isomorphism from $B$, the algebra of fixed points of $\iota$ in $B_\CC$, onto $B_Q$, the algebra of fixed points of $\iota$ in $S_\CC\oplus M_\CC$.

If $B'$ is another basic real $*$-algebra with radical square zero and $Q(B')\cong Q$, then in the same way $B'\cong B_Q$ for some admissible signs, after transport along the isomorphism. It remains to show that $B_Q$ does not depend on the admissible signs. Let $(\kappa,\mu)$ and $(\kappa',\mu')$ be admissible, and put $a_\alpha=\kappa_\alpha\kappa'_\alpha$ and $b_\alpha=\mu_\alpha\mu'_\alpha$. By \eqref{eq:sign-conditions} and the conditions on fixed arrows,
\[ a_{\alpha^{\top}}=a_\alpha,\qquad b_{\bar\alpha}=b_\alpha,\qquad a_\alpha b_{\alpha^{\top}}=b_\alpha a_{\bar\alpha}, \]
and $a_\alpha=1$ if $\alpha^{\top}=\alpha$, $b_\alpha=1$ if $\bar\alpha=\alpha$. In each orbit of arrows under the two involutions, choose an arrow $\alpha$. Put $z_\alpha=\sqrt{-1}$ if $\bar\alpha=\alpha^{\top}\neq\alpha$ and $a_\alpha b_\alpha=-1$, and $z_\alpha=1$ otherwise; then put $z_{\alpha^{\top}}=a_\alpha z_\alpha$, $z_{\bar\alpha}=b_\alpha\bar z_\alpha$ and $z_{\bar\alpha^{\top}}=a_{\bar\alpha}b_\alpha\bar z_\alpha$. A check of the five kinds of orbits, using the identities above, shows that these values are consistent and that $z_{\beta^{\top}}=a_\beta z_\beta$ and $z_{\bar\beta}=b_\beta\bar z_\beta$ for every arrow $\beta$. Consider the map of $S_\CC\oplus M_\CC$ that fixes $S_\CC$ and multiplies the copy of $M_{d_t\times d_s}(\CC)$ that belongs to each arrow $\beta$ by $z_\beta$. By \eqref{eq:block-star}, it carries the involution and the conjugation defined by $(\kappa,\mu)$ to those defined by $(\kappa',\mu')$. It is an algebra automorphism, because $M_\CC^2=0$, so it restricts to a $*$-isomorphism between the two algebras $B_Q$. Hence two basic real $*$-algebras with radical square zero and isomorphic doubly-signed quivers are $*$-isomorphic.

Every $B$ has a basic $*$-algebra. Let $f\in S_B$ be the sum of the matrix units $E_{11}$ of the factors $M_{r,\eta}(K)$ in Step~1. Then $f$ is self-adjoint and $BfB=B$. Moreover $fBf/fNf\cong fS_Bf$ is a product of the algebras $K$, so $fBf$ is basic by Lemma~\ref{lem:star-simple}(3). By (4), $Q(fBf)\cong Q(B)$. So two real $*$-algebras with radical square zero and isomorphic doubly-signed quivers have $*$-isomorphic basic $*$-algebras.

\emph{Signed quivers.} Let $A$ be a finite-dimensional complex algebra with a $\CC$-linear involution and $\rad(A)^2=0$. The same steps apply without $\iota$. Coordinates as in Step~1 come from (2) at the vertices with $i^{\top}=i$. For $i^{\top}\neq i$, they come from coordinates of $S_{i^{\top}}$ in which $*$ maps $X\in S_i$ to $X^T$. Then only $\tau$, the first identity of (R), and cases (a) and (b) of Step~4 without the conditions on the conjugation are needed, and they yield a signed quiver.
\end{proof}

\section{\texorpdfstring{$*$}{*}-Simple Real \texorpdfstring{$*$}{*}-Algebras and Matrices with Signs}
\label{sec:matrices-with-signs}
A $*$-algebra is \emph{$*$-simple} if it is non-zero and has no $*$-stable two-sided ideals other than $0$ and itself; a semisimple $*$-algebra is the product of its $*$-simple factors. For a doubly-signed quiver $Q$, let $\mathcal O(Q)$ be the set of orbits of the group generated by $(\cdot)^{\top}$ and $\overline{(\cdot)}$ on the vertices of $Q$. If $B$ is a finite-dimensional real $*$-algebra with $\rad(B)^2=0$, then by Theorem~\ref{thm:classification}(2) and its proof the orbits of $Q(B)$ correspond to the $*$-simple factors of $B/\rad(B)$. Each orbit $c$ has one of the ten types of Table~\ref{tab:vertex-types}, and we write $K_c$ for the $*$-algebra of minimal size listed there. We say that $c$, or $K_c$, is of \emph{signature type} if $K_c$ is $\RR$ with the trivial involution, $\CC$ with complex conjugation or $\HH$ with quaternion conjugation, that is, if the unitary groups $\Un_n$ listed in Table~\ref{tab:vertex-types} are compact for every $n$.

For a $*$-algebra $B$ and signs $\eta=(\eta_1,\dots,\eta_n)\in\{\pm1\}^n$, let $M_{n,\eta}(B)$ denote $M_n(B)$ with the involution $X\mapsto\eta X^*\eta$, where $X^*$ is the conjugate transpose and $\eta$ also denotes $\diag(\eta_1,\dots,\eta_n)$. Thus the $(i,j)$ entry of the adjoint of $X$ is $\eta_i\eta_j(x_{ji})^*$. This is the $*$-algebra of $B$-linear endomorphisms of $B^n$ with the adjoint involution of the Hermitian form $h(x,y)=\sum_ix_i\eta_iy_i^*$ (Remark~\ref{rem:indefinite}).

\begin{lemma}\label{lem:star-simple}
Let $S$ be a finite-dimensional real $*$-algebra that is $*$-simple. Then $S\cong M_{n,\eta}(K)$ for some $n\geq1$, some signs $\eta$ and one of the ten $*$-algebras $K$ of Table~\ref{tab:vertex-types}. Here $K$ and $n$ are determined by $S$, and the type of $K$ is the type of the doubly-signed quiver of $S$. Moreover:
\begin{enumerate}
  \item every invertible self-adjoint element of $K$ has the form $\pm uu^*$ with $u\in K^\times$, and the sign can be taken to be $+$ if $K$ is not of signature type; in that case one can take $\eta=(1,\dots,1)$;
  \item if $K$ is of signature type and $M_{n,\eta}(K)\cong M_{n,\eta'}(K)$, then the numbers of entries $+1$ and $-1$ in $\eta$ agree with those in $\eta'$, up to exchanging the two numbers;
  \item every self-adjoint idempotent of $S$ is central if and only if $n=1$.
\end{enumerate}
\end{lemma}

\begin{proof}
The radical of $S$ is a $*$-stable ideal different from $S$ (Lemma~\ref{lem:radical-stable}), so it is zero and $S$ is semisimple. Suppose first that $S$ is not simple, and let $S_1$ be a simple two-sided ideal. Then $S_1\cap S_1^*$ is a $*$-stable ideal different from $S$, hence zero, and $S_1+S_1^*=S$. So $S=S_1\times S_1^*$ with $S_1\cong M_n(D)$ for some $D\in\{\RR,\CC,\HH\}$, and the involution identifies $S_1^*$ with $S_1^{\mathrm{op}}$. Hence $S\cong M_n(D)\times M_n(D)^{\mathrm{op}}$ with the involution that exchanges the factors. This is $M_n(K)$ with the conjugate transpose, where $K=D\times D^{\mathrm{op}}$ carries the exchange involution: the pair $(X,Y^{\mathrm{op}})$ corresponds to the matrix with entries $(x_{ij},y_{ji}^{\mathrm{op}})$. It is one of the three exchange types, with $\eta=(1,\dots,1)$.

Suppose now that $S$ is simple, $S\cong M_m(D)$. Its involution is the adjoint involution of a non-degenerate $\epsilon$-Hermitian form $h$ on $D^m$ with respect to an involution $\theta$ of $D$~\cite[Ch.~I, \S4]{KMRT1998}. We may take $\theta$ to be the identity of $\RR$, the identity or the conjugation of $\CC$, or quaternion conjugation; when $\theta$ is the conjugation of $\CC$, we may take $\epsilon=1$, replacing $h$ by $ih$. If $\theta=\mathrm{id}$ and $\epsilon=-1$, then $h$ is alternating, so $m$ is even and $h$ has a symplectic basis. Grouping that basis into hyperbolic pairs identifies $S$ with $M_{m/2}(K)$, where $K=M_2(D)$ carries the symplectic involution and $\eta=(1,\dots,1)$. In all other cases $h$ has an orthogonal basis (step~4 of the discussion of Theorem~\ref{thm:hermitian-ks}), $h\cong\langle\lambda_1,\dots,\lambda_m\rangle$, and $\langle\lambda\rangle\cong\langle u\lambda\theta(u)\rangle$ for every $u\in D^\times$. This normalises each $\lambda_k$ to $\pm1$ if $(D,\theta,\epsilon)$ is $(\RR,\mathrm{id},1)$, $(\CC,\bar{\ },1)$ or $(\HH,\bar{\ },1)$, and to $1$ if it is $(\CC,\mathrm{id},1)$; then $S\cong M_{m,\eta}(K)$ with $K=(D,\theta)$ and $\eta_k=\lambda_k$. If $(D,\theta,\epsilon)=(\HH,\bar{\ },-1)$, the $\lambda_k$ are non-zero pure quaternions. The elements $u\lambda\bar u$ with $u\in\HH^\times$ are all the non-zero pure quaternions, so the $\lambda_k$ can be normalised to $j$. The adjoint involution is then $X\mapsto j^{-1}X^*j$, and $S\cong M_m(K)$, where $K=\HH$ carries the involution $x\mapsto j^{-1}\bar xj$, which is of orthogonal type, and $\eta=(1,\dots,1)$. In every case $K$ is one of the ten algebras of Table~\ref{tab:vertex-types}. By Theorem~\ref{thm:classification}(4), applied to the corner $E_{11}SE_{11}\cong K$, the doubly-signed quivers of $S$ and $K$ agree. The ten algebras have different doubly-signed quivers, so $K$ is determined by $S$, and then so is $n$, because $\dim S=n^2\dim K$.

(1) Let $\lambda\in K$ be invertible and self-adjoint. If $K$ is of signature type, its self-adjoint elements are real numbers, and $\lambda=\pm uu^*$ with $u=|\lambda|^{1/2}$. Otherwise we find $u$ with $\lambda=uu^*$ as follows. For $\CC$ with the trivial involution, take $u=\lambda^{1/2}$. For $\HH$ with the involution $x\mapsto j^{-1}\bar xj$, whose self-adjoint elements are the combinations of $1$, $i$ and $k$, the element $\lambda j$ is a non-zero pure quaternion; so $\lambda j=uj\bar u$ for some $u$, and then $\lambda=uj\bar uj^{-1}=uu^*$. For $M_2(D)$ with the symplectic involution, $u^*$ is the adjugate of $u$, so $u^*=\tr(u)I-u$; the self-adjoint elements are therefore the scalars $\mu$, and $\mu=uu^*$ for $u=\diag(\mu,1)$. For $D\times D^{\mathrm{op}}$, the self-adjoint elements are the pairs $(d,d^{\mathrm{op}})$, and $(d,d^{\mathrm{op}})=uu^*$ for $u=(d,1)$. In particular $-1=uu^*$ for some $u$ when $K$ is not of signature type. Replacing the $k$-th standard basis vector $b_k$ of $K^n$ by $ub_k$ then changes $h(b_k,b_k)=\eta_k$ into $u\eta_ku^*=-\eta_k$. So all signs can be made $+1$ by an isometry, which induces a $*$-isomorphism of the corresponding algebras $M_{n,\eta}(K)$.

(2) Let $\psi:M_{n,\eta}(K)\to M_{n,\eta'}(K)$ be a $*$-isomorphism. If $K=\CC$ and $\psi$ is not $\CC$-linear, we compose it with entrywise complex conjugation, which is a $*$-automorphism of $M_{n,\eta'}(\CC)$ because $\eta'$ is real. Then $\psi$ is linear over the centre ($\RR$ or $\CC$) of the central simple algebra $M_n(K)$, so $\psi(X)=gXg^{-1}$ for an invertible $g$ by the Skolem--Noether theorem (Theorem~\ref{thm:skolem-noether}). Comparing $\psi(\eta X^*\eta)$ with $\eta'\psi(X)^*\eta'$ shows that the matrix $\eta g^*\eta'g$ commutes with every matrix, so $g^*\eta'g=\lambda\eta$ for a central scalar $\lambda$, which is real because both sides are Hermitian. Thus the Hermitian matrices $\eta'$ and $\lambda\eta$ are congruent, and the claim follows from Sylvester's law of inertia and its Hermitian and quaternionic analogues (Theorem~\ref{thm:hermitian-ks}).

(3) If $n\geq2$, then $E_{11}$ is a self-adjoint idempotent that is not central. If $n=1$, then $S\cong K$, and the only self-adjoint idempotents of $K$ are $0$ and $1$. Indeed, $\RR$, $\CC$ and $\HH$ have no other idempotents; the self-adjoint elements of $M_2(D)$ with the symplectic involution are the scalars, as shown in (1); and those of $D\times D^{\mathrm{op}}$ are the pairs $(d,d^{\mathrm{op}})$, which are idempotent only for $d\in\{0,1\}$.
\end{proof}

By part~(3), a real $*$-algebra $A$ is basic as a $*$-algebra (Definition~\ref{def:basic-star}) exactly when every $*$-simple factor of $A/\rad(A)$ is one of the ten algebras of Table~\ref{tab:vertex-types}.

\section{Path \texorpdfstring{$*$}{*}-Algebras}
\label{sec:path-star-algebras}
We now construct the path $*$-algebra $\RR Q$ of a doubly-signed quiver $Q$. Its complexification is built explicitly, as the path algebra of $Q$ with small matrices as coefficients, and $\RR Q$ is the set of fixed points of a conjugation. This is the form in which Chapter~\ref{division-problem} localises and completes $\RR Q$, and in which the implementation of Chapter~\ref{chap:implementation} stores it.

Let $Q$ be a doubly-signed quiver; we work with its complexified vertices and arrows, on which the transpose and the conjugation act. Give each vertex $v$ the block size $d_v$ and the matrices $G_v,C_v\in\GL_{d_v}(\CC)$ of Table~\ref{tab:vertex-blocks}, where $J=\bigl(\begin{smallmatrix}0&1\\-1&0\end{smallmatrix}\bigr)$. Thus $d_v=2$ if $\epsilondiv(v)=-1$ or $\epsilonform(v)=-1$, and $d_v=1$ otherwise; $G_v=J$ if $\epsilonform(v)=-1$, $C_v=J$ if $\epsilondiv(v)=-1$, and otherwise they are identity matrices. Put $g_v=\epsilonform(v)$ if $v^{\top}=v$ and $g_v=1$ otherwise, and $h_v=\epsilondiv(v)$ if $\bar v=v$ and $h_v=1$ otherwise. By rules (2) and (3) of Definition~\ref{def:doubly-signed}, $g$ and $h$ are constant on the orbits of the transpose and the conjugation. Moreover $G_{v^{\top}}=G_v$ and $C_{\bar v}=C_v$, and the table shows that
\begin{equation}\label{eq:vertex-data}
G_v^T=g_vG_v,\qquad C_v\overline{C_v}=h_vI,\qquad C_{v^{\top}}\overline{G_v}\,C_v^T=G_{\bar v}.
\end{equation}

\begin{table}[htbp]
\centering
\caption{The block model of the ten vertex types of Table~\ref{tab:vertex-types}: the names used in Chapters~\ref{division-problem} and~\ref{chap:implementation}, the complex vertices of an orbit, and the block size, vertex form and real structure at each of them.}
\label{tab:vertex-blocks}
\small
\renewcommand{\arraystretch}{1.2}
\begin{tabular}{@{}lllccc@{}}
\toprule
\textbf{Name} & \textbf{Vertex data} & \textbf{Complex vertices} & $d_v$ & $G_v$ & $C_v$ \\
\midrule
\texttt{R}   & $\epsilondiv=+1$, $\epsilonform=+1$ & $b$ & 1 & $I$ & $I$ \\
\texttt{Rsp} & $\epsilondiv=+1$, $\epsilonform=-1$ & $b$ & 2 & $J$ & $I$ \\
\texttt{Hsp} & $\epsilondiv=-1$, $\epsilonform=-1$ & $b$ & 2 & $J$ & $J$ \\
\texttt{Hor} & $\epsilondiv=-1$, $\epsilonform=+1$ & $b$ & 2 & $I$ & $J$ \\
\texttt{Cor} & $\epsilondiv=0$, $\epsilonform=+1$ & $b,\bar b$ & 1 & $I$ & $I$ \\
\texttt{Csp} & $\epsilondiv=0$, $\epsilonform=-1$ & $b,\bar b$ & 2 & $J$ & $I$ \\
\texttt{Cu}  & $\epsilondiv=0$, $b^{\top}=\bar b$ & $b,\bar b$ & 1 & $I$ & $I$ \\
\texttt{RR}  & $\epsilondiv=+1$, $b^{\top}\neq b$ & $b,b^{\top}$ & 1 & $I$ & $I$ \\
\texttt{HH}  & $\epsilondiv=-1$, $b^{\top}\neq b$ & $b,b^{\top}$ & 2 & $I$ & $J$ \\
\texttt{CC}  & $\epsilondiv=0$, $b^{\top}\notin\{b,\bar b\}$ & $b,\bar b,b^{\top},\bar b^{\top}$ & 1 & $I$ & $I$ \\
\bottomrule
\end{tabular}
\end{table}

Let $T$ be the complex vector space of finite sums $\sum_wX_w\,w$ over the paths $w$ of $Q$, with $X_w\in M_{d_t\times d_s}(\CC)$ for $w:s\to t$, and with the product $(X\,w)(Y\,w')=(XY)\,(ww')$, which is zero when $ww'$ is not defined. Since $M_{d_t\times d_m}(\CC)\otimes_{M_{d_m}(\CC)}M_{d_m\times d_s}(\CC)\cong M_{d_t\times d_s}(\CC)$, $T$ is the tensor algebra $T_{S_\CC}(M_\CC)$ of the bimodule $M_\CC=\bigoplus_{\alpha:s\to t}M_{d_t\times d_s}(\CC)$ over $S_\CC=\prod_vM_{d_v}(\CC)$, which contains one copy of the simple bimodule $M_{d_t\times d_s}(\CC)$ for each arrow $\alpha:s\to t$. Given signs $\kappa_\alpha,\mu_\alpha\in\{\pm1\}$ for the arrows, extended multiplicatively to paths, define, for a path $w:s\to t$,
\begin{equation}\label{eq:block-star}
(X\,w)^*=\kappa(w)\,G_sX^TG_t^{-1}\,w^{\top},\qquad \iota(X\,w)=\mu(w)\,C_t\overline XC_s^{-1}\,\bar w ,
\end{equation}
where $w^{\top}$ is the path of the transposed letters of $w$ in reverse order and $\bar w$ is the path of its conjugated letters. We call the signs \emph{admissible} if, for every arrow $\alpha:s\to t$,
\begin{equation}\label{eq:sign-conditions}
\kappa_\alpha\kappa_{\alpha^{\top}}=g_sg_t,\qquad \mu_\alpha\mu_{\bar\alpha}=h_sh_t,\qquad \kappa_\alpha\mu_{\alpha^{\top}}=\mu_\alpha\kappa_{\bar\alpha},
\end{equation}
and moreover $\kappa_\alpha=\epsilonform(\alpha)$ whenever $\alpha^{\top}=\alpha$, and $\mu_\alpha=1$ whenever $\bar\alpha=\alpha$.

\begin{definition}\label{def:path-star-algebra}
Let $Q$ be a doubly-signed quiver, and fix admissible signs; they exist by Proposition~\ref{prop:block-model}(2). The \emph{path $*$-algebra} of $Q$ is the real algebra
\[ \RR Q:=T^\iota=\{x\in T:\iota(x)=x\}, \]
with the involution $*$ of \eqref{eq:block-star}. It is graded by the length of paths. We write $S_Q$ and $M_Q$ for its parts of degree $0$ and $1$, $\RR Q_{\geq\ell}$ for the ideal of elements of degree at least $\ell$, and
\[ B_Q:=\RR Q/\RR Q_{\geq2}=S_Q\oplus M_Q . \]
$\RR Q_{\geq1}$ is the \emph{arrow ideal}. An ideal $I$ of $\RR Q$ is \emph{admissible} if $I^*=I$ and $\RR Q_{\geq m}\subseteq I\subseteq\RR Q_{\geq2}$ for some $m\geq2$.
\end{definition}

\begin{proposition}[The block model]\label{prop:block-model}
Let $Q$ be a doubly-signed quiver.
\begin{enumerate}
  \item The formulas \eqref{eq:block-star} define a $\CC$-linear anti-automorphism $*$ and a conjugate-linear automorphism $\iota$ of $T$. They satisfy $x^{**}=x$, $\iota(\iota(x))=x$ and $\iota(x^*)=\iota(x)^*$ for all $x\in T$ if and only if the signs satisfy \eqref{eq:sign-conditions} for every arrow.
  \item There are admissible signs.
  \item For admissible signs, $\RR Q$ is a real $*$-algebra with $T=\RR Q\oplus i\RR Q$, so that $\CC\otimes_\RR\RR Q\cong T$. It is the tensor algebra $T_{S_Q}(M_Q)$, with the involution that reverses tensors, and $S_Q\cong\prod_{c\in\mathcal O(Q)}K_c$, with the algebras $K_c$ of Table~\ref{tab:vertex-types}. The algebra $B_Q$ is a basic real $*$-algebra with radical $M_Q$, $M_Q^2=0$ and doubly-signed quiver $Q$.
\end{enumerate}
\end{proposition}
\begin{proof}
(1) For paths $w:s\to t$ and $w':r\to s$, the coefficient of $(ww')^{\top}=w'^{\top}w^{\top}$ in $\bigl((X\,w)(Y\,w')\bigr)^*$ is
\[ \kappa(w)\kappa(w')\,G_rY^TX^TG_t^{-1}=\bigl(\kappa(w')\,G_rY^TG_s^{-1}\bigr)\bigl(\kappa(w)\,G_sX^TG_t^{-1}\bigr), \]
which is the coefficient of $(Y\,w')^*(X\,w)^*$. So $*$ reverses products; in the same way $\iota$ preserves them, because $C_s^{-1}C_s=I$. Using \eqref{eq:vertex-data}, $G_{v^{\top}}=G_v$ and $C_{\bar v}=C_v$, one computes
\[ (X\,w)^{**}=\kappa(w)\kappa(w^{\top})\,g_sg_t\;X\,w,\qquad \iota(\iota(X\,w))=\mu(w)\mu(\bar w)\,h_sh_t\;X\,w , \]
and finds that $\iota\bigl((X\,w)^*\bigr)$ and $\iota(X\,w)^*$ are $\kappa(w)\mu(w^{\top})$ and $\mu(w)\kappa(\bar w)$ times the same element, because the third identity of \eqref{eq:vertex-data} gives $C_{s^{\top}}\overline{G_s}=G_{\bar s}C_s^{-T}$ and $\overline{G_t}^{\,-1}C_{t^{\top}}^{-1}=C_t^TG_{\bar t}^{-1}$. Each of the three scalar factors is the product of the corresponding factors of the letters of $w$, because the inner vertices of $w$ contribute $g_v^2=h_v^2=1$. Hence the three identities hold on $T$ if and only if they hold on the arrows, which is \eqref{eq:sign-conditions}.

(2) For each orbit of arrows under the group generated by the transpose and the conjugation, choose a representative $\alpha:s\to t$, put $\kappa_\alpha=\epsilonform(\alpha)$ if $\alpha^{\top}=\alpha$ and $\kappa_\alpha=1$ otherwise, and define, for $\beta\in\{\alpha,\bar\alpha\}$ and $\beta'\in\{\alpha,\alpha^{\top}\}$,
\[ \kappa_\beta=\kappa_\alpha,\qquad \kappa_{\beta^{\top}}=g_sg_t\,\kappa_\alpha,\qquad \mu_{\beta'}=1,\qquad \mu_{\bar\beta'}=h_sh_t . \]
These prescriptions agree where the orbit has fewer than four elements. If $\alpha^{\top}=\alpha$, then $t=s^{\top}$, so $g_sg_t=g_s^2=1$. If $\bar\alpha=\alpha$, then $\bar s=s$ and $\bar t=t$, and rule~(4) of Definition~\ref{def:doubly-signed} gives $\epsilondiv(s)=\epsilondiv(t)$, so $h_sh_t=1$. If $\bar\alpha=\alpha^{\top}$, then $t=\bar s^{\top}$ lies in the orbit of $s$, so $g_sg_t=h_sh_t=1$. Every arrow $\beta$ of the orbit has $g_{s(\beta)}g_{t(\beta)}=g_sg_t$ and $h_{s(\beta)}h_{t(\beta)}=h_sh_t$, and by construction $\kappa_{\beta^{\top}}=g_sg_t\kappa_\beta$, $\kappa_{\bar\beta}=\kappa_\beta$, $\mu_{\bar\beta}=h_sh_t\mu_\beta$ and $\mu_{\beta^{\top}}=\mu_\beta$. This gives \eqref{eq:sign-conditions}. Finally, if $\beta^{\top}=\beta$ for some $\beta$ in the orbit, then $\alpha^{\top}=\alpha$, the orbit is contained in $\{\alpha,\bar\alpha\}$, and $\kappa_\beta=\epsilonform(\alpha)=\epsilonform(\beta)$ by Definition~\ref{def:doubly-signed}(3); and if $\bar\beta=\beta$ for some $\beta$, then $\bar\alpha=\alpha$ and every $\mu_\beta$ equals $1$.

(3) Since $\iota$ is a conjugate-linear automorphism of order two, $T=T^\iota\oplus iT^\iota$, and $T^\iota$ is a real subalgebra, stable under $*$ because $*$ commutes with $\iota$. Since $\iota$ is multiplicative and preserves the path degree, the multiplication maps identify $T^\iota$ with the tensor algebra of its part $M_Q=M_\CC^c$ of degree $1$ over its part $S_Q=S_\CC^c$ of degree $0$, because they do so after extending scalars to $\CC$. The quotient $B_Q=S_Q\oplus M_Q$ is a real $*$-algebra with radical $M_Q$ and $M_Q^2=0$. Its complexification is $S_\CC\oplus M_\CC$ with the maps \eqref{eq:block-star}, and comparison with Theorem~\ref{thm:classification} shows that its doubly-signed quiver is $Q$: the involution acts on a factor $M_{d_v}(\CC)$ with $v^{\top}=v$ by $X\mapsto G_vX^TG_v^{-1}$, where $G_v^T=\epsilonform(v)G_v$; the conjugation acts on a factor with $\bar v=v$ by $X\mapsto C_v\overline XC_v^{-1}$, where $C_v\overline{C_v}=\epsilondiv(v)I$; on the copy of $M_{d_t\times d_s}(\CC)$ that belongs to an arrow $\alpha$ with $\alpha^{\top}=\alpha$, the involution is $\epsilonform(\alpha)$ times the map $X\mapsto G_sX^TG_t^{-1}$, which is the reference map of Theorem~\ref{thm:classification}(3); and the conjugation maps the copy that belongs to $\alpha$ onto the copy that belongs to $\bar\alpha$, so that $\iota\circ*$ maps the copy that belongs to an arrow with $\bar\alpha=\alpha^{\top}$ onto itself. Finally, the part of $S_Q$ at an orbit $c$ of vertices is the algebra of fixed points of $\iota$ in $\prod_{v\in c}M_{d_v}(\CC)$, with the involution of \eqref{eq:block-star}, and inspection of the ten rows of Table~\ref{tab:vertex-blocks} identifies it with the algebra $K_c$ of Table~\ref{tab:vertex-types}. For one vertex $b$ with $C_b=I$ the fixed points form $M_{d_b}(\RR)$, with the transpose (type \texttt{R}) or the symplectic involution $X\mapsto JX^TJ^{-1}$ (\texttt{Rsp}). With $C_b=J$ they form the quaternions, as computed after this proof, with quaternion conjugation (\texttt{Hsp}) or with the transpose, which is $x\mapsto j^{-1}\bar xj$ and whose fixed points have dimension $3$, so that it is of orthogonal type (\texttt{Hor}). For two vertices $b,\bar b$ exchanged by $\iota$, the fixed points are the pairs $(X,\overline X)$, a copy of $M_{d_b}(\CC)$, with the transpose (\texttt{Cor}), the symplectic involution (\texttt{Csp}) or, when $b^{\top}=\bar b$, the conjugate transpose (\texttt{Cu}). For two vertices $b,b^{\top}$ fixed by $\iota$ (\texttt{RR}, \texttt{HH}), and for the four vertices of type \texttt{CC}, the fixed points form $D\times D$ with $D=\RR$, $\HH$ or $\CC$, and the involution exchanges the two factors. So $S_Q\cong\prod_cK_c$, and $B_Q$ is basic as a $*$-algebra (Lemma~\ref{lem:star-simple}(3)).
\end{proof}

By part~(3), $B_Q$ realises the converse statement of Theorem~\ref{thm:classification}: every doubly-signed quiver is the doubly-signed quiver of a basic real $*$-algebra with radical square zero. By the uniqueness statement of that theorem, $B_Q$ is determined by $Q$ up to $*$-isomorphism, whatever admissible signs are chosen. Hence so is $\RR Q=T_{S_Q}(M_Q)$: since $M_Q^2=0$, a $*$-isomorphism between two such algebras $B_Q$ induces compatible isomorphisms of their parts $S_Q$ and $M_Q$, and hence of the tensor algebras. The results of Chapters~\ref{chap:indexing} to~\ref{chap:implementation} about $\RR Q$ and its quotients do not use this uniqueness; they hold for any fixed choice of admissible signs.

When every $K_c$ is a division algebra, $\RR Q$ without its involution is the tensor algebra of a species in the sense of Dlab and Ringel~\cite{DlabRingel1976}. In general it is Morita equivalent to one, because $T_{S}(M)$ is Morita equivalent to $T_{eSe}(eMe)$ for every idempotent $e\in S$ with $SeS=S$. If every vertex of $Q$ has type $(\RR,+1)$, then all blocks have size $1$, $S_Q=\RR^{Q_0}$, the arrows form a basis of $M_Q$, and $\RR Q$ is the path algebra of Definition~\ref{def:path-algebra}; for the signs of Proposition~\ref{prop:block-model}(2), its involution is that of Section~\ref{sec:symmetric-quivers}, multiplied by $\epsilonform(\alpha)$ on each arrow $\alpha$ fixed by the transpose. In general, if $e\in S_\CC$ is the sum of one diagonal matrix unit from each factor $M_{d_v}(\CC)$, then $eTe$ is the complex path algebra of the quiver underlying $Q$. So $T$ is Morita equivalent to that path algebra, which justifies the name.

For the real singular value decomposition (Figure~\ref{fig:dsq-examples}(a)), all blocks have size $1$, $G=C=I$, and the signs are $1$, so that $(X\,a)^*=X\,a^{\top}$; this is the algebra of Section~\ref{sec:unification-example}, with $b=a^{\top}$. For the quaternionic singular value decomposition (Figure~\ref{fig:dsq-examples}(b)), $d_u=d_v=2$ and all the matrices $G$ and $C$ equal $J$. The fixed points of $X\mapsto J\overline XJ^{-1}$ in $M_2(\CC)$ are the blocks $\bigl(\begin{smallmatrix}p&q\\-\bar q&\bar p\end{smallmatrix}\bigr)$, which represent the quaternions $p+qj$, and $X\mapsto JX^TJ^{-1}$ maps such a block to the block of the conjugate quaternion. So the involution of the block model multiplies by conjugate quaternions, as in Section~\ref{subsec:cat-quaternion-svd}. In the arithmetic of $T$, the matrices $G_v$ and $C_v$ are the only places where quaternions and symplectic forms enter.

\section{The Structure Theorem: An Analogue of Gabriel's Theorem}
\label{sec:structure}
Theorem~\ref{thm:classification} describes the real $*$-algebras whose radical squares to zero. We now show that every finite-dimensional real $*$-algebra is a matrix algebra, with sizes and signs at the orbits of vertices, over a quotient $\RR Q/I$ of the path $*$-algebra of its doubly-signed quiver (Theorem~\ref{thm:structure}). We first describe these matrix algebras. To assemble matrices over a basic $*$-algebra whose vertices have different sizes, we use the following tensor product.

\begin{definition}\label{def:matrix-like}
Let $\mathcal O$ be a finite set, and let $\RR^{\mathcal O}$ carry the trivial involution, with minimal idempotents $1_c$ for $c\in\mathcal O$. An \emph{$\RR^{\mathcal O}$-$*$-algebra} is a $*$-algebra $B$ with a unital $*$-homomorphism $\RR^{\mathcal O}\to B$, that is, with pairwise orthogonal self-adjoint idempotents $1_c\in B$ such that $\sum_c1_c=1$. We write $B_{cd}=1_cB1_d$.
\begin{enumerate}
  \item The \emph{tensor product over $\RR^{\mathcal O}$} of two $\RR^{\mathcal O}$-$*$-algebras $E$ and $B$ is
  \[ E\otimes_{\RR^{\mathcal O}}B:=\bigoplus_{c,d\in\mathcal O}E_{cd}\otimes_\RR B_{cd}=e\,(E\otimes_\RR B)\,e,\qquad e=\sum_{c\in\mathcal O}1_c\otimes1_c . \]
  Since $e$ is a self-adjoint idempotent of the $*$-algebra $E\otimes_\RR B$ (Definition~\ref{def:tensor}), this is a $*$-algebra with unit $e$. When the $1_c$ are central in $E$ and in $B$, it is the usual tensor product over $\RR^{\mathcal O}$.
  \item For each $c\in\mathcal O$, let $V_c$ be a finite-dimensional real vector space with a non-degenerate symmetric bilinear form $\eta_c$, and let $(V,\eta)$ be their orthogonal sum. We write $\End(V,\eta)$ for $\End_\RR(V)$ with the adjoint involution of $\eta$. It is an $\RR^{\mathcal O}$-$*$-algebra, $1_c$ being the orthogonal projection onto $V_c$.
\end{enumerate}
\end{definition}

Choose bases $(v_i)$ of the spaces $V_c$ that are orthogonal for $\eta$, with $\eta(v_i,v_i)=\eta_i\in\{\pm1\}$, and write $c(i)$ for the index of the space containing $v_i$. The corresponding matrix units $e_{ij}\in\End(V)$ satisfy $e_{ij}^*=\eta_i\eta_je_{ji}$. Therefore $\sum_{i,j}e_{ij}\otimes x_{ij}\mapsto(x_{ij})$ identifies $\End(V,\eta)\otimes_{\RR^{\mathcal O}}B$ with the $*$-algebra of generalised matrices $X=(x_{ij})$ with entries $x_{ij}\in B_{c(i)c(j)}$, multiplied as matrices, with the involution
\begin{equation}\label{eq:signed-adjoint}
  X\longmapsto\eta X^*\eta,\qquad\text{that is,}\qquad (x_{ij})\longmapsto\bigl(\eta_i\eta_j\,(x_{ji})^*\bigr).
\end{equation}
These are matrices over $B$ indexed by idempotents, as in Section~\ref{sec:indexing}, with $n_c=\dim V_c$ rows and columns for each $c$, and with signs. Special cases are $M_n(B)$ with the conjugate transpose (all $n_c=n$ and all $\eta_i=1$), the corner algebras $eM_N(B)e$ at diagonal idempotents $e$ (all $\eta_i=1$), and $M_{n,\eta}(K)$ (when $\mathcal O$ is a single point). In the language of Chapter~\ref{chap:staralg}, $\End(V,\eta)\otimes_{\RR^{\mathcal O}}B$ is the $*$-algebra of $B$-linear endomorphisms, with the adjoint involution, of the Hermitian $B$-module $\bigoplus_iB1_{c(i)}$ with the form $\sum_ix_i\eta_iy_i^*$. The theorem below is therefore a statement of Hermitian Morita theory~\cite[Ch.~I, \S9]{Knus1991}.

Through $\RR^{\mathcal O(Q)}\to S_Q$, the quotients $\RR Q/I$ of the path $*$-algebra of a doubly-signed quiver $Q$ are $\RR^{\mathcal O(Q)}$-$*$-algebras.

Throughout the rest of this section, $A$ is a finite-dimensional real $*$-algebra with radical $N$. Since $N^2$ is stable under the involution, $A/N^2$ is a $*$-algebra whose radical squares to zero, and we call $Q(A):=Q(A/N^2)$ the \emph{doubly-signed quiver of $A$}. Since $(A/N^2)/(N/N^2)=A/N$, the orbits of $Q(A)$ correspond to the $*$-simple factors of $A/N$.

\begin{theorem}[Structure theorem; an analogue of Gabriel's theorem]\label{thm:structure}
Let $A$ be a finite-dimensional real $*$-algebra, $Q=Q(A)$ its doubly-signed quiver, and $\mathcal O=\mathcal O(Q)$.
\begin{enumerate}
  \item There are non-zero finite-dimensional real spaces $V_c$ with non-degenerate symmetric bilinear forms $\eta_c$ ($c\in\mathcal O$), and an admissible ideal $I$ of $\RR Q$, such that
  \[ A\;\cong\;\End(V,\eta)\otimes_{\RR^{\mathcal O}}\RR Q/I \]
  as $*$-algebras. Here $\RR Q/I$ is basic as a $*$-algebra, and $\RR Q_{\geq m}\subseteq I$ for every $m$ with $\rad(A)^m=0$. Moreover, $\eta_c$ can be taken positive definite if $c$ is not of signature type, and of signature $(p_c,q_c)$ with $p_c\geq q_c$ if it is; we then call the decomposition \emph{normalised}.
  \item Let $A\cong\End(V',\eta')\otimes_{\RR^{\mathcal O'}}\RR Q'/I'$ be another normalised decomposition, for a doubly-signed quiver $Q'$ with $\mathcal O'=\mathcal O(Q')$, non-zero spaces $V'_c$ and an admissible ideal $I'$. Then there is an isomorphism $Q'\cong Q$ under which $\dim V'_c=\dim V_c$ and the signatures of $\eta'_c$ and $\eta_c$ agree. Thus $Q$, the dimensions $n_c=\dim V_c$ and the signatures $(p_c,q_c)$ are invariants of $A$.
  \item In the situation of (2), let $B=\RR Q/I$ and $B'=\RR Q'/I'$. Then $B'$ is $*$-isomorphic to $B^t$, the algebra $B$ with the involution $x\mapsto x^{*_t}:=tx^*t$, for some $t=\sum_c\tau_c1_c$ with $\tau_c\in\{\pm1\}$ and $\tau_c=1$ unless $c$ is of signature type and $p_c=q_c$. Conversely, $\End(V,\eta)\otimes_{\RR^{\mathcal O}}B^t\cong\End(V,\eta)\otimes_{\RR^{\mathcal O}}B$ for every such $t$.
\end{enumerate}
\end{theorem}
Since $t$ is self-adjoint with $t^2=1$, the map $x\mapsto tx^*t$ in (3) is again an involution. The condition $V'_c\neq0$ in (2) cannot be dropped: $\RR\cong\End(V,\eta)\otimes_{\RR^2}(\RR\times\RR)$ with $V_1=\RR$ and $V_2=0$.

\begin{remark}\label{rem:structure}
By Lemma~\ref{lem:star-simple}(3), $A$ is basic as a $*$-algebra exactly when all $n_c=1$, and then Theorem~\ref{thm:structure} says that $A\cong\RR Q/I$. This is the analogue for real $*$-algebras of Gabriel's theorem~\cite[Theorem~II.3.7]{AssemSimsonSkowronski2006}. In general, the theorem separates three kinds of data. The quiver $Q$, the sizes $n_c$ and the signatures $(p_c,q_c)$ are discrete invariants of $A$. The twist $t$ in (3) is a finite ambiguity; it occurs only at orbits of split signature $p_c=q_c$, and it is genuine (Example~\ref{ex:structure-twist}). All continuous parameters lie in the ideal $I$, which is determined by $\RR Q/I$ up to $*$-automorphisms of $\RR Q/\RR Q_{\geq m}$ (Proposition~\ref{prop:ideal-uniqueness}). When $\rad(A)^2=0$, we have $I=\RR Q_{\geq2}$, and $B^t\cong B$ for every $t$, via the map that fixes $S_Q$ and sends $r\in M_Q$ to $tr$. In that case $A\cong\End(V,\eta)\otimes_{\RR^{\mathcal O}}B_Q$ is determined by $Q$, the sizes and the signatures, which refines Theorem~\ref{thm:classification}.
\end{remark}

The structure theorem also accounts for the indexing of Section~\ref{sec:indexing}.

\begin{corollary}[Positions]\label{cor:positions}
In the situation of Theorem~\ref{thm:structure}(1), put $B=\RR Q/I$ and $n_c=\dim V_c$, and call the pairs $c_k$ with $c\in\mathcal O$ and $1\leq k\leq n_c$ the \emph{positions}.
\begin{enumerate}
  \item $A$ is $*$-isomorphic to the algebra of matrices $X=(x_{c_k,d_l})$ indexed by the positions, with entries $x_{c_k,d_l}\in1_cB1_d$, multiplied as matrices, with the involution $X\mapsto(\eta_{c_k}\eta_{d_l}\,x_{d_l,c_k}^*)$ of \eqref{eq:signed-adjoint}.
  \item For $n\geq1$, $M_n(A)$ with the conjugate transpose is $*$-isomorphic to $\End(V^n,\eta^n)\otimes_{\RR^{\mathcal O}}B$, where $(V^n,\eta^n)$ is the orthogonal sum of $n$ copies of $(V,\eta)$. Thus $M_n(A)$ has the doubly-signed quiver $Q$ and the relations $I$, with $n\,n_c$ positions and the signature $(np_c,nq_c)$ at the orbit $c$.
\end{enumerate}
\end{corollary}
\begin{proof}
Part (1) is the identification \eqref{eq:signed-adjoint}, applied to Theorem~\ref{thm:structure}(1). For part (2), recall from Section~\ref{subsec:tensor} that $M_n(A)\cong M_n(\RR)\otimes_\RR A$, with the transpose on $M_n(\RR)$. Moreover $M_n(\RR)\otimes_\RR\End(V,\eta)=\End(\RR^n\otimes V,I_n\otimes\eta)$, and both $M_n(\RR)\otimes_\RR\bigl(\End(V,\eta)\otimes_{\RR^{\mathcal O}}B\bigr)$ and $\End(\RR^n\otimes V,I_n\otimes\eta)\otimes_{\RR^{\mathcal O}}B$ are the corner of $M_n(\RR)\otimes_\RR\End(V)\otimes_\RR B$ at the idempotent $\sum_cI_n\otimes1_c\otimes1_c$. The form $I_n\otimes\eta$ is the orthogonal sum of $n$ copies of $\eta$. On $\RR^n\otimes V_c$ it is positive definite where $\eta_c$ is, and it has signature $(np_c,nq_c)$, so the decomposition is normalised, and Theorem~\ref{thm:structure}(2) gives the invariants.
\end{proof}

So the positions of Section~\ref{sec:indexing} are intrinsic to $A$. Their orbits are the orbits of the vertices of $Q$ under the transpose and the conjugation, the number of positions in each orbit is an invariant, and the entries between the positions of the orbits $c$ and $d$ lie in the corner $1_c(\RR Q/I)1_d$. Passing from $A$ to matrices over $A$ only multiplies the numbers of positions. The order in which the positions are listed is not intrinsic. It is chosen by the algorithm, and it matters (Proposition~\ref{prop:golub-kahan}).

\subsection{Examples}
\label{subsec:structure-examples}
The following examples illustrate the three parts of Theorem~\ref{thm:structure}; its proof follows them, in Section~\ref{subsec:structure-proof}.

\begin{example}[Truncated polynomial rings]\label{ex:structure-loop}
Let $Q$ have one vertex of type $(\RR,+1)$ and one loop $x$ fixed by the transpose, with sign $\pm1$. Then $\RR Q=\RR[x]$ with $x^*=\pm x$, and the admissible ideals are the ideals $(x^m)$ with $m\geq2$. By Theorem~\ref{thm:structure}, the real $*$-algebras with doubly-signed quiver $Q$ are the algebras $M_{n,\eta}(\RR[x]/(x^m))$ with $\eta=(1,\dots,1,-1,\dots,-1)$, containing $p\geq q$ entries $+1$ and $-1$ respectively, and these are pairwise non-isomorphic. For $m=2$ these are the algebras of Sections~\ref{subsec:cat-symmetric} and~\ref{subsec:cat-skew}, with an indefinite form when $q>0$. For $m\geq3$ they are matrices over truncated power series, as in Chapter~\ref{chap:infinitesimals}; for instance, $\RR[\varepsilon]/(\varepsilon^3)$ with $\varepsilon^*=\varepsilon$ is $\RR Q/(x^3)$ for the sign $+1$.
\end{example}

\begin{example}[Degenerating diagonal forms]\label{ex:structure-degenerate}
Let $G=\diag(s_1\varepsilon^{w_1},\dots,s_d\varepsilon^{w_d})$ with signs $s_i=\pm1$ and integers $w_i\geq0$, and consider
\[ \Lambda_G=\{X\in M_d(\RR[\varepsilon]):G^{-1}X^TG\in M_d(\RR[\varepsilon])\}, \]
where $G^{-1}X^TG$ is computed in $M_d(\RR[\varepsilon,\varepsilon^{-1}])$, with the involution $X\mapsto G^{-1}X^TG$, the adjoint involution of the form $G$. The reduction $A_G=\Lambda_G/\varepsilon\Lambda_G$ is a real $*$-algebra of dimension $d^2$. It describes the adjoint involution ``at $\varepsilon=0$'', where $G$ itself degenerates. Let $\omega_1<\dots<\omega_r$ be the distinct values of the $w_i$. Then
\[ A_G\cong\End(V,\eta)\otimes_{\RR^r}\RR Q/I . \]
Here $Q$ is the double of the quiver $1\to2\to\cdots\to r$: all vertices have type $(\RR,+1)$, and there are arrows $a_j:j\to j+1$ and $a_j^{\top}:j+1\to j$. The ideal $I$ is generated by the $2$-cycles $a_j^{\top}a_j$ and $a_ja_j^{\top}$, and $(V_j,\eta_j)$ is the form $\langle s_i:w_i=\omega_j\rangle$. Indeed, $\Lambda_G$ consists of the matrices whose $(k,l)$ entry is divisible by $\varepsilon^{\max(0,w_l-w_k)}$. The elements $\hat e_{kl}=\varepsilon^{\max(0,w_l-w_k)}E_{kl}$ form a basis of $A_G$, with $\hat e_{kl}^*=s_ks_l\hat e_{lk}$, and $\hat e_{kl}\hat e_{lm}=\hat e_{km}$ if $w_l$ lies between $w_k$ and $w_m$, and $0$ otherwise. On the other side, $\RR Q/I$ has a basis of monotone paths, which multiply by the same rule. So only the order of the weights and the signs matter. Because $I$ is spanned by paths, rescaling the arrows $a_j^{\top}$ by signs shows that $B^t\cong B$ for every $t$, so the signs can be normalised as in Theorem~\ref{thm:structure}(1) without changing $\RR Q/I$. For $r=2$, $\RR Q/I$ is the four-dimensional algebra of the singular value decomposition (Section~\ref{sec:unification-example}). For $r\geq3$, its radical does not square to zero, since $a_2a_1\neq0$. As a check, we also verified this description by exact computation for several choices of weights and signs.
\end{example}

\begin{example}[The twist in Theorem~\ref{thm:structure}(3) is necessary]\label{ex:structure-twist}
Let $Q$ have two vertices $1,2$ of type $(\RR,+1)$, a loop $\ell$ at $1$ fixed by the transpose with sign $+1$, and arrows $a:1\to2$ and $a^{\top}:2\to1$. Write $b=a^{\top}=a^*$. The ideal $I=(\ell^2-ba,\;ab,\;a\ell,\;\ell b)$ is admissible: it is $*$-stable, since the involution fixes $\ell^2-ba$ and exchanges $a\ell$ with $\ell b$, and it contains every path of length three. The algebra $B=\RR Q/I$ has the basis $e_1,e_2,\ell,a,b,\ell^2$. Let $t=-e_1+e_2$, so that the involution of $B^t$ maps $a$ to $a^{*_t}=-b$. In $B$ there are $x\in e_2\rad(B)e_1$ and a self-adjoint $y\in e_1\rad(B)e_1$ with $x^*x=y^2\neq0$, namely $x=a$ and $y=\ell$. In $B^t$ there are no $x$ and $y$ with $x^{*_t}x=y^2\neq0$: $x=\gamma a$ gives $x^{*_t}x=-\gamma^2\ell^2$, whereas $y=\alpha\ell+\beta\ell^2$ gives $y^2=\alpha^2\ell^2$. A $*$-isomorphism $B\to B^t$ may be assumed to fix $e_1$ and $e_2$, by Lemma~\ref{lem:unitary-lifting}(2) and because $e_1Be_1$ and $e_2Be_2$ have different dimensions; it would then preserve this property. Hence $B\not\cong B^t$. Nevertheless, for $V_1=\RR^2$ with $\eta_1=\diag(1,-1)$ and $V_2=\RR$, the $17$-dimensional algebras $\End(V,\eta)\otimes_{\RR^2}B$ and $\End(V,\eta)\otimes_{\RR^2}B^t$ are isomorphic, as Theorem~\ref{thm:structure}(3) predicts; an isomorphism exchanges the two basis vectors of $V_1$.
\end{example}

\subsection{Proof of the Structure Theorem}
\label{subsec:structure-proof}
\begin{proof}[Proof of Theorem~\ref{thm:structure}]
Let $N=\rad(A)$.

\emph{Step 1: matrix units.} By Theorem~\ref{thm:wedderburn-malcev}, $A=S\oplus N$ for a $*$-stable subalgebra $S\cong A/N$; let $S=\prod_{c\in\mathcal O}S_c$ be its decomposition into $*$-simple factors. By Lemma~\ref{lem:star-simple}, with $K=K_c$ because the type of $S_c$ is that of the orbit $c$, and since $M_{n,\eta}(K)=M_{n,-\eta}(K)$ and the signs may be permuted by conjugating with a permutation matrix, there are $*$-isomorphisms $S_c\cong M_{n_c,\eta_c}(K_c)$ with normalised signs $\eta_c=(\eta_{c,1},\dots,\eta_{c,n_c})$: all equal to $+1$ if $c$ is not of signature type, and otherwise $p_c$ entries $+1$ followed by $q_c\leq p_c$ entries $-1$. In particular $\eta_{c,1}=1$. Let $E^c_{kl}\in S_c$ correspond to the matrix with entry $1$ at position $(k,l)$. Then $E^c_{kl}E^c_{lm}=E^c_{km}$ and $(E^c_{kl})^*=\eta_{c,k}\eta_{c,l}E^c_{lk}$.

\emph{Step 2: a basic $*$-algebra.} Put $f_c=E^c_{11}$, $f=\sum_cf_c$ and $u_{c,k}=E^c_{k1}$. We write $i,j$ for pairs $(c,k)$, and $c(i)=c$, $\eta_i=\eta_{c,k}$ and $u_i=u_{c,k}$ for $i=(c,k)$. Then $f$ is a self-adjoint idempotent, $u_i=u_if_{c(i)}$, and
\begin{equation}\label{eq:structure-units}
  u_i^*u_j=\delta_{ij}\,\eta_i\,f_{c(i)},\qquad\sum_i\eta_i\,u_iu_i^*=1 .
\end{equation}
Let $B=fAf$, a $*$-algebra with unit $f$ and an $\RR^{\mathcal O}$-$*$-algebra with $1_c=f_c$. We have $B=fSf\oplus fNf$, where $fSf=\prod_cf_cS_cf_c\cong\prod_cK_c$ is semisimple and $fNf$ is a nilpotent ideal of $B$. Hence $\rad(B)=fNf$ and $B/\rad(B)\cong\prod_cK_c$, and $B$ is basic as a $*$-algebra by Lemma~\ref{lem:star-simple}(3).

\emph{Step 3: $A$ is a matrix algebra over $B$.} Let $V_c=\RR^{n_c}$ with the form $\diag(\eta_c)$, and define $\Phi:A\to\End(V,\eta)\otimes_{\RR^{\mathcal O}}B$, in the coordinates of \eqref{eq:signed-adjoint}, by
\[ \Phi(a)_{ij}=\eta_i\,u_i^*\,a\,u_j\;\in\;f_{c(i)}Af_{c(j)}=B_{c(i)c(j)} . \]
By \eqref{eq:structure-units},
\[ (\Phi(a)\Phi(b))_{ij}=\eta_iu_i^*a\Bigl(\sum_k\eta_ku_ku_k^*\Bigr)bu_j=\Phi(ab)_{ij},\qquad\Phi(1)_{ij}=\eta_iu_i^*u_j=\delta_{ij}f_{c(i)}, \]
and $\Phi(a^*)_{ij}=\eta_i(u_j^*au_i)^*=\eta_i\eta_j(\Phi(a)_{ji})^*$, which is the involution \eqref{eq:signed-adjoint}. The map $X\mapsto\sum_{i,j}u_i\,x_{ij}\,\eta_ju_j^*$ is inverse to $\Phi$, again by \eqref{eq:structure-units}. So $\Phi$ is a $*$-isomorphism. For a single orbit and positive signs this is Proposition~\ref{prop:matrix-structure}.

\emph{Step 4: $B$ has doubly-signed quiver $Q$.} Let $N_B=\rad(B)=fNf$. By \eqref{eq:structure-units},
\[ fN^2f=\sum_i\eta_i\,(fNu_i)(u_i^*Nf)\subseteq(fNf)^2=N_B^2 , \]
since $fNu_i=fNu_if$ and $u_i^*Nf=fu_i^*Nf$; the reverse inclusion is clear. So $B/N_B^2=\bar f(A/N^2)\bar f$, where $\bar f$ is the image of $f$. Moreover $(A/N^2)\bar f(A/N^2)=A/N^2$, because $1=\sum_i\eta_iu_ifu_i^*$. Thus $B/N_B^2$ is a basic $*$-algebra of $A/N^2$, and $Q(B)=Q(A/N^2)=Q$ by Theorem~\ref{thm:classification}(4).

\emph{Step 5: a presentation of $B$.} Let $\bar S=B/N_B$ and $\bar M=N_B/N_B^2$. Then $\bar M$ is an $\bar S$-bimodule, and the involution of $B$ induces involutions of $\bar S$ and of $\bar M$ with $(srs')^*=s'^*r^*s^*$. So the tensor algebra $T=T_{\bar S}(\bar M)$, with the involution that reverses tensors, is a $*$-algebra. The subalgebra $S'=fSf$ of $B$ is $*$-stable with $B=S'\oplus N_B$; let $\sigma:\bar S\to S'$ be the inverse of the quotient map, a $*$-isomorphism with $\sigma(\bar f_c)=f_c$. Since $S'$ is semisimple and $\RR$ has characteristic $0$, $S'$ is separable, and every $S'$-bimodule is semisimple~\cite[Ch.~10]{Pierce1982}. So $N_B$ has an $S'$-bimodule projection $p$ with kernel $N_B^2$. Then $p'=\frac12(p+{*}\circ p\circ{*})$ is an $S'$-bimodule map, because $S'^*=S'$. It vanishes on $N_B^2$ and induces the identity on $N_B/N_B^2$, because $(N_B^2)^*=N_B^2$. So it is a projection with kernel $N_B^2$, and it commutes with the involution. Its image $W$ is a $*$-stable sub-bimodule with $N_B=W\oplus N_B^2$. The inverse of the quotient map $W\to\bar M$ is a map $\rho:\bar M\to W$ with $\rho(srs')=\sigma(s)\rho(r)\sigma(s')$ and $\rho(r^*)=\rho(r)^*$. By the universal property of the tensor algebra, there is an algebra homomorphism $\varphi:T\to B$ with $\varphi|_{\bar S}=\sigma$ and $\varphi(r_1\otimes\cdots\otimes r_\ell)=\rho(r_1)\cdots\rho(r_\ell)$. It commutes with the involutions, because it does so on $\bar S$ and $\bar M$ and both involutions reverse products.

The homomorphism $\varphi$ is surjective. From $N_B=W+N_B^2$ we obtain, by induction, $N_B=W+W^2+\cdots+W^{m-1}+N_B^m$ for all $m$, and $N_B^m\subseteq fN^mf=0$ if $N^m=0$; so $\varphi(T)\supseteq S'+N_B=B$. Its kernel $J$ is $*$-stable. It contains $T_{\geq m}$, because $\varphi(T_{\geq m})\subseteq N_B^m$. It is contained in $T_{\geq2}$, because $\varphi$ maps $\bar S\oplus\bar M$ isomorphically onto $S'\oplus W$, $\varphi(T_{\geq2})\subseteq N_B^2$, and $B=S'\oplus W\oplus N_B^2$.

Finally, $B/N_B^2$ and $B_Q$ are basic radical-square-zero $*$-algebras with the same doubly-signed quiver (Step~4). So by Theorem~\ref{thm:classification} there is a $*$-isomorphism $\chi:B/N_B^2\to B_Q$. It maps the radical $\bar M$ onto $M_Q$ and induces a $*$-isomorphism $\bar S\to B_Q/M_Q\cong S_Q$, compatibly with the bimodule structures. Hence it induces a degree-preserving $*$-isomorphism $T\to\RR Q$. The image $I$ of $J$ is an admissible ideal, $I\supseteq\RR Q_{\geq m}$ if $N^m=0$, and $\RR Q/I\cong B$.

\emph{Step 6: proof of (1).} In degree $0$, the isomorphism $B\cong\RR Q/I$ of Step~5 restricts to a $*$-isomorphism $fSf\cong S_Q$. This maps each $f_c$ to the unit $1_{\pi(c)}$ of a factor $K_{\pi(c)}\cong K_c$, for a type-preserving bijection $\pi$ of $\mathcal O$. Combining this with Step~3, we get
\[ A\cong\End(V,\eta)\otimes_{\RR^{\mathcal O}}B\cong\End(V',\eta')\otimes_{\RR^{\mathcal O}}\RR Q/I,\qquad (V'_{\pi(c)},\eta'_{\pi(c)})=(V_c,\eta_c). \]
This proves (1); the remaining claims of (1) were shown in Steps~2 and~5.

\emph{Proof of (2).} Consider any normalised decomposition $A\cong\End(V',\eta')\otimes_{\RR^{\mathcal O'}}B'$ with $B'=\RR Q'/I'$ and all $V'_c\neq0$, and write $n'_c=\dim V'_c$. The ideal $\End(V',\eta')\otimes_{\RR^{\mathcal O'}}\rad(B')$ is nilpotent. Its square is $\End(V',\eta')\otimes_{\RR^{\mathcal O'}}\rad(B')^2$, and the quotient by it is $\End(V',\eta')\otimes_{\RR^{\mathcal O'}}S_{Q'}=\prod_cM_{n'_c,\eta'_c}(K_c)$, which is semisimple. Hence $A/N\cong\prod_cM_{n'_c,\eta'_c}(K_c)$, and $A/N^2\cong\End(V',\eta')\otimes_{\RR^{\mathcal O'}}B_{Q'}$ because $B'/\rad(B')^2=\RR Q'/\RR Q'_{\geq2}=B_{Q'}$. Since every $V'_c$ is non-zero, the idempotent $\sum_ce_{(c,1)(c,1)}\otimes1_c$ is full, and the corner at it is $B_{Q'}$, because $\eta'_{c,1}=1$. So Theorem~\ref{thm:classification}(4) gives an isomorphism $\theta':Q'\to Q(A/N^2)=Q$, under which the $*$-simple factor of $A/N$ belonging to $\theta'(c)$ is $M_{n'_c,\eta'_c}(K_c)$. Applying this to both decompositions in (2), and comparing the two descriptions of each factor, Lemma~\ref{lem:star-simple} shows that $n'_c=n_{\theta(c)}$ for the resulting isomorphism $\theta:Q'\to Q$. It also shows that, if $c$ is of signature type, the signatures of $\eta'_c$ and $\eta_{\theta(c)}$ agree up to order, hence agree, because both are normalised.

\emph{Proof of (3).} Let $f_c\in A$ ($c\in\mathcal O$) and $f'_{c'}\in A$ ($c'\in\mathcal O'$) be the self-adjoint idempotents corresponding to $e_{(c,1)(c,1)}\otimes1_c$ under the two isomorphisms, and put $f=\sum_cf_c$ and $f'=\sum_{c'}f'_{c'}$. Since $\eta_{c,1}=\eta'_{c',1}=1$, the corners of the two decompositions at these idempotents are $B$ and $B'$, so $fAf\cong B$ and $f'Af'\cong B'$. By the proof of (2), each $*$-simple factor $F$ of $A/N$ contains the image of exactly one $f_c$ and of exactly one $f'_{c'}$, with $\theta(c')=c$, and $F\cong M_{n_c,\eta_c}(K_c)$ is the $*$-algebra of the Hermitian $K_c$-module $(K_c^{n_c},\eta_c)$. The images $\bar f_c$ and $\bar f'_{c'}$ are self-adjoint idempotents of $F$ with corners isomorphic to $K_c$, that is, orthogonal projections onto free submodules of rank one. A Hermitian form on a free module of rank one is $\langle\lambda\rangle$ for an invertible self-adjoint $\lambda\in K_c$, so by Lemma~\ref{lem:star-simple}(1) it is isometric to $\langle1\rangle$ if $c$ is not of signature type, and to $\langle1\rangle$ or $\langle-1\rangle$ if it is. By Witt's cancellation theorem, a consequence of Theorem~\ref{thm:hermitian-ks}, two such projections are conjugate by a unitary of $F$ if and only if the forms on their images are isometric. So if $c$ is not of signature type, $\bar f_c$ and $\bar f'_{c'}$ are conjugate. If it is, the projections fall into at most two conjugacy classes, and by Sylvester's law of inertia every family of $n_c$ pairwise orthogonal ones contains $p_c$ members of one class and $q_c$ of the other. Both $\bar f_c$ and $\bar f'_{c'}$ belong to a class that contributes $p_c=\max(p_c,q_c)$ members, because $\eta_{c,1}=\eta'_{c',1}=1$ and $(p'_{c'},q'_{c'})=(p_c,q_c)$. Hence $\bar f_c$ and $\bar f'_{c'}$ are conjugate unless $c$ is of signature type, $p_c=q_c$, and they lie in different classes. In that case put $\tau_c=-1$, and let $f''_c$ correspond to $e_{(c,k)(c,k)}\otimes1_c$ under the first isomorphism, for some $k$ with $\eta_{c,k}=-1$. Otherwise put $\tau_c=1$ and $f''_c=f_c$. By \eqref{eq:signed-adjoint}, the corner of $\End(V,\eta)\otimes_{\RR^{\mathcal O}}B$ at the idempotent corresponding to $f''=\sum_cf''_c$ is $B$ with the involution multiplied by $\tau_c\tau_d$ on $B_{cd}$; that is, $f''Af''\cong B^t$. The images of $f''$ and $f'$ are conjugate in every factor $F$, hence by a unitary of $A/N$, and Lemma~\ref{lem:unitary-lifting} gives $B'\cong f'Af'\cong f''Af''\cong B^t$.

Conversely, $\End(V,\eta)\otimes_{\RR^{\mathcal O}}B^t$ and $\End(V,t\eta)\otimes_{\RR^{\mathcal O}}B$, where $t\eta$ multiplies $\eta_c$ by $\tau_c$, have the same underlying algebra and, by \eqref{eq:signed-adjoint}, the same involution. If $\tau_c=-1$ only where $p_c=q_c$, there is an isometry $(V,t\eta)\to(V,\eta)$ that preserves each $V_c$. Conjugation by it is a $*$-isomorphism $\End(V,t\eta)\otimes_{\RR^{\mathcal O}}B\to\End(V,\eta)\otimes_{\RR^{\mathcal O}}B$.
\end{proof}

\subsection{Uniqueness of the Relations}
\label{subsec:relations-uniqueness}
The ideal $I$ in Theorem~\ref{thm:structure} is not unique. For example, let $Q$ have one vertex of type $(\RR,+1)$ and two loops $x,y$ fixed by the transpose with sign $+1$, so that $\RR Q=\RR\langle x,y\rangle$ with $x^*=x$ and $y^*=y$. The $*$-automorphism of $\RR Q$ that fixes $x$ and maps $y$ to $x+y$ maps the admissible ideal $I=(x^2-y^2,\,xy,\,yx)$ to an admissible ideal that contains $x^2+xy\notin I$, with an isomorphic quotient. The following proposition shows that this is the only kind of ambiguity. It uses a $*$-stable version of the conjugacy part of the Wedderburn--Malcev theorem.

\begin{lemma}\label{lem:star-conjugacy}
Let $B$ be a finite-dimensional real $*$-algebra, and let $S_1,S_2$ be $*$-stable subalgebras with $B=S_1\oplus\rad(B)=S_2\oplus\rad(B)$. Then $S_2=uS_1u^*$ for a unitary $u\in1+\rad(B)$.
\end{lemma}
\begin{proof}
By Theorem~\ref{thm:wedderburn-malcev}, $S_2=gS_1g^{-1}$ for some $g\in1+\rad(B)$. Applying the involution gives $S_2=g^{-*}S_1g^*$, so $h=g^*g$ satisfies $hS_1h^{-1}=S_1$. For $s\in S_1$, the element $hsh^{-1}-s$ lies in $S_1\cap\rad(B)=0$, because $h\equiv1$ modulo $\rad(B)$. So $h$ commutes with $S_1$, and so does the polynomial $h^{-1/2}$ in $h$. Then $u=gh^{-1/2}$ is unitary, and $uS_1u^*=gS_1g^{-1}=S_2$.
\end{proof}

\begin{proposition}\label{prop:ideal-uniqueness}
Let $I$ and $I'$ be admissible ideals of $\RR Q$ with $\RR Q_{\geq m}\subseteq I\cap I'$. If $\RR Q/I\cong\RR Q/I'$ as $*$-algebras, then $I'/\RR Q_{\geq m}=\alpha(I/\RR Q_{\geq m})$ for a $*$-automorphism $\alpha$ of $\RR Q/\RR Q_{\geq m}$.
\end{proposition}
\begin{proof}
Write $T=\RR Q/\RR Q_{\geq m}$ and $B'=\RR Q/I'$. Let $\pi':T\to B'$ be the quotient map, and let $\pi=\chi\circ\pi_I$, where $\pi_I:T\to\RR Q/I$ is the quotient map and $\chi:\RR Q/I\to B'$ is the given isomorphism. Both $\pi$ and $\pi'$ are surjective $*$-homomorphisms that map $S_Q$ isomorphically onto $*$-stable complements of $\rad(B')$, map $T_{\geq1}$ onto $\rad(B')$, and induce isomorphisms $M_Q\to\rad(B')/\rad(B')^2$. By Lemma~\ref{lem:star-conjugacy}, after replacing $\pi$ by $x\mapsto u\pi(x)u^*$ for a unitary $u$, which does not change its kernel, we may assume $\pi(S_Q)=\pi'(S_Q)$. Then $\beta=(\pi'|_{S_Q})^{-1}\circ\pi|_{S_Q}$ is a $*$-automorphism of $S_Q$, and $\pi(srs')=\pi'(\beta(s))\,\pi(r)\,\pi'(\beta(s'))$ for $r\in M_Q$. Since $S_Q$ is separable, $M_Q$ is projective as an $S_Q$-bimodule, even with the action twisted by $\beta$. So $\pi|_{M_Q}$ lifts through the surjection $\pi':T_{\geq1}\to\rad(B')$ to a map $\lambda_0:M_Q\to T_{\geq1}$ with $\lambda_0(srs')=\beta(s)\lambda_0(r)\beta(s')$ and $\pi'\circ\lambda_0=\pi|_{M_Q}$. The map $\lambda=\frac12(\lambda_0+{*}\circ\lambda_0\circ{*})$ has the same two properties and commutes with the involution. By the universal property of the tensor algebra, $\beta$ and $\lambda$ define a $*$-endomorphism $\alpha$ of $T$; it is well defined on $T$ because $\lambda(M_Q)\subseteq T_{\geq1}$. Moreover $\pi'\circ\alpha=\pi$, because both sides are homomorphisms that agree on $S_Q$ and $M_Q$. On the associated graded algebra, $\alpha$ is induced by $\beta$ and by the degree-one part of $\lambda$. The latter is bijective, because $\pi$ and $\pi'$ induce isomorphisms $M_Q\to\rad(B')/\rad(B')^2$; hence $\alpha$ is bijective. Finally, $\pi'\circ\alpha=\pi$ gives $\alpha(I/\RR Q_{\geq m})\subseteq I'/\RR Q_{\geq m}$, with equality because both quotients have dimension $\dim B'$.
\end{proof}
%======================================================================
%   CHAPTER 6: DOUBLY-SIGNED QUIVERS ENCODE MATRIX DECOMPOSITIONS
%======================================================================
\chapter{Matrix Decompositions from Quivers}
\label{chap:indexing}

This chapter proves the two outcomes of Section~\ref{sec:outcomes} for the encoded matrices: parts~(2) and~(3) of Theorem~\ref{thm:A}, which say what a doubly-signed quiver \emph{encodes}, and, with Proposition~\ref{prop:golub-kahan}, parts~(1) to~(3) of Theorem~\ref{thm:B}, which say what the programs of Section~\ref{sec:outcomes} compute on the encoded matrices.

Let $Q$ be a doubly-signed quiver. Its path $*$-algebra $\RR Q$ is graded by the length of paths: the vertices have degree $0$ and the arrows degree $1$. Consider the Hermitian matrices of degree $1$, that is, the matrices indexed by positions (Corollary~\ref{cor:positions}) whose entries are linear combinations of arrows. Theorem~\ref{thm:encoding} shows that unitary similarity of such matrices is the same as similarity by unitary matrices of degree $0$, which form a product of classical groups. The decomposition indexed by $Q$ is, by definition, a normal form for this problem. Section~\ref{sec:catalogue} works through a catalogue of examples.

Section~\ref{sec:encoded-algorithms} runs the generic algorithms on Hermitian matrices of degree $1$, in the completed localisation $\Lambda$ of Chapter~\ref{division-problem}, and identifies what Householder reduction and the QR algorithm compute there (Propositions~\ref{prop:householder-graded}, \ref{prop:qr-steps} and~\ref{prop:2R-algorithms}). Section~\ref{sec:beyond-radical} treats Hermitian matrices with a part of degree $0$, which leads to perturbation theory, and Section~\ref{sec:wild-behaviour} discusses what representation type says about generic algorithms.

\section{The Encoding Theorem}
\label{sec:encoding}
Throughout this chapter, $Q$ is a doubly-signed quiver, $\mathcal O=\mathcal O(Q)$ is the set of orbits of its vertices under the transpose and the conjugation, and $A=\RR Q/I$, where $I\subseteq\RR Q_{\geq2}$ is an ideal that is stable under the involution and spanned by homogeneous elements, for example $I=0$ or $I=\RR Q_{\geq m}$ with $m\geq2$ (Definition~\ref{def:path-star-algebra}). Then $A=\bigoplus_{d\geq0}A_d$ is graded by the path degree, the involution preserves degrees, $A_0=S_Q=\prod_{c\in\mathcal O}K_c$, and $A_1=M_Q$. For $I=\RR Q_{\geq2}$ we obtain $A=B_Q$, whose radical squares to zero; this is the case of Section~\ref{sec:unification-example}. Fix real spaces $V_c$ with non-degenerate symmetric forms $\eta_c$, put $n_c=\dim V_c$, and let
\[ E=\End(V,\eta)\otimes_{\RR^{\mathcal O}}A \]
be the $*$-algebra of matrices indexed by the positions $c_k$ ($c\in\mathcal O$, $1\leq k\leq n_c$) with entries in $A$ (Definition~\ref{def:matrix-like} and \eqref{eq:signed-adjoint}). It is graded by $E_d=\End(V,\eta)\otimes_{\RR^{\mathcal O}}A_d$, and we write $X=X_0+X_1+\cdots$ for the homogeneous parts of $X\in E$. Let $\Herm(E_1)$ be the set of Hermitian elements of degree $1$, and $\Un(E)$ the unitary group of $E$.

For each orbit $[\alpha]$ of arrows under the transpose and the conjugation, let $M_{[\alpha]}\subseteq M_Q$ be the real form of the sum of the copies of $M_{d_t\times d_s}(\CC)$ that belong to the arrows of $[\alpha]$ in the block model (Proposition~\ref{prop:block-model}).

\begin{theorem}[Encoding theorem]\label{thm:encoding}
In this situation:
\begin{enumerate}
  \item If $U\in\Un(E)$, then $U_0\in\Un(E_0)$. Conversely, $\Un(E_0)\subseteq\Un(E)$.
  \item Let $H,H'\in\Herm(E_1)$. If $H'=UHU^*$ for some $U\in\Un(E)$, then $H'=U_0HU_0^*$. If $A=B_Q$, then $UHU^*=U_0HU_0^*$ for all $U\in\Un(E)$ and all $H\in\Herm(E_1)$. In either case, two elements of $\Herm(E_1)$ are unitarily similar over $A$ if and only if they lie in the same orbit of $\Un(E_0)$.
  \item $\Un(E_0)=\prod_{c\in\mathcal O}\Un_c$, where $\Un_c$ is the unitary group of $M_{n_c,\eta_c}(K_c)$. If $\eta_c$ is positive definite, or $c$ is not of signature type, then $\Un_c$ is isomorphic to the classical group of Table~\ref{tab:vertex-types}; otherwise it is an indefinite form of it (Remark~\ref{rem:indefinite}). The space $\Herm(E_1)$ is the direct sum, over the orbits $[\alpha]$ of arrows, of $\Un(E_0)$-stable subspaces $\mathcal H_{[\alpha]}$. On $\mathcal H_{[\alpha]}$ only the factors $\Un_c$ and $\Un_{c'}$ act, where $c$ and $c'$ are the orbits of the source and the target of $\alpha$:
  \begin{enumerate}
    \item if $c\neq c'$, then $\mathcal H_{[\alpha]}\cong M_{n_{c'}\times n_c}(P_{[\alpha]})$, where $P_{[\alpha]}=1_{c'}M_{[\alpha]}1_c$ is a $K_{c'}$--$K_c$-bimodule, and $(g_c,g_{c'})$ acts by $W\mapsto g_{c'}Wg_c^{-1}$;
    \item if $c=c'$, then $\mathcal H_{[\alpha]}$ is the space of matrices $W=(w_{kl})$ of size $n_c$ with entries in $M_{[\alpha]}$ and $w_{kl}=\eta_k\eta_l\,w_{lk}^*$, and $g\in\Un_c$ acts by $W\mapsto gWg^{-1}$.
  \end{enumerate}
\end{enumerate}
\end{theorem}

\begin{proof}
(1) A product of homogeneous elements is homogeneous of the sum of their degrees, and all degrees are non-negative, so $(XY)_0=X_0Y_0$. Moreover $(X^*)_0=(X_0)^*$, because the involution preserves degrees. Applying this to $U^*U=UU^*=1$ gives $U_0^*U_0=U_0U_0^*=1$. The converse is clear.

(2) For $H\in E_1$ and $U\in E$, the part of degree $1$ of $UHU^*$ is $U_0HU_0^*$, because the parts of $U$ and $U^*$ have non-negative degrees, and only their parts of degree $0$ combine with $H$ to degree $1$. Hence $H'=(UHU^*)_1=U_0HU_0^*$. If $A=B_Q$, then $E_d=0$ for $d\geq2$, and $UHU^*$ has no part of degree $0$, so $UHU^*=(UHU^*)_1$. The last statement follows, by (1).

(3) The orbit idempotents $1_c$ are central in $S_Q$, so $E_0$ is the product of the algebras $\End(V_c,\eta_c)\otimes_\RR K_c=M_{n_c,\eta_c}(K_c)$, and its unitary group is the product of theirs. In the block model, $M_Q\otimes_\RR\CC=M_\CC$ is the direct sum of copies of $M_{d_t\times d_s}(\CC)$, one for each arrow, and by \eqref{eq:block-star} the involution and the conjugation map the copy of $\alpha$ onto the copies of $\alpha^{\top}$ and $\bar\alpha$. So the sum of the copies over an orbit is stable under both, and it is the complexification of $M_{[\alpha]}$. Hence $M_Q=\bigoplus_{[\alpha]}M_{[\alpha]}$ is a direct sum of $S_Q$-sub-bimodules that are stable under the involution, $E_1$ is the direct sum of the spaces $\End(V,\eta)\otimes_{\RR^{\mathcal O}}M_{[\alpha]}$, each stable under the involution and under multiplication by $E_0$ on both sides, and $\Herm(E_1)$ is the direct sum of their Hermitian parts $\mathcal H_{[\alpha]}$, each stable under $X\mapsto UXU^*$ for $U\in\Un(E_0)$. The arrows of $[\alpha]$ go from vertices of $c$ to vertices of $c'$ or, after transposition, from $c'$ to $c$, so $M_{[\alpha]}\subseteq1_{c'}M_Q1_c+1_cM_Q1_{c'}$, and only $\Un_c$ and $\Un_{c'}$ act on $\mathcal H_{[\alpha]}$. If $c\neq c'$, an element $X$ of $\mathcal H_{[\alpha]}$ has the block $W=1_{c'}X1_c\in M_{n_{c'}\times n_c}(P_{[\alpha]})$, its block $1_cX1_{c'}=W^*$ is determined by $W$, and its other blocks vanish; every $W$ occurs, since $W+W^*$ is Hermitian. For $U=(g_c,g_{c'})$ the block $1_{c'}(UXU^*)1_c$ is $g_{c'}Wg_c^*$, and $g_c^*=g_c^{-1}$. If $c=c'$, then $\mathcal H_{[\alpha]}$ consists of the Hermitian elements of $\End(V_c,\eta_c)\otimes_\RR M_{[\alpha]}$, which by \eqref{eq:signed-adjoint} form the stated space, and $g$ acts by $X\mapsto gXg^*=gXg^{-1}$.
\end{proof}

Theorem~\ref{thm:encoding} answers the question of what a doubly-signed quiver encodes. \emph{Given $Q$ and the spaces $(V_c,\eta_c)$, the decomposition indexed by $Q$ is a normal form for the Hermitian matrices of degree $1$ under unitary similarity.} By part~(2), this problem does not depend on the homogeneous relations: it is the same over $B_Q$, over $\RR Q/\RR Q_{\geq m}$ and over $\RR Q$ itself, and it is the action of $\Un(E_0)$ on $\Herm(E_1)$. By part~(3), it is a direct sum of classical problems, one for each orbit of arrows, coupled through the groups at the vertices. The numbers $n_c$ are the sizes of the matrices; for example, $n_u\neq n_v$ gives rectangular matrices in the singular value decomposition. In the language of Section~\ref{sec:star-reps}, $\Herm(E_1)$ is a space of $*$-representations of $Q$, and unitary similarity by $E_0$ is their isomorphism. By part~(2), every algorithm that transforms Hermitian matrices by unitary similarity over $A$, such as a generic Jacobi algorithm, or the QR algorithm taken in double steps (over the quivers of the singular value decomposition a single step exchanges the two orbits, Proposition~\ref{prop:qr-steps}), acts on Hermitian matrices of degree $1$ through the classical groups at the vertices; Section~\ref{sec:encoded-algorithms} makes this concrete for algorithms run in the completed localisation.

\begin{remark}[Indefinite forms]\label{rem:indefinite}
By Lemma~\ref{lem:star-simple}(1), $\eta_c$ may be taken positive definite at every orbit $c$ that is not of signature type, since there every invertible self-adjoint element of $K_c$ has the form $uu^*$. At an orbit of signature type, the signature of $\eta_c$ is part of the data: $\Un_c$ is the indefinite group $\Or(p,q)$, $\Un(p,q)$ or $\Sp(p,q)$, and the problems of Theorem~\ref{thm:encoding}(3) become their indefinite analogues. By Theorem~\ref{thm:hermitian-ks}, the signature is an invariant of the Hermitian module and belongs to the data of the representation rather than to the quiver.
\end{remark}

\section{A Catalogue of Decompositions and Their Classical Analogues}
\label{sec:catalogue}
We now apply Theorem~\ref{thm:encoding} to a series of doubly-signed quivers. In each case we describe the $*$-algebra $B=B_Q$, compute the Hermitian matrices of degree $1$ and the action of the unitary matrices of degree $0$, and identify the classical decomposition. Each of these decompositions can also be derived classically, case by case, often as the Hermitian eigendecomposition of a larger block matrix: for example, the singular value decomposition of $W$ from the eigendecomposition of the symmetric Jordan--Wielandt matrix $\bigl(\begin{smallmatrix}0&W^T\\W&0\end{smallmatrix}\bigr)$~\cite{GolubVanLoan2013}, and the quaternionic decompositions from the representation of a quaternionic matrix by a complex matrix of twice the size~\cite{Zhang1997}. The encoding gives a uniform derivation, not the only one. By Theorem~\ref{thm:encoding}(2), the same problem arises over $\RR Q/I$ for every homogeneous $I$. Unless stated otherwise, every orbit has $n$ positions and the forms $\eta_c$ are positive definite, so that $E=M_n(B)$ with the conjugate transpose. We write $B=S\oplus R$ with $S=S_Q$ and $R=M_Q=\rad(B)$, and $\Herm_n(R)$ for the Hermitian matrices with entries in $R$; these are the Hermitian matrices of degree $1$. Diagrams of the quivers are in Figures~\ref{fig:signed-quiver-example} and~\ref{fig:dsq-examples}, and the results are summarised in Table~\ref{tab:catalogue}. Following the convention of Chapter~\ref{chap:quivers}, in labels we write the division algebra $\RR$, $\CC$ or $\HH$ in place of the sign $\epsilondiv=+1,0,-1$. Section~\ref{sec:encoded-algorithms} describes what the generic algorithms compute for each entry (Table~\ref{tab:algorithms}).

\subsection{The Symmetric Eigendecomposition}
\label{subsec:cat-symmetric}
\emph{Quiver.} One vertex $v$ of type $(\RR,+1)$ with one loop $x$ fixed by the transpose with sign $+1$ (Figure~\ref{fig:signed-quiver-example}(a)).

\emph{Algebra.} $B=\RR[x]/(x^2)$, the dual numbers, with the trivial involution $x^*=x$.

\emph{Encoding.} A matrix $H=xW\in M_n(R)$ is Hermitian if and only if $W^T=W$. The unitary group of $S=\RR$ is $\Or(n)$, and $U_0(xW)U_0^*=x\,(PWP^T)$ for $U_0=P$.

\emph{Classical analogue.} Orthogonal similarity of real symmetric matrices, whose normal form is the spectral theorem $W=P\Lambda P^T$~\cite{HornJohnson2013}. The diagonal entries of the normal form, $x\lambda_i$, are the ``eigenvalues'' of $H$ over $B$. Section~\ref{sec:beyond-radical} shows that the full Hermitian matrices over this algebra, $H=X+xW$, are the first-order perturbations of symmetric matrices, recovering~\cite{Gutin2022DualNumbers} and Theorem~\ref{thm:spectral-power-series} modulo $\varepsilon^2$.

\subsection{The Normal Form of Skew-Symmetric Matrices}
\label{subsec:cat-skew}
\emph{Quiver.} One vertex of type $(\RR,+1)$ with one loop $h$ fixed by the transpose with sign $-1$ (Figure~\ref{fig:signed-quiver-example}(b)).

\emph{Algebra.} The dual numbers $B=\RR[h]/(h^2)$ with the involution $h^*=-h$.

\emph{Encoding.} Since $(hW)^*=W^Th^*=-hW^T$, the matrix $hW$ is Hermitian if and only if $W$ is skew-symmetric. The action is again $W\mapsto PWP^T$ with $P\in\Or(n)$.

\emph{Classical analogue.} Orthogonal similarity of real skew-symmetric matrices. Its normal form is
\[ W=P\Bigl(\bigoplus_j\begin{pmatrix}0&\sigma_j\\-\sigma_j&0\end{pmatrix}\oplus0\Bigr)P^T,\qquad\sigma_j>0, \]
see~\cite[Corollary~2.5.11]{HornJohnson2013}. Here diagonalisation is impossible: the only Hermitian $1\times1$ matrix with radical entry is $0$, because $(hw)^*=-hw$. Apart from the zero block, the indecomposable summands therefore have dimension two, which explains the $2\times2$ blocks. These are the parity blocks of Corollary~\ref{cor:twisted-hermitian} reduced modulo $\varepsilon^2$.

\subsection{The Real Singular Value Decomposition}
\label{subsec:cat-svd}
\emph{Quiver.} Two vertices $u,v$ of type $(\RR,+1)$ and arrows $a:u\to v$ and $a^{\top}:v\to u$ exchanged by the transpose (Figures~\ref{fig:signed-quiver-example}(d) and~\ref{fig:dsq-examples}(a)).

\emph{Algebra.} The four-dimensional algebra $B$ of Section~\ref{sec:unification-example}, with basis $e_u,e_v,a,b$, where $b=a^*=a^{\top}$.

\emph{Encoding.} $\Herm_n(R)=\{aW+bW^T:W\in M_n(\RR)\}$, $\Un_n(S)=\Or(n)\times\Or(n)$, and the action is $W\mapsto RWP^T$ (equation \eqref{eq:svd-action}). With $n$ positions for $u$ and $m$ for $v$, $W\in M_{m\times n}(\RR)$. This is case~(a) of Theorem~\ref{thm:encoding}(3), with $P_{[a]}=\RR a$.

\emph{Classical analogue.} The singular value decomposition $W=R^T\Sigma P$~\cite{GolubVanLoan2013}. The vertices contribute the two orthogonal factors, and the pair of arrows contributes the matrix being decomposed.

\subsection{The Complex Singular Value Decomposition}
\label{subsec:cat-complex-svd}
\emph{Quiver.} As for the real SVD, but with two complex vertices of unitary type; after complexification each becomes a pair exchanged both by conjugation and by the transpose (Figure~\ref{fig:dsq-examples}(c)).

\emph{Algebra.} $S=\CC\times\CC$, each factor with complex conjugation as involution, and $R=\CC a\oplus\CC a^*$, where $a$ goes from $u$ to $v$, $za=az$ for $z\in\CC$, and $(za)^*=\bar z\,a^*$.

\emph{Encoding.} $\Herm_n(R)=\{aW+a^*W^H:W\in M_n(\CC)\}$, $\Un_n(S)=\Un(n)\times\Un(n)$, and the action is $W\mapsto RWP^H$.

\emph{Classical analogue.} The complex singular value decomposition $W=R^H\Sigma P$. The only change compared with the real SVD is the vertex type, which replaces $\Or(n)$ by $\Un(n)$; this is exactly what the conjugation in a doubly-signed quiver records. In the diagram, the conjugate pair $(v,\bar v)$ indicates that one must keep track of $\bar U=(U^H)^T$ as well as $U$.

\subsection{The Quaternionic Singular Value Decomposition}
\label{subsec:cat-quaternion-svd}
\emph{Quiver.} As for the real SVD, but with two quaternionic vertices with $\epsilondiv=\epsilonform=-1$ (Figure~\ref{fig:dsq-examples}(b)).

\emph{Algebra.} $S=\HH\times\HH$ with quaternion conjugation, and $R=\HH a\oplus\HH a^*$ with $qa=aq$ for $q\in\HH$ and $(qa)^*=a^*\bar q$.

\emph{Encoding.} $\Herm_n(R)=\{aW+a^*W^*:W\in M_n(\HH)\}$, where $W^*$ is the quaternionic conjugate transpose. The unitary group of $M_n(\HH)$ with this involution is the compact symplectic group $\Sp(n)$, so the action is $W\mapsto RWP^*$ with $P,R\in\Sp(n)$.

\emph{Classical analogue.} The singular value decomposition of quaternionic matrices~\cite{Zhang1997}. That quaternions carry, by default, a \emph{symplectic} Hermitian form is visible in the quiver: after complexification, quaternion conjugation becomes the symplectic involution of $M_2(\CC)$, so $\epsilonform=-1$. In the block model it is recorded by the matrix $G_v=J$ at the vertex (Table~\ref{tab:vertex-blocks}).

\subsection{Hermitian Eigendecompositions over \texorpdfstring{$\CC$}{C} and \texorpdfstring{$\HH$}{H}}
\label{subsec:cat-hermitian}
\emph{Quiver.} One vertex of unitary type (complex) or of type $(\HH,-1)$ (quaternionic), with a loop $x$. In the complex case, the transpose of $x$ is its conjugate partner, and after complexification there are two vertices $v,\bar v$ with one loop each, exchanged by both involutions.

\emph{Algebra.} $S=\CC$ with conjugation (respectively $\HH$ with quaternion conjugation), and $R=\CC x$ (respectively $\HH x$) with $zx=xz$ and $x^*=x$.

\emph{Encoding.} $(wx)^*=x\bar w=\bar wx$, so $xW$ is Hermitian if and only if $W^H=W$ (respectively $W^*=W$). The action is $W\mapsto PWP^H$ with $P\in\Un(n)$ (respectively $P\in\Sp(n)$). In the complex case the orbit of arrows is $\{x,\bar x\}$, with $x^{\top}=\bar x\neq x$. Both arrows are loops at the same orbit, so this is case~(b) of Theorem~\ref{thm:encoding}(3), in which the Hermitian condition is a symmetry condition on a single matrix.

\emph{Classical analogue.} The spectral theorem for complex Hermitian and for quaternionic Hermitian matrices~\cite{HornJohnson2013,Zhang1997}.
\subsection{The Takagi Factorisation}
\label{subsec:cat-takagi}
\emph{Quiver.} One complex vertex of unitary type, with an arrow $\alpha$ from $v$ to its conjugate $\bar v=v^{\top}$ that is fixed by the transpose with sign $+1$ (Figure~\ref{fig:dsq-examples}(d)).

\emph{Algebra.} $S=\CC$ with complex conjugation, and $R=\CC a$ with the \emph{twisted} bimodule structure $za=a\bar z$, and $a^*=a$. The involution is consistent with this rule: $(za)^*=a^*\bar z=a\bar z=za$.

\emph{Encoding.} For $H=aW$, with entries $w_{ij}a$, we get $(H^*)_{ij}=(w_{ji}a)^*=w_{ji}a$, so $H$ is Hermitian if and only if $W^T=W$, that is, $W$ is complex symmetric. For $U_0=P$ with $P\in\Un(n)$,
\[ (U_0HU_0^*)_{ij}=\sum_{k,l}p_{ik}\,w_{kl}\,a\,\bar p_{jl}=\sum_{k,l}p_{ik}\,w_{kl}\,p_{jl}\,a=(PWP^T)_{ij}\,a, \]
using $a\bar p=pa$. So the action is \emph{unitary congruence} $W\mapsto PWP^T$.

\emph{Classical analogue.} The Takagi factorisation: every complex symmetric matrix is $W=P\Sigma P^T$ with $P$ unitary and $\Sigma$ diagonal with non-negative entries~\cite{Takagi1925},~\cite[Corollary~4.4.4]{HornJohnson2013}. Comparing with the complex Hermitian eigendecomposition of Section~\ref{subsec:cat-hermitian}, the only difference in the quiver is that the arrow joins $v$ to its conjugate $\bar v$ instead of being a loop at $v$; algebraically, this is the difference between the untwisted rule $zx=xz$ and the twisted rule $za=a\bar z$. If instead $a^*=-a$ (sign $-1$), $W$ is complex skew-symmetric and the normal form under unitary congruence is $\bigoplus_j\sigma_j\bigl(\begin{smallmatrix}0&1\\-1&0\end{smallmatrix}\bigr)\oplus0$~\cite{Youla1961}.

\subsection{Congruence of Symmetric Matrices and the \texorpdfstring{$LDL^T$}{LDLT} Decomposition}
\label{subsec:cat-ldl}
\emph{Quiver.} A pair of real vertices $u,u^{\top}$ exchanged by the transpose (type $(\RR,0)$, so that the vertex orbit is $\twoR$), with an arrow $a:u\to u^{\top}$ fixed by the transpose with sign $+1$ (Figure~\ref{fig:signed-quiver-example}(c)).

\emph{Algebra.} $S=\twoR$ with idempotents $e=(1,0)$, $e^*=(0,1)$, and $R=\RR a$ with $a=e^*ae$ and $a^*=a$.

\emph{Encoding.} $aW$ is Hermitian if and only if $W^T=W$. By Section~\ref{subsec:2R}, $\Un_n(\twoR)=\{eP+e^*P^{-T}:P\in\GL_n(\RR)\}$, and a computation as in Section~\ref{sec:unification-example} gives $U_0(aW)U_0^*=a\,(P^{-T}WP^{-1})$. Writing $g=P^{-T}$, the action is congruence $W\mapsto gWg^T$ by arbitrary invertible $g$.

\emph{Classical analogue.} Sylvester's law of inertia: $W$ is congruent to $\diag(I_p,-I_q,0)$, and $(p,q)$ is a complete invariant. Numerically, a congruence normal form is computed from the symmetric indefinite factorisation $\Pi W\Pi^T=LDL^T$, in which $\Pi$ is a permutation chosen by symmetric pivoting, $L$ is unit lower triangular and $D$ is block diagonal with $1\times1$ and $2\times2$ blocks (Bunch--Kaufman pivoting)~\cite{GolubVanLoan2013}; a final congruence diagonalises $D$. With sign $-1$ one obtains congruence of skew-symmetric matrices, whose normal form is a symplectic (Darboux) basis.

\emph{Completing the cycle.} In this example the path algebra is acyclic. Every column $x$ with entries in $R$ satisfies $x^*x=0$ already in $\RR Q$, because $a^*a=aa$ is not a path; so the algorithms that need norms do not apply, and the generic algorithm for this quiver is the $LDL^*$ factorisation, which needs none (Section~\ref{sec:encoded-algorithms}). The stabiliser of the normal form is large (an indefinite orthogonal group), and correspondingly there are many normal-form algorithms: $g$ may be chosen lower triangular, as in $LDL^T$, or composed of an orthogonal eigenvector matrix and a diagonal scaling. Adding a reverse arrow $b:u^{\top}\to u$, also fixed with sign $+1$, removes this ambiguity. A Hermitian matrix with radical entries is then $aW+bZ$ with $W,Z$ symmetric, and $\GL_n(\RR)$ acts by $(W,Z)\mapsto(gWg^T,g^{-T}Zg^{-1})$. If we fix $Z=\varepsilon I$ with $\varepsilon>0$, for instance a formal infinitesimal, then the admissible $g$ are exactly those with $g^Tg=I$, and the problem becomes the orthogonal eigendecomposition of $W$. In Chapter~\ref{division-problem} the reverse arrow reappears as the inverse $a^{-1}$ in the localisation of the path algebra (Section~\ref{sec:localisation}).

\subsection{Similarity of Real Matrices}
\label{subsec:cat-similarity}
\emph{Quiver.} A pair of real vertices $v,v^{\top}$ exchanged by the transpose, with a loop $x$ at $v$ and its partner $x^{\top}$ at $v^{\top}$.

\emph{Algebra.} $S=\twoR$, $R=\RR x\oplus\RR x^*$ with $x=exe$.

\emph{Encoding.} $\Herm_n(R)=\{xX+x^*X^T:X\in M_n(\RR)\}$ and, as in Section~\ref{subsec:2R}, the action of $U_0=eP+e^*P^{-T}$ is $X\mapsto PXP^{-1}$.

\emph{Classical analogue.} Similarity of real matrices, whose normal form is the real Jordan form. As a Hermitian eigendecomposition, it is a direct sum of blocks $xJ+x^*J^T$, one for each real Jordan block $J$, none of which is unitarily similar to a direct sum of smaller blocks; each can be taken tridiagonal, since every real matrix that is similar to a single real Jordan block is similar to a tridiagonal one. This is the problem of Section~\ref{subsec:2R} again, now in degree~$1$; there it arose in degree~$0$, from the semisimple algebra $\twoR$ itself. The underlying quiver (a loop) is tame (Section~\ref{sec:tits_form}): the normal form has continuous parameters, the eigenvalues, and a discrete part, the Jordan structure, which depends discontinuously on the matrix (Section~\ref{sec:limitations}).

\begin{table}[htbp]
\centering
\caption{Doubly-signed quivers, the groups and matrix spaces that they encode by Theorem~\ref{thm:encoding}, and the corresponding classical decompositions. Vertex labels give the division algebra and $\epsilonform$; ``u'' denotes a complex vertex of unitary type, and $\rightleftarrows$ a pair of arrows exchanged by the transpose.}
\label{tab:catalogue}
\small
\renewcommand{\arraystretch}{1.25}
\begin{tabular}{@{}>{\raggedright\arraybackslash}p{0.22\linewidth}>{\raggedright\arraybackslash}p{0.16\linewidth}>{\raggedright\arraybackslash}p{0.25\linewidth}>{\raggedright\arraybackslash}p{0.29\linewidth}@{}}
\toprule
\textbf{Quiver} & \textbf{Group} & \textbf{Matrices $W$} & \textbf{Classical decomposition} \\
\midrule
$(\RR,+)$, loop $(+)$ & $\Or(n)$ & real symmetric & spectral theorem \\
$(\RR,+)$, loop $(-)$ & $\Or(n)$ & real skew-symmetric & real normal form, $2\times2$ blocks \\
\rowrule
$(\RR,+)\rightleftarrows(\RR,+)$ & $\Or(n)\times\Or(n)$ & real rectangular & real SVD \\
$(\CC,\mathrm{u})\rightleftarrows(\CC,\mathrm{u})$ & $\Un(n)\times\Un(n)$ & complex rectangular & complex SVD \\
$(\HH,-)\rightleftarrows(\HH,-)$ & $\Sp(n)\times\Sp(n)$ & quaternionic rectangular & quaternionic SVD \\
\rowrule
$(\CC,\mathrm{u})$, loop & $\Un(n)$ & complex Hermitian & Hermitian spectral theorem \\
$(\HH,-)$, loop $(+)$ & $\Sp(n)$ & quaternionic Hermitian & quaternionic spectral theorem \\
\rowrule
$(\CC,\mathrm{u})$, $v\to\bar v$ $(+)$ & $\Un(n)$ & complex symmetric & Takagi factorisation \\
$(\CC,\mathrm{u})$, $v\to\bar v$ $(-)$ & $\Un(n)$ & complex skew-symmetric & Youla normal form \\
\rowrule
$(\RR,0)$, $u\to u^{\top}$ $(+)$ & $\GL_n(\RR)$ & real symmetric & Sylvester inertia, $LDL^T$ \\
$(\RR,0)$, loops $x,x^{\top}$ & $\GL_n(\RR)$ & real square & real Jordan form \\
\bottomrule
\end{tabular}
\end{table}

\subsection{From a Matrix Problem to a Quiver}
\label{subsec:recipe}
Theorem~\ref{thm:encoding} can be read backwards, to encode a new matrix problem and run the generic algorithms on it. The problem must have the form that the theorem describes: a product of classical groups acting on a direct sum of matrix spaces, each of which is moved by one or two of the groups.
\begin{enumerate}
  \item For each group factor, choose the vertex type of Table~\ref{tab:vertex-types} whose unitary group it is, and give the corresponding orbit $c$ as many positions $n_c$ as the index $n$ of the group $\Un_n$ in that table, which is the size of the matrices over $K_c$ on which the factor acts; for an indefinite group, choose a form $\eta_c$ of the right signature (Remark~\ref{rem:indefinite}).
  \item For each matrix space that is moved by two factors, as $W\mapsto g'Wg^{-1}$, add an orbit of arrows between the two orbits; for each space that is moved by one factor, add an orbit of arrows within its orbit (Theorem~\ref{thm:encoding}(3)). Table~\ref{tab:catalogue} shows which arrows give which actions. For example, an arrow from a vertex to its transpose gives congruence, of symmetric or skew-symmetric matrices according to the sign of the arrow, and a loop at a vertex of exchange type gives similarity.
  \item Encode the data as a Hermitian matrix of degree~$1$, indexed by the positions in a cyclic order (Section~\ref{sec:indexing}); its entries between the positions of two orbits are the combinations of arrows of Theorem~\ref{thm:encoding}(3).
  \item Run the generic algorithm over the completed localisation $\Lambda$ of Chapter~\ref{division-problem}, as in Listing~\ref{lst:session}.
  \item Check that the unitary factor and the result contain no inverse letters. Then the identities that certify the result, such as $H=UTU^*$ and $U^*U=I$, hold over the original algebra (Theorem~\ref{thm:reduction}), and the coefficients of the arrows in the result are the classical answer. By Proposition~\ref{prop:householder-graded}, this is automatic for Householder reduction when its steps are defined and every component of the complexified quiver is an oriented cycle or a linearly oriented path.
\end{enumerate}
The generic algorithms report a breakdown when a column is isotropic, as for the congruence quiver (Section~\ref{subsec:cat-ldl}), or when a denominator is not invertible in the corner of its column. Example~\ref{ex:golub-kahan-numbers} carries out steps (3) to (5) by hand for the real singular value decomposition.

\section{Generic Algorithms on Encoded Matrices}
\label{sec:encoded-algorithms}
The generic algorithms of this thesis are run in the completed localisation $\Lambda$ of the complexified path $*$-algebra (Chapter~\ref{division-problem}), on matrices indexed by positions. A Hermitian matrix of degree $1$ over $A$ is lifted to $\Lambda$ by sending each arrow to itself. If an algorithm transforms it only by unitary matrices of degree $0$, and returns a matrix of degree $1$, then by Theorem~\ref{thm:encoding} it has stayed within the classical problem encoded by $Q$, and by Theorem~\ref{thm:reduction} its results reduce to $A$. This section gives a condition on $Q$ under which Householder reduction behaves in this way, and it identifies one step of the QR algorithm for several quivers of the catalogue. Table~\ref{tab:algorithms} summarises the results. They describe what the generic algorithms compute; they do not claim that these algorithms are better in practice than the specialised classical ones, such as the bidiagonalisation of Golub and Kahan~\cite{GolubKahan1965} as implemented in LAPACK, which are faster and better studied.

\paragraph{What this section uses from Chapter~\ref{division-problem}.} The \emph{localisation} $kQ_{\mathrm{loc}}$ adjoins to the path algebra an inverse $\alpha^{-1}$, from $e_w$ to $e_v$, of every arrow $\alpha:v\to w$ (Definition~\ref{def:localisation}). Its elements are linear combinations of \emph{reduced words}: products of letters $\alpha$ and $\alpha^{-1}$ with matching end points and no factor $\alpha\alpha^{-1}$ or $\alpha^{-1}\alpha$. The \emph{degree} of a word is the number of its letters $\alpha$ minus the number of its letters $\alpha^{-1}$, so that a path has its length as degree. The \emph{completed localisation} $\Lambda$ consists of the series $\sum_{d\geq d_0}x_d$ in which each $x_d$ is a finite combination of reduced words of degree $d$ (Definition~\ref{def:completion}); for one vertex with one loop $t$, it is the field of Laurent series $k((t))$. An element whose part of lowest degree is invertible is itself invertible (Lemma~\ref{lem:lambda-units}). When the underlying graph of every component of the quiver has at most one cycle, and that cycle has non-zero degree, every corner of $\Lambda$ is a space of matrices over a field (Proposition~\ref{prop:lambda-complete}); the degree of a cycle is the degree of the reduced word that goes once around it, up to sign. $*$-square roots are computed by the graded rule of Section~\ref{sec:lambda-sqrt}, which, for a homogeneous element of degree $2d$, looks for a root of degree $d$. A result holds in the original algebra when the computed factors contain no inverse letters (Theorem~\ref{thm:reduction}).

\begin{lemma}\label{lem:degree-zero}
Let $Q$ be a connected quiver with $n$ vertices and cycle rank $r$ (Section~\ref{sec:tits_form}), fix a vertex $o$, and for each vertex $v$ let $\ell(v)$ be the degree of a reduced word from $o$ to $v$. The following are equivalent:
\begin{enumerate}
  \item every reduced word of degree $0$ is trivial;
  \item either $r=0$ and the $\ell(v)$ are pairwise distinct, or $r=1$, the cycle of $Q$ has non-zero degree $\delta$, and the $\ell(v)$ are pairwise distinct modulo $\delta$;
  \item $Q$ is linearly oriented of type $A_n$, $1\to2\to\cdots\to n$, or an oriented cycle; for $n=1$, a vertex or a loop.
\end{enumerate}
In that case, the homogeneous elements of degree $0$ of $\Lambda$ are the linear combinations of the vertex idempotents, with coefficient blocks.
\end{lemma}
\begin{proof}
By Proposition~\ref{prop:localisation} and the proof of Proposition~\ref{prop:lambda-complete}, the reduced words from $v$ to $w$ correspond to the elements $g\in F_r$, and the word of $g$ has degree $\ell(w)-\ell(v)+\varphi(g)$ for a homomorphism $\varphi:F_r\to\ZZ$; the trivial word $e_v$ corresponds to $w=v$ and $g=1$. (Another choice of the words that define $\ell$ changes $\ell(v)$ by an element of $\varphi(F_r)$.) Hence (1) holds if and only if $\varphi$ is injective and $\ell(w)-\ell(v)\notin\varphi(F_r)$ for all $v\neq w$. As in the proof of Proposition~\ref{prop:lambda-complete}, $\varphi$ is injective if and only if $r=0$, or $r=1$ and the cycle has non-zero degree $\delta$; then $\varphi(F_r)$ is $0$ or $\delta\ZZ$. This proves that (1) and (2) are equivalent.

Along an arrow $\alpha:v\to w$, $\ell(w)=\ell(v)+1$, modulo $\delta$ if $r=1$. Assume (2). If $r=0$, $Q$ is a tree on which $\ell$ is injective. A vertex with three neighbours would give two of them the same value of $\ell$, so $Q$ is a path, and $\ell$ is monotone along it, so all its arrows point in the same direction. If $r=1$, the values $\ell(v)$ are $n$ distinct residues modulo $\delta$, so $n\leq|\delta|$. But $|\delta|$ is at most the length of the cycle, which is at most $n$. Hence the cycle passes through all the vertices, so that $Q$ is the cycle, and all its arrows point in the same direction. Conversely, for $1\to\cdots\to n$ we may take $\ell(v)=v-1$, and for an oriented cycle with $n$ vertices $\delta=n$ and $\ell$ takes the values $0,\dots,n-1$. The last statement holds because a homogeneous element of $\Lambda$ is a finite linear combination of reduced words of its degree.
\end{proof}

Every connected component of the complexified quiver of every quiver of Table~\ref{tab:catalogue} satisfies the conditions of the lemma. These components are single vertices with one loop, oriented $2$-cycles, for the singular value decompositions and for the Takagi and Youla factorisations, and, for the congruence quiver, a single arrow. The conditions are restrictive: they fail, for example, for the quiver of Example~\ref{ex:elimination-fails}, in which $a^{-1}b$ is a non-trivial word of degree $0$.

\begin{proposition}\label{prop:householder-graded}
Suppose that every connected component of the complexified quiver of $Q$ satisfies the conditions of Lemma~\ref{lem:degree-zero}, and that every form $\eta_c$ is positive definite, so that the involution of $E$ is the conjugate transpose. Let $H\in\Herm(E_1)$, and run Householder tridiagonalisation on $H$ over $\Lambda$, with a square-root rule that returns a real homogeneous root of degree $1$ for a real self-adjoint homogeneous element of degree $2$, such as the graded rule of Section~\ref{sec:lambda-sqrt}. Suppose that every step is defined: every root lies in the corner required by Proposition~\ref{prop:householder-any-root}, and every denominator is invertible in the corner of its column. Then every reflector lies in $\Un(E_0)$, and $H=UTU^*$ with $U\in\Un(E_0)$ and $T\in\Herm(E_1)$ tridiagonal. These identities hold over every homogeneous quotient $A$ of $\RR Q$, so $T$ lies in the orbit of $H$ under the group $\Un(E_0)$ of Theorem~\ref{thm:encoding}.
\end{proposition}
\begin{proof}
Suppose that the current matrix is Hermitian of degree $1$, and consider a step with column $x$, whose entries have degree $1$. Then $x^*x$ is homogeneous of degree $2$, the root $\gamma$ and its negative are homogeneous of degree $1$, so is $w=x-\hat e_1\gamma$, and $s=w^*x$ is homogeneous of degree $2$. The inverse of $s$ from $f$ to $f$ is homogeneous of degree $-2$: if $y$ is the inverse, then comparing the parts of degree $0$ in $ys=f$ and $sy=f$ shows that the part of $y$ of degree $-2$ is also an inverse, and inverses are unique. Hence $P=I-ws^{-1}w^*$ is homogeneous of degree $0$. It is real, because all the quantities are, so by Lemma~\ref{lem:degree-zero} its entries are real linear combinations of vertex idempotents with coefficient blocks, that is, $P\in E_0$; and it is unitary by Proposition~\ref{prop:householder-any-root}. So the next matrix $PHP^*$ is again Hermitian of degree $1$, and the first two claims follow by induction. Since $U$ and $T$ contain no inverse letters, Theorem~\ref{thm:reduction} gives the identities over $A$.
\end{proof}

For the real singular value decomposition, Proposition~\ref{prop:golub-kahan} identifies the reflectors with those of the bidiagonalisation of Golub and Kahan. For the symmetric and skew-symmetric loops, the same computation, with the loop in place of $a$ and $b$, shows that they are the classical Householder reflectors of $W$; for instance, for the skew loop a column $hy$ has $(hy)^*(hy)=-|y|^2h^2$, with the graded root $|y|\,h$. For the other quivers of the catalogue except the congruence quiver, every step was defined in the experiments of Chapter~\ref{chap:implementation}, with $U$ of degree $0$ and $T$ tridiagonal (Table~\ref{tab:experiments}). Over the congruence quiver every column with radical entries is isotropic (Section~\ref{subsec:cat-ldl}), and Householder reduction stops at the first non-zero column.

The QR algorithm can also be identified. For an invertible matrix $W$ over $D\in\{\RR,\CC,\HH\}$, let $W=Q_WR_W$ be its QR factorisation, with $R_W$ upper triangular with positive real diagonal. Both this factorisation and the Cholesky factorisation $G=LL^*$ of a positive definite Hermitian matrix, with $L$ lower triangular with positive real diagonal, are unique. One step of the \emph{Cholesky LR algorithm}, a variant of Rutishauser's LR algorithm~\cite{Rutishauser1958}, maps $G=LL^*$ to $L^*L$.

\begin{proposition}\label{prop:qr-steps}
Run one step $H=QR\mapsto H'=RQ$ of the QR algorithm over $\Lambda$, with $Q$ and $R$ computed by the Gram--Schmidt process with graded square roots whose coefficients are principal roots, such as $\gamma=|\hat y|\,a$ below.
\begin{enumerate}
  \item For the symmetric and skew-symmetric loops and the Hermitian loops over $\CC$ and $\HH$, let $H=xW$ with $W$ invertible, where $x$ is the loop. Then $H'=x\,R_WQ_W$: the step is the classical unshifted QR step on $W$.
  \item For the singular value decomposition over $D$, with the positions in cyclic order, let $H$ encode an invertible square $W$. Then $Q$ is homogeneous of degree $0$ and exchanges the orbits $u$ and $v$, and $H'$ encodes $W_1=R_WQ_{W^*}$, with the roles of $u$ and $v$ exchanged. Moreover $W_1^*W_1=L^*L$ and $W_1W_1^*=L'^*L'$, where $WW^*=LL^*$ and $W^*W=L'L'^*$: the step is one step of the Cholesky LR algorithm on each Gram matrix. Two steps give $W_2$ with $W_2^*W_2=R_1Q_1$ and $W_2W_2^*=R_2Q_2$, where $W^*W=Q_1R_1$ and $WW^*=Q_2R_2$: one step of the unshifted QR algorithm on each Gram matrix.
  \item For the Takagi factorisation, let $H$ encode an invertible complex symmetric $W$. Then $H'$ encodes $W_1=R_W\overline{Q_W}=Q_W^*W\overline{Q_W}$, a unitary congruence of $W$, and $W_1^*W_1=L^*L$, where $WW^*=LL^*$. Two steps give $W_2$ with $W_2W_2^*=R_2Q_2$, where $WW^*=Q_2R_2$.
\end{enumerate}
\end{proposition}
\begin{proof}
We prove (2) for $D=\RR$; over $\CC$ and $\HH$ the transposes become conjugate transposes. The columns of $H$ at the positions $u_j$ and $v_i$ are supported on the positions of $v$ and of $u$ respectively, so they are orthogonal, and the Gram--Schmidt process treats the columns of $u$ and of $v$ separately. The column of $u_j$ is $a\,(We_j)$, with the entries $w_{ij}a$ at the positions $v_i$. For columns $ay$ and $ay'$ with real $y$ and $y'$ we have $(ay')^*(ay)=(y'^Ty)\,ba$. So a residual $a\hat y\neq0$ has the graded root $\gamma=|\hat y|\,a$, with $\gamma^{-1}=|\hat y|^{-1}a^{-1}$, and the normalised column is $(\hat y/|\hat y|)\,e_v$, because $aa^{-1}=e_v$. Hence the Gram--Schmidt process on the columns of $u$ is the classical one on the columns of $W$, with the residuals multiplied by $a$ and the normalised columns by $e_v$: the column of $Q$ at the position $u_j$ has the entries $(Q_W)_{ij}e_v$ at the positions $v_i$, and $R$ has the entries $(R_W)_{kj}a$ at the positions $(u_k,u_j)$. In the same way the columns of $v$ give the entries $(Q_{W^T})_{ji}e_u$ of $Q$ at the positions $(u_j,v_i)$ and the entries $(R_{W^T})_{ki}b$ of $R$ at the positions $(v_k,v_i)$. Thus $Q$ is homogeneous of degree $0$, and after the step the positions $u_k$ carry the idempotent $e_v$ and the positions $v_i$ carry $e_u$. The product $RQ$ has the entries $(R_WQ_{W^T})_{ki}a$ at the positions $(u_k,v_i)$, their adjoints at $(v_i,u_k)$, and no others, which proves the first claim. From $W^T=Q_{W^T}R_{W^T}$ we get $WQ_{W^T}=R_{W^T}^T$, so
\[ W_1^TW_1=Q_{W^T}^TW^TWQ_{W^T}=R_{W^T}R_{W^T}^T,\qquad WW^T=R_{W^T}^TR_{W^T}=LL^T\ \text{with}\ L=R_{W^T}^T . \]
Similarly $W_1W_1^T=R_WR_W^T$ and $W^TW=R_W^TR_W$. A second step applies this to $W_1$, so $W_2^TW_2$ is the result of two Cholesky LR steps on $W^TW$, and $W_2W_2^T$ of two on $WW^T$. Finally, let $G=LL^T$ be positive definite, $L^TL=L_1L_1^T$, and $G=Q_GR_G$. Then $G^2=(LL_1)(LL_1)^T$ and $G^2=G^TG=R_G^TR_G$, so $R_G=(LL_1)^T$ by the uniqueness of the Cholesky factorisation, and
\[ L_1^TL_1=(LL_1)^{-1}G\,(LL_1)=R_G^{-T}GR_G^T=R_GGR_G^{-1}=R_GQ_G , \]
where the third equality holds because $R_G^TR_G=G^2$ commutes with $G$.

(1) This is the same computation, with the loop in place of $a$ and $b$ and without the exchange of orbits. For the skew loop $h$, a residual $h\hat y$ has $(h\hat y)^*(h\hat y)=-|\hat y|^2h^2$ and the graded root $|\hat y|\,h$.

(3) The real elements of degrees $0$ and $1$ are $s(z)=ze_v+\bar ze_{\bar v}$ and $m(z)=z\alpha+\bar z\bar\alpha$ with $z\in\CC$, and $s(z)m(w)=m(\bar zw)$, $m(w)s(z)=m(wz)$, $m(w)^*=m(w)$ and $m(w)^*m(w')=\bar ww'\,\bar\alpha\alpha+w\bar w'\,\alpha\bar\alpha$. Hence the Gram--Schmidt process on the columns of $H=(m(w_{ij}))$ is the classical one on the columns of $W$. A residual $m(\hat y)$ has the graded root $|\hat y|(\alpha+\bar\alpha)$, whose inverse is $|\hat y|^{-1}(\alpha^{-1}+\bar\alpha^{-1})$, and the normalised column is $s(\overline{\hat y}/|\hat y|)$. So $Q=(s(\overline{(Q_W)_{ij}}))$, $R=(m((R_W)_{ij}))$, and $RQ=(m((R_W\overline{Q_W})_{ij}))$. Since $W=W^T=R_W^TQ_W^T$ and $Q_W^T\overline{Q_W}=I$, we get $WW^*=R_W^T\overline{R_W}=LL^*$ with $L=R_W^T$, and $W_1^*W_1=Q_W^T\,\overline WW\,\overline{Q_W}=\overline{R_W}R_W^T=L^*L$. The matrix $W_1$ is again symmetric, so $W_1W_1^*=\overline{W_1^*W_1}=\overline{L^*L}$. Applying the first statement to $W_1$ and conjugating, $W_2W_2^*$ is the result of two Cholesky LR steps on $WW^*$, and the last part of the proof of (2) applies.
\end{proof}

The last step of the proof of (2), that two steps of the Cholesky LR algorithm on a positive definite matrix are one step of the QR algorithm, is classical; see~\cite{Xu1998} and the references there. What the proposition adds is that a single step of the generic QR algorithm, for the quiver of the singular value decomposition or of the Takagi factorisation, is a Cholesky LR step on the Gram matrices.

For the singular value decomposition, part~(2) can also be derived without $*$-algebras. Read as a real or complex matrix, with the positions interleaved, $H$ is the Hermitian matrix $\bigl(\begin{smallmatrix}0&W^*\\W&0\end{smallmatrix}\bigr)$ with its rows and columns interleaved, the generic step is the classical QR step on it, and its square is $\diag(W^*W,WW^*)$, interleaved in the same way; and if two unshifted QR steps take a matrix $G$ to $G_2$, then one QR step takes $G^2$ to $G_2^2$.

Part~(2) also explains the role of the double steps in the experiments of Chapter~\ref{chap:implementation}. A single step exchanges the two orbits: $QQ^*$ is the identity, but $Q^*Q$ is the identity of the swapped indexing, so the step is a unitary transformation from one indexing to the other rather than a unitary similarity over $A$. Only a double step returns to the original positions, where it performs one step of the unshifted QR algorithm on $W^*W$ and on $WW^*$, and so converges to the singular values.

The split-complex numbers $\twoR$ give a different relation between the QR and LR algorithms, over a different $*$-algebra. In degree $0$, the algebra of the real singular value decomposition and $\twoR$ are both $\RR\times\RR$, but the involution of the former fixes the two idempotents, so that its unitary group is $\Or(n)\times\Or(m)$, whereas that of $\twoR$ exchanges them, so that its unitary group is $\GL_n(\RR)$. Over $\twoR$, the algorithms become elimination algorithms, and the QR algorithm becomes the LR algorithm on an invertible real matrix, whenever the steps are defined, rather than on a positive definite Gram matrix. We state this for $\twoR$ itself, which encodes similarity of real matrices in degree $0$ (Section~\ref{subsec:2R}).

\begin{proposition}\label{prop:2R-algorithms}
Let $X\in\GL_n(\RR)$, let $H=(X,X^T)$ be the corresponding Hermitian matrix over $\twoR$, and use the $*$-square roots $(c,1)$ of the self-adjoint elements $(c,c)$.
\begin{enumerate}
  \item The Cholesky factorisation $H=R^*R$, with $R=(R_1,R_2)$ upper triangular, is defined if and only if the leading principal minors of $X$ are non-zero, and then $X=R_2^TR_1$ is the Doolittle factorisation of $X$, with $R_2^T$ unit lower triangular.
  \item One step $H=QR\mapsto RQ$ of the QR algorithm, with the Gram--Schmidt process, is defined if and only if the leading principal minors of $X^2$ are non-zero. It then maps $X$ to $UXU^{-1}$, where $X^2=LU$ is the Doolittle factorisation. If the two steps $X=L_1U_1\mapsto U_1L_1=L_2U_2\mapsto U_2L_2$ of the LR algorithm are defined, then $U_2L_2=UXU^{-1}$.
\end{enumerate}
With the roots $(1,c)$ instead, $R$ changes by a diagonal factor: (1) gives the Crout factorisation, and the result of (2) changes by a diagonal similarity. Statement~(2), with this remark, and the tridiagonalisation by similarity of Section~\ref{subsec:2R} carry over to the similarity quiver of Section~\ref{subsec:cat-similarity}, with $H=xX+x^{\top}X^T$. Statement~(1) carries over only with roots that are not homogeneous: the pivots $c\,(x+x^{\top})$ of this $H$ have odd degree, so the graded rule finds no root, but $c\,x+e_{v^{\top}}$ is a root of $c\,(x+x^{\top})$, and with these roots $H=R^*R$ with $R=xU+e_{v^{\top}}L^T$, where $X=LU$ is the Doolittle factorisation.
\end{proposition}
\begin{proof}
(1) For $R=(R_1,R_2)$ we have $R^*=(R_2^T,R_1^T)$, so $R^*R=(R_2^TR_1,R_1^TR_2)$, and $H=R^*R$ says that $X=R_2^TR_1$. The factorisation takes as diagonal entries of $R$ the roots $(c,1)$ of the pivots $(c,c)$, so $R_2$ has unit diagonal. A pivot can be inverted if and only if $c\neq0$, and the pivots are the ratios of consecutive leading principal minors of $X$.

(2) Since $H$ is invertible, no residual of the Gram--Schmidt process vanishes, and the process is defined if and only if the Cholesky factorisation of $H^*H$ is, with the same roots. Its factor is $R$, because $R^*R=H^*H$. Now $H^*H=H^2=(X^2,(X^2)^T)$, so by (1) the step is defined if and only if the leading principal minors of $X^2$ are non-zero, and then $X^2=R_2^TR_1=LU$. The first component of $RQ=RHR^{-1}$ is $R_1XR_1^{-1}=UXU^{-1}$. If the two LR steps are defined, then $X^2=L_1U_1L_1U_1=L_1L_2U_2U_1$, so $L=L_1L_2$ and $U=U_2U_1$ by the uniqueness of the Doolittle factorisation, and $U_2L_2=U_2U_1XU_1^{-1}U_2^{-1}=UXU^{-1}$.

With the roots $(1,c)$, $R_1$ has unit diagonal instead, and $R$ is multiplied by a diagonal matrix. For the similarity quiver, the real elements $p\,x^k+q\,(x^{\top})^k$ with $p,q\in\RR$ multiply and are adjoint to each other like the elements $(p,q)\,t^k$ of $\twoR((t))$, because $xx^{\top}=x^{\top}x=0$ and $x^*=x^{\top}$. In the QR step and in Householder reduction, the elements whose roots are taken, $v^*v$ for columns $v$ of degree $1$, have degree $2$, and their graded roots have degree $1$. So every computation of the QR step and of Householder reduction is the one over $\twoR$, with the quantities of degree $k$ multiplied by the corresponding powers of $x$ and $x^{\top}$.
\end{proof}

The proposition makes precise the remark of Section~\ref{subsec:2R} that the generic algorithm can break down although the classical problem has a solution: the QR step breaks down exactly at the isotropic vectors, that is, when a leading principal minor of $X^2$ vanishes. It can also succeed where the LR algorithm fails: for $X=\bigl(\begin{smallmatrix}0&1\\1&0\end{smallmatrix}\bigr)$ the LR algorithm stops at once, whereas $X^2=I$, and the QR step is defined.

\begin{table}[htbp]
\centering
\caption{What the generic algorithms compute on the Hermitian matrices of degree $1$ encoded by the quivers of Table~\ref{tab:catalogue}, and over $\twoR$. The statements are proved in Propositions~\ref{prop:golub-kahan}, \ref{prop:householder-graded}, \ref{prop:qr-steps} and~\ref{prop:2R-algorithms} under their hypotheses: every step defined, $W$ square and invertible for the QR steps, and the LR steps defined; for the real loops under Householder reduction, see the remark after Proposition~\ref{prop:householder-graded}. That every step is defined, the bidiagonalisation over $\CC$ and $\HH$, Lanczos and $LDL^*$ were observed in the experiments of Table~\ref{tab:experiments}.}
\label{tab:algorithms}
\small
\renewcommand{\arraystretch}{1.25}
\begin{tabular}{@{}>{\raggedright\arraybackslash}p{0.17\linewidth}>{\raggedright\arraybackslash}p{0.22\linewidth}>{\raggedright\arraybackslash}p{0.28\linewidth}>{\raggedright\arraybackslash}p{0.25\linewidth}@{}}
\toprule
\textbf{Quiver} & \textbf{Householder reduction} & \textbf{QR algorithm} & \textbf{Other} \\
\midrule
symmetric and skew loops & classical Householder tridiagonalisation of $W$ & one step: the classical QR step on $W$ & \\
\rowrule
Hermitian loops, $\CC$ and $\HH$ & tridiagonalisation by unitary similarity & one step: the classical QR step on $W$ & \\
\rowrule
SVD over $\RR$, $\CC$, $\HH$ & bidiagonalisation; over $\RR$, that of Golub and Kahan & two steps: a QR step on $W^*W$ and on $WW^*$; one step: a Cholesky LR step on each & Lanczos: the \mbox{Golub--Kahan--Lanczos} process \\
\rowrule
Takagi, Youla & tridiagonalisation by unitary congruence & Takagi, two steps: a QR step on $WW^*$; one step: $W\mapsto Q_W^*W\overline{Q_W}$, a Cholesky LR step on $WW^*$ & \\
\rowrule
congruence & stops: $x^*x=0$ & stops: $x^*x=0$ & $LDL^*$: $LDL^T$ with the inertia of $W$ \\
\rowrule
similarity, and $\twoR$ & tridiagonalisation by similarity & one step: $X\mapsto UXU^{-1}$ with $X^2=LU$, which is two LR steps & Cholesky over $\twoR$: Doolittle $LU$ \\
\bottomrule
\end{tabular}
\end{table}

\section{Beyond Degree One: Perturbation Theory}
\label{sec:beyond-radical}
Theorem~\ref{thm:encoding} concerns Hermitian matrices of degree $1$. For a general Hermitian matrix $H=H_0+H_1$ over $B=B_Q$, with $H_0\in\Herm_n(S)$ and $H_1\in\Herm_n(R)$, the radical part $U_1$ of a unitary matrix also acts, and the problem becomes one of first-order perturbation theory. We illustrate this with the dual numbers, where the conclusion is the spectral theorem of~\cite{Gutin2022DualNumbers}.

\begin{proposition}\label{prop:dual-spectral}
Let $B=\RR[\varepsilon]/(\varepsilon^2)$ with the trivial involution. Every $H=X+\varepsilon W\in\Herm_n(B)$ is unitarily similar to a diagonal matrix $\diag(\lambda_i+\varepsilon\mu_i)$ with $\lambda_i,\mu_i\in\RR$.
\end{proposition}
\begin{proof}
The unitary matrices are $U=P(I+\varepsilon K)$ with $P\in\Or(n)$ and $K^T=-K$, since $(I+\varepsilon K)^T(I+\varepsilon K)=I+\varepsilon(K^T+K)$. First choose $P$ with $P^TXP=\Lambda=\diag(\lambda_i)$, grouping equal eigenvalues together, and replace $W$ by $P^TWP$. Next, conjugating by $I+\varepsilon K$ gives
\[ (I-\varepsilon K)(\Lambda+\varepsilon W)(I+\varepsilon K)=\Lambda+\varepsilon\,(W+\Lambda K-K\Lambda), \]
and $(\Lambda K-K\Lambda)_{ij}=(\lambda_i-\lambda_j)k_{ij}$. Choosing $k_{ij}=-w_{ij}/(\lambda_i-\lambda_j)$ whenever $\lambda_i\neq\lambda_j$, which defines a skew-symmetric $K$ because $W$ is symmetric, removes all entries of $W$ connecting distinct eigenvalues. Finally, the remaining diagonal blocks of $W$, one for each eigenvalue of $X$, are diagonalised by an orthogonal matrix that is block diagonal and therefore commutes with $\Lambda$.
\end{proof}

The proof shows the general pattern, which the example of Section~\ref{sec:unification-example} also exhibits. The part of degree $0$, $H_0$, is first brought into normal form by the classical groups at the vertices. The radical part of $U$ then removes the components of $H_1$ that connect distinct ``eigenvalues'' of $H_0$, by solving Sylvester equations. What remains is a problem of degree $1$, governed by Theorem~\ref{thm:encoding}, for the stabiliser of the normal form of $H_0$. For the dual numbers this is the familiar statement that the first-order perturbations of a multiple eigenvalue are the eigenvalues of the compression of the perturbation to its eigenspace~\cite{Kato1995}.

Algorithmically, the proof requires knowing which eigenvalues of $H_0$ coincide. Iterative algorithms such as the Jacobi method do not know this in advance. When they meet two diagonal entries with the same part of degree $0$, they must divide the corresponding differences, which are zero divisors, as in the example after Example~\ref{ex:dual-division}. This is the division problem. Chapter~\ref{division-problem} solves it by computing in $\Lambda$. It can also be avoided, as in the proof above, by first computing the eigendecomposition of $H_0$, deciding which eigenvalues are equal, exactly or up to a tolerance, and then applying classical algorithms to the corresponding blocks of $H_1$.

\section{Representation Type and the Behaviour of Generic Algorithms}
\label{sec:wild-behaviour}
Theorem~\ref{thm:encoding} allows a count of degrees of freedom, analogous to the Tits form of Section~\ref{sec:tits_form}, for each decomposition. The number of continuous parameters of a normal form equals
\[ \dim(\text{matrix space})-\dim(\text{group})+\dim(\text{stabiliser of a generic element}). \]
For the real SVD of $n\times n$ matrices the generic stabiliser is finite, and the count is $n^2-2\cdot\tfrac{n(n-1)}{2}=n$, the number of singular values; for the Takagi factorisation it is $n(n+1)-n^2=n$, and for the symmetric eigendecomposition $\tfrac{n(n+1)}{2}-\tfrac{n(n-1)}{2}=n$. For the complex Hermitian eigendecomposition the stabiliser of a generic diagonal matrix is the torus $\Un(1)^n$, so the count is $n^2-n^2+n=n$. In all the examples of Table~\ref{tab:catalogue}, the number of parameters grows at most linearly in $n$.

Consider instead a vertex of type $(\RR,+1)$ with \emph{two} loops $x$ and $y$, both fixed with sign $+1$. Then $B=\RR[x,y]/(x,y)^2$ with $x^*=x$ and $y^*=y$, a Hermitian matrix of degree $1$ is $xW_1+yW_2$ with $W_1,W_2$ symmetric, and the group $\Or(n)$ acts by simultaneous orthogonal similarity. Now
\[ \dim\{(W_1,W_2)\}-\dim\Or(n)=n(n+1)-\tfrac{n(n-1)}{2}=\tfrac{n(n+3)}{2}, \]
which is a lower bound for the number of continuous parameters and grows quadratically in $n$. This quadratic growth is the hallmark of a wild problem, and the underlying quiver (one vertex with two loops) is indeed wild.

What does a generic algorithm do in such a case? A Jacobi-type method for the symmetric eigenproblem, applied to $xW_1+yW_2$, has to divide elements such as $xw_{1,ij}+yw_{2,ij}$ by elements such as $x\alpha+y\beta$. Unlike the quotients needed for the catalogue, these do not exist in the completed localisation $\Lambda$: if $\alpha,\beta\neq0$, then $x\alpha+y\beta$ is not invertible in $\Lambda$. Indeed, it is homogeneous of degree $1$, so an inverse would be homogeneous of degree $-1$, as in the proof of Proposition~\ref{prop:householder-graded}, and would therefore lie in $kQ_{\mathrm{loc}}=k[F_2]$ (Proposition~\ref{prop:localisation}). But the homomorphism from $k[F_2]$ to the ring of Laurent polynomials $k[s^{\pm1},t^{\pm1}]$ that abelianises $F_2$ maps $x\alpha+y\beta$ to $s\alpha+t\beta$, which is not a unit, since the units of that ring are the non-zero multiples of monomials. Here the condition of Proposition~\ref{prop:lambda-complete} fails, since $r=2$, and the elimination of Section~\ref{sec:lambda-inverses} stops, as the implementation confirms.

The quotients do exist in every skew field that contains the path algebra $k\langle x,y\rangle$, for example in Cohn's free field, and no such skew field is commutative (Section~\ref{sec:embedding}). Whether the Jacobi method can then be completed is a further question: its rotations also need $*$-square roots, of elements such as $1+t^*t$, and we have not shown that these exist in the free field. If they do, the output, a ``diagonalisation'', has entries in the skew field that need not lie in $\widehat{kQ}$, so Theorem~\ref{thm:reduction} gives no reduction to $B$, and the output lies outside the domain of the original decomposition. We have not run such an algorithm, and we do not attempt a general theorem about wild cases.

The following proposition shows that a vertex on two distinct simple oriented cycles, the situation in which a non-commutative skew field is unavoidable, always gives a wild quiver.

\begin{proposition}\label{prop:two-cycles-wild}
If some vertex of a connected quiver $Q$ lies on two distinct simple oriented cycles, then $Q$ is wild.
\end{proposition}
\begin{proof}
Two distinct simple oriented cycles have distinct sets of arrows, so the underlying graph of $Q$ contains two independent cycles, and its cycle rank is at least $2$. By Section~\ref{sec:tits_form}, $Q$ is wild.
\end{proof}

The converse fails: the Kronecker quiver with three arrows is acyclic and wild. Thus the quivers for which a commutative coefficient field suffices in Chapter~\ref{division-problem} include all tame quivers but also some wild ones. The completed localisation is a narrower tool. Its elimination is guaranteed to be complete only under the condition of Proposition~\ref{prop:lambda-complete}, and without it the elimination can stop on an invertible element, even for a tame quiver: the quiver of Example~\ref{ex:elimination-fails} is the Kronecker quiver with two arrows.
%======================================================================
%   CHAPTER 7: LOCALISATION AND COMPLETION OF PATH ALGEBRAS
%======================================================================
\chapter{Localisation (Adjoining Inverse Arrows) and Completion of Path Algebras}
\label{division-problem}

Chapters~\ref{chap:matrixops} and~\ref{chap:indexing} showed that generic algorithms must divide by zero divisors and normalise vectors whose norms vanish. This chapter constructs the algebra in which they do so. The encoded decompositions do not require it: for each quiver of the catalogue, the encoded problem is a classical matrix problem (Theorem~\ref{thm:encoding}), which classical algorithms solve directly. The algebra constructed here serves to run one program, written once, over the $*$-algebras of all these quivers.

Starting from a quotient $kQ/I$ of a path algebra, it is obtained in three steps: lift to the path algebra $kQ$ (Section~\ref{sec:division-problem}), invert every arrow between its end points (Section~\ref{sec:localisation}), and complete in the path degree (Section~\ref{sec:completion}). The result, the \emph{completed localisation} $\Lambda$, keeps the involution and the conjugation of a doubly-signed quiver. Sections~\ref{sec:lambda-inverses} and~\ref{sec:lambda-sqrt} show how to invert elements and take $*$-square roots in $\Lambda$, and Section~\ref{sec:reduction} shows when a result holds in $kQ/I$. Section~\ref{sec:embedding} looks beyond $\Lambda$, at embeddings of $kQ$ in simple Artinian rings. Throughout, $k$ is a field. For a doubly-signed quiver, $k=\CC$, and the constructions are applied to its complexified path $*$-algebra (Definition~\ref{def:path-star-algebra}).

\section{The Division Problem}
\label{sec:division-problem}
Let $Q$ be a quiver and $I$ an ideal of $kQ$ with $R_Q^m\subseteq I\subseteq R_Q^2$ for some $m$, so that $kQ/I$ is a finite-dimensional basic algebra with quiver $Q$ (Section~\ref{sec:quivers-basic}). By Theorem~\ref{thm:structure}, every finite-dimensional real $*$-algebra is a matrix algebra, indexed by positions (Section~\ref{sec:indexing}), over a quotient $\RR Q/I$ of the path $*$-algebra of its doubly-signed quiver, so it suffices to compute in such quotients. We describe the constructions for path algebras $kQ$; for a path $*$-algebra they are applied to its complexification, which is a block-matrix version of $\CC Q$ (see the remarks after Proposition~\ref{prop:localisation}). In $kQ/I$ every arrow is nilpotent, and no non-zero element of the radical is invertible between idempotents (Section~\ref{sec:division}). We therefore first \emph{lift} elements to the path algebra $kQ$, by a linear section of the quotient map; for instance, choose paths of length less than $m$ whose residues form a basis of $kQ/I$, and send each residue to its path. In $kQ$, products of arrows no longer vanish: the norm $x^*x$ of Section~\ref{sec:unification-example} becomes a non-zero element of degree $2$. Arrows are still not invertible in $kQ$, so we make each of them invertible between its own end points, in the sense of Definition~\ref{def:invertible-between}.

This chapter considers two ways to do this; Leavitt path algebras and graph $C^*$-algebras adjoin partial inverses of the arrows instead (Section~\ref{sec:related}). The first adjoins inverses of the arrows and completes the result (Sections~\ref{sec:localisation} and~\ref{sec:completion}). It is universal, it keeps the involution, its arithmetic is effective, and it is what the implementation of Chapter~\ref{chap:implementation} uses. When every connected component of $Q$ has at most one cycle and that cycle has non-zero degree, it already yields matrix rings over $k$ or over fields of Laurent series (Proposition~\ref{prop:lambda-complete}). The second embeds $kQ$ in a simple Artinian ring $M_n(D)$, in which every non-zero element between two vertices is invertible (Section~\ref{sec:embedding}); for general quivers this requires a non-commutative skew field $D$, and no algorithm of this thesis uses it.

\section{Localisation at the Arrows}
\label{sec:localisation}
\begin{definition}\label{def:localisation}
Let $Q$ be a quiver. The \emph{localisation of $kQ$ at the arrows} is the $k$-algebra $kQ_{\mathrm{loc}}$ obtained from $kQ$ by adjoining, for every arrow $\alpha:v\to w$, an element $\alpha^{-1}$ subject to the relations
\[ \alpha^{-1}=e_v\,\alpha^{-1}e_w,\qquad \alpha^{-1}\alpha=e_v,\qquad \alpha\alpha^{-1}=e_w . \]
\end{definition}

Thus every arrow $\alpha:v\to w$ becomes invertible from $e_v$ to $e_w$, with inverse $\alpha^{-1}$, and nothing else is inverted. By construction, $kQ_{\mathrm{loc}}$ has the following universal property: every homomorphism $\varphi$ from $kQ$ to a ring $R$ such that $\varphi(\alpha)$ is invertible from $\varphi(e_v)$ to $\varphi(e_w)$ for every arrow $\alpha:v\to w$ extends uniquely to a homomorphism $kQ_{\mathrm{loc}}\to R$. In the language of Cohn~\cite{Cohn2006}, $kQ_{\mathrm{loc}}$ is the universal localisation of $kQ$ at the homomorphisms of projective modules $kQe_w\to kQe_v$, $y\mapsto y\alpha$, given by the arrows. In the language of Gabriel and Zisman~\cite{GabrielZisman1967}, it is the algebra of the category obtained from the path category of $Q$ by inverting all arrows, that is, of the free groupoid on $Q$. Each arrow is inverted only between its own end points: if two arrows $\alpha,\beta$ end at the same vertex $w$, then $\alpha\alpha^{-1}=\beta\beta^{-1}=e_w$.

A \emph{word} is a product of letters $\alpha$ and $\alpha^{-1}$ ($\alpha\in Q_1$) whose end points match, or a trivial word $e_v$. It is \emph{reduced} if it contains no factor $\alpha\alpha^{-1}$ or $\alpha^{-1}\alpha$, and its \emph{degree} is the number of its letters $\alpha$ minus the number of its letters $\alpha^{-1}$. The paths are the reduced words without letters $\alpha^{-1}$, and their degree is their length. If $Q$ is connected and has cycle rank $1$, the \emph{degree of its cycle} is the degree of a reduced word that goes once around the cycle of its underlying graph; it is defined up to sign. For an oriented cycle with $n$ arrows it is $\pm n$, and for the Kronecker quiver $u\rightrightarrows u'$ with arrows $a$ and $b$ it is $0$, the degree of $a^{-1}b$.

\begin{proposition}\label{prop:localisation}
Let $Q$ be a connected quiver with $n$ vertices and cycle rank $r$ (Section~\ref{sec:tits_form}). Fix a vertex $o$ and a spanning tree $T$ of the underlying graph, and let $\tau_v$ be the reduced word along $T$ from $o$ to $v$. Let $F_r$ be the free group on symbols $\xi_\alpha$, one for each of the $r$ arrows $\alpha$ not in $T$.
\begin{enumerate}
  \item There is an isomorphism $\Phi:kQ_{\mathrm{loc}}\to M_n(k[F_r])$ with $\Phi(e_v)=E_{vv}$, $\Phi(\alpha)=E_{wv}$ for $\alpha:v\to w$ in $T$, and $\Phi(\alpha)=\xi_\alpha E_{wv}$ for $\alpha:v\to w$ not in $T$.
  \item The reduced words form a basis of $kQ_{\mathrm{loc}}$, and $\Phi$ maps a reduced word $\omega:v\to w$ to $gE_{wv}$, where $g\in F_r$ is obtained from $\omega$ by deleting the letters of $T$ and replacing $\alpha^{\pm1}$ by $\xi_\alpha^{\pm1}$. In particular, the product of two reduced words is the reduced form of their concatenation, or zero if their end points do not match, and $kQ\to kQ_{\mathrm{loc}}$ is injective.
\end{enumerate}
\end{proposition}
\begin{proof}
(1) The proposed images of the $e_v$ and of the arrows satisfy the relations of $kQ$, and each $\Phi(\alpha)$ is invertible from $E_{vv}$ to $E_{ww}$, with inverse $E_{vw}$ or $\xi_\alpha^{-1}E_{vw}$. By the universal properties of $kQ$ and of $kQ_{\mathrm{loc}}$, they define a homomorphism $\Phi$. Since $\tau_v$ consists of letters of $T$, $\Phi(\tau_v)=E_{vo}$. For an arrow $\alpha:v\to w$ not in $T$, the element $\tau_w^{-1}\alpha\tau_v$ of $e_okQ_{\mathrm{loc}}e_o$ is invertible there, with inverse $\tau_v^{-1}\alpha^{-1}\tau_w$. Hence $\xi_\alpha\mapsto\tau_w^{-1}\alpha\tau_v$ extends to a homomorphism from $F_r$ to the group of units of $e_okQ_{\mathrm{loc}}e_o$, and to an algebra homomorphism $\psi:k[F_r]\to e_okQ_{\mathrm{loc}}e_o$. Put $\Psi(gE_{wv})=\tau_w\psi(g)\tau_v^{-1}$. Since $\tau_v^{-1}\tau_{v'}=\delta_{vv'}e_o$, $\Psi$ is an algebra homomorphism. For an arrow $\alpha:v\to w$ in $T$ we have $\tau_w=\alpha\tau_v$ or $\tau_v=\alpha^{-1}\tau_w$, so $\Psi\Phi(\alpha)=\tau_w\tau_v^{-1}=\alpha$; for $\alpha$ not in $T$, $\Psi\Phi(\alpha)=\tau_w\tau_w^{-1}\alpha\tau_v\tau_v^{-1}=\alpha$. Conversely, $\Phi\psi(\xi_\alpha)=\xi_\alpha E_{oo}$, so $\Phi\Psi(gE_{wv})=E_{wo}\,gE_{oo}\,E_{ov}=gE_{wv}$. Hence $\Phi$ and $\Psi$ are inverse isomorphisms.

(2) Let $\omega:v\to w$ be a word. Letters of $T$ map to matrix units and the other letters to $\xi_\alpha^{\pm1}$ times matrix units, so $\Phi(\omega)=gE_{wv}$ with $g$ as described. If $\omega$ is reduced, then so is $g$: a factor $\xi_\alpha^{\varepsilon}\xi_\alpha^{-\varepsilon}$ of $g$ would come from a factor $\alpha^{\varepsilon}\tau\alpha^{-\varepsilon}$ of $\omega$ in which $\tau$ is a reduced word in the letters of $T$ from a vertex to itself; such a word is trivial because $T$ is a tree, so $\omega$ would not be reduced. Conversely, given $v$, $w$ and a reduced $g\in F_r$, there is exactly one reduced word $\omega:v\to w$ with $\Phi(\omega)=gE_{wv}$: its letters outside $T$ are prescribed by $g$, and between them, and at both ends, it must follow the unique reduced paths in $T$. Hence $\Phi$ maps the reduced words bijectively onto the basis $\{gE_{wv}\}$ of $M_n(k[F_r])$. The remaining statements follow, since paths are reduced words.
\end{proof}

For one vertex with one loop $x$, $kQ_{\mathrm{loc}}=k[x,x^{-1}]$. For the quiver of Section~\ref{sec:unification-example}, with $T=\{a\}$ and $o=u$, the cycle gives $\tau_u^{-1}b\tau_v=ba$, and $kQ_{\mathrm{loc}}\cong M_2(k[g,g^{-1}])$ with $g=ba$; this is the algebra $\Lambda_0$ used there. For an oriented cycle with $n$ vertices, $kQ_{\mathrm{loc}}\cong M_n(k[x,x^{-1}])$, and for one vertex with two loops, $kQ_{\mathrm{loc}}$ is the group algebra $k[F_2]$ of the free group.

For the path $*$-algebra of a doubly-signed quiver we apply the construction to the complexification $T_{S_\CC}(M_\CC)$ of Definition~\ref{def:path-star-algebra}. In $T_{S_\CC}(M_\CC)$ a path from $v$ to $w$ carries a single coefficient block in $M_{d_w\times d_v}(\CC)$, so that $T_{S_\CC}(M_\CC)$ is the algebra of block matrices with entries in $M_{d_w\times d_v}(\CC)\otimes_\CC e_w\CC Qe_v$. Replacing $\CC Q$ by $\CC Q_{\mathrm{loc}}$ gives its localisation: its elements are finite sums $\sum_\omega X_\omega\,\omega$ over reduced words $\omega:v\to w$ with $X_\omega\in M_{d_w\times d_v}(\CC)$, multiplied by $(X\omega)(Y\omega')=(XY)(\omega\omega')$ followed by cancellation. This is how the implementation stores elements (Chapter~\ref{chap:implementation}).

\begin{lemma}\label{lem:loc-involution}
Let $\theta$ be an involution of $kQ$, linear or conjugate-linear, that maps every vertex idempotent to a vertex idempotent and every arrow to a non-zero multiple of an arrow. Then $\theta$ extends uniquely to an involution of $kQ_{\mathrm{loc}}$, with $\theta(\alpha^{-1})=\theta(\alpha)^{-1}$. The same holds for an automorphism of order two, linear or conjugate-linear, that maps vertices to vertices and arrows to non-zero multiples of arrows.
\end{lemma}
\begin{proof}
Consider an anti-automorphism $\theta$; the other case is the same without passing to opposite rings. Regard $\theta$ as a homomorphism from $kQ$ to the opposite ring of $kQ_{\mathrm{loc}}$, which is conjugate-linear if $\theta$ is; the universal property applies equally to such homomorphisms. For an arrow $\alpha:v\to w$, $\theta(\alpha)=\lambda\beta$ with $\lambda\neq0$ and $\beta$ an arrow from $\theta(e_w)$ to $\theta(e_v)$, which is invertible in $kQ_{\mathrm{loc}}$ with inverse $\lambda^{-1}\beta^{-1}$. In the opposite ring, $\lambda\beta$ is therefore invertible from $\theta(e_v)$ to $\theta(e_w)$. By the universal property, $\theta$ extends uniquely to an anti-homomorphism $\tilde\theta$ of $kQ_{\mathrm{loc}}$, and $\tilde\theta(\alpha^{-1})=\theta(\alpha)^{-1}$. The composite $\tilde\theta^2$ is a homomorphism that extends the identity of $kQ$, so it is the identity by the uniqueness in the universal property.
\end{proof}

For the complexified path $*$-algebra of a doubly-signed quiver, with coefficient blocks, the involution and the conjugation are given on paths by the formulas~\eqref{eq:block-star} of the block model (Proposition~\ref{prop:block-model}). We extend them to reduced words by $(\alpha^{-1})^{\top}=(\alpha^{\top})^{-1}$, $\overline{\alpha^{-1}}=\bar\alpha^{-1}$, $\kappa(\alpha^{-1})=\kappa(\alpha)$ and $\mu(\alpha^{-1})=\mu(\alpha)$. The extended involution reverses products, because $\kappa$ is multiplicative and cancellation removes pairs of letters $\alpha^{\pm1}\alpha^{\mp1}$, for which $\kappa(\alpha)\kappa(\alpha^{-1})=1$. It squares to the identity, because on a letter $\alpha^{-1}$ its square is multiplication by the same factor as on $\alpha$, namely $1$. The same holds for the conjugation.

\section{Completion in the Path Degree}
\label{sec:completion}
By Proposition~\ref{prop:localisation}(2), applied to each connected component, $kQ_{\mathrm{loc}}$ is graded by the degree of reduced words, since cancellation removes letters in pairs of total degree $0$:
\[ kQ_{\mathrm{loc}}=\bigoplus_{d\in\ZZ}(kQ_{\mathrm{loc}})_d,\qquad (kQ_{\mathrm{loc}})_d\,(kQ_{\mathrm{loc}})_{d'}\subseteq(kQ_{\mathrm{loc}})_{d+d'} . \]

\begin{definition}\label{def:completion}
The \emph{completed localisation} $\Lambda$ of $kQ$ is the set of formal sums $x=\sum_{d\geq d_0}x_d$ with $x_d\in(kQ_{\mathrm{loc}})_d$ and $d_0\in\ZZ$, with componentwise addition and the product $(xy)_d=\sum_{i+j=d}x_iy_j$. The \emph{valuation} $v(x)$ is the least $d$ with $x_d\neq0$, with $v(0)=\infty$, and $x^\circ=x_{v(x)}$ is the \emph{lowest-degree part} of $x$.
\end{definition}

The sum defining $(xy)_d$ is finite, because $x_i=0$ for $i<v(x)$ and $y_j=0$ for $j<v(y)$. Hence $\Lambda$ is an associative algebra that contains $kQ_{\mathrm{loc}}$, $v(xy)\geq v(x)+v(y)$, and $\Lambda$ is the direct sum of its corners $e_w\Lambda e_v$. Each homogeneous part $x_d$ is a finite linear combination of reduced words, but $x$ may involve infinitely many degrees. If an involution or a conjugation maps arrows to multiples of arrows, it preserves degrees and extends termwise to $\Lambda$.

For one vertex with one loop $t$, $\Lambda$ is the field $k((t))$ of formal Laurent series, and for the real path $*$-algebra with $t^*=t$ it is the field $\RR((t))$ of Section~\ref{sec:pythagorean}, in which sums of squares have square roots (Proposition~\ref{prop:laurent-pythagorean}). Thus the same construction supplies the field of Laurent series in which Theorem~\ref{thm:bruhat} computes the Bruhat decomposition, and the inverses of the present chapter; Chapter~\ref{chap:implementation} uses it for both. It is not real closed; the limits of Proposition~\ref{prop:autolimit} are taken in the real closed field of algebraic Puiseux series, which contains every element of $\RR((t))$ obtained from polynomial inputs by field operations and square roots. For the quiver of Section~\ref{sec:unification-example}, $\Lambda\cong M_2(k((g)))$ with $g=ba$ (Proposition~\ref{prop:lambda-complete}).

In computations, an element of $\Lambda$ is represented by finitely many homogeneous parts together with a \emph{precision} $p$: it is known modulo terms of degree at least $p$, as a power series is known to finitely many terms. The product of $x$ and $y$ is then known to precision $\min(p_x+v(y),\,p_y+v(x))$, and the operations of the next two sections are carried out to a prescribed number of degrees. This makes the arithmetic of $\Lambda$ effective.

\section{Inverses in the Completed Localisation}
\label{sec:lambda-inverses}
\begin{lemma}\label{lem:lambda-units}
Let $f$ and $g$ be sums of vertex idempotents, let $x\in g\Lambda f$, and suppose that the lowest-degree part $x^\circ$ of $x$ is invertible from $f$ to $g$ in $kQ_{\mathrm{loc}}$. Then $x$ is invertible from $f$ to $g$ in $\Lambda$, with inverse
\[ x^{-1}=\sum_{k\geq0}u^k(x^\circ)^{-1},\qquad u=f-(x^\circ)^{-1}x, \]
where $v(u)\geq1$.
\end{lemma}
\begin{proof}
Let $d=v(x)$ and $y=(x^\circ)^{-1}$. The inverse $y$ of the homogeneous element $x^\circ$ is homogeneous of degree $-d$: its part of degree $-d$ is also an inverse, as one sees by comparing the parts of degree $0$ in $yx^\circ=f$ and $x^\circ y=g$, and inverses are unique. Hence $u=-y(x-x^\circ)$ lies in $f\Lambda f$ and has valuation at least $-d+(d+1)=1$. The series $\sum_ku^k$ converges in $\Lambda$, since $v(u^k)\geq k$, and $(f-u)\sum_{k\geq0}u^k=\sum_{k\geq0}u^k(f-u)=f$. Since $x=x^\circ yx=x^\circ(f-u)$, we get $x^{-1}x=\sum_ku^k(f-u)=f$ and $xx^{-1}=x^\circ(f-u)\sum_ku^ky=x^\circ y=g$.
\end{proof}

The lemma applies, for example, to the arrows between vertices of equal block size, to $e_v+\alpha$ for a loop $\alpha$ at $v$, and to Laurent-type elements such as $2\alpha+\alpha\alpha^{\top}\alpha$, whose lowest-degree part $2\alpha$ is a single word.

\paragraph{Elimination.} Other elements are inverted by Gauss--Jordan elimination. With block sizes $1$, an element $x\in g\Lambda f$ is a matrix $(x_{wv})$ with entries $x_{wv}\in e_w\Lambda e_v$, indexed by the vertices $w$ in $g$ and $v$ in $f$; with larger block sizes, rows and columns are indexed by a vertex together with a row of its block. An entry can serve as a pivot if its lowest-degree part is a single word with a non-zero coefficient, since it is then invertible by Lemma~\ref{lem:lambda-units}. Among such entries the elimination chooses one of lowest degree and, among those, one whose leading coefficient has the largest modulus. It either inverts $x$ from $f$ to $g$, or it stops because no remaining entry of a column can serve as a pivot. The following proposition says when it cannot stop on an invertible element.

\begin{proposition}\label{prop:lambda-complete}
Let $Q$ be a connected quiver with cycle rank $r$, and let $\Lambda$ be the completed localisation of $kQ$, with arbitrary block sizes. The following are equivalent:
\begin{enumerate}
  \item between any two vertices there is at most one reduced word of each degree;
  \item $r=0$, or $r=1$ and the cycle of $Q$ has non-zero degree.
\end{enumerate}
If they hold, there is a field $K$, namely $k$ if $r=0$ and a field of Laurent series over $k$ if $r=1$, such that each corner $e_w\Lambda e_v$ is isomorphic to $M_{d_w\times d_v}(K)$, compatibly with multiplication. Moreover, in exact arithmetic, the elimination just described inverts every element of $\Lambda$ that is invertible between two given sums of vertex idempotents.
\end{proposition}
\begin{proof}
By Proposition~\ref{prop:localisation}, the reduced words $\omega:v\to w$ correspond bijectively to the elements $g\in F_r$, via $\Phi(\omega)=gE_{wv}$, and $\deg\omega=\deg\tau_w-\deg\tau_v+\varphi(g)$, where $\varphi:F_r\to\ZZ$ is the homomorphism with $\varphi(\xi_\alpha)=\deg(\tau_w^{-1}\alpha\tau_v)$ for $\alpha:v\to w$. So (1) says exactly that $\varphi$ is injective. This is impossible for $r\geq2$, because $F_r$ is not abelian. For $r=1$ it means that the generator, which corresponds to the cycle, has non-zero degree, and for $r=0$ it is automatic.

Under (2) with $r=1$, let $\zeta$ be the generator of $F_1$ of positive degree $\delta$. A reduced word $v\to w$ corresponds to $\zeta^mE_{wv}$ for a unique $m\in\ZZ$, and its degree $\deg\tau_w-\deg\tau_v+m\delta$ increases with $m$. Completing in the degree therefore turns $e_w\Lambda e_v$ into the set of series $\sum_{m\geq m_0}X_m\zeta^m$ with $X_m\in M_{d_w\times d_v}(k)$, that is, into $M_{d_w\times d_v}(K)$ with $K=k((\zeta))$, compatibly with multiplication. For $r=0$, take $K=k$. A non-zero scalar entry from $v$ to $w$ then corresponds to a non-zero element of $K$, and its lowest-degree part is a single word. So the possible pivots are exactly the non-zero entries, and the elimination is Gauss--Jordan elimination over the field $K$, which inverts every invertible matrix.
\end{proof}

For every quiver of Table~\ref{tab:catalogue}, each connected component of the complexified quiver satisfies condition (2); the implementation reports the components and their cycle ranks. So does an oriented cycle with $n$ vertices, whose completed localisation is $M_n(k((x)))$. When the condition fails, the elimination can stop on an invertible element.

\begin{example}\label{ex:elimination-fails}
Let $Q$ have two vertices $u,u'$ and two arrows $a,b:u\to u'$. This is the complexified quiver of the doubly-signed quiver of Section~\ref{subsec:cat-ldl} with a second arrow fixed by the transpose, which encodes pairs of real symmetric matrices under simultaneous congruence. Here $r=1$, but the cycle $a^{-1}b$ has degree $0$, and $e_u\Lambda e_u=k[g,g^{-1}]$ with $g=a^{-1}b$. Put $p=e_u+g$ and
\[ M=\begin{pmatrix}e_u+p^2&p^3+2p\\ p&e_u+p^2\end{pmatrix}\in M_2(k[g,g^{-1}]),\qquad \det M=(e_u+p^2)^2-(p^3+2p)p=e_u . \]
Conjugating $M$ by $\diag(e_u,a)$ transports it to the vertices $u$ and $u'$ and gives
\[ x=(e_u+p^2)+(p^3+2p)\,a^{-1}+a\,p+a(e_u+p^2)a^{-1}\in(e_u+e_{u'})\Lambda(e_u+e_{u'}), \]
which is invertible from $e_u+e_{u'}$ to itself; its inverse is obtained in the same way from $M^{-1}$. But each of the four entries of $x$ is homogeneous and has at least two words, for instance $ap=a+b$, so no entry can serve as a pivot. The implementation confirms that the elimination stops on $x$, and that the transported inverse of $M$ is an inverse of $x$.
\end{example}

\section{\texorpdfstring{$*$}{*}-Square Roots in the Completed Localisation}
\label{sec:lambda-sqrt}
Let $\Lambda$ be the completed localisation of the complexified path $*$-algebra of a doubly-signed quiver. A word $\omega:v\to w$ carries a coefficient block in $M_{d_w\times d_v}(\CC)$, and the involution has the form $(X\omega)^*=\kappa(\omega)\,G_vX^TG_w^{-1}\,\omega^{\top}$, with invertible matrices $G_v$ and signs $\kappa(\omega)=\pm1$ (Proposition~\ref{prop:block-model}). A $*$-square root of a self-adjoint $s\in\Lambda$ is any $x\in\Lambda$ with $x^*x=s$. Roots are not unique (Lemma~\ref{lem:sqrt-unique}), and Householder reflectors work with any of them that lies in the corner of the next position and for which the denominator of Proposition~\ref{prop:householder-any-root} is invertible. We compute roots by three rules. Each rule checks its result word by word, below the precision of $s$, and reports a failure instead of returning a wrong root.
\begin{itemize}
  \item \emph{The graded rule.} Let $s$ be homogeneous, with words of even length. Each word is cut in the middle, $\omega=\nu'^{\top}\nu$ with half-words $\nu,\nu'$ of equal length, and $x=\sum_\nu Z_\nu\nu$ is built from the half-words $\nu$, so that $x$ has half the degree of $s$. The coefficients of $s$ whose words pass through a middle vertex $m$ form a block matrix $N$, and the coefficients $Z$ of $x$ at $m$ solve a congruence $Z^TG_m^{-1}Z=N$ if $m^{\top}=m$, a factorisation $Z'^TZ=N$ if $m^{\top}\neq m$, and a Hermitian congruence at vertices of type \texttt{Cu}. For real $s$ the solution is made real. A non-homogeneous $s$ is treated like $\sqrt{t^{2k}(a+\cdots)}=t^k\sqrt{a+\cdots}$: its lowest-degree part receives a graded root $x_0$, and $x=rx_0$ for a principal root $r$ of degree $0$.
  \item \emph{The principal rule.} If the lowest-degree part $s_0$ of $s$ has degree $0$, is vertex-diagonal, and has blocks without eigenvalues in $(-\infty,0]$, then $x=x_0+x_1+\cdots$ is a self-adjoint square root, where $x_0$ is the principal square root of $s_0$ at each vertex and the part $x_d$ of degree $d$ solves the Sylvester equation $x_0x_d+x_dx_0=s_d-\sum_{0<i<d}x_ix_{d-i}$, which has a unique solution because the spectra of the blocks of $x_0$ lie in the open right half-plane.
  \item \emph{The arbitrary rule.} The graded root, multiplied by a random unitary that has terms of positive degree. It serves to test that an algorithm does not depend on the choice of root.
\end{itemize}
The refusals of the graded rule match the theory. For $s=-1_c$, where $1_c$ is the idempotent of an orbit $c$, it looks for a root whose words are trivial, that is, for an element $y$ of the algebra $K$ of the orbit with $y^*y=-1$. Such $y$ exists exactly when $K$ is not of signature type: if it is not, by Lemma~\ref{lem:star-simple}(1), and if $K$ is $\RR$, $\CC$ or $\HH$ with its standard involution, $y^*y=|y|^2\geq0$. Accordingly, the graded rule refuses a root of $-1$ exactly at the vertex types \texttt{R}, \texttt{Hsp} and \texttt{Cu} (Table~\ref{tab:vertex-blocks}), whose unitary groups are compact.

Over $\twoR$ it offers the roots $(c,1)$ and $(1,c)$ of a self-adjoint element $(c,c)$, and also $(\sqrt c,\sqrt c)$ if $c>0$ and $(-\sqrt{|c|},\sqrt{|c|})$ if $c<0$; the root $(c,1)$ is the one used in Section~\ref{subsec:2R} and in Proposition~\ref{prop:2R-algorithms}. The graded rule also refuses when $N$ has rank greater than $d_m$, since $Z$ has only $d_m$ rows.

It has two known gaps, in which a root exists but is not found: words with inverse letters, where cancellation can hide the half-words, and elements such as $(e_v+\beta)^*(e_v+\beta)$ for an arrow $\beta:w\to v$ with $w\neq v$, whose lowest-degree part does not reach every vertex that the element reaches. We have not proved that the graded rule finds a root of every element $z^*z$ with $z$ homogeneous; in the tests of Chapter~\ref{chap:implementation} it did so for the Gram elements of every quiver of the catalogue and for $306$ elements $z^*z$ with $z$ real and supported on a random path of a random doubly-signed quiver.

\section{Reducing Back to the Original Algebra}
\label{sec:reduction}
After an algorithm has been executed in $\Lambda$, its results must be brought back to $kQ/I$.

\begin{theorem}\label{thm:reduction}
Let $\widehat{kQ}\subseteq\Lambda$ consist of the elements whose words are all paths, that is, contain no letter $\alpha^{-1}$.
\begin{enumerate}
  \item $\widehat{kQ}$ is a subring of $\Lambda$ that contains $kQ$ and is stable under every involution and conjugation that maps arrows to multiples of arrows. An element of $\widehat{kQ}$ is invertible in $\widehat{kQ}$ if and only if its part of degree $0$, which lies in the semisimple algebra $S$ spanned by the vertices (with coefficient blocks, $S=\prod_vM_{d_v}(k)$), is invertible in $S$.
  \item Let $I$ be an ideal of $kQ$ with $R_Q^m\subseteq I$ for some $m$. Then the quotient map $kQ\to kQ/I$ extends to a ring homomorphism $\rho_I:\widehat{kQ}\to kQ/I$ that vanishes on the elements of valuation at least $m$. If $I$ is stable under the involution and the conjugation, then $\rho_I$ commutes with them.
\end{enumerate}
\end{theorem}
\begin{proof}
(1) A product of two paths is a path or zero, without cancellation, so $\widehat{kQ}$ is closed under multiplication, and it is closed under the involution and the conjugation because these map paths to multiples of paths. Elements of $\widehat{kQ}$ have non-negative valuation. If the part $x_0$ of degree $0$ of $x\in\widehat{kQ}$ is invertible in $S$, then $x^\circ=x_0$, and Lemma~\ref{lem:lambda-units}, with $f=g=1$, gives the inverse $\sum_ku^kx_0^{-1}$, whose words are paths because those of $u=1-x_0^{-1}x$ are. Conversely, if $xy=yx=1$ in $\widehat{kQ}$, then $x_0y_0=y_0x_0=1$.

(2) The elements of valuation at least $m$ form an ideal $\widehat{kQ}_{\geq m}$ of $\widehat{kQ}$, and truncation identifies $\widehat{kQ}/\widehat{kQ}_{\geq m}$ with $kQ/R_Q^m$, since the paths of length less than $m$ form a basis of both. Composing with the quotient map $kQ/R_Q^m\to kQ/I$ gives $\rho_I$. The involution and the conjugation preserve degrees, hence $\widehat{kQ}_{\geq m}$; if they also preserve $I$, they commute with $\rho_I$.
\end{proof}

\begin{corollary}\label{cor:transfer}
Let $I$ be as in Theorem~\ref{thm:reduction}(2) and stable under the involution and the conjugation. Suppose that an algorithm is run in $\Lambda$ on elements $y_1,\dots,y_k\in\widehat{kQ}$, which lift $a_i=\rho_I(y_i)\in kQ/I$, and that some of the quantities it computes, $z_1,\dots,z_l$, lie in $\widehat{kQ}$. Then every identity among the $y_i$ and the $z_j$ that is built from addition, multiplication, the involution and the conjugation, and that holds in $\Lambda$, also holds among the $a_i$ and the $\rho_I(z_j)$ in $kQ/I$. The same applies entrywise to matrices indexed by positions.
\end{corollary}

The corollary follows from Theorem~\ref{thm:reduction}(2). Intermediate quantities need not lie in $\widehat{kQ}$. In the Householder steps of Proposition~\ref{prop:golub-kahan}, $s^{-1}=c'^{-1}a^{-1}b^{-1}$ has inverse letters, but the column $x$, the root $\gamma$ and the reflector $P$ lie in $\widehat{kQ}$, so the identities $P^*P=I$ and $Px=e_1\gamma$, which certify each step, hold in $B=kQ/R_Q^2$. Similarly, the quotient $(xb)/(xc)=b/c$ of Example~\ref{ex:dual-division}, computed in $k((x))$, lies in $\widehat{kQ}=k[[x]]$, and its image $b/c$ in the dual numbers is the intended value. A quantity that does not lie in $\widehat{kQ}$, such as $x/x^2=x^{-1}$, has no image in $kQ/I$: if it is an output, the generic algorithm has produced a result outside the original algebra, which Section~\ref{sec:wild-behaviour} discusses for wild quivers.

For the loop quiver, Theorem~\ref{thm:reduction} takes the following classical form.

\begin{proposition}\label{prop:reduction-loop}
Let $\mathcal O\subseteq k(x)$ be the ring of rational functions without a pole at $x=0$. The map
\[ \mathcal O\to k[x]/(x^2),\qquad f\mapsto f(0)+f'(0)\,x, \]
is a surjective ring homomorphism whose restriction to $k[x]$ is the quotient map.
\end{proposition}
\begin{proof}
Taylor expansion at $0$ embeds $\mathcal O$ in $k[[x]]=\widehat{kQ}$, and $\rho_{(x^2)}$ maps a power series $f$ to $f(0)+f'(0)\,x$.
\end{proof}

Proposition~\ref{prop:reduction-loop} is Proposition~\ref{prop:dual-quotient} in another guise: if $x$ is a positive infinitesimal in a real closed field, then reduction is first-order Taylor expansion. In this sense, the completed localisation replaces nilpotent infinitesimals by invertible ones: the computation is carried out in a ring in which they are invertible, and the result is truncated at the end.

\section{Beyond the Completed Localisation: Simple Artinian Rings}
\label{sec:embedding}
When the condition of Proposition~\ref{prop:lambda-complete} holds and all block sizes are $1$, every non-zero element of $\Lambda$ between two vertices is invertible. Otherwise this can fail. In Example~\ref{ex:elimination-fails}, the element $e_u+g$ of $e_u\Lambda e_u=k[g,g^{-1}]$ is not invertible, since the units of that ring are the non-zero multiples of powers of $g$; for a vertex with two loops $x$ and $y$, neither is $x\alpha+y\beta$ with non-zero scalars $\alpha,\beta$ (Section~\ref{sec:wild-behaviour}). Rings in which every non-zero element between two vertices is invertible are provided by simple Artinian rings (Definition~\ref{def:artinian}). By the Wedderburn--Artin theorem (Theorem~\ref{thm:wedderburn-artin}), such a ring is $M_n(D)$ for a skew field $D$, and in $M_n(D)$ a non-zero element $dE_{ji}$ is invertible from $E_{ii}$ to $E_{jj}$, with inverse $d^{-1}E_{ij}$. This section shows that the path algebra of every finite quiver embeds in such a ring, and when $D$ can be taken commutative. No algorithm of this thesis uses these embeddings; they show what lies beyond $\Lambda$.

Let $Q$ be a quiver with vertices $1,\dots,n$ and cycle rank $r$ (Section~\ref{sec:tits_form}). Choose a spanning forest $T$ of its underlying graph, and label the $r$ arrows not in $T$ by symbols $x_\alpha$. Let $D_r$ be $k$ if $r=0$, the field of rational functions $k(x)$ if $r=1$, and Cohn's \emph{free field} $k\{x_\alpha:\alpha\notin T\}$ if $r\geq2$. The free field $k\{X\}$ is the universal skew field of fractions of the free algebra $k\langle X\rangle$: an $m\times m$ matrix over $k\langle X\rangle$ is invertible over $k\{X\}$ if and only if it is \emph{full}, that is, not a product of an $m\times s$ and an $s\times m$ matrix over $k\langle X\rangle$ with $s<m$, and every skew field generated by $k\langle X\rangle$ is a specialisation of $k\{X\}$~\cite[Chapter~7]{Cohn2006}. The group algebra of the free group on $X$ embeds in $k\{X\}$~\cite{Lewin1974}, and the word problem for $k\{X\}$ is decidable~\cite{CohnReutenauer1999}, so its arithmetic is effective, although much more expensive than that of rational functions.

\begin{theorem}[Embedding theorem]\label{thm:embedding}
\begin{enumerate}
  \item The assignments
  \[ e_v\mapsto E_{vv},\qquad \alpha\mapsto E_{t(\alpha)s(\alpha)}\ (\alpha\in T),\qquad \alpha\mapsto x_\alpha\,E_{t(\alpha)s(\alpha)}\ (\alpha\notin T) \]
  extend to an injective $k$-algebra homomorphism $\Phi:kQ\to M_n(D_r)$. Consequently, every non-zero element of $e_jkQe_i$ becomes invertible from $E_{ii}$ to $E_{jj}$ in $M_n(D_r)$.
  \item If no vertex of $Q$ lies on two distinct simple oriented cycles, then the same assignments, with the $x_\alpha$ interpreted as \emph{commuting} indeterminates, give an injective homomorphism $kQ\to M_n(k(x_\alpha:\alpha\notin T))$ into a matrix ring over a field of rational functions.
  \item If some vertex of $Q$ lies on two distinct simple oriented cycles, then there is no injective homomorphism from $kQ$ to a matrix ring $M_N(F)$ over a commutative ring $F$, for any $N$.
\end{enumerate}
The homomorphisms of (1) and (2) are given by explicit formulas, and the arithmetic in their targets is effective; in this sense the embedding is constructive.
\end{theorem}

\begin{proof}
(1) For connected $Q$, $\Phi$ is the composite of the injective map $kQ\to kQ_{\mathrm{loc}}\cong M_n(k[F_r])$ of Proposition~\ref{prop:localisation} with the inclusion of $k[F_r]$ in $D_r$: for $r=1$ this is $k[x,x^{-1}]\subseteq k(x)$, and for $r\geq2$ it is Lewin's embedding~\cite{Lewin1974}. For a quiver with several components, apply this to each component, with distinct symbols, and place the blocks along the diagonal. A non-zero element of $e_jkQe_i$ maps to $dE_{ji}$ with $d\neq0$, which is invertible from $E_{ii}$ to $E_{jj}$.

(2) Suppose that no vertex lies on two distinct simple oriented cycles. Then each strongly connected component of $Q$ is a single vertex or a single simple cycle. (Otherwise the component contains a simple cycle $Z$ and, by the ear decomposition of strongly connected digraphs~\cite{BangJensenGutin2009}, a path $P$ outside $Z$ whose endpoints $u,w$ lie on $Z$; then $u$ lies on $Z$ and on the simple cycle formed by $P$ and the part of $Z$ from $w$ to $u$.) The argument below works for every choice of the spanning forest $T$. With commuting variables, a path $p$ maps to a commutative monomial, which records how often each non-tree arrow is used. We show that this determines $p$. First, the number of times $n_\alpha$ that $p$ uses each arrow $\alpha$ is determined: two count vectors with the same endpoints and the same non-tree counts differ by a vector supported on $T$ whose boundary is zero, that is, by a cycle in a forest, which must vanish. Second, the counts determine the path: the arrows between different strongly connected components are used at most once each, and the ones used form a chain in the (acyclic) condensation of $Q$, which determines the sequence of components visited and the entry and exit vertices; within a cyclic component, the walk from the entry vertex to the exit vertex is determined by the number of arrows it uses. Distinct paths therefore give distinct monomials, which are linearly independent in the field of rational functions.

(3) Let $v$ lie on two distinct simple oriented cycles $c_1$ and $c_2$, regarded as elements of $e_vkQe_v$. A closed path at $v$ factors uniquely into paths that return to $v$ only at their end, and $c_1$ and $c_2$ are two different such paths, so distinct words in $c_1$ and $c_2$ are distinct paths. Put $y_i=c_1c_2^{\,i}$ for $i\geq1$. A word in the $y_i$ is a word in $c_1$ and $c_2$ from which the $y_i$ can be read off, since each $y_i$ begins with $c_1$, so distinct words in the $y_i$ are distinct, hence linearly independent, paths. In particular, for every $N$ the standard polynomial
\[ s_{2N}(y_1,\dots,y_{2N})=\sum_{\sigma\in S_{2N}}\operatorname{sgn}(\sigma)\,y_{\sigma(1)}y_{\sigma(2)}\cdots y_{\sigma(2N)} \]
is non-zero in $kQ$. By the theorem of Amitsur and Levitzki~\cite{AmitsurLevitzki1950}, $s_{2N}(z_1,\dots,z_{2N})=0$ for all $z_1,\dots,z_{2N}$ in $M_N(F)$, for every commutative ring $F$. An injective homomorphism $kQ\to M_N(F)$ would therefore map a non-zero element to zero, which is impossible.
\end{proof}

By (2) and (3), a commutative coefficient field suffices exactly when no vertex lies on two distinct simple oriented cycles; otherwise the quiver has cycle rank at least $2$ and is wild (Proposition~\ref{prop:two-cycles-wild}), and a non-commutative skew field, such as the free field, is unavoidable. When the condition of Proposition~\ref{prop:lambda-complete} holds and all block sizes are $1$, $D_r$ is $k$ or $k(\zeta)$, which lies in the field $K$ of that proposition, and the canonical map $kQ\to\Lambda\cong M_n(K)$ is $\Phi$ followed by the inclusion $M_n(D_r)\subseteq M_n(K)$. For the Kronecker quiver of Example~\ref{ex:elimination-fails}, with $T=\{a\}$, we get $\Phi(b)=xE_{u'u}$, and the extension of $\Phi$ to $kQ_{\mathrm{loc}}$ maps $g=a^{-1}b$ to $xE_{uu}$, so the element $e_u+g$, which is not invertible in $\Lambda$, becomes the invertible element $(1+x)E_{uu}$ of $M_2(k(x))$. Here $k(x)$ inverts an element that $\Lambda$ cannot, because $g=a^{-1}b$ has degree $0$, and the completion in the path degree does not affect the part of degree $0$.

\paragraph{Involutions.} The embedding $\Phi$ ignores the involution and the conjugation of a doubly-signed quiver. For the quivers of the catalogue this does not matter: by Proposition~\ref{prop:lambda-complete}, $\Lambda$ is then a finite product of algebras $M_N(K)$ over fields, one for each component of the complexified quiver, and the involution and the conjugation permute the factors. For one vertex with block size $1$, the path algebra is the free algebra $k\langle X\rangle$ on the loops, and the free field carries the involution: an involution of $k\langle X\rangle$ is an isomorphism $k\langle X\rangle\to k\langle X\rangle^{\mathrm{op}}$, so, by the uniqueness of the universal skew field of fractions, it extends to an anti-automorphism of $k\{X\}$, whose square is the identity because it fixes the generating subring $k\langle X\rangle$; the same applies to the conjugation.

The same argument applies to every doubly-signed quiver, with the localisation in place of the path algebra. By Proposition~\ref{prop:localisation}, each connected component of the complexified localisation, with its blocks, is a matrix ring $M_N(\CC[F_r])$. The group algebra $\CC[F_r]$ is a fir, whose universal skew field of fractions $U$ is $\CC$, $\CC(x)$ or the free field according as $r=0$, $r=1$ or $r\geq2$~\cite{Cohn2006,Lewin1974}, and $M_N(U)$ is the universal localisation of $M_N(\CC[F_r])$ at the matrices that are full over $\CC[F_r]$~\cite[Chapter~7]{Cohn2006}.

Under the isomorphisms of Proposition~\ref{prop:localisation}, the involution maps the component of a vertex $v$ onto that of $v^{\top}$ by a map of the form $X\mapsto D\,\psi(X)^TD^{-1}$, where $\psi$ is the anti-isomorphism of group algebras induced by an anti-isomorphism of the free groups, which reverses words and transposes their letters, and $D$ is an invertible matrix whose entries are products of signs, the blocks $G_v$ and elements of the free group. The conjugation has the same form, with an isomorphism of the free groups, complex conjugation of the coefficients and no transpose. Maps of this form send matrices that are full over $\CC[F_r]$ to matrices that are full over $\CC[F_{r'}]$, since they preserve factorisations through matrices of smaller size, and a matrix over the opposite of a ring is invertible exactly when its transpose is invertible over the ring.

So, by the universal property of the localisations $M_N(U)$, the involution and the conjugation extend uniquely to the product of the rings $M_N(U)$, and the identities $*^2=\iota^2=1$ and $\iota*=*\iota$ extend by uniqueness. Hence, after complexification, the localisation of every doubly-signed quiver, and with it its path algebra, embeds compatibly with the involution and the conjugation in a finite product of simple Artinian rings that they permute. Whether the $*$-square roots that a Jacobi-type method needs exist there is open (Section~\ref{sec:open-problems}).
%======================================================================
%   CHAPTER 8: IMPLEMENTATION
%======================================================================
\chapter{Implementation}
\label{chap:implementation}

The seven operations of Section~\ref{sec:proposal}, and the algorithms written against them, have been implemented as a Python library, \texttt{dsqlib}. This chapter describes what the library computes, how its parts correspond to the theory of the preceding chapters, and what it has been used to check.

The library computes in the completed localisation of the path $*$-algebra of a doubly-signed quiver (Chapter~\ref{division-problem}): every arrow is invertible, and elements may be series in the path degree. Its matrix algorithms (Householder reduction, the QR algorithm, the Lanczos process, and the Cholesky, $LDL^*$ and $LU$ factorisations) are written once, against the seven operations. No algorithm in the library is specific to a decomposition; only the translation of a classical matrix into its encoding, and back, is.

In the library, the seven operations are written \texttt{x + y}, \texttt{x - y} and \texttt{x * y} for elements, with \texttt{@} for the product of matrices; \texttt{x * y.inverse()}, where \texttt{y.inverse()} inverts $y$ between the vertices that it touches; \texttt{x.H}, the adjoint of an element and the conjugate transpose of a matrix; \texttt{s.sqrt()}; and \texttt{H[i, j]}, the entry of a matrix at its $i$-th and $j$-th positions.

The library consists of about 3700 lines of documented Python; its source code is printed in Appendix~\ref{chap:code}, and its test suite, from which the measurements of this chapter are taken, is part of the source distribution submitted with the thesis. It uses floating-point arithmetic and depends only on NumPy and SciPy~\cite{Harris2020NumPy,Virtanen2020SciPy}. Its test suite of 65 tests records every quantity that it measures, and the errors quoted in this chapter are taken from that record. The results of the thesis do not depend on the implementation; the implementation shows that the constructions are effective, and it allows the encoded decompositions to be checked numerically.

\section{Overview}
\label{sec:impl-overview}
Table~\ref{tab:implementation} lists the components of the library; each builds on those above it. Only \texttt{quiver} refers to the vertex types. The other components see only the block data that it produces (Section~\ref{sec:impl-quivers}), and only the square-root rules use separate formulas for different block data.

\begin{table}[htbp]
\centering
\caption{The components of \texttt{dsqlib} and the parts of the thesis that they implement.}
\label{tab:implementation}
\small
\renewcommand{\arraystretch}{1.25}
\begin{tabular}{@{}>{\raggedright\arraybackslash}p{0.2\linewidth}>{\raggedright\arraybackslash}p{0.42\linewidth}>{\raggedright\arraybackslash}p{0.3\linewidth}@{}}
\toprule
\textbf{Module} & \textbf{Contents} & \textbf{Theory} \\
\midrule
\texttt{quiver} & doubly-signed quivers built from the ten vertex types; the block model; the signs of the arrows & Definition~\ref{def:doubly-signed}, Proposition~\ref{prop:block-model} \\
\rowrule
\texttt{algebra} & the path $*$-algebra, localised at the arrows and completed in path degree; inverses between vertices & Definition~\ref{def:path-star-algebra}, Sections~\ref{sec:localisation}--\ref{sec:lambda-inverses} \\
\rowrule
\texttt{sqrt} & $*$-square roots by three rules & Sections~\ref{sec:hermitian-sqrt} and~\ref{sec:lambda-sqrt} \\
\rowrule
\texttt{matrix} & matrices indexed by positions with self-adjoint idempotents & Section~\ref{sec:indexing}, Definition~\ref{def:matrix-like} \\
\rowrule
\texttt{linalg} & Householder, QR, Lanczos, Cholesky, $LDL^*$ and $LU$ & Propositions~\ref{prop:householder-any-root}, \ref{prop:golub-kahan}, \ref{prop:householder-graded}, \ref{prop:qr-steps} and~\ref{prop:2R-algorithms} \\
\rowrule
\texttt{catalogue}, \texttt{problems} & the quivers of the catalogue and the matrix problems that they encode & Theorem~\ref{thm:encoding}, Table~\ref{tab:catalogue} \\
\bottomrule
\end{tabular}
\end{table}

Listing~\ref{lst:session} builds the doubly-signed quiver of the complex singular value decomposition (Section~\ref{subsec:cat-complex-svd}) and computes in its algebra. It then runs Example~\ref{ex:golub-kahan-numbers}: it encodes the real $3\times2$ matrix $W$ of that example as a Hermitian matrix $H$ (Theorem~\ref{thm:encoding}), reduces $H$ by the generic Householder algorithm, checks that the unitary factor has degree $0$, and decodes the result, the bidiagonal matrix computed by hand in the example. The comments show the values computed; the library numbers positions from $0$.

\begin{lstlisting}[float=tbp,label={lst:session},caption={A session with \texttt{dsqlib}: arithmetic in the algebra of the complex SVD quiver, and Example~\ref{ex:golub-kahan-numbers} run by the generic Householder algorithm. The comments show the values computed.}]
import numpy as np
from dsqlib import DSQ, Algebra, problems
from dsqlib.linalg import householder

Q = DSQ.build([("u", "Cu"), ("v", "Cu")], [("a", "u", "v")])
A = Algebra(Q)
a = A.arrow("a")
A.word(["a^-1"]) * a                       # 1·e_u
z = A.real_element(A.arrow("a", 1 + 2j))  # (1+2j)·a + (1-2j)·a~
y = (z.H * z).sqrt()                       # 2.236·a + 2.236·a~

W = np.array([[3, 1], [4, 2], [0, 0.3]])
H = problems.svd_problem(W, "R")    # positions u0, v0, u1, v1, v2
T, U, _ = householder(H)            # H = U T U*
U.degree_range()                    # (0, 0)
problems.extract_svd(T, "R").real   # [[-5, 2.2], [0, 0.5], [0, 0]]
\end{lstlisting}

\section{Doubly-Signed Quivers}
\label{sec:impl-quivers}
A doubly-signed quiver is entered by giving its vertices, each with one of the ten types of Table~\ref{tab:vertex-types}, named as in Table~\ref{tab:vertex-blocks}, and some of its arrows. The library works with the complexified quiver, in which each vertex type expands into one, two or four complex vertices: $b$, its conjugate $\bar b$ (written \texttt{b\textasciitilde}) and its transpose $b^{\top}$ (written \texttt{b'}). It completes each arrow to its orbit under the transpose and the conjugation, asks for the sign $\epsilonform(\alpha)$ of each arrow $\alpha:v\to v^{\top}$, which it takes to be fixed by the transpose, takes each arrow from $v$ to $v^H$ to be fixed by the conjugate transpose, and checks the rules of Definition~\ref{def:doubly-signed}. It then sets up the block model of Proposition~\ref{prop:block-model}: the block sizes and the matrices $G_v$ and $C_v$ of Table~\ref{tab:vertex-blocks}, and signs $\kappa_\alpha$ and $\mu_\alpha$ that satisfy \eqref{eq:sign-conditions}, which it finds by a search over each orbit of arrows; by Proposition~\ref{prop:block-model}(2), the search always succeeds. On random elements, for every quiver of the catalogue and for 120 random doubly-signed quivers covering all ten vertex types, the axioms of a $*$-algebra with a commuting conjugation hold to a relative error of $8\cdot10^{-17}$.

\section{The Algebra and Its Operations}
\label{sec:impl-algebra}
The module \texttt{algebra} implements the completed localisation $\Lambda$ of Chapter~\ref{division-problem} for the complexified path $*$-algebra of a doubly-signed quiver, with the block data of Section~\ref{sec:impl-quivers}. An element is stored as a dictionary from reduced words to coefficient blocks, together with its precision (Section~\ref{sec:completion}). Series inverses and square roots of exact elements are computed to a fixed number of degrees beyond the lowest one. A block in which a sum cancels to within a relative tolerance ($10^{-12}$ by default) of its summands is set to zero, so that zeros of the algebra appear as exact zeros at any scale.

Inverses are computed as in Section~\ref{sec:lambda-inverses}: by the series of Lemma~\ref{lem:lambda-units} when an element lies in a single corner and its lowest-degree part is a single word with an invertible block, and otherwise by elimination with the pivoting described there, which is complete under the condition of Proposition~\ref{prop:lambda-complete}. When cancellation uses up the precision of an exact input, the series are recomputed with more terms. By default an element is inverted between the vertices that it touches (Section~\ref{sec:division}); inversion between the orbit idempotents that it touches, or in the whole algebra, is also available. Every quotient $x/y$ computed by an algorithm is checked to satisfy $(x/y)\,y=x$, and a breakdown is reported when it does not. The $*$-square roots are computed by the three rules of Section~\ref{sec:lambda-sqrt}. The test suite checks the graded rule on the Gram elements of every quiver of the catalogue and on $306$ random paths, its refusals at all ten vertex types, and the principal and arbitrary rules on random elements.

Matrices are indexed by positions (Section~\ref{sec:indexing}). The \emph{plain} indexing gives every position the idempotent $1$. The \emph{orbit} indexing gives the position $c_k$ the orbit idempotent $1_c$ and lists the positions in a cyclic order, which can be chosen. Products and sums check that the algebras and the indexings agree, and a method reports the entries that do not lie in their corners.

\section{Algorithms Written Once}
\label{sec:impl-algorithms}
The algorithms use only the seven operations of Section~\ref{sec:proposal}, with division by inverses between vertices and indexing by positions, together with comparisons to choose signs and pivots. The square-root rule and the mode of inversion are parameters of each algorithm. A zero divisor that has to be inverted, or an isotropic vector, that is a vector $v\neq0$ with $v^*v=0$, stops an algorithm with an error that reports the breakdown. What the algorithms compute on the encoded matrices is stated, and mostly proved, in Section~\ref{sec:encoded-algorithms} (Table~\ref{tab:algorithms}); here we describe how they are implemented.

\paragraph{Householder reflections.} The reflectors are those of Proposition~\ref{prop:householder-any-root}, with $\gamma$ taken from the square-root rule and replaced by the root $-\gamma$ when this avoids cancellation in $x_1-\gamma$. The denominator $s=w^*x$ is inverted in the corner of its column, and a breakdown is reported when it is not invertible there; over $\twoR$ this depends on the root (Section~\ref{sec:hermitian-sqrt}). With the orbit indexing, the library checks that $\gamma$ lies in the corner of the position to which it is assigned, and it refuses an indexing for which this fails (Section~\ref{sec:unification-example}).

\paragraph{The QR algorithm and other algorithms.} A QR step orthogonalises the columns by the Gram--Schmidt process and normalises each residual column $v_j$ by a $*$-square root of $v_j^*v_j$. If $H$ is rank-deficient, the default inversion maps zero to zero, so that $Q$ is a partial isometry and $H=QR$ still holds. The Lanczos process uses full reorthogonalisation. The Cholesky, $LDL^*$ and $LU$ factorisations are the textbook algorithms, with every pivot inverted between vertices and every quotient checked.

\section{Experiments}
\label{sec:impl-experiments}
Table~\ref{tab:experiments} summarises what the generic algorithms compute on the quivers of the catalogue and on $\twoR$, for random inputs. Each result was compared with an independent computation (singular values, eigenvalues and QR factorisations from NumPy, $LU$ factorisations and LR steps from a separate implementation of Gaussian elimination, quaternionic QR factorisations from a separate Gram--Schmidt process) or checked by reconstructing the input from the computed factors. The error shown is the largest one observed over all variants tested. For the Householder reduction on the SVD quiver, the variants are real, complex and quaternionic matrices, square and rectangular shapes, both cyclic orders, and the graded root with and without a random unitary factor of degree $0$. For the quivers of Table~\ref{tab:catalogue} other than the congruence quiver, on which Householder reduction stops (Section~\ref{subsec:cat-ldl}), and the similarity quiver, whose computation was compared with the one over $\twoR$, the unitary factors of the Householder reduction were also checked to be real, unitary and of degree $0$, as Proposition~\ref{prop:householder-graded} predicts, and for the SVD quivers those of the Householder reduction and of the Lanczos process were checked to have every entry in its corner. With the arbitrary rule instead, the unitary factor of the Householder reduction has terms of degrees $0$ to $4$, and the reduced matrix is still tridiagonal.

\begin{table}[htbp]
\centering
\caption{The classical results reproduced by the generic algorithms of \texttt{dsqlib}, with the largest error observed. Numbers in parentheses refer to the sections in which the quivers are described.}
\label{tab:experiments}
\small
\renewcommand{\arraystretch}{1.25}
\begin{tabular}{@{}>{\raggedright\arraybackslash}p{0.24\linewidth}>{\raggedright\arraybackslash}p{0.16\linewidth}>{\raggedright\arraybackslash}p{0.39\linewidth}>{\raggedright\arraybackslash}p{0.13\linewidth}@{}}
\toprule
\textbf{Quiver} & \textbf{Algorithm} & \textbf{Classical result} & \textbf{Error} \\
\midrule
SVD over $\RR$, $\CC$, $\HH$ (\ref{subsec:cat-svd}--\ref{subsec:cat-quaternion-svd}) & Householder & bidiagonalisation~\cite{GolubKahan1965}; singular values kept & $5\cdot10^{-15}$ \\ \addlinespace[3pt]
 & Lanczos & \mbox{Golub--Kahan--Lanczos} bidiagonalisation & $2\cdot10^{-15}$ \\ \addlinespace[3pt]
 & QR step & two steps: a QR step on $W^*W$ and $WW^*$; one step: a Cholesky LR step on each (Proposition~\ref{prop:qr-steps}) & $9\cdot10^{-15}$ \\ \addlinespace[3pt]
 & double QR steps & singular values, after 20--28 double steps & $4\cdot10^{-15}$ \\
\rowrule
symmetric, skew (\ref{subsec:cat-symmetric}, \ref{subsec:cat-skew}) & Householder & tridiagonalisation; eigenvalues kept & $9\cdot10^{-16}$ \\ \addlinespace[3pt]
Hermitian over $\CC$, $\HH$ (\ref{subsec:cat-hermitian}) & Householder & tridiagonalisation; eigenvalues kept & $7\cdot10^{-15}$ \\ \addlinespace[3pt]
symmetric, skew, Hermitian over $\CC$, $\HH$ & QR step & classical QR step (Proposition~\ref{prop:qr-steps}) & $2\cdot10^{-15}$ \\ \addlinespace[3pt]
Hermitian over $\CC$ (\ref{subsec:cat-hermitian}) & QR iteration & eigenvalues, after 40 steps & $2\cdot10^{-14}$ \\
\rowrule
Takagi, Youla (\ref{subsec:cat-takagi}) & Householder & tridiagonalisation; Takagi or singular values kept & $5\cdot10^{-15}$ \\ \addlinespace[3pt]
Takagi (\ref{subsec:cat-takagi}) & QR step & two steps: a QR step on $WW^*$; one step: $Q_W^*W\overline{Q_W}$, a Cholesky LR step (Proposition~\ref{prop:qr-steps}) & $1\cdot10^{-15}$ \\
\rowrule
$\twoR$ (\ref{subsec:2R}) & Householder & tridiagonalisation of a real matrix by similarity & $1\cdot10^{-13}$ \\ \addlinespace[3pt]
 & QR step & two LR steps~\cite{Rutishauser1958}, relative error (Proposition~\ref{prop:2R-algorithms}) & $6\cdot10^{-14}$ \\ \addlinespace[3pt]
 & Cholesky & Doolittle $LU$ decomposition (Proposition~\ref{prop:2R-algorithms}) & $2\cdot10^{-15}$ \\
\rowrule
similarity (\ref{subsec:cat-similarity}) & Householder, QR step & the computations over $\twoR$ (Proposition~\ref{prop:2R-algorithms}) & $0$ \\
\rowrule
congruence (\ref{subsec:cat-ldl}) & $LDL^*$ & $LDL^T$ with the inertia of $W$ & $3\cdot10^{-16}$ \\
\rowrule
one loop, $\RR((t))$ (\ref{sec:bruhat}) & Gram--Schmidt & Bruhat decomposition, all 24 cells for $n=4$ (Theorem~\ref{thm:bruhat}) & $6\cdot10^{-16}$ \\
\bottomrule
\end{tabular}
\end{table}

The last row checks Theorem~\ref{thm:bruhat} on one matrix $M=LwU$ with random triangular factors in each of the $24$ Bruhat cells of $\GL_4(\RR)$. In every case the orders of the norms give $w$, the standard parts give a Bruhat decomposition $M=L'wU'$, and the orthogonal factor of $\diag(1,t,t^2,t^3)M$ tends to the signed permutation matrix of $w$, with the signs of Theorem~\ref{thm:bruhat}(3). The tests also check the explicit signs of Proposition~\ref{prop:block-model}(2) on $300$ random doubly-signed quivers, the conditions of Lemma~\ref{lem:degree-zero} for the catalogue, and that the elimination refuses to invert $x+2y$ for the quiver with two loops, as Section~\ref{sec:wild-behaviour} predicts.

\section{Limitations}
\label{sec:impl-limitations}
All arithmetic is in floating point, and series are truncated after a fixed number of degrees; exact arithmetic over Puiseux series, in which the auto-limit of Proposition~\ref{prop:autolimit} is exact, is not provided. Quotients $\RR Q/I$ by admissible ideals (Theorem~\ref{thm:structure}) are available only for $I=\RR Q_{\geq m}$, and the indefinite forms of Remark~\ref{rem:indefinite} are not implemented. Inversion is complete under the condition of Proposition~\ref{prop:lambda-complete} and can fail without it (Example~\ref{ex:elimination-fails}), and the graded rule has the two gaps described in Section~\ref{sec:lambda-sqrt}. The QR algorithm is unshifted and has no deflation, and although the experiments show errors at the level of rounding, no analysis of stability has been made (Section~\ref{sec:open-problems}). Finally, the library is written in pure Python over dictionaries of small blocks; a Householder reduction for a $10\times10$ complex SVD, with $20$ positions, takes about one second. Specialised classical routines, such as those of LAPACK that NumPy and SciPy call, compute these decompositions much faster and are better studied; the library is not meant to replace them.
%======================================================================
%   CHAPTER 9: CONCLUSION
%======================================================================
\chapter{Conclusion and Future Work}
\label{chap:conclusion}

\section{Summary of Contributions}
\label{sec:summary}
This thesis has overloaded the seven operations against which the algorithms of numerical linear algebra are written (addition, subtraction, multiplication, division, the adjoint and the $*$-square root of scalars, and $(i,j)$ indexing of matrices), so that the scalars of an algorithm may come from a $*$-algebra, and it has developed the algebraic theory that makes this work. The outcome is that decompositions are unified, as instances of unitary similarity of Hermitian matrices over the $*$-algebras of doubly-signed quivers (Theorem~\ref{thm:A}), and so are several of the algorithms that compute them: one program becomes different classical algorithms over different $*$-algebras (Theorem~\ref{thm:B}). Most of the results of the thesis are thus about one problem, the eigendecomposition of Hermitian matrices over $*$-algebras, and about the algorithms that compute it (Chapter~\ref{chap:intro} and Section~\ref{sec:chain}). Each classical result can also be reached without this reduction. The encoded problems have scalars in finite-dimensional $*$-algebras, but the generic algorithms run in the completed localisation of a path $*$-algebra, which is in general infinite-dimensional; their results return to the finite-dimensional algebra when they contain no inverted arrows. We summarise the results under the two outcomes of Section~\ref{sec:outcomes}; Table~\ref{tab:status} lists the main results with their status.

\begin{table}[htbp]
\centering
\caption{The main results of the thesis and their status: proved, observed in the experiments of Chapter~\ref{chap:implementation}, or open.}
\label{tab:status}
\small
\renewcommand{\arraystretch}{1.2}
\begin{tabular}{@{}>{\raggedright\arraybackslash}p{0.37\linewidth}>{\raggedright\arraybackslash}p{0.18\linewidth}>{\raggedright\arraybackslash}p{0.37\linewidth}@{}}
\toprule
\textbf{Result} & \textbf{Where} & \textbf{Status} \\
\midrule
Classification of basic real $*$-algebras with radical square zero & Theorem~\ref{thm:classification} & proved; existence also by an explicit construction (Proposition~\ref{prop:block-model}(3)) \\
\rowrule
Analogue of Gabriel's theorem (Theorem~\ref{thm:A}(1)) & Theorem~\ref{thm:structure} & proved \\
\rowrule
Encoding theorem and catalogue (Theorem~\ref{thm:A}(2), (3)) & Theorem~\ref{thm:encoding}, Section~\ref{sec:catalogue} & proved \\
\rowrule
Householder reflectors with any $*$-square root & Proposition~\ref{prop:householder-any-root} & proved \\
\rowrule
Golub--Kahan bidiagonalisation (Theorem~\ref{thm:B}(1)) & Proposition~\ref{prop:golub-kahan} & proved for real $W$ with $m\geq n$; the other cases observed \\
\rowrule
Householder reduction in degree $0$ (Theorem~\ref{thm:B}(1)) & Proposition~\ref{prop:householder-graded} & proved when every step is defined; that the steps are defined is observed for most quivers of the catalogue \\
\rowrule
QR steps for the SVD and Takagi quivers (Theorem~\ref{thm:B}(2)) & Proposition~\ref{prop:qr-steps} & proved \\
\rowrule
QR and Cholesky over $\twoR$ (Theorem~\ref{thm:B}(3)) & Proposition~\ref{prop:2R-algorithms} & proved \\
\rowrule
Bruhat decomposition (Theorem~\ref{thm:B}(4)) & Theorem~\ref{thm:bruhat} & proved \\
\rowrule
Dijkstra's algorithm as an auto-limit & Proposition~\ref{prop:dijkstra} & proved \\
\rowrule
QR factorisation as a limit of the singular value decomposition & Proposition~\ref{prop:qr-from-svd} & proved \\
\rowrule
Inverses and reduction in $\Lambda$ & Proposition~\ref{prop:lambda-complete}, Theorem~\ref{thm:reduction} & proved \\
\rowrule
Embedding in simple Artinian rings, and compatibility with the involution & Theorem~\ref{thm:embedding}, Section~\ref{sec:embedding} & proved \\
\rowrule
Lanczos as \mbox{Golub--Kahan--Lanczos}; $LDL^*$ with the inertia & Table~\ref{tab:experiments} & observed \\
\rowrule
Jacobi-type methods on wild quivers & Section~\ref{sec:wild-behaviour} & quotients outside $\Lambda$ proved for two loops; the rest open \\
\bottomrule
\end{tabular}
\end{table}

\paragraph{The seven operations.} Addition, subtraction, multiplication and the adjoint are those of the path $*$-algebra of a doubly-signed quiver (Proposition~\ref{prop:block-model}). Indexing is by positions, which the structure theorem shows to be intrinsic (Corollary~\ref{cor:positions}). Division is defined between idempotents. No non-zero element of the radical is invertible there, so the algorithms divide in the completed localisation $\Lambda$, where inverses are computed by series and by elimination (Proposition~\ref{prop:lambda-complete}). A $*$-square root is not unique, but Householder reflectors work with any root for which the denominator of the reflector is invertible (Proposition~\ref{prop:householder-any-root}). A result holds in the original algebra whenever the computed factors contain no inverse letters (Theorem~\ref{thm:reduction}). Over a field with an infinitesimal, running an algorithm once and taking standard parts computes its limit (Proposition~\ref{prop:autolimit}).

\paragraph{Unified decompositions.} Every finite-dimensional real $*$-algebra is a matrix algebra over a quotient of the path $*$-algebra of its doubly-signed quiver, and the quiver, the sizes and the signatures are invariants (Theorem~\ref{thm:structure}). Unitary similarity of Hermitian matrices of degree one over such an algebra is a classical matrix problem, whose group is a product of classical groups (Theorem~\ref{thm:encoding}). For the quivers of the catalogue its normal forms are the spectral theorems, the singular value decompositions, the Takagi factorisation, Sylvester's law of inertia and the real Jordan form (Table~\ref{tab:catalogue}). Matrices with a part of degree zero lead to first-order perturbation theory (Proposition~\ref{prop:dual-spectral}). Power series lead to spectral theorems over formal power series rings (Theorem~\ref{thm:spectral-power-series}), which generalise the dual-number results of~\cite{Gutin2022DualNumbers}, and to a parity-equivariant spectral theorem (Theorem~\ref{thm:parity}), which makes precise an observation of~\cite{Gutin2024Clifford}.

\paragraph{Unified algorithms.} Theorem~\ref{thm:B} says what the generic algorithms compute (Table~\ref{tab:algorithms}): the bidiagonalisation of Golub and Kahan, QR and Cholesky LR steps on Gram matrices, Rutishauser's LR algorithm and Doolittle's $LU$ factorisation over $\twoR$, and the Bruhat decomposition over a field of Laurent series. Beyond Theorem~\ref{thm:B}, Gaussian elimination over a field of Puiseux series is Dijkstra's shortest-path algorithm in leading terms (Proposition~\ref{prop:dijkstra}). A Python library implements the seven operations and these algorithms, and reproduces the classical algorithms to within rounding errors (Table~\ref{tab:experiments}).

\paragraph{Beyond the completed localisation.} When the quiver is wild, there is no normal form of the kind found in the catalogue, and a generic algorithm may need quotients that do not exist in $\Lambda$ (Section~\ref{sec:wild-behaviour}). They exist in a matrix ring over a skew field, in which the path algebra of every finite quiver embeds (Theorem~\ref{thm:embedding}); the skew field can be taken commutative if and only if no vertex lies on two distinct simple oriented cycles.

Taken together, these results show that $*$-algebras, combined with the operator overloading facilities of programming languages, unify the decompositions of the catalogue and several classical algorithms. Each of these decompositions can also be computed without $*$-algebras, by the classical algorithms for it, which are faster and better studied than the generic algorithms.

\section{Finite-Dimensional versus Operator-Algebraic Settings}
\label{sec:operator}
The thesis mentions Hilbert spaces and $C^*$-algebras, so it is natural to ask whether the framework is inherently finite-dimensional, and whether it could be generalised to operator algebras.

The generic algorithms run in an algebra that is, in general, not finite-dimensional. Finite-dimensional $*$-algebras are where the proposal starts, and Theorem~\ref{thm:structure} describes each of them as a matrix algebra over a quotient $\RR Q/I$ of a path $*$-algebra. The algorithms, however, run in the completed localisation $\Lambda$ of Chapter~\ref{division-problem}, which is infinite-dimensional as soon as the underlying graph of $Q$ contains a cycle (Proposition~\ref{prop:localisation}); for one vertex with one loop, it is the field $\RR((t))$ of formal Laurent series. Finite dimension is needed to describe the algebras from which the framework starts (through the Wedderburn--Malcev decomposition, the nilpotence of the radical and the description of basic algebras by quivers; Sections~\ref{subsec:radicals} and~\ref{sec:quivers-basic}), and to bring results back to them, which Theorem~\ref{thm:reduction} can do because $I$ contains all sufficiently long paths. What is finite throughout is the discrete data: the quiver, the sizes of the matrices and the blocks at the vertices. Accordingly, the transfer principle is applied to matrices of each fixed size separately (Section~\ref{subsec:fol}).

Two features, however, separate the framework from operator algebras. First, the finite-dimensional algebras from which it starts are in general not semisimple, whereas a $C^*$-algebra has zero Jacobson radical and its finite-dimensional examples are direct sums of full matrix algebras; the dual numbers, for example, admit no $C^*$-norm, since $\varepsilon^*\varepsilon=0$ (Remark~\ref{rem:cstar}). Second, $\Lambda$ deals with the radical algebraically, not analytically. For every quiver of the catalogue, $\Lambda$ is, after complexification, a finite product of matrix algebras over fields (Proposition~\ref{prop:lambda-complete}), as a finite-dimensional complex $C^*$-algebra is; but where $Q$ has a cycle these fields are fields of Laurent series, and the arrows are invertible elements of $\Lambda$ of positive degree, that is, invertible infinitesimals in place of nilpotent ones (Section~\ref{sec:reduction}). $\Lambda$ is complete for the path degree, not for a norm. Indeed, already for one vertex with one loop, the completed localisation of the real path $*$-algebra is $\RR((t))$, which admits no submultiplicative norm, because by the theorem of Gelfand and Mazur every real division algebra with a submultiplicative norm is isomorphic to $\RR$, $\CC$ or $\HH$~\cite{Zelazko1973}. Graph $C^*$-algebras, the operator algebras built from directed graphs, invert arrows differently (Section~\ref{sec:related}): there an arrow is a partial isometry, and it is invertible between vertices, with its adjoint as inverse, only when it is the only arrow leaving its source.

Nevertheless, there are several plausible directions of generalisation.
\begin{itemize}
  \item \emph{Semisimple parts.} The vertex groups of Table~\ref{tab:vertex-types} have infinite-dimensional analogues (unitary groups of real, complex and quaternionic Hilbert spaces), and Hermitian modules generalise to Hilbert $C^*$-modules, which carry inner products with values in a $C^*$-algebra~\cite{Lance1995}. Approximately finite-dimensional $C^*$-algebras, which are inductive limits of finite-dimensional ones~\cite{Davidson1996}, are a natural first case.
  \item \emph{Infinitesimals.} The elements of $\Lambda$ are formal series in the path degree, and Theorem~\ref{thm:spectral-power-series} is an algebraic form of Rellich's theorem. Analytic perturbation theory extends it to holomorphic families of self-adjoint operators on Hilbert spaces, at least for isolated eigenvalues of finite multiplicity~\cite{Kato1995}. Convergent power series with operator coefficients, or Banach algebras with non-zero radical, are thus natural analytic counterparts of $\Lambda$ and of the radical, rather than $C^*$-algebras.
  \item \emph{Sequences of matrices.} The $C^*$-algebraic approach to numerical analysis of Hagen, Roch and Silbermann~\cite{HagenRochSilbermann2001} studies sequences $(A_n)$ of matrices, such as finite sections of Toeplitz operators, in the $C^*$-algebra of bounded matrix sequences modulo null sequences. Running a generic algorithm uniformly along such a sequence amounts to replacing the scalars by a $*$-algebra built from these sequences, so their framework provides a possible bridge from the matrices of finite size treated here to operators on Hilbert spaces.
\end{itemize}
We conclude that the framework is not inherently finite-dimensional. Its algorithms run in $\Lambda$, which is in general infinite-dimensional, and finite dimension is needed at the two ends, to describe the algebras from which the framework starts and to bring results back to them. As it stands, however, the framework is algebraic rather than analytic, and its discrete data are finite: $\Lambda$ is complete for the path degree rather than for a norm, the radicals of the finite-dimensional algebras have no counterpart in $C^*$-algebras, and quivers and matrices are finite. Extensions to operator algebras appear feasible for the semisimple parts and for analytic perturbations, and we leave them to future work.

\section{Limitations and Open Problems}
\label{sec:open-problems}
\begin{itemize}
  \item \emph{Numerical behaviour.} We have not studied the numerical stability or conditioning of generic algorithms. The thesis addresses the algebraic obstacles first, because they turned out to be more substantial than expected. The experiments of Section~\ref{sec:impl-experiments} show errors at the level of rounding for the algorithms and quivers tested, but they do not replace an analysis.
  \item \emph{Relations.} Theorem~\ref{thm:structure} describes every finite-dimensional real $*$-algebra by its doubly-signed quiver together with an admissible ideal of relations. The encoding theorem (Theorem~\ref{thm:encoding}) applies to Hermitian matrices of degree one when the relations are homogeneous, as for $\RR[\varepsilon]/(\varepsilon^3)=\RR Q/(x^3)$. When they are not, the path degree is not defined on the quotient. For matrices with a part of degree zero, we only have the perturbation theory of Section~\ref{sec:beyond-radical} and the power series results of Chapter~\ref{chap:infinitesimals}. A general encoding theorem for path $*$-algebras with relations is open.
  \item \emph{Involutions and square roots.} The completed localisation carries the involution, and for every quiver of the catalogue it is a product of simple Artinian rings that the involution permutes. Beyond $\Lambda$, every doubly-signed quiver, including one with a vertex on two distinct simple oriented cycles, embeds after complexification, compatibly with the involution and the conjugation, in a finite product of simple Artinian rings that they permute (Section~\ref{sec:embedding}). It is open whether the $*$-square roots that a Jacobi-type method needs exist there, and whether the graded rule of Section~\ref{sec:lambda-sqrt} finds a $*$-square root of every element $z^*z$ with $z$ homogeneous.
  \item \emph{Reducing results.} Theorem~\ref{thm:reduction} brings results back from the completed localisation whenever the computed factors contain no inverse letters. Results computed in the free field, which quivers with a vertex on two distinct simple oriented cycles require, have no such reduction yet.
  \item \emph{Wild quivers.} Our explanation of the behaviour of generic algorithms on wild quivers (Section~\ref{sec:wild-behaviour}) is supported by examples, not by a general theorem: for two loops we showed that the quotients needed by a Jacobi-type method do not exist in the completed localisation, but not whether the method can be completed in a skew field, or what it returns. Arithmetic in the free field, which these cases require, is also expensive.
\end{itemize}

\section{Future Research Directions}
\label{sec:future}
\begin{itemize}
    \item Develop the library of Chapter~\ref{chap:implementation} into a tool or compiler that applies these transformations automatically to existing code bases, for example in Fortran, C++ or Python with NumPy.
    \item Explore the connection between the representation type of a quiver and the numerical stability and convergence of the corresponding generic algorithms.
    \item Extend the framework to operator algebras along the lines of Section~\ref{sec:operator}, beginning with sequences of matrices in the sense of~\cite{HagenRochSilbermann2001}.
    \item Extend the visualisation of Appendix~\ref{chap:gutin-2022} to $3\times3$ matrices; we are not aware of an analogous construction.
\end{itemize}

\section{The Role of the Visualisation Appendix}
\label{sec:appendix-summary}
Appendix~\ref{chap:gutin-2022} draws the dynamics of the QR algorithm on $2\times2$ real matrices in the hyperbolic plane of inner products on $\RR^2$, and defines the geometric objects that it uses. Two things connect it with the rest of the thesis.

The first is involutions. A point of that plane is an involution of $M_2(\RR)$, and the picture of $M$ shows for which of these involutions the Hermitian part $\frac12(M+M^*)$ is definite. The QR algorithm can be run under each of them, and it turns the picture about the corresponding point (Proposition~\ref{prop:qr-lr}). When that point moves to the boundary along the weighted inner products~\eqref{eq:weighted-ip}, the QR algorithm becomes the LR algorithm, which is the $2\times2$ case of Proposition~\ref{prop:gs-limit} and Theorem~\ref{thm:bruhat}; by Proposition~\ref{prop:2R-algorithms} the LR algorithm is also the QR algorithm over $\twoR$, two steps at a time. So the appendix shows, for $2\times2$ matrices, the QR algorithm changing into the LR algorithm as the involution degenerates, where the main text obtains the same change by changing the scalars.

The second is the two-dimensional $*$-algebras $\CC$, the dual numbers and $\twoR$. A real $2\times2$ matrix that is not scalar generates an algebra isomorphic to one of them, and this decides the kind of its picture (Proposition~\ref{prop:shapes} and Table~\ref{tab:kinds}); Kisil~\cite{Kisil2012} develops the geometry of M\"obius transformations from these three algebras. The projective lines over them are the M\"obius, Laguerre and Minkowski planes, and planar Lie sphere geometry contains all three~\cite{Benz1973,Yaglom1968-ComplexNumbers,cecil1992lie}. The group of Lie sphere geometry has identity component $\mathrm{PSp}(4,\RR)\cong\SO^+(3,2)$. This is the projective unitary group of the $2\times2$ matrices over $M_2(\RR)$ with its symplectic involution. The points of Lie sphere geometry are the Lagrangian planes in $\RR^4$; they form a three-dimensional part of the projective line over $M_2(\RR)$, which is the Grassmannian of planes in $\RR^4$.
%======================================================================
%   APPENDIX A: VISUALISATION
%======================================================================
%======================================================================
%   BIBLIOGRAPHY (before the appendices)
%======================================================================
\cleardoublepage
\phantomsection
\addcontentsline{toc}{chapter}{Bibliography}
\bibliographystyle{plainnat}
\bibliography{thesis}

% Bibliography recovered from the author's matching thesis PDF.
% Source: thesisonelasttime_revised (20).pdf, pages 200--210.
% Exact cite.<key> destinations supply keys; printed labels supply ordering.
% No bibliographic records were added or updated.
\begin{thebibliography}{99}
\providecommand{\url}[1]{\texttt{#1}}
\providecommand{\doi}[1]{\begingroup\urlstyle{rm}\nolinkurl{#1}\endgroup}

% Original reference [1]
\bibitem{AbramsAraSilesMolina2017}
Gene Abrams, Pere Ara, and Mercedes Siles Molina. \emph{Leavitt Path Algebras}, volume 2191 of \emph{Lecture Notes in Mathematics}. Springer, London, 2017. doi: \doi{10.1007/978-1-4471-7344-1}.

% Original reference [2]
\bibitem{AmitsurLevitzki1950}
Shimshon A. Amitsur and Jacob Levitzki. Minimal identities for algebras. \emph{Proceedings of the American Mathematical Society}, 1(4):449--463, 1950. doi: \doi{10.1090/S0002-9939-1950-0036751-9}.

% Original reference [3]
\bibitem{AssemSimsonSkowronski2006}
Ibrahim Assem, Daniel Simson, and Andrzej Skowro\'{n}ski. \emph{Elements of the Representation Theory of Associative Algebras, Volume 1: Techniques of Representation} \emph{Theory}, volume 65 of \emph{London Mathematical Society Student Texts}. Cambridge University Press, Cambridge, 2006. doi: \doi{10.1017/CBO9780511614309}.

% Original reference [4]
\bibitem{AuslanderReitenSmalo1995}
Maurice Auslander, Idun Reiten, and Sverre O. Smal\o{}. \emph{Representation Theory of} \emph{Artin Algebras}, volume 36 of \emph{Cambridge Studies in Advanced Mathematics}. Cambridge University Press, Cambridge, 1995. doi: \doi{10.1017/CBO9780511623608}.

% Original reference [5]
\bibitem{BangJensenGutin2009}
J\o{}rgen Bang-Jensen and Gregory Gutin. \emph{Digraphs: Theory, Algorithms and Applications}. Springer Monographs in Mathematics. Springer, London, second edition, 2009. doi: \doi{10.1007/978-1-84800-998-1}.

% Original reference [6]
\bibitem{goodshift2022}
Jess Banks, Jorge Garza-Vargas, and Nikhil Srivastava. Global convergence of Hessenberg shifted QR I: Exact arithmetic. \emph{Foundations of Computational Mathematics}, 25(5):1561--1594, 2025. doi: \doi{10.1007/s10208-024-09658-7}.

% Original reference [7]
\bibitem{BasuPollackRoy2006}
Saugata Basu, Richard Pollack, and Marie-Fran\c{c}oise Roy. \emph{Algorithms in Real Algebraic Geometry}, volume 10 of \emph{Algorithms and Computation in Mathematics}. Springer, Berlin, Heidelberg, second edition, 2006. doi: \doi{10.1007/3-540-33099-2}.

\bibitem{Beardon1983}
Alan F. Beardon. \emph{The Geometry of Discrete Groups}, volume 91 of \emph{Graduate Texts in Mathematics}. Springer, New York, 1983. doi: \doi{10.1007/978-1-4612-1146-4}.

% Original reference [8]
\bibitem{Benz1973}
Walter Benz. \emph{Vorlesungen \"{u}ber Geometrie der Algebren: Geometrien von M\"{o}bius,} \emph{Laguerre-Lie, Minkowski in einheitlicher und grundlagengeometrischer Behandlung}, volume 197 of \emph{Grundlehren der mathematischen Wissenschaften}. Springer, Berlin, Heidelberg, 1973. doi: \doi{10.1007/978-3-642-88670-6}.

% Original reference [9]
\bibitem{BochnakCosteRoy1998}
Jacek Bochnak, Michel Coste, and Marie-Fran\c{c}oise Roy. \emph{Real Algebraic Geometry}, volume 36 of \emph{Ergebnisse der Mathematik und ihrer Grenzgebiete (3)}. Springer, Berlin, Heidelberg, 1998. doi: \doi{10.1007/978-3-662-03718-8}.

% Original reference [10]
\bibitem{BoettcherSilbermann1999}
Albrecht B\"{o}ttcher and Bernd Silbermann. \emph{Introduction to Large Truncated Toeplitz} \emph{Matrices}. Universitext. Springer, New York, 1999. doi: \doi{10.1007/978-1-4612-1426-7}.

% Original reference [11]
\bibitem{BunseGerstnerByersMehrmannNichols1991}
Angelika Bunse-Gerstner, Ralph Byers, Volker Mehrmann, and Nancy K. Nichols. Numerical computation of an analytic singular value decomposition of a matrix valued function. \emph{Numerische Mathematik}, 60(1):1--39, 1991. doi: \doi{10.1007/BF01385712}.

% Original reference [12]
\bibitem{BunseGerstnerByersMehrmann1992}
Angelika Bunse-Gerstner, Ralph Byers, and Volker Mehrmann. A chart of numerical methods for structured eigenvalue problems. \emph{SIAM Journal on Matrix Analysis and} \emph{Applications}, 13(2):419--453, 1992. doi: \doi{10.1137/0613028}.

% Original reference [13]
\bibitem{buurema1958geometric}
Hendrik Jan Buurema. \emph{A Geometric Proof of Convergence for the QR Method}. PhD thesis, Rijksuniversiteit te Groningen, Groningen, 1970.

% Original reference [14]
\bibitem{Carre1971}
Bernard A. Carr\'{e}. An algebraic approach to network routing problems. \emph{IMA Journal} \emph{of Applied Mathematics}, 7(3):273--294, 1971. doi: \doi{10.1093/imamat/7.3.273}.

% Original reference [15]
\bibitem{cecil1992lie}
Thomas E. Cecil. \emph{Lie Sphere Geometry: With Applications to Submanifolds}. Universitext. Springer, New York, 1992. doi: \doi{10.1007/978-1-4757-4096-7}.

% Original reference [16]
\bibitem{Cohn2006}
Paul M. Cohn. \emph{Free Ideal Rings and Localization in General Rings}, volume 3 of \emph{New Mathematical Monographs}. Cambridge University Press, Cambridge, 2006. doi: \doi{10.1017/CBO9780511542794}.

% Original reference [17]
\bibitem{CohnReutenauer1999}
Paul M. Cohn and Christophe Reutenauer. On the construction of the free field. \emph{International Journal of Algebra and Computation}, 9(3--4):307--323, 1999. doi: \doi{10.1142/S0218196799000205}.

% Original reference [18]
\bibitem{CoxHigham1999}
Anthony J. Cox and Nicholas J. Higham. Row-wise backward stable elimination methods for the equality constrained least squares problem. \emph{SIAM Journal on Matrix} \emph{Analysis and Applications}, 21(1):313--326, 1999. doi: \doi{10.1137/S0895479898335957}.

% Original reference [19]
\bibitem{Davidson1996}
Kenneth R. Davidson. $C$${}^{*}$\emph{-Algebras by Example}, volume 6 of \emph{Fields Institute Monographs}. American Mathematical Society, Providence, RI, 1996. doi: \doi{10.1090/fim/006}.

% Original reference [20]
\bibitem{DemmelKahan1990}
James Demmel and William Kahan. Accurate singular values of bidiagonal matrices. \emph{SIAM Journal on Scientific and Statistical Computing}, 11(5):873--912, 1990. doi: \doi{10.1137/0911052}.

% Original reference [21]
\bibitem{DerksenWeyman2002}
Harm Derksen and Jerzy Weyman. Generalized quivers associated to reductive groups. \emph{Colloquium Mathematicum}, 94(2):151--173, 2002. doi: \doi{10.4064/cm94-2-1}.

% Original reference [22]
\bibitem{DerksenWeyman2017}
Harm Derksen and Jerzy Weyman. \emph{An Introduction to Quiver Representations}, volume 184 of \emph{Graduate Studies in Mathematics}. American Mathematical Society, Providence, RI, 2017. doi: \doi{10.1090/gsm/184}.

% Original reference [23]
\bibitem{Dijkstra1959}
Edsger W. Dijkstra. A note on two problems in connexion with graphs. \emph{Numerische} \emph{Mathematik}, 1:269--271, 1959. doi: \doi{10.1007/BF01386390}.

% Original reference [24]
\bibitem{DlabRingel1976}
Vlastimil Dlab and Claus Michael Ringel. Indecomposable representations of graphs and algebras. \emph{Memoirs of the American Mathematical Society}, 6(173), 1976. doi: \doi{10.1090/memo/0173}.

% Original reference [25]
\bibitem{Drozd1980}
Yuriy A. Drozd. Tame and wild matrix problems. In Vlastimil Dlab and Peter Gabriel, editors, \emph{Representation Theory II: Proceedings of the Second International} \emph{Conference on Representations of Algebras (Ottawa, 1979)}, volume 832 of \emph{Lecture} \emph{Notes in Mathematics}, pages 242--258. Springer, Berlin, Heidelberg, 1980. doi: \doi{10.1007/BFb0088467}.

% Original reference [26]
\bibitem{DummitFoote2003}
David S. Dummit and Richard M. Foote. \emph{Abstract Algebra}. John Wiley \& Sons, Hoboken, NJ, third edition, 2004.

% Original reference [27]
\bibitem{EdelmanJeong2023}
Alan Edelman and Sungwoo Jeong. Fifty three matrix factorizations: A systematic approach. \emph{SIAM Journal on Matrix Analysis and Applications}, 44(2):415--480, 2023. doi: \doi{10.1137/21M1416035}.

% Original reference [28]
\bibitem{FassbenderMackeyMackeyXu1999}
Heike Fa\ss{}bender, D. Steven Mackey, Niloufer Mackey, and Hongguo Xu. Hamiltonian square roots of skew-Hamiltonian matrices. \emph{Linear Algebra and its Applications}, 287(1--3):125--159, 1999. doi: \doi{10.1016/S0024-3795(98)10137-4}.

% Original reference [29]
\bibitem{Francis1961}
John G. F. Francis. The QR transformation: A unitary analogue to the LR transformation---Part 1. \emph{The Computer Journal}, 4(3):265--271, 1961. doi: \doi{10.1093/comjnl/4.3.265}.

% Original reference [30]
\bibitem{Francis1962}
John G. F. Francis. The QR transformation---Part 2. \emph{The Computer Journal}, 4(4):332--345, 1962. doi: \doi{10.1093/comjnl/4.4.332}.

% Original reference [31]
\bibitem{Gabriel1972}
Peter Gabriel. Unzerlegbare Darstellungen I. \emph{Manuscripta Mathematica}, 6(1):71--103, 1972. doi: \doi{10.1007/BF01298413}.

% Original reference [32]
\bibitem{GabrielZisman1967}
Peter Gabriel and Michel Zisman. \emph{Calculus of Fractions and Homotopy Theory}, volume 35 of \emph{Ergebnisse der Mathematik und ihrer Grenzgebiete}. Springer, Berlin, Heidelberg, 1967. doi: \doi{10.1007/978-3-642-85844-4}.

% Original reference [33]
\bibitem{GolubKahan1965}
Gene H. Golub and William Kahan. Calculating the singular values and pseudo-inverse of a matrix. \emph{Journal of the Society for Industrial and Applied Mathematics,} \emph{Series B: Numerical Analysis}, 2(2):205--224, 1965. doi: \doi{10.1137/0702016}.

% Original reference [34]
\bibitem{GolubVanLoan2013}
Gene H. Golub and Charles F. Van Loan. \emph{Matrix Computations}. Johns Hopkins Studies in the Mathematical Sciences. Johns Hopkins University Press, Baltimore, MD, fourth edition, 2013.

% Original reference [35]
\bibitem{GondranMinoux2008}
Michel Gondran and Michel Minoux. \emph{Graphs, Dioids and Semirings: New Models and} \emph{Algorithms}, volume 41 of \emph{Operations Research/Computer Science Interfaces Series}. Springer, New York, 2008. doi: \doi{10.1007/978-0-387-75450-5}.

% Original reference [36]
\bibitem{GriewankWalther2008}
Andreas Griewank and Andrea Walther. \emph{Evaluating Derivatives: Principles and} \emph{Techniques of Algorithmic Differentiation}. Society for Industrial and Applied Mathematics, Philadelphia, PA, second edition, 2008. doi: \doi{10.1137/1.9780898717761}.

% Original reference [37]
\bibitem{Gutin2021AllRings}
Ran Gutin. Towards a singular value decomposition and spectral theory for all rings. arXiv preprint arXiv:2112.10247, 2021.

% Original reference [38]
\bibitem{Gutin2021DoubleNumbers}
Ran Gutin. Matrix decompositions over the double numbers. arXiv preprint arXiv:2105.08047, 2021.

% Original reference [39]
\bibitem{Gutin2022DualNumbers}
Ran Gutin. Generalizations of singular value decomposition to dual-numbered matrices. \emph{Linear and Multilinear Algebra}, 70(20):5107--5114, 2022. doi: \doi{10.1080/03081087.2021.1903830}.

% Original reference [40]
\bibitem{Gutin2022PSD}
Ran Gutin. A visualisation for conveying the dynamics of iterative eigenvalue algorithms over PSD matrices. arXiv preprint arXiv:2204.00383, 2022.

% Original reference [41]
\bibitem{Gutin2022Visualisation}
Ran Gutin. An interactive visualisation for all $2 \times 2$ real matrices, with applications to conveying the dynamics of iterative eigenvalue algorithms. arXiv preprint arXiv:2211.02075, 2022.

% Original reference [42]
\bibitem{Gutin2023Conjecture}
Ran Gutin. A conjecture concerning $\ast$-algebras that unifies some matrix decompositions. arXiv preprint arXiv:2308.04513, 2023.

% Original reference [43]
\bibitem{Gutin2023DoubleComplex}
Ran Gutin. An analogue of the relationship between SVD and pseudoinverse over double-complex matrices. \emph{Linear and Multilinear Algebra}, 71(11):1775--1793, 2023. doi: \doi{10.1080/03081087.2022.2076798}.

% Original reference [44]
\bibitem{Gutin2024Clifford}
Ran Gutin. Unitary canonical forms over Clifford algebras, and an observed unification of some real-matrix decompositions. \emph{Linear and Multilinear Algebra}, 72(14):2288--2311, 2024. doi: \doi{10.1080/03081087.2023.2252979}.

% Original reference [45]
\bibitem{HagenRochSilbermann2001}
Roland Hagen, Steffen Roch, and Bernd Silbermann. $C$${}^{*}$\emph{-Algebras and Numerical} \emph{Analysis}, volume 236 of \emph{Pure and Applied Mathematics}. Marcel Dekker, New York, 2001.

% Original reference [46]
\bibitem{Harris2020NumPy}
Charles R. Harris, K. Jarrod Millman, St\'{e}fan J. van der Walt, Ralf Gommers, et al. Array programming with NumPy. \emph{Nature}, 585(7825):357--362, 2020. doi: \doi{10.1038/s41586-020-2649-2}.

% Original reference [47]
\bibitem{Helgason1978}
Sigurdur Helgason. \emph{Differential Geometry, Lie Groups, and Symmetric Spaces}, volume 80 of \emph{Pure and Applied Mathematics}. Academic Press, New York, 1978.

% Original reference [48]
\bibitem{Higham2002}
Nicholas J. Higham. \emph{Accuracy and Stability of Numerical Algorithms}. Society for Industrial and Applied Mathematics, Philadelphia, PA, second edition, 2002. doi: \doi{10.1137/1.9780898718027}.

% Original reference [49]
\bibitem{Higham2008}
Nicholas J. Higham. \emph{Functions of Matrices: Theory and Computation}, volume 104 of \emph{Other Titles in Applied Mathematics}. Society for Industrial and Applied Mathematics, Philadelphia, PA, 2008. doi: \doi{10.1137/1.9780898717778}.

% Original reference [50]
\bibitem{HornJohnson2013}
Roger A. Horn and Charles R. Johnson. \emph{Matrix Analysis}. Cambridge University Press, Cambridge, second edition, 2013. doi: \doi{10.1017/CBO9781139020411}.

% Original reference [51]
\bibitem{HornSergeichuk2006}
Roger A. Horn and Vladimir V. Sergeichuk. Canonical forms for complex matrix congruence and ${}^{*}$congruence. \emph{Linear Algebra and its Applications}, 416(2--3):1010--1032, 2006. doi: \doi{10.1016/j.laa.2006.01.005}.

% Original reference [52]
\bibitem{Humphreys1975}
James E. Humphreys. \emph{Linear Algebraic Groups}, volume 21 of \emph{Graduate Texts in} \emph{Mathematics}. Springer-Verlag, New York, 1975.

% Original reference [53]
\bibitem{Jacobson1989}
Nathan Jacobson. \emph{Basic Algebra II}. W. H. Freeman, New York, second edition, 1989.

% Original reference [54]
\bibitem{Kac1980}
Victor G. Kac. Infinite root systems, representations of graphs and invariant theory. \emph{Inventiones Mathematicae}, 56(1):57--92, 1980. doi: \doi{10.1007/BF01403155}.

% Original reference [55]
\bibitem{Kato1995}
Tosio Kato. \emph{Perturbation Theory for Linear Operators}. Classics in Mathematics. Springer, Berlin, Heidelberg, second edition, 1995. doi: \doi{10.1007/978-3-642-66282-9}.

\bibitem{Kisil2012}
Vladimir V. Kisil. \emph{Geometry of M\"{o}bius Transformations: Elliptic, Parabolic and Hyperbolic Actions of $\mathrm{SL}_2(\mathbb{R})$}. Imperial College Press, London, 2012.

% Original reference [56]
\bibitem{KleinsteuberHelmkeHuper2004}
Martin Kleinsteuber, Uwe Helmke, and Knut H\"{u}per. Jacobi's algorithm on compact Lie algebras. \emph{SIAM Journal on Matrix Analysis and Applications}, 26(1):42--69, 2004. doi: \doi{10.1137/S0895479802420069}.

% Original reference [57]
\bibitem{Knus1991}
Max-Albert Knus. \emph{Quadratic and Hermitian Forms over Rings}, volume 294 of \emph{Grundlehren der mathematischen Wissenschaften}. Springer, Berlin, Heidelberg, 1991. doi: \doi{10.1007/978-3-642-75401-2}.

% Original reference [58]
\bibitem{KMRT1998}
Max-Albert Knus, Alexander Merkurjev, Markus Rost, and Jean-Pierre Tignol. \emph{The Book of Involutions}, volume 44 of \emph{American Mathematical Society Colloquium Publications}. American Mathematical Society, Providence, RI, 1998. doi: \doi{10.1090/coll/044}.

% Original reference [59]
\bibitem{Kock2006}
Anders Kock. \emph{Synthetic Differential Geometry}, volume 333 of \emph{London Mathematical} \emph{Society Lecture Note Series}. Cambridge University Press, Cambridge, second edition, 2006. doi: \doi{10.1017/CBO9780511550812}.

% Original reference [60]
\bibitem{Lam2001}
Tsit-Yuen Lam. \emph{A First Course in Noncommutative Rings}, volume 131 of \emph{Graduate} \emph{Texts in Mathematics}. Springer, New York, second edition, 2001. doi: \doi{10.1007/978-1-4419-8616-0}.

% Original reference [61]
\bibitem{Lam2005}
Tsit-Yuen Lam. \emph{Introduction to Quadratic Forms over Fields}, volume 67 of \emph{Graduate} \emph{Studies in Mathematics}. American Mathematical Society, Providence, RI, 2005. doi: \doi{10.1090/gsm/067}.

% Original reference [62]
\bibitem{Lance1995}
E. Christopher Lance. \emph{Hilbert} $C$${}^{*}$\emph{-Modules: A Toolkit for Operator Algebraists}, volume 210 of \emph{London Mathematical Society Lecture Note Series}. Cambridge University Press, Cambridge, 1995. doi: \doi{10.1017/CBO9780511526206}.

% Original reference [63]
\bibitem{LeiteSaldanhaTomei2013IMRN}
Ricardo S. Leite, Nicolau C. Saldanha, and Carlos Tomei. Dynamics of the symmetric eigenvalue problem with shift strategies. \emph{International Mathematics Research Notices}, 2013(19):4382--4412, 2013. doi: \doi{10.1093/imrn/rns186}.

% Original reference [64]
\bibitem{Lewin1974}
Jacques Lewin. Fields of fractions for group algebras of free groups. \emph{Transactions of the American Mathematical Society}, 192:339--346, 1974. doi: \doi{10.1090/S0002-9947-1974-0338055-4}.

% Original reference [65]
\bibitem{Loos1969}
Ottmar Loos. \emph{Symmetric Spaces. I: General Theory}. Mathematics Lecture Note Series. W. A. Benjamin, New York, 1969.

% Original reference [66]
\bibitem{Lorand2020}
Jonathan Lorand. \emph{Duality Involutions, Representations, and Geometry}. PhD thesis, Universit\"{a}t Z\"{u}rich, Z\"{u}rich, 2020. URL \url{https://www.zora.uzh.ch/id/eprint/197366/}.

% Original reference [67]
\bibitem{MackeyMackeyTisseur2006}
D. Steven Mackey, Niloufer Mackey, and Fran\c{c}oise Tisseur. Structured factorizations in scalar product spaces. \emph{SIAM Journal on Matrix Analysis and Applications}, 27(3):821--850, 2006. doi: \doi{10.1137/040619363}.

% Original reference [68]
\bibitem{MaclaganSturmfels2015}
Diane Maclagan and Bernd Sturmfels. \emph{Introduction to Tropical Geometry}, volume 161 of \emph{Graduate Studies in Mathematics}. American Mathematical Society, Providence, RI, 2015. doi: \doi{10.1090/gsm/161}.

% Original reference [69]
\bibitem{Marker2002}
David Marker. \emph{Model Theory: An Introduction}, volume 217 of \emph{Graduate Texts in} \emph{Mathematics}. Springer, New York, 2002. doi: \doi{10.1007/b98860}.

% Original reference [70]
\bibitem{Pierce1982}
Richard S. Pierce. \emph{Associative Algebras}, volume 88 of \emph{Graduate Texts in Mathematics}. Springer, New York, 1982. doi: \doi{10.1007/978-1-4757-0163-0}.

% Original reference [71]
\bibitem{QuebbemannScharlauSchulte1979}
Heinz-Georg Quebbemann, Winfried Scharlau, and Manfred Schulte. Quadratic and Hermitian forms in additive and abelian categories. \emph{Journal of Algebra}, 59(2):264--289, 1979. doi: \doi{10.1016/0021-8693(79)90126-1}.

% Original reference [72]
\bibitem{Raeburn2005}
Iain Raeburn. \emph{Graph Algebras}, volume 103 of \emph{CBMS Regional Conference Series in} \emph{Mathematics}. American Mathematical Society, Providence, RI, 2005. doi: \doi{10.1090/cbms/103}.

% Original reference [73]
\bibitem{Ribenboim1999}
Paulo Ribenboim. \emph{The Theory of Classical Valuations}. Springer Monographs in Mathematics. Springer, New York, 1999. doi: \doi{10.1007/978-1-4612-0551-7}.

% Original reference [74]
\bibitem{Robinson1966}
Abraham Robinson. \emph{Non-standard Analysis}. Studies in Logic and the Foundations of Mathematics. North-Holland, Amsterdam, 1966.

% Original reference [75]
\bibitem{Rutishauser1958}
Heinz Rutishauser. Solution of eigenvalue problems with the LR-transformation. In \emph{Further Contributions to the Solution of Simultaneous Linear Equations and the} \emph{Determination of Eigenvalues}, volume 49 of \emph{National Bureau of Standards Applied} \emph{Mathematics Series}, pages 47--81. National Bureau of Standards, Washington, DC, 1958.

% Original reference [76]
\bibitem{Scharlau1985}
Winfried Scharlau. \emph{Quadratic and Hermitian Forms}, volume 270 of \emph{Grundlehren der} \emph{mathematischen Wissenschaften}. Springer, Berlin, Heidelberg, 1985. doi: \doi{10.1007/978-3-642-69971-9}.

% Original reference [77]
\bibitem{Schiffler2014}
Ralf Schiffler. \emph{Quiver Representations}. CMS Books in Mathematics. Springer, Cham, 2014. doi: \doi{10.1007/978-3-319-09204-1}.

% Original reference [78]
\bibitem{Schofield1985}
Aidan H. Schofield. \emph{Representations of Rings over Skew Fields}, volume 92 of \emph{London} \emph{Mathematical Society Lecture Note Series}. Cambridge University Press, Cambridge, 1985. doi: \doi{10.1017/CBO9780511661914}.

% Original reference [79]
\bibitem{Sergeichuk1988}
Vladimir V. Sergeichuk. Classification problems for systems of forms and linear mappings. \emph{Mathematics of the USSR-Izvestiya}, 31(3):481--501, 1988. doi: \doi{10.1070/IM1988v031n03ABEH001086}.

% Original reference [80]
\bibitem{Shmelkin2006}
Dmitri A. Shmelkin. Signed quivers, symmetric quivers and root systems. \emph{Journal of} \emph{the London Mathematical Society. Second Series}, 73(3):586--606, 2006. doi: \doi{10.1112/S0024610706022812}.

% Original reference [81]
\bibitem{Taft1957}
Earl J. Taft. Invariant Wedderburn factors. \emph{Illinois Journal of Mathematics}, 1(4):565--573, 1957. doi: \doi{10.1215/ijm/1255380679}.

% Original reference [82]
\bibitem{Takagi1925}
Teiji Takagi. On an algebraic problem related to an analytic theorem of Carath\'{e}odory and Fej\'{e}r and on an allied theorem of Landau. \emph{Japanese Journal of Mathematics}, 1:83--93, 1924. doi: \doi{10.4099/jjm1924.1.0_83}.

% Original reference [83]
\bibitem{Tao2007}
Terence Tao. Ultrafilters, nonstandard analysis, and epsilon management. \emph{What's New} (blog), \url{https://terrytao.wordpress.com/2007/06/25/ultrafilters-nonstandard-analysis-and-epsilon-management/}, 2007.

% Original reference [84]
\bibitem{Tarski1951}
Alfred Tarski. \emph{A Decision Method for Elementary Algebra and Geometry}. University of California Press, Berkeley, CA, second edition, 1951.

% Original reference [85]
\bibitem{VanLoan1984}
Charles F. Van Loan. A symplectic method for approximating all the eigenvalues of a Hamiltonian matrix. \emph{Linear Algebra and its Applications}, 61:233--251, 1984. doi: \doi{10.1016/0024-3795(84)90034-X}.

% Original reference [86]
\bibitem{Virtanen2020SciPy}
Pauli Virtanen, Ralf Gommers, Travis E. Oliphant, Matt Haberland, et al. SciPy 1.0: Fundamental algorithms for scientific computing in Python. \emph{Nature Methods}, 17(3):261--272, 2020. doi: \doi{10.1038/s41592-019-0686-2}.

% Original reference [87]
\bibitem{Wall1976}
Charles T. C. Wall. Classification of Hermitian forms. VI. Group rings. \emph{Annals of} \emph{Mathematics. Second Series}, 103(1):1--80, 1976. doi: \doi{10.2307/1971019}.

% Original reference [88]
\bibitem{watkins2008qr}
David S. Watkins. The QR algorithm revisited. \emph{SIAM Review}, 50(1):133--145, 2008. doi: \doi{10.1137/060659454}.

\bibitem{WatkinsElsner1991}
David S. Watkins and Ludwig Elsner. Convergence of algorithms of decomposition type for the eigenvalue problem. \emph{Linear Algebra and its Applications}, 143:19--47, 1991. doi: \doi{10.1016/0024-3795(91)90004-G}.

% Original reference [89]
\bibitem{Weihrauch2000}
Klaus Weihrauch. \emph{Computable Analysis: An Introduction}. Texts in Theoretical Computer Science. An EATCS Series. Springer, Berlin, Heidelberg, 2000. doi: \doi{10.1007/978-3-642-56999-9}.

% Original reference [90]
\bibitem{Weyl1939}
Hermann Weyl. \emph{The Classical Groups: Their Invariants and Representations}, volume 1 of \emph{Princeton Mathematical Series}. Princeton University Press, Princeton, NJ, 1939.

% Original reference [91]
\bibitem{Xu1998}
Hongguo Xu. The relation between the QR and LR algorithms. \emph{SIAM Journal on Matrix Analysis and Applications}, 19(2):551--555, 1998. doi: \doi{10.1137/S0895479896299937}.

% Original reference [92]
\bibitem{Yaglom1968-ComplexNumbers}
Isaak M. Yaglom. \emph{Complex Numbers in Geometry}. Academic Press, New York and London, 1968.

% Original reference [93]
\bibitem{Yaglom1979-SimpleNonEuclid}
Isaak M. Yaglom. \emph{A Simple Non-Euclidean Geometry and Its Physical Basis: An Elementary Account of Galilean Geometry and the Galilean Principle of Relativity}. Heidelberg Science Library. Springer, New York, 1979. doi: \doi{10.1007/978-1-4612-6135-3}.

% Original reference [94]
\bibitem{Youla1961}
Dante C. Youla. A normal form for a matrix under the unitary congruence group. \emph{Canadian Journal of Mathematics}, 13:694--704, 1961. doi: \doi{10.4153/CJM-1961-059-8}.

% Original reference [95]
\bibitem{Zelazko1973}
Wies\l{}aw \.{Z}elazko. \emph{Banach Algebras}. Elsevier, Amsterdam, 1973.

% Original reference [96]
\bibitem{Zhang1997}
Fuzhen Zhang. Quaternions and matrices of quaternions. \emph{Linear Algebra and its} \emph{Applications}, 251:21--57, 1997. doi: \doi{10.1016/0024-3795(95)00543-9}.

\end{thebibliography}

\appendix
\chapter{Visualising QR Dynamics on \texorpdfstring{$2\times2$}{2x2} Real Matrices}
\label{chap:gutin-2022}
\addtocontents{toc}{\protect\setcounter{tocdepth}{1}}

\paragraph{Provenance.} This appendix is based on a preprint written by the author during his doctoral studies~\cite{Gutin2022Visualisation} (arXiv:\allowbreak 2211.02075); its warm-up visualisation for positive semidefinite matrices first appeared in~\cite{Gutin2022PSD}. The pictures are those of the preprint, redrawn. The way in which they are explained and derived is new. The preprint obtained them from Lie sphere geometry; here they are derived from the space of inner products on $\RR^2$ (Section~\ref{sec:disc}), and Lie sphere geometry remains as a recipe for drawing (Section~\ref{sec:recipe}). The new derivation proves the formula for the Lie quadric representation (Lemma~\ref{lem:pencil}), corrects the preprint's formulas for the envelopes and its description of the QR step in visualisation v1.5, and makes precise its statements about rotations and about the speed of convergence. Section~\ref{sec:qr-lr}, on the LR algorithm, is also new. No result of the thesis depends on this appendix; Section~\ref{sec:qr-lr} says what it illustrates of the thesis.

\section{Overview}
This appendix shows what the QR algorithm does to a $2\times2$ real matrix, in pictures. It is meant for self-study or teaching, and it is written so that the pictures come first: Sections~\ref{sec:directions} to~\ref{sec:qr-nonneg} need only eigenvectors, the QR factorisation and the power method, and every proof is postponed to Section~\ref{sec:recipe}. One matrix is followed throughout,
\[ M=\bmat{3&1\\2&2},\qquad\text{with eigenvalues $4$ and $1$ and eigenvectors $(1,1)^T$ and $(1,-2)^T$,} \]
which we call \emph{the running example}; its numbers can be checked by hand. Section~\ref{sec:exercises} has exercises with answers.

There are two pictures, called v1 and v2 as in the preprint. Visualisation v1 (Section~\ref{sec:v1}) is a warm-up for positive semidefinite matrices, which are drawn as ellipses. Visualisation v2 covers every invertible $2\times2$ real matrix: a pair of circles when the determinant is positive, which is the main case, and an arc crossed by arrows when it is negative (Section~\ref{sec:qr-negdet}). It is built in three steps.
\begin{itemizedot}
  \item The directions of the plane form a circle, and a matrix moves the points of this circle (Section~\ref{sec:directions}).
  \item The QR algorithm is the power method on this circle, watched from a frame that turns with it (Section~\ref{sec:turning}).
  \item Joining each direction to its image by an arrow gives a picture of the matrix itself: the arrows touch two circles, and these two circles are the picture (Section~\ref{sec:picture}).
\end{itemizedot}
Section~\ref{sec:qr-nonneg} then watches the QR algorithm in this picture: why it converges, when it is fast or slow, where it fails, and what a shift does.

Up to that point the inside of the disc is only drawing paper. Section~\ref{sec:disc} explains what it is: its points are the inner products on $\RR^2$, with the standard inner product at the centre. Each inner product gives an adjoint $X\mapsto X^*$ of $2\times2$ matrices, and the two circles of a matrix $M$ mark where its Hermitian part $\frac12(M+M^*)$ stops being definite (Theorem~\ref{thm:picture}). With this, many properties of $M$ can be read off from the picture:
\begin{itemizedot}
  \item its Jordan normal form;
  \item which matrices it is orthogonally similar to;
  \item its condition number;
  \item whether it is upper triangular, lower triangular, diagonal, symmetric or orthogonal;
  \item a geometric construction for the matrix-vector product $Mv$;
  \item a geometric construction for the matrix-matrix product $MK$.
\end{itemizedot}
Table~\ref{tab:reading} in Section~\ref{sec:reading} says how each of these is read. Section~\ref{sec:qr-lr} changes the inner product, which turns the QR algorithm into the LR algorithm; this is where the appendix meets the main text of the thesis. Section~\ref{sec:qr-negdet} treats matrices with negative determinant.

Our aim is to convey the dynamics of the \emph{unshifted} QR algorithm, $M=QR\mapsto RQ$~\cite{Francis1961,Francis1962}, with as few formulas as possible. Pedagogical treatments of the QR algorithm exist~\cite{watkins2008qr,buurema1958geometric}, but they are more formula-driven and less visual than ours. For $3 \times 3$ matrices we have not found anything analogous. There are ways to extend v1 to all $2 \times 2$ real matrices, for instance by drawing the polar decomposition of $M$, but then most of the properties listed above cannot be read off; such a picture might still be useful for showing the dynamics of algorithms that compute SVDs or polar decompositions, and we do not pursue this here.

Every figure in this appendix is drawn from the formulas given here, and Table~\ref{tab:matrices} at the end lists the matrices that the figures show.

\section{Visualisation v1: Ellipses (Warm-Up)}
\label{sec:v1}
\begin{figure}[ht]
\centering
\includegraphics[width=\linewidth]{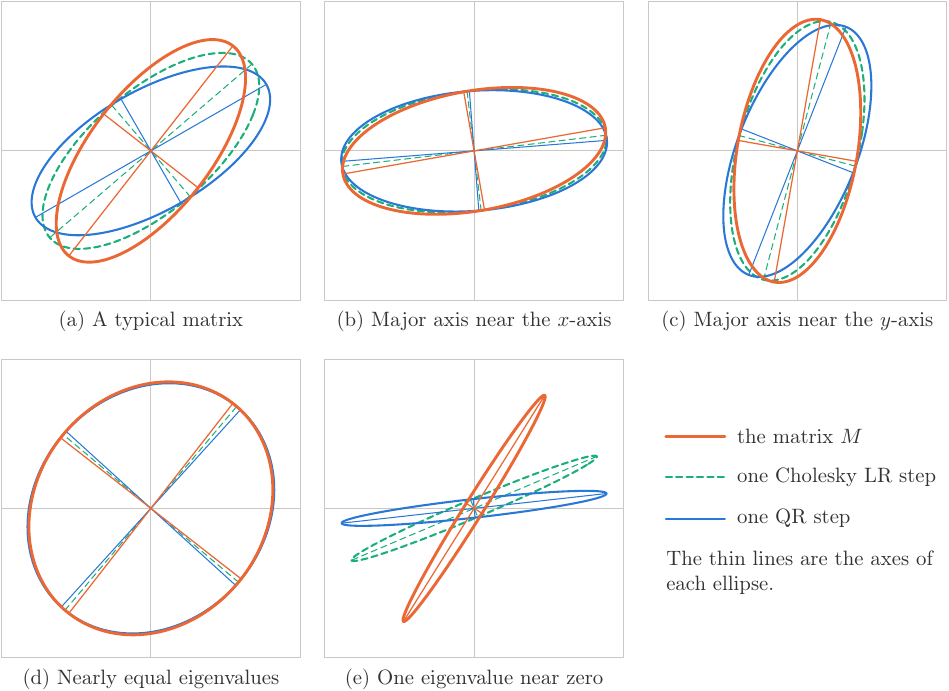}
\caption{Visualisation v1. Each panel shows a $2\times2$ PSD matrix $M$ as the ellipse $\{Mv:\|v\|=1\}$, together with one step of the Cholesky LR algorithm and one step of the QR algorithm. An ellipse with its major axis along the $x$-axis is a diagonal matrix with decreasing diagonal entries, the limit of both algorithms.}
\label{fig:v1}
\end{figure}

This picture is \textbf{valid only for PSD matrices}, unlike v2, which covers every invertible $2\times2$ real matrix. A $2\times 2$ PSD matrix $M$ is drawn as the set $\{ Mv : v \in \RR^2,\ \| v\|_2 = 1 \}$, which by the spectral theorem is an ellipse whose semi-axes are the eigenvectors of $M$ scaled by the eigenvalues. Figure~\ref{fig:v1} shows five examples. In each panel the ellipse of $M$ is drawn with a heavy line, together with two other ellipses. The dashed one is the result of one step of Rutishauser's LR algorithm~\cite{Rutishauser1958} in its Cholesky form, $M=LL^T\mapsto L^TL$, which keeps $M$ symmetric and positive semidefinite. The thin solid one is the result of one step of the QR algorithm, which is the same as two Cholesky LR steps, as recalled after Proposition~\ref{prop:qr-steps}.

The five panels illustrate the following observations.
\begin{itemizedot}
  \item Convergence happens in every case, though not necessarily quickly (panel~a).
  \item Ellipses whose major semi-axis is parallel to the $x$-axis are attracting fixed points, and convergence near them is fast (panel~b).
  \item Ellipses whose major semi-axis is parallel to the $y$-axis are repelling fixed points, and convergence near them is slow (panel~c).
  \item When the eigenvalues are close to each other the ellipse is nearly a circle, and convergence is slow (panel~d).
  \item When one eigenvalue is close to zero and the other is not, the ellipse is very flat, and convergence is fast (panel~e). A \emph{shift} $M\mapsto M-\mu I$, with $\mu$ just below the smaller eigenvalue, flattens the ellipse.
\end{itemizedot}
The closer a matrix is to a repelling fixed point, the more steps it needs to get away (panel~c), so no number of steps is enough for all matrices at once; Section~\ref{sec:qr-nonneg} returns to this.

\section{A Matrix Moves Directions}
\label{sec:directions}
\begin{figure}[ht]
\centering
\includegraphics{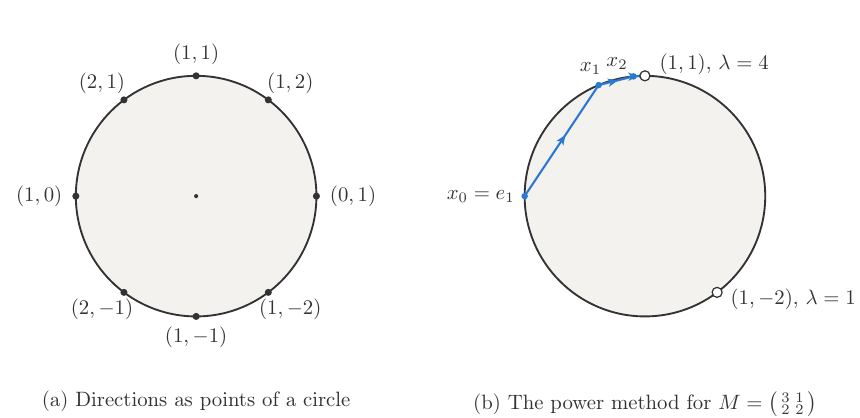}
\caption{(a)~The directions of the plane as points of a circle: the direction of $(v_1,v_2)^T$ is drawn at the point~\eqref{eq:direction-point}, so perpendicular directions are opposite. (b)~The power method for the running example, started at $e_1$. The points $x_k$, the directions of $M^ke_1$, run to the eigenvector $(1,1)^T$. The rings mark the two eigenvectors.}
\label{fig:directions}
\end{figure}

An eigenvector matters only up to scale, so what a matrix does to eigenvectors is best seen on \emph{directions}: lines through the origin of $\RR^2$. A direction is spanned by a vector $v=(v_1,v_2)^T$, and $v$ and $-v$ span the same direction. So if $v$ makes the angle $\psi$ with $e_1$, the direction depends only on $\psi$ modulo $\pi$, and doubling the angle places the directions on a circle. We draw the direction of $v$ at the point
\begin{equation}\label{eq:direction-point}
-e^{-2i\psi}=\Bigl(\frac{v_2^2-v_1^2}{v_1^2+v_2^2},\ \frac{2v_1v_2}{v_1^2+v_2^2}\Bigr)
\end{equation}
of the unit circle: $e_1$ is leftmost, and turning $v$ anticlockwise by $\psi$ moves its point clockwise by $2\psi$. Section~\ref{sec:disc} gives the reason for this choice. Figure~\ref{fig:directions}(a) shows some examples. The direction of $e_1$ is the point $(-1,0)$ on the left, the direction of $e_2$ is the point $(1,0)$ on the right, and perpendicular directions are opposite points. The circle of directions is the real projective line $\RP^1$; we call it \emph{the boundary}, because later it will bound a disc.

An invertible matrix $M$ maps lines through the origin to lines through the origin, so it moves the points of the boundary: the direction of $v$ goes to the direction of $Mv$. We write $f$ for this map of the boundary to itself. When a formula is needed, we label the direction of $(v_1,v_2)^T$ by the number $x=v_1/v_2$, so that $e_1$ has the label $\infty$ and $e_2$ the label $0$; then
\[ f(x)=\frac{ax+b}{cx+d}\qquad\text{for } M=\bmat{a&b\\ c&d}, \]
because $M(x,1)^T=(ax+b,\,cx+d)^T$. This is a M\"obius transformation, and the point~\eqref{eq:direction-point} is $\frac{ix+1}{-ix+1}$.

A direction is fixed by $f$ exactly when it is an eigenvector of $M$. For the running example the fixed points are the directions of $(1,1)^T$ and $(1,-2)^T$, which by~\eqref{eq:direction-point} are the points $(0,1)$ and $(0.6,-0.8)$. Every other direction is moved away from the second and towards the first, because the first has the larger eigenvalue: this is the power method. Starting from $e_1$,
\[ e_1=\begin{pmatrix}1\\0\end{pmatrix},\quad Me_1=\begin{pmatrix}3\\2\end{pmatrix},\quad M^2e_1=\begin{pmatrix}11\\10\end{pmatrix},\quad M^3e_1=\begin{pmatrix}43\\42\end{pmatrix},\ \dots \]
and the corresponding points $x_0,x_1,x_2,\dots$ of the boundary are $(-1,0)$, $(-0.385,0.923)$, $(-0.095,0.995)$, $(-0.024,1.000)$, which run to $(0,1)$ (Figure~\ref{fig:directions}(b)). The distance to the limit shrinks by about a factor $4$ at each step, the ratio of the eigenvalues. We call $(0,1)$ the \emph{attracting} fixed point of $f$ and $(0.6,-0.8)$ the \emph{repelling} one.

\section{The QR Algorithm Is the Power Method in a Turning Frame}
\label{sec:turning}
Run the unshifted QR algorithm $M_0=M$, $M_k=Q_kR_k$, $M_{k+1}=R_kQ_k$, where $Q_k$ is orthogonal and $R_k$ is upper triangular with positive diagonal entries. Two observations explain what a step does.

First, $M_1=R_0Q_0=Q_0^{-1}MQ_0$, so a step replaces $M$ by an orthogonally similar matrix. A similarity by a rotation $R_\theta=\bigl(\begin{smallmatrix}\cos\theta&-\sin\theta\\ \sin\theta&\cos\theta\end{smallmatrix}\bigr)$ does something simple to the boundary. The matrix $R_\theta MR_\theta^{-1}$ moves the direction of $R_\theta v$ to the direction of $R_\theta Mv$, so whatever $M$ does on the boundary, $R_\theta MR_\theta^{-1}$ does the same after every point has been turned. Turning the plane by $\theta$ turns the boundary by $-2\theta$, by~\eqref{eq:direction-point}. Call the boundary with its motion $f$ \emph{the picture of $M$}; Section~\ref{sec:picture} draws it. Then \emph{similarity by a rotation is rotation of the picture}. Similarity by a reflection of the plane reflects the picture in a diameter.

Second, the first column of $M=Q_0R_0$ reads $Me_1=r_{11}Q_0e_1$ with $r_{11}>0$: the orthogonal matrix $Q_0$ takes $e_1$ to the direction of $Me_1$, so $Q_0^{-1}$ takes the direction of $Me_1$ back to $e_1$. When $\det M>0$, $Q_0$ is a rotation, and putting the two together: \emph{one QR step turns the picture of $M$ so that the direction of $Me_1$ comes to the point $(-1,0)$.} Repeating this gives the following statement, which is proved in Section~\ref{sec:recipe}.

\begin{proposition}\label{prop:qr-picture}
Let $M$ be an invertible $2\times2$ real matrix, and let $x_k$ be the direction of $M^ke_1$. Then $M_k=\tilde Q_k^TM\tilde Q_k$, where $\tilde Q_k=Q_0Q_1\cdots Q_{k-1}$ is the orthogonal matrix whose first column is the unit vector along $M^ke_1$. If $\det M>0$, the picture of $M_k$ is the picture of $M$ turned by the rotation that takes the point $x_k$ to $(-1,0)$.
\end{proposition}

So the QR algorithm is the power method $x_0,x_1,x_2,\dots$ of Section~\ref{sec:directions}, watched from a frame that turns so as to keep the current point $x_k$ at $(-1,0)$. For the running example the points $x_k$ run to the eigenvector $(1,1)^T$. In the turning frame this eigenvector therefore comes to rest at $(-1,0)$, which means that $e_1$ is an eigenvector of the limit: the limit is upper triangular. Indeed
\[ M_1=\begin{pmatrix}4.077&-0.385\\0.615&0.923\end{pmatrix},\quad
M_2=\begin{pmatrix}4.041&-0.855\\0.145&0.959\end{pmatrix},\quad
M_3=\begin{pmatrix}4.011&-0.965\\0.035&0.989\end{pmatrix},\quad\dots, \]
which converge to $\bigl(\begin{smallmatrix}4&-1\\0&1\end{smallmatrix}\bigr)$,
and the entry below the diagonal shrinks by about a factor $4$ at each step, as the points $x_k$ did. Three facts are worth remembering.

\emph{The algorithm can only rotate.} Each QR step is an orthogonal similarity, and when $\det M>0$ it is a similarity by a rotation, which is a rotation of the picture.

\emph{The target is the point $(-1,0)$.} A matrix is upper triangular exactly when $e_1$ is an eigenvector, that is, when $(-1,0)$ is a fixed point of $f$. The whole purpose of the QR algorithm is to turn the picture until a fixed point sits at $(-1,0)$.

\emph{The inverse is the reversed picture.} The matrix $M^{-1}$ moves every direction back, so it exchanges the attracting and the repelling fixed point.

\section{The Picture of a Matrix}
\label{sec:picture}
So far a matrix is a motion of the points of a circle. To see the whole motion at a glance, we join each direction $x$ to its image $f(x)$ by an arrow. There are two natural ways of drawing the arrows, and both lead to the same picture.

\subsection{Straight Arrows: Visualisation v1.5}
\begin{figure}[ht]
\centering
\includegraphics{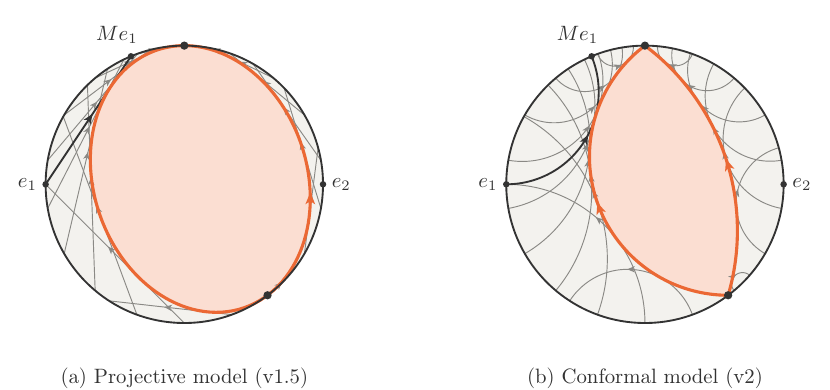}
\caption{The picture of the running example. The circle is the boundary, and the dots on it are the two eigenvectors. Thin lines: canonical arrows from a direction to its image, drawn as chords in (a) and as arcs of circles orthogonal to the boundary in (b); the arrow of $e_1$ is drawn heavier. Heavy line: the envelope of the arrows. Shaded: the region that no arrow enters; Section~\ref{sec:disc} explains what it is.}
\label{fig:models}
\end{figure}

Draw the chord from $x$ to $f(x)$ for every direction $x$ (Figure~\ref{fig:models}(a)). For a matrix with positive determinant, such as the running example, the chords do not fill the disc: they all touch one curve, their \emph{envelope}, and they avoid the region inside it. The envelope is a conic, and it touches the boundary at the eigenvectors, where the chords shrink to points. We give the envelope the direction of the chords that touch it; this direction distinguishes $M$ from $M^{-1}$, which has the same chords, reversed.

The conic already is a picture of the matrix, and one step of the QR algorithm has a brief description in it: draw the chord from $e_1$ to $Me_1$, and turn the picture about the centre of the circle by the rotation that takes the point of $Me_1$ to the point of $e_1$ (Figure~\ref{fig:qrstep}(a)).

\subsection{Round Arrows: Visualisation v2}
\label{subsec:v2}
A conic is awkward to draw and to recognise. The envelope becomes simpler if the arrows are drawn round. Through two points $x$ and $f(x)$ of the boundary there is exactly one circle that crosses the boundary at right angles (a diameter, if the two points are opposite), and we now call this circle, directed from $x$ to $f(x)$ through the disc, the \emph{canonical arrow} $x\to f(x)$ (Figure~\ref{fig:models}(b)). This construction is based on ideas from Yaglom's books~\cite{Yaglom1968-ComplexNumbers,Yaglom1979-SimpleNonEuclid}.

The canonical arrows of a matrix with positive determinant again have an envelope, and now it consists of circles. It is convenient to continue the envelope outside the disc, as the preprint and its interactive program do. The envelope is then a pair of whole oriented circles, which we call the \emph{cycles} of $M$. The inversion $w\mapsto1/\bar w$ in the boundary maps every canonical arrow to itself with its direction reversed, so it exchanges the two cycles: each cycle is the image of the other under the inversion, with its direction reversed. We call such a pair \emph{mutually inverse}. \textbf{The pair of cycles is the picture of the matrix}; it determines the matrix up to a non-zero scalar factor.

The cycles are computed from the entries. For $M=\bmat{a&b\\ c&d}$ with $\Delta=\det M>0$, put
\[ X=-(b+c),\qquad Y=d-a,\qquad C=b-c,\qquad S=a+d,\qquad Z=\pm2\sqrt\Delta . \]
These are the coordinates of Section~\ref{sec:recipe}; here they are only abbreviations. The two signs of $Z$ give the two cycles, and each is the circle with
\begin{equation}\label{eq:cycles}
\text{centre }\Bigl(\frac{X}{C-Z},\ \frac{Y}{C-Z}\Bigr)\qquad\text{and signed radius }\frac{S}{C-Z},
\end{equation}
where a positive radius means anticlockwise and a negative one clockwise (Lemma~\ref{lem:pencil}); if $C=Z$ the cycle is a straight line. Conversely, $a=\frac{S-Y}2$, $d=\frac{S+Y}2$, $b=\frac{C-X}2$ and $c=-\frac{C+X}2$.

For the running example $X=-3$, $Y=-1$, $C=-1$, $S=5$ and $Z=\pm4$. So $Z=4$ gives the circle with centre $(0.6,0.2)$ and radius $1$, run clockwise, and $Z=-4$ the circle with centre $(-1,-\frac13)$ and radius $\frac53$, run anticlockwise. Both pass through the eigenvectors $(0,1)$ and $(0.6,-0.8)$; for instance, $(0,1)$ is at distance $\sqrt{0.6^2+0.8^2}=1$ from $(0.6,0.2)$. They are the heavy circles of Figure~\ref{fig:shapes}(c).

\subsection{A First Reading of the Picture}
\begin{figure}[ht]
\centering
\includegraphics[width=\linewidth]{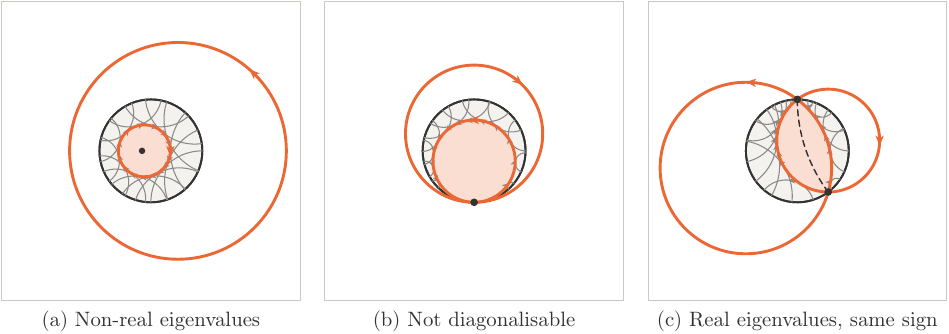}
\caption{The three kinds of picture for $\det M>0$. Heavy lines: the two cycles. Thin arcs: canonical arrows, drawn inside the disc only. Shaded: the part of the disc that no arrow enters. (a)~$M=\bigl(\begin{smallmatrix}2&-2\\1&2\end{smallmatrix}\bigr)$, with non-real eigenvalues: one cycle lies inside the disc and the other outside. (b)~$M=\bigl(\begin{smallmatrix}5&1\\-1&3\end{smallmatrix}\bigr)$, which is not diagonalisable: both cycles touch the boundary at the eigenvector (dot). (c)~The running example: both cycles pass through the two eigenvectors (dots). The dot in (a) and the dashed arc in (c) are explained in Section~\ref{sec:disc}.}
\label{fig:shapes}
\end{figure}

Three things can be read from the cycles before any theory is needed.
\begin{itemizedot}
  \item \emph{The eigenvectors.} The cycles meet the boundary exactly at the fixed points of $f$, the real eigenvectors of $M$. So $M$ is upper triangular exactly when a cycle passes through $(-1,0)$, lower triangular exactly when a cycle passes through $(1,0)$, and diagonal when both happen.
  \item \emph{The kind of matrix.} There are three kinds of picture (Figure~\ref{fig:shapes}). If the eigenvalues are not real there is no fixed point on the boundary: one cycle lies inside the disc and the other outside. If $M$ is not diagonalisable there is one fixed point, and both cycles touch the boundary there. If the eigenvalues are real and distinct, both cycles cross the boundary at the two eigenvectors.
  \item \emph{The inverse.} The picture of $M^{-1}$ is the picture of $M$ with both cycles reversed.
\end{itemizedot}
The scalar matrices are the exception: they move no direction, and their picture is the boundary circle itself.

\section{Watching the QR Algorithm}
\label{sec:qr-nonneg}
\begin{figure}[ht]
\centering
\includegraphics{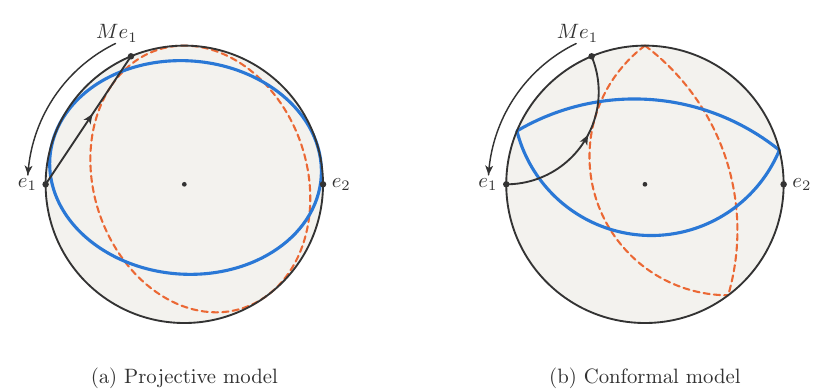}
\caption{One step of the QR algorithm for the running example, in the two ways of drawing the arrows. Dashed: the envelope of $M$. Solid: the envelope of the next iterate $M_1$. The step turns the picture about the centre (dot) by the rotation that takes the point of $Me_1$ to the point of $e_1$; the chord and the arc from $e_1$ to $Me_1$ are the canonical arrow of $e_1$.}
\label{fig:qrstep}
\end{figure}

This section treats matrices with positive determinant; Section~\ref{sec:qr-negdet} treats the others. By Proposition~\ref{prop:qr-picture}, every iterate is the input picture turned about the centre, by the rotation that brings $x_k$ to $(-1,0)$; the cycles are drawn from the motion of the boundary, so they turn with it. Figure~\ref{fig:qrstep} shows one step for the running example.

The remaining figures of this section show the cycles as whole circles, cropped by the frame where they are large; the parts outside the disc mirror the parts inside. An interactive program draws the same pictures. In it the user controls one cycle of the pair by translating it, scaling it, or reversing its orientation; the other cycle follows, and may become a straight line. In this way the user can explore all $2\times2$ real matrices with non-negative determinant. The conventions in the figures are as follows.
\begin{itemizedot}
  \item The shaded disc is bounded by the boundary, the unit circle, which is the picture of the identity matrix and its scalar multiples. The small ring on it marks the point $(-1,0)$.
  \item The input matrix $M_0$ is the pair of oriented cycles drawn with a heavy line; the arrowheads give the orientations.
  \item Of the iterates $M_1,M_2,\dots$ only the last is drawn, with a thinner line. The steps in between are shown by dots. Each dot is the place of one \emph{marked point} of the picture after a step, and the larger dot is its place in the input.
  \item When $M$ has a real eigenvector, the marked point is the eigenvector whose eigenvalue has the larger modulus. This is the point that the algorithm brings to $(-1,0)$, and a thin arc outside the boundary shows which way it travels. When $M$ has no real eigenvector, the marked point is the dot of Figure~\ref{fig:shapes}(a), which lies inside the smaller cycle and travels on a circle about the centre of the disc.
\end{itemizedot}

\subsection{The Three Kinds of Picture}

\begin{figure}[ht]
\centering
\includegraphics[width=\linewidth]{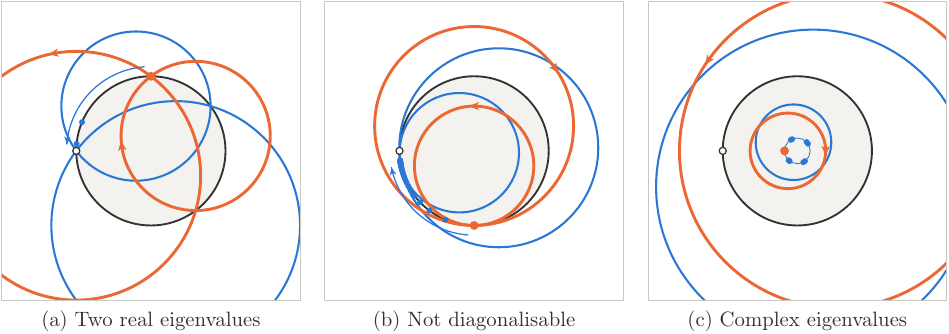}
\caption{The unshifted QR algorithm for the three matrices of Figure~\ref{fig:shapes}: (a)~the running example, (b)~$\bigl(\begin{smallmatrix}5&1\\-1&3\end{smallmatrix}\bigr)$ and (c)~$\bigl(\begin{smallmatrix}2&-2\\1&2\end{smallmatrix}\bigr)$. Heavy line: the input. Thinner line: the last iterate, after six steps in panel~(a), thirty in panel~(b) and nine in panel~(c). Dots: the marked point after each step. The ring marks $(-1,0)$.}
\label{fig:cases123}
\end{figure}

Consider the input cycles, drawn with a heavy line, in Figure~\ref{fig:cases123}.
\begin{itemizedot}
  \item In panel (a), the boundary meets each input cycle at two distinct points. These correspond to two \emph{real eigenvectors}. Since $\det M>0$, the two eigenvalues have the same sign, so $M$ is \emph{diagonalisable with real eigenvalues of the same sign}. The map $f$ moves every other direction away from one eigenvector and towards the other, so the points $x_k$ converge, and the error is multiplied by about $\lambda_2/\lambda_1$ at each step, where $|\lambda_1|>|\lambda_2|$.
  \item In panel (b), the boundary meets each input cycle at a single point, which corresponds to the only real eigenvector of $M$. Therefore $M$ is \emph{not diagonalisable, and its only eigenvalue is real}. The points $x_k$ still converge to the eigenvector, but the error decreases only like $1/k$.
  \item In panel (c), the input cycles do not meet the boundary, so $M$ has no real eigenvector. Its eigenvectors are complex, with a complex eigenvalue $\lambda$; since $\bar\lambda$ is also an eigenvalue, $M$ is \emph{diagonalisable over $\CC$ with non-real eigenvalues}. The points $x_k$ go round the circle for ever (Section~\ref{subsec:complex}) and do not converge.
\end{itemizedot}
These are the three kinds of Section~\ref{sec:picture}. The non-zero scalar matrices are represented by the boundary circle itself and are fixed by the algorithm, the non-zero singular matrices are those for which the two cycles coincide, and the matrices whose eigenvalues are real with opposite signs are treated in Section~\ref{sec:qr-negdet}.

Each iterate is obtained from the input by a rotation about the centre, which reflects the fact that each iterate is orthogonally similar to the previous ones. In panels (a) and (b) the iterates converge, and in panel (c) they do not. When the algorithm converges, the last iterates pass through $(-1, 0)$, so the limit is upper triangular. The non-diagonalisable matrix of panel (b) converges more slowly than the matrix of panel (a), but it still converges. Moving from (a) to (b) to (c), convergence slows down until it stops.

\subsection{Behaviour near the Upper Triangular Limit}

\begin{figure}[ht]
\centering
\includegraphics{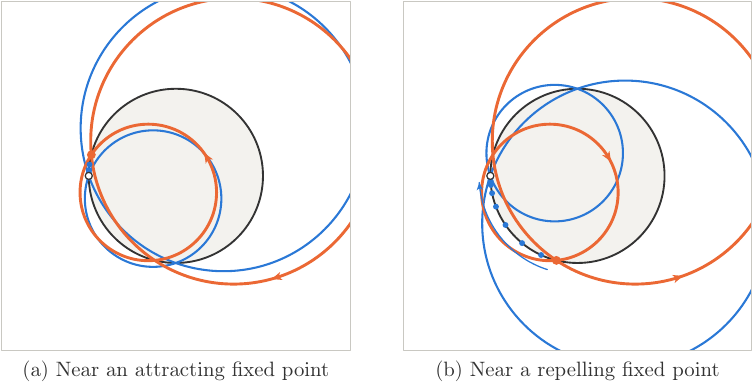}
\caption{Fixed points of the QR algorithm. (a)~The attracting fixed point $\bigl(\begin{smallmatrix}2&1\\0&1\end{smallmatrix}\bigr)$, turned a little: the iterates return to it. (b)~Its inverse, the same pair with the orientations reversed, which is a repelling fixed point: the iterates move away and settle at another fixed point.}
\label{fig:cases45}
\end{figure}

An upper triangular matrix $T$ is a fixed point of the algorithm, and its cycles pass through $(-1,0)$, because $e_1$ is an eigenvector of $T$. Let its diagonal entries differ. If $e_1$ belongs to the eigenvalue of larger modulus, a matrix that starts close to $T$ stays close and converges to an upper triangular matrix close to $T$; we call $T$ \emph{attracting}. If $e_1$ belongs to the smaller one, a matrix that starts close to $T$ is carried away, however close it started (unless $e_1$ is an eigenvector of it too), and converges to an upper triangular matrix with the diagonal entries in the other order; we call $T$ \emph{repelling}. Figure~\ref{fig:cases45} starts the algorithm close to one fixed point of each kind, of which the second is obtained from the first by reversing the orientation, that is, by taking the matrix inverse. \emph{The matrix inverse of an attracting fixed point is a repelling fixed point.}

The closer a matrix starts to a repelling fixed point, the more steps it needs to get away, so no number of steps is enough for all matrices at once. A continuous step cannot make every fixed point attracting. Turning the picture of $M$ through all angles gives a closed curve of matrices, two of which are upper triangular; if both attracted their neighbours and every matrix converged, the curve would split into two disjoint open sets, which a connected curve cannot do. For symmetric tridiagonal matrices, a continuous shift strategy of the kind studied in~\cite{LeiteSaldanhaTomei2013IMRN} fails to converge for some matrices. The Wilkinson shift is discontinuous, and the strategy of~\cite{goodshift2022} adds exceptional shifts.

\subsection{Behaviour near the Identity Matrix and near Singular Matrices}

\begin{figure}[ht]
\centering
\includegraphics{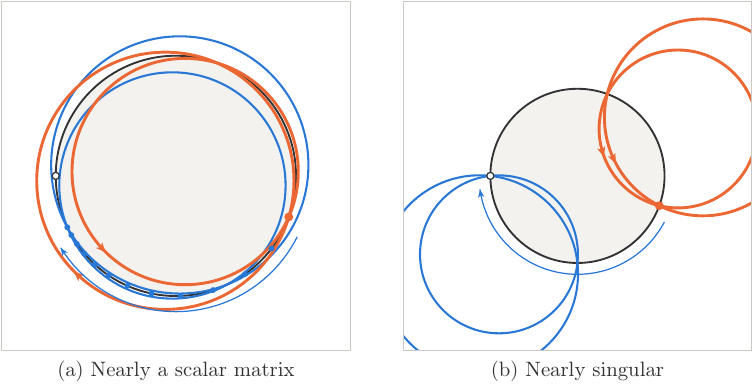}
\caption{Twelve QR steps for two matrices with the same eigenvectors. (a)~A matrix close to the identity, with eigenvalues $1$ and $0.88$, drawn enlarged: the input cycles nearly coincide with the boundary, and convergence is very slow. (b)~A matrix close to a singular one, with eigenvalues $1$ and $0.003$: the two input cycles nearly coincide with each other, and convergence is very fast. After the first step the marked point is already at $(-1,0)$.}
\label{fig:cases678}
\end{figure}

In Figure~\ref{fig:cases678}(a), the two input cycles are very close to the boundary, so the matrix is close to a multiple of the identity. The algorithm converges very slowly: after twelve steps the iterates are still on their way to $(-1,0)$. The error shrinks at each step by about the ratio of the eigenvalues, here $0.88$, and $0.88^{12}\approx0.2$. Moreover, when two eigenvalues nearly coincide, a small change of the matrix can turn its eigenvectors by a large angle, whatever the algorithm.

Figure~\ref{fig:cases678}(b) is the opposite. The two mutually inverse cycles almost coincide with each other, as they do when one eigenvalue is much smaller in modulus than the other, $|\lambda_2|\ll|\lambda_1|$. Convergence is \textbf{very fast}, since the subdiagonal entry is multiplied by about $\lambda_2/\lambda_1$ at each step. The non-diagonalisable matrix of Figure~\ref{fig:cases123}(b) converges, but slowly.

\subsection{The Role of Shifts}
The speed is governed by the ratio $\lambda_2/\lambda_1$ of the eigenvalues, where $|\lambda_1|\geq|\lambda_2|$.
\begin{itemizedot}
  \item Ratio close to $1$ (Figure~\ref{fig:cases678}(a)): slow.
  \item Equal eigenvalues, not diagonalisable (Figure~\ref{fig:cases123}(b)): slower still, like $1/k$. Beyond this lie the non-real eigenvalues, where there is no convergence.
  \item Ratio close to $0$ (Figure~\ref{fig:cases678}(b)): fast.
\end{itemizedot}
A \emph{shift} moves a matrix from the first two cases towards the third~\cite{goodshift2022}: choose $\mu$ close to an eigenvalue, take a QR step of the nearly singular $M-\mu I$, and add $\mu I$ back.

\subsection{The Complex Eigenvalue Case}
\label{subsec:complex}

\begin{figure}[ht]
\centering
\includegraphics{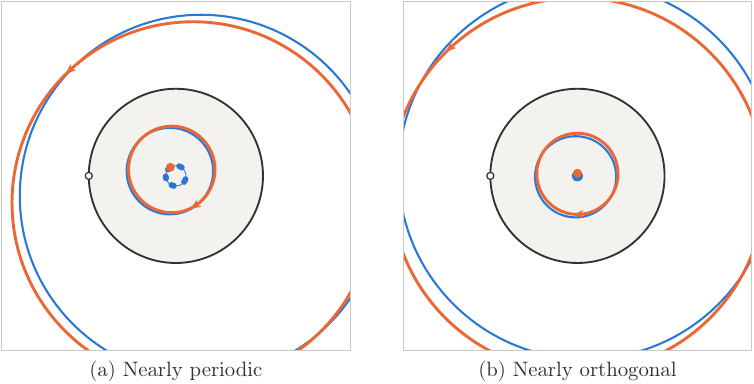}
\caption{Fifteen QR steps for two matrices with complex eigenvalues. (a)~Eigenvalues $e^{\pm i\theta}$ with $\theta/2\pi=0.102$, close to $\frac1{10}$: the iterates are nearly periodic with period five, and the dots fall into five groups. (b)~A matrix close to an orthogonal matrix, whose two cycles are nearly concentric with the boundary: it is nearly a fixed point, and the dots stay together.}
\label{fig:cases910}
\end{figure}

For complex eigenvalues the iterates are nearly periodic, except when $M$ is close to an orthogonal matrix (Figure~\ref{fig:cases910}). If the eigenvalues are $re^{\pm i\theta}$, then $M=g\,(rR_\theta)\,g^{-1}$ for some real matrix $g$, so in suitable coordinates $M$ turns every direction by the angle $\theta$. If $\theta=2\pi m/n$ with $m/n$ in lowest terms, the iterates return to the start after $n$ steps if $n$ is odd and after $n/2$ steps if $n$ is even, because $M^{n/2}$ is then a scalar matrix. Otherwise they never return exactly, but nearly do whenever $\theta/2\pi$ is close to such a fraction. Panel (a) shows the near periodicity. The matrix of panel (b) is close to an orthogonal matrix: its two cycles are nearly concentric with the boundary circle, and it is nearly a fixed point of the algorithm. Indeed, a scalar multiple $rQ$ of an orthogonal matrix has the QR factorisation $Q\cdot(rI)$, so one QR step returns $rQ$ itself.

\section{What the Inside of the Disc Means}
\label{sec:disc}
Until now the disc inside the boundary was drawing paper. This section says what its points are, and why the arrows of a matrix leave a region of it empty.

\subsection{Points of the Disc Are Inner Products}
The QR algorithm depends on a choice that has been invisible so far: the inner product with respect to which $Q$ is orthogonal. An inner product on $\RR^2$ is $\langle u,v\rangle_G=u^TGv$ for a positive definite symmetric matrix $G$; the standard one has $G=I$. Multiplying $G$ by a positive scalar changes no angle and no orthogonality, so we consider inner products up to positive scalars, and write $\mathcal P$ for the set of them. To the inner product $G=\bigl(\begin{smallmatrix}g_{11}&g_{12}\\ g_{12}&g_{22}\end{smallmatrix}\bigr)$ we assign the point
\begin{equation}\label{eq:charts}
w(G)=\frac{g_{11}-g_{22}-2i\,g_{12}}{g_{11}+g_{22}+2\sqrt{\det G}}
\end{equation}
of the plane. This is a bijection from $\mathcal P$ onto the open unit disc, which we call the \emph{conformal chart}. Three examples show how it works.
\begin{itemizedot}
  \item The standard inner product $I$ is the centre of the disc.
  \item The inner product $\diag(1,\varepsilon)$, which counts the second coordinate with the weight $\varepsilon$, is the point $\bigl(\frac{1-\sqrt\varepsilon}{1+\sqrt\varepsilon},0\bigr)$ of the horizontal diameter. As $\varepsilon\to0$ it runs to $(1,0)$, the direction of $e_2$, and as $\varepsilon\to\infty$ it runs to $(-1,0)$, the direction of $e_1$.
  \item In general, as $G$ degenerates to a matrix of rank one, its point runs to the boundary, and it arrives at the direction that the limit matrix sends to zero. This is the reason for the choice~\eqref{eq:direction-point}: the boundary of the disc of inner products is the circle of directions.
\end{itemizedot}
There is a second chart,
\[ k(G)=\frac{g_{11}-g_{22}-2i\,g_{12}}{g_{11}+g_{22}}=\frac{2w}{1+|w|^2}, \]
the \emph{projective chart}, which also maps $\mathcal P$ onto the open unit disc, with the same boundary. It moves each point of the disc outwards along its radius. Visualisation v1.5 is drawn in the projective chart and visualisation v2 in the conformal chart.

Now the arrows acquire a meaning. For two different directions $u$ and $v$, the inner products for which $u\perp v$ form a curve in the disc, which joins the boundary points $u$ and $v$. In the projective chart this curve is the chord from $u$ to $v$, because $u^TGv=0$ is one linear condition on $G$; in the conformal chart it is the arc of the circle through $u$ and $v$ orthogonal to the boundary. So in both pictures
\[ \text{\emph{the canonical arrow of $v$ is the set of inner products for which $Mv$ is perpendicular to $v$.}} \]
For the running example and $v=e_1$ we have $Me_1=(3,2)^T$, and for instance $G=\bigl(\begin{smallmatrix}2&-3\\-3&5\end{smallmatrix}\bigr)$ lies on the arrow of $e_1$, since $e_1^TG\,(3,2)^T=6-6=0$.

\subsection{The Region Where the Hermitian Part Is Definite}
An inner product $G$ gives every $2\times2$ matrix an adjoint,
\[ X^*=G^{-1}X^TG,\qquad\text{so that}\qquad\langle Xu,v\rangle_G=\langle u,X^*v\rangle_G . \]
The map $X\mapsto X^*$ is an involution of the algebra $M_2(\RR)$, so each point of the disc is a way of making $M_2(\RR)$ into a \staralg{}. The \emph{Hermitian part} of $M$ for $G$ is $\frac12(M+M^*)$, and
\[ \langle v,(M+M^*)v\rangle_G=2\langle v,Mv\rangle_G . \]
So an arrow of $M$ passes through the point $G$ exactly when $\langle v,Mv\rangle_G=0$ for some direction $v$, that is, when the quadratic form of $M+M^*$ vanishes along a direction. If $M+M^*$ is positive or negative definite, this never happens. That is the empty region of Figure~\ref{fig:models}, and the envelope is where definiteness is lost.

\begin{theorem}\label{thm:picture}
Let $M$ be an invertible $2\times2$ real matrix, and let $G$ be an inner product on $\RR^2$.
\begin{enumerate}
\item The number of canonical arrows of $M$ through the point $G$ is $0$ if $M+M^*$ is definite, $1$ if it has rank one, and $2$ if it is indefinite. If $M+M^*=0$, every canonical arrow passes through $G$.
\item If $\det M<0$, then $M+M^*$ is indefinite for every $G$, so exactly two canonical arrows pass through every point.
\item If $\det M>0$, the inner products for which $M+M^*$ is singular form the envelope of the canonical arrows. In the projective chart the envelope is the conic $\det(GM+M^TG)=0$, which meets the boundary only at the real eigenvectors of $M$ and touches it there. The envelope separates the region
\[ \mathcal D(M)=\{\,G\in\mathcal P:\ M+M^*\ \text{is definite}\,\}, \]
through whose points no arrow passes, from the region through whose points two arrows pass.
\item For $g\in\GL_2(\RR)$, the arrows, the envelope and the region $\mathcal D$ of $gMg^{-1}$ are the images under $g$ of those of $M$.
\end{enumerate}
\end{theorem}

For the running example and the standard inner product, $M+M^T=\bigl(\begin{smallmatrix}6&3\\3&4\end{smallmatrix}\bigr)$ is positive definite, so the centre of the disc lies in $\mathcal D(M)$, and indeed it lies inside both cycles. Exercise~\ref{ex:diameter} finds the part of the horizontal diameter that lies in $\mathcal D(M)$.

Part (4) of the theorem is the rule by which the pictures move. A matrix $g$ acts on inner products by $g\cdot G=g^{-T}Gg^{-1}$, which is chosen so that $\langle gu,gv\rangle_{g\cdot G}=\langle u,v\rangle_G$, and conjugation by $g$ carries the adjoint for $G$ to the adjoint for $g\cdot G$:
\begin{equation}\label{eq:equivariance}
\bigl(gXg^{-1}\bigr)^{*_{g\cdot G}}=g\,X^{*_G}\,g^{-1}.
\end{equation}
A rotation $R_\theta$ of the plane fixes the standard inner product, and it turns the disc about its centre by $-2\theta$. This is the rule of Section~\ref{sec:turning}, now for the whole disc and not only for its boundary.

\subsection{The Disc Is the Hyperbolic Plane}
To say what shape the region $\mathcal D(M)$ has, we need distances in the disc. The set $\mathcal P$ carries a natural distance,
\[ d(G,H)=\tfrac12\bigl|\ln(\mu_1/\mu_2)\bigr|,\qquad\text{where $\mu_1,\mu_2$ are the eigenvalues of $G^{-1}H$,} \]
which measures how far $H$ is from being a multiple of $G$, and which every $g$ preserves. For example, $\diag(1,\varepsilon)$ is at distance $\frac12|\ln\varepsilon|$ from the centre. With this distance $\mathcal P$ is the hyperbolic plane; the conformal chart is its Poincar\'e disc model and the projective chart is its projective model~\cite{Beardon1983,Yaglom1968-ComplexNumbers}. If $z$ is the point of the upper half-plane with $w=\frac{iz+1}{-iz+1}$, then $G$ is proportional to $\bigl(\begin{smallmatrix}1&-\Re z\\ -\Re z&|z|^2\end{smallmatrix}\bigr)$, and a matrix $g$ with M\"obius transformation $f_g$ acts on $\mathcal P$ by $z\mapsto f_g(z)$ if $\det g>0$ and by $z\mapsto f_g(\bar z)$ if $\det g<0$. We use three standard facts.
\begin{itemizedot}
  \item The shortest paths, the \emph{geodesics}, are the curves described above: the geodesic joining two boundary points $u$ and $v$ is the set of inner products for which $u\perp v$. So a canonical arrow is a geodesic.
  \item The circles about a point, the horocycles based at a boundary point, and the curves at a fixed distance from a geodesic, which we call \emph{equidistant curves}, all appear in the conformal chart as arcs of Euclidean circles: a circle inside the disc, a circle inside the disc touching the boundary at the base point, and an arc joining the two ends of the geodesic.
  \item A matrix that is orthogonal for $G$ fixes the point $G$. It acts as a rotation about $G$ if its determinant is positive, and as a reflection in a geodesic through $G$ if it is negative.
\end{itemizedot}

\subsection{The Three Kinds of Picture, Again}
The three kinds of picture of Section~\ref{sec:picture} can now be described exactly. The kind is decided by the algebra $\RR[M]$ of polynomials in $M$. If $M$ is not a scalar matrix, $\RR[M]$ is a commutative algebra of dimension two, and so it is isomorphic to $\CC$, to the dual numbers $\RR[\varepsilon]/(\varepsilon^2)$ or to $\twoR=\RR\oplus\RR$, according as the eigenvalues of $M$ are non-real, equal, or real and distinct. These are the three kinds of M\"obius transformation, called elliptic, parabolic and hyperbolic; Kisil~\cite{Kisil2012} develops the geometry of all three from the same three algebras.

\begin{proposition}\label{prop:shapes}
Let $M$ be a $2\times2$ real matrix with $\det M>0$ that is not a scalar matrix.
\begin{enumerate}
\item \emph{Elliptic case.} Let the eigenvalues be $re^{\pm i\theta}$ with $0<\theta<\pi$. There is exactly one inner product $p$ for which $M$ is $r$ times a rotation, and it is the only point of $\mathcal P$ fixed by $M$. The region $\mathcal D(M)$ is the open disc about $p$ of radius $d$, where
\[ \cosh d=\frac1{\sin\theta}. \]
It is empty when $\theta=\frac\pi2$, in which case every arrow passes through $p$.
\item \emph{Parabolic case.} Let $M=\lambda(I+N)$ with $N^2=0$ and $N\neq0$. Then $\mathcal D(M)$ is the region inside a horocycle based at the eigenvector of $M$: it consists of the inner products for which the operator norm of $N$ is less than $2$.
\item \emph{Hyperbolic case.} Let the eigenvalues $\lambda_1,\lambda_2$ be real and distinct. The inner products for which $M$ is symmetric form the geodesic joining the two eigenvectors, which we call the \emph{axis} of $M$. The region $\mathcal D(M)$ consists of the points at distance less than $\delta$ from the axis, where
\[ \sinh\delta\,\sinh\frac\ell2=1,\qquad \ell=\Bigl|\ln\frac{\lambda_1}{\lambda_2}\Bigr|. \]
\end{enumerate}
\end{proposition}

Table~\ref{tab:kinds} lists the three cases together with the remaining kinds of matrix. A non-zero singular matrix is a limit of the hyperbolic case: as $\ell\to\infty$ the region shrinks to the axis, which now joins the kernel to the range, and $M+M^*$ is singular exactly there.

For the running example $\ell=\ln4$ and $\sinh\frac\ell2=\frac34$, so $\sinh\delta=\frac43$ and $\delta=\ln3$: the region $\mathcal D(M)$ consists of the inner products within the distance $\ln3$ of the axis. For the matrix $\bigl(\begin{smallmatrix}2&-2\\1&2\end{smallmatrix}\bigr)$ of Figure~\ref{fig:shapes}(a), Exercise~\ref{ex:elliptic} finds the point $p$ and the radius $d$.

\begin{table}[ht]
\centering
\caption{The kinds of invertible $2\times2$ real matrix $M$, the algebra $\RR[M]$ that each generates, and its picture.}
\label{tab:kinds}
\small
\begin{tabular}{@{}>{\raggedright\arraybackslash}p{0.27\linewidth}>{\raggedright\arraybackslash}p{0.16\linewidth}>{\raggedright\arraybackslash}p{0.23\linewidth}>{\raggedright\arraybackslash}p{0.25\linewidth}@{}}
\toprule
Eigenvalues of $M$ & $\RR[M]$ & Fixed by $M$ & Region $\mathcal D(M)$ \\
\midrule
Non-real & $\CC$ & a point $p$ of the disc & a disc about $p$ \\\addlinespace[3pt]
Equal, $M$ not diagonalisable & dual numbers & a boundary point & the inside of a horocycle based there \\\addlinespace[3pt]
Real and distinct, of the same sign & $\twoR$ & two boundary points and the axis joining them & the points within $\delta$ of the axis \\\addlinespace[3pt]
Real, of opposite signs & $\twoR$ & the same & empty \\\addlinespace[3pt]
$M$ scalar & $\RR$ & every point & every point \\
\bottomrule
\end{tabular}
\end{table}

\section{Reading a Matrix from Its Picture}
\label{sec:reading}
\begin{figure}[ht]
\centering
\includegraphics[width=\linewidth]{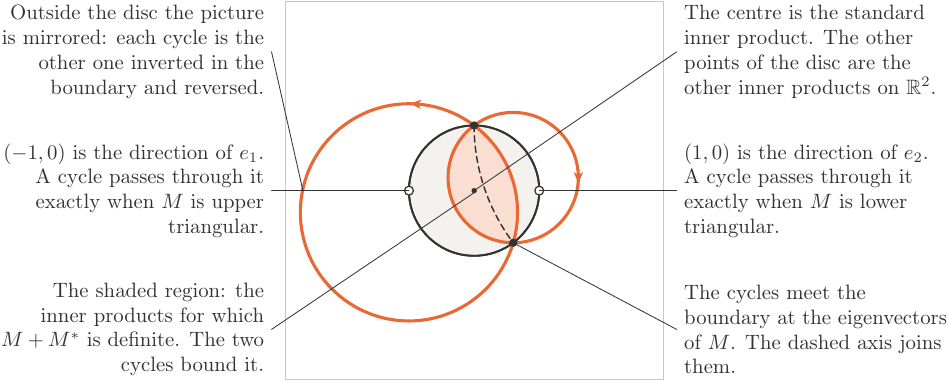}
\caption{The parts of visualisation v2, for the running example.}
\label{fig:key}
\end{figure}

Figure~\ref{fig:key} labels the parts of the picture of the running example, and Table~\ref{tab:reading} collects what the picture says about a matrix. The table covers the properties promised in the Overview. Its rows follow from Sections~\ref{sec:picture} and~\ref{sec:disc}. For instance, $M$ moves the centre $I$ to the inner product $M^{-T}M^{-1}$, whose distance from $I$ is $\ln(\sigma_1/\sigma_2)$, where $\sigma_1\geq\sigma_2$ are the singular values of $M$; this is the row for the condition number. For the running example $\kappa(M)=4.27$, so $M$ moves the centre by the distance $1.45$, and the matrix $V$ of unit eigenvectors has $\kappa(V)=(\sqrt{10}+1)/3=1.39$, so the centre is at the distance $0.33$ from the axis.

\begin{table}[tp]
\centering
\caption{Reading a matrix $M$ with $\det M>0$ from its picture in the conformal chart. Distances are hyperbolic distances in the disc.}
\label{tab:reading}
\small
\begin{tabular}{@{}>{\raggedright\arraybackslash}p{0.33\linewidth}>{\raggedright\arraybackslash}p{0.63\linewidth}@{}}
\toprule
Property of $M$ & In the picture \\
\midrule
$M+M^*$ is definite for the inner product $G$ & The point $G$ lies in the region bounded by the cycles. \\\addlinespace[3pt]
Real eigenvectors & The points where the cycles meet the boundary. \\\addlinespace[3pt]
Jordan normal form & The cycles cross the boundary at two points: real eigenvalues, distinct and of the same sign. They touch it at one point: not diagonalisable. They do not meet it: non-real eigenvalues. Both are the boundary itself: a scalar matrix. \\\addlinespace[3pt]
Orthogonally similar matrices & The pictures obtained by turning the pair about the centre, for similarity by a rotation, or by reflecting it in a diameter, for similarity by a reflection. \\\addlinespace[3pt]
Condition number $\kappa(M)$ & $\ln\kappa(M)$ is the distance from the centre to its image under $M$. For a cycle with centre $e$ and radius $r$, $\sinh\bigl(\tfrac12\ln\kappa(M)\bigr)=2|e|\,/\,\bigl|1+r^2-|e|^2\bigr|$. \\\addlinespace[3pt]
Condition number $\kappa(V)$ of the matrix $V$ of unit eigenvectors & $\ln\kappa(V)$ is the distance from the centre to the axis. \\\addlinespace[3pt]
Upper triangular & A cycle passes through $(-1,0)$. \\\addlinespace[3pt]
Lower triangular & A cycle passes through $(1,0)$. \\\addlinespace[3pt]
Diagonal & A cycle passes through both $(-1,0)$ and $(1,0)$. \\\addlinespace[3pt]
Symmetric & The axis passes through the centre: the cycles meet the boundary at opposite points. \\\addlinespace[3pt]
Orthogonal, up to a scalar & The centre is the fixed point: both cycles are centred at the origin. \\\addlinespace[3pt]
Inverse $M^{-1}$ & The same pair with both orientations reversed. \\\addlinespace[3pt]
Direction of $Mv$ & Two circles through the point of $v$ cross the boundary at right angles and touch the cycles. Directed like the cycle where they touch, one leaves the point of $v$ and ends at the point of $Mv$; the other arrives from $M^{-1}v$. \\\addlinespace[3pt]
Product $MK$ & Follow a canonical arrow of $K$, then a canonical arrow of $M$. The arrows $x\to f_M(f_K(x))$ obtained in this way are the canonical arrows of $MK$. \\
\bottomrule
\end{tabular}
\end{table}

\section{From the QR Algorithm to the LR Algorithm}
\label{sec:qr-lr}
The QR algorithm of Sections~\ref{sec:turning} and~\ref{sec:qr-nonneg} uses the standard inner product: the factors $Q_k$ are orthogonal for it, and the pictures turn about the centre of the disc, which is its point. The same program can be run in the \staralg{} $M_2(\RR)$ with the adjoint of any other inner product $G$. Factorise $M_k=Q_kR_k$, where $Q_k^*Q_k=I$ and $R_k$ is upper triangular with positive diagonal entries, and put $M_{k+1}=R_kQ_k$. If $G=D^TD$ with $D$ upper triangular, this is the standard algorithm applied to $DMD^{-1}$ and conjugated back by $D$. What changes in the picture is the centre of rotation.

\begin{proposition}\label{prop:qr-lr}
Let $M$ be a $2\times2$ real matrix with $\det M>0$, and run the QR algorithm for the inner product $G$.
\begin{enumerate}
\item The picture of $M_k$ is the image of the picture of $M$ under the rotation about the point $G$ that takes the point $x_k$, the direction of $M^ke_1$, to $(-1,0)$.
\item Let $G_\varepsilon=\diag(1,\varepsilon)$, the weighted inner product~\eqref{eq:weighted-ip} for $n=2$. Its point in the conformal chart is $\bigl(\frac{1-\sqrt\varepsilon}{1+\sqrt\varepsilon},0\bigr)$, which moves along a diameter from the centre to $(1,0)$, the direction of $e_2$, as $\varepsilon$ decreases from $1$ to $0$. If $m_{11}\neq0$, then as $\varepsilon\to0^+$ the factor $Q_0$ converges, up to sign, to the unit lower triangular factor $L$ of the Doolittle factorisation $M=LU$, and the step $M\mapsto Q_0^{-1}MQ_0$ converges to the step $M=LU\mapsto UL$ of the LR algorithm.
\item The matrix $L$ fixes the boundary point $(1,0)$ and moves the other points of the disc along the horocycles based there. No such motion moves $(1,0)$ itself, and accordingly the LR step is undefined exactly when the direction of $Me_1$ is the point $(1,0)$, that is, when $m_{11}=0$.
\end{enumerate}
\end{proposition}

\begin{figure}[ht]
\centering
\includegraphics[width=\linewidth]{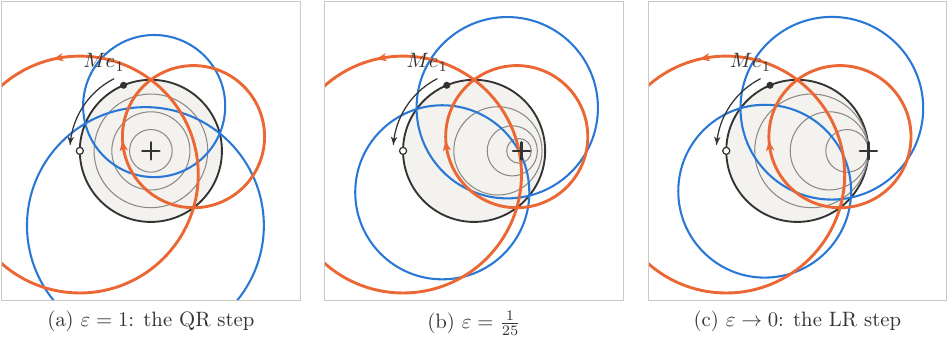}
\caption{One step of the QR algorithm for the inner product $\diag(1,\varepsilon)$, for the running example. Heavy line: the input. Thinner line: the result of the step. The cross marks the point of the inner product, and the thin grey curves are the paths along which the step moves the points of the disc. (a)~The standard inner product: the paths are the circles about the centre. (b)~The cross has moved towards $(1,0)$, and the paths are the circles about it. (c)~The limit, the LR step: the cross is the boundary point $(1,0)$, and the paths are the horocycles based there. In each panel the step carries the point of $Me_1$ to $(-1,0)$.}
\label{fig:qr-to-lr}
\end{figure}

Figure~\ref{fig:qr-to-lr} shows one step of each algorithm on the way from QR to LR. In all three panels the step brings the point of $Me_1$ to $(-1,0)$, and $k$ steps bring the same point $x_k$ there, so all three algorithms converge to an upper triangular matrix. Only the motion that brings the point there differs: a rotation about the centre of the disc, a rotation about another of its points, or a motion along horocycles.

This is the case $n=2$ of two results of the thesis. Proposition~\ref{prop:gs-limit} and Theorem~\ref{thm:bruhat} show that the Gram--Schmidt process for the inner product~\eqref{eq:weighted-ip} tends to Gaussian elimination as $\varepsilon\to0^+$, and that the exceptional matrices, here those with $m_{11}=0$, are the ones whose Bruhat decomposition has a non-trivial permutation. Proposition~\ref{prop:2R-algorithms} shows that the LR algorithm is what the QR algorithm becomes over $\twoR$, where one QR step is two LR steps. So the LR algorithm is reached from the QR algorithm in two ways: by changing the scalars, and, as here, by letting the involution of $M_2(\RR)$ degenerate. That the QR and LR algorithms are one iteration with different changes of coordinates is the starting point of the theory of GR algorithms~\cite{WatkinsElsner1991}; the picture shows what the change of coordinates is for $2\times2$ matrices.

For the running example the Doolittle factorisation is $M=\bigl(\begin{smallmatrix}1&0\\ 2/3&1\end{smallmatrix}\bigr)\bigl(\begin{smallmatrix}3&1\\ 0&4/3\end{smallmatrix}\bigr)$, and the LR step gives $UL=\bigl(\begin{smallmatrix}11/3&1\\ 8/9&4/3\end{smallmatrix}\bigr)=\bigl(\begin{smallmatrix}3.667&1\\0.889&1.333\end{smallmatrix}\bigr)$. The QR step for the inner product $\diag(1,\frac1{25})$ gives $\bigl(\begin{smallmatrix}3.690&0.955\\0.873&1.310\end{smallmatrix}\bigr)$, which is already close to it, while the standard QR step gave $M_1=\bigl(\begin{smallmatrix}4.077&-0.385\\0.615&0.923\end{smallmatrix}\bigr)$.

\section{Matrices with Negative Determinant}
\label{sec:qr-negdet}
\subsection{The Picture}
\label{subsec:negdet}
Let $\det M<0$. By Theorem~\ref{thm:picture}(2) there is \textbf{no envelope}: two arrows pass through every point of the disc. The eigenvalues are real with opposite signs, and as in Proposition~\ref{prop:shapes}(3) the inner products for which $M$ is symmetric form the axis joining the two eigenvectors, which we denote by $\Gamma$. The matrix $M$ maps $\Gamma$ to itself and exchanges its two sides, so every arrow crosses $\Gamma$.

\begin{figure}[ht]
\centering
\includegraphics[width=\linewidth]{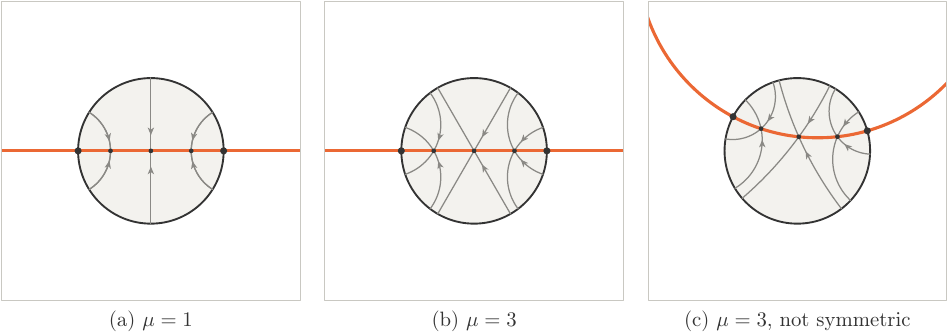}
\caption{Matrices with negative determinant, in the conformal chart. There is no envelope. Instead the axis $\Gamma$ (heavy line; a diameter in panels~a and~b), which consists of the inner products for which $M$ is symmetric, is crossed by the canonical arrows at a fixed angle. Each panel shows the two arrows through each of three points of $\Gamma$. In panel~(a), where $\mu=1$, the two arrows through a point are the same circle with opposite directions, and they cross $\Gamma$ at a right angle.}
\label{fig:negdet}
\end{figure}

To find the angle of crossing, use the upper half-plane and conjugate $M$ so that $f(u)=-\mu u$ with $\mu>0$, where $\mu$ is the ratio of the moduli of the two eigenvalues; then $\Gamma$ is the imaginary axis. The canonical arrow through $x$ and $-\mu x$ is the circle with centre $\frac{1-\mu}2x$ and radius $\frac{1+\mu}2|x|$. It crosses the imaginary axis at height $\sqrt\mu\,|x|$, and at the crossing point its tangent makes an angle $\varphi$ with the horizontal, that is, with the perpendicular to $\Gamma$, where
\begin{equation}\label{eq:theta}
\tan\varphi=\frac{|1-\mu|}{2\sqrt{\mu}}=\sinh\frac\ell2,\qquad \ell=|\ln\mu|,\qquad 0\leq\varphi<\frac{\pi}{2}.
\end{equation}
This angle does not depend on $x$. So every arrow crosses $\Gamma$ at the same angle, which is $\varphi$ away from a right angle, and every point of $\Gamma$ lies on two arrows, those of $x$ and of $-x$, which lean to opposite sides. Figure~\ref{fig:negdet} shows three examples. When $\tr M=0$ we have $\mu=1$: then $M$ acts on the disc as the reflection in $\Gamma$, and the arrows cross $\Gamma$ at right angles.

\subsection{The QR Algorithm}
For matrices with negative determinant, the picture of Section~\ref{subsec:negdet} is used. In the interactive program two arrowheads on the boundary mark the ends of the axis $\Gamma$, which are the eigenvectors; the user turns one of them, and the other moves with it.

\begin{figure}[ht]
\centering
\includegraphics{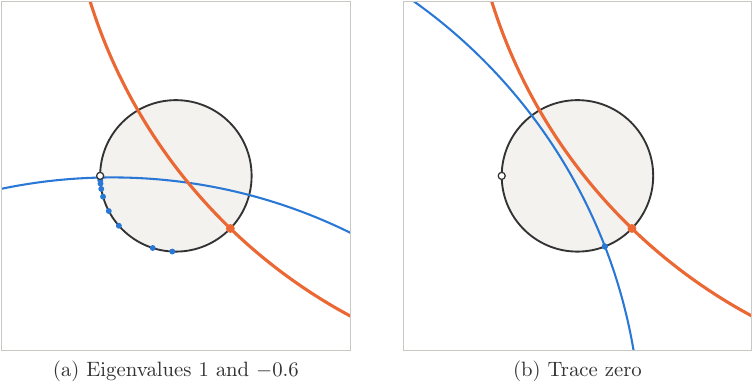}
\caption{The QR algorithm for two matrices with negative determinant. Heavy line: the axis $\Gamma$ of the input. Thinner line: the axis of the last iterate. Dots: the eigenvector whose eigenvalue has the larger modulus, after each step; the larger dot is its place in the input. (a)~Ten steps for eigenvalues $1$ and $-0.6$: the dots run to $(-1,0)$, from alternate sides, and the axis converges to a circle through that point. (b)~Trace zero, five steps, with the eigenvector of the positive eigenvalue marked: the axis alternates between two circles.}
\label{fig:qr-negdet}
\end{figure}

The algorithm acts on this picture much as before (Figure~\ref{fig:qr-negdet}). Since $\det M<0$ and $R_k$ has positive diagonal entries, $Q_k$ is now a reflection, so each step reflects the picture in a diameter instead of turning it. When the two eigenvalues differ in modulus, $\Gamma$ converges to a circle through $(-1,0)$. When they have the same modulus, that is, when the trace is zero, $M^2$ is a scalar matrix, and the algorithm returns to $M$ after two steps.

\section{Proofs, and Formulas for Drawing}
\label{sec:recipe}
\subsection{Proofs of the Statements of Sections~\ref{sec:turning} to~\ref{sec:qr-lr}}
\begin{proof}[Proof of Proposition~\ref{prop:qr-picture}]
By induction on $k$ we have $M_k=\tilde Q_k^TM\tilde Q_k$ and also $M^k=\tilde Q_k\tilde R_k$, where $\tilde R_k=R_{k-1}\cdots R_0$ is upper triangular with positive diagonal entries~\cite{watkins2008qr}; so $M^ke_1$ is a positive multiple of $\tilde Q_ke_1$. If $\det M>0$ then $\det\tilde Q_k>0$, so $\tilde Q_k$ is a rotation. As in Section~\ref{sec:turning}, the matrix $M_k=\tilde Q_k^{-1}M\tilde Q_k$ moves the direction of $\tilde Q_k^{-1}v$ to the direction of $\tilde Q_k^{-1}Mv$, and $\tilde Q_k^{-1}$ turns the boundary about its centre and takes the direction $x_k$ of $\tilde Q_ke_1$ to the direction of $e_1$.
\end{proof}

\begin{proof}[Proof of Theorem~\ref{thm:picture}]
(1) An arrow through $G$ is a direction along which the form $q(v)=v^T(GM+M^TG)v$ vanishes. Such a direction is never an eigenvector, since $q(v)=2\lambda\langle v,v\rangle_G\neq0$ when $Mv=\lambda v$. A non-zero quadratic form in two variables vanishes along $0$, $1$ or $2$ directions according as it is definite, of rank one or indefinite.

(2) Every $2\times2$ matrix $A$ satisfies $\det(A+A^T)=4\det A-(a_{12}-a_{21})^2$. With $A=GM$ this gives $\det(GM+M^TG)\leq4\det G\det M$, which is negative when $\det M<0$.

(3) Write $v=(t,1)^T$ and $GM+M^TG=\bigl(\begin{smallmatrix}\alpha&\beta\\ \beta&\gamma\end{smallmatrix}\bigr)$, where $\alpha,\beta,\gamma$ are linear in $G$. In the projective chart the arrow of $v$ lies on the line $\alpha t^2+2\beta t+\gamma=0$. The envelope of a family of lines that depends quadratically on a parameter $t$ is the discriminant conic $\beta^2=\alpha\gamma$, that is, $\det(GM+M^TG)=0$. At a boundary point $G=\xi\xi^T$ we have $GM+M^TG=\xi\eta^T+\eta\xi^T$ with $\eta=M^T\xi$, and its determinant is $-(\xi_1\eta_2-\xi_2\eta_1)^2\leq0$. It vanishes exactly when $\eta$ is a multiple of $\xi$, that is, when the kernel of $\xi^T$ is an eigenvector of $M$. Since the determinant does not change sign on the boundary, the conic touches the boundary at these points. The last statement follows from~(1).

(4) By~\eqref{eq:equivariance}, the Hermitian part of $gMg^{-1}$ for $g\cdot G$ is the conjugate by $g$ of the Hermitian part of $M$ for $G$; and $g$ maps geodesics to geodesics and the direction $x$ to $f_g(x)$.
\end{proof}

\begin{proof}[Proof of Proposition~\ref{prop:shapes}]
By Theorem~\ref{thm:picture}(4) we may replace $M$ by a conjugate, and $\mathcal D(M)$ is mapped to itself by every matrix that commutes with $M$.

(1) Let $M=rR_\theta$. It is $r$ times a rotation for the standard inner product, and it turns the disc about the centre by $-2\theta$, which is not a multiple of $2\pi$, so the centre is its only fixed point; an inner product for which $M$ is $r$ times a rotation is fixed by $M$. Since $M$ commutes with the rotations, $\mathcal D(M)$ is a union of circles about the centre, and it suffices to examine $G=\diag(1,\varepsilon)$ with $0<\varepsilon\leq1$, at distance $d=\frac12\ln(1/\varepsilon)$ from the centre. There
\[ \det(GM+M^TG)=r^2\bigl(4\varepsilon\cos^2\theta-(1-\varepsilon)^2\sin^2\theta\bigr), \]
which is positive exactly when $\sinh d=\frac{1-\varepsilon}{2\sqrt\varepsilon}<|\cot\theta|$, that is, when $\cosh d<1/\sin\theta$.

(2) Let $N=\bigl(\begin{smallmatrix}0&1\\0&0\end{smallmatrix}\bigr)$. Then $M$ commutes with the matrices $\bigl(\begin{smallmatrix}1&t\\0&1\end{smallmatrix}\bigr)$, which fix the boundary point $(-1,0)$ and move the other points along the horocycles based there, so $\mathcal D(M)$ is a union of such horocycles. For $G=\diag(1,\varepsilon)$ we find $\det(GM+M^TG)=\lambda^2(4\varepsilon-1)$ and $\tr(N^*N)=1/\varepsilon$, which is the square of the operator norm because $N$ has rank one.

(3) Let $M=\diag(\lambda_1,\lambda_2)$. It is symmetric for $G$ exactly when $G$ is diagonal, and the diagonal inner products form the geodesic from $(-1,0)$ to $(1,0)$. Since $M$ commutes with the diagonal matrices, which map this axis to itself, $\mathcal D(M)$ is a union of equidistant curves. The inner products $G=\bigl(\begin{smallmatrix}1&g\\ g&1\end{smallmatrix}\bigr)$, $|g|<1$, form the geodesic through the centre perpendicular to the axis, and $G$ is at distance $\delta$ from the axis, where $\tanh\delta=|g|$. Now $\det(GM+M^TG)=4\lambda_1\lambda_2-g^2(\lambda_1+\lambda_2)^2$ is positive exactly when
\[ \sinh^2\delta=\frac{g^2}{1-g^2}<\frac{4\lambda_1\lambda_2}{(\lambda_1-\lambda_2)^2}=\frac{1}{\sinh^2(\ell/2)}. \qedhere \]
\end{proof}

\begin{proof}[Proof of Proposition~\ref{prop:qr-lr}]
(1) As in Proposition~\ref{prop:qr-picture}, $M_k=\tilde Q_k^{-1}M\tilde Q_k$ and $M^k=\tilde Q_k\tilde R_k$, where now $\tilde Q_k$ is orthogonal for $G$, with positive determinant, so it acts as a rotation about the point $G$; and $\tilde Q_k^{-1}$ takes the direction $x_k$ of $M^ke_1$ to that of $e_1$.

(2) The value of $w(G_\varepsilon)$ follows from~\eqref{eq:charts}. Let $D=\diag(1,\sqrt\varepsilon)$, so that $G_\varepsilon=D^TD$ and $Q_0=D^{-1}Q'D$, where $Q'$ is the orthogonal factor of $DMD^{-1}$. The first column of $Q'$ is the unit vector along $DMe_1=(m_{11},\sqrt\varepsilon\,m_{21})^T$, so $Q'=\pm R_\varphi$ with $\tan\varphi=\sqrt\varepsilon\,m_{21}/m_{11}$, and
\[ Q_0=\pm\begin{pmatrix}\cos\varphi&-\sqrt\varepsilon\,\sin\varphi\\[2pt] \sin\varphi/\sqrt\varepsilon&\cos\varphi\end{pmatrix}\ \longrightarrow\ \pm\begin{pmatrix}1&0\\ m_{21}/m_{11}&1\end{pmatrix}=\pm L\qquad(\varepsilon\to0^+). \]
Then $Q_0^{-1}MQ_0\to L^{-1}ML=UL$.

(3) The M\"obius transformation of $L$ is $x\mapsto x/(lx+1)$ with $l=m_{21}/m_{11}$. It fixes only $x=0$, the direction of $e_2$, so it is parabolic and preserves the horocycles based at $(1,0)$. Finally, the direction of $Me_1$ has the label $m_{11}/m_{21}$, which is $0$ exactly when $m_{11}=0$.
\end{proof}

\subsection{The Projective and the Conformal Chart as Group Actions}
Visualisation v1.5 is the picture that comes from the accidental isomorphism, one special to small dimensions,
\[\PSL_2(\RR)\;\cong\;\SO^+(2,1).\]
Here $\PSL_2(\RR)=\SL_2(\RR)/\{\pm I\}$ is the group of M\"obius transformations of $\RP^1$, and a matrix of positive determinant determines such a transformation up to a non-zero scalar factor. $\SO^+(2,1)$ is the identity component of the group of linear maps of $\RR^3$ preserving a quadratic form of signature $(2,1)$; here $\RR^3$ is the space of symmetric $2\times2$ matrices and the form is $-\det$. The group acts on the cone $\det G=0$, whose rays form the boundary circle, and on the projective chart. Visualisation v2 comes in the same way from the accidental isomorphism
\[ \PSL_2(\RR)\;\cong\;\PU(1,1), \]
where $\PU(1,1)$ is the group of complex $2\times2$ matrices preserving the Hermitian form $|z_1|^2-|z_2|^2$, modulo scalars; it acts by M\"obius transformations on the unit disc, the conformal chart, and on its boundary circle. By Proposition~\ref{prop:shapes} and the second of the facts in Section~\ref{sec:disc}, the envelope in the conformal chart consists of arcs of circles; continued outside the disc, they are the two cycles of Section~\ref{subsec:v2}.

\subsection{The Cycles in Lie Sphere Geometry}
Formula~\eqref{eq:cycles} comes from \emph{Lie sphere geometry}~\cite{cecil1992lie}, which was the starting point of the preprint. Lie sphere geometry also helps to understand the simple Lie group $\Sp(4,\RR)$, through the accidental isomorphism $\Sp(4,\RR)/\{\pm I\}\cong\SO^+(3,2)$, in which the right-hand side is the identity component of the group of transformations of Lie sphere geometry; we do not use this.

\begin{definition}[Lie quadric]
\emph{Oriented cycles}, that is, oriented circles, oriented lines and points (including $\infty$) on the sphere $S^2\cong\widehat\CC$, are represented by the points of the quadric
\[
\LieQuad:\quad X^2+Y^2+Z^2=C^2+S^2 \qquad\text{in } \RP^4,
\]
with homogeneous coordinates $[X:Y:Z:C:S]$~\cite{cecil1992lie}. If the left-hand side vanishes then so does the right-hand side, and $[0:0:0:0:0]$ is not a point of $\RP^4$; so we may normalise $X^2 + Y^2 + Z^2 = C^2 + S^2=1$. Then $(X,Y,Z)\in S^2$ is the spherical centre of the cycle, and $(C,S)=(\cos r,\sin r)$ encodes its signed spherical radius $r$, the sign of $r$ giving the orientation.
\end{definition}

Let
\[
Q(\xi)=X^2+Y^2+Z^2-C^2-S^2,\qquad
\mathcal B(\xi,\eta)=XX'+YY'+ZZ'-CC'-SS'
\]
be the quadratic form, of signature $(3,2)$, and its polarisation. Two cycles $p,q\in\LieQuad$ are \emph{orientedly tangent} if and only if $\mathcal B(p,q)=0$.

\begin{remark}
The condition $\mathcal B(p,q)=0$ is a reasonable definition of oriented tangency. It says that $XX'+YY'+ZZ'=CC'+SS'$, an equality of two dot products of unit vectors, which by $u\cdot v=\cos\angle(u,v)$ means that the spherical distance between the centres equals the difference of the signed radii. A point is a cycle of radius zero, and it is orientedly tangent to a cycle exactly when it lies on it.
\end{remark}

Let $M=\bmat{a & b\\ c & d}$ with $\Delta=\det M>0$ and $\tau=\tr M$. We identify the sphere $S^2$ with $\widehat\CC$, in which the conformal chart is drawn, by the stereographic projection $w=(X+iY)/(1-Z)$; the boundary $|w|=1$ is then the equator $Z=0$.

\begin{lemma}[The cycles of a matrix]\label{lem:pencil}
The two cycles of $M$ are the points
\begin{equation}\label{eq:lie-pair}
[X:Y:Z:C:S]=\bigl[-(b+c):\; d-a:\; \pm\,2\sqrt{\Delta}:\; b-c:\; a+d\bigr]
\end{equation}
of the Lie quadric. In the conformal chart, the cycle $[X:Y:Z:C:S]$ is the circle with
\begin{equation}\label{eq:disc-circle}
\text{centre }\frac{X+iY}{C-Z}\qquad\text{and signed radius }\frac{S}{C-Z},
\end{equation}
which is a straight line when $C=Z$; the radius is positive when the cycle runs anticlockwise.
\end{lemma}

\begin{proof}
Let $\det G=1$ and put $J_G=G^{-1}J$ with $J=R_{\pi/2}$; then $J_G^2=-I$, and $J_G$ is the rotation by a right angle for $G$. Since $a_{12}-a_{21}=\tr(JA)$ for every $2\times2$ matrix $A$, and $JG=G^{-1}J$, the identity in the proof of Theorem~\ref{thm:picture}(2) becomes
\begin{equation}\label{eq:trace-form}
\det(GM+M^TG)=4\Delta-\tr(J_GM)^2 .
\end{equation}
So the envelope consists of the two curves $\tr(J_GM)=\mp2\sqrt\Delta$. Write $G=\bigl(\begin{smallmatrix}T+\xi&-\eta\\ -\eta&T-\xi\end{smallmatrix}\bigr)$ with $T^2-\xi^2-\eta^2=1$. By~\eqref{eq:charts}, $w(G)=(\xi+i\eta)/(1+T)$, which is the stereographic projection of the point $(\xi,\eta,-1)/T$ of the sphere, the point cycle $[\xi:\eta:-1:T:0]$. A direct computation gives $\tr(J_GM)=T(b-c)+\xi(b+c)-\eta(d-a)$. Hence, with $X,Y,Z,C$ as in~\eqref{eq:lie-pair}, the equation $\tr(J_GM)=-Z$ reads $X\xi+Y\eta-Z-CT=0$, which says that the point $w(G)$ lies on the cycle $[X:Y:Z:C:S]$. The equation of the quadric determines $S$ up to sign, $S^2=X^2+Y^2+Z^2-C^2=(a+d)^2$, and the sign of $S$ is the orientation. Formula~\eqref{eq:disc-circle} is the stereographic projection of the circle on the sphere. That $S=a+d$, with the sign convention of~\eqref{eq:disc-circle}, is the direction of the arrows that touch the cycle was verified numerically for many random matrices; for the rotation $M=R_\theta$ the two cycles are the circles $|w|=|\tan(\frac\pi4\pm\frac\theta2)|$, in agreement with Proposition~\ref{prop:elliptic} below.
\end{proof}

The coordinates in~\eqref{eq:lie-pair} are simply the trace-free part of $M$, in $X$, $Y$ and $C$, together with $\pm2\sqrt{\det M}$ and $\tr M$. The heavy circles in the figures of this appendix are drawn from \eqref{eq:lie-pair} and \eqref{eq:disc-circle}, which together are formula~\eqref{eq:cycles}.

\subsection{The Two Cycles as Level Sets}
\label{subsec:levels}
The cycles can also be described in the upper half-plane, where the boundary is $\RR$, as level sets of a M\"obius transformation $h$ adapted to the fixed points of $f$. Since M\"obius transformations map circles and lines to circles and lines, every level set $\{|h(z)|=\rho\}$ and $\{\Arg h(z)=\kappa\}$ is a circle or a line.

\begin{proposition}[Hyperbolic case]\label{prop:hyperbolic}
Let $\Delta>0$ and $\tau^2>4\Delta$, let $p,q\in\RR$ be the fixed points of $f$, and let $h(z)=\frac{z-p}{z-q}$. The two cycles of $M$ are the level sets $\{\Arg h(z)=\pm\kappa_H\}$ (modulo $\pi$), where
\[ \tan\kappa_H=\sinh\frac\ell2=\frac{\sqrt{\tau^2-4\Delta}}{2\sqrt\Delta}; \]
they pass through $p$ and $q$ and meet the boundary at the angle $\kappa_H$.
\end{proposition}
\begin{proof}
The map $h$ sends $p$ to $0$ and $q$ to $\infty$, so it sends the axis to the imaginary axis. The curves equidistant from the imaginary axis are the rays from $0$, and the rays at the angle $\kappa$ to the real axis are at distance $\delta$ with $\sinh\delta=\cot\kappa$. So Proposition~\ref{prop:shapes}(3) gives $\tan\kappa_H=\sinh\frac\ell2$, and $\sinh\frac\ell2=\frac{|\lambda_1-\lambda_2|}{2\sqrt{\lambda_1\lambda_2}}$.

This can be checked directly. In the chart given by $h$ we have $f(x) = \lambda x$, where $\lambda>0$ is the ratio of the eigenvalues. The canonical arrow from $x$ to $\lambda x$ is the circle with centre $\frac{1+\lambda}{2}x$ and radius $\frac{|\lambda-1|}{2}|x|$. A line $y=mx$ is tangent to it if and only if its distance from the centre equals the radius, that is, $m^2(1+\lambda)^2=(\lambda-1)^2(1+m^2)$, which gives $m=\pm\frac{\lambda-1}{2\sqrt\lambda}$, independently of $x$. Hence the envelope is the pair of lines $y = \pm\frac{\lambda - 1}{2\sqrt\lambda} x$, as in the following diagram, in which the circle is the canonical arrow from $1$ to $f(1) = \lambda$ and the heavy lines are the envelope:
\begin{center}%
    \includegraphics{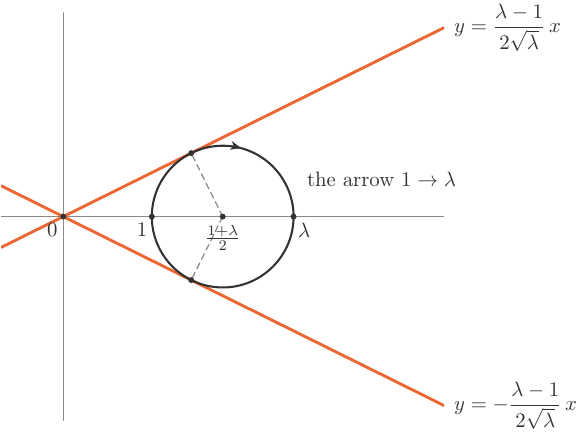}%
\end{center}
One of the two lines, directed as the arrows that touch it, is obtained from the other by reflection in the real axis followed by reversal of direction.
\end{proof}

\begin{proposition}[Elliptic case]\label{prop:elliptic}
Let $\Delta>0$ and $\tau^2<4\Delta$, write $\cos\theta=\tau/(2\sqrt{\Delta})$ with $0<\theta<\pi$, let $p$ be the fixed point of $f$ in the upper half-plane, and let $h(z)=\frac{z-p}{z-\bar p}$. Then the two cycles of $M$ are the circles
\[
|h(z)|=\rho^{\pm1},\qquad \rho=\Bigl|\tan\Bigl(\frac\pi4-\frac\theta2\Bigr)\Bigr| .
\]
\end{proposition}
\begin{proof}
The map $h$ sends the upper half-plane to the unit disc and $p$ to $0$, so it sends the circle of radius $d$ about $p$ to the circle $|h|=\tanh\frac d2$. By Proposition~\ref{prop:shapes}(1), $\cosh d=1/\sin\theta$, and then $\tanh\frac d2=\sqrt{(1-\sin\theta)/(1+\sin\theta)}=|\tan(\frac\pi4-\frac\theta2)|$. The other cycle is the inversion $|h|=1/\rho$. Directly: $h$ conjugates $f$ to a rotation about $0$ by an angle $\pm2\theta$, for which the canonical arrow joining $e^{-i\theta}$ to $e^{i\theta}$ is the circle with centre $\sec\theta$ and radius $|\tan\theta|$, whose closest point to the origin is at distance $|\sec\theta-\tan\theta|=\rho$.
\end{proof}
The identity matrix corresponds to $\theta\to0$, $\rho\to1$, where both cycles approach the boundary. An orthogonal matrix with $\theta=\pi/2$, which is a rotation by a right angle up to scaling, gives $\rho=0$, where the inner cycle degenerates to the point $p$.

\subsection{Practical Recipe}
\begin{enumerate}
\item Given $M=\bmat{a&b\\ c&d}$, compute $\Delta=ad-bc$, $\tau=a+d$, and the fixed points of $f(z)=\frac{az+b}{cz+d}$, which are the roots of $cz^2+(d-a)z-b=0$.
\item If $\Delta>0$, the two cycles are $[-(b+c):\,d-a:\,\pm2\sqrt{\Delta}:\,b-c:\,a+d]$, and \eqref{eq:disc-circle} gives them directly in the conformal chart. Alternatively:
  \begin{itemize}
  \item \textbf{Elliptic:} with the fixed point $p$ in the upper half-plane, compute $\theta$ from $\cos\theta=\tau/(2\sqrt{\Delta})$, set $\rho=|\tan(\frac\pi4-\frac\theta2)|$, and plot $\bigl|\frac{z-p}{z-\bar p}\bigr|=\rho^{\pm1}$.
  \item \textbf{Hyperbolic:} with the real fixed points $p,q$, compute $\kappa_H$ from $\tan\kappa_H=\sqrt{\tau^2-4\Delta}/(2\sqrt{\Delta})$ and plot $\Arg\frac{z-p}{z-q}=\pm\kappa_H$.
  \end{itemize}
  Then map the picture to the disc by $w=\frac{iz+1}{-iz+1}$.
\item If $\Delta<0$, plot the axis $\Gamma$, which is the circle through the two fixed points orthogonal to the boundary, and through sample points of $\Gamma$ draw the two canonical arrows, which cross $\Gamma$ symmetrically at the angle determined by \eqref{eq:theta}.
\end{enumerate}

\section{Exercises}
\label{sec:exercises}
The answers follow the exercises.
\begin{enumerate}
\item\label{ex:point} Where on the boundary is the direction of $(1,3)^T$? Which direction is drawn at the point $(0.8,-0.6)$?
\item\label{ex:power} For the running example, compute $M^4e_1$ and the point $x_4$. By what factor does the entry below the diagonal of $M_k$ shrink from $M_1$ to $M_2$ and from $M_2$ to $M_3$?
\item\label{ex:shear} Compute the two cycles of the shear $\bigl(\begin{smallmatrix}1&1\\0&1\end{smallmatrix}\bigr)$ from~\eqref{eq:cycles}. Where do they meet the boundary, and what does the QR algorithm do to this matrix?
\item\label{ex:transpose} Compute the cycles of the transpose of the running example, and compare them with the cycles of the running example. What is the rule?
\item\label{ex:elliptic} For $M=\bigl(\begin{smallmatrix}2&-2\\1&2\end{smallmatrix}\bigr)$, find an inner product for which $M$ is a multiple of a rotation, its point $p$ in the disc, and the radius $d$ of the region $\mathcal D(M)$.
\item\label{ex:diameter} For the running example, for which $\varepsilon>0$ is $M+M^*$ definite when the inner product is $\diag(1,\varepsilon)$? Compare the answer with the points where the two cycles cross the horizontal diameter.
\item\label{ex:lr} For the running example, carry out a second step of the LR algorithm. Which entry shows that the iterates converge, and at what rate?
\end{enumerate}

\paragraph{Answers.}
(\ref{ex:point}) By~\eqref{eq:direction-point}, at $(0.8,0.6)$. The point $(0.8,-0.6)$ is the direction of $(1,-3)^T$.
(\ref{ex:power}) $M^4e_1=(171,170)^T$, and $x_4=(-0.006,1.000)$. The entries are $0.615$, $0.145$, $0.035$, so the factors are $0.24$ and $0.24$, close to the ratio $\frac14$ of the eigenvalues.
(\ref{ex:shear}) $X=-1$, $Y=0$, $C=1$, $S=2$, $Z=\pm2$: the circle with centre $(1,0)$ and radius $2$, and the circle with centre $(-\frac13,0)$ and radius $\frac23$. Both touch the boundary at $(-1,0)$, the direction of the only eigenvector $e_1$. The matrix is upper triangular, so it is a fixed point of the QR algorithm.
(\ref{ex:transpose}) The cycles of $M^T$ are the circle with centre $(-0.6,-0.2)$ and radius $1$, anticlockwise, and the circle with centre $(1,\frac13)$ and radius $\frac53$, clockwise. The picture of $M^T$ is the picture of $M$ turned by a half-turn about the centre, with both cycles reversed. The same holds for every $M$, because transposition changes the sign of $C$ in~\eqref{eq:cycles} and nothing else.
(\ref{ex:elliptic}) The fixed point of $f$ in the upper half-plane is $z=i\sqrt2$, which gives $G=\diag(1,2)$; indeed $M^TGM=6G$. Its point is $p=\bigl(\frac{1-\sqrt2}{1+\sqrt2},0\bigr)=(-0.172,0)$. The eigenvalues are $2\pm i\sqrt2=\sqrt6\,e^{\pm i\theta}$ with $\sin\theta=1/\sqrt3$, so $\cosh d=\sqrt3$ and $d=1.146$.
(\ref{ex:diameter}) $GM+M^TG=\bigl(\begin{smallmatrix}6&1+2\varepsilon\\1+2\varepsilon&4\varepsilon\end{smallmatrix}\bigr)$ has determinant $-(4\varepsilon^2-20\varepsilon+1)$, which is positive for $0.0505<\varepsilon<4.95$. The corresponding points $\bigl(\frac{1-\sqrt\varepsilon}{1+\sqrt\varepsilon},0\bigr)$ run from $(0.633,0)$ to $(-0.380,0)$, and these are the points where the cycle with centre $(-1,-\frac13)$ and the cycle with centre $(0.6,0.2)$ cross the diameter inside the disc.
(\ref{ex:lr}) The second step gives $\bigl(\begin{smallmatrix}3.909&1\\0.264&1.091\end{smallmatrix}\bigr)$. The entry below the diagonal goes $2$, $0.889$, $0.264$, and in the limit it shrinks by the factor $\frac14$ at each step; the limit is $\bigl(\begin{smallmatrix}4&1\\0&1\end{smallmatrix}\bigr)$.

\begin{table}[ht]
\centering
\caption{The matrices drawn in the figures of this appendix, with entries rounded to three decimal places, and their eigenvalues. Figures~\ref{fig:directions}(b) and~\ref{fig:qrstep} show the running example.}
\label{tab:matrices}
\footnotesize
\begin{tabular}{@{}lll@{\hspace{3em}}lll@{}}
\toprule
Figure & Matrix & Eigenvalues & Figure & Matrix & Eigenvalues \\
\midrule\addlinespace[6pt]
\ref{fig:v1}(a) & $\begin{pmatrix}0.658&0.267\\0.267&0.792\end{pmatrix}$ & $1,\ 0.45$ & \ref{fig:cases123}(a) & $\begin{pmatrix}3.000&1.000\\2.000&2.000\end{pmatrix}$ & $4,\ 1$ \\\addlinespace[4pt]
\ref{fig:v1}(b) & $\begin{pmatrix}0.983&0.094\\0.094&0.467\end{pmatrix}$ & $1,\ 0.45$ & \ref{fig:cases123}(b) & $\begin{pmatrix}5.000&1.000\\{-}1.000&3.000\end{pmatrix}$ & $4,\ 4$ \\\addlinespace[4pt]
\ref{fig:v1}(c) & $\begin{pmatrix}0.467&0.094\\0.094&0.983\end{pmatrix}$ & $1,\ 0.45$ & \ref{fig:cases123}(c) & $\begin{pmatrix}2.000&{-}2.000\\1.000&2.000\end{pmatrix}$ & $2\pm1.414\,i$ \\\addlinespace[4pt]
\ref{fig:v1}(d) & $\begin{pmatrix}0.913&0.068\\0.068&0.947\end{pmatrix}$ & $1,\ 0.86$ & \ref{fig:cases45}(a) & $\begin{pmatrix}1.864&1.106\\0.106&1.136\end{pmatrix}$ & $2,\ 1$ \\\addlinespace[4pt]
\ref{fig:v1}(e) & $\begin{pmatrix}0.331&0.418\\0.418&0.739\end{pmatrix}$ & $1,\ 0.07$ & \ref{fig:cases45}(b) & $\begin{pmatrix}0.568&{-}0.553\\{-}0.053&0.932\end{pmatrix}$ & $1,\ 0.5$ \\\addlinespace[4pt]
\ref{fig:models} & $\begin{pmatrix}3.000&1.000\\2.000&2.000\end{pmatrix}$ & $4,\ 1$ & \ref{fig:cases678}(a) & $\begin{pmatrix}0.904&{-}0.017\\{-}0.137&0.976\end{pmatrix}$ & $1,\ 0.88$ \\\addlinespace[4pt]
\ref{fig:shapes}(a) & $\begin{pmatrix}2.000&{-}2.000\\1.000&2.000\end{pmatrix}$ & $2\pm1.414\,i$ & \ref{fig:cases678}(b) & $\begin{pmatrix}0.204&{-}0.140\\{-}1.137&0.799\end{pmatrix}$ & $1,\ 0.003$ \\\addlinespace[4pt]
\ref{fig:shapes}(b) & $\begin{pmatrix}5.000&1.000\\{-}1.000&3.000\end{pmatrix}$ & $4,\ 4$ & \ref{fig:cases910}(a) & $\begin{pmatrix}0.921&{-}0.691\\0.538&0.682\end{pmatrix}$ & $0.802\pm0.598\,i$ \\\addlinespace[4pt]
\ref{fig:shapes}(c) & $\begin{pmatrix}3.000&1.000\\2.000&2.000\end{pmatrix}$ & $4,\ 1$ & \ref{fig:cases910}(b) & $\begin{pmatrix}0.805&{-}0.645\\0.643&0.727\end{pmatrix}$ & $0.766\pm0.643\,i$ \\\addlinespace[4pt]
\ref{fig:negdet}(a) & $\begin{pmatrix}1.000&0\\0&{-}1.000\end{pmatrix}$ & $1,\ {-}1$ & \ref{fig:qr-to-lr} & $\begin{pmatrix}3.000&1.000\\2.000&2.000\end{pmatrix}$ & $4,\ 1$ \\\addlinespace[4pt]
\ref{fig:negdet}(b) & $\begin{pmatrix}1.500&0\\0&{-}0.500\end{pmatrix}$ & $1.5,\ {-}0.5$ & \ref{fig:qr-negdet}(a) & $\begin{pmatrix}{-}0.297&{-}0.524\\{-}0.749&0.697\end{pmatrix}$ & $1,\ {-}0.6$ \\\addlinespace[4pt]
\ref{fig:negdet}(c) & $\begin{pmatrix}1.573&{-}0.291\\0.517&{-}0.573\end{pmatrix}$ & $1.5,\ {-}0.5$ & \ref{fig:qr-negdet}(b) & $\begin{pmatrix}{-}0.622&{-}0.655\\{-}0.936&0.622\end{pmatrix}$ & $1,\ {-}1$ \\\addlinespace[4pt]
\ref{fig:key} & $\begin{pmatrix}3.000&1.000\\2.000&2.000\end{pmatrix}$ & $4,\ 1$ &  &  &  \\\addlinespace[4pt]
\bottomrule
\end{tabular}
\end{table}

\paragraph{Acknowledgements for this appendix.} For the work in this appendix, the author would like to thank Jan Techter, Gregory Gutin and Abbas Edalat for helpful discussions.
%======================================================================
%   APPENDIX B: SOURCE CODE OF THE LIBRARY
%   Generated by gen_code_appendix.py from code_appendix.tmpl and the
%   source files of dsqlib 0.2.0. Edit the template, not this file.
%======================================================================
\chapter{Source Code of the Library}
\label{chap:code}

This appendix prints the source code of the package \texttt{dsqlib}, version~0.2.0, the library described in Chapter~\ref{chap:implementation}. The test suite, from which the measurements reported in Chapter~\ref{chap:implementation} are taken, is not printed: it is part of the source distribution of the library, which is submitted with this thesis as supplementary material, together with a record of a run of the tests. The modules of the package are printed in the order of Table~\ref{tab:implementation}, preceded by \texttt{\_\_init\_\_.py}, which defines the public interface, and followed by \texttt{errors.py}, which defines the exceptions.

Each heading gives the path of a file relative to the root of the library's source tree. The library requires Python~3.10 or later, NumPy~1.24 or later and SciPy~1.10 or later. In the source distribution, the tests are run with pytest, with the directories \texttt{src} and \texttt{tests} on the Python path. Line numbers, shown on every fifth line, are those of the file; a line too long for the page is broken, and its continuation is marked by~$\hookrightarrow$.

The code is released under the MIT licence, which reads as follows.
\begin{lstlisting}[style=dsqtext]
MIT License

Copyright (c) 2026 Ran Gutin

Permission is hereby granted, free of charge, to any person obtaining a copy
of this software and associated documentation files (the "Software"), to deal
in the Software without restriction, including without limitation the rights
to use, copy, modify, merge, publish, distribute, sublicense, and/or sell
copies of the Software, and to permit persons to whom the Software is
furnished to do so, subject to the following conditions:

The above copyright notice and this permission notice shall be included in all
copies or substantial portions of the Software.

THE SOFTWARE IS PROVIDED "AS IS", WITHOUT WARRANTY OF ANY KIND, EXPRESS OR
IMPLIED, INCLUDING BUT NOT LIMITED TO THE WARRANTIES OF MERCHANTABILITY,
FITNESS FOR A PARTICULAR PURPOSE AND NONINFRINGEMENT. IN NO EVENT SHALL THE
AUTHORS OR COPYRIGHT HOLDERS BE LIABLE FOR ANY CLAIM, DAMAGES OR OTHER
LIABILITY, WHETHER IN AN ACTION OF CONTRACT, TORT OR OTHERWISE, ARISING FROM,
OUT OF OR IN CONNECTION WITH THE SOFTWARE OR THE USE OR OTHER DEALINGS IN THE
SOFTWARE.
\end{lstlisting}

\addtocontents{toc}{\protect\setcounter{tocdepth}{1}}
\section{The Package}
\label{sec:code-package}
\subsection{\texorpdfstring{\texttt{src/dsqlib/\_\_init\_\_.py}}{src/dsqlib/\_\_init\_\_.py}}
\begin{lstlisting}[style=dsqcode]
"""dsqlib: *-algebras of doubly-signed quivers, localised at the arrows and completed in path degree.

The main objects:

* :class:`DSQ` -- a doubly-signed quiver (``DSQ.build`` from the vertex kinds of Table 5.1);
* :class:`Algebra` and :class:`Element` -- the completed localisation of its path *-algebra;
* :func:`star_sqrt` -- *-square roots by the graded, principal and arbitrary rules;
* :class:`Index` and :class:`Matrix` -- matrices indexed by self-adjoint idempotents;
* :mod:`dsqlib.linalg` -- Householder, QR, Lanczos, Cholesky, LDL and LU, written once for every DSQ;
* :mod:`dsqlib.catalogue` and :mod:`dsqlib.problems` -- the quivers of the thesis catalogue and the
  classical matrix problems they encode.

Example
-------
>>> import numpy as np
>>> from dsqlib import problems, linalg
>>> H = problems.hermitian_problem(np.array([[2.0, 1j], [-1j, 3.0]]))
>>> T, U, _ = linalg.householder(H)
"""

from . import catalogue, linalg, problems
from .algebra import INF, Algebra, Element, quaternion
from .errors import BreakdownError, DSQError, NoSquareRootError, NotInvertibleError, PrecisionError, QuiverError
from .matrix import Index, Matrix, cyclic_positions, orbit_index, plain_index
from .quiver import DSQ, KINDS, SIGNATURE_KINDS
from .sqrt import random_semisimple_unitary, random_unitary, sqrt_graded, sqrt_principal, sqrt_rule, star_sqrt

__version__ = "0.2.0"

__all__ = [
    "DSQ",
    "INF",
    "KINDS",
    "SIGNATURE_KINDS",
    "Algebra",
    "BreakdownError",
    "DSQError",
    "Element",
    "Index",
    "Matrix",
    "NoSquareRootError",
    "NotInvertibleError",
    "PrecisionError",
    "QuiverError",
    "__version__",
    "catalogue",
    "cyclic_positions",
    "linalg",
    "orbit_index",
    "plain_index",
    "problems",
    "quaternion",
    "random_semisimple_unitary",
    "random_unitary",
    "sqrt_graded",
    "sqrt_principal",
    "sqrt_rule",
    "star_sqrt",
]
\end{lstlisting}
\needspace{6\baselineskip}
\subsection{\texorpdfstring{\texttt{src/dsqlib/quiver.py}}{src/dsqlib/quiver.py}}
\begin{lstlisting}[style=dsqcode]
"""Doubly-signed quivers (thesis, Definition 5.7), stored at the complexified level.

A doubly-signed quiver (DSQ) has vertices and arrows with two commuting involutions:

* the transpose ``T``, which reverses arrows (``a: i -> j`` gives ``a^T: j^T -> i^T``);
* the conjugation ``B``, which preserves them (``a: i -> j`` gives ``a~: i~ -> j~``);

and signs:

* ``eps_div`` = +-1 on the vertices fixed by ``B`` (real or quaternionic), 0 elsewhere;
* ``eps_form`` = +-1 on the vertices and arrows fixed by ``T``, 0 elsewhere.

Most users build a DSQ from the ten real vertex kinds of Table 5.1 with :meth:`DSQ.build`. Each kind
expands into one, two or four complexified vertices. Each arrow given is closed under the Klein
four-group {1, T, B, TB}; its partners are named ``a'`` (transpose), ``a~`` (conjugate) and ``a~'``
(both), unless they coincide with ``a``.

Block model
-----------
Each complexified vertex ``v`` carries a block size ``d_v`` (1 or 2) and two ``d_v x d_v`` matrices:

* ``G_v``, the vertex form: on ``S_v = M_{d_v}(C)`` the involution is ``X -> G_v X^T G_v^-1``;
* ``C_v``, the real structure: the conjugation is ``X -> C_v conj(X) C_v^-1``.

Each arrow carries signs ``kappa``, ``mu`` = +-1 that fix the involution and the conjugation on its
bimodule: ``(Y a)* = kappa_a G_s Y^T G_t^-1 a^T`` and ``c(Y a) = mu_a C_t conj(Y) C_s^-1 a~``. They are
found by a finite search such that ``** = 1``, ``cc = 1`` and ``c* = *c``; thesis Proposition 5.14(2)
gives an explicit solution, so the search always succeeds.
"""

from __future__ import annotations

import itertools
from collections.abc import Iterable, Mapping, Sequence
from typing import Any, NamedTuple

import numpy as np

from .errors import QuiverError

__all__ = ["DSQ", "KINDS", "SIGNATURE_KINDS", "Arrow", "Component", "Vertex", "VertexKind"]

_J2 = np.array([[0, 1], [-1, 0]], dtype=complex)


class VertexKind(NamedTuple):
    """One of the ten real vertex kinds of Table 5.1."""

    pattern: str
    """How the kind expands: 'fixed', 'cfix', 'cu', 'exch' or 'cc'."""
    eps_div: int
    eps_form: int
    algebra: str
    """The simple *-algebra of minimal size."""
    unitary_group: str
    """The unitary group of ``M_n`` of that algebra."""


KINDS: dict[str, VertexKind] = {
    "R": VertexKind("fixed", +1, +1, "R, trivial involution", "O(n)"),
    "Rsp": VertexKind("fixed", +1, -1, "M2(R), symplectic involution", "Sp(2n,R)"),
    "Hsp": VertexKind("fixed", -1, -1, "H, quaternion conjugation", "Sp(n)"),
    "Hor": VertexKind("fixed", -1, +1, "H, involution of orthogonal type", "O*(2n)"),
    "Cor": VertexKind("cfix", 0, +1, "C, trivial involution", "O(n,C)"),
    "Csp": VertexKind("cfix", 0, -1, "M2(C), symplectic involution", "Sp(2n,C)"),
    "Cu": VertexKind("cu", 0, 0, "C, complex conjugation", "U(n)"),
    "RR": VertexKind("exch", +1, 0, "R x R^op, exchange involution", "GL_n(R)"),
    "HH": VertexKind("exch", -1, 0, "H x H^op, exchange involution", "GL_n(H)"),
    "CC": VertexKind("cc", 0, 0, "C x C^op, exchange involution", "GL_n(C)"),
}
"""The ten vertex kinds of Table 5.1, by their short names."""

SIGNATURE_KINDS: tuple[str, ...] = ("R", "Hsp", "Cu")
"""The vertex kinds with a signature (compact unitary group)."""

_GROUP = ("1", "T", "B", "TB")  # the Klein four-group generated by transpose and conjugation
_SUFFIX = {"1": "", "T": "'", "B": "~", "TB": "~'"}


def _group_mul(g: str, h: str) -> str:
    """Product in the Klein group (it is abelian)."""
    s = set(g.replace("1", "")) ^ set(h.replace("1", ""))
    return {frozenset(): "1", frozenset("T"): "T", frozenset("B"): "B", frozenset("TB"): "TB"}[frozenset(s)]


class Vertex:
    """A complexified vertex.

    Attributes
    ----------
    name, top, bar : str
        The vertex and its partners under the transpose and the conjugation.
    eps_div, eps_form : int
        The signs of Definition 5.7 (0 where they are not defined).
    orbit : str or None
        The base name given to :meth:`DSQ.build` (names the STAR-BAR orbit).
    kind : str
        The vertex kind (a key of :data:`KINDS`).
    block_size : int
        ``d_v``, 1 or 2 (set when the vertex joins a :class:`DSQ`).
    G, C : numpy.ndarray
        The vertex form and the real structure (set when the vertex joins a :class:`DSQ`).
    """

    def __init__(
        self,
        name: str,
        top: str,
        bar: str,
        eps_div: int,
        eps_form: int,
        orbit: str | None = None,
        kind: str | None = None,
    ) -> None:
        self.name, self.top, self.bar = name, top, bar
        self.eps_div, self.eps_form = eps_div, eps_form
        self.orbit = orbit
        self.kind: str = kind or ""
        self.block_size = 1
        self.G = self.C = self.G_inv = self.C_inv = np.eye(1, dtype=complex)
        self._g = self._h = 1  # the scalars G_{v^T} G_v^-T and C_{v~} conj(C_v)

    def __repr__(self) -> str:
        return f"Vertex({self.name!r}, kind={self.kind!r})"


class Arrow:
    """A complexified arrow ``source -> target``.

    Attributes
    ----------
    name, source, target : str
    top, bar : str
        The names of the transposed and the conjugate arrow.
    eps_form : int
        +-1 for an arrow fixed by the transpose, 0 otherwise.
    kappa, mu : int
        The signs of the involution and the conjugation on the arrow's bimodule.
    square : bool
        Whether the blocks are square (``d_source == d_target``), so that a single element ``X a`` can
        be invertible.
    """

    def __init__(self, name: str, source: str, target: str, top: str, bar: str, eps_form: int = 0) -> None:
        self.name, self.source, self.target = name, source, target
        self.top, self.bar, self.eps_form = top, bar, eps_form
        self.kappa = self.mu = 0
        self.square = True

    def __repr__(self) -> str:
        return f"Arrow({self.name!r}: {self.source} -> {self.target})"


class Component(NamedTuple):
    """A connected component of a quiver, as reported by :meth:`DSQ.localisation_structure`."""

    vertices: list[str]
    arrows: list[str]
    cycle_rank: int


class DSQ:
    """A doubly-signed quiver.

    Usually constructed with :meth:`DSQ.build`; the constructor takes complete vertex and arrow
    dictionaries, validates them (raising :class:`~dsqlib.errors.QuiverError`), sets up the block
    model and solves for the arrow signs.

    Parameters
    ----------
    vertices : mapping of str to Vertex
    arrows : mapping of str to Arrow
    """

    def __init__(self, vertices: Mapping[str, Vertex], arrows: Mapping[str, Arrow]) -> None:
        self.vertices: dict[str, Vertex] = dict(vertices)
        self.arrows: dict[str, Arrow] = dict(arrows)
        self.validate()
        self._vertex_blocks()
        self._solve_signs()
        for a in self.arrows.values():
            a.square = self.vertices[a.source].block_size == self.vertices[a.target].block_size
        self._orbit_cache: dict[str, tuple[str, ...]] | None = None

    # ------------------------------------------------------------------ construction
    @classmethod
    def build(cls, orbits: Iterable[tuple[str, str]], arrows: Iterable[Sequence[Any]] = ()) -> DSQ:
        """Build a DSQ from real vertex kinds and arrows.

        Parameters
        ----------
        orbits : iterable of (base_name, kind)
            ``kind`` is a key of :data:`KINDS`. A kind with one complexified vertex keeps the base
            name ``b``; conjugate vertices are named ``b~``, transposed ones ``b'``, and the fourth
            vertex of ``CC`` is ``b~'``.
        arrows : iterable of (name, source, target) or (name, source, target, options)
            ``source`` and ``target`` are complexified vertex names. ``options`` is a dict with

            * ``sign`` (+-1): required for an arrow fixed by the transpose, i.e. ``i -> i^T``;
            * ``top_is_bar`` (bool): declares ``a^T = a~``, as for the loop of the Hermitian
              eigenproblem over C. This holds for every arrow ``i -> i^H`` (rule (5) of
              Definition 5.7), so the option only checks that the arrow has this form.

        Returns
        -------
        DSQ

        Examples
        --------
        >>> Q = DSQ.build([("u", "Cu"), ("v", "Cu")], [("a", "u", "v")])   # the complex SVD quiver
        >>> sorted(Q.arrows)
        ['a', "a'", 'a~', "a~'"]
        """
        V: dict[str, Vertex] = {}
        for b, kind in orbits:
            if kind not in KINDS:
                raise QuiverError(f"unknown vertex kind {kind!r}; the kinds are {list(KINDS)}")
            if any(v.orbit == b for v in V.values()):
                raise QuiverError(f"orbit name {b!r} used twice")
            pattern, ed, ef = KINDS[kind][:3]
            if pattern == "fixed":
                V[b] = Vertex(b, b, b, ed, ef, b, kind)
            elif pattern == "cfix":
                V[b] = Vertex(b, b, b + "~", 0, ef, b, kind)
                V[b + "~"] = Vertex(b + "~", b + "~", b, 0, ef, b, kind)
            elif pattern == "cu":
                V[b] = Vertex(b, b + "~", b + "~", 0, 0, b, kind)
                V[b + "~"] = Vertex(b + "~", b, b, 0, 0, b, kind)
            elif pattern == "exch":
                V[b] = Vertex(b, b + "'", b, ed, 0, b, kind)
                V[b + "'"] = Vertex(b + "'", b, b + "'", ed, 0, b, kind)
            else:  # 'cc'
                names = {"1": b, "B": b + "~", "T": b + "'", "TB": b + "~'"}
                for g, nm in names.items():
                    V[nm] = Vertex(nm, names[_group_mul(g, "T")], names[_group_mul(g, "B")], 0, 0, b, kind)
        A: dict[str, Arrow] = {}
        for spec in arrows:
            name, s, t = spec[:3]
            opts = spec[3] if len(spec) > 3 else {}
            unknown = set(opts) - {"sign", "top_is_bar"}
            if unknown:
                raise QuiverError(
                    f"arrow {name}: unknown options {sorted(unknown)} (the options are 'sign', 'top_is_bar')"
                )
            for v in (s, t):
                if v not in V:
                    raise QuiverError(f"arrow {name}: unknown vertex {v!r}; the vertices are {list(V)}")
            for a in cls._arrow_orbit(V, name, s, t, opts):
                if a.name in A:
                    raise QuiverError(f"arrow name {a.name!r} used twice")
                A[a.name] = a
        return cls(V, A)

    @staticmethod
    def _act(V: Mapping[str, Vertex], g: str, v: str) -> str:
        if g in ("T", "TB"):
            v = V[v].top
        if g in ("B", "TB"):
            v = V[v].bar
        return v

    @classmethod
    def _arrow_orbit(cls, V: Mapping[str, Vertex], name: str, s: str, t: str, opts: Mapping[str, Any]) -> list[Arrow]:
        """Close one arrow under the Klein group, respecting the normalisations of Definitions 5.4
        and 5.7."""

        def act(g: str, v: str) -> str:
            return cls._act(V, g, v)

        ends = {
            "1": (s, t),
            "T": (act("T", t), act("T", s)),
            "B": (act("B", s), act("B", t)),
            "TB": (act("TB", t), act("TB", s)),
        }
        stab = {"1"}
        sign = opts.get("sign")
        if t == V[s].top:  # arrows i -> i^T are fixed by T
            stab.add("T")
            if sign is None:
                raise QuiverError(f"arrow {name}: {s} -> {t} = {s}^T is fixed by the transpose; give its sign")
        elif sign is not None:
            raise QuiverError(f"arrow {name}: a sign is only allowed for arrows i -> i^T")
        if V[s].bar == s and V[t].bar == t and V[s].eps_div == V[t].eps_div:
            stab.add("B")  # real arrows between conjugation-fixed vertices
        if ends["TB"] == (s, t):
            stab.add("TB")  # arrows i -> i^H are fixed by the conjugate transpose (rule (5))
        elif opts.get("top_is_bar"):
            raise QuiverError(f"arrow {name}: top_is_bar needs {t} = ({s})^T~")
        for g, h in itertools.product(list(stab), repeat=2):
            stab.add(_group_mul(g, h))
        for g in stab:
            if ends[g] != (s, t):
                raise QuiverError(f"arrow {name}: inconsistent with the vertex involutions")
        cosets: dict[frozenset[str], str] = {}
        for g in _GROUP:
            key = frozenset(_group_mul(g, h) for h in stab)
            if key not in cosets:
                cosets[key] = name + _SUFFIX[g]

        def name_of(g: str) -> str:
            return cosets[frozenset(_group_mul(g, h) for h in stab)]

        out = []
        for key, aname in cosets.items():
            g = next(iter(key))
            ss, tt = ends[g]
            eps = int(sign) if (sign is not None and "T" in stab) else 0
            out.append(Arrow(aname, ss, tt, name_of(_group_mul("T", g)), name_of(_group_mul("B", g)), eps))
        return out

    # ------------------------------------------------------------------ checks
    def validate(self) -> None:
        """Check the axioms of Definition 5.7; raise :class:`~dsqlib.errors.QuiverError` listing every
        violation."""
        V, A = self.vertices, self.arrows
        err = []
        for v in V.values():
            if v.top not in V or v.bar not in V:
                err.append(f"{v.name}: unknown partner")
            elif V[v.top].top != v.name or V[v.bar].bar != v.name:
                err.append(f"{v.name}: maps are not involutions")
            elif V[v.top].bar != V[v.bar].top:
                err.append(f"{v.name}: transpose and conjugation do not commute")
            if (v.bar == v.name) != (v.eps_div in (1, -1)):
                err.append(f"{v.name}: eps_div must be +-1 exactly on conjugation-fixed vertices")
            if (v.top == v.name) != (v.eps_form in (1, -1)):
                err.append(f"{v.name}: eps_form must be +-1 exactly on transpose-fixed vertices")
        for v in V.values():
            if v.top in V and V[v.top].eps_div != v.eps_div:
                err.append(f"{v.name}: eps_div(v^T) != eps_div(v)")
            if v.bar in V and V[v.bar].eps_form != v.eps_form:
                err.append(f"{v.name}: eps_form(v~) != eps_form(v)")
        for a in A.values():
            if a.source not in V or a.target not in V:
                err.append(f"{a.name}: unknown endpoint")
                continue
            if a.top not in A or a.bar not in A:
                err.append(f"{a.name}: unknown partner")
                continue
            at, ab = A[a.top], A[a.bar]
            if at.top != a.name or ab.bar != a.name:
                err.append(f"{a.name}: maps are not involutions")
            if (at.source, at.target) != (V[a.target].top, V[a.source].top):
                err.append(f"{a.name}: transpose must reverse the arrow")
            if (ab.source, ab.target) != (V[a.source].bar, V[a.target].bar):
                err.append(f"{a.name}: conjugation must preserve the arrow")
            if A[at.bar].name != A[ab.top].name:
                err.append(f"{a.name}: transpose and conjugation do not commute")
            if (a.top == a.name) != (a.eps_form in (1, -1)):
                err.append(f"{a.name}: eps_form must be +-1 exactly on transpose-fixed arrows")
            if ab.eps_form != a.eps_form:
                err.append(f"{a.name}: eps_form(a~) != eps_form(a)")
            if a.target == V[a.source].top and a.top != a.name:
                err.append(f"{a.name}: arrows i -> i^T must be fixed by the transpose")
            conj_fixed_ends = V[a.source].bar == a.source and V[a.target].bar == a.target
            if conj_fixed_ends and (a.bar == a.name) != (V[a.source].eps_div == V[a.target].eps_div):
                err.append(f"{a.name}: between conjugation-fixed vertices, a~ = a exactly when eps_div agree")
            if a.target == V[V[a.source].bar].top and A[ab.top].name != a.name:
                err.append(f"{a.name}: arrows i -> i^H must be fixed by the conjugate transpose")
        if err:
            raise QuiverError("not a doubly-signed quiver:\n  " + "\n  ".join(err))

    # ------------------------------------------------------------------ vertex blocks and signs
    def _vertex_blocks(self) -> None:
        for v in self.vertices.values():
            transpose_fixed, conj_fixed = v.top == v.name, v.bar == v.name
            v.block_size = 2 if (v.eps_div == -1 or v.eps_form == -1) else 1
            eye = np.eye(v.block_size, dtype=complex)
            v.G = (_J2.copy() if v.eps_form == -1 else eye.copy()) if transpose_fixed else eye.copy()
            v.C = (_J2.copy() if v.eps_div == -1 else eye.copy()) if conj_fixed else eye.copy()
            v.G_inv, v.C_inv = np.linalg.inv(v.G), np.linalg.inv(v.C)
            v._g = v.eps_form if transpose_fixed else 1
            v._h = v.eps_div if conj_fixed else 1
            if not v.kind:
                v.kind = self._detect_kind(v)
        for v in self.vertices.values():  # the vertex data must make * and c commute
            vt, vb = self.vertices[v.top], self.vertices[v.bar]
            if not np.allclose(vt.C @ np.conj(v.G) @ v.C.T, vb.G):
                raise QuiverError(f"vertex {v.name}: C G C^T != G~ (internal error)")

    @staticmethod
    def _detect_kind(v: Vertex) -> str:
        transpose_fixed, conj_fixed = v.top == v.name, v.bar == v.name
        if transpose_fixed and conj_fixed:
            pattern = "fixed"
        elif transpose_fixed:
            pattern = "cfix"
        elif conj_fixed:
            pattern = "exch"
        elif v.top == v.bar:
            pattern = "cu"
        else:
            pattern = "cc"
        for k, kind in KINDS.items():
            p = kind.pattern
            if (
                p == pattern
                and (kind.eps_div == v.eps_div or p in ("cfix", "cu", "cc"))
                and (kind.eps_form == v.eps_form or p in ("cu", "exch", "cc"))
            ):
                return k
        raise QuiverError(f"vertex {v.name}: no vertex kind matches its data")

    def _solve_signs(self) -> None:
        """kappa, mu in {+-1} on each arrow orbit with

        (I)   kappa_a kappa_{a^T} = g_s g_t          (** = 1)
        (II)  mu_a mu_{a~}        = h_s h_t          (cc = 1)
        (III) kappa_a mu_{a^T}    = mu_a kappa_{a~}  (c* = *c)

        and the normalisations kappa_a = eps_form(a) when a^T = a, mu_a = +1 when a~ = a."""
        A = self.arrows
        done: set[str] = set()
        for a0 in A.values():
            if a0.name in done:
                continue
            orb = sorted({a0.name, a0.top, a0.bar, A[a0.top].bar})
            done |= set(orb)
            for bits in itertools.product((1, -1), repeat=2 * len(orb)):
                kap = dict(zip(orb, bits[: len(orb)], strict=True))
                mu = dict(zip(orb, bits[len(orb) :], strict=True))
                if all(self._signs_ok(A[n], kap, mu) for n in orb):
                    for n in orb:
                        A[n].kappa, A[n].mu = kap[n], mu[n]
                    break
            else:
                raise QuiverError(f"no consistent signs for the arrow orbit {orb}")

    def _signs_ok(self, a: Arrow, kap: Mapping[str, int], mu: Mapping[str, int]) -> bool:
        vs, vt = self.vertices[a.source], self.vertices[a.target]
        n = a.name
        if a.top == n and kap[n] != a.eps_form:
            return False
        if a.bar == n and mu[n] != 1:
            return False
        return (
            kap[n] * kap[a.top] == vs._g * vt._g
            and mu[n] * mu[a.bar] == vs._h * vt._h
            and kap[n] * mu[a.top] == mu[n] * kap[a.bar]
        )

    # ------------------------------------------------------------------ structure
    def star_bar_orbits(self) -> list[list[str]]:
        """The orbits of <T, B> on the vertices (they correspond to the *-simple factors of
        ``B/rad B``), each in the order of :attr:`vertices`."""
        order = list(self.vertices).index
        seen: set[str] = set()
        out = []
        for v in self.vertices:
            if v in seen:
                continue
            orb: set[str] = set()
            todo = [v]
            while todo:
                x = todo.pop()
                if x in orb:
                    continue
                orb.add(x)
                todo += [self.vertices[x].top, self.vertices[x].bar]
            seen |= orb
            out.append(sorted(orb, key=order))
        return out

    def orbit_of(self, v: str) -> tuple[str, ...]:
        """The STAR-BAR orbit of the vertex ``v``."""
        if self._orbit_cache is None:
            self._orbit_cache = {u: tuple(o) for o in self.star_bar_orbits() for u in o}
        return self._orbit_cache[v]

    def orbit_names(self) -> dict[str, tuple[str, ...]]:
        """The STAR-BAR orbits by name: the base name given to :meth:`build` (or the first vertex)."""
        out = {}
        for o in self.star_bar_orbits():
            name = self.vertices[o[0]].orbit or o[0]
            out[name] = tuple(o)
        return out

    def localisation_structure(self) -> list[Component]:
        """The connected components, with their cycle ranks. Inverting every arrow of a component with
        n vertices and cycle rank r turns its basic algebra into ``M_n(C[F_r])`` (``F_r`` free of
        rank r), before completion."""
        parent = {v: v for v in self.vertices}

        def find(x: str) -> str:
            while parent[x] != x:
                parent[x] = parent[parent[x]]
                x = parent[x]
            return x

        for a in self.arrows.values():
            parent[find(a.source)] = find(a.target)
        comps: dict[str, list[str]] = {}
        for v in self.vertices:
            comps.setdefault(find(v), []).append(v)
        out = []
        for vs in comps.values():
            names = [a.name for a in self.arrows.values() if a.source in vs]
            out.append(Component(vs, names, len(names) - len(vs) + 1))
        return out

    def describe(self) -> str:
        """A multi-line summary: orbits with their kinds, arrows with their partners, components."""
        lines = ["STAR-BAR orbits (Table 5.1 kinds):"]
        for o in self.star_bar_orbits():
            v = self.vertices[o[0]]
            k = KINDS[v.kind]
            lines.append(
                f"  {{{', '.join(o)}}}: {v.kind:3s}  {k.algebra}, unitary group {k.unitary_group}, "
                f"block size {v.block_size}"
            )
        lines.append("arrows:")
        for a in self.arrows.values():
            extra = f"  eps_form {a.eps_form:+d}" if a.eps_form else ""
            lines.append(
                f"  {a.name}: {a.source} -> {a.target}   T: {a.top}   bar: {a.bar}{extra}   "
                f"blocks {self.vertices[a.target].block_size}x{self.vertices[a.source].block_size}"
            )
        lines.append("localisation (every arrow inverted), per connected component:")
        for c in self.localisation_structure():
            lines.append(f"  vertices {c.vertices}, cycle rank {c.cycle_rank}")
        return "\n".join(lines)

    def __repr__(self) -> str:
        kinds = ", ".join(f"{name}:{self.vertices[vs[0]].kind}" for name, vs in self.orbit_names().items())
        return f"DSQ(orbits=[{kinds}], arrows={list(self.arrows)})"
\end{lstlisting}
\needspace{6\baselineskip}
\subsection{\texorpdfstring{\texttt{src/dsqlib/algebra.py}}{src/dsqlib/algebra.py}}
\begin{lstlisting}[style=dsqcode]
"""The completed localisation of the path *-algebra of a doubly-signed quiver.

Localisation
    Every arrow ``a: i -> j`` receives a formal *local* inverse ``a^-1: j -> i`` with ``a^-1 a = e_i``
    and ``a a^-1 = e_j`` (Gabriel-Zisman / Cohn localisation: only the arrow is inverted, between its
    two vertices). The result is spanned by reduced words in the arrows and their inverses, and the
    product is concatenation with cancellation; this works for every quiver, wild ones included. The
    involution extends by ``(a^-1)* = (a*)^-1`` and the conjugation by ``c(a^-1) = c(a)^-1``. When
    ``d_i != d_j`` no single element ``X a`` is invertible, but sums over several arrows can be
    (Morita equivalence), which the block model represents exactly.

Completion
    The path degree (+1 per arrow, -1 per inverse letter) is a grading. An element is a series
    ``sum_d x_d`` of homogeneous parts, bounded below and known modulo degree ``>= prec``
    (``prec = inf`` for an exact element), like a Laurent series known to some order. An element
    whose lowest-degree part is a single word is a unit, inverted by a geometric series.

Coefficients
    Everything is computed in the complexification with complex floating-point blocks: a word from
    ``s`` to ``t`` carries a ``d_t x d_s`` block. A block in which a sum cancels to within ``tol``
    (relative, default 1e-12) of the size of its summands is set to zero, so exact zeros of the
    algebra come out as zeros whatever the scale of the data. The real algebra is the fixed set of
    the conjugation; real inputs give real outputs.

Words
    A word is a reduced tuple of letters ``(arrow_name, +1 or -1)``, written as for composition of
    functions (the last letter is applied first). Functions that accept words from users take a list
    of arrow names, with ``'^-1'`` appended for an inverse letter, e.g. ``['a', "a'"]`` or ``['b^-1']``.
"""

from __future__ import annotations

import numbers
from collections import defaultdict
from collections.abc import Iterable, Mapping, Sequence
from functools import cache
from typing import TYPE_CHECKING, Any

import numpy as np

from .errors import NotInvertibleError, PrecisionError

if TYPE_CHECKING:
    from .quiver import DSQ, Arrow, Vertex

__all__ = ["INF", "Algebra", "Element", "Key", "Scalar", "Word", "quaternion"]

INF = float("inf")
"""The precision of an exact element."""

Word = tuple[tuple[str, int], ...]
"""A reduced word: a tuple of letters ``(arrow_name, +1 | -1)``; the last letter is applied first."""

Key = tuple[str, str, Word]
"""The key of a term: ``(target, source, word)``."""

_WITHIN = ("support", "orbits", "all")

Scalar = int | float | complex | np.number
"""A number that can scale an element."""
_SCALARS = (int, float, complex, np.number)


class Element:
    """An element of an :class:`Algebra`: terms ``{(target, source, word): block}``, known modulo
    degree ``>= prec``.

    Elements support ``+``, ``-``, ``*`` (the algebra product, or scaling by a number), ``/`` by a
    number and ``**`` with a non-negative integer exponent; a number added to an element stands for
    that multiple of the unit. ``x.H`` is the involution ``x*``.

    Elements are normally made by the constructors of :class:`Algebra` (``vertex``, ``word``,
    ``arrow``, ``one``, ``element``) and by arithmetic, not directly.
    """

    __slots__ = ("algebra", "prec", "terms")
    __array_ufunc__ = None  # numpy operands defer to our reflected operators

    def __init__(self, algebra: Algebra, terms: dict[Key, np.ndarray] | None = None, prec: float = INF) -> None:
        self.algebra = algebra
        self.terms: dict[Key, np.ndarray] = terms if terms is not None else {}
        self.prec = prec

    # ------------------------------------------------------------------ arithmetic
    def _coerce(self, other: Any) -> Element:
        if isinstance(other, Element):
            if other.algebra is not self.algebra:
                raise ValueError("elements of different algebras")
            return other
        if isinstance(other, _SCALARS):
            return self.algebra.scale(other, self.algebra.one())
        return NotImplemented

    def __add__(self, other: Any) -> Element:
        o = self._coerce(other)
        return NotImplemented if o is NotImplemented else self.algebra.add(self, o)

    __radd__ = __add__

    def __sub__(self, other: Any) -> Element:
        o = self._coerce(other)
        return NotImplemented if o is NotImplemented else self.algebra.sub(self, o)

    def __rsub__(self, other: Any) -> Element:
        o = self._coerce(other)
        return NotImplemented if o is NotImplemented else self.algebra.sub(o, self)

    def __neg__(self) -> Element:
        return self.algebra.scale(-1.0, self)

    def __pos__(self) -> Element:
        return self

    def __mul__(self, other: Any) -> Element:
        if isinstance(other, Element):
            return self.algebra.mul(self, self._coerce(other))
        if isinstance(other, _SCALARS):
            return self.algebra.scale(other, self)
        return NotImplemented

    def __rmul__(self, other: Any) -> Element:
        if isinstance(other, _SCALARS):
            return self.algebra.scale(other, self)
        return NotImplemented

    def __truediv__(self, other: Any) -> Element:
        if isinstance(other, _SCALARS):
            return self.algebra.scale(1.0 / other, self)
        return NotImplemented

    def __pow__(self, n: int) -> Element:
        if not isinstance(n, numbers.Integral) or n < 0:
            raise ValueError("only non-negative integer powers are defined; use inverse() for inverses")
        out = self.algebra.one()
        for _ in range(int(n)):
            out = self.algebra.mul(out, self)
        return out

    # ------------------------------------------------------------------ structure maps
    @property
    def H(self) -> Element:
        """The involution ``x*``."""
        return self.algebra.star(self)

    def conj(self) -> Element:
        """The conjugation ``c(x)`` (a conjugate-linear automorphism)."""
        return self.algebra.conj(self)

    def real(self) -> Element:
        """The real part ``(x + c(x)) / 2``."""
        return self.algebra.real(self)

    def inverse(self, within: str = "support") -> Element:
        """The inverse; see :meth:`Algebra.inverse`."""
        return self.algebra.inverse(self, within)

    def sqrt(self, rule: str = "graded", choice: str = "positive", split: str = "left", rng: Any = None) -> Element:
        """A *-square root ``y`` with ``y* y = x``; see :func:`dsqlib.sqrt.star_sqrt`."""
        from .sqrt import star_sqrt

        return star_sqrt(self, rule, choice, split, rng)

    # ------------------------------------------------------------------ degrees and precision
    @property
    def valuation(self) -> float:
        """The lowest degree present (or ``prec``, whichever is smaller); ``inf`` for an exact zero."""
        return self.algebra.valuation(self)

    @property
    def degrees(self) -> list[int]:
        """The degrees of the words present, in increasing order."""
        return self.algebra.degrees(self)

    @property
    def is_zero(self) -> bool:
        """True if no term is present (the element may still be inexact: ``O(deg p)``)."""
        return not self.terms

    @property
    def is_exact(self) -> bool:
        """True if the element is known exactly (``prec = inf``)."""
        return self.prec == INF

    def part(self, d: int) -> Element:
        """The homogeneous part of degree ``d`` (``O(deg prec)`` if the element is not known at
        degree ``d``)."""
        return self.algebra.part(self, d)

    def truncate(self, prec: float) -> Element:
        """The element known modulo degree ``>= prec`` only."""
        return self.algebra.truncate(self, prec)

    # ------------------------------------------------------------------ sizes and tests
    def norm(self) -> float:
        """The sum of the absolute values of all coefficients."""
        return self.algebra.norm(self)

    def dist(self, other: Element) -> float:
        """The norm of ``self - other`` without any rounding to zero (see :meth:`Algebra.dist`)."""
        return self.algebra.dist(self, other)

    def is_real(self, tol: float = 1e-9) -> bool:
        """Whether ``c(x) = x`` to relative tolerance ``tol``."""
        return self.algebra.is_real(self, tol)

    def is_selfadjoint(self, tol: float = 1e-9) -> bool:
        """Whether ``x* = x`` to relative tolerance ``tol``."""
        return self.algebra.is_selfadjoint(self, tol)

    @property
    def support(self) -> list[tuple[str, str]]:
        """The corners ``(target, source)`` in which the element has terms."""
        return self.algebra.support(self)

    def corner(self, target: str, source: str) -> Element:
        """The part of the element in the corner ``e_target A e_source``."""
        return self.algebra.corner(self, target, source)

    def coefficient(
        self, letters: Sequence[str] = (), target: str | None = None, source: str | None = None
    ) -> np.ndarray:
        """The block of the given word.

        Parameters
        ----------
        letters : sequence of str
            The word as arrow names (``'^-1'`` appended for inverses); empty for a vertex.
        target, source : str, optional
            The end vertices: required for the empty word (and then equal); for a non-empty word they
            default to its ends, and must agree with them if given.

        Returns
        -------
        numpy.ndarray
            The ``d_target x d_source`` block (zeros if the word is absent).

        Raises
        ------
        ValueError
            If the word is invalid, does not have the given ends, or has a degree at which the element
            is not known (``>= prec``).
        """
        return self.algebra.coefficient(self, letters, target, source)

    # ------------------------------------------------------------------ display
    def format(self, digits: int = 4) -> str:
        """The element as a sum of ``coefficient·word`` terms."""
        return self.algebra.format(self, digits)

    def __repr__(self) -> str:
        return self.algebra.format(self)


class Algebra:
    """The completed localisation of the path *-algebra of a doubly-signed quiver.

    Parameters
    ----------
    quiver : DSQ
    tol : float, optional
        Relative cancellation tolerance: a block in which a sum cancels to within ``tol`` of the size
        of its summands is set to zero.
    relprec : int, optional
        The number of degrees (relative to the lowest one) to which series inverses and square roots
        of exact elements are computed.

    Notes
    -----
    Besides the element constructors, the class offers every operation in functional form
    (``add``, ``mul``, ``star``, ...); the operators of :class:`Element` call these.

    Examples
    --------
    >>> from dsqlib import Algebra, catalogue
    >>> A = Algebra(catalogue.svd("R"), relprec=6)
    >>> a = A.arrow("a")
    >>> A.word(["a^-1"]) * a  # a^-1 is a local inverse
    1·e_u
    >>> (1 + 0.5 * a).inverse()  # a a = 0, so the series stops
    1·e_u + 1·e_v + -0.5·a
    >>> (1 + 0.5 * (a + A.arrow("a'"))).inverse().prec  # an infinite series, known to degree relprec
    6
    """

    def __init__(self, quiver: DSQ, tol: float = 1e-12, relprec: int = 16) -> None:
        self.quiver = quiver
        self.tol = tol
        self.relprec = relprec
        self.vertices: dict[str, Vertex] = quiver.vertices
        self.arrows: dict[str, Arrow] = quiver.arrows

    def __repr__(self) -> str:
        return f"Algebra({self.quiver!r}, tol={self.tol:g}, relprec={self.relprec})"

    def block_size(self, v: str) -> int:
        """``d_v``, the block size of the vertex ``v``."""
        return self.vertices[v].block_size

    # ------------------------------------------------------------------ words
    @staticmethod
    @cache
    def degree(w: Word) -> int:
        """The path degree of a word (+1 per arrow, -1 per inverse letter)."""
        return sum(e for _, e in w)

    @staticmethod
    def concat(w1: Word, w2: Word) -> Word:
        """The reduced product ``w1 w2`` (``w2`` applied first)."""
        i, j = len(w1), 0
        while i > 0 and j < len(w2) and w1[i - 1][0] == w2[j][0] and w1[i - 1][1] == -w2[j][1]:
            i -= 1
            j += 1
        return w1[:i] + w2[j:]

    @staticmethod
    def word_inverse(w: Word) -> Word:
        """The inverse word ``w^-1``."""
        return tuple((n, -e) for n, e in reversed(w))

    def word_transpose(self, w: Word) -> Word:
        """The transposed word ``w^T``."""
        return tuple((self.arrows[n].top, e) for n, e in reversed(w))

    def word_conjugate(self, w: Word) -> Word:
        """The conjugate word ``w~``."""
        return tuple((self.arrows[n].bar, e) for n, e in w)

    def word_kappa(self, w: Word) -> int:
        """The product of the involution signs ``kappa`` of the letters of ``w``."""
        k = 1
        for n, _ in w:
            k *= self.arrows[n].kappa
        return k

    def word_mu(self, w: Word) -> int:
        """The product of the conjugation signs ``mu`` of the letters of ``w``."""
        m = 1
        for n, _ in w:
            m *= self.arrows[n].mu
        return m

    def word_ends(self, w: Word) -> tuple[str, str]:
        """``(target, source)`` of a non-empty word."""
        n, e = w[-1]
        a = self.arrows[n]
        src = a.source if e > 0 else a.target
        n, e = w[0]
        a = self.arrows[n]
        tgt = a.target if e > 0 else a.source
        return tgt, src

    @staticmethod
    def word_str(w: Word, vertex: str | None = None) -> str:
        """A readable form of a word (``e_v`` for the empty word at ``v``)."""
        if not w:
            return f"e_{vertex}"
        return "·".join(n if e > 0 else n + "⁻¹" for n, e in w)

    def parse_word(self, letters: Iterable[str], check: bool = True) -> Word:
        """The reduced word of a list of arrow names (``'^-1'`` appended for an inverse letter).

        With ``check`` (the default) the arrows must exist and consecutive letters must compose
        (the last letter is applied first); a ``ValueError`` says which do not.
        """
        w = tuple((L[:-3], -1) if L.endswith("^-1") else (L, 1) for L in letters)
        if check:
            for n, _ in w:
                if n not in self.arrows:
                    raise ValueError(f"unknown arrow {n!r}; the arrows are {list(self.arrows)}")

            def src(letter: tuple[str, int]) -> str:
                a = self.arrows[letter[0]]
                return a.source if letter[1] > 0 else a.target

            def tgt(letter: tuple[str, int]) -> str:
                a = self.arrows[letter[0]]
                return a.target if letter[1] > 0 else a.source

            for i in range(len(w) - 1):
                if tgt(w[i + 1]) != src(w[i]):
                    raise ValueError(
                        f"the letters {w[i + 1][0]} and {w[i][0]} do not compose "
                        f"({tgt(w[i + 1])} != {src(w[i])}); the last letter is applied first"
                    )
        red: Word = ()
        for L in w:
            red = self.concat(red, (L,))
        return red

    # ------------------------------------------------------------------ construction
    def element(self, terms: Mapping[Key, np.ndarray], prec: float = INF) -> Element:
        """An element from ``{(target, source, word): block}``, dropping words of degree ``>= prec``
        and zero blocks."""
        out = {}
        for k, X in terms.items():
            if self.degree(k[2]) >= prec:
                continue
            if np.any(X != 0):
                out[k] = X
        return Element(self, out, prec)

    def zero(self, prec: float = INF) -> Element:
        """Zero; with a finite ``prec``, the element ``O(deg prec)``."""
        return Element(self, {}, prec)

    def idempotent(self, v: str) -> Element:
        """The vertex idempotent ``e_v``."""
        if v not in self.vertices:
            raise ValueError(f"unknown vertex {v!r}; the vertices are {list(self.vertices)}")
        return Element(self, {(v, v, ()): np.eye(self.block_size(v), dtype=complex)})

    def one(self) -> Element:
        """The unit ``1 = sum of e_v``."""
        return Element(self, {(v, v, ()): np.eye(x.block_size, dtype=complex) for v, x in self.vertices.items()})

    def _block(self, X: Any, rows: int, cols: int) -> np.ndarray:
        X = np.atleast_2d(np.asarray(X, dtype=complex))
        if X.shape == (1, 1) and (rows, cols) != (1, 1):
            if rows != cols:
                raise ValueError("a scalar block needs equal block sizes")
            X = X[0, 0] * np.eye(rows, dtype=complex)
        if X.shape != (rows, cols):
            raise ValueError(f"block must be {rows}x{cols}")
        return X

    def vertex(self, v: str, X: Any = 1.0) -> Element:
        """``X e_v`` with ``X`` a ``d_v x d_v`` block (a number means a multiple of the identity)."""
        if v not in self.vertices:
            raise ValueError(f"unknown vertex {v!r}; the vertices are {list(self.vertices)}")
        d = self.block_size(v)
        return self.element({(v, v, ()): self._block(X, d, d)})

    def word(self, letters: Iterable[str], X: Any = 1.0) -> Element:
        """``X w`` for the word ``w`` given as a list of arrow names (``'^-1'`` for a local inverse;
        the last letter is applied first). ``X`` is a ``d_target x d_source`` block or a number."""
        red = self.parse_word(letters)
        if not red:
            raise ValueError("the word reduces to a vertex; use vertex()")
        t, s = self.word_ends(red)
        return self.element({(t, s, red): self._block(X, self.block_size(t), self.block_size(s))})

    def arrow(self, name: str, X: Any = 1.0) -> Element:
        """``X a`` for the arrow ``a`` (or ``X a^-1`` for ``name = 'a^-1'``)."""
        return self.word([name], X)

    def real_element(self, x: Element) -> Element:
        """The real element whose coefficients at the words of ``x`` are those of ``x``.

        A term ``X w`` with ``c(w) != w`` contributes ``X w + c(X w)``; a term with ``c(w) = w`` must
        have a block of real type and is kept (projected onto the real blocks). So for a complex
        vertex ``z e_v`` gives ``z e_v + conj(z) e_v~``, and for the complex SVD quiver ``z a`` gives
        ``z a + conj(z) a~``.
        """
        moved = Element(self, {k: X for k, X in x.terms.items() if self._conj_key(k) != k}, x.prec)
        fixed = Element(self, {k: X for k, X in x.terms.items() if self._conj_key(k) == k}, x.prec)
        cm = self.conj(moved)
        if set(cm.terms) & set(moved.terms):
            raise ValueError("real_element: x contains a term and its conjugate; give only one of them")
        return self.add(self.add(moved, cm), self.real(fixed))

    def _conj_key(self, k: Key) -> Key:
        t, s, w = k
        return (self.vertices[t].bar, self.vertices[s].bar, self.word_conjugate(w))

    # ------------------------------------------------------------------ valuation and precision
    def valuation(self, x: Element) -> float:
        """The lowest degree present (or ``prec``, whichever is smaller); ``inf`` for an exact zero."""
        d = min((self.degree(k[2]) for k in x.terms), default=INF)
        return min(d, x.prec)

    def truncate(self, x: Element, prec: float) -> Element:
        """``x`` known modulo degree ``>= prec`` only."""
        return self.element(dict(x.terms), min(prec, x.prec))

    def part(self, x: Element, d: int) -> Element:
        """The homogeneous part of degree ``d`` (``O(deg prec)`` if ``x`` is not known at degree ``d``)."""
        if d >= x.prec:
            return Element(self, {}, x.prec)
        return Element(self, {k: X for k, X in x.terms.items() if self.degree(k[2]) == d}, INF)

    def degrees(self, x: Element) -> list[int]:
        """The degrees of the words present, in increasing order."""
        return sorted({self.degree(k[2]) for k in x.terms})

    # ------------------------------------------------------------------ ring operations
    def add(self, x: Element, y: Element) -> Element:
        """``x + y``. A block in which the two summands cancel to within ``tol`` of their size is set to
        zero (relative cancellation, so the result does not depend on the scale of the data)."""
        prec = min(x.prec, y.prec)
        deg, tol = self.degree, self.tol
        t = {k: X for k, X in x.terms.items() if deg(k[2]) < prec}
        for k, Y in y.terms.items():
            if deg(k[2]) >= prec:
                continue
            X = t.get(k)
            if X is None:
                t[k] = Y
            else:
                S = X + Y
                if np.abs(S).max() <= tol * (np.abs(X).max() + np.abs(Y).max()):
                    del t[k]
                else:
                    t[k] = S
        return Element(self, t, prec)

    def sub(self, x: Element, y: Element) -> Element:
        """``x - y``."""
        return self.add(x, self.scale(-1.0, y))

    def scale(self, c: Scalar, x: Element) -> Element:
        """``c x`` for a number ``c``."""
        if c == 0:
            return Element(self, {}, x.prec)
        return Element(self, {k: c * X for k, X in x.terms.items()}, x.prec)

    def mul(self, x: Element, y: Element) -> Element:
        """The product ``x y``, with precision ``min(prec x + val y, prec y + val x)``."""
        prec = min(x.prec + self.valuation(y), y.prec + self.valuation(x))
        by_target = defaultdict(list)
        for (t2, s2, w2), B in y.terms.items():
            by_target[t2].append((s2, w2, self.degree(w2), B))
        out: dict[Key, np.ndarray] = {}
        mag: dict[Key, float] = {}
        for (t1, s1, w1), A in x.terms.items():
            d1, nA = self.degree(w1), np.abs(A).max() * A.shape[1]
            for s2, w2, d2, B in by_target.get(s1, ()):
                if d1 + d2 >= prec:
                    continue
                k = (t1, s2, self.concat(w1, w2))
                P = A @ B
                m = nA * np.abs(B).max()
                if k in out:
                    out[k] = out[k] + P
                    mag[k] += m
                else:
                    out[k] = P
                    mag[k] = m
        tol = self.tol  # drop the blocks where the contributions cancel
        return Element(self, {k: P for k, P in out.items() if np.abs(P).max() > tol * mag[k]}, prec)

    def star(self, x: Element) -> Element:
        """The involution: ``(X w)* = kappa(w) G_s X^T G_t^-1 w^T`` for ``w: s -> t``."""
        out = {}
        for (t, s, w), X in x.terms.items():
            Vs, Vt = self.vertices[s], self.vertices[t]
            out[(Vs.top, Vt.top, self.word_transpose(w))] = self.word_kappa(w) * (Vs.G @ X.T @ Vt.G_inv)
        return Element(self, out, x.prec)

    def conj(self, x: Element) -> Element:
        """The conjugation (conjugate-linear): ``c(X w) = mu(w) C_t conj(X) C_s^-1 w~``."""
        out = {}
        for (t, s, w), X in x.terms.items():
            Vs, Vt = self.vertices[s], self.vertices[t]
            out[(Vt.bar, Vs.bar, self.word_conjugate(w))] = self.word_mu(w) * (Vt.C @ np.conj(X) @ Vs.C_inv)
        return Element(self, out, x.prec)

    def real(self, x: Element) -> Element:
        """The projection ``(x + c(x)) / 2`` onto the real algebra."""
        return self.scale(0.5, self.add(x, self.conj(x)))

    def norm(self, x: Element) -> float:
        """The sum of the absolute values of all coefficients."""
        return float(sum(np.abs(X).sum() for X in x.terms.values()))

    def dist(self, x: Element, y: Element) -> float:
        """The sum of the absolute values of the coefficients of ``x - y`` over the words known in
        both, without any rounding to zero."""
        prec = min(x.prec, y.prec)
        t = {k: X for k, X in x.terms.items() if self.degree(k[2]) < prec}
        for k, Y in y.terms.items():
            if self.degree(k[2]) >= prec:
                continue
            t[k] = t[k] - Y if k in t else -Y
        return float(sum(np.abs(X).sum() for X in t.values()))

    def is_real(self, x: Element, tol: float = 1e-9) -> bool:
        """Whether ``c(x) = x`` to relative tolerance ``tol``."""
        return self.dist(self.conj(x), x) <= tol * self.norm(x)

    def is_selfadjoint(self, x: Element, tol: float = 1e-9) -> bool:
        """Whether ``x* = x`` to relative tolerance ``tol``."""
        return self.dist(self.star(x), x) <= tol * self.norm(x)

    # ------------------------------------------------------------------ corners and coefficients
    def corner(self, x: Element, target: str, source: str) -> Element:
        """The part of ``x`` in ``e_target A e_source``."""
        return Element(self, {k: X for k, X in x.terms.items() if k[0] == target and k[1] == source}, x.prec)

    def support(self, x: Element) -> list[tuple[str, str]]:
        """The corners ``(target, source)`` in which ``x`` has terms."""
        return sorted({(k[0], k[1]) for k in x.terms})

    def corner_project(self, x: Element, rows: Iterable[str], cols: Iterable[str]) -> Element:
        """The part of ``x`` with target in ``rows`` and source in ``cols``."""
        rows, cols = set(rows), set(cols)
        return Element(self, {k: X for k, X in x.terms.items() if k[0] in rows and k[1] in cols}, x.prec)

    def coefficient(
        self, x: Element, letters: Sequence[str] = (), target: str | None = None, source: str | None = None
    ) -> np.ndarray:
        """The block of a word in ``x``; see :meth:`Element.coefficient`."""
        w = self.parse_word(letters)
        if w:
            t, s = self.word_ends(w)
            if (target or t, source or s) != (t, s):
                raise ValueError(
                    f"the word {self.word_str(w)} goes from {s} to {t}, not from {source or s} to {target or t}"
                )
            target, source = t, s
        elif target is None or source is None or target != source:
            raise ValueError("the empty word needs equal target and source (a vertex)")
        for v in (target, source):
            if v not in self.vertices:
                raise ValueError(f"unknown vertex {v!r}; the vertices are {list(self.vertices)}")
        if self.degree(w) >= x.prec:
            raise ValueError(f"degree {self.degree(w)} is not known: the element is known below degree {x.prec}")
        X = x.terms.get((target, source, w))
        if X is not None:
            return X
        return np.zeros((self.block_size(target), self.block_size(source)), dtype=complex)

    # ------------------------------------------------------------------ inverses
    def inverse(self, x: Element, within: str = "support") -> Element:
        """The inverse of ``x``.

        Parameters
        ----------
        x : Element
        within : {'support', 'orbits', 'all'}
            Between which idempotents the inverse is taken:

            * ``'support'``: the vertices ``x`` touches. If ``x`` maps ``e_C`` to ``e_R`` (``C``, ``R``
              the sets of source and target vertices of its terms), the ``y`` in ``e_C A e_R`` with
              ``y x = e_C`` and ``x y = e_R`` (for a single corner this is Definition 3.5). Zero maps to
              zero.
            * ``'orbits'``: the STAR-BAR orbit idempotents ``x`` touches (real and self-adjoint), so a
              zero divisor of an orbit corner (such as ``(c, 0)`` in the split-complex numbers) has no
              inverse. Zero maps to zero.
            * ``'all'``: the two-sided inverse in the whole algebra.

        Returns
        -------
        Element

        Raises
        ------
        NotInvertibleError
            If there is no such inverse, or ``x`` is only known to vanish below some degree.
        PrecisionError
            If the elimination runs out of precision (a subclass of ``NotInvertibleError``).

        Notes
        -----
        A bijection between single vertices whose corners have unit leading terms is inverted corner
        by corner, by geometric series. Otherwise Gauss-Jordan elimination runs on the scalar entries,
        with pivots of lowest degree first and then of largest coefficient. This is complete when
        between any two vertices there is at most one reduced word of each degree (true for every quiver
        of the catalogue); otherwise it can miss units whose elimination needs a pivot made of several
        words of the same degree.
        """
        if within not in _WITHIN:
            raise ValueError(f"within must be one of {_WITHIN}, not {within!r}")
        if within == "all":
            verts = list(self.vertices)
            return self._gauss_jordan(x, verts, verts)
        if not x.terms:
            if x.prec == INF:
                return self.zero()
            raise NotInvertibleError(f"O(deg {x.prec}) is not known well enough to invert")
        sup = self.support(x)
        if within == "support":
            return self.inverse_between(x, {t for t, _ in sup}, {s for _, s in sup})
        rows = [v for t, _ in sup for v in self.quiver.orbit_of(t)]
        cols = [v for _, s in sup for v in self.quiver.orbit_of(s)]
        return self.inverse_between(x, rows, cols)

    def inverse_between(self, x: Element, rows: Iterable[str], cols: Iterable[str]) -> Element:
        """The ``y`` in ``e_cols A e_rows`` with ``y x = e_cols`` and ``x y = e_rows``, for ``x`` in
        ``e_rows A e_cols`` (``rows``, ``cols``: vertex sets).

        Raises
        ------
        NotInvertibleError
            If there is none.
        """
        order = list(self.vertices).index
        rows, cols = sorted(set(rows), key=order), sorted(set(cols), key=order)
        sup = self.support(x)
        if not sup:
            raise NotInvertibleError("zero is not invertible")
        if {t for t, _ in sup} == set(rows) and {s for _, s in sup} == set(cols) and len(sup) == len(rows) == len(cols):
            leads = [self._leading_unit(self.corner(x, t, s)) for t, s in sup]  # a bijection between
            if all(lu is not None for lu in leads):  # vertices: invert corner by corner
                out = self.zero()
                for (t, s), lu in zip(sup, leads, strict=True):
                    assert lu is not None
                    out = self.add(out, self._series_inverse(self.corner(x, t, s), lu))
                return out
        return self._gauss_jordan(x, rows, cols)

    def corner_inverse(self, z: Element) -> Element:
        """The inverse between vertices (Definition 3.5) of ``z`` in a single corner ``e_t A e_s``: the
        element of ``e_s A e_t`` with ``z^-1 z = e_s`` and ``z z^-1 = e_t``."""
        sup = self.support(z)
        if not sup:
            raise NotInvertibleError("zero is not invertible")
        if len(sup) != 1:
            raise ValueError(f"corner_inverse: the element is not in a single corner (support {sup})")
        lu = self._leading_unit(z)
        if lu is not None:
            return self._series_inverse(z, lu)
        ((t, s),) = sup
        return self._gauss_jordan(z, [t], [s])

    def _leading_unit(self, z: Element) -> tuple[Key, np.ndarray] | None:
        """For ``z`` in one corner: (key, block) of its lowest-degree part if that part is a unit (a
        single word with an invertible square block), else None."""
        if not z.terms:
            return None
        d0 = min(self.degree(k[2]) for k in z.terms)
        lead = [(k, X) for k, X in z.terms.items() if self.degree(k[2]) == d0]
        if len(lead) != 1:
            return None
        key, X = lead[0]
        if X.shape[0] != X.shape[1] or np.linalg.cond(X) > 1e12:
            return None
        return key, X

    def _series_inverse(self, z: Element, lead: tuple[Key, np.ndarray]) -> Element:
        """The geometric series for the inverse of ``z = z0 + r``, ``z0 = X w`` the unit leading part."""
        (t, s, w), X = lead
        z0inv = Element(self, {(s, t, self.word_inverse(w)): np.linalg.inv(X)}, INF)
        d0 = self.degree(w)
        r = self.sub(z, Element(self, {(t, s, w): X}, INF))
        if not r.terms and r.prec == INF:
            return z0inv
        prec = -d0 + ((z.prec - d0) if z.prec < INF else self.relprec)
        u = self.scale(-1.0, self.mul(z0inv, r))
        total, term = z0inv, z0inv
        while True:
            term = self.truncate(self.mul(u, term), prec)
            if not term.terms:
                break
            total = self.add(total, term)
        return self.truncate(total, prec)

    # Gauss-Jordan elimination on scalar entries. A vertex v of block size d_v contributes the d_v
    # indices (v, 0), ..., (v, d_v - 1); entry ((t, i), (s, j)) of x is the scalar series
    # sum_w X_w[i, j] w over the words w: s -> t. A scalar series is stored as (dict word -> complex,
    # prec). An entry can be a pivot when its lowest-degree part is a single word (a unit); the pivot
    # of lowest degree is taken, and among those the largest coefficient.
    def _sval(self, a: tuple[dict, float]) -> float:
        return min(min((self.degree(w) for w in a[0]), default=INF), a[1])

    def _strunc(self, a: tuple[dict, float], prec: float) -> tuple[dict, float]:
        return ({w: c for w, c in a[0].items() if self.degree(w) < prec}, min(a[1], prec))

    def _sadd(self, a: tuple[dict, float], b: tuple[dict, float], sign: float = 1.0) -> tuple[dict, float]:
        prec = min(a[1], b[1])
        t = {w: c for w, c in a[0].items() if self.degree(w) < prec}
        for w, c in b[0].items():
            if self.degree(w) >= prec:
                continue
            if w in t:
                S = t[w] + sign * c
                if abs(S) <= self.tol * (abs(t[w]) + abs(c)):
                    del t[w]
                else:
                    t[w] = S
            else:
                t[w] = sign * c
        return (t, prec)

    def _smul(self, a: tuple[dict, float], b: tuple[dict, float]) -> tuple[dict, float]:
        prec = min(a[1] + self._sval(b), b[1] + self._sval(a))
        t: dict = {}
        mag: dict = {}
        for w1, c1 in a[0].items():
            d1 = self.degree(w1)
            for w2, c2 in b[0].items():
                if d1 + self.degree(w2) < prec:
                    w = self.concat(w1, w2)
                    t[w] = t.get(w, 0) + c1 * c2
                    mag[w] = mag.get(w, 0) + abs(c1 * c2)
        return ({w: c for w, c in t.items() if abs(c) > self.tol * mag[w]}, prec)

    def _sleading(self, a: tuple[dict, float]) -> tuple[Word, complex] | None:
        if not a[0]:
            return None
        d0 = min(self.degree(w) for w in a[0])
        lead = [(w, c) for w, c in a[0].items() if self.degree(w) == d0]
        return lead[0] if len(lead) == 1 else None

    def _sinv(self, a: tuple[dict, float], relprec: float) -> tuple[dict, float]:
        lead = self._sleading(a)
        assert lead is not None
        w0, c0 = lead
        d0 = self.degree(w0)
        inv0 = ({self.word_inverse(w0): 1.0 / c0}, INF)
        r = ({w: c for w, c in a[0].items() if w != w0}, a[1])
        if not r[0] and r[1] == INF:
            return inv0
        prec = -d0 + ((a[1] - d0) if a[1] < INF else relprec)
        u = self._smul(inv0, r)
        u = ({w: -c for w, c in u[0].items()}, u[1])
        total = term = inv0
        while True:
            term = self._strunc(self._smul(u, term), prec)
            if not term[0]:
                break
            total = self._sadd(total, term)
        return self._strunc(total, prec)

    def _scalar_index(self, verts: Sequence[str]) -> list[tuple[str, int]]:
        return [(v, i) for v in verts for i in range(self.block_size(v))]

    def _gauss_jordan(self, x: Element, rverts: Sequence[str], cverts: Sequence[str]) -> Element:
        """Inverse of ``x`` from ``e_rverts`` to ``e_cverts``. For an exact ``x`` the series are
        recomputed with more precision (up to 8 relprec) when cancellations during the elimination
        leave ``x x^-1 = 1`` known to less than degree relprec; the result carries its true precision."""
        rp = self.relprec
        out = None
        for attempt in range(4):
            try:
                out = self._gj(x, rverts, cverts, rp)
            except PrecisionError:
                if x.prec < INF or attempt == 3:
                    raise
                rp *= 2
                continue
            if x.prec < INF or out.prec == INF or out.prec + self.valuation(x) >= self.relprec or attempt == 3:
                return out
            rp *= 2
        assert out is not None
        return out

    def _gj(self, x: Element, rverts: Sequence[str], cverts: Sequence[str], relprec: float) -> Element:
        rows, cols = self._scalar_index(rverts), self._scalar_index(cverts)
        if len(rows) != len(cols):
            raise NotInvertibleError(
                f"no inverse between vertices: {list(cverts)} and {list(rverts)} have different total "
                f"block sizes ({len(cols)} and {len(rows)})"
            )
        n = len(rows)
        ri = {r: k for k, r in enumerate(rows)}
        ci = {c: k for k, c in enumerate(cols)}
        L: list[list[tuple[dict, float]]] = [[({}, x.prec) for _ in range(n)] for _ in range(n)]
        for (t, s, w), X in x.terms.items():
            for i in range(X.shape[0]):
                for j in range(X.shape[1]):
                    if X[i, j] != 0:
                        if (t, i) not in ri or (s, j) not in ci:
                            raise ValueError("the element has terms outside the given corners")
                        L[ri[(t, i)]][ci[(s, j)]][0][w] = X[i, j]
        R: list[list[tuple[dict, float]]] = [
            [({(): 1.0}, INF) if i == j else ({}, INF) for j in range(n)] for i in range(n)
        ]
        used, pivrow = [False] * n, {}
        for k in range(n):
            best = None
            for i in range(n):
                if used[i]:
                    continue
                lu = self._sleading(L[i][k])
                if lu is None:
                    continue
                key = (self._sval(L[i][k]), -abs(lu[1]))
                if best is None or key < best[0]:
                    best = (key, i)
            if best is None:
                msg = (
                    f"no pivot for index {cols[k]}: no remaining entry of that column has a single word as "
                    "its lowest-degree part"
                )
                if any(not L[i][k][0] and L[i][k][1] < INF for i in range(n) if not used[i]):
                    raise PrecisionError(msg + " (some entries are only known to be O(deg ...): precision lost)")
                raise NotInvertibleError(msg)
            p = best[1]
            zi = self._sinv(L[p][k], relprec)
            L[p] = [self._smul(zi, v) for v in L[p]]
            R[p] = [self._smul(zi, v) for v in R[p]]
            for i in range(n):
                f = L[i][k]
                if i == p or (not f[0] and f[1] == INF):
                    continue
                L[i] = [self._sadd(L[i][j], self._smul(f, L[p][j]), -1.0) for j in range(n)]
                R[i] = [self._sadd(R[i][j], self._smul(f, R[p][j]), -1.0) for j in range(n)]
            used[p] = True
            pivrow[k] = p
        out: dict[Key, np.ndarray] = {}
        prec = INF
        for k in range(n):  # row k of the inverse maps index rows[j] to cols[k]
            s, i = cols[k]
            for j in range(n):
                t, jj = rows[j]
                ser, pr = R[pivrow[k]][j]
                prec = min(prec, pr)
                for w, c in ser.items():
                    term = (s, t, w)
                    if term not in out:
                        out[term] = np.zeros((self.block_size(s), self.block_size(t)), dtype=complex)
                    out[term][i, jj] += c
        return self.element(out, prec)

    # ------------------------------------------------------------------ square roots
    def sqrt(
        self, x: Element, rule: str = "graded", choice: str = "positive", split: str = "left", rng: Any = None
    ) -> Element:
        """A *-square root of ``x`` (some ``y`` with ``y* y = x``); see :func:`dsqlib.sqrt.star_sqrt`.

        With the default arguments, ``algebra.sqrt`` can be passed directly as the ``sqrt`` argument of
        the functions in :mod:`dsqlib.linalg`.
        """
        from .sqrt import star_sqrt

        return star_sqrt(x, rule, choice, split, rng)

    # ------------------------------------------------------------------ back to the path algebra
    def is_path_element(self, x: Element) -> bool:
        """Whether ``x`` has no inverse letters."""
        return all(e > 0 for k in x.terms for _, e in k[2])

    def reduce(self, x: Element, m: int = 2) -> Element:
        """The image of ``x`` in ``RQ / RQ_{>=m}`` (``m = 2`` gives the radical-square-zero algebra
        ``B_Q``). ``x`` must lie in the completed path algebra (no inverse letters) and be known up to
        degree ``m``."""
        bad = {k: X for k, X in x.terms.items() if any(e < 0 for _, e in k[2])}
        if bad:
            raise ValueError("the element has poles (terms with local inverses): " + self.format(Element(self, bad)))
        if x.prec < m:
            raise ValueError(f"the element is only known below degree {x.prec}; cannot reduce modulo degree {m}")
        return self.element({k: X for k, X in x.terms.items() if len(k[2]) < m}, INF)

    # ------------------------------------------------------------------ display
    def format(self, x: Element, digits: int = 4) -> str:
        """``x`` as a sum of ``coefficient·word`` terms, by degree."""
        if not x.terms:
            return "0" if x.prec == INF else f"O(deg {x.prec})"
        parts = []
        for (_, s, w), X in sorted(x.terms.items(), key=lambda kv: (self.degree(kv[0][2]), str(kv[0]))):
            if X.size == 1:
                c = X[0, 0]
                cs = (
                    f"{c.real:.{digits}g}"
                    if abs(c.imag) <= 1e-14 * abs(c)
                    else f"({c.real:.{digits}g}{c.imag:+.{digits}g}j)"
                )
            else:
                cs = np.array2string(np.round(X, digits), separator=",").replace("\n", "")
            parts.append(f"{cs}·{self.word_str(w, s)}")
        tail = "" if x.prec == INF else f" + O(deg {x.prec})"
        return " + ".join(parts) + tail


def quaternion(a: float, b: float = 0.0, c: float = 0.0, d: float = 0.0) -> np.ndarray:
    """The 2x2 complex block of the quaternion ``a + bi + cj + dk`` (the real form of a block at a
    vertex with real structure ``C = J``)."""
    p, q = complex(a, b), complex(c, d)
    return np.array([[p, q], [-np.conj(q), np.conj(p)]], dtype=complex)
\end{lstlisting}
\needspace{6\baselineskip}
\subsection{\texorpdfstring{\texttt{src/dsqlib/sqrt.py}}{src/dsqlib/sqrt.py}}
\begin{lstlisting}[style=dsqcode]
"""*-square roots: given ``x = x*``, some ``y`` with ``y* y = x``.

If ``y`` and ``y'`` are both *-square roots of ``x`` and ``y`` is invertible between the vertices it
touches, then ``y' = U y`` with ``U = y' y^+`` unitary there; so a *-square root is a choice, and the
rules below make that choice consistently. Every rule checks its result (``y* y = x`` word by word,
below the precision of both sides, and ``c(y) = y`` when ``x`` is real) and raises
:class:`~dsqlib.errors.NoSquareRootError` when it cannot produce one.

Rules
-----
``'graded'`` (the default)
    For homogeneous ``x`` whose words all have even length. Each word is cut in the middle,
    ``w = u'^T u``, and ``y = sum Z_u u`` is built from the half-words ``u``; for paths ``y`` is
    homogeneous of half the degree. The words of ``x`` fall into independent groups, one for each
    pair ``{m, m^T}`` of middle vertices:

    * ``m = m^T``: a congruence ``Z^T G_m^-1 Z = N`` through the ``d_m``-dimensional middle;
    * ``m != m^T``: a split ``Z'^T Z = N`` between the half-words at ``m`` and at ``m^T``.

    For real ``x`` the conjugation decides how the solution is made real: a Hermitian congruence when
    ``m~ = m^T != m`` (vertex kind Cu); explicit formulas in the real middle algebras R, M2(R), H and the
    exchanged pairs when ``m~ = m``; otherwise the solution at ``m`` is transported to ``m~``. Where
    the signature of a vertex of kind R, Hsp or Cu forbids a root, or ``x`` has rank greater than
    ``d_m`` at a middle vertex, there is no graded root and the rule says so.

    ``choice='positive'`` takes principal square roots of the blocks; ``choice='random'`` then
    multiplies by a random unitary of the semisimple part. For an exchanged pair ``(s, p)`` (``s`` the
    vertex listed first), ``split='left'`` puts ``x`` on ``s`` and 1 on ``p``, ``'right'`` the reverse
    and ``'balanced'`` equal square roots on both: over the split-complex numbers these are the roots
    ``(c, 1)``, ``(1, c)`` and ``(sqrt c, sqrt c)`` of ``(c, c)``.

    For ``x`` that is not homogeneous, the lowest-degree part ``x0`` gets the graded root ``y0`` and
    ``y = r y0`` with ``r`` the principal square root of ``(y0^+)* x y0^+``, as
    ``sqrt(t^2k (c + ...)) = t^k sqrt(c + ...)``.

``'principal'``
    For ``x`` with an invertible, vertex-diagonal lowest part of degree 0 in the semisimple part: the
    self-adjoint ``y`` with ``y^2 = x`` whose degree-0 part is the principal square root, found degree by
    degree from Sylvester equations.

``'arbitrary'``
    ``U y`` for the graded (or principal) root ``y`` and a random unitary ``U`` of the orbit corner where
    ``y`` ends: a Cayley transform with terms of positive degree, times a random unitary of the
    semisimple part. Not homogeneous; matrix entries stay in their orbit corners.

Known gaps of the graded rule (a root exists but no rule finds it): words with inverse letters, where
cancellation can hide the half-words; and an ``x`` whose lowest-degree part does not reach every vertex
that ``x`` reaches (e.g. ``(e_v + a)*(e_v + a)``, which has the root ``e_v + a``).
"""

from __future__ import annotations

import itertools
from collections.abc import Callable, Sequence
from typing import Any

import numpy as np
from scipy.linalg import solve_sylvester, sqrtm

from .algebra import INF, Algebra, Element, Key
from .errors import NoSquareRootError, NotInvertibleError

__all__ = [
    "CHOICES",
    "RULES",
    "SPLITS",
    "random_semisimple_unitary",
    "random_unitary",
    "sqrt_graded",
    "sqrt_principal",
    "sqrt_rule",
    "star_sqrt",
]

RULES = ("graded", "principal", "arbitrary")
CHOICES = ("positive", "random")
SPLITS = ("left", "right", "balanced")


def _psqrt(N: np.ndarray) -> np.ndarray:
    return np.asarray(sqrtm(np.asarray(N, dtype=complex)), dtype=complex)


def _rank(N: np.ndarray, rtol: float = 1e-9) -> int:
    if N.size == 0:
        return 0
    s = np.linalg.svd(N, compute_uv=False)
    return int(np.sum(s > rtol * s[0])) if s[0] > 0 else 0


def _close(A: np.ndarray, B: np.ndarray, rtol: float = 1e-9) -> bool:
    return bool(np.abs(A - B).max(initial=0) <= rtol * max(np.abs(A).max(initial=0), np.abs(B).max(initial=0)))


# ---------------------------------------------------------------------- block solvers
def _bilinear_candidates(M: np.ndarray, N: np.ndarray) -> list[np.ndarray]:
    """Solutions Z of Z^T M Z = N for square blocks (M invertible, N of the same symmetry as M), the
    principal one first."""
    d = M.shape[0]
    if d == 1:
        z = np.sqrt(complex(N[0, 0] / M[0, 0]))
        return [np.array([[z]]), np.array([[-z]])]
    if _close(M, M.T) and _close(N, N.T):
        return [np.linalg.solve(_psqrt(M), _psqrt(N))]
    if _close(M, -M.T) and _close(N, -N.T):  # 2x2 skew:  Z^T M Z = det(Z) M
        r = N[0, 1] / M[0, 1]
        return [np.sqrt(complex(r)) * np.eye(2, dtype=complex), np.diag([r, 1.0]).astype(complex)]
    return []


def _quaternion_symmetric_sqrt(N: np.ndarray) -> np.ndarray:
    """A quaternion block Z = [[a, c], [-c~, a~]] with Z^T Z = N, for N = [[p, i b], [i b, p~]] (b real);
    it always exists for N != 0 (the involution of orthogonal type on H has no signature)."""
    p, beta = N[0, 0], (N[0, 1] / 1j).real
    s = np.hypot(beta, abs(p))
    rho, sig = np.sqrt(max(beta + s, 0.0)), np.sqrt(max(s - beta, 0.0))
    ph = np.exp(0.5j * np.angle(p)) if abs(p) > 0 else 1.0
    P, Q = rho * ph, sig * ph  # P Q = p, |P|^2 - |Q|^2 = 2 beta
    a, c = (P + Q) / 2, np.conj((P - Q) / 2j)
    return np.array([[a, c], [-np.conj(c), np.conj(a)]], dtype=complex)


def _real_square(M: np.ndarray, N: np.ndarray, C: np.ndarray) -> np.ndarray:
    """Z^T M Z = N with Z of real type for the structure C (C = 1: real matrices, C = J: quaternions),
    square blocks of size 1 or 2."""
    d = M.shape[0]
    if not np.any(N):
        return np.zeros((d, d), dtype=complex)
    if d == 1:
        r = N[0, 0] / M[0, 0]
        if abs(r.imag) <= 1e-9 * (1 + abs(r)) and r.real >= 0:
            return np.array([[np.sqrt(r.real)]], dtype=complex)
        raise NoSquareRootError(f"no real solution of z^2 = {r:.4g} (signature)")
    quaternionic = not _close(C, np.eye(2))
    if _close(M, -M.T):  # Rsp or Hsp middle: det(Z) = r
        r = N[0, 1] / M[0, 1]
        if abs(r.imag) > 1e-9 * (1 + abs(r)):
            raise NoSquareRootError("block congruence: N is not a real multiple of M")
        r = r.real
        if r > 0:
            return np.sqrt(r) * np.eye(2, dtype=complex)
        if quaternionic:
            raise NoSquareRootError(f"no quaternion of norm {r:.4g} (signature)")
        return np.diag([r, 1.0]).astype(complex)
    if quaternionic:  # Hor middle: G = M = 1
        if not _close(M, np.eye(2)):
            raise NoSquareRootError("block congruence: unexpected form at a quaternionic vertex")
        if not (_close(N, N.T) and _close(N, np.array([[N[0, 0], N[0, 1]], [N[0, 1], np.conj(N[0, 0])]]))):
            raise NoSquareRootError("block congruence: N is not a symmetric quaternion")
        return _quaternion_symmetric_sqrt(N)
    Z = np.linalg.solve(_psqrt(M), _psqrt(N))
    if np.abs(Z.imag).max() > 1e-9 * (1 + np.abs(Z).max()):
        raise NoSquareRootError("no real solution of the block congruence")
    return Z.real.astype(complex)


def _pivots(N: np.ndarray, k: int, rows: bool | None = None) -> tuple[list[int], list[int]]:
    """Index sets (I, J) of size k maximising |det N[I, J]| (I = J for the symmetric problems)."""
    R, Cn = N.shape
    best = -1.0
    arg: tuple[list[int], list[int]] = ([], [])
    for J in itertools.combinations(range(Cn), k):
        Is = [J] if rows is None else itertools.combinations(range(R), k)
        for I in Is:
            v = abs(np.linalg.det(N[np.ix_(I, J)]))
            if v > best:
                best, arg = v, (list(I), list(J))
    return arg


def _sym_factor(M: np.ndarray, N: np.ndarray, d: int) -> np.ndarray:
    """Some complex Z (d x D) with Z^T M Z = N, N of rank <= d."""
    D = N.shape[0]
    r = _rank(N)
    if r == 0:
        return np.zeros((d, D), dtype=complex)
    if r > d:
        raise NoSquareRootError(f"x has rank {r} at a middle vertex of size {d}: no graded *-square root")
    if d == D and r == d:
        Zs = _bilinear_candidates(M, N)
        if Zs and _close(Zs[0].T @ M @ Zs[0], N):
            return Zs[0]
        raise NoSquareRootError("the block congruence has no solution")
    k = r
    Mk = M[:k, :k] if k < d else M
    if k < d and not _close(M, np.diag(np.diag(M))):
        raise NoSquareRootError("x is degenerate at a middle vertex with a skew form")
    I, _ = _pivots(N, k)
    ZJ = _bilinear_candidates(Mk, N[np.ix_(I, I)])
    if not ZJ:
        raise NoSquareRootError("the block congruence has no solution")
    Zk = np.linalg.solve(Mk, np.linalg.solve(ZJ[0].T, N[I, :]))
    Z = np.vstack([Zk, np.zeros((d - k, D), dtype=complex)])
    if not _close(Z.T @ M @ Z, N):
        raise NoSquareRootError("graded rule: x does not factor through the middle vertex")
    return Z


def _split_factor(N: np.ndarray, d: int, split: str) -> tuple[np.ndarray, np.ndarray]:
    """Some (Z', Z), each with d rows, with Z'^T Z = N (rank <= d)."""
    R, Cn = N.shape
    r = _rank(N)
    if r == 0:
        return np.zeros((d, R), dtype=complex), np.zeros((d, Cn), dtype=complex)
    if r > d:
        raise NoSquareRootError(f"x has rank {r} at a middle vertex of size {d}: no graded *-square root")
    if R == Cn == r:
        I = J = list(range(r))
    else:
        I, J = _pivots(N, r, rows=True)
    P = N[np.ix_(I, J)]
    if split == "left":
        Zk, Zpk = N[I, :], np.linalg.solve(P.T, N[:, J].T)
    elif split == "right":
        Zk, Zpk = np.linalg.solve(P, N[I, :]), N[:, J].T
    elif split == "balanced":
        if r == 1:
            S = np.sqrt(abs(P[0, 0])) * np.ones((1, 1))
            Zk, Zpk = N[I, :] / S[0, 0], (N[:, J] * S[0, 0] / P[0, 0]).T
        else:
            S = _psqrt(P)
            Zk, Zpk = np.linalg.solve(S, N[I, :]), np.linalg.solve(S.T, N[:, J].T)
    else:
        raise ValueError(f"split must be one of {SPLITS}, not {split!r}")

    def pad(Z: np.ndarray) -> np.ndarray:
        return np.vstack([Z, np.zeros((d - r, Z.shape[1]), dtype=complex)])

    Zp, Z = pad(Zpk), pad(Zk)
    if not _close(Zp.T @ Z, N):
        raise NoSquareRootError("graded rule: x does not factor through the middle vertices")
    return Zp, Z


def _phase_fix(v: np.ndarray) -> np.ndarray:
    k = int(np.argmax(np.abs(v)))
    return v * (np.conj(v[k]) / abs(v[k])) if abs(v[k]) > 0 else v


def _hermitian_factor(H: np.ndarray, N: np.ndarray, d: int) -> np.ndarray:
    """Some Z (d x D) with Z^H H Z = N (H invertible; H and N both Hermitian or both anti-Hermitian).
    For H = +-1 and a square block this is the principal square root of +-N."""
    if not _close(H, H.conj().T):
        H, N = -1j * H, -1j * N
    if not (_close(H, H.conj().T) and _close(N, N.conj().T)):
        raise NoSquareRootError("graded rule: x is not Hermitian at a middle vertex")
    N = (N + N.conj().T) / 2
    D = N.shape[0]
    if _rank(N) > d:
        raise NoSquareRootError(f"x has rank {_rank(N)} at a middle vertex of size {d}: no graded *-square root")
    for h in (1.0, -1.0):
        if d == D and _close(H, h * np.eye(d)):
            ev = np.linalg.eigvalsh(h * N)
            if ev.min() < -1e-12 * abs(ev).max():
                raise NoSquareRootError(
                    "graded rule: x is not of the sign of the form at a Cu-type middle vertex (signature)"
                )
            return _psqrt(h * N)
    eh, Wh = np.linalg.eigh(H)
    en, Wn = np.linalg.eigh(N)
    tol = 1e-10 * abs(en).max()
    Z = np.zeros((d, D), dtype=complex)
    free = list(range(d))
    for j in np.argsort(-abs(en)):
        if abs(en[j]) <= tol:
            continue
        match = [i for i in free if np.sign(eh[i]) == np.sign(en[j])]
        if not match:
            raise NoSquareRootError("graded rule: x does not fit the signature of the middle vertex")
        i = match[0]
        free.remove(i)
        Z += np.sqrt(en[j] / eh[i]) * np.outer(Wh[:, i], _phase_fix(Wn[:, j]).conj())
    if not _close(Z.conj().T @ H @ Z, N):
        raise NoSquareRootError("graded rule: x has too large a rank at a middle vertex")
    return Z


# ---------------------------------------------------------------------- the graded rule
_HalfKey = tuple[str, Any]  # (source vertex, half-word)


def _order(alg: Algebra) -> Callable[[_HalfKey], tuple[int, int, str]]:
    pos = {v: i for i, v in enumerate(alg.vertices)}
    return lambda k: (pos[k[0]], len(k[1]), str(k[1]))


def _layout(alg: Algebra, keys: Sequence[_HalfKey]) -> tuple[dict[_HalfKey, tuple[int, int]], int]:
    """Column offsets for a list of (source, half-word) keys."""
    off, out = 0, {}
    for k in keys:
        out[k] = (off, off + alg.block_size(k[0]))
        off += alg.block_size(k[0])
    return out, off


def _matrix(
    alg: Algebra, blocks: dict, rows: Sequence[_HalfKey], cols: Sequence[_HalfKey], transform: Callable
) -> np.ndarray:
    ro, R = _layout(alg, rows)
    co, C = _layout(alg, cols)
    N = np.zeros((R, C), dtype=complex)
    for (r, c), X in blocks.items():
        a, b = ro[r]
        e, f = co[c]
        N[a:b, e:f] = transform(r, c, X)
    return N


def _terms(alg: Algebra, m: str, keys: Sequence[_HalfKey], Z: np.ndarray) -> Element:
    lay, _ = _layout(alg, keys)
    t: dict[Key, np.ndarray] = {}
    for s, u in keys:
        a, b = lay[(s, u)]
        if np.abs(Z[:, a:b]).max(initial=0) > 0:
            t[(m, s, u)] = Z[:, a:b]
    return Element(alg, t)


def _blockrow(alg: Algebra, y: Element, v: str, keys: Sequence[_HalfKey]) -> np.ndarray:
    lay, D = _layout(alg, keys)
    Z = np.zeros((alg.block_size(v), D), dtype=complex)
    for (t, s, u), X in y.terms.items():
        if t == v:
            a, b = lay[(s, u)]
            Z[:, a:b] = X
    return Z


def _realify(alg: Algebra, y0: Element, split: str) -> Element:
    """A real *-square root of y0* y0 from a complex one: y = r y1 with y1 the real part of B y0 (B = 1,
    then fixed other choices) and r a degree-0 *-square root of K = g* g, where g y1 = y0."""
    mids = sorted({t for (t, _, _) in y0.terms}, key=list(alg.vertices).index)
    rng = np.random.default_rng(12345)
    for attempt in range(6):
        if attempt == 0:
            Bel = alg.one()
        else:
            Bel = alg.zero()
            for v in mids:
                d = alg.block_size(v)
                Bel = alg.add(Bel, alg.vertex(v, rng.normal(size=(d, d)) + 1j * rng.normal(size=(d, d))))
        y1 = alg.real(alg.mul(Bel, y0))
        keys = sorted({(s, u) for (_, s, u) in list(y0.terms) + list(y1.terms)}, key=_order(alg))
        G, ok = alg.zero(), True
        for v in mids:
            Z0, Z1 = _blockrow(alg, y0, v, keys), _blockrow(alg, y1, v, keys)
            g = Z0 @ np.linalg.pinv(Z1)
            if not _close(g @ Z1, Z0):
                ok = False
                break
            G = alg.add(G, Element(alg, {(v, v, ()): g}))
        if not ok:
            continue
        K = alg.real(alg.mul(alg.star(G), G))
        r = _graded(alg, K, split, base=True)
        return alg.mul(r, y1)
    raise NoSquareRootError("graded rule: could not find a real *-square root")


def _graded(alg: Algebra, x: Element, split: str, base: bool = False) -> Element:
    if not x.terms:
        return alg.zero()
    degs = alg.degrees(x)
    if len(degs) != 1:
        raise NoSquareRootError(f"graded rule: x is not homogeneous (degrees {degs})")
    V = alg.vertices
    real = alg.is_real(x)
    groups: dict[str, dict] = {}
    for (t, s, w), X in x.terms.items():
        if len(w) % 2:
            raise NoSquareRootError(f"graded rule: the word {alg.word_str(w, s)} has odd length")
        h = len(w) // 2
        u, up = w[h:], alg.word_transpose(w[:h])
        m = alg.word_ends(u)[0] if u else s
        groups.setdefault(m, {})[((V[t].top, up), (s, u))] = X
    key = _order(alg)

    def norm(r: _HalfKey, c: _HalfKey, X: np.ndarray) -> np.ndarray:  # N[u', u] = kappa(u') G_s'^-1 X
        return alg.word_kappa(r[1]) * V[r[0]].G_inv @ X

    y, done = alg.zero(), set()
    for m in sorted(groups, key=list(V).index):
        if m in done:
            continue
        blocks, mT, mb, d = groups[m], V[m].top, V[m].bar, V[m].block_size
        cols = sorted({c for (_, c) in blocks}, key=key)
        rows = sorted({r for (r, _) in blocks}, key=key)
        if m == mT:  # Z^T G_m^-1 Z = N
            idx = sorted(set(cols) | set(rows), key=key)
            N, M = _matrix(alg, blocks, idx, idx, norm), V[m].G_inv
            if real and mb == m:
                single = len(idx) == 1 and V[idx[0][0]].block_size == d and alg.word_conjugate(idx[0][1]) == idx[0][1]
                if single:
                    yg = _terms(alg, m, idx, _real_square(M, N, V[m].C))
                elif base:
                    raise NoSquareRootError("graded rule: internal error (base case with several half-words)")
                else:
                    y0 = _terms(alg, m, idx, _sym_factor(M, N, d))
                    yg = y0 if alg.is_real(y0) else _realify(alg, y0, split)
            else:
                yg = _terms(alg, m, idx, _sym_factor(M, N, d))
                if real:
                    yg = alg.add(yg, alg.conj(yg))
                    done.add(mb)
            done.add(m)
        else:
            N = _matrix(alg, blocks, rows, cols, norm)  # Z'^T Z = N  (G = 1 off the T-fixed vertices)
            if real and mb == mT:  # Hermitian: Z' = c(Z)
                conj_cols = [(V[s].bar, alg.word_conjugate(u)) for (s, u) in cols]
                if set(conj_cols) != set(rows):
                    raise NoSquareRootError("graded rule: the half-words at m^T are not the conjugates of those at m")
                Nc = _matrix(alg, blocks, conj_cols, cols, norm)
                lay, _ = _layout(alg, cols)
                Ns = Nc.copy()
                for s, u in cols:  # N'[v, u] = mu(v) C_s(v)^T N[v~, u]
                    a, b = lay[(s, u)]
                    Ns[a:b, :] = alg.word_mu(u) * V[s].C.T @ Nc[a:b, :]
                Hm = V[m].C.T @ V[mT].G_inv
                Z = _hermitian_factor(Hm, Ns, d)
                yg = _terms(alg, m, cols, Z)
                yg = alg.add(yg, alg.conj(yg))
            else:
                Zp, Z = _split_factor(N, d, split)
                yg = alg.add(_terms(alg, m, cols, Z), _terms(alg, mT, rows, Zp))
                if real and mb == m:
                    if not alg.is_real(yg):
                        if base:
                            found = None
                            for sp in SPLITS:
                                Zp, Z = _split_factor(N, d, sp)
                                cand = alg.add(_terms(alg, m, cols, Z), _terms(alg, mT, rows, Zp))
                                if alg.is_real(cand):
                                    found = cand
                                    break
                            if found is None:
                                raise NoSquareRootError("graded rule: no real split at an exchanged pair")
                            yg = found
                        else:
                            yg = _realify(alg, yg, split)
                elif real:
                    yg = alg.add(yg, alg.conj(yg))
                    done |= {mb, V[mb].top}
            done |= {m, mT}
        y = alg.add(y, yg)
    return y


def sqrt_graded(
    x: Element, choice: str = "positive", split: str = "left", rng: np.random.Generator | None = None
) -> Element:
    """The graded *-square root of ``x`` (see the module documentation).

    For ``x`` that is not homogeneous, ``y = r y0`` where ``y0`` is the graded root of the lowest-degree
    part ``x0`` and ``r`` the principal square root of ``(y0^+)* x y0^+``, whose degree-0 part is the
    identity of the vertices where ``y0`` ends.

    Raises
    ------
    NoSquareRootError
        If the rule cannot produce a root.
    """
    if choice not in CHOICES:
        raise ValueError(f"choice must be one of {CHOICES}, not {choice!r}")
    if split not in SPLITS:
        raise ValueError(f"split must be one of {SPLITS}, not {split!r}")
    alg = x.algebra
    if not x.terms:  # x = O(deg p): y = O(deg floor(p/2)) is consistent
        y = alg.zero(prec=(x.prec // 2) if x.prec < INF else INF)
    elif len(alg.degrees(x)) == 1:
        y = _graded(alg, x, split)
        if x.prec < INF and y.terms:  # y is known modulo degree prec(x) - val(y)
            y = alg.truncate(y, x.prec - alg.valuation(y))
    else:
        y0 = _graded(alg, alg.part(x, int(alg.valuation(x))), split)
        C = {s for (_, s, _) in y0.terms}
        if any(t not in C or s not in C for (t, s, _) in x.terms):
            raise NoSquareRootError(
                "graded rule: x reaches vertices where its lowest-degree part vanishes; the graded and "
                "principal rules do not cover this case (e.g. (e_v + a)*(e_v + a) has the root e_v + a)"
            )
        try:
            y0i = alg.inverse(y0)
        except NotInvertibleError as e:
            raise NoSquareRootError(f"graded rule: the root of the lowest-degree part is not invertible ({e})") from e
        m = alg.mul(alg.mul(alg.star(y0i), x), y0i)
        y = alg.mul(sqrt_principal(m), y0)
    if choice == "random":
        y = alg.mul(random_semisimple_unitary(alg, rng or np.random.default_rng()), y)
    return _check(alg, x, y, "graded rule")


# ---------------------------------------------------------------------- the principal rule
def sqrt_principal(x: Element) -> Element:
    """The principal *-square root: ``y = y*`` with ``y^2 = x``, for ``x`` with an invertible,
    vertex-diagonal degree-0 lowest part in the semisimple part.

    Raises
    ------
    NoSquareRootError
        If ``x`` is not of that form, or the result fails the check.
    """
    alg = x.algebra
    if not x.terms:
        return alg.zero()
    if alg.valuation(x) < 0 or 0 not in alg.degrees(x):
        raise NoSquareRootError("principal rule: needs a lowest part of degree 0 (use rule='graded')")
    x0 = alg.part(x, 0)
    if any(w for (_, _, w) in x0.terms):
        raise NoSquareRootError("principal rule: the degree-0 part must lie in the semisimple part S")
    if any(t != s for (t, s, _) in x0.terms):
        raise NoSquareRootError("principal rule: the degree-0 part must be vertex-diagonal")
    Y0: dict[str, np.ndarray] = {}
    y = alg.zero()
    for (_t, s, _w), X in x0.terms.items():
        if np.linalg.cond(X) > 1e12:
            raise NoSquareRootError(f"principal rule: the degree-0 part of x is not invertible at vertex {s}")
        Y0[s] = _psqrt(X)
        y = alg.add(y, Element(alg, {(s, s, ()): Y0[s]}))
    missing = ({k[0] for k in x.terms} | {k[1] for k in x.terms}) - set(Y0)
    if missing:
        raise NoSquareRootError(
            f"principal rule: x reaches the vertices {sorted(missing)}, where its degree-0 part vanishes"
        )
    prec = min(x.prec, alg.relprec)
    parts = {0: y}
    for d in range(1, int(prec)):
        rhs = alg.part(x, d)
        for i in range(1, d):
            rhs = alg.sub(rhs, alg.mul(parts[i], parts[d - i]))
        yd = alg.zero()
        for (t, s, w), Rb in rhs.terms.items():
            yd = alg.add(yd, Element(alg, {(t, s, w): solve_sylvester(Y0[t], Y0[s], Rb)}))
        parts[d] = yd
        y = alg.add(y, yd)
    y = alg.truncate(y, prec)
    xx = alg.truncate(x, prec)
    if not alg.is_selfadjoint(y):
        raise NoSquareRootError("principal rule: the result is not self-adjoint (is the degree-0 part positive?)")
    return _check(alg, xx, y, "principal rule")


# ---------------------------------------------------------------------- unitaries and the arbitrary rule
def random_semisimple_unitary(alg: Algebra, rng: np.random.Generator, signs: bool = True) -> Element:
    """A random real unitary of the semisimple part: for each STAR-BAR orbit ``O``, the Cayley
    transform ``(e_O - K)(e_O + K)^-1`` of a random skew-adjoint ``K`` in ``e_O S e_O``, times a random
    sign (if ``signs``)."""
    U = alg.zero()
    for orb in alg.quiver.star_bar_orbits():
        R, eo = alg.zero(), alg.zero()
        for v in orb:
            d = alg.block_size(v)
            R = alg.add(R, alg.vertex(v, rng.normal(size=(d, d)) + 1j * rng.normal(size=(d, d))))
            eo = alg.add(eo, alg.idempotent(v))
        R = alg.real(R)
        K = alg.sub(R, alg.star(R))
        Uc = alg.mul(alg.sub(eo, K), alg.inverse(alg.add(eo, K)))
        if signs and rng.random() < 0.5:
            Uc = alg.scale(-1.0, Uc)
        U = alg.add(U, Uc)
    return U


def _paths(alg: Algebra, max_len: int) -> list[Key]:
    out: list[Key] = [(v, v, ()) for v in alg.vertices]
    frontier = list(out)
    for _ in range(max_len):
        new = []
        for t, s, w in frontier:
            for n, a in alg.arrows.items():
                if a.source == t:
                    new.append((a.target, s, ((n, 1), *w)))
        out += new
        frontier = new
    return out


def random_unitary(
    alg: Algebra,
    rng: np.random.Generator,
    vertices: Sequence[str] | set[str] | None = None,
    max_len: int = 2,
    scale: float = 0.5,
) -> Element:
    """A random real unitary of the corner ``e_T A e_T``, ``T = vertices`` (default: all; ``T`` should
    be a union of STAR-BAR orbits).

    It is the Cayley transform ``(e_T - K)(e_T + K)^-1`` of ``K = R - R*``, where ``R`` is a random real
    combination of the idempotents ``e_v`` (``v`` in ``T``) and of the paths of length ``<= max_len``
    between vertices of ``T``. It has terms of positive degree: a unitary of the completion.
    """
    T = set(alg.vertices) if vertices is None else set(vertices)
    R = alg.zero()
    for t, s, w in _paths(alg, max_len):
        if t in T and s in T:
            dt, ds = alg.block_size(t), alg.block_size(s)
            X = scale * (rng.normal(size=(dt, ds)) + 1j * rng.normal(size=(dt, ds)))
            R = alg.add(R, Element(alg, {(t, s, w): X}))
    R = alg.real(R)
    K = alg.sub(R, alg.star(R))
    eT = alg.zero()
    for v in T:
        eT = alg.add(eT, alg.idempotent(v))
    return alg.mul(alg.sub(eT, K), alg.inverse_between(alg.add(eT, K), T, T))


def _magnitudes(alg: Algebra, x: Element, y: Element) -> dict[Key, float]:
    """For each word of the product x y, the sum over its contributions of |X| |Y| (block maxima times the
    inner dimension): the scale on which rounding errors in that coefficient live."""
    by_target: dict[str, list] = {}
    for (t2, s2, w2), B in y.terms.items():
        by_target.setdefault(t2, []).append((s2, w2, np.abs(B).max()))
    prec = min(x.prec + alg.valuation(y), y.prec + alg.valuation(x))
    mag: dict[Key, float] = {}
    for (t1, s1, w1), A in x.terms.items():
        a = np.abs(A).max() * A.shape[1]
        for s2, w2, b in by_target.get(s1, ()):
            w = alg.concat(w1, w2)
            if alg.degree(w) < prec:
                k = (t1, s2, w)
                mag[k] = mag.get(k, 0.0) + a * b
    return mag


def _check(alg: Algebra, x: Element, y: Element, what: str) -> Element:
    """y* y = x word by word, below the precision of both, and c(y) = y when x is real. At each word the
    error may be 1e-8 of the larger of the coefficient of x and the size of the contributions to that
    coefficient of y* y (rounding lives on that scale), plus 1e-13 of the largest coefficient of x; so a
    word of x that y* y misses is always caught."""
    ys = alg.star(y)
    yy = alg.mul(ys, y)
    mag = _magnitudes(alg, ys, y)
    prec = min(x.prec, yy.prec)
    big = max((np.abs(X).max() for X in x.terms.values()), default=0.0)
    for k in set(x.terms) | set(yy.terms):
        if alg.degree(k[2]) >= prec:
            continue
        X, Y = x.terms.get(k), yy.terms.get(k)
        a = np.abs(X).max() if X is not None else 0.0
        e = np.abs((Y if Y is not None else 0) - (X if X is not None else 0)).max()
        if e > 1e-8 * max(a, mag.get(k, 0.0)) + 1e-13 * big:
            raise NoSquareRootError(
                f"{what}: y* y differs from x at the word {alg.word_str(k[2], k[1])} by {e:.2e} "
                f"(coefficient of x there: {a:.2e})"
            )
    if alg.is_real(x) and not alg.is_real(y):
        raise NoSquareRootError(f"{what}: could not find a real *-square root")
    return y


def star_sqrt(
    x: Element,
    rule: str = "graded",
    choice: str = "positive",
    split: str = "left",
    rng: np.random.Generator | None = None,
) -> Element:
    """A *-square root of ``x``: some ``y`` with ``y* y = x``.

    Parameters
    ----------
    x : Element
        A self-adjoint element.
    rule : {'graded', 'principal', 'arbitrary'}
    choice : {'positive', 'random'}
        For the graded rule: principal square roots of the blocks, or those times a random unitary
        of the semisimple part.
    split : {'left', 'right', 'balanced'}
        For the graded rule: where ``x`` goes at an exchanged pair of vertices.
    rng : numpy.random.Generator, optional
        For the random choices.

    Returns
    -------
    Element

    Raises
    ------
    NoSquareRootError
        If the rule cannot produce a root.

    Examples
    --------
    >>> from dsqlib import Algebra, catalogue
    >>> A = Algebra(catalogue.svd("C"))
    >>> z = A.real_element(A.arrow("a", 1 + 2j))      # (1+2i) a + (1-2i) a~
    >>> star_sqrt(z.H * z)
    2.236·a + 2.236·a~
    """
    if rule == "graded":
        return sqrt_graded(x, choice, split, rng)
    if rule == "principal":
        return sqrt_principal(x)
    if rule == "arbitrary":
        alg = x.algebra
        rng = rng or np.random.default_rng()
        try:
            y = sqrt_graded(x, "positive", split)
        except NoSquareRootError:
            y = sqrt_principal(x)
        if y.terms:  # a random unitary of the orbit corner where y ends
            T = {v for (t, _, _) in y.terms for v in alg.quiver.orbit_of(t)}
            y = alg.mul(alg.mul(random_unitary(alg, rng, T), random_semisimple_unitary(alg, rng)), y)
        return _check(alg, alg.truncate(x, y.prec), y, "arbitrary rule")
    raise ValueError(f"rule must be one of {RULES}, not {rule!r}")


def sqrt_rule(
    rule: str = "graded", choice: str = "positive", split: str = "left", seed: int | None = None
) -> Callable[[Element], Element]:
    """A one-argument *-square root ``x -> y`` for the ``sqrt`` argument of :mod:`dsqlib.linalg`.

    Examples
    --------
    >>> sq = sqrt_rule("graded", "random", seed=1)   # reproducible random choices
    """
    if rule not in RULES:
        raise ValueError(f"rule must be one of {RULES}, not {rule!r}")
    if choice not in CHOICES:
        raise ValueError(f"choice must be one of {CHOICES}, not {choice!r}")
    if split not in SPLITS:
        raise ValueError(f"split must be one of {SPLITS}, not {split!r}")
    rng = np.random.default_rng(seed)

    def sqrt(x: Element) -> Element:
        return star_sqrt(x, rule, choice, split, rng)

    sqrt.__doc__ = f"*-square root by the {rule!r} rule (choice={choice!r}, split={split!r})"
    return sqrt
\end{lstlisting}
\needspace{6\baselineskip}
\subsection{\texorpdfstring{\texttt{src/dsqlib/matrix.py}}{src/dsqlib/matrix.py}}
\begin{lstlisting}[style=dsqcode]
"""Matrices over the completed localisation, indexed by self-adjoint idempotents.

A matrix index is a list of positions. Position ``i`` stands for a real self-adjoint idempotent
``E_i`` of the algebra, and entry ``(i, j)`` of a matrix ``H`` lives in the corner ``E_i A E_j``. This
generalises the usual indexing, where every ``E_i = 1``. Two indexings are provided:

:func:`plain_index`
    every ``E_i = 1``: ``n x n`` matrices over the algebra in the usual sense;
:func:`orbit_index`
    one counter per STAR-BAR orbit ``O``, standing for ``e_O`` = the sum of ``e_v`` over the vertices
    of ``O`` (real and self-adjoint, because ``O`` is closed under the transpose and the conjugation).
    The positions visit the orbits in a cyclic order and come back to the first, skipping an orbit
    whose counter is exhausted: for ``(n_u, n_v) = (3, 5)`` they are ``u0 v0 u1 v1 u2 v2 v3 v4``.

A :class:`Matrix` holds its entries together with a row index and a column index.
"""

from __future__ import annotations

import numbers
from collections.abc import Callable, Iterator, Mapping, Sequence
from typing import Any, overload

import numpy as np

from .algebra import INF, Algebra, Element

__all__ = ["Index", "Matrix", "cyclic_positions", "orbit_index", "plain_index"]

Position = tuple[str, int]


def cyclic_positions(dims: Mapping[str, int], order: Sequence[str]) -> list[Position]:
    """The positions ``(orbit, k)``, cycling through the orbits in the given order and skipping an
    orbit whose counter has reached ``dims[orbit]``.

    Examples
    --------
    >>> cyclic_positions({"u": 2, "v": 3}, ["u", "v"])
    [('u', 0), ('v', 0), ('u', 1), ('v', 1), ('v', 2)]
    """
    idx, k = [], 0
    while any(k < dims[o] for o in order):
        for o in order:
            if k < dims[o]:
                idx.append((o, k))
        k += 1
    return idx


class Index:
    """A list of positions ``(label, k)``; each label names a set of vertices whose idempotents add up
    to the idempotent ``E`` of the positions with that label.

    Parameters
    ----------
    algebra : Algebra
    positions : sequence of (label, k)
    groups : mapping of label to a tuple of vertices
    """

    def __init__(self, algebra: Algebra, positions: Sequence[Position], groups: Mapping[str, Sequence[str]]) -> None:
        self.algebra = algebra
        self.positions: list[Position] = [tuple(p) for p in positions]  # type: ignore[misc]
        self.groups: dict[str, tuple[str, ...]] = {k: tuple(v) for k, v in groups.items()}
        for label, _ in self.positions:
            if label not in self.groups:
                raise ValueError(f"position label {label!r} has no vertex group")
        self._E: dict[str, Element] = {}
        for label, vs in self.groups.items():
            e = algebra.zero()
            for v in vs:
                e = algebra.add(e, algebra.idempotent(v))
            self._E[label] = e

    def __len__(self) -> int:
        return len(self.positions)

    @property
    def n(self) -> int:
        """The number of positions."""
        return len(self.positions)

    @overload
    def __getitem__(self, i: int) -> Position: ...
    @overload
    def __getitem__(self, i: slice) -> Index: ...
    def __getitem__(self, i: int | slice) -> Position | Index:
        if isinstance(i, slice):
            return Index(self.algebra, self.positions[i], self.groups)
        return self.positions[i]

    def __iter__(self) -> Iterator[Position]:
        return iter(self.positions)

    def __eq__(self, other: object) -> bool:
        if not isinstance(other, Index):
            return NotImplemented
        return (
            self.algebra is other.algebra
            and self.positions == other.positions
            and all(self.groups[p[0]] == other.groups[p[0]] for p in self.positions)
        )

    def __hash__(self) -> int:
        return hash(tuple(self.positions))

    def __repr__(self) -> str:
        return "Index(" + ", ".join(str(k) if lab == "1" else f"{lab}{k}" for lab, k in self.positions) + ")"

    def idempotent(self, i: int) -> Element:
        """``E_i``, the idempotent of position ``i``."""
        return self._E[self.positions[i][0]]

    def label(self, i: int) -> str:
        """The label of position ``i``."""
        return self.positions[i][0]

    def vertices(self, i: int) -> tuple[str, ...]:
        """The vertices whose idempotents add up to ``E_i``."""
        return self.groups[self.positions[i][0]]

    def position(self, label: str, k: int) -> int:
        """The position number of ``(label, k)``."""
        return self.positions.index((label, k))

    def subset(self, positions: Sequence[int]) -> Index:
        """The index made of the given positions, in that order."""
        return Index(self.algebra, [self.positions[i] for i in positions], self.groups)

    def relabel(self, mapping: Mapping[str, str]) -> Index:
        """The same positions with the labels renamed by ``mapping`` (for example after a QR step on
        the SVD indexing, which maps each orbit to the other end of the arrows)."""
        return Index(self.algebra, [(mapping.get(lab, lab), k) for lab, k in self.positions], self.groups)


def plain_index(algebra: Algebra, n: int) -> Index:
    """``n`` positions with every ``E_i = 1`` (the label ``'1'`` stands for all vertices)."""
    return Index(algebra, [("1", k) for k in range(n)], {"1": tuple(algebra.vertices)})


def orbit_index(algebra: Algebra, dims: Mapping[str, int], order: Sequence[str] | None = None) -> Index:
    """One counter per STAR-BAR orbit, cycling through the orbits in ``order``.

    Parameters
    ----------
    algebra : Algebra
    dims : mapping of orbit name to the number of its positions
        Orbits are named by the base names given to :meth:`DSQ.build`.
    order : sequence of orbit names, optional
        The cyclic order (default: the order of the quiver).

    Examples
    --------
    >>> from dsqlib import Algebra, catalogue
    >>> A = Algebra(catalogue.svd("R"))
    >>> orbit_index(A, {"u": 3, "v": 5}, ["u", "v"])
    Index(u0, v0, u1, v1, u2, v2, v3, v4)
    """
    names = algebra.quiver.orbit_names()
    for o in dims:
        if o not in names:
            raise ValueError(f"unknown orbit {o!r}; the orbits are {list(names)}")
    order = list(order) if order else [o for o in names if o in dims]
    if sorted(order) != sorted(dims):
        raise ValueError(f"the cyclic order {order} must list each orbit of dims {list(dims)} once")
    return Index(algebra, cyclic_positions(dims, order), {o: names[o] for o in order})


def _exact_zero(x: Element) -> bool:
    return not x.terms and x.prec == INF


class Matrix:
    """A matrix over an :class:`~dsqlib.algebra.Algebra`, with a row and a column :class:`Index`.

    Entry ``(i, j)`` lives in the corner ``E_i A E_j`` (``E_i`` from the row index, ``E_j`` from the
    column index); :meth:`corner_violation` measures how far the entries are from that.

    Supports ``H[i, j]`` (an :class:`~dsqlib.algebra.Element`), slicing ``H[r0:r1, c0:c1]``, ``A @ B``,
    ``A + B``, ``A - B``, ``-A``, ``c * A`` for a number ``c``, and ``A.H`` (the adjoint: transpose with
    the involution applied entrywise).

    Parameters
    ----------
    entries : list of lists of Element
    rows : Index
    cols : Index, optional
        Defaults to ``rows``.
    """

    __slots__ = ("_e", "algebra", "cols", "rows")
    __array_ufunc__ = None

    def __init__(self, entries: Sequence[Sequence[Element]], rows: Index, cols: Index | None = None) -> None:
        cols = rows if cols is None else cols
        if rows.algebra is not cols.algebra:
            raise ValueError("the row and column indices belong to different algebras")
        self.algebra: Algebra = rows.algebra
        self.rows, self.cols = rows, cols
        self._e: list[list[Element]] = [list(r) for r in entries]
        if len(self._e) != len(rows) or any(len(r) != len(cols) for r in self._e):
            raise ValueError(f"entries are not {len(rows)} x {len(cols)}")
        for r in self._e:
            for x in r:
                if not isinstance(x, Element) or x.algebra is not self.algebra:
                    raise TypeError("Matrix entries must be Elements of the indices' algebra")

    # ------------------------------------------------------------------ construction
    @classmethod
    def zeros(cls, rows: Index, cols: Index | None = None) -> Matrix:
        """The zero matrix."""
        cols = rows if cols is None else cols
        z = rows.algebra.zero()
        return cls([[z for _ in range(len(cols))] for _ in range(len(rows))], rows, cols)

    @classmethod
    def identity(cls, index: Index) -> Matrix:
        """``diag(E_i)``, the identity of the matrices with this index."""
        M = cls.zeros(index)
        for i in range(len(index)):
            M._e[i][i] = index.idempotent(i)
        return M

    @classmethod
    def from_function(cls, rows: Index, f: Callable[[int, int], Element | None], cols: Index | None = None) -> Matrix:
        """The matrix with entries ``f(i, j)`` (``None`` for zero)."""
        cols = rows if cols is None else cols
        z = rows.algebra.zero()
        return cls([[(f(i, j) or z) for j in range(len(cols))] for i in range(len(rows))], rows, cols)

    @classmethod
    def unit_vector(cls, index: Index, i: int) -> Matrix:
        """The column vector with ``E_i`` in position ``i`` (its column index is position ``i``)."""
        i = range(len(index))[i]  # normalises negative positions, raises IndexError when out of range
        z = index.algebra.zero()
        return cls([[index.idempotent(i) if k == i else z] for k in range(len(index))], index, index.subset([i]))

    @classmethod
    def from_columns(cls, columns: Sequence[Sequence[Element]], rows: Index, cols: Index | None = None) -> Matrix:
        """The matrix with the given columns (each a list of entries). ``cols`` defaults to the first
        ``len(columns)`` positions of ``rows``."""
        cols = rows[: len(columns)] if cols is None else cols
        return cls([[c[i] for c in columns] for i in range(len(rows))], rows, cols)

    # ------------------------------------------------------------------ access
    @property
    def shape(self) -> tuple[int, int]:
        """``(number of rows, number of columns)``."""
        return (len(self.rows), len(self.cols))

    def __getitem__(self, key: tuple[Any, Any]) -> Any:
        if not (isinstance(key, tuple) and len(key) == 2):
            raise TypeError("index a Matrix with H[i, j] or H[r0:r1, c0:c1]")
        i, j = key
        if isinstance(i, numbers.Integral) and isinstance(j, numbers.Integral):
            return self._e[i][j]
        ri = range(len(self.rows))[i] if isinstance(i, slice) else [range(len(self.rows))[i]]
        cj = range(len(self.cols))[j] if isinstance(j, slice) else [range(len(self.cols))[j]]
        return Matrix([[self._e[a][b] for b in cj] for a in ri], self.rows.subset(list(ri)), self.cols.subset(list(cj)))

    def __setitem__(self, key: tuple[int, int], value: Element) -> None:
        if not (
            isinstance(key, tuple)
            and len(key) == 2
            and isinstance(key[0], numbers.Integral)
            and isinstance(key[1], numbers.Integral)
        ):
            raise TypeError("assign single entries: H[i, j] = x")
        if not isinstance(value, Element) or value.algebra is not self.algebra:
            raise TypeError("a Matrix entry must be an Element of the matrix's algebra")
        i, j = key
        self._e[i][j] = value

    def row(self, i: int) -> list[Element]:
        """Row ``i`` as a list of entries."""
        return list(self._e[i])

    def column(self, j: int) -> list[Element]:
        """Column ``j`` as a list of entries."""
        return [r[j] for r in self._e]

    def tolist(self) -> list[list[Element]]:
        """The entries as a list of lists (a copy of the lists; the entries are shared)."""
        return [list(r) for r in self._e]

    def copy(self) -> Matrix:
        """A copy (the entries, which are immutable in practice, are shared)."""
        return Matrix(self._e, self.rows, self.cols)

    def with_index(self, rows: Index, cols: Index | None = None) -> Matrix:
        """The same entries with another row (and column) index."""
        return Matrix(self._e, rows, rows if cols is None else cols)

    # ------------------------------------------------------------------ arithmetic
    def __matmul__(self, other: Matrix) -> Matrix:
        if not isinstance(other, Matrix):
            return NotImplemented
        if self.shape[1] != other.shape[0]:
            raise ValueError(f"shapes {self.shape} and {other.shape} do not compose")
        if self.algebra is not other.algebra:
            raise ValueError("matrices over different algebras")
        if self.cols != other.rows:
            raise ValueError(
                f"the column index {self.cols} of the left factor is not the row index {other.rows} of the "
                "right one (relabel with with_index if that is intended)"
            )
        alg = self.algebra
        out = []
        for i in range(self.shape[0]):
            row = []
            for j in range(other.shape[1]):
                acc = alg.zero()
                for k in range(self.shape[1]):
                    a, b = self._e[i][k], other._e[k][j]
                    if _exact_zero(a) or _exact_zero(b):
                        continue
                    acc = alg.add(acc, alg.mul(a, b))
                row.append(acc)
            out.append(row)
        return Matrix(out, self.rows, other.cols)

    def _check_same_shape(self, other: Matrix) -> None:
        if self.shape != other.shape:
            raise ValueError(f"shapes {self.shape} and {other.shape} differ")
        if self.algebra is not other.algebra:
            raise ValueError("matrices over different algebras")
        if self.rows != other.rows or self.cols != other.cols:
            raise ValueError("the matrices have different indices (relabel with with_index if that is intended)")

    def __add__(self, other: Matrix) -> Matrix:
        if not isinstance(other, Matrix):
            return NotImplemented
        self._check_same_shape(other)
        add = self.algebra.add
        return Matrix(
            [[add(a, b) for a, b in zip(ra, rb, strict=True)] for ra, rb in zip(self._e, other._e, strict=True)],
            self.rows,
            self.cols,
        )

    def __sub__(self, other: Matrix) -> Matrix:
        if not isinstance(other, Matrix):
            return NotImplemented
        self._check_same_shape(other)
        sub = self.algebra.sub
        return Matrix(
            [[sub(a, b) for a, b in zip(ra, rb, strict=True)] for ra, rb in zip(self._e, other._e, strict=True)],
            self.rows,
            self.cols,
        )

    def __neg__(self) -> Matrix:
        return (-1.0) * self

    def __mul__(self, c: Any) -> Matrix:
        if not isinstance(c, (int, float, complex, np.number)):
            return NotImplemented
        scale = self.algebra.scale
        return Matrix([[scale(c, a) for a in r] for r in self._e], self.rows, self.cols)

    __rmul__ = __mul__

    @property
    def H(self) -> Matrix:
        """The adjoint: ``(H*)_ij = (H_ji)*``."""
        star = self.algebra.star
        m, n = self.shape
        return Matrix([[star(self._e[j][i]) for j in range(m)] for i in range(n)], self.cols, self.rows)

    def conj(self) -> Matrix:
        """The entrywise conjugation."""
        conj = self.algebra.conj
        return Matrix([[conj(a) for a in r] for r in self._e], self.rows, self.cols)

    def real(self) -> Matrix:
        """The entrywise real part ``(x + c(x)) / 2``."""
        real = self.algebra.real
        return Matrix([[real(a) for a in r] for r in self._e], self.rows, self.cols)

    # ------------------------------------------------------------------ checks
    def corner_violation(self, rows: Index | None = None, cols: Index | None = None) -> float:
        """The total norm of the parts of the entries outside their corners ``E_i A E_j`` (``rows`` and
        ``cols`` override the matrix's own indices)."""
        rows, cols = rows or self.rows, cols or self.cols
        alg, bad = self.algebra, 0.0
        for i, r in enumerate(self._e):
            for j, x in enumerate(r):
                inside = alg.corner_project(x, rows.vertices(i), cols.vertices(j))
                bad += alg.dist(x, inside)
        return bad

    def unitary_defect(self, rows: Index | None = None, cols: Index | None = None) -> float:
        """``|U* U - diag(E_cols)| + |U U* - diag(E_rows)|``: zero exactly when ``U`` is a unitary from
        the corners of the column index to those of the row index."""
        rows, cols = rows or self.rows, cols or self.cols
        alg = self.algebra
        P1, P2 = self.H @ self, self @ self.H
        bad = 0.0
        for i in range(P1.shape[0]):
            for j in range(P1.shape[1]):
                bad += alg.dist(P1._e[i][j], cols.idempotent(i) if i == j else alg.zero())
        for i in range(P2.shape[0]):
            for j in range(P2.shape[1]):
                bad += alg.dist(P2._e[i][j], rows.idempotent(i) if i == j else alg.zero())
        return bad

    def off_band(self, band: int = 0) -> float:
        """The total norm of the entries with ``|i - j| > band``."""
        norm = self.algebra.norm
        return sum(norm(x) for i, r in enumerate(self._e) for j, x in enumerate(r) if abs(i - j) > band)

    def is_real(self) -> bool:
        """Whether every entry is real."""
        return all(self.algebra.is_real(x) for r in self._e for x in r)

    def is_selfadjoint(self, tol: float = 1e-9) -> bool:
        """Whether ``H* = H`` entrywise, to relative tolerance ``tol``."""
        alg = self.algebra
        m, n = self.shape
        return m == n and all(
            alg.dist(alg.star(self._e[j][i]), self._e[i][j]) <= tol * alg.norm(self._e[i][j])
            for i in range(m)
            for j in range(n)
        )

    def degree_range(self) -> tuple[int, int] | tuple[None, None]:
        """``(lowest, highest)`` degree among the entries (``(None, None)`` for the zero matrix)."""
        degs = [d for r in self._e for x in r for d in self.algebra.degrees(x)]
        return (min(degs), max(degs)) if degs else (None, None)

    def coefficients(
        self, letters: Sequence[str] = (), target: str | None = None, source: str | None = None
    ) -> np.ndarray:
        """The coefficient blocks of one word in every entry.

        Returns an array of shape ``(m, n)`` when the blocks are 1 x 1, else ``(m, n, d_t, d_s)``;
        see :meth:`Element.coefficient` for the arguments.
        """
        blocks = [[x.coefficient(letters, target, source) for x in r] for r in self._e]
        arr = np.array(blocks, dtype=complex)
        return arr[:, :, 0, 0] if arr.shape[2:] == (1, 1) else arr

    # ------------------------------------------------------------------ display
    def format(self, digits: int = 3, width: int = 40) -> str:
        """A grid of the formatted entries, each cut to ``width`` characters."""
        cells = [[self.algebra.format(x, digits) for x in r] for r in self._e]
        cells = [[c if len(c) <= width else c[: width - 1] + "…" for c in r] for r in cells]
        w = max((len(c) for r in cells for c in r), default=1)
        return "\n".join("[ " + "  ".join(c.ljust(w) for c in r) + " ]" for r in cells)

    def __repr__(self) -> str:
        m, n = self.shape
        return f"Matrix {m}x{n}, rows {self.rows}, cols {self.cols}\n{self.format()}"
\end{lstlisting}
\needspace{6\baselineskip}
\subsection{\texorpdfstring{\texttt{src/dsqlib/linalg.py}}{src/dsqlib/linalg.py}}
\begin{lstlisting}[style=dsqcode]
"""Matrix algorithms over the *-algebra of any doubly-signed quiver, written once.

Every function takes a :class:`~dsqlib.matrix.Matrix` and, where it needs them:

``sqrt``
    a *-square root ``x -> y`` with ``y* y = x``: a callable (for example ``algebra.sqrt`` or
    :func:`~dsqlib.sqrt.sqrt_rule`), a rule name (``'graded'``, ``'principal'``, ``'arbitrary'``), or
    ``None`` for the graded rule with principal choices;
``inverse``
    how pivots are inverted: ``'support'`` (the default: between the vertices an element touches),
    ``'orbits'`` (between STAR-BAR orbit idempotents, so that a zero divisor of an orbit corner is a
    breakdown), ``'all'``, or a callable. Every quotient ``x / y`` is checked: ``(x / y) y`` (or
    ``y (y^-1 x)``) must equal ``x``, which fails when ``y`` is invertible only on part of the corner of
    ``x``; that is reported as a breakdown. Householder reflectors do not use ``inverse``: their
    denominator ``w* x`` must be invertible in the corner ``f A f`` of its column (thesis, Proposition
    'householder-any-root'), and it is inverted there.

Nothing here depends on the quiver. :class:`~dsqlib.errors.BreakdownError` is raised when a pivot
has no inverse or does not divide exactly, at an isotropic vector (``v != 0`` with ``v* v = 0``) in
Gram-Schmidt, Householder and Lanczos, when the denominator of a Householder reflector is not invertible
in the corner of its column, at a zero pivot in :func:`cholesky` and :func:`ldl`, and at a zero pivot
with nonzero entries below it in :func:`lu`. With the default ``inverse='support'`` (which maps zero to
zero) a residual column that is exactly zero, i.e. a rank-deficient ``H``, gives a zero column of ``Q``
in Gram-Schmidt: ``Q`` is then a partial isometry, ``H = Q R`` still holds, and for a self-adjoint ``H``
the QR step ``Q* H Q`` keeps the nonzero spectrum and replaces the kernel by zeros.

Examples
--------
>>> import numpy as np
>>> from dsqlib import problems
>>> from dsqlib.linalg import householder
>>> A = np.array([[2.0, 1.0, 0.5], [1.0, 3.0, 1.0], [0.5, 1.0, 4.0]])
>>> T, U, _ = householder(problems.symmetric_problem(A))
>>> T.off_band(1) == 0  # tridiagonal
True
>>> B = problems.extract_loop(T, ["x"]).real
>>> bool(np.allclose(np.linalg.eigvalsh(B), np.linalg.eigvalsh(A)))
True
"""

from __future__ import annotations

import functools
from collections.abc import Callable, Sequence
from typing import NamedTuple

from .algebra import INF, Algebra, Element
from .errors import BreakdownError, NotInvertibleError
from .matrix import Matrix
from .sqrt import sqrt_rule

__all__ = [
    "HouseholderResult",
    "LDLResult",
    "LUResult",
    "QRResult",
    "QRStepResult",
    "cholesky",
    "householder",
    "householder_qr",
    "lanczos",
    "ldl",
    "lu",
    "qr",
    "qr_step",
]

SqrtArg = Callable[[Element], Element] | str | None
InverseArg = Callable[[Element], Element] | str | None


class HouseholderResult(NamedTuple):
    """``H = U T U*`` with ``T`` tridiagonal."""

    T: Matrix
    U: Matrix
    residuals: list[float]
    """The size of what is left below the subdiagonal of column ``k`` after step ``k``."""


class QRResult(NamedTuple):
    """``H = Q R`` with ``Q`` unitary and ``R`` upper triangular."""

    Q: Matrix
    R: Matrix


class QRStepResult(NamedTuple):
    """One QR step: ``H = Q R`` and ``H_next = R Q = Q* H Q``."""

    H_next: Matrix
    Q: Matrix
    R: Matrix


class LDLResult(NamedTuple):
    """``H = L D L*`` with ``L`` unit lower triangular and ``D`` diagonal."""

    L: Matrix
    D: Matrix


class LUResult(NamedTuple):
    """``H = L U`` with ``L`` unit lower triangular and ``U`` upper triangular."""

    L: Matrix
    U: Matrix


# ---------------------------------------------------------------------- helpers
def _sqrt_fn(sqrt: SqrtArg) -> Callable[[Element], Element]:
    if sqrt is None:
        return sqrt_rule()
    if isinstance(sqrt, str):
        return sqrt_rule(sqrt)
    return sqrt


def _inverse_fn(alg: Algebra, inverse: InverseArg, what: str) -> Callable[[Element], Element]:
    base: Callable[[Element], Element]
    if inverse is None or isinstance(inverse, str):
        within = inverse or "support"
        if within not in ("support", "orbits", "all"):
            raise ValueError(f"inverse must be 'support', 'orbits', 'all' or a callable, not {within!r}")
        base = functools.partial(alg.inverse, within=within)
    else:
        base = inverse

    def guarded(x: Element) -> Element:
        try:
            return base(x)
        except NotInvertibleError as e:
            raise BreakdownError(f"{what}: a pivot has no inverse ({e})") from e

    return guarded


def _exact_zero(x: Element) -> bool:
    return not x.terms and x.prec == INF


def _dot(alg: Algebra, x: Sequence[Element], y: Sequence[Element]) -> Element:
    """``sum_i x_i* y_i``."""
    acc = alg.zero()
    for a, b in zip(x, y, strict=True):
        if _exact_zero(a) or _exact_zero(b):
            continue
        acc = alg.add(acc, alg.mul(alg.star(a), b))
    return acc


def _square(H: Matrix, what: str) -> None:
    if H.shape[0] != H.shape[1]:
        raise ValueError(f"{what} needs a square matrix, not {H.shape[0]} x {H.shape[1]}")


def _check_sign(sign: str | None) -> None:
    if sign not in ("auto", None):
        raise ValueError(f"sign must be 'auto' or None, not {sign!r}")


def _in_corner(x: Element, rows: Sequence[str], cols: Sequence[str]) -> bool:
    return all(k[0] in rows and k[1] in cols for k in x.terms)


_QUOTIENT_RTOL = 1e-9


def _checked(alg: Algebra, x: Element, q: Element, y: Element, what: str, left: bool = False) -> Element:
    """Return the quotient ``q`` (``x y^-1``, or ``y^-1 x`` if ``left``) after checking that it is exact:
    ``q y = x`` (``y q = x``). This fails when ``y`` is invertible only on part of the corner of ``x``,
    as for a pivot ``(c, 0)`` of the split-complex numbers dividing an entry ``(p, q)`` with ``q != 0``."""
    back = alg.mul(y, q) if left else alg.mul(q, y)
    if alg.dist(back, x) > _QUOTIENT_RTOL * (alg.norm(x) + alg.norm(back)):
        raise BreakdownError(f"{what}: a pivot is invertible only on part of the corner of the entry it divides")
    return q


def _corner_inverse(alg: Algebra, s: Element, verts: Sequence[str], what: str) -> Element:
    """The inverse of ``s`` from ``f`` to ``f``, ``f`` the sum of the idempotents of ``verts`` (the
    corner of a column), as a Householder reflector requires."""
    try:
        return alg.inverse_between(s, verts, verts)
    except NotInvertibleError as e:
        raise BreakdownError(
            f"{what}: the denominator w* x is not invertible in the corner of its column ({e}); this root of "
            "x* x gives no reflector"
        ) from e


# ---------------------------------------------------------------------- Householder
def householder(
    H: Matrix, sqrt: SqrtArg = None, inverse: InverseArg = None, sign: str | None = "auto"
) -> HouseholderResult:
    """Tridiagonalise a self-adjoint ``H`` by Householder reflectors: ``H = U T U*``.

    The reflector for the column part ``x`` is ``P = I - w (w* x)^+ w*`` with ``w = x - e_1 gamma`` and
    ``gamma = +-sqrt(x* x)``; ``P`` is unitary for *any* *-square root ``gamma``, so no phase condition is
    needed.

    Parameters
    ----------
    H : Matrix
        Self-adjoint.
    sqrt, inverse
        See the module documentation.
    sign : {'auto', None}
        ``'auto'`` uses ``-gamma`` (also a *-square root) when that avoids cancellation in
        ``x_0 - gamma``, as the classical algorithm does; ``None`` uses ``gamma`` as returned (the same
        in exact arithmetic, but a nearly reduced column then loses digits in floating point).

    Returns
    -------
    HouseholderResult
        ``(T, U, residuals)``.
    """
    _square(H, "householder")
    _check_sign(sign)
    alg, n = H.algebra, H.shape[0]
    sq = _sqrt_fn(sqrt)
    _inverse_fn(alg, inverse, "householder")  # validates the argument; the reflector uses its own inverse
    rows = H.rows
    T = H.tolist()
    U = Matrix.identity(rows)
    resid = []
    for k in range(n - 2):
        x = [T[i][k] for i in range(k + 1, n)]
        g = _dot(alg, x, x)
        if not g.terms and any(xi.terms for xi in x):
            raise BreakdownError(f"householder: column {k} is isotropic (x != 0, x* x = 0)")
        gamma = sq(g)
        if not _in_corner(gamma, rows.vertices(k + 1), H.cols.vertices(k)):
            raise ValueError(
                f"householder: the *-square root of x* x for column {k} does not lie in the corner of position "
                f"({k + 1}, {k}); Householder reflectors do not suit this indexing"
            )
        if sign == "auto" and alg.norm(alg.add(x[0], gamma)) > alg.norm(alg.sub(x[0], gamma)):
            gamma = alg.scale(-1.0, gamma)
        w = [alg.sub(x[0], gamma), *x[1:]]
        if all(not wi.terms for wi in w):  # x = e_1 gamma already (to the known precision)
            resid.append(sum(alg.norm(T[i][k]) for i in range(k + 2, n)))
            continue
        si = _corner_inverse(alg, _dot(alg, w, x), H.cols.vertices(k), f"householder: column {k}")
        P = Matrix.identity(rows).tolist()
        for i in range(k + 1, n):
            for j in range(k + 1, n):
                P[i][j] = alg.sub(P[i][j], alg.mul(alg.mul(w[i - k - 1], si), alg.star(w[j - k - 1])))
        Pm = Matrix(P, rows)
        Ps = Pm.H
        T = (Pm @ Matrix(T, rows) @ Ps).tolist()
        U = U @ Ps
        resid.append(sum(alg.norm(T[i][k]) for i in range(k + 2, n)))
    return HouseholderResult(Matrix(T, rows), U, resid)


def householder_qr(H: Matrix, sqrt: SqrtArg = None, inverse: InverseArg = None, sign: str | None = "auto") -> QRResult:
    """``H = Q R`` by Householder reflectors (the reflectors of :func:`householder`).

    The reflector for column ``k`` puts the *-square root ``gamma`` of ``x* x`` at position ``(k, k)``,
    so it needs ``gamma`` in the corner ``E_k A E_k``: true on the plain indexing and on one-orbit
    problems, but not on the SVD orbit indexing (where the root of a degree-2 Gram element maps one
    orbit to the other); there it raises ``ValueError`` and Gram-Schmidt (:func:`qr`) is the method.

    Returns
    -------
    QRResult
        The full factorisation: ``Q`` (``m x m``) unitary, ``Q* Q = diag(E_i)``, and ``R`` (``m x n``)
        upper triangular.

    Raises
    ------
    BreakdownError
        At an isotropic column or a pivot with no inverse.
    """
    _check_sign(sign)
    alg, (m, n) = H.algebra, H.shape
    sq = _sqrt_fn(sqrt)
    _inverse_fn(alg, inverse, "householder_qr")  # validates the argument; the reflector uses its own inverse
    R = H.tolist()
    Q = Matrix.identity(H.rows).tolist()
    for k in range(min(m - 1, n)):
        x = [R[i][k] for i in range(k, m)]
        g = _dot(alg, x, x)
        if not g.terms:
            if any(xi.terms for xi in x):
                raise BreakdownError(f"householder_qr: column {k} is isotropic (x != 0, x* x = 0)")
            continue
        gamma = sq(g)
        if not _in_corner(gamma, H.rows.vertices(k), H.cols.vertices(k)):
            raise ValueError(
                f"householder_qr: the *-square root of x* x for column {k} does not lie in the corner of "
                f"position ({k}, {k}); use qr(H, method='gram-schmidt') on this indexing"
            )
        if sign == "auto" and alg.norm(alg.add(x[0], gamma)) > alg.norm(alg.sub(x[0], gamma)):
            gamma = alg.scale(-1.0, gamma)
        w = [alg.sub(x[0], gamma), *x[1:]]
        if all(not wi.terms for wi in w):
            continue
        si = _corner_inverse(alg, _dot(alg, w, x), H.cols.vertices(k), f"householder_qr: column {k}")
        for j in range(n):  # rows k..m-1 of R:  R <- P R
            c = _dot(alg, w, [R[k + i][j] for i in range(m - k)])
            if _exact_zero(c):
                continue
            f = alg.mul(si, c)
            for i in range(m - k):
                R[k + i][j] = alg.sub(R[k + i][j], alg.mul(w[i], f))
        sis = alg.star(si)  # Q <- Q P*,  P* = I - w si* w*
        for r in range(m):
            c = alg.zero()
            for i in range(m - k):
                if not _exact_zero(Q[r][k + i]) and not _exact_zero(w[i]):
                    c = alg.add(c, alg.mul(Q[r][k + i], w[i]))
            if _exact_zero(c):
                continue
            f = alg.mul(c, sis)
            for i in range(m - k):
                Q[r][k + i] = alg.sub(Q[r][k + i], alg.mul(f, alg.star(w[i])))
    return QRResult(Matrix(Q, H.rows), Matrix(R, H.rows, H.cols))


def _gram_schmidt(H: Matrix, sq: Callable[[Element], Element], inv: Callable[[Element], Element]) -> Matrix:
    alg, (m, n) = H.algebra, H.shape
    cols: list[list[Element]] = []
    for j in range(n):
        v = H.column(j)
        for q in cols:
            c = _dot(alg, q, v)
            v = [alg.sub(v[i], alg.mul(q[i], c)) for i in range(m)]
        g = _dot(alg, v, v)
        if not g.terms and any(x.terms for x in v):
            raise BreakdownError(f"Gram-Schmidt: column {j} is isotropic (v != 0, v* v = 0)")
        r = sq(g)
        yi = inv(r)  # a zero residual (H rank-deficient) gives a zero column of Q
        cols.append([_checked(alg, v[i], alg.mul(v[i], yi), r, "Gram-Schmidt") for i in range(m)])
    return Matrix.from_columns(cols, H.rows, H.cols)


def qr(H: Matrix, sqrt: SqrtArg = None, inverse: InverseArg = None, method: str = "gram-schmidt") -> QRResult:
    """``H = Q R`` with ``R`` upper triangular.

    Parameters
    ----------
    H : Matrix
        ``m x n``.
    sqrt, inverse
        See the module documentation.
    method : {'gram-schmidt', 'householder'}
        Gram-Schmidt, with ``R_jj`` = the *-square root of the Gram element of the residual column and
        ``R = Q* H``: the thin factorisation, ``Q`` (``m x n``) with ``Q* Q = diag(E_j)`` and ``R``
        (``n x n``); it works on every indexing. Or Householder reflectors (:func:`householder_qr`): the
        full factorisation, ``Q`` (``m x m``) unitary and ``R`` (``m x n``), on indexings where the
        square roots stay in their corners.

    Returns
    -------
    QRResult

    Raises
    ------
    BreakdownError
        See the module documentation.
    """
    if method == "householder":
        return householder_qr(H, sqrt, inverse)
    if method != "gram-schmidt":
        raise ValueError(f"method must be 'householder' or 'gram-schmidt', not {method!r}")
    Q = _gram_schmidt(H, _sqrt_fn(sqrt), _inverse_fn(H.algebra, inverse, "qr"))
    return QRResult(Q, Q.H @ H)


def qr_step(H: Matrix, sqrt: SqrtArg = None, inverse: InverseArg = None, method: str = "gram-schmidt") -> QRStepResult:
    """One step of the QR algorithm: ``H = Q R`` and ``H_next = Q* H Q`` (``= R Q``).

    Parameters
    ----------
    H : Matrix
        Square.
    sqrt, inverse
        See the module documentation.
    method : {'gram-schmidt', 'householder'}
        How the factorisation is computed (see :func:`qr`).

    Returns
    -------
    QRStepResult
        ``(H_next, Q, R)``.

    Notes
    -----
    On an indexing where the square roots change corners (the SVD orbit indexing, where the root of
    a degree-2 Gram element has degree 1), one step maps each orbit to the other end of the arrows;
    use :meth:`Index.relabel` to describe the result, or take double steps.
    """
    _square(H, "qr_step")
    if method == "gram-schmidt":
        Q = _gram_schmidt(H, _sqrt_fn(sqrt), _inverse_fn(H.algebra, inverse, "qr_step"))
        R = Q.H @ H
        return QRStepResult(R @ Q, Q, R)
    Q, R = qr(H, sqrt, inverse, method)
    return QRStepResult(Q.H @ H @ Q, Q, R)


# ---------------------------------------------------------------------- Lanczos
def lanczos(
    H: Matrix,
    q1: Matrix | Sequence[Element],
    steps: int | None = None,
    sqrt: SqrtArg = None,
    inverse: InverseArg = None,
) -> Matrix:
    """Lanczos vectors ``q_1, ..., q_steps`` of a self-adjoint ``H``, with full reorthogonalisation (twice).

    Parameters
    ----------
    H : Matrix
        Self-adjoint.
    q1 : Matrix or sequence of Element
        The starting vector, with ``q1* q1`` an idempotent: an ``n x 1`` matrix such as
        ``Matrix.unit_vector(H.rows, 0)``, or a list of ``n`` entries.
    steps : int, optional
        The number of vectors, at least 1 (default ``n``).

    Returns
    -------
    Matrix
        The ``n x k`` matrix of Lanczos vectors (``k <= steps``: it stops early when the residual
        vanishes, i.e. at an invariant subspace). Its column index is the first ``k`` positions of
        ``H.rows``, which describes the corners of the vectors when ``q1`` lies in the corner of
        position 0; otherwise relabel it with :meth:`Matrix.with_index`.
    """
    _square(H, "lanczos")
    alg, n = H.algebra, H.shape[0]
    steps = n if steps is None else steps
    if steps < 1:
        raise ValueError(f"steps must be at least 1, not {steps}")
    if isinstance(q1, Matrix) and q1.shape != (n, 1):
        raise ValueError(f"q1 must be an {n} x 1 matrix, not {q1.shape[0]} x {q1.shape[1]}")
    sq, inv = _sqrt_fn(sqrt), _inverse_fn(alg, inverse, "lanczos")
    q = q1.column(0) if isinstance(q1, Matrix) else list(q1)
    if len(q) != n:
        raise ValueError(f"q1 must have {n} entries, not {len(q)}")
    Q = [q]
    Ht = H.tolist()
    for _ in range(steps - 1):
        r = []
        for i in range(n):
            acc = alg.zero()
            for k in range(n):
                if _exact_zero(Ht[i][k]) or _exact_zero(Q[-1][k]):
                    continue
                acc = alg.add(acc, alg.mul(Ht[i][k], Q[-1][k]))
            r.append(acc)
        for _ in range(2):
            for qv in Q:
                c = _dot(alg, qv, r)
                r = [alg.sub(r[i], alg.mul(qv[i], c)) for i in range(n)]
        if not any(x.terms for x in r):
            break
        g = _dot(alg, r, r)
        if not g.terms:
            raise BreakdownError("lanczos: isotropic residual (r != 0, r* r = 0)")
        beta = sq(g)
        bi = inv(beta)
        Q.append([_checked(alg, r[i], alg.mul(r[i], bi), beta, "lanczos") for i in range(n)])
    return Matrix.from_columns(Q, H.rows)


# ---------------------------------------------------------------------- triangular factorisations
def cholesky(H: Matrix, sqrt: SqrtArg = None, inverse: InverseArg = None) -> Matrix:
    """``H = R* R`` with ``R`` upper triangular and ``R_jj`` the *-square root of the pivot.

    Raises
    ------
    BreakdownError
        At a zero or non-invertible pivot.
    NoSquareRootError
        When a pivot has no *-square root by the rule used (e.g. a negative pivot at a real vertex).
    """
    _square(H, "cholesky")
    alg, n = H.algebra, H.shape[0]
    sq, inv = _sqrt_fn(sqrt), _inverse_fn(alg, inverse, "cholesky")
    R = Matrix.zeros(H.rows).tolist()
    for j in range(n):
        s = H[j, j]
        for k in range(j):
            s = alg.sub(s, alg.mul(alg.star(R[k][j]), R[k][j]))
        if not s.terms:
            raise BreakdownError(f"cholesky: zero pivot at {j}")
        R[j][j] = sq(s)
        ris = alg.star(inv(R[j][j]))
        for col in range(j + 1, n):
            s = H[j, col]
            for k in range(j):
                s = alg.sub(s, alg.mul(alg.star(R[k][j]), R[k][col]))
            R[j][col] = _checked(alg, s, alg.mul(ris, s), alg.star(R[j][j]), "cholesky", left=True)
    return Matrix(R, H.rows)


def ldl(H: Matrix, inverse: InverseArg = None) -> LDLResult:
    """``H = L D L*`` with ``L`` unit lower triangular (diagonal ``E_i``) and ``D`` diagonal; no square
    roots.

    Raises
    ------
    BreakdownError
        At a zero or non-invertible pivot.
    """
    _square(H, "ldl")
    alg, n = H.algebra, H.shape[0]
    inv = _inverse_fn(alg, inverse, "ldl")
    L = Matrix.identity(H.rows).tolist()
    D = [alg.zero() for _ in range(n)]
    for j in range(n):
        s = H[j, j]
        for k in range(j):
            s = alg.sub(s, alg.mul(alg.mul(L[j][k], D[k]), alg.star(L[j][k])))
        if not s.terms:
            raise BreakdownError(f"ldl: zero pivot at {j}")
        D[j] = s
        Di = inv(s)
        for i in range(j + 1, n):
            s = H[i, j]
            for k in range(j):
                s = alg.sub(s, alg.mul(alg.mul(L[i][k], D[k]), alg.star(L[j][k])))
            L[i][j] = _checked(alg, s, alg.mul(s, Di), D[j], "ldl")
    Dm = Matrix.zeros(H.rows)
    for j in range(n):
        Dm[j, j] = D[j]
    return LDLResult(Matrix(L, H.rows), Dm)


def lu(H: Matrix, inverse: InverseArg = None) -> LUResult:
    """``H = L U`` by Gaussian elimination without pivoting: ``L`` unit lower triangular (diagonal
    ``E_i``), ``U`` upper triangular.

    Raises
    ------
    BreakdownError
        At a zero pivot with nonzero entries below it, or a pivot with no inverse.
    """
    _square(H, "lu")
    alg, n = H.algebra, H.shape[0]
    inv = _inverse_fn(alg, inverse, "lu")
    U = H.tolist()
    L = Matrix.identity(H.rows).tolist()
    for k in range(n):
        if not U[k][k].terms:
            if any(U[i][k].terms for i in range(k + 1, n)):
                raise BreakdownError(f"lu: zero pivot at {k}")
            continue
        pi = inv(U[k][k])
        for i in range(k + 1, n):
            if not U[i][k].terms:
                continue
            m = _checked(alg, U[i][k], alg.mul(U[i][k], pi), U[k][k], "lu")
            L[i][k] = m
            for j in range(k, n):
                U[i][j] = alg.sub(U[i][j], alg.mul(m, U[k][j]))
    return LUResult(Matrix(L, H.rows), Matrix(U, H.rows))
\end{lstlisting}
\needspace{6\baselineskip}
\subsection{\texorpdfstring{\texttt{src/dsqlib/catalogue.py}}{src/dsqlib/catalogue.py}}
\begin{lstlisting}[style=dsqcode]
"""The doubly-signed quivers of the thesis catalogue (Chapter 6) and the double of ``A_r``
(Example 5.20)."""

from __future__ import annotations

from collections.abc import Callable

from .quiver import DSQ

__all__ = [
    "CATALOGUE",
    "double_A",
    "hermitian",
    "ldl",
    "similarity",
    "skew_loop",
    "split_complex",
    "svd",
    "symmetric_loop",
    "takagi",
    "youla",
]


def symmetric_loop() -> DSQ:
    """Symmetric eigendecomposition: one (R, +1) vertex ``v`` with a loop ``x``, ``x* = x``.

    Completed, its algebra is the field of real Laurent series ``R((t))`` with ``t = x``."""
    return DSQ.build([("v", "R")], [("x", "v", "v", {"sign": +1})])


def skew_loop() -> DSQ:
    """Normal form of skew-symmetric matrices: one (R, +1) vertex ``v`` with a loop ``h``, ``h* = -h``."""
    return DSQ.build([("v", "R")], [("h", "v", "v", {"sign": -1})])


def svd(field: str = "R") -> DSQ:
    """Singular value decomposition over ``field`` = 'R', 'C' or 'H': vertices ``u``, ``v`` and an arrow
    ``a: u -> v`` (with its transpose ``a'``)."""
    kinds = {"R": "R", "C": "Cu", "H": "Hsp"}
    if field not in kinds:
        raise ValueError(f"field must be 'R', 'C' or 'H', not {field!r}")
    return DSQ.build([("u", kinds[field]), ("v", kinds[field])], [("a", "u", "v")])


def hermitian(field: str = "C") -> DSQ:
    """Hermitian eigendecomposition over C (a loop ``x`` at ``v`` with ``x^T = x~``) or over H (a loop
    with ``x* = x``)."""
    if field == "C":
        return DSQ.build([("v", "Cu")], [("x", "v", "v", {"top_is_bar": True})])
    if field == "H":
        return DSQ.build([("v", "Hsp")], [("x", "v", "v", {"sign": +1})])
    raise ValueError(f"field must be 'C' or 'H', not {field!r}")


def takagi() -> DSQ:
    """Takagi factorisation: a unitary-type complex vertex ``v`` and an arrow ``al: v -> v~ = v^T`` with
    sign +1."""
    return DSQ.build([("v", "Cu")], [("al", "v", "v~", {"sign": +1})])


def youla() -> DSQ:
    """Youla normal form of complex skew-symmetric matrices: ``al: v -> v~ = v^T`` with sign -1."""
    return DSQ.build([("v", "Cu")], [("al", "v", "v~", {"sign": -1})])


def ldl() -> DSQ:
    """Congruence of real symmetric matrices (``LDL^T``): an exchanged real pair ``u, u'`` and
    ``a: u -> u'`` with sign +1."""
    return DSQ.build([("u", "RR")], [("a", "u", "u'", {"sign": +1})])


def similarity() -> DSQ:
    """Similarity of real matrices: an exchanged real pair ``v, v'`` with a loop ``x`` at ``v`` (and
    ``x'`` at ``v'``)."""
    return DSQ.build([("v", "RR")], [("x", "v", "v")])


def split_complex() -> DSQ:
    """The split-complex numbers ``2R = R x R`` with ``(a, b)* = (b, a)``: one exchanged real pair,
    no arrows."""
    return DSQ.build([("u", "RR")], [])


def double_A(r: int) -> DSQ:
    """The double of ``1 -> 2 -> ... -> r`` with all vertices of kind (R, +1) (Example 5.20)."""
    return DSQ.build([(str(j), "R") for j in range(1, r + 1)], [(f"a{j}", str(j), str(j + 1)) for j in range(1, r)])


CATALOGUE: dict[str, Callable[[], DSQ]] = {
    "symmetric eigendecomposition": symmetric_loop,
    "skew-symmetric normal form": skew_loop,
    "real SVD": lambda: svd("R"),
    "complex SVD": lambda: svd("C"),
    "quaternionic SVD": lambda: svd("H"),
    "Hermitian eigendecomposition over C": lambda: hermitian("C"),
    "Hermitian eigendecomposition over H": lambda: hermitian("H"),
    "Takagi factorisation": takagi,
    "Youla normal form": youla,
    "LDL^T (congruence)": ldl,
    "similarity of real matrices": similarity,
    "split-complex numbers": split_complex,
}
"""The catalogue quivers by problem name (each value builds the quiver)."""
\end{lstlisting}
\needspace{6\baselineskip}
\subsection{\texorpdfstring{\texttt{src/dsqlib/problems.py}}{src/dsqlib/problems.py}}
\begin{lstlisting}[style=dsqcode]
"""Classical matrix problems as matrices over the *-algebras of the catalogue.

Each ``*_problem`` function returns a self-adjoint :class:`~dsqlib.matrix.Matrix` ``H`` that encodes
the classical data (its algebra is ``H.algebra`` and its index ``H.rows``); each ``extract_*``
function reads the classical matrix back from a (transformed) ``H``. Quaternions are 2 x 2 complex
blocks ``[[p, q], [-conj(q), conj(p)]]`` (see :func:`~dsqlib.algebra.quaternion`); a quaternion matrix
is an array of shape ``(m, n, 2, 2)``.

==========================  ==========================================  ===========================
problem                     classical data                              read back with
==========================  ==========================================  ===========================
:func:`svd_problem`         ``W`` over R, C or H (orbit indexing)       :func:`extract_svd`
:func:`svd_problem_plain`   ``W`` over R, C or H (plain indexing)       :func:`extract_svd_plain`
:func:`hermitian_problem`   Hermitian ``A`` over C or H                 :func:`extract_loop`
:func:`symmetric_problem`   real symmetric ``A``                        :func:`extract_loop`
:func:`skew_problem`        real skew-symmetric ``A``                   :func:`extract_loop`
:func:`takagi_problem`      complex symmetric ``A``                     :func:`extract_loop`
:func:`youla_problem`       complex skew-symmetric ``A``                :func:`extract_loop`
:func:`split_problem`       real ``X`` over the split-complex numbers   :func:`extract_split`
:func:`similarity_problem`  real ``X`` on the similarity quiver         :func:`extract_loop`
:func:`congruence_problem`  real symmetric ``S`` (congruence, LDL)      ``D[i, i].coefficient``
:func:`series_problem`      real ``A`` over ``R((t))``                  :func:`standard_part`
==========================  ==========================================  ===========================
"""

from __future__ import annotations

from collections.abc import Sequence

import numpy as np

from . import catalogue
from .algebra import Algebra, Element, quaternion
from .matrix import Matrix, orbit_index, plain_index
from .quiver import DSQ

__all__ = [
    "complex_representation",
    "congruence_problem",
    "extract_loop",
    "extract_split",
    "extract_svd",
    "extract_svd_plain",
    "hermitian_problem",
    "quaternion",
    "random_quaternion_matrix",
    "series_problem",
    "similarity_problem",
    "skew_problem",
    "split_problem",
    "standard_part",
    "svd_problem",
    "svd_problem_plain",
    "symmetric_problem",
    "t_power",
    "takagi_problem",
    "youla_problem",
]


def _block(z: object, field: str) -> np.ndarray:
    return np.asarray(z, dtype=complex).reshape(2, 2) if field == "H" else np.array([[z]], dtype=complex)


def _scalar(X: np.ndarray, field: str) -> np.ndarray | complex:
    return X if field == "H" else X[0, 0]


def _check_field(field: str, allowed: str = "RCH") -> None:
    if field not in allowed:
        raise ValueError(f"field must be one of {list(allowed)}, not {field!r}")


def _check_square(A: np.ndarray) -> None:
    if A.ndim < 2 or A.shape[0] != A.shape[1]:
        raise ValueError(f"expected a square matrix, not shape {A.shape}")


def complex_representation(W: np.ndarray) -> np.ndarray:
    """The complex ``2m x 2n`` representation of a quaternion matrix given as ``(m, n, 2, 2)`` blocks."""
    m, n = W.shape[:2]
    return W.transpose(0, 2, 1, 3).reshape(2 * m, 2 * n)


def random_quaternion_matrix(rng: np.random.Generator, m: int, n: int) -> np.ndarray:
    """A random ``m x n`` quaternion matrix (normal coefficients), as ``(m, n, 2, 2)`` blocks."""
    W = np.zeros((m, n, 2, 2), dtype=complex)
    for i in range(m):
        for j in range(n):
            W[i, j] = quaternion(*rng.normal(size=4))
    return W


# ---------------------------------------------------------------- SVD over R, C, H
def svd_problem(W: np.ndarray, field: str = "R", order: Sequence[str] = ("u", "v"), relprec: int = 16) -> Matrix:
    """The SVD of ``W`` (``m x n``; quaternions as ``(m, n, 2, 2)`` blocks) as a self-adjoint matrix on
    the orbit-cyclic indexing of the SVD quiver ``u -> v``: ``H[(v, i), (u, j)]`` is the real element
    with coefficient ``W_ij`` at the arrow ``a``, and ``H = H*``.

    Parameters
    ----------
    W : numpy.ndarray
    field : {'R', 'C', 'H'}
    order : sequence of str
        The cyclic order of the two orbits; ``('u', 'v')`` gives positions ``u0, v0, u1, v1, ...``.
    relprec : int
        Passed to the :class:`~dsqlib.algebra.Algebra`.
    """
    _check_field(field)
    alg = Algebra(catalogue.svd(field), relprec=relprec)
    m, n = W.shape[:2]
    index = orbit_index(alg, {"u": n, "v": m}, order)
    H = Matrix.zeros(index)
    for I, (oi, ki) in enumerate(index):
        for J, (oj, kj) in enumerate(index):
            if oi == "v" and oj == "u":
                H[I, J] = alg.real_element(alg.arrow("a", _block(W[ki, kj], field)))
                H[J, I] = alg.star(H[I, J])
    return H


def extract_svd(H: Matrix, field: str = "R") -> np.ndarray:
    """The classical matrix of an SVD problem: the coefficients of ``a`` in the ``(v, u)`` entries."""
    _check_field(field)
    pos = H.rows.positions
    m = sum(1 for o, _ in pos if o == "v")
    n = sum(1 for o, _ in pos if o == "u")
    W = np.zeros((m, n, 2, 2), dtype=complex) if field == "H" else np.zeros((m, n), dtype=complex)
    for I, (oi, ki) in enumerate(pos):
        for J, (oj, kj) in enumerate(pos):
            if oi == "v" and oj == "u":
                W[ki, kj] = _scalar(H[I, J].coefficient(["a"], "v", "u"), field)
    return W


def svd_problem_plain(W: np.ndarray, field: str = "R", relprec: int = 16) -> Matrix:
    """The SVD problem on the plain indexing (every ``E_i = 1``; ``N x N`` with ``N = max(m, n)``,
    padded by zeros): ``H_ij = W_ij a + conj(W_ji) a'`` (and conjugates), i.e. ``H = W (x) a + W^H (x) a*``."""
    _check_field(field)
    alg = Algebra(catalogue.svd(field), relprec=relprec)
    m, n = W.shape[:2]
    index = plain_index(alg, max(m, n))
    H = Matrix.zeros(index)
    for i in range(m):
        for j in range(n):
            x = alg.real_element(alg.arrow("a", _block(W[i, j], field)))
            H[i, j] = alg.add(H[i, j], x)
            H[j, i] = alg.add(H[j, i], alg.star(x))
    return H


def extract_svd_plain(H: Matrix, m: int, n: int, field: str = "R") -> np.ndarray:
    """The classical ``m x n`` matrix of :func:`svd_problem_plain`."""
    _check_field(field)
    W = np.zeros((m, n, 2, 2), dtype=complex) if field == "H" else np.zeros((m, n), dtype=complex)
    for i in range(m):
        for j in range(n):
            W[i, j] = _scalar(H[i, j].coefficient(["a"], "v", "u"), field)
    return W


# ---------------------------------------------------------------- one-orbit problems
def _one_orbit(quiver: DSQ, A: np.ndarray, word: Sequence[str], field: str, relprec: int) -> Matrix:
    _check_square(A)
    alg = Algebra(quiver, relprec=relprec)
    n = A.shape[0]
    index = orbit_index(alg, {"v": n})
    H = Matrix.zeros(index)
    for i in range(n):
        for j in range(i, n):
            x = alg.real_element(alg.word(word, _block(A[i, j], field)))
            if i == j:
                x = alg.scale(0.5, alg.add(x, alg.star(x)))
            H[i, j] = x
            H[j, i] = alg.star(x)
    return H


def hermitian_problem(A: np.ndarray, field: str = "C", relprec: int = 16) -> Matrix:
    """A Hermitian matrix over C (or H, as ``(n, n, 2, 2)`` blocks): ``H_ij`` has coefficient ``A_ij`` at
    the loop ``x``."""
    _check_field(field, "CH")
    return _one_orbit(catalogue.hermitian(field), A, ["x"], field, relprec)


def symmetric_problem(A: np.ndarray, relprec: int = 16) -> Matrix:
    """A real symmetric matrix: ``H_ij = A_ij x`` with the loop ``x* = x``."""
    return _one_orbit(catalogue.symmetric_loop(), A, ["x"], "R", relprec)


def skew_problem(A: np.ndarray, relprec: int = 16) -> Matrix:
    """A real skew-symmetric matrix: ``H_ij = A_ij h`` with the loop ``h* = -h``."""
    return _one_orbit(catalogue.skew_loop(), A, ["h"], "R", relprec)


def takagi_problem(A: np.ndarray, relprec: int = 16) -> Matrix:
    """A complex symmetric matrix: ``H_ij`` has coefficient ``A_ij`` at the arrow ``al: v -> v~``."""
    return _one_orbit(catalogue.takagi(), A, ["al"], "C", relprec)


def youla_problem(A: np.ndarray, relprec: int = 16) -> Matrix:
    """A complex skew-symmetric matrix: ``H_ij`` has coefficient ``A_ij`` at ``al: v -> v~`` (sign -1)."""
    return _one_orbit(catalogue.youla(), A, ["al"], "C", relprec)


def extract_loop(
    H: Matrix, word: Sequence[str], target: str | None = None, source: str | None = None, field: str = "C"
) -> np.ndarray:
    """The classical matrix of a one-orbit problem: the coefficients of ``word`` in every entry
    (``['x']`` for the Hermitian and symmetric problems, ``['h']`` for the skew one, ``['al']`` for Takagi
    and Youla). ``target`` and ``source`` default to the ends of the word."""
    _check_field(field)
    n = H.shape[0]
    A = np.zeros((n, n, 2, 2), dtype=complex) if field == "H" else np.zeros((n, n), dtype=complex)
    for i in range(n):
        for j in range(n):
            A[i, j] = _scalar(H[i, j].coefficient(word, target, source), field)
    return A


# ---------------------------------------------------------------- split-complex numbers and congruence
def split_problem(X: np.ndarray, relprec: int = 16) -> Matrix:
    """A real matrix over the split-complex numbers ``2R``: ``H_ij = X_ij e_u + X_ji e_u'``, so that
    ``H = (X, X^T)`` is self-adjoint."""
    _check_square(X)
    alg = Algebra(catalogue.split_complex(), relprec=relprec)
    index = orbit_index(alg, {"u": X.shape[0]})
    return Matrix.from_function(index, lambda i, j: alg.add(alg.vertex("u", X[i, j]), alg.vertex("u'", X[j, i])))


def extract_split(H: Matrix) -> np.ndarray:
    """The first component ``X`` of a matrix ``(X, Y)`` over the split-complex numbers."""
    n = H.shape[0]
    return np.array([[H[i, j].coefficient((), "u", "u")[0, 0] for j in range(n)] for i in range(n)])


def similarity_problem(X: np.ndarray, relprec: int = 16) -> Matrix:
    """Similarity of a real square ``X`` on the similarity quiver: ``H_ij = X_ij x + X_ji x'``, with the
    loop ``x`` at ``v`` and its transpose ``x'`` at ``v'``, so that ``H`` is self-adjoint. Read ``X`` back
    with ``extract_loop(H, ['x'], 'v', 'v', 'R')``."""
    _check_square(X)
    alg = Algebra(catalogue.similarity(), relprec=relprec)
    index = orbit_index(alg, {"v": X.shape[0]})
    return Matrix.from_function(index, lambda i, j: alg.add(alg.arrow("x", X[i, j]), alg.arrow("x'", X[j, i])))


def congruence_problem(S: np.ndarray, relprec: int = 16) -> Matrix:
    """Congruence of a real symmetric ``S``: ``H_ij = S_ij a`` with ``a: u -> u'`` (the LDL^T quiver)."""
    _check_square(S)
    alg = Algebra(catalogue.ldl(), relprec=relprec)
    index = orbit_index(alg, {"u": S.shape[0]})
    return Matrix.from_function(index, lambda i, j: alg.arrow("a", S[i, j]))


# ---------------------------------------------------------------- real Laurent series R((t))
def series_problem(A: np.ndarray, relprec: int = 24) -> Matrix:
    """A real matrix over the field of real Laurent series ``R((t))``, realised as the completed algebra
    of the one-loop quiver :func:`~dsqlib.catalogue.symmetric_loop` (``t`` = the loop ``x``, with its
    local inverse ``x^-1 = t^-1``); plain indexing."""
    _check_square(A)
    alg = Algebra(catalogue.symmetric_loop(), relprec=relprec)
    index = plain_index(alg, A.shape[0])
    return Matrix.from_function(index, lambda i, j: alg.vertex("v", float(A[i, j])) if A[i, j] != 0 else None)


def t_power(alg: Algebra, k: int) -> Element:
    """``t^k`` (``k`` in Z) in ``R((t))`` = the completed one-loop algebra."""
    k = int(k)
    if k == 0:
        return alg.idempotent("v")
    return alg.word(["x"] * k) if k > 0 else alg.word(["x^-1"] * (-k))


def standard_part(H: Matrix) -> np.ndarray:
    """The ``t^0`` coefficients of a matrix over ``R((t))``, as a float array (NaN where an entry has a
    pole or is not known at degree 0)."""
    m, n = H.shape
    S = np.zeros((m, n))
    for i in range(m):
        for j in range(n):
            x = H[i, j]
            if (x.terms and x.valuation < 0) or x.prec <= 0:
                S[i, j] = np.nan
            else:
                S[i, j] = x.coefficient((), "v", "v")[0, 0].real
    return S
\end{lstlisting}
\needspace{6\baselineskip}
\subsection{\texorpdfstring{\texttt{src/dsqlib/errors.py}}{src/dsqlib/errors.py}}
\begin{lstlisting}[style=dsqcode]
"""Exceptions raised by dsqlib.

Every exception derives from :class:`DSQError` and also from the built-in exception it refines, so
code that catches ``ValueError`` or ``ZeroDivisionError`` keeps working.
"""

from __future__ import annotations

__all__ = [
    "BreakdownError",
    "DSQError",
    "NoSquareRootError",
    "NotInvertibleError",
    "PrecisionError",
    "QuiverError",
]


class DSQError(Exception):
    """Base class of all dsqlib exceptions."""


class QuiverError(DSQError, ValueError):
    """The data passed to :meth:`DSQ.build` or :class:`DSQ` do not define a doubly-signed quiver.

    (Invalid arguments elsewhere raise the usual ``ValueError``, ``KeyError`` or ``TypeError``.)"""


class NotInvertibleError(DSQError, ZeroDivisionError):
    """No inverse was found: the element is zero, is known only to vanish below some degree, or is
    not a unit between the requested idempotents."""


class PrecisionError(NotInvertibleError):
    """Gauss-Jordan elimination ran out of precision: every remaining entry of a pivot column is
    only known to be O(deg p)."""


class NoSquareRootError(DSQError, ValueError):
    """A *-square-root rule cannot produce some y with y* y = x (for example because of the
    signature of a vertex, or because x has too large a rank at a middle vertex)."""


class BreakdownError(DSQError, ZeroDivisionError):
    """A matrix algorithm met a zero or non-invertible pivot, or an isotropic vector
    (v != 0 with v* v = 0)."""
\end{lstlisting}

\addtocontents{toc}{\protect\setcounter{tocdepth}{2}}

%======================================================================
%   BIBLIOGRAPHY
%======================================================================
% (the bibliography now precedes the appendices; see part09_appendix.tex)

%======================================================================
%   DOCUMENT END
%======================================================================
\end{document}